\PassOptionsToPackage{dvipsnames}{xcolor}
\documentclass[12pt]{amsart}   \textwidth=13.5cm \oddsidemargin=1cm \evensidemargin=1cm
\usepackage{amsmath}
\usepackage{amsxtra}
\usepackage{amscd}
\usepackage{amsthm}
\usepackage{amsfonts}
\usepackage{amssymb}\usepackage{eucal}
\usepackage{verbatim}
\usepackage[all]{xy}
\usepackage{tikz}
\usepackage[dvipsnames]{xcolor}
\definecolor{darkspringgreen}{rgb}{0.09, 0.45, 0.27}
\definecolor{Red}{rgb}{0.7,0,0}
\usepackage{stmaryrd}

\usepackage[bookmarks=false,colorlinks=true,linkcolor=black,urlcolor=darkspringgreen,citecolor=darkspringgreen,pdfnewwindow=true]{hyperref}
\usepackage{geometry}
\usepackage[utf8]{inputenc}
\usepackage[T1]{fontenc}
\usepackage[english]{babel}
\usepackage{microtype} \usepackage{tabularx} \usepackage{booktabs} \usepackage{amsmath, amsfonts,amssymb,amsthm,mathbbol, textcomp, yfonts}
\usepackage[all]{xy}
\usepackage{tikz,nicematrix}
\usetikzlibrary{decorations.pathreplacing}
\usepackage{appendix}
\usepackage{etoolbox}\usepackage[T1]{fontenc}
\usepackage{lmodern}
 \patchcmd{\footnotemark}{\stepcounter{footnote}}{\refstepcounter{footnote}}{}{}
\usepackage{blkarray}
\usepackage{multirow}
\usepackage{mathtools}
\usepackage{graphicx}
\usepackage{placeins}

\newtheorem{thm}{Theorem}[subsection]
\newtheorem{cor}[thm]{Corollary}

\newtheorem{lem}[thm]{Lemma}
\newtheorem{prop}[thm]{Proposition}
\newtheorem{conj}[thm]{Conjecture}

\newtheorem*{thmrestateinner}{\thmrestateheading}
\newcommand{\thmrestateheading}{}

\makeatletter
\newenvironment{thmrestate}[2][]{\renewcommand{\thmrestateheading}{Theorem~\ref{#2}$'$\if\relax\detokenize{#1}\relax \else\space(#1)\fi }\begin{thmrestateinner} \phantomsection \protected@edef\@currentlabel{\getrefnumber{#2}\noexpand\ensuremath{'}}}{\end{thmrestateinner} }
\makeatother
\usepackage{refcount}

\makeatletter

\newtheorem*{proprestateinner}{\proprestateheading}
\newcommand{\proprestateheading}{}

\newenvironment{proprestate}[2][]{\renewcommand{\proprestateheading}{Proposition~\getrefnumber{#2}$'$\if\relax\detokenize{#1}\relax \else\space(#1)\fi }\phantomsection
\edef\@currentlabel{\getrefnumber{#2}$'$}\begin{proprestateinner} }{\end{proprestateinner} }

\makeatother

\newtheorem{defn}[thm]{Definition}

\theoremstyle{definition}

\newtheorem{rem}[thm]{Remark}

\theoremstyle{definition}
\newtheorem{condition}[thm]{Condition}

\newcommand{\nc}{\newcommand}
\nc{\renc}{\renewcommand} \nc{\ssec}{\subsection}
\nc{\sssec}{\subsubsection}
\newcommand{\case}[1]{\smallskip\noindent\textbf{#1}}
\nc{\wh}{\widehat}
\nc\ol{\overline} \nc\ul{\underline} \nc\wt{\widetilde}

\nc{\BA}{{\mathbb{A}}} \nc{\BC}{{\mathbb{C}}} \nc{\BF}{{\mathbb{F}}}
\nc{\BD}{{\mathbb{D}}} \nc{\BG}{{\mathbb{G}}} \nc{\BQ}{{\mathbb{Q}}}
\nc{\BM}{{\mathbb{M}}} \nc{\BN}{{\mathbb{N}}} \nc{\BO}{{\mathbb{O}}}
\nc{\BP}{{\mathbb{P}}} \nc{\BR}{{\mathbb{R}}}
\nc{\BZ}{{\mathbb{Z}}} \nc{\BS}{{\mathbb{S}}} \nc{\BW}{{\mathbb{W}}}
\nc{\BL}{{\mathbb{L}}}

\nc{\CA}{{\mathcal{A}}} \nc{\CB}{{\mathcal{B}}} \nc{\CalD}{{\mathcal{D}}}
\nc{\CE}{{\mathcal{E}}} \nc{\CF}{{\mathcal{F}}}
\nc{\CG}{{\mathcal{G}}} \nc{\CH}{{\mathcal{H}}}
\nc{\CI}{{\mathcal{I}}} \nc{\CK}{{\mathcal{K}}} \nc{\CL}{{\mathcal{L}}}
\nc{\CM}{{\mathcal{M}}} \nc{\CN}{{\mathcal{N}}}
\nc{\CO}{{\mathcal{\bfO}}} \nc{\CP}{{\mathcal{P}}}
\nc{\CQ}{{\mathcal{Q}}} \nc{\CR}{{\mathcal{R}}}
\nc{\CS}{{\mathcal{S}}} \nc{\CT}{{\mathcal{T}}}
\nc{\CU}{{\mathcal{U}}} \nc{\CV}{{\mathcal{V}}}  \nc{\CY}{{\mathcal Y}}
\nc{\CW}{{\mathcal{W}}} \nc{\CZ}{{\mathcal{Z}}}

\nc{\cM}{{\check{\mathcal M}}{}} \nc{\csM}{{\check{\mathcal A}}{}}
\nc{\oM}{{\overset{\circ}{\mathcal M}}{}}
\nc{\obM}{{\overset{\circ}{\mathbf M}}{}}
\nc{\oCA}{{\overset{\circ}{\mathcal A}}{}}
\nc{\obA}{{\overset{\circ}{\mathbf A}}{}}
\nc{\ooM}{{\overset{\circ}{M}}{}}
\nc{\osM}{{\overset{\circ}{\mathsf M}}{}}
\nc{\vM}{{\overset{\bullet}{\mathcal M}}{}}
\nc{\nM}{{\underset{\bullet}{\mathcal M}}{}}
\nc{\oD}{{\overset{\circ}{\mathcal D}}{}}
\nc{\obD}{{\overset{\circ}{\mathbf D}}{}}
\nc{\oA}{{\overset{\circ}{\mathbb A}}{}}
\nc{\op}{{\overset{\bullet}{\mathbf p}}{}}
\nc{\cp}{{\overset{\circ}{\mathbf p}}{}}
\nc{\oU}{{\overset{\bullet}{\mathcal U}}{}}
\nc{\ofZ}{{\overset{\circ}{\mathfrak Z}}{}}

\nc{\ff}{{\mathfrak{f}}} \nc{\fv}{{\mathfrak{v}}}
\nc{\fa}{{\mathfrak{a}}} \nc{\fb}{{\mathfrak{b}}}
\nc{\fd}{{\mathfrak{d}}} \nc{\fe}{{\mathfrak{e}}}
\nc{\fg}{{\mathfrak{g}}} \nc{\fgl}{{\mathfrak{gl}}}
\nc{\fh}{{\mathfrak{h}}} \nc{\fri}{{\mathfrak{i}}}
\nc{\fj}{{\mathfrak{j}}} \nc{\fk}{{\mathfrak{k}}} \nc{\fl}{{\mathfrak{l}}}
\nc{\fm}{{\mathfrak{m}}} \nc{\fn}{{\mathfrak{n}}}
\nc{\ft}{{\mathfrak{t}}} \nc{\fu}{{\mathfrak{u}}}
\nc{\fw}{{\mathfrak{w}}} \nc{\fz}{{\mathfrak{z}}}
\nc{\fp}{{\mathfrak{p}}} \nc{\fq}{{\mathfrak{q}}} \nc{\frr}{{\mathfrak{r}}}
\nc{\fs}{{\mathfrak{s}}} \nc{\fsl}{{\mathfrak{sl}}}
\nc{\hsl}{{\widehat{\mathfrak{sl}}}}
\nc{\hgl}{{\widehat{\mathfrak{gl}}}}
\nc{\hg}{{\widehat{\mathfrak{g}}}}
\nc{\chg}{{\widehat{\mathfrak{g}}}{}^\vee}
\nc{\hn}{{\widehat{\mathfrak{n}}}}
\nc{\chn}{{\widehat{\mathfrak{n}}}{}^\vee}

\nc{\fA}{{\mathfrak{A}}} \nc{\fB}{{\mathfrak{B}}} \nc{\fC}{{\mathfrak{C}}}
\nc{\fD}{{\mathfrak{D}}} \nc{\fE}{{\mathfrak{E}}}
\nc{\fF}{{\mathfrak{F}}} \nc{\fG}{{\mathfrak{G}}} \nc{\fH}{{\mathfrak{H}}}
\nc{\fI}{{\mathfrak{I}}} \nc{\fJ}{{\mathfrak{J}}}
\nc{\fK}{{\mathfrak{K}}} \nc{\fL}{{\mathfrak{L}}}
\nc{\fM}{{\mathfrak{M}}} \nc{\fN}{{\mathfrak{N}}}
\nc{\frP}{{\mathfrak{P}}} \nc{\fQ}{{\mathfrak{Q}}}
\nc{\fS}{{\mathfrak{S}}} \nc{\fT}{{\mathfrak{T}}} \nc{\fU}{{\mathfrak{U}}}
\nc{\fV}{{\mathfrak{V}}} \nc{\fW}{{\mathfrak{W}}}
\nc{\fX}{{\mathfrak{X}}} \nc{\fY}{{\mathfrak{Y}}}
\nc{\fZ}{{\mathfrak{Z}}}

\nc{\ba}{{\mathbf{a}}}
\nc{\bb}{{\mathbf{b}}} \nc{\bc}{{\mathbf{c}}}
\nc{\be}{{\mathbf{e}}}
\nc{\bE}{{\mathbf{E}}}
\nc{\bj}{{\mathbf{j}}} \nc{\bm}{{\mathbf{m}}}
\nc{\bn}{{\mathbf{n}}} \nc{\bp}{{\mathbf{p}}}
\nc{\bq}{{\mathbf{q}}} \nc{\br}{{\mathbf{r}}} \nc{\bt}{{\mathbf{t}}}
\nc{\bfu}{{\mathbf{u}}} \nc{\bv}{{\mathbf{v}}}
\nc{\bx}{{\mathbf{x}}} \nc{\by}{{\mathbf{y}}} \nc{\bz}{{\mathbf{z}}}
\nc{\bw}{{\mathbf{w}}} \nc{\bA}{{\mathbf{A}}}
\nc{\bB}{{\mathbf{B}}} \nc{\bC}{{\mathbf{C}}}
\nc{\bD}{{\mathbf{D}}} \nc{\bF}{{\mathbf{F}}} \nc{\bG}{{\mathbf{G}}}
\nc{\bH}{{\mathbf{H}}} \nc{\bI}{{\mathbf{I}}} \nc{\bJ}{{\mathbf{J}}}
\nc{\bK}{{\mathbf{K}}} \nc{\bM}{{\mathbf{M}}} \nc{\bN}{{\mathbf{N}}}
\nc{\bO}{{\mathbf{\bfO}}} \nc{\bS}{{\mathbf{S}}} \nc{\bT}{{\mathbf{T}}}
\nc{\bU}{{\mathbf{U}}} \nc{\bV}{{\mathbf{V}}} \nc{\bW}{{\mathbf{W}}}
\nc{\bX}{{\mathbf{X}}}
\nc{\bY}{{\mathbf{Y}}} \nc{\bP}{{\mathbf{P}}}\nc{\bQ}{{\mathbf{Q}}}
\nc{\bR}{{\mathbf{R}}}
\nc{\bZ}{{\mathbf{Z}}} \nc{\bh}{{\mathbf{h}}}

\nc{\sA}{{\mathsf{A}}} \nc{\sB}{{\mathsf{B}}}
\nc{\sC}{{\mathsf{C}}} \nc{\sD}{{\mathsf{D}}}
\nc{\sE}{{\mathsf{E}}} \nc{\sF}{{\mathsf{F}}} \nc{\sG}{{\mathsf{G}}}
\nc{\sI}{{\mathsf{I}}} \nc{\sK}{{\mathsf{K}}} \nc{\sL}{{\mathsf{L}}}
\nc{\sfm}{{\mathsf{m}}} \nc{\sM}{{\mathsf{M}}} \nc{\sO}{{\mathsf{\bfO}}}
\nc{\sQ}{{\mathsf{Q}}} \nc{\sP}{{\mathsf{P}}}
\nc{\sT}{{\mathsf{T}}} \nc{\sZ}{{\mathsf{Z}}}
\nc{\sV}{{\mathsf{V}}} \nc{\sW}{{\mathsf{W}}}
\nc{\sfp}{{\mathsf{p}}} \nc{\sq}{{\mathsf{q}}} \nc{\sr}{{\mathsf{r}}}
\nc{\st}{{\mathsf{t}}} \nc{\sfb}{{\mathsf{b}}}
\nc{\sfc}{{\mathsf{c}}} \nc{\sd}{{\mathsf{d}}}
\nc{\sz}{{\mathsf{z}}}
\nc{\si}{{\mathsf{i}}}
\nc{\sj}{{\mathsf{j}}}
\nc{\sk}{{\mathsf{k}}}
\nc{\se}{{\mathsf{e}}}
\nc{\sm}{{\mathrm{sm}}}
\nc{\sinf}{{\frac{\infty}{2}}}

\nc{\tA}{{\widetilde{\mathbf{A}}}}
\nc{\tB}{{\widetilde{\mathcal{B}}}}
\nc{\tg}{{\widetilde{\mathfrak{g}}}} \nc{\tG}{{\widetilde{G}}}
\nc{\TM}{{\widetilde{\mathbb{M}}}{}}
\nc{\tO}{{\widetilde{\mathsf{\bfO}}}{}}
\nc{\tU}{{\widetilde{\mathfrak{U}}}{}} \nc{\TZ}{{\tilde{Z}}}
\nc{\tx}{{\tilde{x}}} \nc{\tbv}{{\tilde{\bv}}}
\nc{\tfP}{{\widetilde{\mathfrak{P}}}{}} \nc{\tz}{{\tilde{\zeta}}}
\nc{\tmu}{{\tilde{\xi}}}

\nc{\urho}{\underline{\rho}} \nc{\uB}{\underline{B}}
\nc{\uC}{{\underline{\mathbb{C}}}} \nc{\ui}{\underline{i}}
\nc{\uj}{\underline{j}} \nc{\ofP}{{\overline{\mathfrak{P}}}}
\nc{\oB}{{\overline{\mathcal{B}}}}
\nc{\og}{{\overline{\mathfrak{g}}}} \nc{\oI}{{\overline{I}}}

\nc{\eps}{\varepsilon} \nc{\hrho}{{\hat{\rho}}}
\nc{\blambda}{{\boldsymbol{\lambda}}} \nc{\bxi}{{\boldsymbol{\xi}}} \nc{\bnu}{{\boldsymbol{\eta}}}\nc{\bmu}{{\boldsymbol{\mu}}}

\usepackage{tikz}

\usepackage{adjustbox}
\newcommand{\dashbackslash}{\mathbin{\adjustbox{valign=m}{\tikz[scale=0.30,baseline]\draw[dashed,dash pattern=on 1pt off 1pt,line width=0.3pt] (0.5,0) -- (0,1);}}}
\newcommand{\dashslash}{\mathbin{\adjustbox{valign=m}{\tikz[scale=0.30,baseline]\draw[dashed,dash pattern=on 1pt off 1pt,line width=0.3pt] (0,0) -- (0.5,1);}}}

\nc{\one}{{\mathbf{1}}} \nc{\two}{{\mathbf{t}}}
\nc{\Cat}{\mathop{\operatorname{\rm 1-Cat}}}
\nc{\Sym}{\mathop{\operatorname{\rm Sym}}}
\nc{\Tot}{{\mathop{\operatorname{\rm Tot}}}}
\nc{\sprd}{{\mathop{\operatorname{\rm sprd}}}}
\nc{\Spec}{\mathop{\operatorname{\rm Spec}}}
\nc{\Ker}{{\mathop{\operatorname{\rm Ker}}}}
\nc{\Isom}{{\mathop{\operatorname{\rm Isom}}}}
\nc{\Hilb}{{\mathop{\operatorname{\rm Hilb}}}}
\nc{\Conf}{{\mathop{\operatorname{\rm Conf}}}}
\nc{\deeq}{{\mathop{\operatorname{\rm deeq}}}}
\nc{\FS}{{\mathop{\operatorname{\rm -FactMod}}}}
\nc{\End}{{\mathop{\operatorname{\rm End}}}}
\nc{\Ran}{{\mathop{\operatorname{\rm Ran}}}}
\nc{\Ext}{{\mathop{\operatorname{\rm Ext}}}}
\nc{\Hom}{{\mathop{\operatorname{\rm Hom}}}}
\nc{\CHom}{{\mathop{\operatorname{{\mathcal{H}}\it om}}}}
\nc{\GL}{{\mathop{\operatorname{\rm GL}}}}
\nc{\PGL}{{\mathop{\operatorname{\rm PGL}}}}
\nc{\PSp}{{\mathop{\operatorname{\rm PSp}}}}
\nc{\SOSp}{{\mathop{\operatorname{\rm SOSp}}}}
\nc{\SO}{{\mathop{\operatorname{\rm SO}}}}
\nc{\Sp}{{\mathop{\operatorname{\rm Sp}}}}
\nc{\unit}{{\mathop{\operatorname{\rm unit}}}}
\nc{\Mir}{{\mathop{\operatorname{\rm Mir}}}}
\nc{\St}{{\mathop{\operatorname{\rm St}}}}
\nc{\oblv}{{\mathop{\operatorname{\rm oblv}}}}
\nc{\gr}{{\mathop{\operatorname{\rm gr}}}}
\nc{\Id}{{\mathop{\operatorname{\rm Id}}}}
\nc{\perf}{{\mathop{\operatorname{\rm perf}}}}
\nc{\defi}{{\mathop{\operatorname{\rm def}}}}
\nc{\length}{{\mathop{\operatorname{\rm length}}}}
\nc{\supp}{{\mathop{\operatorname{\rm supp}}}}
\nc{\colim}{{\mathop{\operatorname{\rm colim}}}}
\nc{\Fun}{{\mathop{\operatorname{\rm Funct}}}}
\nc{\Hei}{{\mathop{\operatorname{\rm Heis}}}}
\nc{\Vac}{{\mathop{\operatorname{\rm Vac}}}}
\nc{\HC}{{\mathcal H}{\mathcal C}}
\nc{\ren}{{\mathsf{ren}}}
\nc{\locc}{{\mathsf{loc.c}}}
\nc{\pr}{{\operatorname{pr}}}
\nc{\Cliff}{{\mathsf{Cliff}}}
\nc{\loc}{{\operatorname{loc}}}
\nc{\Fl}{{\mathbf{Fl}}} \nc{\Ffl}{{\mathcal{F}\ell}}
\nc{\Fib}{{\mathsf{Fib}}}
\nc{\Coh}{{\mathsf{Coh}}} \nc{\FCoh}{{\mathsf{FCoh}}}
\nc{\Perf}{{\mathsf{Perf}}}

\nc{\reg}{{\text{\rm reg}}}
\nc{\gvee}{{\mathfrak g}^{\!\scriptscriptstyle\vee}}
\nc{\tvee}{{\mathfrak t}^{\!\scriptscriptstyle\vee}}
\nc{\nvee}{{\mathfrak n}^{\!\scriptscriptstyle\vee}}
\nc{\bvee}{{\mathfrak b}^{\!\scriptscriptstyle\vee}}
\nc{\rhovee}{\rh\bfO^{\!\scriptscriptstyle\vee}}

\nc{\cplus}{{\mathbf{C}_+}} \nc{\cminus}{{\mathbf{C}_-}}
\nc{\cthree}{{\mathbf{C}_*}} \nc{\Qbar}{{\bar{Q}}}

\nc{\Gtimes}{\vphantom{j^{X^2}}\smash{\overset{G}{\vphantom{\rule{0pt}{0.30em}}\smash{\times}}}}
\nc{\sGtimes}{\vphantom{j^{X^2}}\smash{\overset{\mathsf G}{\vphantom{\rule{0pt}{0.30em}}\smash{\times}}}}

\nc{\bOmega}{{\overline{\Omega}}}

\nc{\seq}[1]{\stackrel{#1}{\sim}}
\nc{\nilp}{{\operatorname{Nilp}}}
\nc{\Bun}{{\operatorname{Bun}}}
\nc{\mxd}{{\operatorname{mxd}}}
\nc{\Lus}{{\mathrm{Lus}}}
\nc{\free}{{\mathrm{free}}}
\nc{\cofr}{{\mathrm{cofr}}}
\nc{\DK}{{\mathrm{DK}}}
\nc{\funshf}{{la correspondance faisceaux-fonctions}}
\newcommand{\Weil}{\mathrm{Weil}}
\nc{\aff}{{\operatorname{aff}}}
\nc{\rot}{{\mathrm{rot}}}
\nc{\niso}{{\mathrm{N-Iso}}}
\nc{\iiso}{{\mathrm{Iso}}}
\nc{\Irr}{{\mathrm{Irr}}}
\nc{\Verma}{{\mathrm{Verma}}}
\nc{\res}{{\mathrm{res}}}
\nc{\Sph}{{\mathrm{Sph}}}
\nc{\Sat}{{\mathrm{Sat}}}
\nc{\super}{{\mathrm{super}}}
\nc{\pos}{{\mathrm{pos}}}
\nc{\fr}{{\mathrm{fr}}}
\nc{\locnilp}{{\mathrm{loc.nil}}}
\nc{\Dr}{{\mathrm{Dr}}}
\nc{\SVect}{{\mathrm{SVect}}}
\nc{\ratio}{{\mathrm{ratio}}}
\nc{\quant}{{\mathrm{quant}}}
\nc{\Kir}{{\mathrm{Kir}}}
\nc{\Cone}{{\mathrm{Cone}}}
\nc{\Supp}{{\mathrm{Supp}}}
\nc{\geo}{{\mathrm{geo}}}
\nc{\glob}{{\mathrm{glob}}}
\nc{\stan}{{\mathrm{st}}}
\nc{\Sch}{{\mathrm{Sch}}}
\nc{\Sect}{{\mathrm{Sect}}}
\nc{\disj}{{\mathrm{disj}}}
\nc{\gen}{{\operatorname{gen}}}
\nc{\fin}{{\operatorname{fin}}}
\nc{\mir}{{\operatorname{mir}}}
\nc{\triv}{{\operatorname{triv}}}
\nc{\ext}{{\operatorname{ext}}}
\nc{\righ}{{\operatorname{right}}}
\nc{\lef}{{\operatorname{left}}}
\nc{\forg}{{\operatorname{forg}}}
\nc{\fid}{{\operatorname{fd}}}
\nc{\modu}{{\operatorname{-mod}}}
\nc{\Mor}{{\operatorname{Mor}}}
\nc{\Gr}{{\mathbf{Gr}}}
\nc{\ord}{{\mathrm{ord}}}
\nc{\rel}{{\mathrm{rel}}}
\nc{\DrPl}{{\mathrm{DrPl}}}
\nc{\FT}{{\operatorname{FT}}}
\nc{\Mat}{{\operatorname{Mat}}}
\nc{\MSt}{{\operatorname{MSt}}}
\nc{\sph}{{\operatorname{sph}}}
\nc{\GR}{{\mathbf{Gr}}}
\nc{\Perv}{{\operatorname{Perv}}}
\nc{\Rep}{{\operatorname{Rep}}}
\nc{\Ind}{{\operatorname{Ind}}}
\nc{\coInd}{{\operatorname{coInd}}}
\nc{\IC}{{\operatorname{IC}}}
\nc{\Proj}{{\operatorname{Proj}}}
\nc{\colimc}{{\operatorname{colim}^{co}}}
\nc{\Stab}{{\operatorname{Stab}}}
\nc{\pt}{{\operatorname{pt}}}
\nc{\bfmu}{{\boldsymbol{\xi}}}
\nc{\bfomega}{{\boldsymbol{\omega}}}
\nc{\calM}{\mathcal M}
\nc{\calA}{\mathcal A}
\nc{\calO}{\mathcal O}
\nc{\cC}{\mathcal C}
\nc{\CC}{\mathbb C}
\nc{\calN}{\mathcal N}
\nc{\grg}{\mathfrak g}
\nc{\tslash}{/\!\!/\!\!/}
\nc\grt{\mathfrak t}
\nc\bfM{\mathbf M}
\nc\bfN{\mathbf N}

\nc\ZZ{\mathbb{Z}}
\nc\calC{\mathcal C}
\nc\calF{\mathcal F}
\nc\calX{\mathcal X}
\nc\calY{\mathcal Y}
\nc\seminf{\frac{\infty}{2}}
\nc\QCoh{\operatorname{QCoh}}
\nc\Spf{\operatorname{Spf}}
\nc\SD{\operatorname{Shv}}
\nc\SPerv{\operatorname{SPerv}}
\nc\FM{\operatorname{-FactMod}}
\nc\fact{\operatorname{fact}}
\nc\diag{\operatorname{diag}}
\nc\IndCoh{\operatorname{IndCoh}}
\nc\Maps{\operatorname{Maps}}
\nc\Dmod{D-\operatorname{mod}}
\newcommand\Hecke{\operatorname{Hecke}}
\newcommand\Hec{\mathrm{Hecke}}

\newcommand\DotHec{\overset{\bullet}{\mathrm{Hecke}}}
\nc{\calD}{\mathcal D}
\nc\bfO{\mathbf O}
\nc\bfF{\mathbf F}

\nc\GG{\mathbb G}
\nc\calK{\mathcal K}
\nc{\calG}{\mathcal G}
\nc\RHom{\operatorname{RHom}}
\nc\Res{\operatorname{Res}}
\nc\Av{\operatorname{Av}}
\nc\grs{\mathfrak s}
\nc{\tilX}{\widetilde X}
\nc\calB{\mathcal B}
\nc\calS{\mathcal S}
\nc\calT{\mathcal T}
\nc\calZ{\mathcal Z}
\nc\LS{\operatorname{LocSys}}
\nc\bfL{\on{\mathbf L}}

\makeatletter
\newcommand{\raisemath}[1]{\mathpalette{\raisem@th{#1}}}
\newcommand{\raisem@th}[3]{\raisebox{#1}{$#2#3$}}
\nc{\binlim}[2][]{\def\@tempa{#1}\@ifnextchar^{\@binlim{#2}}{\@binlim{#2}^{}}}
\def\@binlim#1^#2{\mathbin{\@ifempty{#2}{\mathop{#1}}{\mathop{#1}\@xp\displaylimits\@tempa^{#2}}}}

\makeatother

\nc\cX{{\mathcal X}}

\nc\Gm{{\mathbb G_m}}
\nc{\isoto}{\stackrel{\displaystyle\raisemath{-.2ex}{\sim}}\to}
\nc\cD\CalD
\nc\setm\setminus
\nc\cE\CE \renc\Hecke{\mathit{\CH\kern-.2ex ecke}}
\nc\bsl\backslash
\nc\Fq{\mathbb F_q}
\nc\x\times
\nc\bGO{{\bG_\bO}}
\nc\bla\blambda
\nc\blt\bullet
\nc\opp{{\on{op}}}
\nc\srel\stackrel
\nc\tbx{\binlim{\widetilde\boxtimes{}}}
\nc\chkbx{\binlim{\check\boxtimes}}
\nc\into\hookrightarrow
\nc\otto\leftrightarrow \renc\setminus\smallsetminus
\nc\phitau{\varphi\tau}
\newenvironment{i-ii-iii}{\begin{enumerate}}{\end{enumerate}}
\nc\ceil[1]{\lceil#1\rceil}  \nc\floor[1]{\lfloor#1\rfloor}
\nc\Lie{\on{Lie}}
 \let\arXiv\arxiv
\nc\on\operatorname
\nc\kap{\kappa}
\nc\gra{\mathfrak a}
\nc\gl{\mathfrak{gl}}
\nc\sTr{\operatorname{sTr}}
\nc\hatG{\widehat{G}}
\nc\calL{\mathcal L}
\nc\Whit{\operatorname{Whit}}
\nc\KL{\operatorname{KL}}

\makeatletter

\allowdisplaybreaks
\nc\mto{\mapsto }
\nc\en{\enspace }

\numberwithin{equation}{section}
\newtheorem*{rep@theorem}{\rep@title}
\newcommand{\newreptheorem}[2]{\newenvironment{rep#1}[1]{\def\rep@title{#2 \ref{##1}}\begin{rep@theorem}}{\end{rep@theorem}}}
\makeatother
\usepackage{tikz-cd}
\usepackage{wasysym, stackengine, makebox, graphicx}

\usepackage{xcolor,verbatim}
 \newcommand{\ncmd}{\newcommand*}
\newcommand{\rncmd}{\renewcommand*}
\ncmd{\DMO}{\DeclareMathOperator}
\ncmd{\ncmdd}[2]{\ncmd{#1}{{#2}}}

\makeatletter
\ncmd{\DefOps}[1]{\def\OPERATOR@NAME##1{\DeclareMathOperator{##1}{\expandafter\@gobble\string##1}}\def\OPERATOR@LIST##1{\ifcat\noexpand\relax\noexpand##1\OPERATOR@NAME##1\expandafter\OPERATOR@LIST\fi}\OPERATOR@LIST#1.}

\ncmd{\DefRm}[1]{\def\OPERATOR@NAME##1{\ncmd{##1}{\mathrm{\expandafter\@gobble\string##1}}}\def\OPERATOR@LIST##1{\ifcat\noexpand\relax\noexpand##1\OPERATOR@NAME##1\expandafter\OPERATOR@LIST\fi}\OPERATOR@LIST#1.}

\nc\negquad{\mkern-18mu}
\nc\lhs{\quad&\negquad}
\ncmd{\phtr}[2]{\lefteqn{#1{\phantom{#2}}}#2}

\def\Alphabet{ABCDEFGHIJKLMNOPQRSTUVWXYZ}
\def\newalph#1#2{\begingroup
\def\procL@tt@r##1{\@xp\gdef\csname#1\endcsname{#2}}\proc@lph@bet\endgroup}
\def\proc@lph@bet{\@xp\prlist@\Alphabet\relax}
\def\prlist@#1{\ifx#1\relax\else\procL@tt@r{#1}\@xp\prlist@\fi}

\ncmd{\SmSub}[2][]{_\bgroup #2\smsub@{#1}}
\def\smsub@#1{\@ifnextchar_{#1\@smsb}{\egroup}}
\def\@smsb_#1{#1\smsub@{}}

\ncmd{\SmSup}[2][]{^\bgroup #2\smsup@{#1}}
\def\smsup@#1{\@ifnextchar^{#1\@smsp}{\@ifnextchar'{\prime\@xp\smsup@\@xp{\@xp}\@gobble}{\egroup}}}
\def\@smsp^#1{#1\smsup@{}}

\def\@binlim#1_#2{\mathbin{\@ifempty{#2}{\mathop{#1}}{\mathop{#1}\@xp\displaylimits\@tempa_{#2}}}}

\catcode`\=\active
\ncmd{\RedNote}[1]{\Text{\textcolor{Red}{\sffamily#1}}}
\ncmd\Text[1]{\ifmmode\text{#1}\else#1\fi}
\ncmd\todo[1][todo]{\RedNote{#1}\ifx\undefined\@todoflag \global\let\@todoflag\relax \AtEndDocument{\par\vspace{3em}\noindent \RedNote{TO DO:}}\fi
\edef\pagenmbr{\thepage}\expandafter\redendnote \expandafter{\pagenmbr}{#1}}
\ncmd\redendnote[2]{\AtEndDocument{\\$\bullet$\ \RedNote{page #1: #2}}}
\def{\todo}
\ncmd\hidetodos{\rncmd\todo\relax}

\theoremstyle{remark} \newalph{c#1}{\ensuremath{{\mathcal #1}}} \newalph{s#1}{{\mathsf #1}} \newalph{#1#1}{{\mathbb #1}} \newalph{e#1}{{\EuScript #1}} \newalph{g#1}{{\mathfrak #1}} \newalph{r#1}{{\ensuremath{\mathscr #1}}}
   
\ncmd{\Fp}{{\FF_p}}

\DMO{\cHom}{\text{\textrm{\itshape{\cH}\kern-.2ex{}om}}}
\DMO{\cEnd}{\text{\textrm{\itshape{\cE}\kern-.2ex{}nd}}}    \DMO{\cExt}{\text{\textrm{\itshape{\cE}\kern-.2ex{}xt}}}

\ncmd{\sff}{\mathsf f}

\let\Im\undefined   \let\det\undefined
\ncmd{\young}[1]{\vcenter{\begin{Young}#1\crcr\end{Young}}}
\DefOps   { \REnd  \Aut \det \tr \rk \curv  \Coker \Im  \ad \Ad \add}
\DefRm{ \id \Vect \Pic \Loc \Hitch \Higgs \Eu \Fr}

\ncmd{\RGam}{\text{\upshape R}\Gamma}
\DMO{\chr}{char}
\ncmd{\angs}[1]{\langle#1\rangle}
\ncmd{\hmod}{\text{\upshape-mod}}
\ncmd{\cxym}[1]{\ensuremath{\vcenter{\xymatrix{#1}}}}
 \let\arXiv\arxiv
\makeatother

\title{Iwahori Gaiotto conjecture for cotangent type}
\dedicatory{To Michael Finkelberg, with admiration and gratitude.}
\date{}
\author[R.Travkin]{Roman Travkin}
\address{Skolkovo Institute of Science and Technology, Moscow, Russia}
\email{roman.travkin2012@gmail.com}
\author[R.Yang]{Ruotao Yang}
\address{Institute of Mathematics, Chinese Academy of Sciences, Beijing, 100190, China, and School of Mathematical Sciences, University of Chinese Academy of Sciences, Beijing, China}
\email{yruotao@gmail.com}

\newcommand{\fitmath}[1]{\begingroup \setbox0=\hbox{$\displaystyle #1$}\ifdim\wd0>.96\linewidth \resizebox{.96\linewidth}{!}{\box0}\else\box0\fi\endgroup}
\begin{document}

\begin{abstract}
We study the Iwahori Gaiotto conjecture for cotangent types and its application to the representation theory of the (mixed, small, and Lusztig) quantum supergroups. Using local and global models and Sakellaridis--Wang Zastava spaces, we construct factorization functors and compare the Iwahori Gaiotto category with the category $\cO$ of the mixed quantum supergroup. We then study the images of intersection cohomology sheaves and the metaplectic Hecke action. In type A, we describe the finite and affine monodromic categories, the corresponding Yokonuma--Hecke actions, and their role in the Kazhdan--Lusztig algorithm. The last section treats the pointwise purity.
\end{abstract}
\maketitle \setcounter{tocdepth}{2}
\tableofcontents

\setcounter{tocdepth}{2}

\section{Background and motivation}\label{sec:introduction} The Kazhdan--Lusztig theory, or the Kazhdan--Lusztig type geometric character formula, not only gives rise to a beautiful formula but also provides a deep connection that translates important representation theoretical information into geometric information. Usually, it translates irreducible objects into intermediate extension sheaves, and the characters are determined by the corresponding $*$-stalks or $!$-stalks, etc.

In the classical case, Beilinson--Bernstein and Brylinski--Kashiwara proved the Kazhdan--Lusztig conjecture by realizing the irreducible modules and Verma modules in the usual category $\cO$ as the corresponding objects on the flag varieties. For the quantum groups case, in particular, the root of unity case, there are also several well-developed geometric models, which are able to translate the representation information into geometric information, cf. \cite{[ABG],[ABBGM],[AJS],[KT],[CDR]}, etc.

However, in the super case, although there are already character formulas for the basic classical Lie superalgebras, cf. \cite{[Ser1],[Bru],[CLW],[CLW2],[BLW],[Bao]}, there is a lack of a uniform geometrization as in the non-super case, which translates the irreducible objects, standard objects, and the multiplicities into the corresponding intermediate extension sheaves, $!$-extension sheaves, and the $!$-stalks, etc.

Recently, motivated by the study of the boundary conditions in the supersymmetric Yang--Mills theory, D.~Gaiotto proposed a conjectural geometrization of the category of representations of the quantum supergroups associated with basic classical Lie superalgebras, which claims a braided monoidal equivalence between the category of representations of the quantum supergroup and the category of certain twisted sheaves on the affine Grassmannian, cf. \cite{[BFT1],[TY]}. We should remark here that it can be regarded as the quantum version of the relative Langlands program proposed by Ben-Zvi--Sakellaridis--Venkatesh, \cite{[BZSV]}, and also as the super generalization of the Fundamental Local Equivalence in the quantum Langlands program proposed by Gaitsgory--Lurie.

In the works \cite{[BFT],[BFT1],[BFT2]}, Braverman--Finkelberg--Travkin initiated the study of the Gaiotto conjecture, and proved the corresponding conjecture in the $q$ generic cases for the supergroup $\GL(N-1|N)$, and the same proof also works for $\GL(N|N)$, $\SOSp(2k|2k)$ and $\SOSp(2k|2k-2)$. Later, for $q$ generic, the cases for $\GL(M|N)$ and $\on{F}(4)$ were also proved in \cite{[TY]} and \cite{[FTY]}. However, a Kazhdan--Lusztig theory is still missing, in particular, in the root of unity case.

The goal of this paper is to study (the Iwahori extension of) the Gaiotto conjecture in the root of unity case, and establish the Kazhdan--Lusztig theory for the quantum supergroups. We will divide the basic classical Lie superalgebras into two families, cotangent type and non-cotangent type, depending on whether the corresponding geometric model can be polarized. In this paper, we work with the cotangent type case, and the non-cotangent type will be studied in future works joint with M.~Finkelberg.

\subsection{The Gaiotto conjecture}
\subsubsection{Spherical case} As we mentioned above, the Gaiotto conjecture is the quantum version of the relative Langlands conjecture. Roughly speaking, the latter claims that there is a bijection between certain coisotropic Hamiltonian varieties with certain properties, cf. \cite{[BZSV]}. In the anomaly-free case, the BZSV duality predicts that there is a canonical bijection between the hyperspherical varieties for $G$ and the hyperspherical varieties for the Langlands dual group ${G^\vee}$.
\begin{equation}
\begin{split}
\{G-\text{hyperspherical variety}\}&\Leftrightarrow \{{G^\vee}-\text{hyperspherical variety}\}\\
\bM&\leftrightsquigarrow {\bM^\vee}
\end{split}
\end{equation}
The spherical version of the relative Langlands conjecture claims that there is an equivalence between the category of $G(\bO)$-equivariant modules over an algebra $Q\bM$ arising via the deformation quantization of the structure sheaf $O_{\bM}$, and the category of ${G^\vee}$-equivariant quasi-coherent sheaves on ${\bM^\vee}$ up to a shearing, cf. \cite{[BZSV]}.

If ${\bM}$ is a cotangent bundle $T^*(G/H)$, or more generally, the $\psi$-twisted cotangent bundle $T_\psi^*(G/H)$ of a homogeneous spherical variety, then the algebra $Q\bM$ is given by the ring of $\psi$-twisted differential operators on $(G/H)(\bF)$, and the  conjecture says
\begin{conj}
\begin{equation}
D^{H(\bF),\psi}(\Gr_G)\simeq' \QCoh^{{G^\vee}}({(\bM^\vee)^\bullet}).
\end{equation}
\end{conj}
Here, the superscript $\bullet$ indicates the shearing, and $\simeq'$ means an equivalence up to a detail of renormalization of the categories, which we prefer not to discuss in detail here.

Given a basic classical supergroup $G_s$, it is known that the odd part $\fg_{\bar{1}}$ of $\fg_s$ is a coisotropic Hamiltonian variety of the even part, cf. \cite{[FU]}. Take ${\bM^\vee}$ to be $\fg_{\bar{1}}$ and ${G^\vee}=G_{\bar{0}}$, and up to a Koszul duality, \cite{[BFGT]}, the relative Langlands conjecture can be reformulated as
\begin{conj}\label{conj deg}
\begin{equation}
D^{H(\bF),\psi}(\Gr_G)\simeq' \Rep(\underline{G}_s).
\end{equation}
\end{conj}
Here, $\underline{G}_s$ denotes the corresponding degenerate supergroup. The data on the left-hand side of Conjecture \ref{conj deg} are classified in \cite{[FU]}.

The Gaiotto conjecture, in the cotangent type case, then says that
\begin{conj}\label{conj:spherical-gaiotto}
\begin{equation}\label{sph gai int}
D_\kappa^{H(\bF),\psi}(\Gr_G)\simeq' \Rep^\Lus_q({G}_s),
\end{equation}
where $D_\kappa^{H(\bF),\psi}(\Gr_G)$ denotes the category of $\kappa$-twisted $(H(\bF),\psi)$-equivariant $D$-modules on the affine Grassmannian, and $\Rep^\Lus_q({G}_s)$ denotes the category of representations of the Lusztig quantum supergroup. The quantum parameter and the twisting are related by $q=\exp(\pi\cdot i\cdot \kappa)$.
\end{conj}
An important caveat is that the full Lusztig quantum supergroup, with the precise properties needed here at a root of unity, is not constructed in this paper. Thus, the spherical statements involving it should be understood under the hypotheses in Remark~\ref{rem:spherical-gaiotto-conjecture}.

\subsubsection{Iwahori case} In the main body of the paper, we will study the Iwahori extension of Conjecture \ref{conj:spherical-gaiotto}. First, note that the relative Langlands conjecture admits a tamely ramified extension
\begin{conj}
\begin{equation}
\begin{split}
D^{H(\bF),\psi}(\Fl_G)\simeq' \QCoh^{{G^\vee}}({(\bM^\vee)^\bullet}\underset{{(\fg^\vee)^*}}{\times}{\tilde{\fg^\vee}^{*,\bullet}}),
\end{split}
\end{equation}
where ${\tilde{\fg^\vee}^{*}}$ denotes the Grothendieck--Springer resolution of ${(\fg^\vee)^*}$.
\end{conj}
We note that the spherical generated subcategory of the whole tamely ramified category is studied in \cite{[LPY]}.

Similar to Conjecture \ref{conj deg}, one can also rewrite the spectral side of the above conjecture in terms of the representation-theoretical category. Namely, there is
\begin{conj}
\begin{equation}
D^{H(\bF),\psi}(\Fl_G)\simeq' \cO^{\mxd}(\underline{U}(\fg_s)),
\end{equation}
where $\cO^{\mxd}(\underline{U}(\fg_s))$ denotes the category $\cO$ of a mixed Lie superalgebra, which can be defined as the relative Drinfeld center of the module category over the degenerate supergroup of the positive part with respect to the torus.
\end{conj}

In the quantum setting, the Iwahori extension of the Gaiotto conjecture says
\begin{conj}\label{iwa q}
\begin{equation}
D_\kappa^{H(\bF),\psi}(\Fl_G)\simeq' \cO_q^{\mxd}(\fg_s).
\end{equation}
Here, $\cO_q^{\mxd}(\fg_s)$ denotes the category $\cO$ of the mixed quantum supergroup, whose positive part is given by the Lusztig quantum supergroup, and the negative part is given by the De Concini--Kac quantum supergroup.
\end{conj}

An important difference with the spherical Gaiotto conjecture is that although the (whole) tamely ramified Gaiotto conjecture should be strictly stronger, the category on the representation side admits a definition without the use of the full Lusztig quantum supergroup. Thus, the above conjecture can be stated precisely. Namely, following Gaitsgory's definition of the mixed quantum group via the relative Drinfeld center, cf. \cite[Sections 4.4 and 5.3.1]{[G8]}, the right-hand side is exactly the relative Drinfeld center of the category of representations of the positive part of the Lusztig quantum supergroup with respect to the category of representations of the quantum torus, see \eqref{eq:mixed-quantum-supergroup-center}. Here, the positive part is defined as the graded Hopf dual of the negative part of the De Concini--Kac quantum supergroup, which has already been described explicitly, cf., \cite{[Y1]}.

\subsection{What do we do in the paper?} This paper is divided into five parts, with different subtopics on the same theme about describing the Kazhdan--Lusztig theory of the quantum supergroups at the root of unity.

\subsubsection{Equivalence of categories} In the first part, we generalize the main result in \cite{[Yan]} to the super setting and prove Conjecture \ref{iwa q} under a mild assumption on $q$, see Conditions~\ref{av s o} and~\ref{cond:star}.

First, we should make the statement of Conjecture \ref{iwa q} more precise by explaining the expected renormalizations. Let us note that the category $\cO$ of the mixed quantum supergroup admits a highest weight structure with the standard objects given by the Verma modules. Thus, one expects that the geometric side should contain the objects that correspond to the standard objects on the representation side, i.e., Verma modules; and it should be generated by the standard objects.

Motivated by \cite{[Yan]}, which in turn learned from D.~Gaitsgory and \cite{[AB]}, we construct the objects \[\Delta_\lambda:=\Av_!^{H(\bF),\psi}(J_\lambda)[-\delta(\lambda)],\] where $J_\lambda$ denotes the Wakimoto sheaf in \eqref{eq:twisted-bmw} indexed by the coweight $\lambda$, and $\delta$ is defined in \eqref{eq:local-degree-normalization}. Equivalently, one can also write $\Delta_\lambda=\Delta_0\star_{\rm H}J_\lambda$, where $\star_{\rm H}$ denotes the normalized action defined later. We denote by $\cC^\Delta$ the ind-complete category generated by the standards $\Delta_\lambda$, which is a full subcategory of $D_\kappa^{H(\bF),\psi}(\Fl_G)$. We note that it is quite different from the spherical generated part of $D_\kappa^{H(\bF),\psi}(\Fl_G)$ in \cite{[LPY]}; the latter does not have the Verma objects.

In the first part of the paper, we prove
\begin{thm}\label{thm 1.1.7}
For the homogeneous cases in Section~\ref{pol des}, assume that $q$ satisfies Condition~\ref{cond:star}.
\begin{enumerate}
\item There is a functor
\begin{equation}
F_\Fl^\Lus\colon D_\kappa^{H(\bF),\psi}(\Fl_G)\rightarrow  \cO_q^{\mxd}(\fg_s).
\end{equation}
\item The restriction of $F_\Fl^\Lus$ to $\cC^\Delta$ induces a $t$-exact equivalence of highest weight categories
\begin{equation}
F_\Fl^\Lus\colon \cC^\Delta\simeq \cO_q^{\mxd}(\fg_s),
\end{equation}
which sends the standard object $\Delta_\lambda$ to the Verma module $\Pi^{d(\lambda)}\Verma_\lambda$. Here, $\Pi$ denotes the parity shift, $d$ is defined in \eqref{eq:local-degree-normalization}, and we always take even highest weight lines for Verma modules and irreducible modules.
\item For any object $\cF\in D_\kappa^{H(\bF),\psi}(\Fl_G)$, there is a canonical isomorphism
\begin{equation}
\RHom_{D_\kappa^{H(\bF),\psi}(\Fl_G)}(\Delta_\lambda, \cF)\simeq \Pi^{d(\lambda)}\RHom_{\cO_q^{\mxd}(\fg_s)}(\Verma_\lambda, F_\Fl^\Lus(\cF)).
\end{equation}
\end{enumerate}
\end{thm}
A similar proof also works for $\GL(N|N)$, if we use the mirabolic model.

In fact, similar to \cite{[Yan]}, the proof uses the factorization method, and we prove a stronger version of the equivalence. Namely, one can consider general sheaf theory (e.g., de Rham version, $\ell$-adic version, Betti version, positive characteristic version, etc.) rather than just the $D$-modules. Theorem \ref{main} establishes a $t$-exact equivalence between the highest weight category $\cC^\Delta$ and the category of factorization modules over a certain factorization algebra. In the $D$-module setting, we use the argument in \cite{[G7]} to identify this factorization algebra with the factorization algebra arising from the quantum group side. In general, it is interesting to determine the Hopf (sub)-algebra and module category corresponding to the factorization algebra and the category of factorization modules in other sheaf settings.

In Section~\ref{sec:positive-characteristic}, we establish a bridge from the $D$-module setting to the mixed sheaf setting. Namely, we prove that the $!$-stalks of the intermediate extension in the $D$-module setting have the same dimension as the $!$-stalks of the intermediate extension in the $\ell$-adic setting when the characteristic is big enough. Thus, we want to use the mixed sheaves to calculate the canonical basis.

\subsubsection{Irreducible objects} In order to apply the result to the Kazhdan--Lusztig theory of the quantum supergroup at the root of unity, one has to determine the objects corresponding to the irreducible objects. Once this is done, then Theorem \ref{thm 1.1.7} (3) would give us a geometric description of the $\RHom$ of Vermas and irreducibles.

Naively, one may wonder if they are the intermediate extension sheaves on the relevant $H(\bF)$-orbits. However, this is not true: since the $t$-structure on the affine flags, which corresponds to the $t$-structure on the representation side, is the (stable) $t$-structure defined by the standard objects, the irreducible objects do not coincide with the irreducible objects in the perverse $t$-structure.

Motivated by \cite{[Lys],[GL]}, we consider the pullback to the affine flags of the intermediate extension sheaves on the affine Grassmannian indexed by \textit{restricted} dominant coweights. In the non-super/Whittaker setting, they correspond to irreducible representations of the graded small (or, Lusztig) quantum group.

In the super/Gaiotto setting, the situation is a little different: even in the generic case, \cite{[BFT]} noticed that, unlike the Whittaker case, not all orbits indexed by dominant coweights are relevant. Thus, the number of dominant coweights that are relevant and satisfy the naive restricted condition may be strictly smaller than the expected number of restricted weights on the representation side.

The reason is that in the super setting, the labelling of irreducible finite-dimensional representations by highest weights depends on the choice of the super Borel. The highest weight $\lambda$ may change by a simple odd super root $\alpha$ after the corresponding odd reflection, and the condition for whether it is preserved or changed by adding $\alpha$ depends on $\lambda,\alpha$, and $q$. Thus, we should consider the dominant coweights which correspond to the restricted dominant weights on the representation side with respect to the chosen super Borel, see Definition~\ref{def:restricted-coweight}. Later, when we talk about a restricted dominant weight, we mean such a coweight.

In Section~\ref{sec:irreducibles}, we prove that
\begin{thm}\label{thm 2}
Assume Condition~\ref{cond:star}. For $\lambda$ restricted dominant and relevant,
\begin{enumerate}
\item there are $t$-exact functors
\begin{equation}
\begin{split}
F_{\Gr}^{\Lus}\colon D_\kappa^{H(\bF),\psi,\ren}(\Gr_G)\rightarrow \cO_q^{\mxd}(\fg_s),\\
F_{\Gr}^{\sm}\colon D_\kappa^{H(\bF),\psi,\ren}(\Gr_G)\rightarrow \Rep^\sm_q(\fg_s).
\end{split}
\end{equation}
Here, $D_\kappa^{H(\bF),\psi,\ren}(\Gr_G)$ denotes the ind-completion of the locally compact Gaiotto category defined in Section~\ref{sec:global-models}, and its $t$-structure is the extension of the Gaiotto perverse $t$-structure, which is transported from the global model with the normalization \eqref{eq:local-degree-normalization} on each component. Also, $\Rep^\sm_q(\fg_s)$ denotes the category of representations of the graded small quantum supergroup, see Section~\ref{sec:mixed-quantum-supergroups}.
\item $F_{\Gr}^{\Lus}(\IC_{\Gr}^\lambda)$ is the irreducible representation $\Pi^{d(\lambda)}\Irr_\lambda^\mxd$;
\item $F_{\Gr}^{\sm}(\IC_{\Gr}^\lambda)$ is the irreducible representation $\Pi^{d(\lambda)}\Irr_\lambda^\sm$.
\end{enumerate}
\end{thm}
In particular, (2) implies that the pullback $p^!\IC_{\Gr}^\lambda$ along the projection $p\colon\Fl_G\rightarrow\Gr_G$ (with the convention in Remark~\ref{rem 3.5.6}) will go to $\Pi^{d(\lambda)}\Irr_\lambda^\mxd$ under $F^\Lus_{\Fl}$.

A conceptual explanation was learned from G.~Dhillon, assuming the expected Lusztig quantum supergroup and its properties in Remark~\ref{rem:spherical-gaiotto-conjecture}. The irreducible objects in the spherical Gaiotto category should be the intermediate extension sheaves on the orbits. The pullback functor from the spherical Gaiotto category to the Iwahori Gaiotto category should correspond to the forgetful functor (or the restriction functor, which means pullback along the algebra morphism from the mixed quantum supergroup to the Lusztig quantum supergroup) from the category of representations of the Lusztig quantum supergroup to the category $\cO$ of the mixed quantum supergroup. Note that if a weight is restricted, then the corresponding irreducible of the Lusztig quantum supergroup is also an irreducible of the small quantum supergroup, the De Concini--Kac quantum supergroup, as well as the mixed quantum supergroup, so the pullback of the intermediate extension on the affine Grassmannian is still irreducible. We will prove the above theorem directly in Section~\ref{sec:irreducibles}, without using this explanation/expectation.

It is expected that there is also a Steinberg--Lusztig isomorphism for the quantum supergroup. In Theorem~\ref{thm:metaplectic-hecke-action}, we prove the geometric side of the Steinberg--Lusztig theorem that
\begin{thm}
\begin{enumerate}
\item There is a $t$-exact normalized metaplectic Satake action $\star_{\rm H}$ on $D_\kappa^{H(\bF),\psi,\ren}(\Gr_G)$, as in \eqref{eq:normalized-Gaiotto-Hecke-action}.
\item For $\lambda$ restricted dominant relevant and $\lambda'\in (\Lambda^\sharp)^+$ such that $\lambda+\lambda'$ is relevant, we have
\begin{equation}
\IC_{\Gr}^{\lambda}\star_{\rm H}\IC_{\Gr}^{\Sat,\lambda'}\simeq \IC_{\Gr}^{\lambda+\lambda'}.
\end{equation}
\end{enumerate}
\end{thm}

The above theorem is a super analogue of the corresponding result in \cite{[GL]}. With the above theorem, and repeating almost verbatim the constructions and arguments as in \cite{[GL]}, one can prove the analogous theorem in \cite{[GL]} in the super setting, see Remark~\ref{rem:graded-small-quantum-group}. Namely, there is an equivalence
\begin{equation}
D_\kappa^{H(\bF),\psi,\ren}(\Gr_G)\underset{\Rep({G_{q}^\vee})}{\otimes}\Rep({T_{q}^\vee})\simeq \Rep_q^\sm(\fg_s).
\end{equation}
Here, ${G_{q}^\vee}$ is the metaplectic dual group associated with $q$, and ${T_{q}^\vee}$ is the corresponding torus.

Then, we extend the objects corresponding to irreducible modules from restricted dominant weights to other cases by convolving with the twisted Wakimoto sheaves. Namely, for the object $L_\mu$ defined before Corollary~\ref{cor:iwahori-irreducible-image}, we have an isomorphism
\begin{equation}
\RHom_{D_\kappa^{H(\bF),\psi}(\Fl_G)}(\Delta_\lambda, L_\mu)\simeq \Pi^{d(\lambda)+d(\mu)}\RHom_{\cO_q^{\mxd}(\fg_s)}(\Verma_\lambda, \Irr^\mxd_{\mu}).
\end{equation}

Actually, as we explain in Remark~\ref{rem:spherical-gaiotto-conjecture}, with the well-expected assumption of the conjectural Lusztig quantum supergroup, the above theorems also imply the comparison of the transition matrix coefficients of the canonical basis and the costandard basis on both sides in the spherical setting
\begin{equation}
[\RHom_{D_\kappa^{H(\bF),\psi}(\Gr_G)}(\Delta_{\lambda,\Gr}, \IC_{\Gr}^\mu)]= [\Pi^{d(\lambda)+d(\mu)}\RHom_{\Rep_q^\Lus(G_s)}(\on{Kac}_\lambda, \Irr^\Lus_\mu)],
\end{equation}
in the Grothendieck group.

The remaining part of the paper is to study how to get a Kazhdan--Lusztig algorithm to calculate the canonical basis from the geometric side.

\subsubsection{Finite type monodromic category} The third part is a study of the endoscopy theory of the Koszul duality proposed in \cite[Section 1.7]{[BFT3]}.

Unlike the non-monodromic case, where one can always perform the Kazhdan--Lusztig algorithm, a naive attempt to calculate the canonical basis of $D_\kappa^{H(\bF),\psi}(\Fl_G)$ fails immediately, even before entering the affine case and the special case $\GL(N|N)$.

Namely, in the non-monodromic case, \cite{[T1]} gives a collection of generators and then determines the canonical basis in the $B$-equivariant sheaves on the mirabolic flag variety. In the monodromic case, it is noticed in \cite{[And]} that the similarly defined generators in \cite{[T1]} cannot generate the category of $(B_N,\chi)$-equivariant sheaves on the mirabolic flag $(\GL_N\times \BA^N)/(B_N,\chi')$ with right twisting by $\chi'$. Thus, it seems very difficult to determine the duality of the Grothendieck group in the monodromic setting.

Since the problem already appears in the finite type case, let us temporarily only consider the finite type category, rather than consider the category of sheaves on the affine Grassmannian and the affine flags.

The key idea to save the game is the endoscopy theory of the relative Koszul duality in \cite{[BFT3]}. Namely, in loc.\ cit., the authors proposed \cite[Conjecture 1.7.2]{[BFT3]}, which claims that, in the case of general $\lambda$-monodromy (in the notation of loc.\ cit.),
\[
D(\mathcal Y^{T_\lambda})^{B_\lambda} \quad\text{and}\quad \bigl(Q_{\hbar=1}(\mathcal X^\vee)\text{-}\mathrm{mod}\bigr)^{B^\vee_{\mathrm{mon}},\lambda}
\]
should be Koszul dual to each other.

In other words, the category of $B_\lambda$-equivariant $D$-modules on the $T_\lambda$-fixed locus of $\mathcal Y$ is expected to be equivalent, via Koszul duality, to the category of $B^\vee$-monodromic modules over the quantization of the dual variety $\mathcal X^\vee$ with the monodromy $\lambda$.

In the $\GL(M|N)$ case, it says that
\begin{conj}
There is an equivalence
\begin{equation}
D((B_M,\chi)\ltimes (U_{M,N},\psi)\backslash \GL_N/(B_N,\chi'))\simeq' \bigotimes_i D(B_{k_i}\dashbackslash T^*\Hom(\BC^{k_i},\BC^{l_i})\dashslash B_{l_i})
\end{equation}
Here, $k_i$ and $l_i$ are the numbers of the components of the left and the right characters, respectively, that are equal to a given character, and $\sum k_i=M, \sum l_i=N$.
\end{conj}

For our use, we care more about the duality and the $t$-structure, and we expect that the conjecture above, after composing with the Koszul duality for each component, will be compatible with respect to the $t$-structure and the duality, i.e.,
\begin{conj}
There is a $t$-exact equivalence of mixed sheaves,
\begin{equation}\label{1.16 conj}
\cC^{\fin}_{\chi,\chi'}(M|N)\simeq \bigotimes_i \cC^{\fin}_{1,1}(k_i|l_i),
\end{equation}
which is compatible with respect to the Verdier duality. Here, $\cC^{\fin}_{\chi,\chi'}(M|N):=D((B_M,\chi)\ltimes (U_{M,N},\psi)\backslash\GL_N/(B_N,\chi'))$ if $M<N$, we use the mirabolic flag model if $M=N$, and $\cC^{\fin}_{\chi,\chi'}(M|N):=D((B_M,\chi)\backslash\GL_M/((B_N,\chi')\ltimes (U_{N,M},\psi)))$ if $M>N$.
\end{conj}

In Theorem~\ref{Thm 9.2.5}, we prove the above conjecture in the $\GL(M|N)$ case.

\subsubsection{Kazhdan--Lusztig calculations} With the preparations in the former part, we are able to go back to the affine setting, and consider the Kazhdan--Lusztig polynomials in the Iwahori Gaiotto category.

The fourth part is about the calculations. The calculations involve:
\begin{enumerate}
\item The explicit formula to translate the basis in \eqref{1.16 conj};
\item The action formula of Yokonuma's algebra in the finite type monodromic case;
\item The starting generators for the Kazhdan--Lusztig algorithm in the finite non-monodromic case;
\item The action formula of the affine Yokonuma's algebra in the affine monodromic case;
\item The generators for the Kazhdan--Lusztig algorithm in the affine monodromic case.
\end{enumerate}
These formulas determine the Kazhdan--Lusztig polynomials of the IC sheaves on the orbits in the twisted Iwahori Gaiotto category in the $\GL(M|N)$ case, for arbitrary $M,N$, twisting and monodromy, in the finite field setting of Section~\ref{sec:type-a-conventions}. In the case $M=N$, this extends the untwisted calculation in \cite[\S 2]{[FGT]}.

\subsubsection{Pointwise purity} In the fifth part, we deal with the pointwise purity property of the intermediate extension sheaf on the affine flags. In the case of $M=N$ and the non-twisted non-monodromic case, it is proposed in \cite{[FGT]} that the intermediate extension sheaves are pointwise pure, which could imply the parity and the positivity properties. We consider the case for general $M,N$ and the twisted monodromic case, and prove
\begin{thm}
Over a finite field, for any $M,N$ and any twisting and monodromy of finite order, the twisted $(H(\bF),\psi)$-equivariant intermediate extension sheaf on any relevant orbit is pointwise pure.
\end{thm}

\subsection{Notations}\label{sec:notation} Let $T$ and $T_s$ be the maximal tori of the geometric group $G$ and of the even part $G_{\bar0}$ of the supergroup $G_s$, respectively. The roots of $G_s$ (both even and odd) are called super roots, and the roots and coroots of $G$ are called geometric roots and geometric coroots. We will use the following notations.

We denote by $G=G_{\bar0}^\vee$ the group on the geometric side. Its coweight lattice $\Lambda=X_*(T)=X^*(T_s)$ is the weight lattice of the supergroup. Usually, $\lambda,\mu,\nu\in\Lambda$ denote coweights, $\check\lambda\in X^*(T)$ denotes a geometric weight, and $\alpha$ denotes a super root.

If $G=G_1\times G_2$, we denote the components of a coweight by $\lambda=(\lambda^{(1)},\lambda^{(2)})$ and $\lambda^{(a)}=(\lambda_i^{(a)})_i$, and we denote by $\alpha_{a,i}$ (resp. $\check\alpha_{a,i}$) a simple coroot (resp. the corresponding simple root) of $G_a$. Let $k$ be the ground field. We denote by $\bF=k((t))$ the field of Laurent series and by $\bO=k[[t]]$ the ring of formal power series.

We denote by $\se$ the coefficient field of sheaves. As in \cite[\S 1.6.1]{[TY]}, the notation $\SD$ means that we consider sheaves with coefficients in complexes of super vector spaces. (In the introduction, we suppressed it and simply wrote $D$.) Also, $\Pi$ denotes the parity shift, and $\RHom$ denotes the derived Hom. We always take even coefficients for the geometric IC sheaves, the standard and costandard objects, and the BMW sheaves, both locally and on the global model. For the standard, costandard and irreducible objects on the configuration side and on the representation side, we choose even generating lines (resp. even highest weight lines), and the additional parity shifts will be written explicitly as $\Pi^\epsilon$. Note that the uniqueness of irreducible objects is always understood up to parity shift, and the generating families include their parity shifts. For the type A calculations, we use coefficients in complexes of vector spaces, see Section~\ref{sec:type-a-conventions}.

\section{Configuration spaces and factorization modules}\label{sec:configuration-modules} It is known that the category of modules over a quantum group can be realized as the category of factorization modules over the corresponding factorization algebra, cf. \cite{[BFS],[GL],[CF],[L]}. In this section, we briefly recall the notion of factorization algebra and the category of factorization modules, following \cite[Sections 4--5]{[GL]}.

\subsection{Configuration spaces} Let $C$ be a smooth curve over $\BC$. Recall that $\Lambda=X_*(T)$ is the weight lattice of the even torus of the supergroup, and we denote by $\alpha_i$ the simple super roots of the chosen super Borel. We denote by $\Lambda^{\text{neg}}=\sum_i\BZ_{\geq0}(-\alpha_i)$ the monoid generated by the negative simple super roots. For $\alpha\in\Lambda^{\text{neg}}\setminus \{0\}$, we denote by $\Conf^\alpha$ the scheme which classifies $\Lambda^{\textnormal{neg}}$-colored divisors on $C$ of total degree $\alpha$. It is called the configuration space of degree $\alpha$. We denote by $\Conf$ the disjoint union of all $\Conf^\alpha$.

In particular, if $\alpha=-\sum_i n_i\alpha_i$, then $\Conf^\alpha\simeq\prod_i C^{(n_i)}$, where $C^{(n)}$ is the $n$-th symmetric power of $C$. On $\Lambda$, we denote by $\leq_{\super}$ the partial order given by the super roots, i.e., $\mu\leq_{\super}\lambda$ if $\mu-\lambda\in\Lambda^{\mathrm{neg}}$.

The configuration space acquires a non-unital monoid structure with respect to taking addition:
\begin{equation}
\begin{split}
\add: \Conf\times \Conf&\longrightarrow \Conf \\ D, D'&\mapsto D+D'.
\end{split}
\end{equation}
Let $(\Conf\times \Conf)_{\disj}$ be the open subscheme of $\Conf\times \Conf$, such that the supports of $D, D'$ are disjoint.
\subsection{Configuration spaces with a marked point} Let $c\in C$ be a fixed point and $\lambda\in \Lambda$. We denote by $\Conf^{\lambda}_{\infty \cdot c}$ the space which classifies $\Lambda$-colored divisors on $C$ such that the total degree is $\lambda$ and the coefficient of $x\neq c$ is an element of $\Lambda^{\textnormal{neg}}$.  Any point of $\Conf^{\lambda}_{\infty\cdot c}$ can be written as
\[D= \lambda_c\cdot c+ \sum \alpha_x \cdot x,\] where $x\in C-c$, and $\alpha_x\in \Lambda^{\textnormal{neg}}$.

We denote $\Conf_{\infty\cdot c}:=\bigsqcup_{\lambda\in\Lambda}\Conf^\lambda_{\infty\cdot c}$.

Note that taking addition defines a degree-graded map
\begin{equation}
\begin{split}
\add_c: \Conf\times \Conf_{\infty \cdot c}&\longrightarrow \Conf_{\infty \cdot c} \\ D, D'&\mapsto D+D'.
\end{split}
\end{equation}
Let $(\Conf\times \Conf_{\infty \cdot c})_{\disj}$ be the open subspace where the first divisor is disjoint from both $c$ and the support of the second divisor.

\subsection{Factorization algebras and modules} Let us fix a factorization gerbe $\cG^\Lambda$ on $\Conf$ and $\Conf_{\infty\cdot c}$, i.e., there are compatible isomorphisms
\[
\add^*(\cG^\Lambda)|_{(\Conf\times \Conf)_{\disj}} \simeq (\cG^\Lambda\boxtimes\cG^\Lambda)|_{(\Conf\times \Conf)_{\disj}},
\]
and
\[
\add_c^*(\cG^\Lambda)|_{(\Conf\times \Conf_{\infty\cdot c})_{\disj}} \simeq (\cG^\Lambda\boxtimes\cG^\Lambda)|_{(\Conf\times \Conf_{\infty\cdot c})_{\disj}},
\]
with higher homotopy coherence; see \cite[Section 4.3]{[GL]}. We denote the corresponding categories of twisted sheaves by $\SD_{\cG^\Lambda}(\Conf)$ and $\SD_{\cG^\Lambda}(\Conf_{\infty\cdot c})$.

For a $\cG^\Lambda$-twisted sheaf $\cI\in\SD_{\cG^\Lambda}(\Conf)$, we call it a factorization algebra if there is an isomorphism \[\add^!(\cI)|_{(\Conf\times \Conf)_{\disj}}\simeq \cI\boxtimes \cI|_{(\Conf\times \Conf)_{\disj}}{,}\] with higher homotopy coherence.

A factorization module over the factorization algebra $\cI$ is a twisted sheaf $\cF\in\SD_{\cG^\Lambda}(\Conf_{\infty \cdot c})$, such that \[\add_c^!(\cF)|_{(\Conf\times \Conf_{\infty \cdot c})_{\disj}}\simeq \cI\boxtimes \cF|_{(\Conf\times \Conf_{\infty \cdot c})_{\disj}},\] with higher homotopy coherence.

\subsection{The category of factorization modules}
\begin{defn}
Given a perverse factorization algebra $\cI$, assumed holonomic in the $D$-module setting, we denote by $\cI\FS$ the cocomplete category of $\cG^\Lambda$-twisted sheaves on $\Conf_{\infty\cdot c}$ which are factorizable with respect to $\cI$.
\end{defn}

\subsection{Standard and costandard objects} For $\lambda\in \Lambda$, we denote by $\Conf_{=\lambda\cdot c}$ the locally closed subscheme of $\Conf_{\infty\cdot c}$, such that a point belongs to $\Conf_{=\lambda\cdot c}$ if and only if the coefficient of $c$ is exactly $\lambda$. The following lemma is from \cite[Lemma 5.3.5]{[GL]}.

\begin{lem}
There is an equivalence of categories
\begin{equation}
\cI\FS|_{\Conf_{=\lambda\cdot c}}\simeq \SVect_{\lambda\cdot c}\simeq \SVect.
\end{equation}
Here, $\SVect_{\lambda\cdot c}$ is the category of complexes of super vector spaces twisted by the fiber of the gerbe $\cG^\Lambda$ on $\Conf_{\infty\cdot c}$ at $\lambda\cdot c$. It is (non-canonically) equivalent to the (non-twisted) category $\SVect$, by choosing a trivialization of this fiber. We choose the perverse generator corresponding to $\se$, placed in even parity and cohomological degree $0$.
\end{lem}

In the notation for configuration spaces, we write $\leq$ and $\geq$ for $\leq_{\super}$ and $\geq_{\super}$.

For $\lambda\in \Lambda$, we denote by $\Delta^{\fact}_{\lambda}$ (resp. $\nabla^{\fact}_{\lambda}$) the $!$-extension (resp. $*$-extension) of this perverse generator of $ \cI\FS|_{\Conf_{=\lambda\cdot c}}$. Since the embedding $\Conf_{=\lambda\cdot c}\hookrightarrow \Conf_{\infty\cdot c}$ factors as an affine open embedding into $\Conf_{\leq\lambda\cdot c}$ (where the coefficient of $c$ is at most $\lambda$ with respect to $\leq_{\super}$), followed by the closed embedding into $\Conf_{\infty\cdot c}$, the objects $\Delta^{\fact}_{\lambda}$ and $\nabla^{\fact}_{\lambda}$ are perverse, cf. \cite[Proposition 5.4.2]{[GL]}. Furthermore, the perverse heart of $\cI\FS$ has a highest weight structure, such that $\Delta^{\fact}_{\lambda}$ are standards and  $\nabla^{\fact}_{\lambda}$ are costandards. Namely, the objects $\{\Delta^{\fact}_{\lambda},\Pi\Delta^{\fact}_{\lambda}\}$ compactly generate $\cI\FS$, and there are orthogonality relations
\begin{equation}
\RHom_{\cI\FS}(\Delta^{\fact}_{\lambda}, \nabla^{\fact}_{\lambda})=\se,
\end{equation}
and
\begin{equation}
\RHom_{\cI\FS}(\Delta^{\fact}_{\lambda}, \nabla^{\fact}_{\mu})=0,
\end{equation}
if $\lambda\neq \mu$.

We denote the intermediate extension of the chosen generator by $\IC_\lambda^{\fact}$.

For the compact generation and the orthogonality relations, see \cite[Section 5.3.8]{[GL]}, and for the highest weight structure, see \cite[Section 5.4]{[GL]}.

Furthermore, as pointed out in \cite[Proposition 5.4.2]{[GL]}, the perverse $t$-structure of $\cI\FS$ can also be defined by
\begin{defn}\label{2.4.3}
For $\cM\in \cI\FS$, it is coconnective, if and only if
\begin{equation}
H^0\RHom_{\cI\FS}(\Pi^\epsilon\Delta^{\fact}_{\lambda}[k], \cM)=0,
\end{equation}
for any $k>0$, $\epsilon\in\{0,1\}$, and $\lambda\in \Lambda$.
\end{defn}

\section{Twisted Iwahori Gaiotto category}\label{Gaiotto category} In this section, we will introduce the constructible/geometric side of the main theorem in this paper; it will be (the quantum version of) a constructible sheaf category associated with the Gaiotto type/equivariant slice type coisotropic Hamiltonian varieties.

\subsection{Preparation about the associated spherical varieties} All the Gaiotto type coisotropic Hamiltonian varieties (the coisotropic Hamiltonian varieties which are $S$-dual to the ones given by basic classical supergroups $G_s$) have already been classified in \cite{[FU]}, and are presented explicitly by equivariant slices. In this paper, we restrict our attention to the $S$-dual pairs arising from $\GL(M|N)$, $\SOSp(2k|2l)$, and $\on{F}(4)$. In these cases, the corresponding coisotropic Hamiltonian varieties on the constructible side are polarized. In particular, they satisfy the anomaly free condition in \cite{[BZSV]}, and $\bM$ is a $G$-hyperspherical variety. Also, the group $G$ is the Langlands dual group of $G_{\bar{0}}:=(G_s)_{\bar{0}}$, and $G\circlearrowright\bM$ is the $S$-dual or the BZSV-dual of $G_{\bar{0}}\circlearrowright\fg_{\bar{1}}$.

For the cases $\SOSp(2k+1|2l)$, and $\on{G}(3)$, the corresponding coisotropic Hamiltonian varieties on the constructible side are of non-cotangent type, and one has to consider the metaplectic dual group and the twisted $S$-dual, cf., \cite{[BFT2]}; and the case $D(2,1,\alpha)$ is a deformation of the $\SOSp(4|2)$ case, cf., \cite[Remark 3.2.2]{[BFT]}.

\subsubsection{Description of the polarizations}\label{pol des} Since we use the polarization of the corresponding coisotropic varieties in an essential way, we should set specific conventions for each case, and in particular, we have to present the polarization data explicitly. In the homogeneous case, $\bM\simeq T_\psi^*(G/H)$, where $H$ is a spherical subgroup of $G$, and $\psi$ is an additive character of $H$. We denote by $\overset{\circ}{\bX}:=H\backslash G$ the open homogeneous space and by $\bX:=\overline{H\backslash G}^{\aff}$ its affine closure. Note that $G/H$ is identified with $H\backslash G$ via the inversion map $gH\mapsto Hg^{-1}$.

\begin{itemize}
\item $\GL(M|N)$ case: if $M<N$, then $G=\GL_M\times \GL_N$, and we take $H=\GL_M\ltimes U^-_{M,N}$. Here, $U^-_{M,N}$ is the unipotent radical of the lower-triangular parabolic subgroup of $\GL_N$ corresponding to the partition $(M+1, 1,\cdots,1)$. The group $H$ is regarded as a subgroup of $G$ via the following embedding: let $\GL_M$ embed into the upper-left corner of $\GL_N$, and $\GL_M\ltimes U^-_{M,N}$ embeds into $\GL_M\times \GL_N$ via the diagonal embedding
\begin{equation}\label{emd gl}
\begin{split}
H=\GL_M\ltimes U^-_{M,N}&\hookrightarrow G=\GL_M\times \GL_N\\
(g,u)&\mapsto (g, gu).
\end{split}
\end{equation}
The additive character $\psi\colon H\rightarrow \BG_a$ is given by the sum of the sub-diagonal entries in the $U^-_{M,N}$-part, i.e., $\sum_{i=M+2}^N u_{i,i-1}$. Then, the associated hyperspherical variety is $\bM=T_\psi^*(G/H)$; if $M=N$, the corresponding hyperspherical variety is $\bM=T^*(\GL_N\times \BA^N)$;
\item $\SOSp(2k|2l)$ case: let $(M,N)=(2k,2l+1)$ if $k\leq l$, and $(M,N)=(2l+1,2k)$ if $k>l$. Then, $m=(N-M-1)/2\geq0$, $G=\SO(M)\times \SO(N)$, and $H=\SO(M)\ltimes U^-_{M,N}$, where $U^-_{M,N}$ is the unipotent radical of the parabolic subgroup of $\SO(N)$ corresponding to the Levi subgroup $\SO(M+1)\times\BG_m^m$. For simplicity, we take the unipotent subgroup $U^-_{M,N}$ to be the lower triangular subgroup.

We realize $\SO(r)$ as the group of matrices $g$ satisfying $g^tJ_rg=J_r$ and $\det g=1$, where
\[
J_r=\begin{pmatrix}
0&\cdots&0&1\\
\vdots&&1&0\\
0&1&&\vdots\\
1&0&\cdots&0
\end{pmatrix}_{r\times r}.
\]

The group $H$ is regarded as a subgroup of $G$ via the following embedding: we choose an embedding of $\SO(M)$ into the subgroup $\SO(M+1)\subset \SO(N)$. Then, $\SO(M)\ltimes U^-_{M,N}$ embeds into $\SO(M)\times \SO(N)$ via the diagonal embedding as in \eqref{emd gl}. The embedding of the group $\SO(M)$ into $\SO(M+1)$ depends on the choice of an anisotropic vector of $\BC^{M+1}=\oplus_{i=1}^{M+1} \,\BC\cdot v_i$ with the symmetric bilinear form given by $(v_i, v_j)=\delta_{i, M+2-j}$, cf., \cite[\S 3.2]{[BFT]}. For example, if $M=2k$ is even (equivalently, $k\leq l$), then we can just choose the anisotropic vector $v_{k+1}$, and the embedding is given by $\begin{pmatrix} A& B\\
C&D
\end{pmatrix}\mapsto \begin{pmatrix}
A&0& B\\
0&1&0\\
C&0&D
\end{pmatrix}$; if $M=2l+1$ is odd (equivalently, $k>l$), then we take    $v_{\mathrm{an}}:=\frac{1}{\sqrt{2}}\left(v_{l+1}+v_{l+2}\right)$.
Then, $v_{\mathrm{an}}^\perp$ is identified with $\BC^{2l+1}$ by the ordered basis
\[
v_1,\ldots,v_l,\quad w:=\frac{i}{\sqrt{2}}\left(v_{l+1}-v_{l+2}\right),\quad v_{l+3},\ldots,v_{2l+2}.
\]
With respect to this basis, the embedding $\SO(2l+1)\hookrightarrow \SO(2l+2)$ is given by
\[
\begin{pmatrix}
A&\alpha&B\\
\beta&c&\gamma\\
C&\delta&D
\end{pmatrix}
\longmapsto
\begin{pmatrix}
A&-\frac{i}{\sqrt{2}}\alpha& \frac{i}{\sqrt{2}}\alpha&B\\
\frac{i}{\sqrt{2}}\beta& \frac{c+1}{2}&\frac{1-c}{2}& \frac{i}{\sqrt{2}}\gamma\\
-\frac{i}{\sqrt{2}}\beta& \frac{1-c}{2}&\frac{c+1}{2}& -\frac{i}{\sqrt{2}}\gamma\\
C&-\frac{i}{\sqrt{2}}\delta& \frac{i}{\sqrt{2}}\delta&D
\end{pmatrix}.
\]
For the character $\psi$, please see \cite[\S 3.2]{[BFT]};
\item $\on{F}(4)$ case: $G=\PGL_2\times \PSp(6)$, $H=\PGL_2\ltimes \on{U}^-$, where $\on{U}^-$ is the unipotent radical of the lower parabolic subgroup of $\PSp(6)$ corresponding to the Dynkin diagram without the middle vertex, i.e., corresponding to the partition $(2,2,2)$; we regard $\PGL_2$ as a subgroup of $\PSp(6)$ via the diagonal embedding, and the group $H$ is regarded as a subgroup of $G$ via the same formula as in \eqref{emd gl}. The additive character $\psi$ is given by $u_{3,1}+u_{4,2}$. For more details, please see \cite[\S 3.1]{[FTY]}.
\end{itemize}

In all the homogeneous cases above, $G=G_1\times G_2$, and $H$ is the semi-direct product of $H_r=G_1$ and a unipotent subgroup $H_u$. Furthermore, the embedding $H_r\into G$ factors as $H_r=G_1\overset\Delta\into G_1\x G_1\into G_1\x G_2=G$ where the first  map is the diagonal embedding and the second one is induced by an embedding $G_1\into G_2$.  Finally, $H_u$ is the unipotent radical of a parabolic subgroup $P^-=L\ltimes H_u$ of $G_2$ with a Levi subgroup $L$.

Unless otherwise specified, the statements below which are formulated in terms of $H$ are about these homogeneous polarizations. In the case $\GL(N|N)$, we will use the mirabolic model in the finite and affine type A sections later.

\begin{rem}
A more conceptual explanation is via the construction of the equivariant slices in \cite{[FU]}: we take a nilpotent element $e\in\Lie(G_2)$ for a reductive group $G_2$, namely,
\begin{enumerate}
\item in the $\GL(M|N)$ case, we take $e\in \mathfrak{gl}_N$ to be the nilpotent element of the Jordan type $(1^M,N-M)$;
\item in the $\SOSp(2k|2l)$ case, if $k\leq l$, then we take $e\in \mathfrak{so}_{2l+1}$ to be the nilpotent element of the Jordan type $(1^{2k},2l+1-2k)$; if $k>l$, then we take $e\in \mathfrak{so}_{2k}$ to be the nilpotent element of the Jordan type $(1^{2l+1},2k-2l-1)$;
\item in the $\on{F}(4)$ case, we take $e\in \mathfrak{psp}_6$ to be the nilpotent element of the Jordan type $(3,3)$.
\end{enumerate}
We consider the corresponding $\fsl_2$-triple $\{e,f,h\}$. If the grading associated with $h$ is even, then $(\fg_2)_{\leq 0}$ is a parabolic subalgebra which gives the polarization. Namely, $(\fg_2)_{0}$ is the Levi subalgebra, $(\fg_2)_{<0}$ is the nilpotent radical corresponding to $H_u$, and the character $\psi$ is given by pairing with $e$. Also, the group $G_1$ is chosen in the neutral component of the centralizer of the $\fsl_2$-triple.
\end{rem}

\subsubsection{Transversal embedding}\label{trans D} In this paper, we fix the choice of the Borel subgroups to be $B=B_{G_1}^-\times B_{G_2}$, where $B_{G_a}$ (resp. $B_{G_a}^-$) denotes the upper (resp. lower) triangular Borel subgroup of $G_a$. Similarly, let $U_a^-$ and $U_a^+$ be the unipotent radicals of $B_{G_a}^-$ and $B_{G_a}$, respectively, and let $T_a$ be the maximal torus of $G_a$. Here, $a\in\{1,2\}$. We also denote $B_a=B_a^+=B_{G_a}$ and $B_a^-=B_{G_a}^-$. We take $T=T_1\times T_2$, and denote $\Lambda_{G_a}=X_*(T_a)$, thus $\Lambda=\Lambda_{G_1}\oplus\Lambda_{G_2}$. In this paper, when we talk about dominance of coweights, we always mean dominance with respect to the \emph{opposite} Borel subgroup $B^-$ of $B$. We denote by $\Lambda_G^+$ (or simply $\Lambda^+$) the set of these dominant coweights.

The description in Section \ref{pol des} determines a Borel-compatible embedding $\iota\colon G_1\hookrightarrow L$, where the Borel of $L$ is given by the intersection with the (opposite) Borel of $G_2$.

However, from the perspective of the spherical variety, it is natural to consider an embedding $H\hookrightarrow G$ which is transversal to $B$. It is enough to find a representative $\sD\in G_2$ (viewed as $(1,\sD)\in G$), such that $H\sD B$ is open in $G$, or, the $B$-orbit passing through $\sD$ in the corresponding spherical variety is open. Then, the $\sD$-conjugated embedding
\[
\iota'\colon H=G_1\ltimes H_u\hookrightarrow G_1\times G_2,\qquad (g_1,u)\mapsto (g_1,\sD^{-1}\cdot g_1u\cdot \sD). \] of $H$ into $G$ is transversal to $B$. From now on, we denote by $H$ its image under this conjugated embedding, and by $\psi$ the transported character.

We take the representative $\sD$ as follows,
\begin{itemize}
\item In the $\GL(M|N)$ case, if $M<N$, we let
\[
\sD = \left(
\begin{array}{cccccccc}
1 &   &   &   &   &   &   &   \\
& 1 &   &   &   &   &   &   \\
&   & \ddots &   &   &   &   &   \\
&   &   & 1 &   &   &   &   \\
1 & 1 & \cdots & 1 & 1 &   &   &   \\
&   &   &   &   & 1 &   &   \\
&   &   &   &   &   & \ddots &   \\
&   &   &   &   &   &   & 1
\end{array}
\right) \in \GL_{N}
\]
be an $N\times N$-matrix whose $(M+1)$-th row is $(1,1,\cdots, 1, 0,\cdots, 0)$; if $M=N$, the spherical variety is $\GL_N\times \BA^N$, which is not homogeneous, and we let $\sD$ be the point $(\id,(1,1,\cdots, 1))\in \GL_N\times \BA^N$; then we still have that the $B$-orbit passing through $\sD$ is open.

\item In the $\SOSp(2k|2l)$ case, we have a natural embedding $\SO(M+1)\hookrightarrow \SO(N)$, and we can take an element $\sD$ in $\SO(M+1)$. A convenient way of choosing $\sD$ depends on the parity of $M$. When $M=2k$ is even, we can take $\sD$ to be the exponential of the following element in $\mathfrak{so}(M+1)=\mathfrak{so}(2k+1)$,
\[
\sd = \left(
\begin{array}{cccccccc}
0 &   &   &   &   &   &   &   \\
& 0 &   &   &   &   &   &   \\
&   & \ddots &   &   &   &   &   \\
&   &   & 0 &   &   &   &   \\
1 & 1 & \cdots & 1 & 0 &   &   &   \\
&   &   &   &  -1 & 0 &   &   \\
&   &   &   &  \vdots &   & \ddots &   \\
&   &   &   &  -1 &   &   & 0
\end{array}
\right) \in \mathfrak{so}(M+1).
\]
If $M$ is odd, then we can take $\sD=\exp(\sd)$, with $\sd$ given by
\[
\sd= \left(
\begin{array}{cccccccc}
0      &        &        &        &        &        &        &   \\
& \ddots &        &        &        &        &        &   \\
&        & 0      &        &        &        &        &   \\
1      & \cdots & 1      & 0      & 0      &        &        &   \\
1      & \cdots & 1      & 0      & 0      &        &        &   \\
&        &        & -1     & -1     & 0      &        &   \\
&        &        & \vdots & \vdots &        & \ddots &   \\
&        &        & -1     & -1     &        &        & 0
\end{array}
\right) \in \mathfrak{so}(M+1).
\]
\item In the $\on{F}(4)$ case, we take $\sD=1+E_{4,3}\in \PSp(6)$.
\end{itemize}

All the homogeneous cases that we consider here can be written as $\bM=T^*_\psi (G/H)$, see \cite[Section 3]{[BFT2]}, and the restriction of the characters $\psi$ to $H_u(\bF)$ is the most non-degenerate additive character. Namely, it is the $t^{-1}$-coefficient of the sum of the sub-diagonal elements in the $\GL(M|N)$ case, the $t^{-1}$-coefficient of $u_{3,1}+u_{3,2}+u_{4,2}$ in the $\on{F}(4)$ case (after the conjugation by $\sD$), and see \cite[Section 3.2]{[BFT]} for the $\SOSp(2k|2l)$ case. The characters $\psi$ are invariant under the $H_r(\bF)$-action on $H_u(\bF)$, and extend to $H(\bF)$. With some abuse of notation, we denote by $\psi\colon H(\bF)\rightarrow \BG_a$ the additive character on $H(\bF)$ as well as its restriction to its subgroups.

\medskip
\subsection{Definition of the Gaiotto category}
\subsubsection{\texorpdfstring{$\omega^\rho$}{omega-rho}-renormalized spaces and groups} For $a=1,2$, let $I_a\subset G_a(\bO)$ be the preimage of $B\cap G_a$ under the evaluation map at $t=0$, and denote $I=I_1\times I_2$.

The additive character $\psi$ of $H$ determines an additive character of $H(\bF)$. Morally, the Iwahori Gaiotto category is defined as the category of $(H(\bF),\psi)$-equivariant twisted sheaves on the affine flags $\Fl_G=\Fl_{G_1}\times \Fl_{G_2}$. Using the restriction functor to $1\in \Fl_{G_1}$, we can also define it as the category of $I_1\ltimes H_u(\bF)$-equivariant $\kappa$-twisted sheaves on $\Fl_{G_2}$ against the character $\psi$.

However, similar to \cite{[G3]} (for an explicit description in the Gaiotto cases, one can also see \cite[\S 3]{[TY]} and \cite{[FGT]}), we consider the $\omega^\rho$-renormalized affine flags $\Fl_G^\omega$ and renormalized group $H^{\omega}(\bF)$ to have a canonically defined $\psi$ (rather than depending on the choice of a coordinate $t$).

\subsubsection{The grading cocharacter}\label{sec:grading-cocharacter} Denote by $Z_L$ the center of $L$. For example, in the $\GL(M|N)$ case, one has $Z_L\simeq\BG_m^{N-M}$.

We take
\[
\rho=(0,\rho^{(2)}),
\]
where $\rho^{(2)}$ is the difference of the Weyl cocharacters associated with $L$ and $G_2$, up to adding a half-integral cocharacter of $Z_L$. In particular, $\rho$ is only non-trivial on the second factor, and $2\rho^{(2)}$ factors through $Z_L$:
\[
2\rho:\BG_m\to T, \qquad 2\rho^{(2)}:\BG_m\to Z_L\subset T_2.
\]
We choose this central factor such that
\[
\langle\check\alpha,\rho\rangle=1
\]
for any geometric root $\check\alpha$ on whose root subgroup $\psi$ is non-zero.

Note that in the cotangent case, $G$ is the Langlands dual group of the even part of $G_s$, and in Section~\ref{super Borel}, we will find out that $\rho$ corresponds, up to adding a central cocharacter, to the Weyl vector of the supergroup $G_s$.

\subsubsection{Twisted spaces} We take a smooth curve $C$, and fix a square root $\omega^{\otimes \frac12}$ of the canonical line bundle $\omega=\omega_C$. We define the associated $T$-bundle via
\[
\omega^\rho:=(\omega^{\otimes \frac12})^\times\times^{\BG_m,2\rho} T.
\]
Here, $(\omega^{\otimes \frac12})^\times$ denotes the associated $\BG_m$-torsor. Since the first component of $\rho$ is trivial, the induced $T_1$-bundle is trivial, and the twisting occurs only on the $G_2$-factor. Similarly, we denote
\[
\cP_G^\omega:=\omega^{\rho}\times^{T} G=\cP_{G_1}^0\times \cP_{G_2}^\omega, \] and similarly for $\cP_{P^-}^\omega$ and $\cP_L^\omega$.

We denote by $\Ran$ the Ran space of $C$ and by $\Ran_c$ the Ran space with the fixed marked point $c$; see Section~\ref{ranified version}.

The $\omega^\rho$-twisted Beilinson--Drinfeld affine Grassmannian $\Gr_{G,\Ran}^{\omega}$ is defined, as usual, with the only modification of replacing the trivial $G$-bundle with $\cP_G^\omega$. Again, since only the second factor of $\rho$ is non-trivial, one has
\[
\Gr_{G,\Ran}^{\omega} \simeq \Gr_{G_1,\Ran}\times_{\Ran}\Gr_{G_2,\Ran}^{\omega}.
\]
Similarly, one defines the $\omega^\rho$-twisted Beilinson--Drinfeld affine flags $\Fl_{G,\Ran_c}^{\omega}$, $\omega^\rho$-twisted affine Grassmannian $\Gr_{G,c}^{\omega}$ and $\omega^\rho$-twisted affine flags $\Fl_{G,c}^{\omega}$.

Now, let us fix a point $c\in C$, and denote by $\cD_c=\Spec \bO_c, \overset{\circ}{\cD}_c=\Spec \bF_c,$ the formal disc and the punctured disc. Then, we define the $\omega^\rho$-twisted loop group and jet group via
\[
G(\bF)_c^{\omega}:=\Aut(\cP_G^\omega|_{\overset{\circ}{\cD}_c}), \qquad G(\bO)_c^{\omega}:=\Aut(\cP_G^\omega|_{\cD_c}).
\]
Note that there is
\[
G(\bF)_c^{\omega}\simeq G_1(\bF)_c\times G_2(\bF)_c^{\omega}, \qquad G(\bO)_c^{\omega}\simeq G_1(\bO)_c\times G_2(\bO)_c^{\omega}.
\]

Similarly, we have $I^{\omega}:=I_1\times I_2^{\omega}$ and $H(\bF)_c^{\omega}=G_1(\bF)_c\ltimes H_u(\bF)_c^{\omega}$. When the point $c$ is fixed, we omit it from the notation, i.e., we denote $H(\bF)_c^\omega$ and $H_u(\bF)_c^\omega$ by $H^\omega(\bF)$ and $H_u^\omega(\bF)$, respectively.

\subsubsection{Gerbes used in the paper}\label{sec:factorization-gerbes} In this section, we briefly summarize the gerbes used in this paper, cf., \cite{[GL],[GL1],[Zha1],[Zha2],[Zha3]}.

In the Betti and $\ell$-adic settings, let $\cG$ be a factorization $\mathsf e^{\mathsf{torsion},\times}$-gerbe on $\Gr_{G_2,\Ran}$, cf., \cite{[GL1]}. Here, $\mathsf e^{\mathsf{torsion},\times}$ denotes the multiplicative group of torsion elements in $\se$; in the $D$-module setting, we assume that $\cG$ is tame, as in \cite[Section 3]{[Zha1]}.

The $G_2(\bO)$-equivariance of $\cG$ yields a descent gerbe on the Hecke prestack. Taking the pullback along the projection from $\Gr_{G_2,\Ran}^{\omega}$ to the Hecke prestack of $G_2$ gives a factorization gerbe on $\Gr_{G_2,\Ran}^{\omega}$, which we still denote by $\cG$. It is $G_2(\bO)^{\omega}$ (in particular, $G_1(\bO)$)-equivariant, and thus gives rise to a $G_1(\bF)$-equivariant factorization gerbe on $\Gr_{G,\Ran}^{\omega}$, which we denote by $\cG^G$. Its restriction to the first factor is the \textit{inverse} of the naive restriction along $\Gr_{G_1,\Ran}\rightarrow \Gr_{G_2,\Ran}^{\omega}$. By abuse of notation, its pullbacks to
\[
\Fl_{G,\Ran_c}^{\omega},\qquad \Fl_{G,c}^{\omega},\qquad \Gr_{G,c}^{\omega},\qquad G(\bF)_c^{\omega},\qquad H(\bF)_c^{\omega}
\]
will all be denoted by the same symbol.

We note here that the pullback of $\cG^G$ to $G(\bF)_c^{\omega}$ is multiplicative, i.e.,
\[
m^*(\cG^G)\simeq \cG^G\boxtimes \cG^G,
\]
cf., \cite[\S 7.2]{[GL1]} and \cite[Corollary 3.4.7]{[Zha3]}. Hence, its restriction to $H(\bF)_c^\omega$ is also multiplicative. Furthermore, because of the inverse in the definition of $\cG^G$ above, its restriction to $H(\bF)_c^\omega$ admits a multiplicative trivialization. Thus, it makes sense to define $(H(\bF)_c^\omega,\psi)$-equivariant sheaves.

\subsubsection{Gerbe on $\Conf$}Replacing $G$ by $B$ in the definition of $\Gr_{G, \Ran}^{\omega}$, one obtains a factorization prestack $\Gr_{B, \Ran}^{\omega}$. Consider the following diagram of prestacks:
\begin{center}
\begin{equation}
\begin{tikzpicture}[scale=1]
\draw [->] (-1,1) -- (-2.5,0); \draw [->] (1,1) -- (2.5,0); \node at (0,1.2){$\Gr_{B, \Ran}^{\omega}$}; \node [below] at (2.5,0){$\Gr_{T, \Ran}^{\omega}.$}; \node [below] at (-2.5,0){$\Gr_{G, \Ran}^{\omega}$};
\end{tikzpicture}
\end{equation}
\end{center}
The pullback of $\mathcal{G}^G$ on $\Gr_{G,\Ran}^{\omega}$ along the left morphism gives a factorization gerbe on $\Gr_{B, \Ran}^{\omega}$. By \cite[Section 5.1]{[GL1]} and \cite[Lemma 4.3 and Section 4.5]{[Zha1]}, this factorization gerbe descends to a factorization gerbe on $\Gr_{T, \Ran}^{\omega}$. We denote the resulting gerbe by $\mathcal{G}^T$.

We let $(\Gr_{T,\Ran}^{\omega})^{\textnormal{neg}}$ (resp. $(\Gr_{T, \Ran_c}^{\omega})_{\infty\cdot c}^{\textnormal{neg}}$) be the closed sub-prestacks of $\Gr_{T,\Ran}^{\omega}$ (resp. $\Gr_{T,\Ran_c}^{\omega}$) satisfying the regular and non-redundant conditions defined in \cite[\S 4.6.2, 4.6.6]{[GL]}, using the negative simple super roots. By construction, we have maps, cf., \cite[(4.16),(4.18)]{[GL]},
\begin{equation}\label{omega to conf}
(\Gr_{T, \Ran}^{\omega})^{\textnormal{neg}}\to \Conf,
\end{equation}
and
\begin{equation}\label{omega to confx}
(\Gr_{T, \Ran_c}^{\omega})_{\infty\cdot c}^{\textnormal{neg}}\to \Conf_{\infty\cdot c}.
\end{equation}
According to \cite[Lemma 4.6.4]{[GL]}, (\ref{omega to conf}) and (\ref{omega to confx}) induce equivalences of gerbes on $(\Gr_{T,\Ran}^{\omega})^{\textnormal{neg}}$ (resp. $(\Gr_{T, \Ran_c}^{\omega})_{\infty\cdot c}^{\textnormal{neg}}$) and the corresponding configuration spaces.

Hence, one can descend the factorization gerbe $\mathcal{G}^T$ on $(\Gr_{T, \Ran}^{\omega})^{\textnormal{neg}}$ (resp. $(\Gr_{T, \Ran_c}^{\omega})_{\infty\cdot c}^{\textnormal{neg}}$) to a gerbe $\mathcal{G}^\Lambda$ on $\Conf$ (resp. $\Conf_{\infty\cdot c}$). Note that the above maps are compatible with the factorization structures, and thus, the gerbes on $\Conf$ and $\Conf_{\infty\cdot c}$ are factorizable.

\subsubsection{The quadratic form $q$ and the bilinear form $b$}\label{3.2.5 sect} Following \cite[Section 4.2]{[GL1]} and \cite[Section 4.5.3 and Lemma 5.7]{[Zha1]}, we can get a quadratic form
\begin{equation}\label{quad f}
q\colon \Lambda\to \se/\mathbb{Z}
\end{equation}
from a factorization gerbe $\mathcal{G}^\Lambda$ on $\Conf$ in the $D$-module setting. In the Betti setting and the $\ell$-adic setting, $q$ takes values in $\mathsf{e}^{\mathsf{torsion},\times}(-1):= \colim_n \Hom(\mu_n, \mathsf{e}^{\mathsf{torsion},\times})$. Also, we have a bilinear form; in the $D$-module setting, it is given by
\begin{equation}
b(\lambda,\mu)=q(\lambda+\mu)-q(\lambda)-q(\mu),
\end{equation}
and in the Betti and $\ell$-adic settings, it is replaced by
\begin{equation}
b(\lambda,\mu)=q(\lambda+\mu)q^{-1}(\lambda)q^{-1}(\mu).
\end{equation}

\subsubsection{Gerbe on $\Bun_G$} When $C$ is projective, let $\Bun_G$ denote the algebraic stack classifying principal $G$-bundles on the curve $C$. By \cite[Section 2.3.5]{[GL1]} and \cite{[Zha1]}, the gerbe $\cG^G$ on $\Gr_{G,\Ran}^{\omega}$ descends to $\Bun_G$, normalized at the point $\cP_G^\omega$. We will also denote the resulting gerbe on $\Bun_G$ by $\mathcal{G}^G$. Similarly, $\cG^T$ descends to $\Bun_T$, normalized at the point $\omega^\rho$.

\subsubsection{Factorization line bundles}
\begin{rem}\label{rem:determinant-twisting}
For the descriptions below, we assume that $C$ is projective. For the twisting used in this paper (to define those twisted sheaves), one can also think of factorization line bundles $\cP_{\det}$ rather than factorization gerbes, and will not lose the main interest of the paper. Namely,
\begin{enumerate}
\item In the $\GL(M|N)$ case, we consider the determinant line bundle $\cP_{\det,N}$ on $\Fl_{G_2}$, thus, the corresponding bilinear form is given by the trace. Then, after extending to $\Fl_G$, and then taking the $\omega$-renormalization (tensoring with the inverse of the fiber at $(\cP_{G_1}^0, \cP_{G_2}^\omega)$), we get the following relative determinant line bundle over $\Fl^\omega_G$, whose fiber over the point $(\cP_M, \alpha_M, \epsilon_M, \cP_N, \alpha_N, \epsilon_N)$ is
\begin{equation*}
\det^{-1} R\Gamma(C, \cV^M_{\cP_M})\otimes \det R\Gamma(C, \cV^N_{\cP_N})\otimes \det R\Gamma(C, \cV^M_{\cP^{triv}_M})\otimes \det^{-1} R\Gamma(C, \cV^N_{\cP_N^\omega}),
\end{equation*}
where $\cV^M_{\cP_M}$ (resp.\ $\cV^N_{\cP_N}$) denotes the rank $M$ (resp.\ $N$) vector bundle associated with $\cP_M$ (resp.\ $\cP_N$).
\item In the $\SOSp(2k|2l)$ case, like \cite[Section 3.2]{[BFT]}, we can take the determinant line bundle $\cP_{\det}$ on $\Fl_{\SO(N)}$ and the corresponding bilinear form is half of the trace. By the process of extension and $\omega$-renormalization, it gives rise to a relative determinant line bundle over $\Fl^\omega_G$;
\item In the $\on{F}(4)$ case, we can take the renormalized determinant line bundle from \cite[Section 3.1]{[FTY]}.
\end{enumerate}
The reason that we use the terminology of gerbes is just to match the conventions in \cite{[GL],[GL1],[Zha1],[Yan]}, etc.

These determinant line bundles give rise to factorization gerbes, cf., \cite[Section 1.4]{[GL1]} and, for tame parameters, \cite[Section 3.3.4]{[Zha1]}. Then, $\SD_{\cG^G}(\Fl_G^\omega)$ will be identified with the category of sheaves on $\overset{\circ}{\cP}_{\det}$ which are monodromic along the fibers with the given monodromy. Here, $\overset{\circ}{\cP}_{\det}$ denotes the punctured total space of the above determinant line bundle on $\Fl_G^\omega$. 
\end{rem}

\begin{rem}\label{rem:typeA-configuration-twisting}
Correspondingly, if one thinks of the factorization line bundles rather than factorization gerbes, the factorization line bundles on the configuration space can be described similarly, see \cite{[BFT1],[TY],[FTY]}.
\end{rem}

\subsubsection{The canonical character $\psi$} The advantage of considering $\omega$-twisted groups and affine flags is that the character $\psi$ is canonical in this case. Namely, the character $\psi$ of $H^\omega(\bF)$, as well as $H_u^{\omega}(\bF)$ and $I_1\ltimes H_u^{\omega}(\bF)$, is given by the following map
\begin{equation}\label{chi}
\begin{split}
\psi:H^\omega(\bF)&\longrightarrow\BG_a,\\
gu&\longmapsto\operatorname{Res}_c\bigl(\psi_\omega(u)\bigr), \qquad g\in H_r^\omega(\bF),\quad u\in H_u^\omega(\bF).
\end{split}
\end{equation}
Here, $\psi_\omega:H_u^\omega(\bF)\to\omega(\overset{\circ}{\cD}_c)$ is induced by the $H_r$-invariant character chosen in Section~\ref{pol des}.

\begin{defn}
In the homogeneous case, we define the twisted Iwahori Gaiotto category as $\SD_{\cG^G}^{H^\omega(\bF),\psi}(\Fl_G^\omega)$, i.e., the category of $\cG^G$-twisted $(H^\omega(\bF),\psi)$-equivariant sheaves on the renormalized affine flags.

In the Betti setting and in the $\ell$-adic setting in characteristic zero, the $(H^\omega(\bF),\psi)$-equivariance for $\psi\neq0$ is understood via the Kirillov model, see Section~\ref{sec:positive-characteristic}.

Furthermore, via $!$-restriction along the embedding
\begin{equation}
\begin{split}
\Fl_{G_2}^{\omega} &\simeq \{1\}\times \Fl_{G_2}^{\omega} \hookrightarrow \Fl_G^\omega\simeq \Fl_{G_1}\times \Fl_{G_2}^{\omega},
\end{split}
\end{equation}
we have
\begin{equation}
\SD_{\cG^G}^{H^\omega(\bF),\psi}(\Fl_G^\omega) \simeq \SD_{\cG^G}^{I_1\ltimes H_u^\omega(\bF),\psi}(\Fl_{G_2}^{\omega}).
\end{equation}

In the case of $\GL(N|N)$, we use the $D^{!*}$-version of the category of sheaves on $\Fl_N\times\bF^N$ as in \cite[\S 3.2]{[BFGT]}. Namely, in the vector direction, the transition functor between the finite-dimensional quotients along a smooth projection $p$ of relative dimension $d$ is $p^*[d]$ (with the additional Tate twist $(d/2)$ in the mixed setting).

We call an orbit relevant if it supports a non-zero $\cG^G$-twisted $(H^\omega(\bF),\psi)$-equivariant sheaf.

For an object $\cF$ in the twisted Iwahori Gaiotto category, we call $\cF$ a locally compact object if it is compact when forgetting the equivariance of the torus part of $I_1$. In the homogeneous case, it means that $\cF$ is supported on finitely many relevant orbits, and its $!$-stalk at a point of each orbit is bounded and finite-dimensional. In the case of $\GL(N|N)$, it means that $\cF$ is compact when regarded as a sheaf on $\Fl_N\times \bF^N$.

We denote the resulting full subcategory by
\[
\SD_{\cG^G}^{H^\omega(\bF),\psi,\locc}(\Fl_G^\omega) \quad\text{or}\quad \SD_{\cG^G}^{I_1\ltimes H_u^\omega(\bF),\psi,\locc}(\Fl_{G_2}^{\omega}).
\]
Similarly, for an object $\cF$ in the spherical Gaiotto category, we call $\cF$ a locally compact object if its restriction to $\{1\}\times\Gr_{G_2}^{\omega}$ is compact when forgetting the $G_1(\bO)$-equivariance, cf., \cite[Definition 3.2.4]{[TY]}. In the case of $\GL(N|N)$, it means that $\cF$ is compact when regarded as a sheaf on $\Gr_N\times\bF^N$, in the same $D^{!*}$-version as above. We denote the resulting full subcategory by $\SD_{\cG^G}^{H^\omega(\bF),\psi,\locc}(\Gr_G^\omega)$.
\end{defn}

\subsubsection{Basic properties}

We now study some basic properties of the twisted Iwahori Gaiotto category
\[
\SD_{\cG^G}^{H^\omega(\bF),\psi}(\Fl_G^\omega) \simeq \SD_{\cG^G}^{I_1\ltimes H_u^\omega(\bF),\psi}(\Fl_{G_2}^{\omega}).
\]

There is a continuous forgetful functor
\[
\oblv\colon \SD_{\cG^G}^{I_1\ltimes H_u^\omega(\bF),\psi}(\Fl_{G_2}^{\omega}) \longrightarrow \SD_{\cG^G}(\Fl_{G_2}^{\omega}).
\]
By the adjoint functor theorem, it admits a right adjoint $\Av_*^{I_1\ltimes H_u^\omega(\bF),\psi}$, which is not continuous in general.

Since $H_u^\omega(\bF)$ is an ind-pro-group, we may choose a presentation
\[
H_u^\omega(\bF)\simeq \colim_i H_{u,i},
\]
where each $H_{u,i}$ is a pro-finite type group scheme, the transition maps are closed embeddings, and each $H_{u,i}$ is stable under the conjugation action of $I_1$.

The forgetful functor $\oblv$ admits a partially defined left adjoint
\[
\Av_!^{I_1\ltimes H_u^\omega(\bF),\psi} := \colim_i \Av_!^{I_1\ltimes H_{u,i},\psi}.
\]
For an object $\cF\in \SD_{\cG^G}(\Fl_{G_2}^{\omega})$, the object $\Av_!^{I_1\ltimes H_u^\omega(\bF),\psi}(\cF)$ is defined if $\Av_!^{I_1\ltimes H_{u,i},\psi}(\cF)$ is defined for every $i$. In particular, this holds whenever $\cF$ is holonomic.

\begin{lem}
For a spherical variety $\bX$, there are countably many Iwahori orbits in $\bX(\bF)$. Thus, there are countably many irreducible objects in the Iwahori Gaiotto category.
\end{lem}
\begin{proof}
Since there are only finitely many $G$-orbits in $\bX$ and each of them is spherical, it is enough to consider the case $\bX=G/H$. According to \cite[Theorem 3.2.1]{[GN]}, there are countably many $G(\bO)$-orbits in $\bX(\bF)$. For $x\in\bX(\bF)$, we denote by $K_x=\Stab_{G(\bO)}(x)$ its stabilizer, and by $H_x\subset G$ the image of $K_x$ under the evaluation map at $t=0$. Then, $H_x$ is spherical by \cite[Theorem 1.3]{[BMo]}. Since $I$ is the preimage of $B$, there is a bijection
\[
I\backslash G(\bO)/K_x\simeq B\backslash G/H_x.
\]
Note that the right-hand side is finite, which proves the first assertion. For the second assertion, note that the stabilizers have finite component groups, thus there are only finitely many irreducible twisted equivariant local systems on each relevant orbit.
\end{proof}

\begin{lem}
The twisted Iwahori Gaiotto category
\[
\SD_{\cG^G}^{H^\omega(\bF),\psi}(\Fl_G^\omega)\simeq \SD_{\cG^G}^{I_1\ltimes H_u^\omega(\bF),\psi}(\Fl_{G_2}^\omega)
\]
is compactly generated. In particular, it is cocomplete and dualizable.
\end{lem}

\begin{proof}
Let $\{\sO_{\Fl,\tilde w}\}$ be the collection of $(\cG^G,\psi)$-relevant $H^\omega(\bF)$-orbits in $\Fl_G^\omega$, which is a countable set. Let $\{\BO_{\tilde w}\}$ be the corresponding $(\cG^G,\psi)$-relevant $(I_1\ltimes H_u^\omega(\bF))$-orbits in $\Fl_{G_2}^\omega$. For each relevant orbit $\BO_{\tilde w}$, we choose a point $x_{\tilde w}\in \BO_{\tilde w}$ and a trivialization of $\cG^G$ at $x_{\tilde w}$, and let $\delta_{\tilde w}$ be the corresponding skyscraper sheaf. Since $\delta_{\tilde w}$ is holonomic, the object $\Av_!^{I_1\ltimes H_u^\omega(\bF),\psi}(\delta_{\tilde w})$ is well-defined.

We claim that the objects $\Av_!^{I_1\ltimes H_u^\omega(\bF),\psi}(\delta_{\tilde w})$, where $\BO_{\tilde w}$ runs over all relevant orbits, compactly generate $\SD_{\cG^G}^{I_1\ltimes H_u^\omega(\bF),\psi}(\Fl_{G_2}^\omega)$.

Indeed, let $\cF\in \SD_{\cG^G}^{I_1\ltimes H_u^\omega(\bF),\psi}(\Fl_{G_2}^\omega)$, and assume that
\[
\RHom_{\SD_{\cG^G}^{I_1\ltimes H_u^\omega(\bF),\psi}(\Fl_{G_2}^\omega)}\bigl(\Av_!^{I_1\ltimes H_u^\omega(\bF),\psi}(\delta_{\tilde w}),\cF\bigr)=0
\]
for all relevant $\BO_{\tilde w}$. By adjunction, this implies
\[
\RHom_{\SD_{\cG^G}(\Fl_{G_2}^\omega)}\bigl(\delta_{\tilde w},\oblv(\cF)\bigr)=0
\]
for all relevant $\BO_{\tilde w}$. Hence, $\oblv(\cF)$ has zero $!$-stalk at the chosen point $x_{\tilde w}$ of every relevant orbit $\BO_{\tilde w}$. Since $\cF$ is equivariant and its restrictions to the non-relevant orbits vanish, we obtain $\cF=0$. This proves the generation.

To prove the compactness, note that $\delta_{\tilde w}$ is compact in $\SD_{\cG^G}(\Fl_{G_2}^\omega)$. Since $\oblv$ is continuous, for any filtered diagram $\{\cF_\alpha\}$ in $\SD_{\cG^G}^{I_1\ltimes H_u^\omega(\bF),\psi}(\Fl_{G_2}^\omega)$ we have
\begin{equation}
\begin{split}
\RHom_{\SD_{\cG^G}^{I_1\ltimes H_u^\omega(\bF),\psi}(\Fl_{G_2}^\omega)}\bigl(\Av_!^{I_1\ltimes H_u^\omega(\bF),\psi}(\delta_{\tilde w}),\colim_\alpha \cF_\alpha\bigr)\\
\simeq \RHom_{\SD_{\cG^G}(\Fl_{G_2}^\omega)}\bigl(\delta_{\tilde w},\oblv(\colim_\alpha \cF_\alpha)\bigr)\\
\simeq \colim_\alpha \RHom_{{\SD_{\cG^G}(\Fl_{G_2}^\omega)}}\bigl(\delta_{\tilde w},\oblv(\cF_\alpha)\bigr)\\
\simeq \colim_\alpha \RHom_{{\SD_{\cG^G}^{I_1\ltimes H_u^\omega(\bF),\psi}(\Fl_{G_2}^\omega)}}\bigl(\Av_!^{I_1\ltimes H_u^\omega(\bF),\psi}(\delta_{\tilde w}),\cF_\alpha\bigr).
\end{split}
\end{equation}

Thus, $\Av_!^{I_1\ltimes H_u^\omega(\bF),\psi}(\delta_{\tilde w})$ is compact for every relevant $\BO_{\tilde w}$, and the category is compactly generated.
\end{proof}

\section{Local functors}\label{sec:local-functors}
\subsection{Factorization affine Grassmannian and factorization affine flag}

In order to define the desired equivalence functor, let us recall the factorization (or Beilinson--Drinfeld, or Ran-ified) affine Grassmannian $\Gr_{G,\Ran}^{\omega}$ and the factorization affine flag $\Fl_{G,\Ran_c}^{\omega}$.

\subsubsection{Ran-ified version}\label{ranified version} Given a curve $C$ with a marked point, we recall that the Ran space $\Ran$ is the prestack classifying non-empty finite subsets of $C$, and the Ran space $\Ran_c$ with the marked point is the prestack classifying non-empty finite subsets of $C$ containing the marked point $c$. We denote by $\Gr_{G,\Ran}^{\omega}$ the prestack over $\Ran$ classifying triples $(\bar{x},\cP_G,\alpha)$, where $\bar{x}\in \Ran$ is a non-empty finite subset of $C$, $\cP_G$ is a $G$-bundle on $C$, and
\begin{equation}
\alpha:\cP_G|_{C\setminus \bar{x}}\xrightarrow{\sim}\cP_G^\omega|_{C\setminus \bar{x}}=(\omega^\rho\times^T G)|_{C\setminus \bar{x}},
\end{equation}
is an isomorphism away from $\bar{x}$. Similarly, $\Fl_{G,\Ran_c}^{\omega}$ is the prestack over $\Ran_c$ classifying quadruples $(\bar{x},\cP_G,\alpha,\epsilon)$, where $\bar{x}\in \Ran_c$ is a non-empty finite subset of $C$ containing the marked point $c$, $(\cP_G,\alpha)$ is as above, and $\epsilon$ is a $B$-reduction of the fiber $\cP_G|_c$.

They satisfy the usual factorization property:
\begin{equation}
\begin{split}
\Gr_{G,\Ran}^{\omega}\underset{\Ran}{\times}(\Ran\times\Ran)_{\disj} &\simeq \bigl(\Gr_{G,\Ran}^{\omega}\times \Gr_{G,\Ran}^{\omega}\bigr) \underset{\Ran\times\Ran}{\times}(\Ran\times\Ran)_{\disj},\\
\Fl_{G,\Ran_c}^{\omega}\underset{\Ran_c}{\times}(\Ran\times\Ran_c)_{\disj} &\simeq \bigl(\Gr_{G,\Ran}^{\omega}\times \Fl_{G,\Ran_c}^{\omega}\bigr) \underset{\Ran\times\Ran_c}{\times}(\Ran\times\Ran_c)_{\disj}.
\end{split}
\end{equation}
The fiber of $\Gr_{G,\Ran}^{\omega}$ (resp. $\Fl_{G,\Ran_c}^{\omega}$) over the point $\bar{x}=\{x_1,\cdots,x_l\}\in \Ran$ (resp. $\bar{x}=\{c,x_1,\cdots,x_l\}\in \Ran_c$) is given by
\begin{equation}
\bigl(\Gr_{G,\Ran}^{\omega}\bigr)_{\bar{x}}\simeq \Gr_{G,\bar{x}}^{\omega}\simeq \prod_{i=1}^l \Gr_{G,x_i}^\omega
\end{equation}
\begin{samepage}
(resp.
\begin{equation}
\bigl(\Fl_{G,\Ran_c}^{\omega}\bigr)_{\bar{x}}\simeq \Fl_{G,\bar{x}}^{\omega}\simeq\Fl_{G,c}^\omega\times \prod_{i=1}^l \Gr_{G,x_i}^\omega
\end{equation}
).\par
\end{samepage}

We have the sub-prestack ${\sO}_{\Gr,\Ran}$ of $\Gr_{G,\Ran}^{\omega}$, whose fiber over $\bar{x}$ is given by the product of the unital $H^\omega(\bF)$-orbits, i.e., $\prod_{i=1}^l \sO^0_{\Gr,x_i}$. Its closure is denoted by $\overline{\sO}_{\Gr,\Ran}$, whose fiber over $\bar{x}$ is given by the product of the closures of the unital orbits. When we consider the Ran-ified affine Grassmannian with a marked point $\Gr_{G,\Ran_c}^{\omega}$, in addition to ${\sO}_{\Gr,\Ran_c}$ and $\overline{\sO}_{\Gr,\Ran_c}$ which are given by base change, we also have $(\overline{\sO}_{\Gr,\Ran_c})_{\infty\cdot c}$, whose fiber over $\bar{x}\in \Ran_c$ is $\Gr_{G,c}^\omega\times\prod_{i=1}^l \overline{\sO}_{\Gr,x_i}$. For the Ran-ified affine flags $\Fl_{G,\Ran_c}^{\omega}$, we have $(\overline{\sO}_{\Fl,\Ran_c})_{\infty\cdot c}$, whose fiber is given by $\Fl_{G,c}^\omega\times\prod_{i=1}^l \overline{\sO}_{\Gr,x_i}$.

\subsubsection{Configuration version} Similarly, one defines the factorization affine Grassmannian $\Gr_{G,\Conf}^{\omega}$ over $\Conf$ and the factorization affine flag $\Fl_{G,\Conf_{\infty\cdot c}}^{\omega}$ over $\Conf_{\infty\cdot c}$. Their fibers are given by
\begin{equation}
\bigl(\Gr_{G,\Conf}^{\omega}\bigr)_D\simeq \prod_{i=1}^l \Gr_{G,x_i}^{\omega}\qquad\text{for }D=\sum_{i=1}^l \alpha_{x_i}\cdot x_i\in \Conf,
\end{equation}
and
\begin{equation}
\bigl(\Fl_{G,\Conf_{\infty\cdot c}}^{\omega}\bigr)_D\simeq \Fl_{G,c}^{\omega}\times \prod_{i=1}^l \Gr_{G,x_i}^{\omega}\qquad\text{for }D=\lambda_c\cdot c+\sum_{i=1}^l \alpha_{x_i}\cdot x_i\in \Conf_{\infty\cdot c}.
\end{equation}

We denote by $W^{\fin}$ the finite Weyl group of $G$. We will use the factorization sub-prestacks
\begin{equation*}
\fitmath{\overline{\sO}_{\Gr,\Conf},\quad \overline{\sO}_{\Fl,\Conf_{\infty\cdot c}},\quad \overline{\sS}_{\Gr,\Conf},\quad \overline{\sS}_{\Fl,\Conf_{\infty\cdot c}},\quad \sS_{\Gr,\Conf},\quad \sS_{\Fl,\Conf_{\infty\cdot c}},\quad \sS_{\Fl,\Conf_{\infty\cdot c}}^w\ \ (w\in W^{\fin}),}
\end{equation*}
defined in the usual way. We only record their geometric fibers below.

\subsubsection{Sub-prestacks of the factorization affine flags}\label{subpre in fl}

Among the above sub-prestacks, $\overline{\sO}_{\Fl,\Conf_{\infty\cdot c}}$, $\overline{\sS}_{\Fl,\Conf_{\infty\cdot c}}$, $\sS_{\Fl,\Conf_{\infty\cdot c}}$, and $\sS_{\Fl,\Conf_{\infty\cdot c}}^w$ are factorization module sub-prestacks of $\Fl_{G,\Conf_{\infty\cdot c}}^{\omega}$. Their fibers over
\begin{equation}
D=\lambda_c\cdot c+\sum_{i=1}^l \alpha_{x_i}\cdot x_i\in \Conf_{\infty\cdot c}
\end{equation}
are given by
\begin{equation}
\begin{split}
\overline{\sO}_{\Fl,D} &\simeq \Fl_{G,c}^{\omega}\times \prod_{i=1}^l \overline{\sO}_{\Gr,x_i},\qquad \overline{\sS}_{\Fl,D} \simeq \overline{\sS}_{\Fl,c}^{\lambda_c}\times \prod_{i=1}^l \overline{\sS}_{\Gr,x_i}^{\alpha_{x_i}},\\
\sS_{\Fl,D} &\simeq \sS_{\Fl,c}^{\lambda_c}\times \prod_{i=1}^l \sS_{\Gr,x_i}^{\alpha_{x_i}},\qquad \sS_{\Fl,D}^w \simeq \sS_{\Fl,c}^{\lambda_c w}\times \prod_{i=1}^l \sS_{\Gr,x_i}^{\alpha_{x_i}}.
\end{split}
\end{equation}
Here, we denote $\sS_{\Fl,c}^{\lambda_c}=\sS_{\Fl,c}^{\lambda_c e}$, and the semi-infinite orbits in this section are the $U(\bF)^\omega$-orbits, where $U$ denotes the unipotent radical of $B$. Also, $\overline{\sS}_{\Fl,c}^{\lambda_c}$ is the preimage in $\Fl_{G,c}^{\omega}$ of the closure of the semi-infinite orbit $\sS_{\Gr,c}^{\lambda_c}\subset \Gr_{G,c}^{\omega}$ through $t^{\lambda_c}$, and $\sS_{\Fl,c}^{\lambda_c w}$ is the semi-infinite orbit in $\Fl_{G,c}^{\omega}$ through $t^{\lambda_c}w\in \Fl_{G,c}^{\omega}$, where $w\in W^{\fin}$.

The factorization module structure is expressed by canonical isomorphisms
\begin{equation}
\label{3.19 iso}\fitmath{\begin{aligned} &\overline{\sO}_{\Fl,\Conf_{\infty\cdot c}}\underset{\Conf_{\infty\cdot c}}{\times}(\Conf\times \Conf_{\infty\cdot c})_{\disj} \simeq \bigl(\overline{\sO}_{\Gr,\Conf}\times \overline{\sO}_{\Fl,\Conf_{\infty\cdot c}}\bigr) \underset{\Conf\times \Conf_{\infty\cdot c}}{\times}(\Conf\times \Conf_{\infty\cdot c})_{\disj},\\
&\overline{\sS}_{\Fl,\Conf_{\infty\cdot c}}\underset{\Conf_{\infty\cdot c}}{\times}(\Conf\times \Conf_{\infty\cdot c})_{\disj} \simeq \bigl(\overline{\sS}_{\Gr,\Conf}\times \overline{\sS}_{\Fl,\Conf_{\infty\cdot c}}\bigr) \underset{\Conf\times \Conf_{\infty\cdot c}}{\times}(\Conf\times \Conf_{\infty\cdot c})_{\disj},\\
&\sS_{\Fl,\Conf_{\infty\cdot c}}^w\underset{\Conf_{\infty\cdot c}}{\times}(\Conf\times \Conf_{\infty\cdot c})_{\disj} \simeq \bigl(\sS_{\Gr,\Conf}\times \sS_{\Fl,\Conf_{\infty\cdot c}}^w\bigr) \underset{\Conf\times \Conf_{\infty\cdot c}}{\times}(\Conf\times \Conf_{\infty\cdot c})_{\disj}.
\end{aligned}}
\end{equation}
\subsubsection{Sub-prestacks of the factorization affine Grassmannian}

Also, $\overline{\sO}_{\Gr,\Conf}$, $\overline{\sS}_{\Gr,\Conf}$, and $\sS_{\Gr,\Conf}$ are factorization sub-prestacks of $\Gr_{G,\Conf}^{\omega}$. Their fibers over
\begin{equation}
D=\sum_{i=1}^l \alpha_{x_i}\cdot x_i\in \Conf
\end{equation}
are
\begin{equation}
\overline{\sO}_{\Gr,D} \simeq \prod_{i=1}^l \overline{\sO}_{\Gr,x_i},\qquad \overline{\sS}_{\Gr,D} \simeq \prod_{i=1}^l \overline{\sS}_{\Gr,x_i}^{\alpha_{x_i}},\qquad \sS_{\Gr,D} \simeq \prod_{i=1}^l \sS_{\Gr,x_i}^{\alpha_{x_i}}.
\end{equation}
Here, $\overline{\sO}_{\Gr,x_i}$ is the closure of the unital $H(\bF)_{x_i}^{\omega}$-orbit in $\Gr_{G,x_i}^{\omega}$, and $\overline{\sS}_{\Gr,x_i}^{\alpha_{x_i}}$ is the closure of the semi-infinite orbit $\sS_{\Gr,x_i}^{\alpha_{x_i}}\subset \Gr_{G,x_i}^{\omega}$ through $t^{\alpha_{x_i}}\in \Gr_{G,x_i}^{\omega}$.

As in \eqref{3.19 iso}, their factorization property is expressed by canonical isomorphisms
\begin{equation}\label{3.24 iso}
\begin{split}
&\overline{\sO}_{\Gr,\Conf}\underset{\Conf}{\times}(\Conf\times\Conf)_{\disj} \simeq \bigl(\overline{\sO}_{\Gr,\Conf}\times \overline{\sO}_{\Gr,\Conf}\bigr) \underset{\Conf\times\Conf}{\times}(\Conf\times\Conf)_{\disj},\\
&\overline{\sS}_{\Gr,\Conf}\underset{\Conf}{\times}(\Conf\times\Conf)_{\disj} \simeq \bigl(\overline{\sS}_{\Gr,\Conf}\times \overline{\sS}_{\Gr,\Conf}\bigr) \underset{\Conf\times\Conf}{\times}(\Conf\times\Conf)_{\disj},\\
&\sS_{\Gr,\Conf}\underset{\Conf}{\times}(\Conf\times\Conf)_{\disj} \simeq \bigl(\sS_{\Gr,\Conf}\times \sS_{\Gr,\Conf}\bigr) \underset{\Conf\times\Conf}{\times}(\Conf\times\Conf)_{\disj}.
\end{split}
\end{equation}

\subsection{The construction of local functors}
\subsubsection{Zastava spaces} We denote by
\[
v_{\Fl}:\Fl_{G,\Conf_{\infty\cdot c}}^\omega\longrightarrow\Conf_{\infty\cdot c}, \qquad v_{\Gr}:\Gr_{G,\Conf}^\omega\longrightarrow\Conf
\]
the projections to the configuration spaces. With some abuse of notation, we will use the same notations for their restrictions to the semi-infinite strata, their closures, and the Sakellaridis--Wang (SW) Zastava spaces (also in the version with a marked point).

Consider the diagram
\begin{equation}
\fitmath{\xymatrix{ & \Ran_c\times \Fl_{G,c}^{\omega}\ar[ld]_{\pr_{\Fl}}\ar[rd]^{\unit_{\Ran_c,\Fl}}&&\overline{\sO}_{\Fl,\Conf_{\infty\cdot c}}\ar[ld]^{\pr_{\Conf_{\infty\cdot c}\to \Ran_c}}&\overline{\cY}_{\Fl,\infty\cdot c}\ar[l]_{\iota_{\Fl}}\ar[rd]^{v_{\Fl}}&\\
\Fl_{G,c}^{\omega}&&(\overline{\sO}_{\Fl,\Ran_c})_{\infty\cdot c}&&&\Conf_{\infty\cdot c}. }}
\end{equation}
Here,
\begin{equation}\label{3.26 zas}
\begin{split}
\overline{\cY}_{\Fl,\infty\cdot c} &:=\overline{\sO}_{\Fl,\Conf_{\infty\cdot c}}\cap \overline{\sS}_{\Fl,\Conf_{\infty\cdot c}},\\
\cY_{\Fl,\infty\cdot c} &:=\overline{\sO}_{\Fl,\Conf_{\infty\cdot c}}\cap\sS_{\Fl,\Conf_{\infty\cdot c}},
\end{split}
\end{equation}
are the Iwahori versions of the Sakellaridis--Wang Zastava spaces.

\subsubsection{Spreading} We denote by $H(\bF)_{\Ran}^{\omega}$ and $H(\bO)_{\Ran}^{\omega}$ the loop and arc group prestacks over $\Ran$, whose fibers at $\bar x\in\Ran$ are $\prod_{x\in\bar x}H(\bF)_x^\omega$ and $\prod_{x\in\bar x}H(\bO)_x^\omega$, and by $\psi_{\Ran}$ the sum of the local characters $\psi$ at these points. The subscripts $\Ran_c$ and $\Conf$ denote their pullbacks to $\Ran_c$ and to the configuration spaces, respectively.

Given an object $\cF\in \SD_{\cG^G}^{H(\bF)_c^\omega,\psi}(\Fl_{G,c}^{\omega})$ in the twisted Gaiotto category on the affine flags, we first take its $!$-pullback along $\pr_{\Fl}$. This gives an object of $\SD_{\cG^G}^{H(\bF)_c^\omega,\psi}(\Ran_c\times \Fl_{G,c}^{\omega})$, where $H(\bF)_c^\omega$ acts on the second factor.

We denote by $\sprd_{\Ran_c,\Fl}$ the spreading functor constructed as in \cite[Section 6.1]{[Yan]}, cf., \cite[Section 7.2]{[TY]}, with the Iwahori structure at $c$. The object $\sprd_{\Ran_c,\Fl}(\cF)$ satisfies
\begin{equation}
\unit_{\Ran_c,\Fl}^!(\sprd_{\Ran_c,\Fl}(\cF)) \simeq \pr_{\Fl}^!(\cF).
\end{equation}

As in \cite[\S 7.3.1]{[TY]}, we pull back $\sprd_{\Ran_c,\Fl}(\cF)$ to $\overline{\sO}_{\Fl,\Conf_{\infty\cdot c}}$, and denote the composed functor by
\begin{equation}
\label{sprd 3.32}\fitmath{\sprd:=\pr_{\Conf_{\infty\cdot c}\to \Ran_c}^!\circ\sprd_{\Ran_c,\Fl}\colon \SD_{\cG^G}^{H(\bF)_c^\omega,\psi}(\Fl_{G,c}^{\omega})\rightarrow \SD_{\cG^G}^{H(\bF)_{\Conf}^\omega,\psi_{\Conf}}(\overline{\sO}_{\Fl,\Conf_{\infty\cdot c}}).}
\end{equation}
By abuse of notation, we also denote its further pullback along $\iota_{\Fl}$ to $\overline{\cY}_{\Fl,\infty\cdot c}$ by the same symbol.

\subsubsection{Semi-infinite sheaves} Recall the sub-prestacks $\sS_{\Fl,\Conf_{\infty\cdot c}}$ and $\overline{\sS}_{\Fl,\Conf_{\infty\cdot c}}$ in Section \ref{subpre in fl}. Let
\begin{equation}
j_{\Fl}^{\sinf}\colon \sS_{\Fl,\Conf_{\infty\cdot c}}\hookrightarrow \overline{\sS}_{\Fl,\Conf_{\infty\cdot c}}
\end{equation}
be the locally closed embedding. The restriction of $(\cG^G)^{-1}\otimes v_{\Fl}^*(\cG^\Lambda)$ to $\sS_{\Fl,\Conf_{\infty\cdot c}}$ admits a canonical trivialization. Hence, by abuse of notation, we may regard $\omega_{\sS_{\Fl,\Conf_{\infty\cdot c}}}$ as a $(\cG^G)^{-1}\otimes v_{\Fl}^*(\cG^\Lambda)$-twisted sheaf on $\sS_{\Fl,\Conf_{\infty\cdot c}}$. For the $!*$-extension, we use the semi-infinite IC sheaf, cf., \cite{[G5],[G6]} and \cite[Section 13]{[GL]}. We take the corresponding extensions as $(\cG^G)^{-1}\otimes v_{\Fl}^*(\cG^\Lambda)$-twisted sheaves
\begin{equation}\label{3.34 sinf}
j^{\sinf}_{\Fl,!}(\omega_{\sS_{\Fl,\Conf_{\infty\cdot c}}}),\qquad j^{\sinf}_{\Fl,*}(\omega_{\sS_{\Fl,\Conf_{\infty\cdot c}}}),\qquad j^{\sinf}_{\Fl,!*}(\omega_{\sS_{\Fl,\Conf_{\infty\cdot c}}}).
\end{equation}
They are $(\cG^G)^{-1}\otimes v_{\Fl}^*(\cG^\Lambda)$-twisted sheaves on $\overline{\sS}_{\Fl,\Conf_{\infty\cdot c}}$. Since the underlying spaces are different, there is no risk of ambiguity, and we use the same notations for their $!$-restrictions to $\overline{\cY}_{\Fl,\infty\cdot c}$.

\subsubsection{Normalization} Let $\deg$ be the function on the coweight lattice
\begin{equation}
\deg\colon \Lambda\rightarrow \BZ,
\end{equation}
given by $\deg(\lambda)=\langle2\check\rho_G,\lambda\rangle$, where $\check\rho_G$ is the half-sum of the geometric positive roots with respect to $B^-$. It is the geometric degree which will be used in the calculation of the BMW sheaves in Section~\ref{foundations-equivalence}.

We define $\delta$, $d$, and the geometric weight $2\check\rho^\circ$ by
\begin{equation}\label{eq:local-degree-normalization}
\begin{aligned}
\delta(\lambda)&=
\begin{cases}
(N-M-1)\displaystyle\sum_{i=1}^M\lambda_i^{(1)},&\text{in type A},\\
0,&\text{in the other cases},
\end{cases}\\
d(\lambda)&=\deg(\lambda)+\delta(\lambda)=\langle2\check\rho^\circ,\lambda\rangle.
\end{aligned}
\end{equation}
Note that the weight $2\check\rho^\circ$ is obtained from $2\check\rho_G$ by adding a central character,{\interfootnotelinepenalty=10000\footnote{In type A with $M<N$, $\delta$ comes from the character $\det^{-(N-M-1)}$ of $\GL_M$ on $\Lie(H_u)$, see Section~\ref{sec:global-dimension-normalization}. In our geometric coordinates, the pairing with $2\rho^\circ$ in \cite[Section 3.4]{[TY]} is $d(\lambda)+(N-M-1)(\sum_j\lambda_j^{(2)}-\sum_i\lambda_i^{(1)})$, thus the two pairings agree on the lattice of super roots.}} such that $d(\alpha_i)=1$ for any odd isotropic simple super root $\alpha_i$, and $d(\alpha_i)=2$ for any even simple super root $\alpha_i$. In particular, on the lattice of geometric coroots, $\frac12d$ is the usual height function with respect to $B^-$, and $d$ gives rise to a linear function $\Lambda\rightarrow \BZ/2$, which extends the parity on the lattice of super roots. Also, for $\alpha\in\Lambda^{\mathrm{neg}}$, the open SW Zastava space (without the marked point) of degree $\alpha$ has dimension $-d(\alpha)$, see Section~\ref{sec:global-dimension-normalization}.

\begin{rem}\label{rem:target-parity}
On the representation side, we denote by $\se_\lambda$ the one-dimensional even $T_s$-representation of weight $\lambda$, and we take it as the highest weight line of $\Verma_\lambda$ and $\Irr_\lambda$. Similarly, the generators defining $\Delta_\lambda^{\fact}$, $\nabla_\lambda^{\fact}$ and $\IC_\lambda^{\fact}$ are even, see Section~\ref{sec:configuration-modules}.
\end{rem}

We define the functors
\begin{equation}
F_{\Fl}^{\Lus},\ F_{\Fl}^{\DK},\ F_{\Fl}^{\sm}:\SD_{\cG^G}^{H(\bF)_c^\omega,\psi}(\Fl_{G,c}^{\omega})\longrightarrow \SD_{\cG^\Lambda}(\Conf_{\infty\cdot c})
\end{equation}
by
\begin{equation}
\begin{split}
F_{\Fl}^{\Lus}(\cF) &:=\Pi^d v_{\Fl,*}\bigl(\sprd(\cF)\overset{!}{\otimes} j^{\sinf}_{\Fl,!}(\omega_{\sS_{\Fl,\Conf_{\infty\cdot c}}})\bigr)[d],\\
F_{\Fl}^{\DK}(\cF) &:=\Pi^d v_{\Fl,*}\bigl(\sprd(\cF)\overset{!}{\otimes} j^{\sinf}_{\Fl,*}(\omega_{\sS_{\Fl,\Conf_{\infty\cdot c}}})\bigr)[d],\\
F_{\Fl}^{\sm}(\cF) &:=\Pi^d v_{\Fl,*}\bigl(\sprd(\cF)\overset{!}{\otimes} j^{\sinf}_{\Fl,!*}(\omega_{\sS_{\Fl,\Conf_{\infty\cdot c}}})\bigr)[d].
\end{split}
\end{equation}

\begin{rem}
Here, the function $d$ is evaluated on the total degree of the configuration. Note that the cohomological sign and the super Koszul sign of $\Pi^d[d]$ cancel each other, thus this normalization preserves the commutativity constraint of the factorization structure, cf. \cite[\S 7.3.2]{[TY]}.

\end{rem}

\subsubsection{Description of $!$-stalks} Given $\lambda\in \Lambda$, let
\begin{equation}
i_\lambda\colon \{\lambda\cdot c\}\hookrightarrow \Conf_{\infty\cdot c}
\end{equation}
be the closed embedding. The fiber of $\Fl_{G,\Conf_{\infty\cdot c}}^{\omega}$ over $\lambda\cdot c$ is canonically identified with $\Fl_{G,c}^{\omega}$. By base change (for the pushforward $v_{\Fl,*}$ and for the semi-infinite kernel), we obtain
\begin{equation}
i_\lambda^!F_{\Fl}^{\DK}(\cF)\simeq \Pi^{d(\lambda)}H\bigl(\Fl_{G,c}^{\omega},\, \cF\overset{!}{\otimes}j^{\sinf}_{\Fl,\lambda,*}(\omega_{\sS_{\Fl,c}^{\lambda}})\bigr)[d(\lambda)],
\end{equation}
where $j^{\sinf}_{\Fl,\lambda,*}$ is the $*$-extension of the dualizing sheaf on $\sS_{\Fl,c}^{\lambda}$. Here, the base change for the $*$-extension follows from the compatibility of the $!$-pullback with the $*$-extension. For the $!$-extension and the $!*$-extension, we need the corresponding non-trivial base change theorems, i.e., \cite[Corollary A.2.9 and Proposition 6.5.2]{[Yan]} (which are direct analogy of \cite[Section 3.6]{[G6]}), and the local-global comparison of the semi-infinite IC sheaves in \cite[Theorem 3.3.3]{[G6]}, respectively. Thus, one similarly has
\begin{equation}
\begin{split}
i_\lambda^!F_{\Fl}^{\Lus}(\cF) &\simeq \Pi^{d(\lambda)}H\bigl(\Fl_{G,c}^{\omega},\, \cF\overset{!}{\otimes}j^{\sinf}_{\Fl,\lambda,!}(\omega_{\sS_{\Fl,c}^{\lambda}})\bigr)[d(\lambda)],\\
i_\lambda^!F_{\Fl}^{\sm}(\cF) &\simeq \Pi^{d(\lambda)}H\bigl(\Fl_{G,c}^{\omega},\, \cF\overset{!}{\otimes}j^{\sinf}_{\Fl,\lambda,!*}(\omega_{\sS_{\Fl,c}^{\lambda}})\bigr)[d(\lambda)],
\end{split}
\end{equation}
where $j^{\sinf}_{\Fl,\lambda,!}$ and $j^{\sinf}_{\Fl,\lambda,!*}$ are the $!$- and $!*$-extensions of the (twisted) dualizing sheaf on $\sS_{\Fl,c}^{\lambda}$.

\subsection{Factorization property}\label{3.3 fact prop}

Now, let us consider the affine Grassmannian. Let $\overline{\sO}_{\Gr,c}$ be the closure of the unital $H(\bF)_c^\omega$-orbit in $\Gr_{G,c}^{\omega}$, and let $\overline{\sO}_{\Gr,\Ran}$ and $\overline{\sO}_{\Gr,\Conf}$ be the corresponding factorization versions. We denote
\begin{equation}\label{3.44 zas}
\begin{split}
\overline{\cY}_{\Gr} &:=\overline{\sO}_{\Gr,\Conf}\underset{\Gr_{G,\Conf}^{\omega}}{\times}\overline{\sS}_{\Gr,\Conf},\\
\cY_{\Gr} &:=\overline{\sO}_{\Gr,\Conf}\underset{\Gr_{G,\Conf}^{\omega}}{\times}\sS_{\Gr,\Conf}.
\end{split}
\end{equation}
They are the compactified and non-compactified SW Zastava spaces, see \cite{[SW]}, and we denote by $v_{\Gr}$ the natural projection $\overline{\cY}_{\Gr}\to\Conf$.

\subsubsection{The unit object} The following lemma is well-known, and can be obtained by direct calculations.
\begin{lem}\label{minimal-relevance}
In all the homogeneous cases in Section \ref{pol des}, the unital $H^\omega(\bF)$-orbit in $\Gr_G^\omega$ is minimal relevant, i.e., there is no non-zero $(H^\omega(\bF),\psi)$-equivariant sheaf supported on the boundary of its closure.

In particular, the (twisted) rank-one $(H^\omega(\bF),\psi)$-equivariant local system on $\sO_{\Gr}^0$ cleanly extends to $\Gr_{G}^\omega$.
\end{lem}

In type A, it follows from the calculation of relevant orbits in \cite[Section 3.3]{[TY]}. In the $\on{F}(4)$ case, in the coordinates of \cite[Proposition 3.3.1(b)]{[FTY]}, the indices of a $\psi$-relevant orbit satisfy $a_1\geq a_2\geq a_3\geq0$ and $a_4\geq0$. If this orbit is contained in the closure of the unital orbit, then $a_1+a_4\leq0$ by \cite[Lemma 3.4.2]{[FTY]}, thus all four indices are zero. In the orthosymplectic case, one can also obtain the desired statement from Lemmas~\ref{lem:L-orbit-indices} and \ref{lem:spherical-orbit-closure}.

Let $\unit_{\Ran}:\Ran\to\overline{\sO}_{\Gr,\Ran}$ be the unit section. We let $\IC_{\Gr}^0\in \SD_{\cG^G}^{H(\bF)_c^\omega,\psi}(\overline{\sO}_{\Gr,c})$ be the clean extension in Lemma~\ref{minimal-relevance}, corresponding to the unital orbit, and let
\[
\Vac\in \SD_{\cG^G}(\overline{\sO}_{\Gr,\Conf})
\]
be the pullback to $\overline{\sO}_{\Gr,\Conf}$ of
\[
\Av_!^{H(\bF)^\omega_{\Ran}/H(\bO)^\omega_{\Ran},\psi_{\Ran}} \bigl(\unit_{\Ran,!}\omega_{\Ran}\bigr),
\]
cf. \cite[Section 7.2.6]{[TY]}. Note that its $!$-restriction to the fiber over a single point is $\IC_{\Gr}^0$. By abuse of notation, we also denote its pullback along the inclusion $\iota_{\Gr}:\overline{\cY}_{\Gr}\hookrightarrow\overline{\sO}_{\Gr,\Conf}$ by $\Vac$.

The following lemma follows from the construction and the factorization structures in \eqref{3.19 iso} and \eqref{3.24 iso}.
\begin{lem}\label{3.4.1}
The object $\Vac$ carries a canonical factorization structure on $\overline{\sO}_{\Gr,\Conf}$. Furthermore, for any $\cF\in \SD_{\cG^G}^{H(\bF)_c^\omega,\psi}(\Fl_{G,c}^{\omega})$, the object $\sprd(\cF)$ carries a canonical factorization module structure with respect to $\Vac$. More precisely, there are canonical identifications
\begin{equation}
\begin{split}
\Vac\big|_{\overline{\sO}_{\Gr,\Conf}\times_{\Conf}(\Conf\times \Conf)_{\disj}} &\simeq (\Vac\boxtimes \Vac)\big|_{(\overline{\sO}_{\Gr,\Conf}\times \overline{\sO}_{\Gr,\Conf})_{\disj}},\\
\sprd(\cF)\big|_{\overline{\sO}_{\Fl,\Conf_{\infty\cdot c}}\times_{\Conf_{\infty\cdot c}}(\Conf\times \Conf_{\infty\cdot c})_{\disj}} &\simeq (\Vac\boxtimes \sprd(\cF))\big|_{(\overline{\sO}_{\Gr,\Conf}\times \overline{\sO}_{\Fl,\Conf_{\infty\cdot c}})_{\disj}}.
\end{split}
\end{equation}
\end{lem}

\subsubsection{Semi-infinite sheaves} Let
\[
j^{\sinf}_{\Gr}:\sS_{\Gr,\Conf}\hookrightarrow \overline{\sS}_{\Gr,\Conf}
\]
be the open embedding to its closure. Similarly to \eqref{3.34 sinf}, we have the twisted sheaves
\[
j^{\sinf}_{\Gr,!}(\omega_{\sS_{\Gr,\Conf}}),\qquad j^{\sinf}_{\Gr,*}(\omega_{\sS_{\Gr,\Conf}}),\qquad j^{\sinf}_{\Gr,!*}(\omega_{\sS_{\Gr,\Conf}}),
\]
and we also denote their restrictions to $\overline{\cY}_{\Gr}$ by the same notations.

\begin{lem}\label{3.4.2}
The objects $j^{\sinf}_{\Gr,!}(\omega_{\sS_{\Gr,\Conf}})$, $j^{\sinf}_{\Gr,*}(\omega_{\sS_{\Gr,\Conf}})$ and $j^{\sinf}_{\Gr,!*}(\omega_{\sS_{\Gr,\Conf}})$ carry canonical factorization structures. Moreover, $j^{\sinf}_{\Fl,!}(\omega_{\sS_{\Fl,\Conf_{\infty\cdot c}}})$, $j^{\sinf}_{\Fl,*}(\omega_{\sS_{\Fl,\Conf_{\infty\cdot c}}})$, and $j^{\sinf}_{\Fl,!*}(\omega_{\sS_{\Fl,\Conf_{\infty\cdot c}}})$ carry canonical factorization module structures with respect to $j^{\sinf}_{\Gr,!}(\omega_{\sS_{\Gr,\Conf}})$, $j^{\sinf}_{\Gr,*}(\omega_{\sS_{\Gr,\Conf}})$, and $j^{\sinf}_{\Gr,!*}(\omega_{\sS_{\Gr,\Conf}})$, respectively.
\end{lem}

\subsubsection{Factorization algebras and modules}
\begin{defn}\label{def fact alg}
We define
\[
\Omega_q^{\Lus}:=\Pi^d v_{\Gr,*}\bigl(\Vac\overset{!}{\otimes} j^{\sinf}_{\Gr,!}(\omega_{\sS_{\Gr,\Conf}})\bigr)[d]\in \SD_{\cG^\Lambda}(\Conf),
\]
\[
\Omega_q^{\DK}:=\Pi^d v_{\Gr,*}\bigl(\Vac\overset{!}{\otimes} j^{\sinf}_{\Gr,*}(\omega_{\sS_{\Gr,\Conf}})\bigr)[d] \in \SD_{\cG^\Lambda}(\Conf),
\]
and
\[
\Omega_q^{\sm}:=\Pi^d v_{\Gr,*}\bigl(\Vac\overset{!}{\otimes} j^{\sinf}_{\Gr,!*}(\omega_{\sS_{\Gr,\Conf}})\bigr)[d]\in\SD_{\cG^\Lambda}(\Conf).
\]
\end{defn}

Combining Lemma \ref{3.4.1} and Lemma \ref{3.4.2}, we obtain the following

\begin{cor}\label{3.5.4}
The objects $\Omega_q^{\Lus}$, $\Omega_q^{\DK}$, and $\Omega_q^{\sm}$ carry canonical structures of factorization algebras on $\Conf$. Furthermore, for any $\cF\in \SD_{\cG^G}^{H(\bF)_c^\omega,\psi}(\Fl_{G,c}^{\omega})$, the objects $F_{\Fl}^{\Lus}(\cF)$, $F_{\Fl}^{\DK}(\cF)$ and $F_{\Fl}^{\sm}(\cF)$ carry canonical factorization module structures over $\Omega_q^{\Lus}$, $\Omega_q^{\DK}$, and $\Omega_q^{\sm}$, respectively.

In particular, the functors $F_{\Fl}^{\Lus}$, $F_{\Fl}^{\DK}$ and $F_{\Fl}^{\sm}$ factor through $\Omega_q^{\Lus}\FS$, $\Omega_q^{\DK}\FS$, and $\Omega_q^{\sm}\FS$, respectively.
\end{cor}

\begin{rem}\label{rem 3.5.6}
In Theorem~\ref{main}, we mainly care about the factorization algebras obtained from the unit object $\IC_{\Gr}^0$. However, for future use, we remark that the above constructions (with spreading over $\Ran_c$) extend naively to
\begin{equation}\label{3.41+}
F_{\Gr}^{\Lus},\ F_{\Gr}^{\DK}, F_{\Gr}^{\sm}\colon\SD_{\cG^G}^{H(\bF)_c^\omega,\psi}(\Gr_{G,c}^\omega)\longrightarrow \SD_{\cG^\Lambda}(\Conf_{\infty\cdot c}).
\end{equation}
The above functors factor through $\Omega_q^{\Lus}\FS$, $\Omega_q^{\DK}\FS$, and $\Omega_q^{\sm}\FS$, respectively.

Furthermore, if we denote the natural projection by
\begin{equation}\label{fl to gr}
p\colon \Fl_G^\omega\rightarrow \Gr_G^\omega,
\end{equation}
then, we have canonical isomorphisms
\begin{equation}\label{eq:flag-Grassmannian-functors}
\begin{split}
F_{\Fl}^{\Lus}\circ p^!\simeq F_{\Gr}^{\Lus},\\
F_{\Fl}^{\DK}\circ p^!\simeq F_{\Gr}^{\DK}.
\end{split}
\end{equation}
\end{rem}

\section{Global Gaiotto category}\label{sec:global-models} The global model is needed in the case when $H_u$ is non-trivial, and in the corresponding cases listed in Section \ref{pol des}, the spherical variety is homogeneous. In this section, we will assume the curve $C=\BP^1$, and we fix a point $c\in \BP^1$. \par \ifdim\dimexpr\pagegoal-\pagetotal\relax<7\baselineskip\newpage\fi
\subsection{The global Gaiotto category}
\subsubsection{Global models} We denote by $\bX$ the affine closure $\overline{H\backslash G}^\aff$ of the homogeneous $G$-spherical variety, and by $\overset{\circ}{\bX}=H\backslash G$ the open $G$-orbit. We denote $(\bX/G)_{\omega^\rho} :=\omega^\rho\times^{Z_L}  (\bX/G)$; it has an open substack $(\overset\circ\bX/G)_{\omega^\rho}$ defined similarly. We let $\bX_{\omega^\rho}:=  \omega^\rho\times^{Z_L}\bX$.  Following \cite{[GN],[SW]}, our global model stack is defined as the algebraic stack of sections
\[\cM:= \mathrm{Sect}_{\gen}(C,(\bX/G)_{\omega^\rho} \supset (\overset\circ\bX/G)_{\omega^\rho}).\]

Similarly, we denote by $\cM_{\infty\cdot c}$ the algebraic ind-stack which classifies a pair $(\cP_G, \sigma)$, where $\cP_G$ is a $G$-bundle, and $\sigma: C-c\rightarrow \bX_{\omega^\rho}\overset{G}{\times} \cP_G$ is a section which generically lies in  $\overset{\circ}{\bX}_{\omega^\rho}\overset{G}{\times} \cP_G$.

We denote by $\cM'_{\infty\cdot c}:= \cM_{\infty\cdot c}\underset{\pt/G}{\times}\pt/B$ the fiber product, where $\cM_{\infty\cdot c}\rightarrow{\pt/G}$ is given by taking the corresponding $G$-bundle at the marked point $c$. It classifies $(\cP_G, \sigma, \epsilon)$, where $\epsilon$ is a $B$-reduction of $\cP_G$ at $c$.

\subsubsection{Global Gaiotto categories} In this section, we briefly review the global Gaiotto categories in \cite{[TY]}, which are in turn adaptations of \cite{[GN],[G3]}.

Fix an $S$-point $\bar x\in \Ran_{C-c}(S)$. Let $\cM_{\bar x,\infty\cdot c}\subset \cM_{\infty\cdot c}$ and $\cM'_{\bar x,\infty\cdot c}\subset \cM'_{\infty\cdot c}$ be the open substacks which have no defect near $\bar x$, i.e., the restriction of the section $\sigma$ to the formal disc $\cD_{\bar{x}}$ of $\bar{x}$ is contained in $\overset{\circ}{\bX}_{\omega^\rho}\overset{G}{\times} \cP_G|_{\cD_{\bar{x}}}$.

Over the formal neighborhood $\cD_{\bar x}$, a point of $\cM_{\bar x, \infty\cdot c}$ or $\cM'_{\bar x,\infty\cdot c}$ determines an $H^\omega$-bundle $\cP_H$ on $\cD_{\bar x}$. Then, we denote by $\tilde{\cM}_{\bar x, \infty\cdot c}$ and $\tilde{\cM}_{\bar x, \infty\cdot c}'$ the $H(\bO)^\omega_{\bar{x}}$-torsors over  $\cM_{\bar x, \infty\cdot c}$ and $\cM'_{\bar x,\infty\cdot c}$ obtained by adding an isomorphism between the $H^\omega$-bundle $\cP_H$ and $\cP_H^\omega:= \omega^\rho\overset{Z_L}{\times}H$ on $\cD_{\bar x}$, respectively. The action of $H(\bO)^\omega_{\bar{x}}$ extends to an action of $H(\bF)^\omega_{\bar{x}}$.

\begin{defn}
The Iwahori global Gaiotto category is defined as
\begin{equation}
\SD_{\cG^G}^{H,\Hecke,\psi}(\cM'_{\infty\cdot c}):= \lim_{\bar{x}\in \Ran_{C-c}(S)} \SD_{\cG^G}^{H(\bF)^\omega_{\bar{x}},\psi_{\bar{x}}}(\tilde{\cM}'_{\bar x,\infty\cdot c}),
\end{equation}
and the spherical global Gaiotto category is defined as
\begin{equation}
\SD_{\cG^G}^{H,\Hecke,\psi}(\cM_{\infty\cdot c}):= \lim_{\bar{x}\in \Ran_{C-c}(S)} \SD_{\cG^G}^{H(\bF)^\omega_{\bar{x}},\psi_{\bar{x}}}(\tilde{\cM}_{\bar x,\infty\cdot c}).
\end{equation}
Here, the transition functors are defined in \cite[Sections 4.3--4.4]{[TY]}; these categories are functorial both for inclusions $\bar x\subset \bar x'$ and for base change in the test scheme $S$. Therefore, we may take limits over $\Ran_{C-c}$. The forgetful functors use the componentwise shift $[-d_{\Gr}]$, where $d_{\Gr}=\dim\Bun_H^\omega$ as in Remark~\ref{rem:local-functor-degree}.
\end{defn}

The definition above admits a more Hecke-equivariant style description as in \cite{[GN]}. Consider the following universal Hecke correspondence diagram,
\begin{equation}\label{univ heck}
\xymatrix{&\Hecke_{\Gr,\Ran_{C-c}}\ar[ld]_{\overleftarrow h}\ar[rd]^{\overrightarrow h}&\\
\cM_{\infty\cdot c}&&\cM_{\infty\cdot c} }
\end{equation}
Here, an $S$-point of the algebraic stack $\Hecke_{\Gr,\Ran_{C-c}}$ consists of the data
\begin{itemize}
\item an $S$-point $\bar{x}$ in $\Ran_{C-c}$,
\item an $S$-point $(\cP_G,\sigma)$ in $\cM_{\bar{x},\infty\cdot c}$ (in particular, it determines an $H^\omega$-bundle $\cP_H$ on $\cD_{\bar{x}}$),
\item an $H^\omega$-bundle $\tilde{\cP}_H$ on $\cD_{\bar{x}}$
\item an isomorphism $\tilde{\cP}_H|_{\overset{\circ}{\cD}}\simeq {\cP}_H|_{\overset{\circ}{\cD}}$.
\end{itemize}
The projection $\overleftarrow h$ just sends the above data to $(\cP_G,\sigma)$, while the projection $\overrightarrow h$ sends them to $(\cP_G',\sigma')$, where $(\cP_G',\sigma')$ is obtained by gluing $\sigma|_{(C-c)\times S-\bar{x}}$ and $\tilde{\cP}_H$ along the punctured formal neighborhood of $\bar{x}$.

The canonical additive character on $H(\bF)_{\Ran}^\omega$ gives rise to a local system $\psi_{\Ran_{C-c}}$ on the universal Hecke stack over $\Ran_{C-c}$. Also, note that the twisting $\cG^G$ is compatible with the Hecke correspondence.

Then, an object of the generic Hecke equivariant category on $\cM_{\infty\cdot c}$ is just a sheaf $\cF\in\SD_{\cG^G}(\cM_{\infty\cdot c})$, together with the isomorphism
\begin{equation}\label{univ heck iso}
\varepsilon\colon \overleftarrow h^{\,!}(\cF)\xrightarrow{\sim} \bigl(\overrightarrow h^{\,!}(\cF)\overset{!}{\otimes} \psi_{\Ran_{C-c}}\bigr)[\overrightarrow h^*d_{\Gr}-\overleftarrow h^*d_{\Gr}]
\end{equation}
with higher homotopy coherence.

Indeed, it follows from the definition that for any $S$-point $\bar{x}\in \Ran_{C-c}(S)$, the restriction of the Hecke correspondence diagram \eqref{univ heck} gives rise to
\begin{equation}\label{x heck}
\xymatrix{&\Hecke_{\Gr,\bar{x}}\ar[ld]_{\overleftarrow h}\ar[rd]^{\overrightarrow h}&\\
\cM_{\bar{x},\infty\cdot c}&&\cM_{\bar{x},\infty\cdot c}. }
\end{equation}
Then, by \eqref{univ heck iso}, we have
\begin{equation}
\overleftarrow h^{\,!}(\cF)\xrightarrow{\sim} \bigl(\overrightarrow h^{\,!}(\cF)\overset{!}{\otimes} \psi_{\bar x}\bigr)[\overrightarrow h^*d_{\Gr}-\overleftarrow h^*d_{\Gr}].
\end{equation}
Note that an object $\cF\in \SD_{\cG^G}(\cM_{\bar{x},\infty\cdot c})$ together with the above isomorphism is just equivalent to an object in $\SD_{\cG^G}^{H(\bF)^\omega_{\bar{x}},\psi_{\bar{x}}}(\tilde{\cM}_{\bar x,\infty\cdot c})$.

For the Iwahori version of the global category (note that the Iwahori global model differs from the spherical one only at the marked point), one also has a similar description.

We denote by
\[
\IC_{\Gr,\glob}^{0}\in \SD_{\cG^G}^{H,\Hecke,\psi}(\cM)
\]
the perverse generator on the unital stratum; equivalently, it is the clean extension of the twisted generic Hecke equivariant perverse local system from the open stratum $\Bun_H^\omega$ of $\cM$ (the stack of $H^\omega$-bundles). Here, on a smooth component of dimension $r$, the rank-one generator is $\omega[-r]$ with even coefficients.

\par \ifdim\dimexpr\pagegoal-\pagetotal\relax<7\baselineskip\newpage\fi
\subsubsection{Dimensions}\label{sec:global-dimension-normalization}
\begin{rem}\label{rem:local-functor-degree}
Let $g=g(C)$ be the genus of $C$. We denote by $\Bun_B^{\deg=\nu}$ the component on which the associated $T$-bundle has degree $\nu$, and let
\[
\Bun_B^{\omega,\lambda}:=\Bun_B^{\deg=(2g-2)\rho-\lambda}, \qquad d_B(\lambda):=\dim\Bun_B^{\omega,\lambda}.
\]
Since $\check\rho_G$ is taken with respect to $B^-$, by the Riemann--Roch theorem applied to the filtration of $\Lie B$ by root spaces, we have
\[
d_B(\lambda)=(g-1)\dim B+ \langle2\check\rho_G,(2g-2)\rho-\lambda\rangle =d_B(0)-\deg(\lambda).
\]
Note that for the embedding in Section~\ref{trans D}, we have $H\cap B=\{1\}$ by Lemma~\ref{lem:trivial-H-B-intersection}, and the differential of the multiplication map identifies the $\omega^\rho$-twists of $\Lie H\oplus\Lie B$ and $\Lie G$. Taking the Euler characteristics, we obtain
\[
d_H^0+d_B(0)=\dim\Bun_G,
\]
where $d_H^0$ denotes the dimension of the connected component of $\Bun_H^\omega$ containing $\cP_H^\omega$.

In type A, we denote $a(\lambda)=\sum_i\lambda_i^{(1)}$. Since the determinant character of $H_r=\GL_M$ on $\Lie(H_u)$ is $\det^{-(N-M-1)}$, by the Riemann--Roch theorem, we have $\dim\Bun_H^\omega=d_H^0-(N-M-1) m$ on the component where the $H_r$-bundle has determinant degree $-m$. Note that in the fiber of configuration degree $\lambda$, the underlying $G_1$-bundle has determinant degree $-a(\lambda)$, relative to the $G_1$-bundle induced by $\cP_H^\omega$. Hence,
\[
\dim\Bun_G-\dim\Bun_H^\omega-d_B(\lambda) =\deg(\lambda)+(N-M-1)a(\lambda)=d(\lambda),
\]
and thus, on the open SW Zastava space in Section \ref{5.3.1}, we have
\begin{equation}\label{zas dim}
    \dim\bigl(\Bun_H^\omega\underset{\Bun_G}{\times'}\Bun_B^{\omega,\lambda}\bigr)=-d(\lambda).
\end{equation}
Here, the prime denotes the open locus of generically transverse reductions.

In the other homogeneous cases, the adjoint determinant character of $H_r$ is trivial, thus $\dim\Bun_H^\omega=d_H^0$ and $d=\deg$, and the \eqref{zas dim} holds.

\end{rem}

\subsection{Local-global comparison}
\subsubsection{The comparison functors} Let $x_0\in \overset{\circ}{\bX}=H\backslash G$ be a chosen base point.

For a point $(\cP_G,\alpha,\epsilon)\in \Fl_{G,c}^{\omega}$, the trivialization $\alpha$ over $C-c$ identifies $\cP_G|_{C-c}$ with $\cP_G^\omega|_{C-c}$, and hence the point $x_0$ determines a section
\begin{equation}
\sigma:C-c\rightarrow \overset{\circ}{\bX}_{\omega^\rho}\times^G \cP_G.
\end{equation}
This gives a morphism
\begin{equation}
\pi_{\Fl}:\Fl_{G,c}^{\omega}\longrightarrow \cM'_{\infty\cdot c}.
\end{equation}

Similarly, we have
\begin{equation}
\pi_{\Gr}:\Gr_{G,c}^{\omega}\longrightarrow \cM_{\infty\cdot c}.
\end{equation}
Its restriction gives
\begin{equation}
\pi_{\Gr}:\overline{\sO}_{\Gr,c}\longrightarrow \cM.
\end{equation}

As in the Drinfeld compactification case (cf. \cite[7.1.1]{[FGV]}), one can show that the above maps $\pi_\Fl$ and $\pi_{\Gr}$ are formally smooth.

Now, let
\[
d_{\Gr}:=\dim \Bun_H^\omega,\qquad d_{\Fl}:=d_{\Gr}+\dim(G/B).
\]
In type A, on a component of the global model where the underlying $G_1$-bundle has determinant degree $-m$ relative to $\cP_{G_1}^0$, we take $d_{\Gr}=d_H^0-(N-M-1)m$, as in Remark~\ref{rem:local-functor-degree}. In the other homogeneous cases, $d_{\Gr}=d_H^0$. Note that in the Iwahori version, both $\dim\Bun_G'$ and $d_{\Fl}$ acquire the same additional summand $\dim(G/B)$, while $\Bun_B^e\simeq\Bun_B$ (these stacks are defined in the next subsection).

The spherical version of the following theorem is shown in \cite[Theorem 4.5.3]{[TY]}, and the Iwahori version is a direct analogue.
\begin{thm}\label{prop:global-local-strata}
Taking pullback along $\pi_{\Fl}$ induces an equivalence
\begin{equation}
\begin{split}
E_{\Fl}:=\pi_{\Fl}^{!}[d_{\Fl}]: \SD_{\cG^G}^{H,\Hecke,\psi}(\cM'_{\infty\cdot c}) \xrightarrow{\sim} \SD_{\cG^G}^{H(\bF)_c^\omega,\psi}(\Fl_{G,c}^{\omega}).
\end{split}
\end{equation}
Similarly, pullback along $\pi_{\Gr}$ induces an equivalence
\begin{equation}
E_{\Gr}:=\pi_{\Gr}^{!}[d_{\Gr}]: \SD_{\cG^G}^{H,\Hecke,\psi}(\cM_{\infty\cdot c}) \xrightarrow{\sim} \SD_{\cG^G}^{H(\bF)_c^\omega,\psi}(\Gr_{G,c}^{\omega}).
\end{equation}
\end{thm}

\begin{rem}\label{rem:component-normalization}
Via $E_{\Fl}$ and $E_{\Gr}$, we define the perverse (or na\"{i}ve) $t$-structure of the local Gaiotto category by transporting the global perverse $t$-structures. Note that for dominant $\lambda$, the irreducible equivariant ``perverse'' sheaf on the orbit through $t^\lambda$ has $!$-stalk $\se[-d(\lambda)]$ at $t^\lambda$.
\end{rem}
\subsubsection{Defect degree} The global model $\cM_{\infty\cdot c}$ (resp. $\cM'_{\infty\cdot c}$) is stratified according to the defect degree.

(1). At the marked point $c\in \BP^1$. Given any geometric point $(\cP_G,\sigma)$, according to the generic condition, the restriction of $\sigma$ to $\overset{\circ}{\cD}_c$ lands in $\overset{\circ}{\bX}_{\omega^\rho}\overset{G}{\times} \cP_G|_{\overset{\circ}{\cD}_c}$. A trivialization of $\cP_G|_{\cD_c}$ (resp. one compatible with the Iwahori structure at $c$) determines a point in $\overset{\circ}{\bX}_{\omega^\rho}(\bF)_c$. Changing the trivialization changes the resulting point by right multiplication by $G(\bO)^\omega_c$ (resp. $I^\omega_c$), and thus, $(\cP_G,\sigma)$ determines a point in
\begin{equation}
|\overset{\circ}{\bX}_{\omega^\rho}(\bF)_c/G(\bO)^\omega_c|=|H(\bF)_c^\omega\backslash\Gr_{G,c}^\omega|
\end{equation}
(resp.
\begin{equation}
|\overset{\circ}{\bX}_{\omega^\rho}(\bF)_c/I^\omega_c|=|H(\bF)_c^\omega\backslash\Fl_{G,c}^\omega|.
\end{equation})
We will call the corresponding class in $|H(\bF)_c^\omega\backslash\Gr_{G,c}^\omega|$ (resp., $|H(\bF)_c^\omega\backslash\Fl_{G,c}^\omega|$) the defect degree at $c$. Here, in families, a fixed defect degree means that the local datum is fppf locally isomorphic to the chosen representative of this class.

(2). One can also consider the degree of the defect at $x\in C-c$. The corresponding defect degree belongs to $|(\overset{\circ}{\bX}_{\omega^\rho}(\bF)_x\cap \bX_{\omega^\rho}(\bO)_x)/G(\bO)^\omega_x|$, while only the unital stratum is relevant.

\subsubsection{Stratification of the global models} For an $H(\bF)_c^\omega$-orbit $\sO_{\Fl,\tilde w}\subset \Fl_{G,c}^{\omega}$, let
\[
\overset{\circ}{\cM}_{\Fl,\tilde w}\subset \cM'_{\infty\cdot c}
\]
be the corresponding good-elsewhere stratum, i.e., the locally closed substack consisting of points whose defect at $c$ is ${\tilde w}$ and which have no defect away from $c$.  Then, $E_{\Fl}$ induces an equivalence between the subcategories on the strata $\overset{\circ}{\cM}_{\Fl,\tilde w}$ and $\sO_{\Fl,\tilde w}$. Similarly, in the spherical case, for an $H(\bF)_c^\omega$-orbit $\sO_{\Gr}^{\lambda}\subset\Gr_{G,c}^{\omega}$, we denote by $\cM_{=\lambda\cdot c}\subset\cM_{\infty\cdot c}$ the corresponding good-elsewhere stratum. In particular,
\begin{equation}
E_{\Gr}(\IC_{\Gr,\glob}^{0})\simeq \IC_{\Gr}^0.
\end{equation}

\subsection{Global Zastava spaces and global functors}\label{sec:global-zastava-functors}
\subsubsection{Zastava spaces}\label{5.3.1} Let $\Bun_G$ be the stack of principal $G$-bundles on $C$, and let
\begin{equation}
\Bun_G':=\Bun_G\underset{\pt/G}{\times}\pt/B
\end{equation}
be the stack of $G$-bundles equipped with a $B$-reduction at $c$.

Let $\Bun_B$ be the stack of $B$-bundles on $C$. We define $\Bun_B'$ to be the stack classifying a $B$-bundle $\cP_B$ on $C$ together with a second $B$-reduction of the induced $G$-bundle at $c$. Thus, $\Bun_B'$ has two maps to $\pt/B$, and hence there is a natural map
\begin{equation}
\Bun_B'\longrightarrow \pt/B\underset{\pt/G}{\times}\pt/B\simeq B\backslash G/B.
\end{equation}
We denote by
\begin{equation}
\Bun_B^w\subset \Bun_B'
\end{equation}
the preimage of the Bruhat cell corresponding to $w\in W^\fin$. In particular, $\Bun_B^e\simeq\Bun_B$ is the closed stratum in $\Bun_B'$ on which the two reductions agree.  Let $\overline{\Bun}_B$ and $\overline{\Bun}_B'$ be the corresponding Drinfeld compactifications of $\Bun_B$ and $\Bun_B'$, respectively. Also, we denote by $(\Bun_B^w)^{\omega,\lambda}$ the component of degree $(2g-2)\rho-\lambda$.

In \eqref{3.26 zas} and \eqref{3.44 zas}, we defined the SW Zastava spaces in local terms. It is well-known (and more classical) that the SW Zastava space can be described in global terms, see \cite{[GL],[BFT1]}, or \cite{[SW]}.

To be more precise,
\begin{lem}
There are canonical isomorphisms of algebraic stacks (in fact, schemes) of ind-finite type
\begin{equation}
\begin{split}
\overline{\cY}_{\Gr} &\simeq\cM\underset{\Bun_G}{\times'}\overline{\Bun}_B,\\
\overline{\cY}_{\Fl,\infty\cdot c} &\simeq\cM'_{\infty\cdot c}\underset{\Bun_G'}{\times'}\overline{\Bun}_B',
\end{split}
\end{equation}
where $\times'$ means that the generalized $H$-reduction and the $B$-reduction are generically transverse. \footnote{For the first fiber product, we take the components of degree $\alpha\in\Lambda^{\mathrm{neg}}\setminus\{0\}$ to match the convention that the degree of $\Conf$ is non-zero.} These isomorphisms restrict to the following substacks
\begin{equation}
\begin{split}
{\cY}_{\Gr} &\simeq\cM\underset{\Bun_G}{\times'}{\Bun}_B,\\
{\cY}_{\Fl,\infty\cdot c} &\simeq\cM'_{\infty\cdot c}\underset{\Bun_G'}{\times'}{\Bun}_B^e.
\end{split}
\end{equation}
\end{lem}

\subsubsection{Projection maps} We have the following commutative diagrams of natural projections
\begin{equation}
\xymatrix{ \overline{\cY}_{\Gr}\ar[r]^{q_{\Gr}}\ar[d]_{p_{\Gr}} & \overline{\Bun}_B\ar[d]^{\mathfrak p_{\Gr}} \\
\cM\ar[r]_{\mathfrak{q}_{\Gr}} & \Bun_G, } \qquad \xymatrix{ \overline{\cY}_{\Fl,\infty\cdot c}\ar[r]^{q_{\Fl}}\ar[d]_{p_{\Fl}} & \overline{\Bun}_B'\ar[d]^{\mathfrak p_{\Fl}} \\
\cM'_{\infty\cdot c}\ar[r]_{\mathfrak{q}_{\Fl}} & \Bun_G'. }
\end{equation}

\subsubsection{Global IC sheaves} One has the natural open and locally closed embeddings, respectively,
\begin{equation}
j_{\Gr}:\Bun_B\hookrightarrow \overline{\Bun}_B,\qquad j_{\Fl}^w:\Bun_B^w\hookrightarrow \overline{\Bun}_B'\quad(w\in W^{\fin}).
\end{equation}

We consider the ratio of the pullback of the gerbe $\cG^T$ on $\Bun_T$ and the pullback of the gerbe $\cG^G$ on $\Bun_G$. Denote the resulting gerbe on $\overline{\Bun}_B$ and $\overline{\Bun}_B'$ by $\cG^{G,T,\ratio}$. By construction, the restriction of $\cG^{G,T,\ratio}$ to $\Bun_B$ and to $\Bun_B'$ is canonically trivial.

In particular, it makes sense to consider the twisted constant perverse sheaves on $\Bun_B$ and $\Bun_B^w$. For the global semi-infinite sheaves, we let 
\[
\IC_{\Bun_B}|_{\Bun_B^{\omega,\lambda}} =\omega_{\Bun_B^{\omega,\lambda}}[-d_B(\lambda)],
\]
where $d_B(\lambda)=\dim\Bun_B^{\omega,\lambda}$. Similarly, on a component of $\Bun_B^w$ of dimension $r$, we let $\IC_{\Bun_B^w}:=\omega[-r]$. We will use their $!$, $*$, or $!*$-extensions, and denote them by
\[
\fitmath{j_{!,\Gr}(\IC_{\Bun_B}):=(j_{\Gr})_!(\IC_{\Bun_B}), \qquad j_{*,\Gr}(\IC_{\Bun_B}):=(j_{\Gr})_*(\IC_{\Bun_B}), \qquad j_{!*,\Gr}(\IC_{\Bun_B}),}
\]
and
\[
\fitmath{j_{!,\Fl}^w(\IC_{\Bun_B^w}):=(j_{\Fl}^w)_!(\IC_{\Bun_B^w}), \qquad j_{*,\Fl}^w(\IC_{\Bun_B^w}):=(j_{\Fl}^w)_*(\IC_{\Bun_B^w}), \qquad j_{!*,\Fl}^w(\IC_{\Bun_B^w}).}
\]

The following lemma explains the relations between the local semi-infinite sheaves and the global semi-infinite sheaves, and it is (the metaplectic version of) \cite[Theorem 3.3.3]{[G6]}, or can be proved with the same argument.
\begin{lem}\label{semi-inf compatble}
Under the $!$-pullback along the natural projection
\begin{equation}
\overline{\sS}_{\Gr,\Conf}\rightarrow \overline{\Bun}_B,
\end{equation}
the global semi-infinite objects $j_{!,\Gr}(\IC_{\Bun_B})$, $j_{*,\Gr}(\IC_{\Bun_B})$, and $j_{!*,\Gr}(\IC_{\Bun_B})$ go to the local semi-infinite objects $j^{\sinf}_{\Gr,!}(\omega_{\sS_{\Gr,\Conf}})$, $j^{\sinf}_{\Gr,*}(\omega_{\sS_{\Gr,\Conf}})$, $j^{\sinf}_{\Gr,!*}(\omega_{\sS_{\Gr,\Conf}})$, respectively, with the additional shift $[-d_B(\lambda)]$ on configuration degree $\lambda$, where $d_B(\lambda)$ is the dimension of the corresponding component of $\Bun_B$.

Similarly for the Iwahori version.
\end{lem}

\subsubsection{The global functors}
\begin{defn}
We define the global functors
\begin{equation}
\begin{split}
F_{\Gr}^{\Lus,\glob}:\SD_{\cG^G}^{H,\Hecke,\psi}(\cM)\longrightarrow \SD_{\cG^\Lambda}(\Conf),\\
\cF\longmapsto \Pi^d v_{\Gr,*}\bigl(p_{\Gr}^!(\cF)\overset{!}{\otimes}q_{\Gr}^! j_{!,\Gr}(\IC_{\Bun_B})\bigr)[\dim \Bun_G],
\end{split}
\end{equation}
and
\begin{equation}
\begin{split}
F_{\Fl}^{\Lus,\glob}:\SD_{\cG^G}^{H,\Hecke,\psi}(\cM'_{\infty\cdot c})\longrightarrow \SD_{\cG^\Lambda}(\Conf_{\infty\cdot c}),\\
\cF\longmapsto \Pi^d v_{\Fl,*}\bigl(p_{\Fl}^!(\cF)\overset{!}{\otimes}q_{\Fl}^! j_{!,\Fl}^e(\IC_{\Bun_B^e})\bigr)[\dim \Bun_G'].
\end{split}
\end{equation}
Similarly, we define $F_{\Gr}^{\DK,\glob}$ and $F_{\Fl}^{\DK,\glob}$ (resp. $F_{\Gr}^{\sm,\glob}$ and $F_{\Fl}^{\sm,\glob}$) by replacing $j_{!,\Gr}(\IC_{\Bun_B})$ and $j_{!,\Fl}^e(\IC_{\Bun_B^e})$ with $j_{*,\Gr}(\IC_{\Bun_B})$ and $j_{*,\Fl}^e(\IC_{\Bun_B^e})$ (resp. $j_{!*,\Gr}(\IC_{\Bun_B})$ and $j_{!*,\Fl}^e(\IC_{\Bun_B^e})$).
\end{defn}

\begin{rem}
Similarly to Remark \ref{rem 3.5.6}, the global functors of the spherical version also naturally extend to functors from $\SD_{\cG^G}^{H,\Hecke,\psi}(\cM_{\infty\cdot c})$ to $\SD_{\cG^\Lambda}(\Conf_{\infty\cdot c})$.
\end{rem}

Note that the projection $v_{\Fl}$ is ind-proper, cf. \cite[(7.4.2)]{[Yan]}.

\begin{lem}\label{lem:proper-support-global-functors}
Let $\cF$ be a global Gaiotto object supported on finitely many global strata, and let $K$ be one of the three global semi-infinite objects on $\overline{\Bun}_B'$ used in the above functors. Over each component of $\Conf_{\infty\cdot c}$, the support of
\[
p_{\Fl}^!\cF\overset{!}{\otimes}q_{\Fl}^!K
\]
is contained in a closed subscheme of $\overline{\cY}_{\Fl,\infty\cdot c}$ proper and of finite type over that component (of the configuration space).
\end{lem}
\begin{proof}
Since $\cF$ is supported on finitely many global strata, the poles at $c$ are bounded. Then, the claim follows from the same argument as in \cite[Proposition 4.1.1 and Lemma 4.1.2]{[SW]}, with the additional flag at $c$.
\end{proof}

\subsubsection{Comparison with the local functors} According to the constructions and Lemma \ref{semi-inf compatble}, we have the following
\begin{prop}\label{prop:local-global-functor-comparison}
For $?=\Lus,\DK,\sm$, there are canonical isomorphisms of functors
\[
F_{\Gr}^{?}\circ E_{\Gr} \xrightarrow{\sim} F_{\Gr}^{?,\glob}, \qquad F_{\Fl}^{?}\circ E_{\Fl} \xrightarrow{\sim} F_{\Fl}^{?,\glob}.
\]
Here, $F_{\Gr}^{?,\glob}$ denotes the extension to $\cM_{\infty\cdot c}$ described above.
\end{prop}
\begin{proof}
By Theorem~\ref{prop:global-local-strata} and Lemma~\ref{semi-inf compatble}, it is enough to compare the normalizations. On the component of configuration degree $\lambda$, by Remark~\ref{rem:local-functor-degree}, we have
\[
\dim\Bun_G-d_B(\lambda)=d_H^0+\deg(\lambda)=d_{\Gr}+d(\lambda).
\]
Here, $d_{\Gr}=d_H^0-\delta(\lambda)$ on the corresponding component of $\Bun_H^\omega$. Thus, the cohomological shifts agree, and both functors have the same parity factor $\Pi^{d(\lambda)}$. In the case of the affine flags, both sides acquire the additional summand $\dim(G/B)$. Note that the above comparison works for all three kernels, and it is compatible with the spreading and the factorization structures.
\end{proof}

As an immediate consequence, by combining with Corollary \ref{3.5.4}, we have
\begin{cor}\label{cor:omega-from-global-vacuum}
There are canonical isomorphisms
\[
\fitmath{\Omega_q^{\Lus}\simeq F_{\Gr}^{\Lus,\glob}(\IC_{\Gr,\glob}^{0}), \qquad \Omega_q^{\DK}\simeq F_{\Gr}^{\DK,\glob}(\IC_{\Gr,\glob}^{0}), \qquad \Omega_q^{\sm}\simeq F_{\Gr}^{\sm,\glob}(\IC_{\Gr,\glob}^{0}).}
\]

Also, the functors $F_{\Fl}^{\Lus,\glob}$ and $F_{\Gr}^{\Lus,\glob}$, $F_{\Fl}^{\DK,\glob}$ and $F_{\Gr}^{\DK,\glob}$, and $F_{\Fl}^{\sm,\glob}$ and $F_{\Gr}^{\sm,\glob}$ factor through $\Omega_q^{\Lus}\FS$, $\Omega_q^{\DK}\FS$, and $\Omega_q^{\sm}\FS$, respectively.
\end{cor}

\subsubsection{Dimensions of the global strata}
\begin{rem}\label{rem:global-stratum-degree}
Let $\lambda$ be dominant. It follows by a direct calculation that, on each component of $\cM_{=\lambda\cdot c}$, we have
\[
\dim\cM_{=\lambda\cdot c}-d_{\Gr}=d(\lambda),
\]
which is compatible with the normalization in Remark~\ref{rem:component-normalization}. Now, let us consider the SW Zastava space
\[
\cY_{\Gr,=\lambda\cdot c}^{\mu}:=\cM_{=\lambda\cdot c}\times'_{\Bun_G}\Bun_B^{\omega,\mu},
\]
where $\Bun_B^{\omega,\mu}$ is as in Remark~\ref{rem:local-functor-degree}. As in the proof of Proposition~\ref{prop:local-global-functor-comparison}, we have $\dim\Bun_G-d_{\Gr}-d_B(\mu)=d(\mu)$ by Remark~\ref{rem:local-functor-degree}, i.e., the normalized semi-infinite kernel has degree $d(\mu)$ on $\cY_{\Gr,=\lambda\cdot c}^{\mu}$. Hence, over each component of $\cM_{=\lambda\cdot c}$, the corresponding component of $\cY_{\Gr,=\lambda\cdot c}^{\mu}$ has dimension
\[
d(\lambda)-d(\mu)=d(\lambda-\mu).
\]
In particular, for $\lambda=0$, we have $\cM_{=0\cdot c}=\Bun_H^\omega$, and this is compatible with \eqref{zas dim}.
\end{rem}

\subsection{Duality}\label{sec:global-duality} While working with the local model, the definition of the duality of the Gaiotto category is not obvious due to the renormalization and the infinite-dimensional orbits; the duality of the Gaiotto categories is easier to approach, and more convenient to work with, in the global model.

Let
\[
\SD_{\cG^G}^{H,\Hecke,\psi}(\cM'_{\infty\cdot c})^c,\qquad \SD_{\cG^G}^{H,\Hecke,\psi}(\cM_{\infty\cdot c})^c
\]
denote the full subcategories of objects which have support on finitely many strata, and whose $!$-stalks are bounded and finite-dimensional on each stratum. In the $D$-module setting, these are the categories of \textit{coherent} objects in $\SD_{\cG^G}^{H,\Hecke,\psi}(\cM'_{\infty\cdot c})$ and $\SD_{\cG^G}^{H,\Hecke,\psi}(\cM_{\infty\cdot c})$, in the sense of \cite[\S 7.3.1]{[DG]}.

Under the local-global comparison in Theorem \ref{prop:global-local-strata}, they correspond to the full subcategories of locally compact objects $\SD_{\cG^G}^{H^\omega(\bF),\psi,\locc}(\Fl_G^\omega) $ and $\SD_{\cG^G}^{H^\omega(\bF),\psi,\locc}(\Gr_G^\omega)$. See also \cite[Remark 5.3.9 and Proposition 5.3.10]{[TY]}.

With the same argument as in \cite[Remark 4.6.8]{[TY]}, one can show that
\begin{lem}
The Verdier duality on the ambient twisted sheaf categories preserves the generic Hecke equivariant structure, while replacing the twisting by its inverse and $\psi$ by $-\psi$.
\end{lem}
Hence, we have duality functors
\begin{equation}\label{4.26 dual}
\BD: \SD_{\cG^G}^{H,\Hecke,\psi}(\cM'_{\infty\cdot c})^c \longrightarrow \Bigl(\SD_{(\cG^G)^{-1}}^{H,\Hecke,-\psi}(\cM'_{\infty\cdot c})^c\Bigr)^{\opp},
\end{equation}
and
\begin{equation}\label{4.27 dual}
\BD: \SD_{\cG^G}^{H,\Hecke,\psi}(\cM_{\infty\cdot c})^c \longrightarrow \Bigl(\SD_{(\cG^G)^{-1}}^{H,\Hecke,-\psi}(\cM_{\infty\cdot c})^c\Bigr)^{\opp}.
\end{equation}

As a direct corollary,
\begin{cor}
The ind-complete categories $\SD_{\cG^G}^{H^\omega(\bF),\psi,\ren}(\Fl_G^\omega):=\Ind(\SD_{\cG^G}^{H^\omega(\bF),\psi,\locc}(\Fl_G^\omega))$ and $\SD_{\cG^G}^{H^\omega(\bF),\psi,\ren}(\Gr_G^\omega):=\Ind(\SD_{\cG^G}^{H^\omega(\bF),\psi,\locc}(\Gr_G^\omega))$ are dual to the corresponding renormalized categories with twisting $(\cG^G)^{-1}$ and character $-\psi$. On compact objects, the resulting anti-equivalences are induced by \eqref{4.26 dual} and \eqref{4.27 dual}, and we will denote them by the same notation.
\end{cor}

\section{Standard objects in the Iwahori Gaiotto category}\label{sec:gaiotto-standards}

Now, we are going to show that the local functor $F_{\Fl}^{\Lus}$ is an equivalence. We will adapt the strategy of \cite{[Yan]}. In particular, we first recall the twisted BMW sheaves on the enhanced affine flag, and then use them to define the standard objects in the Iwahori Gaiotto category. Recall that, on the constructible side, when we talk about dominance of coweights, we always mean dominance with respect to the \emph{opposite} Borel subgroup $B^-$, that is, a coweight $\lambda$ is dominant if
\begin{equation}U(\bO)\subset t^{\lambda} U(\bO)t^{-\lambda}.\end{equation}
Equivalently,
\begin{equation}
B^\omega(\bO)\subset t^{\lambda} I^\omega t^{-\lambda}.
\end{equation}

\subsection{Twisted BMW sheaves} Let $\widetilde{\Fl}_G:=G(\bF)_c^\omega/I^0$ be the enhanced affine flags, where $I^0\subset I^\omega$ is the pro-unipotent radical of the Iwahori subgroup. The exact sequence
\begin{equation}
1\longrightarrow I^0\longrightarrow I^\omega\longrightarrow T^\omega\longrightarrow 1
\end{equation}
gives a canonical right action of $T^\omega=I^\omega/I^0$ on $\widetilde{\Fl}_G$.

Recall that the quadratic form $q$ on $\Lambda$ determines a symmetric bilinear form $b$. In the $D$-module setting,
\begin{equation}
b(\lambda_1,\lambda_2):=q(\lambda_1+\lambda_2)-q(\lambda_1)-q(\lambda_2),
\end{equation}
while in the Betti and $\ell$-adic settings one uses the corresponding multiplicative formula. We denote by $b_{\lambda}:=b(\lambda,-)$ the associated multiplicative local system on $T^\omega$ (a Kummer sheaf in the corresponding sheaf theory), and by the same symbol its pullback to $I^\omega$.

For $\lambda,\mu\in \Lambda$, we define
\begin{equation}
\SD^{I^\omega,b_{\lambda}}_{\cG^G}(\widetilde{\Fl}_G)_{b_{\mu}} := \SD_{\cG^G}(I^\omega,b_{\lambda}\backslash \widetilde{\Fl}_G/T^\omega,b_{\mu}),
\end{equation}
the category of left $(I^\omega,b_{\lambda})$-equivariant and right $(T^\omega,b_{\mu})$-equivariant $\cG^G$-twisted sheaves on $\widetilde{\Fl}_G$. Similarly, we define
\begin{equation}
\SD_{\cG^G}^{I^0}(\widetilde{\Fl}_G)_{b_{\mu}} := \SD_{\cG^G}(I^0\backslash \widetilde{\Fl}_G/T^\omega,b_{\mu}).
\end{equation}
When ${\mu}=0$, the category is identified with $\SD_{\cG^G}^{I^0}({\Fl}^\omega_G)$.

\subsubsection{} We denote by $W_G^{\mathrm{ext}}$ the extended affine Weyl group of $G$, and by $l$ its length function.

For $w\in W_G^{\mathrm{ext}}$, let $\widetilde{\Fl}_G^{w}\subset \widetilde{\Fl}_G$ be the preimage of the Iwahori orbit
\begin{equation}
\Fl_G^{w}:=I^\omega w I^\omega/I^\omega\subset \Fl_G^\omega.
\end{equation}
Since the group $I^0\times T^\omega$ acts transitively on $\widetilde{\Fl}_G^w$ with pro-unipotent stabilizers, after fixing compatible trivializations of the twisting and the monodromy, and choosing an identification
\begin{equation}\label{5.5 id}
\SD_{\cG^G}^{I^0}(\widetilde{\Fl}_G^{w})_{b_{\mu}}\simeq \SVect
\end{equation}
which preserves the cohomological degree of the $!$-stalk, we denote by $(c_{w})_{\mu}$ the corresponding generator (with even coefficients). Different choices of \eqref{5.5 id} change $(c_{w})_{\mu}$ by tensoring with an even line in cohomological degree $0$. If ${\mu}=0$, then $(c_{w})_0$ is simply the $\cG^G$-twisted dualizing sheaf on this stratum (with even coefficients), together with its $I^0$-equivariant and right $T^\omega$-equivariant structure.

We define $(J_{w,!})_{\mu}$ and $(J_{w,*})_{\mu}$ to be the $!$- and $*$-extensions of $(c_{w})_{\mu}[-l(w)]$ along the locally closed embedding $\widetilde{\Fl}_G^{w}\hookrightarrow \widetilde{\Fl}_G$.

As in \cite[Section 5.3]{[Yan]}, there is a convolution product of the twisted BMW sheaves on $\widetilde{\Fl}_G$, and the same argument as in \cite[Lemma 8]{[AB]} and \cite[Lemmas 3.4, 3.5]{[LY]} gives the twisted analogue of the Wakimoto convolution relation. Namely, if $\lambda=\lambda_1-\lambda_2$ with $\lambda_1,\lambda_2\in \Lambda^+$, we define the twisted BMW sheaf by
\begin{equation}\label{eq:twisted-bmw}
(J_{\lambda})_{\mu}:=(J_{\lambda_1,!})_{-\lambda_2+{\mu}}\star (J_{-\lambda_2,*})_{\mu}.
\end{equation}
Similarly, we define the dual twisted BMW sheaf by
\begin{equation}\label{eq:twisted-bmw-dk}
(J_{\lambda}^{\BD})_{\mu}:=(J_{\lambda_1,*})_{-\lambda_2+{\mu}}\star (J_{-\lambda_2,!})_{\mu}.
\end{equation}
These objects are well-defined up to tensoring by an even line in cohomological degree $0$, and are independent of the choice of $\lambda_1,\lambda_2$. When ${\mu}=0$, we omit the subscript and simply write $J_{\lambda}, J_{\lambda}^{\BD}$.

Moreover, they satisfy the convolution relation
\begin{equation}\label{eq:bmw-convolution}
(J_\eta)_{\lambda+{\mu}}\star (J_{\lambda})_{\mu} \simeq (J_{\eta+\lambda})_{\mu},
\end{equation}
for all $\eta,\lambda,{\mu}\in \Lambda$, and similarly for the dual version.

\subsection{Construction of standard objects}

For ${\mu}\in \Lambda$, we consider the category
\begin{equation}
\SD_{\cG^G}^{H(\bF)_c^\omega,\psi}(\widetilde{\Fl}_G)_{b_{\mu}} := \SD_{\cG^G}(H(\bF)_c^\omega,\psi\backslash \widetilde{\Fl}_G/T^\omega,b_{\mu}).
\end{equation}
If $J\in\SD_{\cG^G}^{I^\omega,b_\nu}(\widetilde{\Fl}_G)_{b_\mu}$, then the right convolution with $J$ gives a functor from the category with right monodromy $b_\nu$ to the category with right monodromy $b_\mu$.

\subsubsection{Dominant orbits} First, we record the following well-known lemma for completeness.
\begin{lem}\label{lem:trivial-H-B-intersection}
In all the (homogeneous) cases in Section \ref{pol des}, we have $H\cap B=\{1\}$.
\end{lem}
\begin{proof}
For $\GL(M|N)$ with $M<N$, the equality $H\cap B=\{1\}$ follows directly from the choice of the matrix $\sD$ in Section~\ref{trans D}. Indeed, if $g\in B_{G_1}^-$, $u\in H_u$, and $\sD^{-1}gu\sD\in B_{G_2}$, then $u=1$ and $g$ is diagonal; the $(M+1,i)$-entry of $\sD^{-1}g\sD$ is $1-g_{ii}$ for $1\leq i\leq M$, hence $g=1$. In the $\on{F}(4)$ case, we have $\sD=1+E_{4,3}$ (cf.\ \cite[\S 3.4.1]{[FTY]}), and the condition $\sD^{-1}gu\sD\in B_{G_2}$ implies $u=1$ and $g\in T_1$. If we take a lift $\diag(a,a^{-1})$ of $g$ to $\on{SL}_2$, then the $(4,3)$-entry of $\sD^{-1}g\sD$ is $a^{-1}-a$, thus $a^2=1$, and $g=1$ in $\PGL_2$. In the orthosymplectic case, it can also be checked directly from the matrices in Section~\ref{trans D}. Namely, if $g\in B_{G_1}^-$, then its image in $\SO(M+1)$ is lower triangular, and so is $\sD$. Thus, the condition $\sD^{-1}gu\sD\in B_{G_2}$ implies that its component in $H_u$ is trivial and its Levi component is diagonal, i.e., $u=1$ and $g=\sD t\sD^{-1}$ for some $t\in T_{\SO(M+1)}$. Since $g$ fixes the anisotropic vector $v_{\mathrm{an}}$ used to define the embedding, $t$ fixes $\sD^{-1}v_{\mathrm{an}}$, which is
\[
v_{k+1}+\sum_{j=k+2}^{2k+1}v_j
\]
for $M=2k$, and is
\[
v_{l+1}+v_{l+2}+2\sum_{j=l+3}^{2l+2}v_j
\]
for $M=2l+1$ (here, we use the multiple $\sqrt2v_{\mathrm{an}}=v_{l+1}+v_{l+2}$). Note that every coordinate character of $T_{\SO(M+1)}$ appears with a non-zero coefficient in these vectors, thus $t=1$, and hence $g=1$.
\end{proof}
For $\lambda\in \Lambda$, let $\sO_{\widetilde{\Fl}}^{\lambda}$ be the preimage of the $H(\bF)_c^\omega$-orbit of $t^{\lambda}\in {\Fl}_G^\omega$ under the natural projection. We have

\begin{prop}\label{4.2.1}
If $\lambda$ is dominant, then taking the $!$-stalk at $t^{\lambda}$ defines a non-canonical equivalence
\begin{equation}
\SD_{\cG^G}^{H(\bF)_c^\omega,\psi}(\sO_{\widetilde{\Fl}}^{\lambda})_{b_{\mu}}\simeq \SVect.
\end{equation}
\end{prop}

\begin{proof}
It is enough to check the stabilizer of $t^\lambda$. We denote
\[
H_1^\omega= \ker\bigl(H(\bO)_c^\omega\longrightarrow H_c^\omega\bigr).
\]
Since $HB$ is open in $G$ and $H\cap B=\{1\}$, we have $I^\omega=H_1^\omega B^\omega(\bO)$. If $\lambda$ is dominant, then by the Iwahori decomposition, we have
\[
t^\lambda I^\omega t^{-\lambda} \subset I^\omega t^\lambda B^\omega(\bO)t^{-\lambda} =H_1^\omega t^\lambda B^\omega(\bO)t^{-\lambda}.
\]
Thus, since $H(\bF)_c^\omega\cap B(\bF)_c^\omega=\{1\}$, we obtain
\[
K_\lambda:=\Stab_{H(\bF)_c^\omega}(t^\lambda I^\omega/I^\omega) =H(\bF)_c^\omega\cap t^\lambda I^\omega t^{-\lambda} \subset H_1^\omega.
\]
Since $H_1^\omega$ is pro-unipotent, the stabilizer $K_\lambda$ is also pro-unipotent. Also, $\psi|_{K_\lambda}=0$, since $H_1^\omega$ is contained in the kernel of $\psi$. Furthermore, the fixed trivialization of the twisting on $H(\bF)_c^\omega$ restricts to $K_\lambda$, and the right $T^\omega$-monodromy $b_\mu$ is also trivial on it. Hence, we obtain the desired equivalence.
\end{proof}

\subsubsection{Standard objects} Let $(\Delta_0)_{\mu}$ denote the $!$-extension of the irreducible object in $\SD_{\cG^G}^{H(\bF)_c^\omega,\psi} (\sO_{\widetilde{\Fl}}^0)_{b_{\mu}}$ whose $!$-stalk at $t^0$ is $\se$ in even parity and cohomological degree $0$, as in Proposition~\ref{4.2.1} and Remark~\ref{rem:component-normalization}. Equivalently, by the fact that the stabilizer $K_0$ of $t^0$ is pro-unipotent, one may write
\begin{equation}
(\Delta_0)_{\mu}\simeq \Av_!^{H(\bF)_c^\omega,\psi}\bigl((J_0)_{\mu}\bigr).
\end{equation}

\begin{defn}
For $\lambda,{\mu}\in \Lambda$, we denote
\begin{equation}\label{eq:standard-object}
(\Delta_{\lambda})_{\mu}:= (\Delta_0)_{\lambda+{\mu}}\star (J_{\lambda})_{\mu}[-\delta(\lambda)].
\end{equation}
In particular, we define the standard object in the Iwahori Gaiotto category $\SD_{\cG^G}^{H(\bF)_c^\omega,\psi}(\Fl_G^\omega) \simeq \SD_{\cG^G}^{H(\bF)_c^\omega,\psi}(\widetilde{\Fl}_G)_{0}$ by    $\Delta_{\lambda}:=(\Delta_{\lambda})_0$.
\end{defn}
\subsubsection{Convolution and duality} By \eqref{eq:standard-object}, \eqref{eq:bmw-convolution}, and the additivity of $\delta$, we have
\begin{equation}\label{eq:normalized-standard-convolution}
(\Delta_\lambda)_{\eta+\mu}\star(J_\eta)_\mu[-\delta(\eta)] \simeq (\Delta_{\lambda+\eta})_\mu.
\end{equation}
Recall that the Verdier dual $\BD$ on the twisted Gaiotto category is the duality transported from the global model via $E_{\Fl}$. (Similarly in the case with the right monodromy.) With the normalization of $E_{\Fl}$ in Theorem~\ref{prop:global-local-strata}, the global Hecke action of $(J_\eta)_\mu$ corresponds to
\[
\cF\longmapsto \cF\star(J_\eta)_\mu[-\delta(\eta)].
\]
Note that the Verdier duality on the global model is compatible with the global Hecke action. Hence, after the transport by $E_{\Fl}$, we have, on locally compact objects,
\begin{equation}\label{eq:bmw-normalized-duality}
\BD\bigl(\cF\star(J_\eta)_\mu[-\delta(\eta)]\bigr) \simeq \bigl(\BD\cF\star(J_\eta^{\BD})_\mu[-\delta(\eta)]\bigr).
\end{equation}
Here, $(J_\eta^{\BD})_\mu$ is taken with the inverse twisting and the inverse monodromies.

\begin{prop}\label{4.2.3}
If $\lambda\in \Lambda$ is dominant, then $(\Delta_{\lambda})_{\mu}$ is the $!$-extension of the unique irreducible object in $\SD_{\cG^G}^{H(\bF)_c^\omega,\psi}(\sO_{\widetilde{\Fl}}^{\lambda})_{b_{\mu}}$ with the normalization of Remark~\ref{rem:component-normalization}, i.e.,
\begin{equation}
(\Delta_{\lambda})_{\mu}\simeq (\Delta^\stan_{\lambda})_{\mu}.
\end{equation}
Here, $\Delta^{\stan}_{\lambda}$ is the $!$-extension of the irreducible perverse sheaf in $\SD_{\cG^G}^{H^\omega(\bF),\psi}(\sO^{\lambda}_{\Fl})$, as in the normalization of Remark~\ref{rem:component-normalization}.
\end{prop}

\begin{proof}
If $\lambda$ is dominant, then the BMW object is $(J_{\lambda,!})_\mu$. Consider the multiplication map
\begin{equation}\label{4.10}
H(\bF)_c^\omega\times^{H_1^\omega} \bigl(I^\omega t^\lambda I^\omega/I^\omega\bigr) \longrightarrow H(\bF)_c^\omega t^\lambda I^\omega/I^\omega.
\end{equation}
Since $I^\omega=H_1^\omega B^\omega(\bO)$, and $B^\omega(\bO)\subset t^\lambda I^\omega t^{-\lambda}$ by the dominance assumption, the group $H_1^\omega$ acts transitively on $I^\omega t^\lambda I^\omega/I^\omega$ with stabilizer $K_\lambda$. Thus, both sides of \eqref{4.10} are identified with $H(\bF)_c^\omega/K_\lambda$, and the multiplication map is an isomorphism. With the same argument, we also have
\begin{equation}\label{4.11}
H(\bF)_c^\omega\times^{H_1^\omega}\widetilde{\Fl}_G^\lambda \xrightarrow{\sim}\sO_{\widetilde{\Fl}}^\lambda,
\end{equation}
and the monodromies match by construction. Therefore, the convolution of the two $!$-extensions is the $!$-extension from this orbit, with the shift $-l(t^\lambda)=-\deg(\lambda)$ coming from the BMW object. Namely, its $!$-stalk is $\se[-\deg(\lambda)]$ with even coefficients, and after the additional shift $[-\delta(\lambda)]$, it becomes $\se[-d(\lambda)]$. By Remark~\ref{rem:component-normalization}, it is exactly the perverse generator, which proves the assertion.
\end{proof}

The following basic property of the duality functor follows from the definition and Propositions \ref{4.2.1} and \ref{4.2.3}.
\begin{lem}\label{3.5.2}
If $\lambda$ is dominant, then $\BD((\Delta_{\lambda})_\mu)\simeq (\nabla^{\stan}_{\lambda})_\mu$ for every $\mu$. Here, $(\nabla^{\stan}_{\lambda})_\mu$ is the $*$-extension of the dual irreducible perverse sheaf on $\sO_{\widetilde{\Fl}}^\lambda$, with inverse twisting, character $-\psi$, and right monodromy $b_\mu^{-1}$. For $\mu=0$, we denote this object by $\nabla^{\stan}_{\lambda}$.
\end{lem}

\section{Equivalence with factorization modules}\label{foundations-equivalence}

In this section, we prove the equivalence between the full subcategory generated by the standard objects in $\SD_{\cG^G}^{H(\bF)_c^\omega,\psi}(\Fl_G^\omega)$ and the category of factorization modules over $\Omega_q^{\Lus}$. Here, $\Omega_q^{\Lus}$ is the factorization algebra in Definition~\ref{def fact alg}. We will prove its perversity in Proposition~\ref{prop:Omega-geo-perverse}, and compare it with the factorization algebra arising from the quantum supergroup in Section~\ref{sec:factorization-comparison}.

\begin{defn}
We denote by
\[
\cC^\Delta\subset \SD_{\cG^G}^{H(\bF)_c^\omega,\psi}(\Fl_G^\omega)
\]
the full DG subcategory compactly generated by $\Delta_\lambda$ and $\Pi\Delta_\lambda$, $\lambda\in\Lambda$. Equivalently, it is the smallest stable cocomplete $\SVect$-linear full subcategory which is closed under colimits and contains the standard objects $\Delta_{\lambda}$.

Similarly, we denote by $\cC^{\Delta,\glob}\subset \SD_{\cG^G}^{H,\Hecke,\psi}(\cM'_{\infty\cdot c})$ the full subcategory consisting of those objects $\cF$ such that $E_{\Fl}(\cF)\in \cC^\Delta$.
\end{defn}

\begin{rem}
Note that the standard objects $\Delta_\lambda$ are compact in $\SD_{\cG^G}^{H(\bF)_c^\omega,\psi}(\Fl_G^\omega)$, thus they compactly generate $\cC^\Delta$.
\end{rem}

Using the standard objects $\Delta_{\lambda}$ defined in the former section, we equip $\cC^\Delta$ with a $t$-structure.

\begin{defn}\label{5.0.3}
For $\cF\in \cC^\Delta$, we say that $\cF$ is coconnective if and only if
\begin{equation}
H^0\RHom_{\cC^\Delta}(\Pi^\epsilon\Delta_{\lambda}[k],\cF)=0
\end{equation}
for every $k>0$, $\epsilon\in\{0,1\}$, and $\lambda\in \Lambda$.
\end{defn}

\subsubsection{The $!$-stalks} We claim:

\begin{prop}\label{5.0.1}
For any object $\cF\in \SD_{\cG^G}^{H(\bF)_c^\omega,\psi}(\Fl_G^\omega)$, we have
\begin{equation}
i_{\lambda}^!\bigl(F_{\Fl}^{\Lus}(\cF)\bigr)\simeq \Pi^{d(\lambda)}\RHom_{\SD_{\cG^G}^{H(\bF)_c^\omega,\psi}(\Fl_G^\omega)}(\Delta_{\lambda},\cF).
\end{equation}
\end{prop}

\begin{proof}
Note that the BMW objects used for the semi-infinite sheaves have the twisting $(\cG^G)^{-1}$, and the inversion of the loop group gives the twisting $\cG^G$ in the convolution formula. For simplicity, we omit the twisting subscripts.

First, we note that for sufficiently anti-dominant $\nu$, the translated orbit $t^{-\nu}\Fl_G^{\nu+\lambda}$ is contained in $\sS_{\Fl,c}^{\lambda}$, and $\sS_{\Fl,c}^{\lambda}$ is the filtered union of these orbits, cf. \cite[Proposition 6.5.6]{[Yan]}. We take the trivializations induced by the equivariant trivialization on $\sS_{\Fl,c}^{\lambda}$, and the transition maps are induced by the adjunction maps
\[
i_!\omega_Y\simeq i_!i^!\omega_{Y'}\longrightarrow\omega_{Y'},
\]
where $i:Y\hookrightarrow Y'$ runs over the closed embeddings between these orbits. Thus, we have
\begin{equation}\label{5.1-colim-j!}
j^{\sinf}_{\Fl,\lambda,!}(\omega_{\sS_{\Fl,c}^{\lambda}})[d(\lambda)] \simeq \colim_{\nu\ \text{anti-dominant}} t^{-\nu}(J^{\BD}_{\nu+\lambda})_{0}[\deg(-\nu)+\delta(\lambda)].
\end{equation}
Here, $\deg(-\nu)=l(t^\nu)$ equals the dimension of the Iwahori orbit passing through $t^\nu$, and the transition maps are the identity on the normalized $!$-stalk at $t^\lambda$.

Taking the $!$-tensor product with $\cF$, by the continuity of the cohomology $H(\Fl_G^\omega,-)$, we obtain
\begin{equation}\label{7.4 eq}
\begin{split}
i_\lambda^!F_{\Fl}^{\Lus}(\cF) \simeq\Pi^{d(\lambda)}\colim_{\nu\ \text{anti-dominant}} H\bigl(\Fl_G^\omega,\cF\overset!\otimes t^{-\nu}(J^{\BD}_{\nu+\lambda})_0[\deg(-\nu)+\delta(\lambda)]\bigr).
\end{split}
\end{equation}
Let $\iota:G(\bF)_c^\omega\to G(\bF)_c^\omega$ be the inversion map, and $p:G(\bF)_c^\omega\to\Fl_G^\omega$ the projection. Then, $\iota^!$ induces an equivalence between the category of right $I^\omega$-equivariant sheaves on $G(\bF)_c^\omega$ and the category of left $I^\omega$-equivariant sheaves on $G(\bF)_c^\omega$. Applying $\iota^!p^!$ to $(J^{\BD}_{\nu+\lambda})_0$ and then translating by $t^{-\nu}$, we identify the cohomology
\[
H\bigl(\Fl_G^\omega,\cF\overset!\otimes t^{-\nu}(J^{\BD}_{\nu+\lambda})_0\bigr)
\]
with the $!$-stalk of $\cF\star(J_{-\nu-\lambda})_{\nu+\lambda}$ at $t^{-\nu}$. Then, by adjunction, it is identified with 
\[
\RHom_{\SD_{\cG^G}^{H(\bF)_c^\omega,\psi}(\widetilde{\Fl}_G)_{b_{\nu+\lambda}}}\bigl( \Av_!^{H(\bF)_c^\omega,\psi}(\delta_{-\nu})[\deg(\nu)], \cF\star(J_{-\nu-\lambda})_{\nu+\lambda}\bigr),
\]
up to a cohomological shift. Here, we regard $\delta_{-\nu}$ as the delta sheaf at $t^{-\nu}\in G(\bF)_c^\omega/I^\omega$, with even coefficients and right monodromy $b_{\nu+\lambda}$. Since $-\nu$ is dominant, by Proposition~\ref{4.2.1} and the orbit calculation in the proof of Proposition~\ref{4.2.3}, we have
\[
\Av_!^{H(\bF)_c^\omega,\psi}(\delta_{-\nu})[\deg(\nu)] \simeq(\Delta_{-\nu})_{\nu+\lambda}[-\delta(\nu)].
\]
Then, by the convolution with the BMW object $(J_{\nu+\lambda})_0$ and \eqref{eq:bmw-convolution}, the above $\RHom$ is identified with
\[
\RHom\bigl((\Delta_{-\nu})_{\nu+\lambda}[-\delta(\nu)], \cF\star(J_{-\nu-\lambda})_{\nu+\lambda}\bigr) \simeq\RHom(\Delta_{\lambda}[\delta(\lambda)],\cF).
\]
Note that these identifications are compatible with the transition maps, which are the identity on the $!$-stalk at $t^\lambda$. Thus, by \eqref{eq:standard-object}, we obtain
\[
\begin{split}
i_\lambda^!F_{\Fl}^{\Lus}(\cF) &\simeq\Pi^{d(\lambda)}\RHom(\Delta_{\lambda}[\delta(\lambda)],\cF)[\delta(\lambda)]\\
&\simeq\Pi^{d(\lambda)}\RHom(\Delta_\lambda,\cF),
\end{split}
\]
as required.
\end{proof}

\subsubsection{The image of standard objects} We next compare the image of standard objects $\Delta_{\lambda}$ under $F_{\Fl}^{\Lus}$ with the standard factorization modules.

\begin{prop}\label{5.0.5}
For every $\lambda\in \Lambda$, we have
\begin{equation}
F_{\Fl}^{\Lus}(\Delta_{\lambda})\simeq \Pi^{d(\lambda)}\Delta_{\lambda}^{\Lus,\fact}.
\end{equation}
\end{prop}

To check the above proposition, one has to check the $*$-stalk of $F_{\Fl}^{\Lus}(\Delta_{\lambda})$, which is quite hard. To achieve this, we use the global model as in \cite{[G]}.

Let $\Delta_{\lambda}^{\glob}\in \cC^{\Delta,\glob}$ be the object corresponding to $\Delta_{\lambda}$, i.e., the object such that $E_{\Fl}(\Delta_{\lambda}^{\glob})\simeq\Delta_{\lambda}$. According to Proposition~\ref{prop:local-global-functor-comparison}, we should show
\[
F_{\Fl}^{\Lus,\glob}(\Delta_\lambda^{\glob})\simeq\Pi^{d(\lambda)}\Delta_\lambda^{\Lus,\fact}.
\]

\subsubsection{Duality} We need the following compatibility of the global functors with the Verdier duality. Note that whenever $\BD$ is applied, the twisting and the character on the target become $(\cG^G)^{-1}$ and $-\psi$, and the factorization twisting is also inverted.

\begin{lem}\label{global-duality-compatibility}
For any locally compact object $\cF\in \SD_{\cG^G}^{H,\Hecke,\psi}(\cM_{\infty\cdot c})$ (resp. $\SD_{\cG^G}^{H,\Hecke,\psi}(\cM'_{\infty\cdot c})$), we have
\begin{equation}
F_{\Gr}^{\sm,\glob}\bigl(\BD(\cF)\bigr)\simeq \BD\bigl(F_{\Gr}^{\sm,\glob}(\cF)\bigr).
\end{equation}
(resp.
\begin{equation}
F_{\Fl}^{\Lus,\glob}\bigl(\BD(\cF)\bigr)\simeq \BD\bigl(F_{\Fl}^{\DK,\glob}(\cF)\bigr).
\end{equation}
)
The same statements hold if we interchange $\Gr$ and $\Fl$.
\end{lem}

\begin{proof}
The first assertion follows by the same argument as the duality in \cite[Section 7]{[G]}, using the universal local acyclicity (ULA) in \cite[Theorem 4.2.1]{[Cam]}.

For the second assertion, we use the same argument as in \cite[Theorem 8.2.2]{[Yan]}, together with Lemma~\ref{lem:proper-support-global-functors} and the ULA property proved in \cite[Theorem 4.2.1]{[Yan2]}, which is a tamely ramified version of \cite{[Cam]}.
\end{proof}

Combining this lemma with Theorem \ref{prop:global-local-strata}, Proposition \ref{prop:local-global-functor-comparison}, and the definition of the duality on the local model, we obtain:

\begin{cor}\label{5.0.8}
For any locally compact object $\cF\in  \SD_{\cG^G}^{H^\omega(\bF),\psi}(\Gr_G^\omega)$ (resp. $\SD_{\cG^G}^{H^\omega(\bF),\psi}(\Fl_G^\omega)$), we have
\begin{equation}
F_{\Gr}^{\sm}\bigl(\BD(\cF)\bigr)\simeq \BD\bigl(F_{\Gr}^{\sm}(\cF)\bigr).
\end{equation}
(resp. 
\begin{equation}
F_{\Fl}^{\Lus}\bigl(\BD(\cF)\bigr)\simeq \BD\bigl(F_{\Fl}^{\DK}(\cF)\bigr).
\end{equation}
)
The same statements hold if we interchange $\Gr$ and $\Fl$.
\end{cor}

\subsubsection{}
\begin{proof}[Proof of Proposition \ref{5.0.5}]
According to Corollary \ref{5.0.8}, it is enough to show
\begin{equation}
F_{\Fl}^{\DK}\bigl(\BD(\Delta_{\lambda})\bigr)\simeq \Pi^{d(\lambda)}\nabla_{\lambda}^{\DK,\fact}.
\end{equation}
Equivalently, we must prove
\begin{equation}
i_{\lambda}^!\bigl(F_{\Fl}^{\DK}(\BD(\Delta_{\lambda}))\bigr)\simeq \Pi^{d(\lambda)}\se,
\end{equation}
and
\begin{equation}
i_{\mu}^!\bigl(F_{\Fl}^{\DK}(\BD(\Delta_{\lambda}))\bigr)\simeq 0 \qquad\text{for }{\mu}\neq \lambda.
\end{equation}

As in \cite[Section 8.3]{[Yan]}, we use the calculation on the translates of Iwahori orbits. Let $\cK=\BD(\Delta_\lambda)$, and note that the standard objects in the first argument of each $\RHom$ below are taken with the inverse twisting and the character $-\psi$. We let $\cK^{\glob}=E_{\Fl}^{-1}(\cK)$, whose support is contained in finitely many global strata. By Lemma~\ref{lem:proper-support-global-functors}, the support of
\[
p_{\Fl}^!\cK^{\glob}\overset!\otimes q_{\Fl}^!j_{*,\Fl}^e(\IC_{\Bun_B^e})
\]
over $\mu\cdot c$ is contained in a scheme of finite type. Under the identification of the local and global Zastava spaces in Section~\ref{sec:global-zastava-functors}, its restriction to the open semi-infinite stratum has support $\operatorname{supp}(\cK)\cap\sS_{\Fl,c}^{\mu}$. Thus, this intersection is of finite type, and it is contained in $t^{-\nu}\Fl_G^{\nu+\mu}$ for sufficiently anti-dominant $\nu$. For such $\nu$, we have
\begin{equation}\label{5.0.5-j*}
\begin{split}
i_\mu^!F_{\Fl}^{\DK}(\cK) &\simeq\Pi^{d(\mu)} H\bigl(\Fl_G^\omega,\cK\overset!\otimes t^{-\nu}(J_{\nu+\mu})_0[\deg(-\nu)+\delta(\mu)]\bigr)\\
&\simeq\Pi^{d(\mu)} \RHom_{\SD_{(\cG^G)^{-1}}^{H(\bF)_c^\omega,-\psi} (\widetilde{\Fl}_G)_{b_{\nu+\mu}^{-1}}} \bigl((\Delta_{-\nu})_{\nu+\mu}[-\delta(\nu)], \cK\star(J^{\BD}_{-\nu-\mu})_{\nu+\mu}\bigr)[\delta(\mu)].
\end{split}
\end{equation}

Let $\eta=-\nu-\mu$. By \eqref{eq:normalized-standard-convolution} and \eqref{eq:bmw-normalized-duality}, we have
\begin{equation}\label{7.13 eq}
\begin{split}
\BD(\Delta_\lambda)\star(J_\eta^{\BD})_{\nu+\mu}[-\delta(\eta)] &\simeq\BD\bigl(\Delta_\lambda\star(J_\eta)_{\nu+\mu}[-\delta(\eta)]\bigr)\\
&\simeq \BD\bigl((\Delta_{\lambda+\eta})_{\nu+\mu}\bigr).
\end{split}
\end{equation}
Note that the cohomological shifts in the resulting $\RHom$ cancel each other, since
\[
\delta(\mu)+\delta(\nu)+\delta(\eta)=0.
\]
Hence,
\begin{equation}\label{5.0.5-rhom}
i_\mu^!F_{\Fl}^{\DK}(\BD\Delta_\lambda) \simeq\Pi^{d(\mu)}\RHom\bigl((\Delta_{-\nu})_{\nu+\mu}, \BD((\Delta_{\lambda-\mu-\nu})_{\nu+\mu})\bigr).
\end{equation}
For sufficiently anti-dominant $\nu$, both $-\nu$ and $\lambda-\mu-\nu$ are dominant, thus by Proposition~\ref{4.2.3} and Lemma~\ref{3.5.2}, the first object is the $!$-extension of the normalized generator on its orbit, and the second one is the $*$-extension of the dual normalized generator. Therefore, if $\lambda=\mu$, the corresponding $\RHom$ is $\se$ by adjunction; if $\lambda\neq\mu$, the generators lie on different orbits, thus the corresponding $\RHom$ is zero.

Thus, we obtain
\begin{equation}
i_{\mu}^!\bigl(F_{\Fl}^{\DK}(\BD(\Delta_{\lambda}))\bigr)\simeq
\begin{cases}
\Pi^{d(\lambda)}\se,& {\mu}=\lambda,\\
0,& {\mu}\neq \lambda.
\end{cases}
\end{equation}
Therefore, by the convention in Remark~\ref{rem:target-parity}, we have
\begin{equation}
F_{\Fl}^{\DK}(\BD(\Delta_{\lambda}))\simeq \Pi^{d(\lambda)}\nabla_{\lambda}^{\DK,\fact}.
\end{equation}
Applying Verdier duality and Corollary \ref{5.0.8}, we conclude that
\begin{equation}
F_{\Fl}^{\Lus}(\Delta_{\lambda})\simeq \Pi^{d(\lambda)}\Delta_{\lambda}^{\Lus,\fact}.
\end{equation}
\end{proof}

\subsubsection{The equivalence} Finally, we prove:

\begin{thm}\label{main}
The functor
\begin{equation}
F_{\Fl}^{\Lus}:\cC^\Delta\longrightarrow \Omega_q^{\Lus}\FS
\end{equation}
induces a $t$-exact equivalence of highest weight categories which sends the standard objects $\Delta_{\lambda}$ to $\Pi^{d(\lambda)}\Delta_{\lambda}^{\Lus,\fact}$ for all $\lambda\in \Lambda$.
\end{thm}

\begin{proof}
To prove full faithfulness, we first note that by Proposition \ref{5.0.1}, for any $\cF\in \cC^\Delta$, we have
\begin{equation}\label{5.20}
\begin{split}
\RHom_{\SD_{\cG^G}^{H(\bF)_c^\omega,\psi}(\Fl_G^\omega)}(\Delta_{\lambda},\cF) &\simeq \Pi^{d(\lambda)}i_{\lambda}^!\bigl(F_{\Fl}^{\Lus}(\cF)\bigr)\\
&\simeq \RHom_{\Omega_q^{\Lus}\FS}\bigl(\Pi^{d(\lambda)}\Delta_{\lambda}^{\Lus,\fact},F_{\Fl}^{\Lus}(\cF)\bigr)\\
&\simeq \RHom_{\Omega_q^{\Lus}\FS}\bigl(F_{\Fl}^{\Lus}(\Delta_{\lambda}),F_{\Fl}^{\Lus}(\cF)\bigr),
\end{split}
\end{equation}
where the last isomorphism follows from Proposition \ref{5.0.5}.

We should explain why the map induced by $F_{\Fl}^{\Lus}$ is an isomorphism for any $\cF\in \cC^\Delta$. Actually, it is true for any $\cF\in \SD_{\cG^G}^{H^\omega(\bF),\psi}(\Fl_G^\omega)$. Namely, since $\RHom_{\SD_{\cG^G}^{H^\omega(\bF),\psi}(\Fl_G^\omega)}(\Delta_\lambda,\Delta_\lambda)\simeq\se$, we can choose a set of generators $\{\Delta_\lambda,c_1,c_2,\ldots\}$ of $\SD_{\cG^G}^{H^\omega(\bF),\psi}(\Fl_G^\omega)$ such that
\[
\RHom_{\SD_{\cG^G}^{H^\omega(\bF),\psi}(\Fl_G^\omega)}(\Delta_\lambda,c_i)=0.
\]
By the compactness of $\Delta_\lambda$ and \eqref{5.20}, both sides commute with colimits, thus it is enough to check that the map of $\RHom$ defined by $F_{\Fl}^{\Lus}$ is an isomorphism for $\cF\in\{\Delta_\lambda,c_1,c_2,\ldots\}$. If $\cF=\Delta_\lambda$, then both sides are $\se$, and since $F_{\Fl}^{\Lus}$ sends $\id$ to $\id$, the above map is an isomorphism. If $\cF=c_i$, then both sides are zero by \eqref{5.20}.

Since the standard objects $\{\Delta_\lambda\}_{\lambda\in\Lambda}$ compactly generate $\cC^\Delta$, we have $F_{\Fl}^{\Lus}$ is fully faithful on $\cC^\Delta$.

For essential surjectivity, we first note that $F_{\Fl}^{\Lus}$ commutes with filtered colimits. Furthermore, $\Omega_q^{\Lus}\FS$ is generated by the standard factorization modules $\Delta_\lambda^{\Lus,\fact}$ under parity shifts, cohomological shifts, extensions, and filtered colimits. Therefore, essential surjectivity follows from Proposition~\ref{5.0.5}.

Finally, $t$-exactness follows from the description of the $t$-structures on $\cC^\Delta$ and $\Omega_q^{\Lus}\FS$ in Definition \ref{5.0.3} and the corresponding definition on the factorization side: both are characterized by vanishing of $\Hom$ from standard objects, and $F_{\Fl}^{\Lus}$ preserves standard objects up to parity shifts by Proposition \ref{5.0.5}.
\end{proof}

\section{Mixed quantum supergroups}\label{sec:mixed-quantum-supergroups} In this section, we recall the quantum supergroups which will be compared with the factorization algebras constructed geometrically. We essentially just copy the relevant part in \cite{[G7]} and formulate it for the supergroups.

In principle, the positive half of the De Concini--Kac form is defined by deforming the braiding parameter over a formal disc and then specializing, while the positive half of the Lusztig form is obtained by duality. In the super case, the role of the Weyl group of the ordinary root system is played by the Weyl groupoid of the chosen super root datum. In this section, we will work over $\BC$, and the algebraic constructions also make sense over any algebraically closed (coefficient) field of characteristic $0$.

The section is organized as follows. First, we recall the quantum torus and the five versions of the positive half of the quantum supergroup, namely, the free, cofree, small, De Concini--Kac and Lusztig positive halves. Secondly, we recall the Weyl groupoid, the Lusztig isomorphisms and the PBW bases, which allow us to describe the De Concini--Kac positive half at roots of unity under two conditions on $q$, and we find a further condition under which its root powers generate a central Hopf subalgebra. Thirdly, we recall the factorization algebras attached to these Hopf superalgebras and the quantum Frobenius map, following \cite{[G7]}. Finally, we define the category $\cO$ of the mixed quantum supergroup, describe it in terms of factorization modules, and fix the super Borels that will be used in the rest of the paper.

\subsection{The quantum torus}

In this subsection, we fix the notations for super roots and recall the braided monoidal category $\Rep_q(T_s)$, in which the positive halves below live. Let $\fg_s$ be one of the Lie superalgebras $\mathfrak{gl}(M|N)$, $\mathfrak{osp}(2k|2l)$, and $\mathfrak f(4)$ attached to the cotangent data fixed in Section~\ref{pol des}, and let $T_s$ be its even maximal torus. We denote its weight lattice $X^*(T_s)=X_*(T)$ by $\Lambda$, and its super root system by $\Phi_s$. For a chosen super Borel $B_s$ indexed by $a$ (an object of the Weyl groupoid recalled below),  we denote its ordered set of simple super roots by $\Pi_a=\{\alpha_i^a\}_{i\in I}$. We also denote by $p_i^a:=p(\alpha_i^a)\in\{0,1\}$ the parity of $\alpha_i^a$ and by $Q_s:=\sum_{i\in I}\BZ\alpha_i^a$ the super root lattice, and we set $\BN=\BZ_{\geq0}$. Then, the positive and negative super root monoids are
\begin{equation*}
\Lambda^{\mathrm{pos}}:=\sum_{i\in I}\BN\alpha_i^a\subset \Lambda, \qquad \Lambda^{\mathrm{neg}}:=-\Lambda^{\mathrm{pos}}.
\end{equation*}

Recall from Section~\ref{3.2.5 sect} that the factorization gerbe $\cG^\Lambda$ determines a multiplicative quadratic form
\begin{equation*}
q:\Lambda\longrightarrow \mathsf e^\times(-1)
\end{equation*}
in the Betti or $\ell$-adic setting, and an additive quadratic form in the $D$-module setting. In the $D$-module setting, we will use the multiplicative convention $q(\lambda)=\exp(2\pi i\,q_{\mathrm{add}}(\lambda))$ for the additive quadratic form $q_{\mathrm{add}}$. Then, we choose a bilinear form
\begin{equation*}
b':\Lambda\times\Lambda\longrightarrow \mathsf e^\times
\end{equation*}
such that, for all $\lambda,\mu\in\Lambda$,
\begin{equation*}
b'(\lambda,\lambda)=q(\lambda), \qquad b'(\lambda,\mu)b'(\mu,\lambda) =q(\lambda+\mu)q(\lambda)^{-1}q(\mu)^{-1}.
\end{equation*}
In other words, $b'$ is a bilinear form refining $q$, i.e., $b'(\lambda,\lambda)=q(\lambda)$, and hence its symmetrization $b'(\lambda,\mu)b'(\mu,\lambda)$ is the symmetric bilinear form $b$ of Section~\ref{3.2.5 sect}, and the corresponding additive formulas are obtained by replacing products by sums.

We denote by $\Rep_q(T_s)$ the braided monoidal category of representations of the quantum torus attached to $q$. Namely, its underlying monoidal category is the category of $\Lambda$-graded super vector spaces, and after choosing $b'$, the braiding is given by
\begin{equation*}
M_{\lambda}\otimes N_{\mu}\longrightarrow N_{\mu}\otimes M_{\lambda}, \qquad m\otimes n\longmapsto (-1)^{p(m)p(n)}b'(\lambda,\mu)\,n\otimes m,
\end{equation*}
for homogeneous elements $m\in M_{\lambda}$ and $n\in N_{\mu}$. Note that, as in \cite[Remark 3.1.3]{[G7]}, the braided category $\Rep_q(T_s)$ depends, up to a braided equivalence, only on $q$ and not on the chosen $b'$.

For super root-lattice degrees $\beta,\beta'$, it is convenient to denote
\begin{equation}\label{eq:full-bicharacter-super}
\chi(\beta,\beta') :=(-1)^{p(\beta)p(\beta')}b'(\beta,\beta').
\end{equation}
Namely, $\chi(\beta,\beta')$ is the scalar occurring in the braiding of homogeneous super root vectors of degrees $\beta$ and $\beta'$, where the sign comes from the super structure.

\subsection{Free, cofree, small, De Concini--Kac, and Lusztig positive halves}

Since $\Rep_q(T_s)$ is braided monoidal, it makes sense to consider braided Hopf superalgebras in $\Rep_q(T_s)$. In this subsection, we recall five such Hopf superalgebras, namely, the five versions of the positive half of the quantum supergroup.

For the super Borel $a$, we set $B_s^+=B_s$, and we denote by $B_s^-$ its opposite super Borel and by $N_s^+$ and $N_s^-$ their unipotent radicals, respectively. Let
\begin{equation*}
V_a^+:=\bigoplus_{i\in I}\Pi^{p_i^a}\BC_{\alpha_i^a}\in\Rep_q(T_s),
\end{equation*}
where $\BC_{\alpha_i^a}$ is one-dimensional, even and of $T_s$-weight $\alpha_i^a$, and we denote by $E_i^a$ the generator of $\Pi^{p_i^a}\BC_{\alpha_i^a}$. The positive half of the free quantum supergroup $U_q^{\fr}(N_s^+)_a$ is the free algebra generated by $V_a^+$ in $\Rep_q(T_s)$, equipped with the structure of a Hopf superalgebra given by
\begin{equation*}
\Delta(E_i^a)=1\otimes E_i^a+E_i^a\otimes 1,\qquad \epsilon(E_i^a)=0.
\end{equation*}

Similarly, we have the cofree Hopf superalgebra $U_q^{\cofr}(N_s^+)_a$, which is cogenerated by $E_i^a$, and one has
\begin{equation*}
U_q^{\cofr}(N_s^+)_a:=\bigl(U_{q^{-1}}^{\fr}(N_s^-)_a\bigr)^\vee,
\end{equation*}
where $(H^\vee)_\lambda=(H_{-\lambda})^*$ with the same parity, and the negative halves are generated in degrees $-\alpha_i^a$. Note that here we use the Hopf duality with the inverse braiding in \cite[Section 3.1.5]{[G7]}, followed by the reversal of the lattice grading, so that the dual above belongs to $\Rep_q(T_s)$. Then, we have the canonical Hopf morphism
\begin{equation*}
U_q^{\fr}(N_s^+)_a\longrightarrow U_q^{\cofr}(N_s^+)_a,
\end{equation*}
which sends the generators $E_i^a$ to the cogenerators $E_i^a$.

The small quantum supergroup is defined by
\begin{equation*}
u_q^{\sm}(N_s^+)_a:=\operatorname{Im}\bigl(U_q^{\fr}(N_s^+)_a\to U_q^{\cofr}(N_s^+)_a\bigr).
\end{equation*}
In other words, $u_q^{\sm}(N_s^+)_a$ is the Nichols algebra of $V_a^+$, and we regard it as the small positive half.

We define $\Rep_q^\sm(\fg_s)$ similarly as in \cite[Sections 27.3 and 28.1]{[GL]}. Namely, we take the relative Drinfeld center of $u_q^{\sm}(N_s^+)_a\modu(\Rep_q(T_s))$ over $\Rep_q(T_s)$, and then the ind-completion of its full subcategory of objects with finite-dimensional total cohomology.

\subsubsection{The deformation}\label{deform algebra} In order to define the De Concini--Kac form, we use the formal deformation, like \cite[Sections 1.7 and 3.3]{[G7]} and \cite{[Y1]}. Roughly speaking, we deform the braiding parameter $q$ over the formal disc $\Spec\BC[[\hbar]]$, so that it becomes generic after inverting $\hbar$. Here, on the super root lattice, the parameter has the form $q(\beta)=\xi^{q_{\BZ}(\beta)}$, where $\xi\in\BC^\times$ is the twisting parameter and $q_{\BZ}$ is an integral normalization of the invariant quadratic form defining the twisting, and since $q_{\BZ}$ comes from the invariant quadratic form, we have $q_{\BZ}(\beta)=0$ for isotropic super roots $\beta$ and $q_{\BZ}(\beta)\neq0$ for non-isotropic super roots. Then, we choose an integral bilinear form $b'_{\BZ}:\Lambda\times\Lambda\to\BZ$ such that
\[
b'_{\BZ}(\lambda,\lambda)=q_{\BZ}(\lambda), \qquad \lambda\in\Lambda.
\]

We set $\zeta_\hbar=\exp(\hbar)$ in the classical topology and $\zeta_\hbar=1+\hbar$ in the $\ell$-adic setting, and define
\begin{equation}\label{eq:formal-quadratic-bilinear-super}
q_\hbar(\lambda) :=q(\lambda)\zeta_\hbar^{q_{\BZ}(\lambda)}, \qquad b'_\hbar(\lambda,\mu) :=b'(\lambda,\mu)\zeta_\hbar^{b'_{\BZ}(\lambda,\mu)}.
\end{equation}

We denote by $\Rep_{q[[\hbar]]}(T_s)$ and $\Rep_{q((\hbar))}(T_s)$ the resulting braided categories over $\BC[[\hbar]]$ and $\BC((\hbar))$, respectively.

\subsubsection{The De Concini--Kac and Lusztig forms} Roughly speaking, the De Concini--Kac positive half is the specialization at $\hbar=0$ of the $\BC[[\hbar]]$-subalgebra generated by the $E_i^a$ in the small positive half over $\BC((\hbar))$, and the Lusztig positive half is obtained by duality. More precisely, let $I_{q((\hbar))}$ be the kernel of the map
\begin{equation*}
U_{q((\hbar))}^{\fr}(N_s^+)_a\longrightarrow U_{q((\hbar))}^{\cofr}(N_s^+)_a,
\end{equation*}
and we define the Hopf ideal and the Hopf superalgebra over $\BC[[\hbar]]$ by
\begin{equation*}
I_{q[[\hbar]]}:=U_{q[[\hbar]]}^{\fr}(N_s^+)_a\cap I_{q((\hbar))}, \qquad U_{q[[\hbar]]}^{\DK}(N_s^+)_a:=U_{q[[\hbar]]}^{\fr}(N_s^+)_a/I_{q[[\hbar]]},
\end{equation*}
where the intersection is taken inside $U_{q((\hbar))}^{\fr}(N_s^+)_a$. Note that this Hopf superalgebra is $\BC[[\hbar]]$-flat, and we define the De Concini--Kac positive half by
\begin{equation*}
U_q^{\DK}(N_s^+)_a:=U_{q[[\hbar]]}^{\DK}(N_s^+)_a\otimes_{\BC[[\hbar]]}\BC.
\end{equation*}
Then, the Lusztig positive half is defined by the componentwise Hopf duality,
\begin{equation*}
U_q^{\Lus}(N_s^+)_a:=\bigl(U_{q^{-1}}^{\DK}(N_s^-)_a\bigr)^\vee.
\end{equation*}
The canonical maps factor as
\begin{equation}\label{eq:five-positive-halves-super}
U_q^{\fr}(N_s^+)_a\twoheadrightarrow U_q^{\DK}(N_s^+)_a\twoheadrightarrow u_q^{\sm}(N_s^+)_a\hookrightarrow U_q^{\Lus}(N_s^+)_a\hookrightarrow U_q^{\cofr}(N_s^+)_a.
\end{equation}
Namely, the De Concini--Kac positive half is a quotient of the free one which surjects onto the small one, and the Lusztig positive half is a subobject of the cofree one which contains the small one.

\subsection{Weyl groupoids, Lusztig isomorphisms, and PBW super root vectors}

In order to describe the De Concini--Kac positive half at roots of unity, we need its PBW basis. In the super case, the simple systems related by odd reflections are not conjugate under the Weyl group, and the role of the Weyl group is played by the Weyl groupoid. In this subsection, we recall the Weyl groupoid and the Lusztig isomorphisms, and use them to construct the super root vectors and the PBW bases. We also introduce two conditions on $q$, i.e., Conditions~\ref{av s o} and~\ref{cond:isotropic-non-Cartan}, under which the De Concini--Kac positive half is described in terms of the distinguished pre-Nichols algebras in the sense of \cite{[A]}, see Lemma~\ref{lem:saturated-DK-pre-Nichols}.

\subsubsection{Weyl groupoids} We recall the Weyl groupoid, following (essentially just copying from) \cite{[CH], [HY], [H], [HY2], [A]} and Yamane's odd-reflection formalism. Let $I$ be a finite index set and let $\{\alpha_i\}_{i\in I}$ be the standard basis of the abstract lattice $\BZ^I$. At an object $a$ of the groupoid, the abstract vector $\alpha_i$ will be identified with the actual simple super root $\alpha_i^a$.

\begin{defn}[\cite{[CH]}, Definition 2.1]
A Cartan scheme is a tuple
\begin{equation*}
\cC=\cC\bigl(I,A,(\rho_i)_{i\in I},(C^a)_{a\in A}\bigr), \qquad C^a=(c_{ij}^a)_{i,j\in I},
\end{equation*}
where $A$ is non-empty, $\rho_i:A\to A$ are involutions, and the generalized Cartan matrices satisfy $c_{ij}^a=c_{ij}^{\rho_i(a)}$. For $a\in A$ and $i\in I$, we define $\sigma_i^a\in\Aut(\BZ^I)$ by $\sigma_i^a(\alpha_j)=\alpha_j-c_{ij}^a\alpha_i$.
\end{defn}

\begin{defn}[\cite{[CH]}]
The Weyl groupoid $\CW(\cC)$ has object set $A$, and its morphisms are generated by $\sigma_i^a:a\to\rho_i(a)$.
\end{defn}

In our situation, $A$ is the orbit of an ordered fundamental system of $\fg_s$ under even and odd reflections. At $a\in A$, we denote by
\begin{equation*}
\Pi_a=\{\alpha_i^a\}_{i\in I}, \qquad \Delta_+^a\subset\Phi_s
\end{equation*}
the corresponding simple and positive super roots, and set $\Delta^a:=\Delta_+^a\sqcup(-\Delta_+^a)$. Then, the Cartan matrices $C^a$ are determined by the super root strings. Namely, for $i\neq j$, we set
\begin{equation}\label{eq:root-string-cartan-integers}
-c_{ij}^a:=\max\left\{n\in\BN\ \middle|\ \alpha_j^a+n\alpha_i^a\in\Delta_+^a\right\}.
\end{equation}
For example, if $\alpha_i^a$ is an isotropic odd simple super root, then, for $j\neq i$, the $\alpha_i^a$-string through $\alpha_j^a$ has length at most one, and hence
\begin{equation*}
c_{ij}^a=
\begin{cases}
-1,& \alpha_j^a+\alpha_i^a\in\Delta^a,\\
0,& \alpha_j^a+\alpha_i^a\notin\Delta^a.
\end{cases}
\end{equation*}

Note that, in the cotangent cases under consideration (not the cases of type $B$ and $\on{G}(3)$), every odd super root is isotropic.

\subsubsection{The braiding and the Cartan integers} On the other hand, the braiding of $\Rep_q(T_s)$ also determines Cartan integers, and we need them to agree, up to sign, with those of the Cartan scheme. Namely, the diagonal braiding at $a$ is encoded by
\begin{equation*}
\widetilde q_{ij}^a :=\chi(\alpha_i^a,\alpha_j^a) =(-1)^{p_i^ap_j^a}b'(\alpha_i^a,\alpha_j^a),
\end{equation*}
and the Cartan integer defined by the braiding is
\begin{equation}\label{eq:Nichols-Cartan-integer}
m_{ij}^a :=\min\left\{n\in\BN\ \middle|\ (n+1)_{\widetilde q_{ii}^a} \left(1-(\widetilde q_{ii}^a)^n\widetilde q_{ij}^a\widetilde q_{ji}^a\right)=0\right\},
\end{equation}
where $(r)_z=1+z+\cdots+z^{r-1}$.

\begin{condition}[avoids small torsion]\label{av s o}
We say that $q$ avoids small torsion if $m_{ij}^a=-c_{ij}^a$ for every $a\in A$ and $i\neq j$. For a purely even isolated $A_1$ component, we also require $q(\alpha_i^a)\neq1$.
\end{condition}

\begin{rem}
The above condition can be understood as follows.
\begin{enumerate}
\item If $\alpha_i^a$ is isotropic odd (so, $q(\alpha_i^a)=1$), then $m_{ij}^a=-c_{ij}^a$ is equivalent to
\begin{equation*}
\alpha_i^a+\alpha_j^a\in\Delta_+^a \quad\Longleftrightarrow\quad \widetilde q_{ij}^a\widetilde q_{ji}^a\neq1.
\end{equation*}
Intuitively, it says that the braiding of $E_i^a$ and $E_j^a$ is symmetric exactly when $\alpha_i^a$ and $\alpha_j^a$ are not adjacent.
\item For a even super root $\alpha_i^a$, since $2(\alpha_i^a,\alpha_j^a)=c_{ij}^a (\alpha_i^a,\alpha_i^a)$, we have
\begin{equation}\label{eq:Cartan-vertex-super}
\widetilde q_{ij}^a\widetilde q_{ji}^a = (\widetilde q_{ii}^a)^{c_{ij}^a}.
\end{equation}
Then,
\begin{equation*}
m_{ij}^a=-c_{ij}^a \quad\Longleftrightarrow\quad \operatorname{ord}(\widetilde q_{ii}^a)>-c_{ij}^a,
\end{equation*}
which is the condition of avoiding small torsion used in \cite{[Lus]} and in \cite[Section 1.1.5]{[G7]}. For an isolated even vertex, the equalities $m_{ij}^a=-c_{ij}^a=0$ impose no condition on $q(\alpha_i^a)$. Thus, we separately require $q(\alpha_i^a)\neq1$, as in the $A_1$ case of \cite[Section 1.1.5]{[G7]}.
\end{enumerate}
\end{rem}

\subsubsection{Cartan roots} Following the terminology for Nichols algebras, a vertex $i$ at $a$ is called a Cartan vertex if \eqref{eq:Cartan-vertex-super} holds for every $j\neq i$, and the Cartan roots are the super roots obtained from Cartan vertices by the Weyl groupoid. Namely, we denote by $\mathcal O^a$ the set of super roots $w(\alpha_i^b)$, where $w:b\to a$ is a morphism of the Weyl groupoid and $i$ is a Cartan vertex at $b$.

Also, we denote
\begin{equation*}
\mathcal O_+^a:=\mathcal O^a\cap\Delta_+^a.
\end{equation*}
These are exactly the positive super roots whose power relations are omitted in the distinguished pre-Nichols algebra at roots of unity, cf., \cite{[A]}. We will see below that, under Conditions~\ref{av s o} and~\ref{cond:isotropic-non-Cartan}, on each connected component other than a single isotropic odd vertex, the positive Cartan roots are exactly the positive non-isotropic super roots, and their root powers survive in the De Concini--Kac positive half, see Theorem~\ref{thm:DK-PBW-super}.

\subsubsection{The isotropic non-Cartan condition} Note that, in the cotangent cases considered here, every non-isotropic super root is even, and for the braiding $\chi$ in \eqref{eq:full-bicharacter-super}, these roots are Cartan by \eqref{eq:Cartan-vertex-super} and the definition of $\mathcal O^a$. However, the converse may not be true, i.e., an isotropic odd super root may also be Cartan, and we impose the following condition to avoid this situation.
\begin{condition}[Isotropic non-Cartan condition]
\label{cond:isotropic-non-Cartan} For every object $a$ of the Weyl groupoid and every non-isolated isotropic odd simple super root $\alpha_i^a$, there exists $j\neq i$ such that
\begin{equation}\label{eq:isotropic-non-Cartan-condition}
c_{ij}^a=-1 \qquad\text{and}\qquad b(\alpha_i^a,\alpha_j^a)\neq -1.
\end{equation}
\end{condition}
Recall that $b(\lambda,\mu):=b'(\lambda,\mu)b'(\mu,\lambda)$ denotes the symmetric bilinear form. The following lemma is easy to check from the definition, and it explains the name of the condition.
\begin{lem}
Assume $q(\alpha_i^a)=1$ for every isotropic simple super root and Condition~\ref{av s o}. Then, Condition~\ref{cond:isotropic-non-Cartan} is equivalent to the assertion that every non-isolated isotropic odd super root is non-Cartan.
\end{lem}
\begin{proof}
Note that an isotropic simple super root has $\widetilde q_{ii}^a=-1$. Then, at a non-adjacent vertex $j$, the Cartan relation is $b(\alpha_i^a,\alpha_j^a)=1$, which follows from Condition~\ref{av s o}, while at an adjacent vertex $j$, it is $b(\alpha_i^a,\alpha_j^a)=-1$. Thus, the vertex fails to be Cartan if and only if the condition in \eqref{eq:isotropic-non-Cartan-condition} holds for some adjacent $j$, and by applying this argument at every object of the Weyl groupoid, we obtain the lemma.
\end{proof}

\begin{rem}
Condition~\ref{cond:isotropic-non-Cartan} is not a weird condition. Just as an example, suppose that $q(\lambda)=\zeta^{q_{\BZ}(\lambda)}$, where $q_{\BZ}$ comes from the invariant form. Then, if $\operatorname{ord}(\zeta)$ is odd, Condition~\ref{cond:isotropic-non-Cartan} is satisfied. Indeed, at each non-isolated isotropic vertex, the value of $b$ on an adjacent vertex has odd order and cannot equal $-1$.
\end{rem}

\subsubsection{Distinguished pre-Nichols algebras}
In this subsubsection, we compare the De Concini--Kac positive half with the distinguished pre-Nichols algebras of \cite{[A]}. Since the isolated odd vertices behave differently, we first treat them separately.
\begin{rem}\label{rem:isolated-odd-factor}
Suppose that an isotropic odd simple super root corresponds to a vertex which is isolated, i.e., not connected to any other vertex of the Dynkin diagram. Then, its generator $E$ has braiding $-1$ throughout the deformation, and the map from the free algebra to the cofree algebra kills $E^2$ and is injective in degrees $0$ and $1$. Thus, the small, De Concini--Kac, and Lusztig positive halves for this component are all $\BC[E]/(E^2)$, i.e., an exterior algebra in one odd variable.
\end{rem}

The following fact is well-expected and known to experts. We thank H.~Yamane for introducing us to the relation between pre-Nichols algebras and quantum supergroups.
\begin{lem}\label{lem:saturated-DK-pre-Nichols}
We consider the deformation in Section~\ref{deform algebra} and assume Conditions~\ref{av s o} and~\ref{cond:isotropic-non-Cartan}.

If $\xi$ has infinite order, then $U_q^{\DK}(N_s^+)_a\simeq u_q^{\sm}(N_s^+)_a$, which is the Nichols algebra of $V_a^+$ in $\Rep_q(T_s)$.

Suppose that $\xi$ is a root of unity. For a connected component $J\subset I$ of the Dynkin diagram, we denote by $V_{a,J}^+\subset V_a^+$ the span of $E_i^a$, $i\in J$, with the braiding inherited from $\Rep_q(T_s)$. If $J$ is not a single isotropic odd vertex, then the subalgebra of $U_q^{\DK}(N_s^+)_a$ generated by $V_{a,J}^+$ is the distinguished pre-Nichols algebra of $V_{a,J}^+$ defined in \cite[Definition 3.1]{[A]}.

If $J=\{i\}$ consists of a single isotropic odd vertex, then this subalgebra is $\BC[E_i^a]/((E_i^a)^2)$, while the distinguished pre-Nichols algebra of $V_{a,J}^+$ in \cite[Definition 3.1]{[A]} is $\BC[E_i^a]$.
\end{lem}
\begin{proof}
The case of a single isotropic odd vertex follows from Remark~\ref{rem:isolated-odd-factor}. For the other connected components, we use the integral form in \cite[Proposition 5.4]{[AAY]}, which, by its PBW basis, embeds into $u_{q((\hbar))}^{\sm}(N_s^+)_a$ after the base change along our deformation, and its image is the $\BC[[\hbar]]$-subalgebra generated by the simple root vectors. Thus, by the definition of $U_{q[[\hbar]]}^{\DK}(N_s^+)_a$, it agrees with the subalgebra of $U_{q[[\hbar]]}^{\DK}(N_s^+)_a$ generated by $E_i^a$, $i\in J$, and the root of unity case follows from \cite[Theorem 5.7]{[AAY]}.

Finally, if $\xi$ has infinite order, then the quantum factorials for even roots are non-zero, and the PBW exponents of odd root vectors are $0$ or $1$. Thus, the pairing in \cite[Proposition 5.5]{[AAY]} remains non-degenerate after replacing the divided powers by the ordinary powers, which identifies the specialization with the Nichols algebra of $V_a^+$ in $\Rep_q(T_s)$.
\end{proof}

Note that, in the generic case, the middle three terms in \eqref{eq:five-positive-halves-super} are isomorphic, while at roots of unity, they are different in general.

From now on, we impose Conditions~\ref{av s o} and~\ref{cond:isotropic-non-Cartan}, and we use the identifications of Lemma~\ref{lem:saturated-DK-pre-Nichols} on each connected component. We denote by $\Delta^a_{+,\niso}$ the set of positive non-isotropic super roots (which is also the set of positive even super roots in the cotangent case), and by $\Delta^a_{+,\iiso}$ the set of positive isotropic super roots (which is also the set of positive odd super roots in the cotangent case). Note that, on each connected component other than a single isotropic odd vertex, the positive Cartan roots are exactly the positive non-isotropic super roots.

\subsubsection{Lusztig isomorphisms and their local form}

In order to define the super root vectors, we need first work with the full De Concini--Kac quantum supergroup, since the reflection of a simple positive generator may produce a negative generator. Let us denote by $U_q^{\DK}(\fg_s)_a$ the full De Concini--Kac quantum supergroup attached to the simple system $a$, defined by the quantum double construction in \cite[Section 3.1]{[A]}. Note that, for an isolated isotropic odd component, both halves are exterior algebras, as in Remark~\ref{rem:isolated-odd-factor}. The following theorem is the Lusztig isomorphism theorem in the super setting.

\begin{thm}[\cite{[Y2], [H], [A]}]\label{thm:lusztig-isomorphism-super}
For every $a\in A$ and $i\in I$, there are algebra isomorphisms
\begin{equation*}
T_i^a:U_q^{\DK}(\fg_s)_a\xrightarrow{\sim}U_q^{\DK}(\fg_s)_{\rho_i(a)}, \qquad (T_i^a)^{-1}:U_q^{\DK}(\fg_s)_{\rho_i(a)}\xrightarrow{\sim}U_q^{\DK}(\fg_s)_a.
\end{equation*}
\end{thm}

Next, we describe $T_i^a(E_j^a)$ for $j\neq i$. Let us fix $i\in I$, and for a homogeneous element $X$ of degree $\gamma$, we define the braided adjoint operator
\begin{equation}\label{eq:braided-adjoint-super}
\operatorname{ad}_c(E_i^a)(X) :=E_i^aX-\chi(\alpha_i^a,\gamma)XE_i^a.
\end{equation}
For $j\neq i$, we define the super root-string vectors inductively by
\begin{equation}\label{eq:recursive-root-string-vector}
E_{j,0}^a:=E_j^a, \qquad E_{j,r+1}^a :=E_i^aE_{j,r}^a -\chi(\alpha_i^a,\alpha_j^a+r\alpha_i^a)E_{j,r}^aE_i^a.
\end{equation}
In other words, $E_{j,r}^a=(\operatorname{ad}_cE_i^a)^r(E_j^a)$. We choose the Lusztig isomorphisms as in \cite[Theorem 2.8 and Proposition 3.4]{[A]}, and then for $j\neq i$, $T_i^a(E_j^a)$ is a non-zero scalar multiple of
\begin{equation}\label{eq:local-Lusztig-root-vector}
E_{j,-c_{ij}^a}^{\rho_i(a)} =\bigl(\operatorname{ad}_cE_i^{\rho_i(a)}\bigr)^{-c_{ij}^a} (E_j^{\rho_i(a)}).
\end{equation}

\subsubsection{Construction of super root vectors}

Roughly speaking, as in the case of Weyl groups, a reduced expression of a morphism of maximal length determines an enumeration of the positive super roots, and by applying the Lusztig isomorphisms, we obtain the corresponding super root vectors. More precisely, let us fix a target object $a\in A$, and let $n=|\Delta_+^a|$. We choose a morphism of maximal length ending at $a$, and a reduced expression
\begin{equation}\label{eq:reduced-maximal-morphism-super}
w= \sigma_{i_n}^{a_{n-1}}\circ\cdots\circ \sigma_{i_2}^{a_1}\circ\sigma_{i_1}^{a_0}:a_0\longrightarrow a_n=a, \qquad a_r=\rho_{i_r}(a_{r-1}).
\end{equation}
For $1\leq j\leq n$, we denote the tail of this word by
\begin{equation*}
w_{j+1,n} :=\sigma_{i_n}^{a_{n-1}}\circ\cdots\circ\sigma_{i_{j+1}}^{a_j}:a_j\longrightarrow a,
\end{equation*}
with $w_{n+1,n}=\id_a$, and we define the $j$-th super root associated with the reduced word by
\begin{equation}\label{eq:positive-root-from-tail}
\beta_j:=w_{j+1,n}(\alpha_{i_j}^{a_j}) =\sigma_{i_n}^{a_{n-1}}\cdots\sigma_{i_{j+1}}^{a_j}(\alpha_{i_j}^{a_j}).
\end{equation}
Then, we have
\begin{equation}\label{eq:positive-roots-enumeration}
\Delta_+^a=\{\beta_1,\ldots,\beta_n\},
\end{equation}
without repetition. Furthermore, the order $\beta_1<\cdots<\beta_n$ is convex, i.e., whenever $\beta_r+\beta_t$ is a super root with $r<t$, its position lies strictly between $r$ and $t$.

Now, we define the super root vectors by applying the same tail of Lusztig isomorphisms as in \eqref{eq:positive-root-from-tail} to the simple generator $E_{i_j}^{a_j}$ of degree $\alpha_{i_j}^{a_j}$, i.e.,
\begin{equation}\label{eq:PBW-root-vector-construction}
E_{\beta_j} :=T_{i_n}^{a_{n-1}}\circ\cdots\circ T_{i_{j+1}}^{a_j} \bigl(E_{i_j}^{a_j}\bigr), \qquad 1\leq j\leq n.
\end{equation}

\subsubsection{PBW bases}

For $\ba=(r_1,\ldots,r_n)\in\BN^n$, we define the ordered PBW monomial
\begin{equation*}
\bE_{\ba}:=E_{\beta_1}^{r_1}E_{\beta_2}^{r_2}\cdots E_{\beta_n}^{r_n}.
\end{equation*}
For every positive super root $\beta\in \Delta_+^a$, we denote
\begin{equation}\label{eq:root-order-N-beta}
N_\beta:=\begin{cases} \ord\bigl(\chi(\beta,\beta)\bigr),&\chi(\beta,\beta)\text{ is a non-trivial root of unity},\\
\infty,&\text{otherwise}.
\end{cases}
\end{equation}

The following PBW theorem for the De Concini--Kac positive half is proved in \cite{[Y1]} in the generic case, and at roots of unity, it follows from Lemma~\ref{lem:saturated-DK-pre-Nichols}, Remark~\ref{rem:isolated-odd-factor}, and \cite[Theorem 3.6]{[A]}.

\begin{thm}\label{thm:DK-PBW-super}
The set
\begin{equation}\label{eq:DK-PBW-basis-super}
\left\{ \bE_{\ba}\ \middle|\ r_j\in\BN,\quad 0\leq r_j<N_{\beta_j} \text{ whenever }\beta_j\in\Delta^a_{+,\iiso} \right\}
\end{equation}
is a basis of $U_q^{\DK}(N_s^+)_a$.
\end{thm}

Similarly, the corresponding PBW basis of a Nichols algebra is given by \cite[Theorem 2.9]{[A]}.
\begin{thm}\label{thm:small-PBW-super}
The set
\begin{equation}\label{eq:small-PBW-basis-super}
\left\{ \bE_{\ba}\ \middle|\ 0\leq r_j<N_{\beta_j}\text{ for every }j \right\}
\end{equation}
is a basis of $u_q^{\sm}(N_s^+)_a$.
\end{thm}

Note that, if $\beta$ is an isotropic odd super root, then $\chi(\beta,\beta)=-1$ for the braiding in \eqref{eq:full-bicharacter-super}, and hence $N_\beta=2$, i.e., its exponent in the small PBW basis is $0$ or $1$. On the other hand, if $\beta$ is a non-isotropic super root with $N_\beta<\infty$, then the element $E_\beta^{N_\beta}$ is non-zero in the De Concini--Kac positive half, and it is zero in the small quotient. Namely, by comparing the two PBW bases, we see that the difference between the De Concini--Kac and the small positive halves comes only from the root powers of the non-isotropic super roots, which we study in the next subsection.

\subsection{Root powers and the De Concini--Kac center}

In this subsection, we study the root powers $E_\beta^{N_\beta}$, and find a condition under which they generate a super-central Hopf subalgebra, which is the analogue of the De Concini--Kac center. We continue to assume that $q$ avoids small torsion and satisfies the isotropic non-Cartan condition, and we only consider the root powers $E_\beta^{N_\beta}$ for the super roots with $N_\beta<\infty$. One has the following commutation property as in the usual quantum group, cf.~\cite{[A]}, which is obtained by reducing to the simple super root case and then calculating directly.

\begin{lem}\label{cor:root-power-commutation}
For every $\beta\in\Delta^a_{+,\niso}$ and every homogeneous $X\in U_q^{\DK}(N_s^+)_a$ of degree $\gamma$,
\begin{equation}\label{eq:root-power-commutation}
E_\beta^{N_\beta}X =\chi(N_\beta\beta,\gamma)XE_\beta^{N_\beta} =(-1)^{p(N_\beta\beta)p(\gamma)} b'(N_\beta\beta,\gamma)XE_\beta^{N_\beta}.
\end{equation}
\end{lem}

In particular, $E_\beta^{N_\beta}$ is super-central in $U_q^{\DK}(N_s^+)_a$ if and only if $b'(N_\beta\beta,\gamma)=1$ for every $\gamma\in\Lambda^{\mathrm{pos}}$, which motivates the following condition.

\begin{condition}\label{cond:star}
The pair $(q,b')$ satisfies condition $(*)$ if:
\begin{enumerate}
\item $q$ avoids small torsion;
\item $q$ satisfies the isotropic non-Cartan condition;
\item for every $\beta\in\Delta^a_{+,\niso}$ with $N_\beta<\infty$ and every $\lambda\in\Lambda$,
\begin{equation}\label{eq:centrality-condition-b-prime}
b'(N_\beta\beta,\lambda) =b'(\lambda,N_\beta\beta)=1.
\end{equation}
\end{enumerate}
\end{condition}

\subsubsection{The root-power subalgebra} We denote by
\begin{equation}\label{eq:root-power-subalgebra}
Z_q^+(N_s^+)_a \subset U_q^{\DK}(N_s^+)_a
\end{equation}
the subalgebra generated by
\begin{equation*}
E_\beta^{N_\beta}, \qquad \beta\in \Delta^a_{+,\niso}.
\end{equation*}
The following proposition says that, under condition $(*)$, it plays the role of the De Concini--Kac center, i.e., the small positive half is the quotient of the De Concini--Kac positive half by the ideal generated by the root powers.

\begin{prop}\label{prop:DK-center-positive}
Assume condition $(*)$. Then, $Z_q^+(N_s^+)_a$ is a super-commutative Hopf subalgebra contained in the super-center of $U_q^{\DK}(N_s^+)_a$, and it is independent of the reduced expression used to define the super root vectors. Furthermore, the kernel of the canonical projection
\begin{equation*}
U_q^{\DK}(N_s^+)_a\twoheadrightarrow u_q^{\sm}(N_s^+)_a
\end{equation*}
is the Hopf ideal generated by the elements $E_\beta^{N_\beta}$ for $\beta\in \Delta^a_{+,\niso}$. Equivalently, there is an exact sequence of braided Hopf superalgebras
\begin{equation}\label{Frob quant}
\BC\longrightarrow Z_q^+(N_s^+)_a \longrightarrow U_q^{\DK}(N_s^+)_a \longrightarrow u_q^{\sm}(N_s^+)_a \longrightarrow\BC.
\end{equation}
\end{prop}
\begin{proof}
At roots of unity, the assertion follows from \cite[Proposition 4.1, Theorems 4.10 and 4.13, and Corollary 4.14]{[A]}, using Lemma~\ref{lem:saturated-DK-pre-Nichols} and Remark~\ref{rem:isolated-odd-factor}. In the generic case, by Lemma~\ref{lem:saturated-DK-pre-Nichols}, we have $U_q^{\DK}(N_s^+)_a\simeq u_q^{\sm}(N_s^+)_a$ and $Z_q^+(N_s^+)_a=\BC$.
\end{proof}

\subsection{Factorization algebras attached to quantum supergroups} In this subsection, we recall the construction of the factorization algebra attached to a Hopf superalgebra in the quantum torus category $\Rep_q(T_s)$. We essentially repeat the content in \cite[Section 3.4]{[G7]} and use the same conventions as in loc.\ cit., see also \cite[Sections 29.3--29.5]{[GL]} and \cite[Section 6]{[CF]}. In particular, the construction is contravariant. Let us take $C=\BP^1$, and we will use the construction on a coordinate open subset $U\simeq\BA^1$ of $C$, with coefficients in $\SVect$ as in \cite[Remark 2.4.6]{[TY]}. Let
\begin{equation*}
H=\bigoplus_{\eta\in\Lambda^{\mathrm{pos}}}H_{\eta}
\end{equation*}
be a $\Lambda^{\mathrm{pos}}$-graded Hopf superalgebra in $\Rep_q(T_s)$, with $H_0=\BC$ and finite-dimensional graded pieces. Then, the construction assigns to $H$ a factorization algebra
\begin{equation*}
\Omega_{H,U}\in\SD_{\cG^\Lambda}(\Conf_U),
\end{equation*}
where $\Conf_U$ is the configuration space of $U$, colored by $\Lambda^{\mathrm{neg}}\setminus\{0\}$.

\subsubsection{} Now, we define the factorization algebras on the configuration space of $C$ used in the paper, following \cite[Sections 1.6--1.7]{[G7]}. Roughly speaking, they all have the same restriction to the open locus $\Conf^\circ$ below, and they differ only in the way they are extended to the whole configuration space. More precisely, let $\Conf^\circ\subset\Conf$ be the open locus of divisors whose coefficients are negative simple super roots, and let $\jmath$ denote its inclusion. We denote by $\Omega_q^\circ$ the common restriction to $\Conf^\circ$ of the factorization algebras in Definition~\ref{def fact alg}.

Under the factorization trivialization of $\cG^\Lambda$ on $\Conf^\circ$, we have
\[
\left.\Omega_q^\circ\right|_{\Conf^{-\alpha_i^a}} \simeq \Pi^{p_i^a}\omega_C[-1].
\]
Note that the restrictions to the other components of $\Conf^\circ$ are determined by factorization, as in \cite[Definition 2.4.1 and Remark 2.4.2]{[TY]}, see also \cite{[BFT1],[G]}. In particular, for $n\geq1$, its restriction to $\Conf^{-n\alpha_i^a}\cap\Conf^\circ$ is $\Pi^n$ times the perverse constant local system if $\alpha_i^a$ is odd and is the perverse sign local system if $\alpha_i^a$ is even. Here, the sign local system corresponds to the sign representation of $S_n$ under the covering which orders the $n$ points. Indeed, when we interchange two factors $\Pi^{p_i^a}\omega_C[-1]$, the cohomological shift contributes $-1$ and the super symmetry contributes $(-1)^{p_i^a}$, and thus, the two signs cancel for an odd simple super root, while the sign $-1$ remains for an even simple super root.

We set
\begin{equation*}
\Omega_{q,\quant}^{\sm}:=\jmath_{!*}\Omega_q^\circ, \qquad \Omega_{q((\hbar)),\quant}^{\DK}:= \jmath_{!*}\Omega_{q((\hbar))}^\circ.
\end{equation*}
Here, $\Omega_{q[[\hbar]]}^\circ$ is defined on $\Conf^{-\alpha_i^a}$ by the same formula and extended to $\Conf^\circ$ by factorization, using the bilinear form $b'_\hbar$ in Section~\ref{deform algebra}, and $\Omega_{q((\hbar))}^\circ$ is obtained by extending the coefficients to $\BC((\hbar))$. Note that $\Omega_{q,\quant}^{\sm}$ corresponds to the small positive half, and the De Concini--Kac factorization algebra is obtained from $\Omega_{q((\hbar)),\quant}^{\DK}$ by a procedure parallel to the definition of $U_q^{\DK}(N_s^+)_a$, as follows.

\subsubsection{\texorpdfstring{$\hbar$}{hbar}-deformation of the De Concini--Kac factorization algebra} Following \cite[Section 1.7.7]{[G7]}, we define the $\hbar$-deformation of the De Concini--Kac factorization algebra by
\begin{equation}\label{def h-deform}
\Omega_{q[[\hbar]],\quant}^{\DK}:= \jmath_*\bigl(\Omega_{q[[\hbar]]}^\circ\bigr) \cap \Omega_{q((\hbar)),\quant}^{\DK} \subset \jmath_*\bigl(\Omega_{q((\hbar))}^\circ\bigr).
\end{equation}
Here, as in \cite[Sections 1.7.2 and 1.7.6--1.7.7]{[G7]}, we regard the twisted perverse sheaves with $\BC((\hbar))$-coefficients as objects in the ind-completion of the category of twisted perverse sheaves with $\BC[[\hbar]]$-coefficients, and the intersection is taken in this category, inside $\jmath_*(\Omega_{q((\hbar))}^\circ)$. According to \cite[Proposition A.2.7]{[G7]}, the intersection is an $\hbar$-flat twisted perverse sheaf with $\BC[[\hbar]]$-coefficients. Then, we set
\begin{equation*}
\Omega_{q,\quant}^{\DK}:= \Omega_{q[[\hbar]],\quant}^{\DK}\otimes_{\BC[[\hbar]]}\BC, \qquad \Omega_{q,\quant}^{\Lus}:= \mathbb D\bigl(\Omega_{q^{-1},\quant}^{\DK}\bigr).
\end{equation*}
Note that these constructions commute with the restriction to open subsets of the curve, and since the external products on disjoint supports are $\hbar$-flat, the intersections in \eqref{def h-deform} commute with factorization, and thus, they are factorization algebras on the configuration space of $C$.

Next, we compare these factorization algebras with the Hopf superalgebras in \eqref{eq:five-positive-halves-super}, and we omit the subscript $a$ when no confusion arises. First, on $U\simeq\BA^1$, we have $\Omega_{q,\quant}^{\sm}|_{\Conf_U}\simeq\Omega_{u_q^{\sm}(N_s^+),U}$ by \cite[Corollary 3.5.8]{[G7]}. Furthermore, by the identification after inverting $\hbar$ and Lemma~\ref{lem:saturated-DK-pre-Nichols}, the factorization algebra attached to the Hopf superalgebra $U_{q[[\hbar]]}^{\DK}(N_s^+)_a$ over $\BC[[\hbar]]$ is identified with \eqref{def h-deform}, inside $\jmath_*\Omega_{q((\hbar))}^\circ$. Indeed, the map from the factorization algebra attached to $U_{q[[\hbar]]}^{\DK}(N_s^+)_a$ to \eqref{def h-deform} is an isomorphism after inverting $\hbar$, and a non-zero $\hbar$-torsion cokernel would produce a kernel after the specialization, which contradicts the injection $\Omega_{U_q^{\DK}(N_s^+),U}\hookrightarrow \jmath_*\Omega_q^\circ$. Here, we follow the argument in \cite[Sections 3.6.3--3.6.7]{[G7]}, with the super PBW basis above in place of the usual one. Then, by Verdier duality, we obtain the identification for the Lusztig form. Thus, we have
\begin{equation*}
\Omega_{q,\quant}^{?}|_{\Conf_U}\simeq\Omega_{H^?,U}, \quad H^{\DK}=U_q^{\DK}(N_s^+),\quad H^{\sm}=u_q^{\sm}(N_s^+),\quad H^{\Lus}=U_q^{\Lus}(N_s^+).
\end{equation*}
Since the construction is contravariant, the Hopf maps in \eqref{eq:five-positive-halves-super} induce
\begin{equation*}
\Omega_{q,\quant}^{\Lus}\longrightarrow\Omega_{q,\quant}^{\sm}\longrightarrow\Omega_{q,\quant}^{\DK},
\end{equation*}
which are the maps induced by the identity on $\Omega_q^\circ$. Note that, in the perverse heart, the first map is an epimorphism and the second one is a monomorphism, and if $\chi(\beta,\beta)$ has infinite order for every non-isotropic super root $\beta$, then both maps are isomorphisms.

\subsection{The classical factorization algebra and quantum Frobenius}\label{sec:classical-frobenius}

In this subsection, we recall the quantum Frobenius map and the action of the classical factorization algebra on $\Omega_{q,\quant}^{\Lus}$, following \cite[Section 4]{[G7]}. We assume that $q$ is torsion valued on $Q_s$ and that condition $(*)$ holds. Then, we define the quantum Frobenius lattice by
\begin{equation}\label{eq:quantum-Frobenius-lattice-super}
\Lambda^\# :=\operatorname{Span}_{\BZ}\left\{N_\beta\beta\ \middle|\ \beta\in \Delta^a_{+,\niso},\ N_\beta<\infty\right\} \subset Q_s\subset\Lambda.
\end{equation}
We denote by $\Lambda^\sharp:=\{\eta\in\Lambda\mid b(\eta,\Lambda)=1\}$ the metaplectic lattice. In particular, we have $\Lambda^\#\subset\Lambda^\sharp$. We denote by $\Conf^\#$ the configuration space with non-zero colors in the monoid generated by $-N_\beta\beta$, $\beta\in\Delta^a_{+,\niso}$, and the inclusion $\Lambda^\#\hookrightarrow\Lambda$ induces a map $F:\Conf^\#\to\Conf$.

We fix a trivialization of $F^*\cG^\Lambda$ compatible with factorization, as in \cite[Section 4.1.2]{[G7]}. Here, for the gerbes considered in this paper, its existence follows from the description of enhanced $\Theta$-data in \cite[Lemma 5.1 and Sections 5.3.5--5.3.6]{[Zha1]}, since $\Lambda^\#\subset Q_G$ (the geometric coroot lattice), and the determinant term $2q(\gamma)+b(\gamma,2\rho)$ vanishes for $\gamma\in\Lambda^\#$ (in additive notation). Note that, by condition $(*)$, we also have $b'(\Lambda^\#,\Lambda)=b'(\Lambda,\Lambda^\#)=1$.

Let $T_s^\#$ be the torus with character lattice $\Lambda^\#$. Then, the inclusion $\Lambda^\#\hookrightarrow\Lambda$ defines a map
\begin{equation*}
\Fr_{q,T_s}:T_s\longrightarrow T_s^\#.
\end{equation*}

\subsubsection{The Frobenius factorization algebra} Following \cite[Theorem 3.9]{[AAR]}, we denote by $U(\fn_s^\#)$ the componentwise graded dual of the central Hopf subalgebra in $U_{q^{-1}}^{\DK}(N_s^-)$ in Proposition~\ref{prop:DK-center-positive}, where $\fn_s^\#$ is the Lie algebra of primitive elements in this dual Hopf algebra. By dualizing \eqref{Frob quant} for the negative half, we obtain
\begin{equation}\label{short ex seq Lus}
u_q^{\sm}(N_s^+) \longrightarrow U_q^{\Lus}(N_s^+) \overset{\Fr_{q,N_s}}{\longrightarrow} \Fr_{q,T_s}^*U(\fn_s^\#).
\end{equation}
Here, $\Fr_{q,N_s}$ is cocentral, and the small positive half is its algebra of coinvariants, see \cite[Remark 4.17]{[A]}. Note that the composition is the augmentation followed by the unit, and by the PBW decomposition, we have
\begin{equation}\label{eq:Frobenius-coinvariants-super}
U_q^{\Lus}(N_s^+)\otimes_{u_q^{\sm}(N_s^+)}\BC \xrightarrow{\sim} \Fr_{q,T_s}^*U(\fn_s^\#).
\end{equation}

For a non-unital commutative algebra $A$ graded by the non-zero colors of $\Conf^\#$, we denote by $\operatorname{Fact}(A)$ the factorization algebra on $\Conf^\#$ constructed in \cite[Section 4.3.1]{[G7]}. We denote by $C^\bullet(\fn_s^\#)$ the Chevalley cochain algebra of $\fn_s^\#$, and by applying $\operatorname{Fact}$ to its augmentation ideal, we define
\begin{equation*}
\Omega^{\#,\mathrm{cl}}:=\operatorname{Fact}\bigl(\ker(C^\bullet(\fn_s^\#)\to\BC)\bigr)
\end{equation*}
on $\Conf^\#$. Note that, on a coordinate open subset $U\simeq\BA^1$, it is identified with $\Omega_{U(\fn_s^\#),U}$.

The following lemma follows from \cite[Theorems 3.7 and 3.9]{[AAR]} and \cite[Sections 4.1.1 and 4.3.3--4.3.4]{[G7]}.
\begin{lem}\label{lem:even-Frobenius-identification}
The Lie algebra $\fn_s^\#$ is purely even and, after choosing simple root vectors, identifies with the positive nilpotent Lie algebra of the metaplectic dual group of $G$. Thus, $\Omega^{\#,\mathrm{cl}}$ is the classical factorization algebra for $G$ defined in \cite[Remark 4.3.4]{[G7]}.
\end{lem}

\subsubsection{The Frobenius action} In order to state the corollary below, we recall the following notions. Let $A$ be a factorization algebra on $\Conf^\#$, equipped with an associative multiplication for convolution by addition of divisors, as in \cite[Sections 4.1.3 and 4.2.1]{[G7]}. Let $K$ be a factorization algebra on $\Conf$ with an action of $A$ in the sense of \cite[Section 4.2.7]{[G7]}, and we denote by $\operatorname{Bar}(A,K)$ the bar construction of \cite[Remark 4.4.4]{[G7]}. Also, for a factorization algebra $K$, we denote
\[
K_C:=\bigoplus_{\gamma\ne0}\iota_\gamma^!K, \qquad \iota_\gamma:C\longrightarrow\Conf^\gamma,\quad x\longmapsto\gamma\cdot x,
\]
where $\gamma$ ranges over its non-zero configuration degrees, and similarly for $A$ on $\Conf^\#$. Here, the tensor products over $A_C$ are taken over $A_C\oplus\omega_C$, as in \cite[Remark 4.2.5]{[G7]}.

Then, we have the following corollary, which is the super analogue of \cite[Theorem 4.4.2]{[G7]}. Roughly speaking, it says that the small factorization algebra is obtained from the Lusztig one by taking the derived coinvariants for the action of the classical factorization algebra.
\begin{cor}\label{cor:frob-action-factorization-super}
There is a canonical action of $\Omega^{\#,\mathrm{cl}}$ on $\Omega_{q,\quant}^{\Lus}$. The morphism $\Omega_{q,\quant}^{\Lus}\to\Omega_{q,\quant}^{\sm}$ is compatible with the trivial, or augmentation, action on $\Omega_{q,\quant}^{\sm}$, and it induces an isomorphism of factorization algebras
\begin{equation}\label{eq:bar-small-factorization-super}
\operatorname{Bar}\bigl(\Omega^{\#,\mathrm{cl}},\Omega_{q,\quant}^{\Lus}\bigr) \xrightarrow{\sim}\Omega_{q,\quant}^{\sm}.
\end{equation}
On each main diagonal, the left-hand side can be written as
\begin{equation}\label{eq:cone-small-factorization-super}
\operatorname{Cone}\bigl((\Omega^{\#,\mathrm{cl}})_C\longrightarrow(\Omega_{q,\quant}^{\Lus})_C\bigr) \otimes_{(\Omega^{\#,\mathrm{cl}})_C}\omega_C \xrightarrow{\sim}(\Omega_{q,\quant}^{\sm})_C.
\end{equation}
\end{cor}
\begin{proof}
We can apply the same proof as in \cite[Theorem 4.4.2 and Sections 4.5--4.6]{[G7]} to \eqref{short ex seq Lus} and \eqref{eq:Frobenius-coinvariants-super}, using Lemma~\ref{lem:even-Frobenius-identification}. Note that $\Lambda^\#$ is even and, by Condition~\ref{av s o}, it contains no simple super roots, as required in the argument of \cite[Sections 4.4.9--4.4.11]{[G7]}. Finally, the last formula follows from the description in \cite[Remark 4.4.4]{[G7]}.
\end{proof}

\subsection{The mixed quantum supergroup}

In this subsection, we define the mixed quantum supergroup category in the sense of Gaitsgory. First, we define
\begin{equation*}
\Rep_q(B_s^+)_a :=U_q^{\Lus}(N_s^+)_a\modu^{\ren}\bigl(\Rep_q(T_s)\bigr)
\end{equation*}
to be the ind-completion of the category of finite-dimensional modules over the Lusztig positive half. Note that the coproduct on $U_q^{\Lus}(N_s^+)_a$ makes $\Rep_q(B_s^+)_a$ monoidal relative to $\Rep_q(T_s)$.

\begin{defn}
The category $\cO$ of the mixed quantum supergroup is the relative Drinfeld center
\begin{equation}\label{eq:mixed-quantum-supergroup-center}
\cO_q^{\mxd}(\fg_s)_a :=Z_{\Dr,\Rep_q(T_s)}\bigl(\Rep_q(B_s^+)_a\bigr).
\end{equation}
\end{defn}

In the non-super case, this definition is due to Gaitsgory, see \cite[Section 5.3.1]{[G8]}, and Situ and Losev studied the blocks of this category at roots of unity and their relation to the affine Hecke category, see \cite{[Sit],[Los3]}.

Roughly speaking, an object of the relative center is a $U_q^{\Lus}(N_s^+)_a$-module together with a compatible half-braiding relative to the torus category, and by unwinding this half-braiding, we obtain an action of the negative De Concini--Kac half. Thus, $\cO_q^{\mxd}(\fg_s)_a$ can be described as the cocomplete category generated by Verma modules for the mixed quantum supergroup, whose positive part is $U_q^{\Lus}(N_s^+)_a$ and whose negative part is $U_q^{\DK}(N_s^-)_a$, and the local nilpotence condition is imposed on the positive Lusztig part. Here, the word ``mixed'' refers to the fact that the two nilpotent halves use different integral forms.

\begin{rem}
We expect that the locally compact Iwahori Gaiotto category corresponds to a full subcategory $\Rep_q^{\mxd}(\fg_s)$ of $\cO_q^{\mxd}(\fg_s)_a$. As the notation suggests, $\cO_q^{\mxd}(\fg_s)_a$ depends on the choice of the super Borel, since the definition of the Verma modules depends on $a$, but $\Rep_q^{\mxd}(\fg_s)$ is independent of the choice. 
\end{rem}

\subsubsection{Factorization modules} By the general Hopf-factorization formalism, we have the following description of the category $\cO$ of the mixed quantum supergroup in terms of factorization modules, which is the super analogue of the constructions in \cite{[BFS],[GL],[G7],[CF],[L]}.
\begin{thm}\label{thm:mixed-O-factorization-modules}
There is a $t$-exact equivalence
\begin{equation}
\cO_q^{\mxd}(\fg_s)_a\simeq\Omega_{q,\quant}^{\Lus}\FS
\end{equation}
whose restriction to the hearts is an equivalence of highest weight categories. Under this equivalence, Verma modules correspond to standard factorization modules, while contragredient Verma modules correspond to costandard factorization modules, with the parity convention of Remark~\ref{rem:target-parity}.
\end{thm}
\subsection{Super Borel}\label{super Borel} In this subsection, we will explicitly describe the super Borels that will correspond to the geometric models chosen in the paper. Note that, under the identifications below, the positive even super roots correspond to the positive geometric coroots for $B^-$ in all three cases.

We take the symmetric bilinear form
\begin{equation*}
(\varepsilon_i,\varepsilon_j)=\delta_{ij},\qquad (\delta_p,\delta_q)=-\delta_{pq},\qquad (\varepsilon_i,\delta_p)=0.
\end{equation*}
For the descriptions of the super root systems and the conventions for arbitrary simple systems, we follow \cite{[Kac],[Y1],[CHW]}. For a chosen simple system $\Pi_a=\{\alpha_i^a\}_{i\in I}$, we denote
\begin{equation*}
A^a:=\bigl((\alpha_i^a,\alpha_j^a)\bigr)_{i,j\in I}.
\end{equation*}
Namely, $A^a$ is the symmetric Gram matrix of the simple super roots, which is not the generalized Cartan matrix of the Cartan scheme (whose diagonal entries are all $2$). In the diagrams below, a white node denotes an even simple super root, a crossed node denotes an isotropic odd simple super root, and a single, double, or triple edge records the absolute value $1$, $2$, or $3$ of the corresponding non-zero off-diagonal entry of $A^a$.

\tikzset{ evenroot/.style={circle,draw,fill=white,minimum size=5.5mm,inner sep=0pt,line width=0.5pt}, isotropicroot/.style={circle,draw,fill=white,minimum size=5.5mm,inner sep=0pt,line width=0.5pt, path picture={ \draw[line width=0.45pt] (path picture bounding box.north west)--(path picture bounding box.south east); \draw[line width=0.45pt] (path picture bounding box.south west)--(path picture bounding box.north east); }}, singlebond/.style={line width=0.55pt}, doublebond/.style={double,double distance=1.4pt,line width=0.45pt}, triplebond/.style={double,double distance=2.4pt,line width=0.45pt} }

\subsubsection{\texorpdfstring{$\mathfrak{gl}(M|N)$ ($1\leq M\leq N$)}{gl(M|N)}}

We denote by $\varepsilon_1,\ldots,\varepsilon_M,\delta_1,\ldots,\delta_N$ the standard basis of weights of the diagonal Cartan subgroup $T_s$, and for $a=1,2$, we denote by $e_i^{(a)}$ the coordinate cocharacters of $G_a$. We fix the identification
\begin{equation}\label{eq:typeA-duality-coordinates}
\varepsilon_i\longmapsto(-e_i^{(1)},0), \qquad \delta_j\longmapsto(0,e_j^{(2)}).
\end{equation}
For a weight $\lambda=(\lambda^{(1)},\lambda^{(2)})$ (on the representation side) with $\lambda^{(a)}=\sum_i\lambda_i^{(a)}e_i^{(a)}$, its expression in the super weight basis is
\begin{equation}\label{eq:typeA-weight-in-geometric-coordinates}
\lambda=-\sum_{i=1}^M\lambda_i^{(1)}\varepsilon_i +\sum_{j=1}^N\lambda_j^{(2)}\delta_j.
\end{equation}
Then, the super root system is given by
\begin{equation*}
\Phi_{\bar 0} = \{\varepsilon_i-\varepsilon_j\mid 1\leq i\neq j\leq M\} \sqcup \{\delta_p-\delta_q\mid 1\leq p\neq q\leq N\},
\end{equation*}
\begin{equation*}
\Phi_{\bar 1} = \{\pm(\varepsilon_i-\delta_p) \mid 1\leq i\leq M,\ 1\leq p\leq N\}.
\end{equation*}
We choose the ordered basis of $\Lambda\otimes_{\BZ}\BC$
\begin{equation*}
(x_1,\dots,x_{M+N}) = (\delta_1,\varepsilon_1,\delta_2,\varepsilon_2,\dots, \delta_M,\varepsilon_M,\delta_{M+1},\dots,\delta_N).
\end{equation*}
The corresponding positive super roots are
\begin{equation*}
\Phi^+ = \{x_b-x_a\mid 1\leq a<b\leq M+N\},
\end{equation*}
where the super root $x_b-x_a$ is even when $x_a$ and $x_b$ are both among the $\varepsilon_i$ or both among the $\delta_p$, and odd otherwise. The simple super roots are
\begin{equation*}
\alpha_r=x_{r+1}-x_r, \qquad 1\leq r\leq M+N-1.
\end{equation*}

In particular, if $M=N$, all the $2M-1$ simple super roots are isotropic odd, and if $M<N$, the first $2M$ simple super roots are isotropic odd, and the remaining $N-M-1$ simple super roots are even of square $-2$. Namely, the chosen super Borel is one with the maximal number of isotropic simple super roots.

Next, we compute the Weyl vectors, which are needed for the $\omega^\rho$-renormalization. The even and odd Weyl vectors are
\begin{equation*}
\rho_{\bar 0} = \frac12\sum_{i=1}^M(-M+2i-1)\,\varepsilon_i + \frac12\sum_{p=1}^N(-N+2p-1)\,\delta_p,
\end{equation*}
\begin{equation*}
\rho_{\bar 1} = \frac12\sum_{i=1}^M(-N+2i)\,\varepsilon_i + \frac12\sum_{p=1}^M(-M+2p-2)\,\delta_p + \frac M2\sum_{p=M+1}^N\delta_p.
\end{equation*}
Hence,
\begin{equation*}
\widetilde\rho = \rho_{\bar 0}-\rho_{\bar 1} = \frac{N-M-1}{2}\sum_{i=1}^M\varepsilon_i - \frac{N-M-1}{2}\sum_{p=1}^M\delta_p + \frac12\sum_{p=M+1}^N(-M-N-1+2p)\,\delta_p.
\end{equation*}
We denote
\begin{equation*}
\omega_{\mathrm{cen}} := \sum_{i=1}^M\varepsilon_i- \sum_{p=1}^N\delta_p.
\end{equation*}
Then, $(\omega_{\mathrm{cen}},\alpha_r)=0$ for every simple super root, and
\begin{equation*}
\widetilde\rho = \rho+ \frac{N-M-1}{2}\,\omega_{\mathrm{cen}}, \qquad \rho := \sum_{p=M+1}^N(-M-1+p)\,\delta_p.
\end{equation*}
Thus, $\rho$ is the representative of $\widetilde\rho$ modulo $\mathbb C\omega_{\mathrm{cen}}$ used for the $\omega^\rho$-renormalization in Section~\ref{Gaiotto category}.

The corresponding super Dynkin diagrams are as follows, where in small ranks, the empty parts of a chain are omitted. In particular, the even part of the second diagram is empty when $N=M+1$, and the first diagram consists of the single vertex $\alpha_1$ when $M=N=1$.
\begin{center}
\begin{tikzpicture}[baseline=(current bounding box.center)]
\node at (-1.1,0) {$M=N\colon$}; \node[isotropicroot,label=below:$\alpha_1$] (a1) at (0,0) {}; \node[isotropicroot,label=below:$\alpha_2$] (a2) at (1.2,0) {}; \node (dots) at (2.4,0) {$\cdots$}; \node[isotropicroot,label=below:$\alpha_{2M-2}$] (ap) at (3.6,0) {}; \node[isotropicroot,label=below:$\alpha_{2M-1}$] (aq) at (4.8,0) {}; \draw[singlebond] (a1)--(a2); \draw[singlebond] (a2)--(dots); \draw[singlebond] (dots)--(ap); \draw[singlebond] (ap)--(aq);
\end{tikzpicture}
\end{center}

\begin{center}
\begin{tikzpicture}[baseline=(current bounding box.center)]
\node at (-1.45,0) {$N\geq M+1\colon$}; \node[isotropicroot,label=below:$\alpha_1$] (a1) at (0,0) {}; \node[isotropicroot,label=below:$\alpha_2$] (a2) at (1.2,0) {}; \node (dots1) at (2.4,0) {$\cdots$}; \node[isotropicroot,label=below:$\alpha_{2M}$] (am) at (3.6,0) {}; \node[evenroot,label=below:$\alpha_{2M+1}$] (b1) at (4.8,0) {}; \node (dots2) at (6.0,0) {$\cdots$}; \node[evenroot,label=below:$\alpha_{M+N-1}$] (bn) at (7.2,0) {}; \draw[singlebond] (a1)--(a2); \draw[singlebond] (a2)--(dots1); \draw[singlebond] (dots1)--(am); \draw[singlebond] (am)--(b1); \draw[singlebond] (b1)--(dots2); \draw[singlebond] (dots2)--(bn);
\end{tikzpicture}
\end{center}

\begin{equation*}
A^a_{\mathfrak{gl}(M|N)} =(a^a_{ij})_{1\leq i,j\leq M+N-1}, \qquad a^a_{rr}=
\begin{cases}
0,&r\leq\min(2M,M+N-1),\\
-2,&r>2M,
\end{cases}
\end{equation*}
\begin{equation*}
a^a_{r,r+1}=a^a_{r+1,r} =
\begin{cases}
(-1)^r,&r<2M,\\
1,&r\geq2M,
\end{cases}
\qquad a^a_{ij}=0\quad\text{if }|i-j|>1.
\end{equation*}

The corresponding Cartan-scheme matrix in \eqref{eq:root-string-cartan-integers} is given by \begin{equation}\label{eq:Cartan-scheme-gl} C^a_{\mathfrak{gl}(M|N)} =
\begin{pmatrix}
2&-1&0&\cdots&0\\
-1&2&-1&\ddots&\vdots\\
0&-1&2&\ddots&0\\
\vdots&\ddots&\ddots&\ddots&-1\\
0&\cdots&0&-1&2
\end{pmatrix}_{r\times r},\qquad r=M+N-1.
\end{equation}

\subsubsection{\texorpdfstring{$\mathfrak{osp}(2k|2l)$}{osp(2k|2l)}}

We assume that $k,l\geq1$, and the geometric factors are ordered as in Section~\ref{pol des}. As above, we denote by $e_i^{(a)}$ the coordinate cocharacters of $G_a$, for $a=1,2$, and we identify
\begin{equation}\label{eq:osp-duality-coordinates}
\varepsilon_i=
\begin{cases}
(-e_i^{(1)},0),&k\leq l,\\
(0,e_i^{(2)}),&k>l,
\end{cases}
\qquad \delta_p=
\begin{cases}
(0,e_p^{(2)}),&k\leq l,\\
(-e_p^{(1)},0),&k>l.
\end{cases}
\end{equation}
Then, the super root system is given by
\begin{equation*}
\Phi_{\bar0} = \{\pm\varepsilon_i\pm\varepsilon_j\mid 1\leq i<j\leq k\} \sqcup \{\pm\delta_p\pm\delta_q\mid 1\leq p<q\leq l\} \sqcup \{\pm2\delta_p\mid 1\leq p\leq l\},
\end{equation*}
\begin{equation*}
\Phi_{\bar1} = \{\pm\varepsilon_i\pm\delta_p \mid 1\leq i\leq k,\ 1\leq p\leq l\}.
\end{equation*}
We choose the ordered basis of $\Lambda\otimes_{\BZ}\BC$
\begin{equation*}
(x_1,\dots,x_{k+l}) =
\begin{cases}
(\varepsilon_1,\dots,\varepsilon_{k-l}, \delta_1,\varepsilon_{k-l+1},\delta_2,\dots, \delta_l,\varepsilon_k),&k>l,\\[1mm]
(\delta_1,\dots,\delta_{l-k+1}, \varepsilon_1,\delta_{l-k+2},\varepsilon_2,\dots, \delta_l,\varepsilon_k),&k\leq l.
\end{cases}
\end{equation*}
Let $r:=k+l$. Then, the positive super roots determined by this order are
\begin{equation*}
\Phi^+ = \{x_j-x_i,-x_i-x_j\mid 1\leq i<j\leq r\} \sqcup \{-2\delta_p\mid 1\leq p\leq l\}.
\end{equation*}
Here, the super roots involving only the $\varepsilon_i$ or only the $\delta_p$ are even, and those involving one $\varepsilon_i$ and one $\delta_p$ are odd. The simple super roots are
\begin{equation*}
\alpha_t=x_{t+1}-x_t, \qquad 1\leq t\leq r-1,
\end{equation*}
\begin{equation*}
\alpha_r=-x_{r-1}-x_r.
\end{equation*}

The even and odd Weyl vectors are
\begin{equation*}
\rho_{\bar 0} = \sum_{i=1}^k(-k+i)\,\varepsilon_i + \sum_{p=1}^{l}(-l+p-1)\,\delta_p,
\end{equation*}
\begin{equation*}
\rho_{\bar 1} = -\sum_{i=1}^k\min(k-i,l)\,\varepsilon_i - \sum_{p=1}^{l}\min(k,l-p+1)\,\delta_p.
\end{equation*}
Hence,
\begin{equation*}
\rho = \rho_{\bar 0}-\rho_{\bar 1} =
\begin{cases}
\displaystyle \sum_{i=1}^{k-l-1}(i-k+l)\,\varepsilon_i,&k>l,\\[3mm]
\displaystyle \sum_{p=1}^{l-k}(p-l+k-1)\,\delta_p,&k\leq l.
\end{cases}
\end{equation*}

We set
\begin{equation*}
s :=
\begin{cases}
k-l-1,&k>l,\\
l-k,&k\leq l.
\end{cases}
\end{equation*}
Then, the first $s$ simple super roots are even, and the remaining $r-s$ simple super roots are isotropic odd. Furthermore, if $k>l$, the first $s$ simple super roots have square $2$, and if $k<l$, they have square $-2$.

For $s>0$, the super Dynkin diagram is
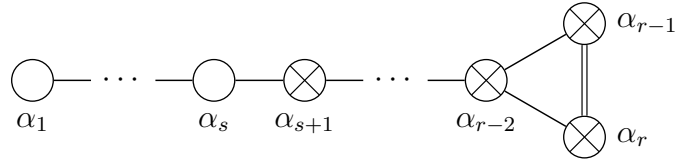
\begin{center}
\begin{tikzpicture}[baseline=(current bounding box.center)]
\node[evenroot,label=below:$\alpha_1$] (a1) at (0,0) {}; \node (dots1) at (1.2,0) {$\cdots$}; \node[evenroot,label=below:$\alpha_s$] (as) at (2.4,0) {}; \node[isotropicroot,label=below:$\alpha_{s+1}$] (b1) at (3.6,0) {}; \node (dots2) at (4.8,0) {$\cdots$}; \node[isotropicroot,label=below:$\alpha_{r-2}$] (bm) at (6.0,0) {}; \node[isotropicroot,label=right:$\alpha_{r-1}$] (bp) at (7.3,0.75) {}; \node[isotropicroot,label=right:$\alpha_r$] (bq) at (7.3,-0.75) {}; \draw[singlebond] (a1)--(dots1); \draw[singlebond] (dots1)--(as); \draw[singlebond] (as)--(b1); \draw[singlebond] (b1)--(dots2); \draw[singlebond] (dots2)--(bm); \draw[singlebond] (bm)--(bp); \draw[singlebond] (bm)--(bq); \draw[doublebond] (bp)--(bq);
\end{tikzpicture}
\end{center}
For $s=0$, the leftmost node is crossed. When $k=1<l$, the chain of crossed nodes before the fork is empty, and the two crossed end nodes are joined to the even node $\alpha_s$. For $k=l=1$, only the two crossed nodes joined by a double edge remain.

Let $A^a_{\mathfrak{osp}(2k|2l)} =(a^a_{ij})_{1\leq i,j\leq r}$, and for $r\geq3$, its entries are as follows.

If $k>l$, so that $s=k-l-1$, then
\begin{equation*}
a^a_{ii} =
\begin{cases}
2,&1\leq i\leq s,\\
0,&s+1\leq i\leq r,
\end{cases}
\end{equation*}
\begin{equation*}
a^a_{i,i+1}=a^a_{i+1,i} =
\begin{cases}
-1,&1\leq i\leq s,\\
(-1)^{i-s+1},&s+1\leq i\leq r-2.
\end{cases}
\end{equation*}

If $k\leq l$, so that $s=l-k$, then
\begin{equation*}
a^a_{ii} =
\begin{cases}
-2,&1\leq i\leq s,\\
0,&s+1\leq i\leq r,
\end{cases}
\end{equation*}
\begin{equation*}
a^a_{i,i+1}=a^a_{i+1,i} =
\begin{cases}
1,&1\leq i\leq s,\\
(-1)^{i-s},&s+1\leq i\leq r-2.
\end{cases}
\end{equation*}

In both cases, $ a^a_{r-2,r}=a^a_{r,r-2}=1$, $a^a_{r-1,r}=a^a_{r,r-1}=-2$, and all entries not specified above are $0$. For $k=l=1$, we have $A^a_{\mathfrak{osp}(2|2)} =
\begin{pmatrix}
0&-2\\
-2&0
\end{pmatrix}$.

The Cartan-scheme matrix $C^a_{\mathfrak{osp}(2k|2l)}=(c^a_{ij})$ is given, for $r\geq3$, by
\begin{equation}
\label{eq:Cartan-scheme-osp}\fitmath{c^a_{ii}=2,\qquad c^a_{i,i+1}=c^a_{i+1,i}=-1\ (1\leq i\leq r-2),\qquad c^a_{r-2,r}=c^a_{r,r-2} =c^a_{r-1,r}=c^a_{r,r-1}=-1,}
\end{equation}
with all other off-diagonal entries zero. For $k=l=1$, it is $\left(\begin{smallmatrix}2&-1\\-1&2\end{smallmatrix}\right)$.

\subsubsection{\texorpdfstring{$\mathfrak f(4)$}{f(4)}}

We follow the description in \cite{[FTY]}, and choose a basis of $\Lambda\otimes_{\BZ}\BC$ such that
\begin{equation*}
(\varepsilon_i,\varepsilon_j)_0=\delta_{ij}, \qquad (\delta,\delta)_0=-3, \qquad (\varepsilon_i,\delta)_0=0.
\end{equation*}
For the geometric group $G=\PGL_2\times\PSp(6)$, we identify
\begin{equation}\label{eq:F4-duality-coordinates}
\delta\longmapsto(-1;0,0,0),\qquad \varepsilon_j\longmapsto(0;e_j).
\end{equation}
Here, the first coordinate $a$ denotes the cocharacter $z\mapsto\diag(z^a,z^{-a})$ of $\PGL_2$ (it may be half-integral), and the last three coordinates are the standard diagonal coordinates for the coweights of $\PSp(6)$. Namely, a geometric coweight $(a;b_1,b_2,b_3)$ is the super weight $-a\delta+\sum_j b_j\varepsilon_j$, and in particular, the positive geometric coroot in the first factor is $-\delta$. Then, the super root system is
\begin{equation*}
\Phi_{\bar0} = \{\pm\varepsilon_i\mid 1\leq i\leq3\} \sqcup \{\pm\varepsilon_i\pm\varepsilon_j\mid 1\leq i<j\leq3\} \sqcup \{\pm\delta\},
\end{equation*}
\begin{equation*}
\Phi_{\bar1} = \left\{ \frac12(\pm\varepsilon_1\pm\varepsilon_2\pm\varepsilon_3\pm\delta) \right\},
\end{equation*}
where all four signs are independent. Note that every odd super root is isotropic. We choose the simple super roots
\begin{equation*}
\alpha_1 = \frac12(-\varepsilon_1+\varepsilon_2+\varepsilon_3-\delta), \qquad \alpha_2 = \frac12(\varepsilon_1-\varepsilon_2-\varepsilon_3-\delta),
\end{equation*}
\begin{equation*}
\alpha_3 = \frac12(-\varepsilon_1+\varepsilon_2-\varepsilon_3+\delta), \qquad \alpha_4=-\varepsilon_2+\varepsilon_3.
\end{equation*}
Then, the first three simple super roots are isotropic odd, and $\alpha_4$ is even. The positive even super roots are
\begin{equation*}
\begin{aligned}
\Phi^+_{\bar0} = \{&-\varepsilon_2+\varepsilon_3, -\varepsilon_3, -\varepsilon_1+\varepsilon_2, -\delta, -\varepsilon_2, -\varepsilon_1+\varepsilon_3,\\
&-\varepsilon_1, -\varepsilon_2-\varepsilon_3, -\varepsilon_1-\varepsilon_3, -\varepsilon_1-\varepsilon_2\},
\end{aligned}
\end{equation*}
and the positive odd super roots are
\begin{equation*}
\begin{aligned}
\Phi^+_{\bar1} = \frac12\{& \varepsilon_1-\varepsilon_2-\varepsilon_3-\delta, -\varepsilon_1+\varepsilon_2+\varepsilon_3-\delta,\\
&-\varepsilon_1+\varepsilon_2-\varepsilon_3+\delta, -\varepsilon_1-\varepsilon_2+\varepsilon_3+\delta,\\
&-\varepsilon_1+\varepsilon_2-\varepsilon_3-\delta, -\varepsilon_1-\varepsilon_2+\varepsilon_3-\delta,\\
&-\varepsilon_1-\varepsilon_2-\varepsilon_3+\delta, -\varepsilon_1-\varepsilon_2-\varepsilon_3-\delta \}.
\end{aligned}
\end{equation*}

The even and odd Weyl vectors are
\begin{equation*}
\rho_{\bar 0} = -\frac12(5\varepsilon_1+3\varepsilon_2+\varepsilon_3+\delta), \qquad \rho_{\bar 1} = -\frac12(3\varepsilon_1+\varepsilon_2+\varepsilon_3+\delta),
\end{equation*}
\begin{equation*}
\rho=\rho_{\bar 0}-\rho_{\bar 1}=-\varepsilon_1-\varepsilon_2.
\end{equation*}

To obtain an integral symmetric Gram matrix, we take $(\ ,\ )_{\mathrm{int}}=2(\ ,\ )_0$, i.e.,
\begin{equation*}
(\varepsilon_i,\varepsilon_j)_{\mathrm{int}}=2\delta_{ij}, \qquad (\delta,\delta)_{\mathrm{int}}=-6, \qquad (\varepsilon_i,\delta)_{\mathrm{int}}=0.
\end{equation*}

With respect to $(\ ,\ )_{\mathrm{int}}$, the super Dynkin diagram and Gram matrix for this non-distinguished simple system are as follows.
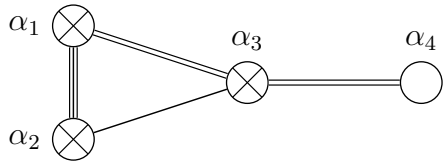
\begin{center}
\begin{tikzpicture}[baseline=(current bounding box.center)]
\node[isotropicroot,label=left:$\alpha_1$] (a1) at (0,1.5) {}; \node[isotropicroot,label=left:$\alpha_2$] (a2) at (0,0) {}; \node[isotropicroot,label=above:$\alpha_3$] (a3) at (2.3,0.75) {}; \node[evenroot,label=above:$\alpha_4$] (a4) at (4.6,0.75) {}; \draw[triplebond] (a1)--(a2); \draw[singlebond] (a1)--(a2); \draw[doublebond] (a1)--(a3); \draw[singlebond] (a2)--(a3); \draw[doublebond] (a3)--(a4);
\end{tikzpicture}
\end{center}

By a direct calculation, we have
\begin{equation*}
A^a_{\mathfrak f(4)} =
\begin{pmatrix}
0&-3&2&0\\
-3&0&1&0\\
2&1&0&-2\\
0&0&-2&4
\end{pmatrix}.
\end{equation*}
In the notation of Section~\ref{deform algebra}, we take
\[
q_{\BZ}(\lambda)=4(\lambda,\lambda)_0 =2(\lambda,\lambda)_{\mathrm{int}} =4\sum_{j=1}^3(\lambda_j^{(2)})^2-12(\lambda^{(1)})^2, \qquad q(\lambda)=\xi^{q_{\BZ}(\lambda)}.
\]
Then, the corresponding polarization is $b_{\BZ}=8(\ ,\ )_0=4(\ ,\ )_{\mathrm{int}}$, and hence $b(\lambda,\mu)=\xi^{b_{\BZ}(\lambda,\mu)}$. Note that this quadratic form is integral on the coweight lattice $\Lambda$, since $\lambda^{(1)}\in\tfrac12\BZ$, and the three $\lambda_j^{(2)}$ are either all integral or all half-integral. In particular,
\[
q(-\varepsilon_1+\varepsilon_2)=\xi^8,\qquad q(-\varepsilon_3)=\xi^4,\qquad q(-\delta)=\xi^{-12}.
\]

The Cartan-scheme matrix is
\begin{equation}\label{eq:Cartan-scheme-f4}
C^a_{\mathfrak f(4)} =
\begin{pmatrix}
2&-1&-1&0\\
-1&2&-1&0\\
-1&-1&2&-1\\
0&0&-1&2
\end{pmatrix}.
\end{equation}

From now on, we will take the super Borel as in Section~\ref{super Borel}, and omit the subscript $a$ in the notation.

\begin{rem}
Although the matrices \eqref{eq:Cartan-scheme-gl}, \eqref{eq:Cartan-scheme-osp}, and \eqref{eq:Cartan-scheme-f4} at the chosen objects are symmetric, a matrix $C^{\rho_i(a)}$ obtained after an odd reflection need not be symmetric.
\end{rem}

\section{Perversity and identification of factorization algebras}\label{sec:factorization-identification}
In this section, we prove that the factorization algebras on the geometric side are perverse, and we identify them with the corresponding factorization algebras on the representation side. First, in Section~\ref{sec:orbit-preliminaries}, we collect some basic results about orbits on the affine Grassmannian, namely, the classification of orbits, the closure relations, the dimension estimates and the ratio local systems. Then, we prove the perversity of factorization algebras in Section~\ref{sec:perversity}, and we compare the factorization algebras on the two sides in Section~\ref{sec:factorization-comparison}. Finally, in Section~\ref{sec:positive-characteristic}, we compare the characteristic zero and the positive characteristic settings.

\subsection{Orbits and dimension estimates}\label{sec:orbit-preliminaries}

\subsubsection{Orbit classifications}

We consider the homogeneous case as in Section~\ref{pol des}. In this section and Section~\ref{sec:irreducibles}, we denote by $\Lambda_{G_i}^{\on{pos}}$ (resp. $\Lambda_{G_i}^{\on{neg}}$) the monoid generated by the positive (resp. negative) simple coroots of $G_i$, for $i=1,2$. We denote by $\leq_G$ the usual order generated by the positive simple coroots of $G=G_1\times G_2$, and $\leq_{\super}$ is the order of Section~\ref{sec:configuration-modules}, generated by the positive simple super roots.

For a coweight $\lambda$ of $G$, we write it in terms of its $G_1$-part and its $G_2$-part, i.e.,
$\lambda=(\lambda^{(1)},\lambda^{(2)})$,
with $\lambda^{(1)}\in\Lambda_{G_1}$ and $\lambda^{(2)}\in\Lambda_{G_2}$. For $\lambda,\theta\in \Lambda$, we denote by $\sO_{\Gr}^{\lambda}$ the corresponding $H(\bF)_c^\omega$-orbit in $\Gr_G^\omega$, and by $\sS_{\Gr}^{\theta}$ the corresponding semi-infinite orbit. For $i=1,2$ and $\nu\in\Lambda_{G_i}$, we let $U_i=U_i^+$ and denote
\[
\sS_{G_i}^{\pm,\nu}:=U_i^\pm(\bF)t^\nu G_i(\bO)/G_i(\bO), \qquad \sS_{\Gr}^{\theta} =\sS_{G_1}^{-,\theta^{(1)}}\times\sS_{G_2}^{+,\theta^{(2)}}.
\]
Note that, by our choice of the Borel subgroup $B=B_{G_1}^-\times B_{G_2}$ in Section~\ref{trans D}, the semi-infinite orbits in $\Gr_{G_1}$ are taken for $U_1^-$, and those in $\Gr_{G_2}$ for $U_2$. Thus, we also simply write $\sS_{G_1}^{\nu}=\sS_{G_1}^{-,\nu}$ and $\sS_{G_2}^{\nu}=\sS_{G_2}^{+,\nu}$. The superscripts $\pm$ here are to remind the readers of the unipotent subgroups used in the definition of the semi-infinite orbits.

Note that by taking the fiber over $t^0\in\Gr_{G_1}$ and the left translation by $\sD$, the $H(\bF)$-orbits in $\Gr_G$ correspond to the $G_1(\bO)H_u(\bF)$-orbits in $\Gr_{G_2}$.

\begin{lem}\label{lem:L-orbit-indices}
The $G_1(\bO)H_u(\bF)$-orbits in $\Gr_{G_2}$ are indexed by pairs $\lambda=(\lambda^{(1)},\lambda^{(2)})\in\Lambda$ satisfying the following conditions.
\begin{itemize}
\item In the $\GL(M|N)$ case, we require
\[
\lambda^{(1)}_1\geq\lambda^{(1)}_2\geq\cdots\geq\lambda^{(1)}_M, \qquad \lambda^{(2)}_1\leq\lambda^{(2)}_2\leq\cdots\leq\lambda^{(2)}_{M+1},
\]
while there is no requirement on $\lambda^{(2)}_{M+2},\ldots,\lambda^{(2)}_N$.
\item In the $\SOSp(2k|2l)$ case, with $M<N$, we require:
\begin{enumerate}
\item if $k\leq l$, i.e., $M$ is even, then
\[
\lambda^{(1)}_1\geq\cdots\geq\lambda^{(1)}_{k-1}\geq|\lambda^{(1)}_k|, \qquad \lambda^{(2)}_{l-k+1}\leq\cdots\leq\lambda^{(2)}_l\leq0;
\]
\item if $k>l$, i.e., $M$ is odd, then
\[
\lambda^{(1)}_1\geq\cdots\geq\lambda^{(1)}_{l-1}\geq\lambda^{(1)}_l\geq0, \qquad \lambda^{(2)}_{k-l}\leq\cdots\leq\lambda^{(2)}_{k-1}\leq-|\lambda^{(2)}_k|.
\]
\end{enumerate}
\item In the $\on{F}(4)$ case,\footnote{Recall that the coweight lattice of $\PGL_2$ is $\frac12\BZ$ in the diagonal coordinate used here, while the simple coroot generates the sublattice $\BZ$, and the coweight lattice of $\PSp(6)$ is given by $\lambda^{(2)}=(\lambda^{(2)}_1,\lambda^{(2)}_2,\lambda^{(2)}_3,-\lambda^{(2)}_3,-\lambda^{(2)}_2,-\lambda^{(2)}_1)\in(\frac12\BZ)^6$, such that $\lambda^{(2)}_1-\lambda^{(2)}_2\in\BZ$ and $\lambda^{(2)}_2-\lambda^{(2)}_3\in\BZ$.} we require
\[
\lambda^{(1)}\in\frac12\BN, \qquad \lambda^{(2)}_2-\lambda^{(2)}_1\in\BN, \qquad 2\lambda^{(2)}_3\in-\BN.
\]
\end{itemize}
In all cases, a representative is given by $t^{-\lambda^{(1)}}\sD t^{\lambda^{(2)}}\in\Gr_{G_2}^\omega$ (i.e., $\iota(t^{-\lambda^{(1)}})\sD t^{\lambda^{(2)}}$). We denote the index set by $\Lambda_G^{L,+}$, since they are the dominant coweights of $G_1\times L$.

Furthermore, if we require the $\psi$-relevance condition, then the orbits that satisfy it are indexed by $\Lambda_G^+$.
\end{lem}

\begin{proof}
For $\GL(M|N)$, see \cite[Lemma 3.3.3 and proof of Proposition 3.3.4]{[TY]}, and for $\on F(4)$, see \cite[Proposition 3.3.1(a,b)]{[FTY]}.

In the orthosymplectic case, by the Iwasawa decomposition for $P^-=LH_u$, it is enough to classify the $G_1(\bO)$-orbits on
\[
\Gr_L=\Gr_{\SO(M+1)}\times \Gr_{\BG_m^{(N-M-1)/2}}.
\]
Thus, the inequalities follow from \cite[Lemma 2.3.2]{[BFT]}, and the remaining coordinates are arbitrary integers.

Now, we check when $\psi$ vanishes on the stabilizer of $t^{-\lambda^{(1)}}\sD t^{\lambda^{(2)}}$. Recall that $\psi$ is given by the sum of residues in \cite[Section 3.2]{[BFT]}. On this stabilizer, the images of the linear terms whose residues define $\psi$ are $t^{\lambda^{(2)}_{r+1}-\lambda^{(2)}_r}\bO$, for $1\leq r\leq l-k$ if $k\leq l$, and $1\leq r\leq k-l-1$ if $k>l$. Note that the $t^{-1}$-coefficient vanishes on $t^a\bO$ if and only if $a\geq0$. Thus, $\psi$ vanishes on the stabilizer if and only if $\lambda^{(2)}_r\leq\lambda^{(2)}_{r+1}$ for these values of $r$, which, together with the inequalities for $L$, are exactly the dominance conditions for $G_2$.
\end{proof}

\subsubsection{Orbit closures} Let us recall some results about the closures of these orbits, cf. \cite{[BFT1],[TY],[FTY]}, etc.

\begin{lem}\label{lem:spherical-orbit-closure}
For $\lambda,\lambda'\in\Lambda_G^{L,+}$, if $\sO_{\Gr}^{\lambda}\subset\overline{\sO}_{\Gr}^{\lambda'}$, then $\lambda\leq_{\super}\lambda'$.
\end{lem}

\begin{proof}
For $\GL(M|N)$, see \cite[Proposition 8.2.1]{[TY]}, and for $\on{F}(4)$, see \cite[Proposition 3.4.1]{[FTY]}. The case of $\SOSp(2k|2l)$ follows from Lemma~\ref{lem:osp-minor-super-order} below, after translation by a sufficiently dominant coweight of $Z_L$.
\end{proof}

In particular, given two dominant coweights $\lambda,\mu$, where dominance is taken with respect to $B^-$, one has $\sO_{\Gr}^{\lambda}\subset \overline{\sO}_{\Gr}^{\mu}$ only if $\lambda\leq_{\super}\mu$.

\begin{lem}\label{lem:orbit-S-nonempty-super-order}
Let $\lambda,\mu\in \Lambda$, and assume that $\lambda$ is dominant with respect to $B^-$. Then, $\overline{\sO}_{\Gr}^{\lambda}\cap \overline{\sS}_{\Gr}^{\mu}\neq \varnothing$ only if $\mu\leq_{\super}\lambda.$
\end{lem}
\begin{proof}
In type A and in the $\on F(4)$ case, it follows from the valuation inequalities in \cite[Section 8.2]{[TY]} and \cite[Section 3.4]{[FTY]}, and in the orthosymplectic case, it is Lemma~\ref{lem:osp-minor-super-order} below.
\end{proof}

For completeness, we also include the proof of the above lemmas in the orthosymplectic case.
\begin{lem}\label{lem:osp-minor-super-order}
In the orthosymplectic cases, let $\lambda,\mu\in\Lambda$, and assume that $\lambda$ is dominant with respect to $B^-$. If $\sO_{\Gr}^{\mu}\subset\overline\sO_{\Gr}^{\lambda}$, then $\mu\leq_{\super}\lambda$. Furthermore, if $\overline\sO_{\Gr}^{\lambda}\cap \overline\sS_{\Gr}^{\mu}\neq\varnothing$, then $\mu\leq_{\super}\lambda$ as well.
\end{lem}
\begin{proof}
For the closure relation, we use the minimal degrees of the determinants of minors, as in \cite[proof of Proposition 8.2.1]{[TY]}, see also \cite[proof of Theorem 3.3.5(a,b)]{[BFT]}. Let us assume that $M=2k$, i.e., $k\leq l$, and let $s=l-k$. As in Lemma~\ref{lem:L-orbit-indices}, we consider the $G_1(\bO)H_u(\bF)$-orbits in $\Gr_{G_2}$, with the representatives $t^{-\lambda^{(1)}}\sD t^{\lambda^{(2)}}$.

First, we consider the minors in the following rows:
\[
\begin{gathered}
1,\ldots,j\quad(1\leq j\leq s),\qquad
1,\ldots,s+i\quad(1\leq i<k),\\
1,\ldots,s+i,l+1\quad(0\leq i<k).
\end{gathered}
\]
For each choice of rows, we take the minimal degree of the determinants over all choices of columns and all representatives of the orbit. Then, the corresponding minimal degrees are
\[
\sum_{h=1}^{j}\lambda_h^{(2)},\qquad
\sum_{h=1}^{s+i}\lambda_h^{(2)}-\sum_{h=1}^{i}\lambda_h^{(1)},\qquad
\sum_{h=1}^{s+i+1}\lambda_h^{(2)}-\sum_{h=1}^{i}\lambda_h^{(1)},
\]
respectively. Indeed, the first $s$ rows remove the contribution of $H_u(\bF)$, and these values follow from the matrix $\sD$ and the dominance conditions in Lemma~\ref{lem:L-orbit-indices}.

Secondly, we consider the $l\times l$-minors in the rows
\[
1,\ldots,l
\qquad\text{and}\qquad
1,\ldots,l-1,l+2,
\]
and their minimal degrees are, respectively,
\[
\sum_{j=1}^{l}\lambda_j^{(2)}-\sum_{i=1}^{k-1}\lambda_i^{(1)}-\lambda_k^{(1)},
\qquad
\sum_{j=1}^{l}\lambda_j^{(2)}-\sum_{i=1}^{k-1}\lambda_i^{(1)}+\lambda_k^{(1)}.
\]
Here, we should consider the two choices of rows separately, since $\SO(2k)$ preserves each of the two components of the variety of maximal isotropic subspaces in its standard representation.

Note that the condition that all these minors belong to $t^n\bO$ is closed. Thus, if $\sO_{\Gr}^{\mu}\subset\overline{\sO}_{\Gr}^{\lambda}$, then each minimal degree for $\mu$ is at least the corresponding one for $\lambda$. Together with the condition
\[
\sum_i(\lambda_i^{(1)}-\mu_i^{(1)})+\sum_j(\lambda_j^{(2)}-\mu_j^{(2)})\in2\BZ
\]
on the connected components, these inequalities imply $\mu\leq_{\super}\lambda$, by the simple super roots in Section~\ref{super Borel}. The case $M=2l+1$, i.e., $k>l$, is similar.

Finally, by the same convolution argument as in \cite[proof of Lemma 3.6.2(b)]{[BFT1]}, applied to the two factors $G_1$ and $G_2$, the closure relation proved above implies
\[
\overline\sO_{\Gr}^{\lambda}\cap\overline\sS_{\Gr}^{\mu}\neq\varnothing
\quad\Longrightarrow\quad
\mu\leq_{\super}\lambda.
\]
\end{proof}

\subsubsection{Dimension estimates}
Recall that, by \eqref{eq:local-degree-normalization}, we have
\begin{equation}\label{eq:super-geometric-dimension-grading}
d\Bigl(\sum_i n_i\alpha_i\Bigr) =\sum_i(2-p(\alpha_i))n_i.
\end{equation}
Here, $p(\alpha_i)\in\{0,1\}$ denotes the parity of the simple super root $\alpha_i$, and $d$ computes the dimensions of the SW Zastava spaces, see Remark~\ref{rem:global-stratum-degree}.

\begin{cor}\label{lem:dimension-estimate-super}
Let $\lambda$ be a dominant coweight (with respect to $B^-$), and assume $\lambda-\mu=\sum_i n_i\alpha_i, n_i\in \BN,$ where $\alpha_i$ are the simple super roots. Then,
\begin{equation}
\begin{split}
&\dim\bigl(\overline{\sO}_{\Gr}^{\lambda}\cap \overline{\sS}_{\Gr}^{\mu}\bigr)\\
&\quad\leq \max_{\substack{m_{1,i},m_{2,j}\in\BN\\
\lambda-\sum_i m_{1,i}\alpha_{1,i} -\sum_j m_{2,j}\alpha_{2,j}\geq_{\super}\mu}} \left(\sum_i m_{1,i}+\sum_j m_{2,j}\right).
\end{split}
\end{equation}
Here, $\alpha_{1,i}$ and $\alpha_{2,j}$ denote the simple coroots of $G_1$ and $G_2$, respectively.

In particular,
\begin{equation}\label{dim}
\dim\bigl(\overline{\sO}_{\Gr}^{\lambda}\cap \overline{\sS}_{\Gr}^{\mu}\bigr) \leq \frac{1}{2}\sum_i (2-p(\alpha_i))\,n_i.
\end{equation}
Note that the right-hand side of \eqref{dim} is $\tfrac12d(\lambda-\mu)$. Moreover, the equality in \eqref{dim} can hold only if $\lambda-\mu$ belongs to the monoid generated by the positive simple coroots of $G_1\times G_2$. Indeed, if the equality holds, then $\lambda-\mu-\sum_i m_{1,i}\alpha_{1,i}-\sum_j m_{2,j}\alpha_{2,j}$ has value zero under $d$, thus it vanishes by \eqref{eq:super-geometric-dimension-grading}. The number $\sum_i(2-p(\alpha_i))\,n_i$ can also be identified with the dimension of the corresponding component of the SW Zastava space $\cY_{\Gr,=\lambda\cdot c}^{\mu}=\cM_{=\lambda\cdot c}\times'_{\Bun_G}\Bun_B^{\omega,\mu}$, see Remark~\ref{rem:global-stratum-degree}.
\end{cor}
\begin{proof}
The same argument as in \cite[Proposition 3.6.1]{[BFT1]} and \cite[\S 8.3]{[TY]} applies here.
\end{proof}

By Corollary~\ref{lem:dimension-estimate-super}, for $\lambda$ dominant with respect to $B^-$, the dimension of $\sO_{\Gr}^{\lambda}\cap\sS_{\Gr}^{\theta}$ is bounded by $\frac12d(\lambda-\theta)$. Here, if $\lambda-\theta =\sum_i m_i\alpha_i$ is written in terms of positive simple super roots, then we have
\[
\frac12d(\lambda-\theta) =\frac12\sum_i\bigl(2-p(\alpha_i)\bigr)m_i.
\]
In particular, if $\lambda-\theta$ is a sum of simple coroots of $G$, then $\frac12d(\lambda-\theta)$ is the number of positive simple coroots in it, counted with multiplicities. We call this upper bound $\frac12d(\lambda-\theta)$ the \emph{critical dimension}, and we call an irreducible component of $\sO_{\Gr}^{\lambda}\cap\sS_{\Gr}^{\theta}$ \emph{critical-dimensional} if its dimension is equal to the critical dimension.

\subsubsection{The parabolic decomposition}

Let $G'$ be a reductive group, and let $P^-\subset G'$ be a parabolic subgroup opposite to a parabolic containing the fixed Borel. Let $U_P^-$ be its unipotent radical, and let $L$ be the Levi subgroup containing the fixed maximal torus. We denote by $\sS_{\Gr_{G'}}^\nu$ and $\sS_{\Gr_{G'}}^{-,\nu}$ the orbits of $t^\nu$ under the loop groups of the unipotent radicals of the fixed Borel and its opposite, respectively, and use the same notation for $L$. Consider the correspondence
\begin{equation}
\xymatrix{ & \Gr_{P^-}\ar[ld]_{\pr_L}\ar[rd]^{\pr_{G'}} &\\
\Gr_L && \Gr_{G'} . }
\end{equation}
Here, $\pr_L$ is induced by the projection $P^-\to L$, and $\pr_{G'}$ is induced by the embedding $P^-\hookrightarrow G'$.

The following lemma is a slight modification of \cite[\S 7.3.1, (7.5), (7.6)]{[G7]}.

\begin{lem}\label{lem:parabolic-product-decomposition}
Let $\lambda,\mu\in \Lambda_{G'}$, where we regard $\mu$ also as a coweight of $L$ via the common maximal torus. Let $\sS_{\Gr_L}^{\mu}\subset \Gr_L, \sS_{\Gr_{G'}}^{\lambda}\subset \Gr_{G'}$ be the corresponding semi-infinite orbits. Then, after choosing coordinates on $\sS_{\Gr_L}^{\mu}$, there is an isomorphism over $\sS_{\Gr_L}^{\mu}$
\begin{equation}
\pr_L^{-1}(\sS_{\Gr_L}^{\mu})\times_{\Gr_{G'}}\sS_{\Gr_{G'}}^{\lambda} \simeq \sS_{\Gr_L}^{\mu}\times \left( \bigl(U_P^-(\bF)\cdot t^\mu G'(\bO)/G'(\bO)\bigr) \cap \sS_{\Gr_{G'}}^{\lambda} \right),
\end{equation}
where the fiber product on the left is taken with respect to $\pr_{G'}:\Gr_{P^-}\to \Gr_{G'}$.

Moreover, under this isomorphism, for any subspace $X\subset \sS_{\Gr_L}^{\mu}$, one has
\begin{equation}
\pr_L^{-1}(X)\times_{\Gr_{G'}}\sS_{\Gr_{G'}}^{\lambda} \simeq X\times \left( \bigl(U_P^-(\bF)\cdot t^\mu G'(\bO)/G'(\bO)\bigr) \cap \sS_{\Gr_{G'}}^{\lambda} \right).
\end{equation}
\end{lem}
In our application below, $G'=G_2$ and $U_P^-=H_u$.

\subsubsection{The sharper dimension estimate} We need a sharper (than Corollary~\ref{lem:dimension-estimate-super}) dimension estimate for the intersection $\sO_{\Gr}^{0}\cap\sS_{\Gr}^{\alpha}$ in the homogeneous cases of Section~\ref{pol des}. Note that in the calculation of the $!$-stalks of $\Omega_q^{\DK}$ below, we may replace $\overline{\sO}_{\Gr}^{0}\cap\sS_{\Gr}^{\alpha}$ by $\sO_{\Gr}^{0}\cap\sS_{\Gr}^{\alpha}$, since the unit object has zero restriction to $\overline{\sO}_{\Gr}^{0}\setminus\sO_{\Gr}^{0}$ by Lemma~\ref{minimal-relevance}. Below is a uniform proof of the critical-dimension arguments that appear in \cite[Proposition 9.4.1 and Remark 9.4.2]{[TY]} and \cite[Proposition 3.9.4]{[FTY]}, and the idea of the proof already appeared in \cite[Remark 9.4.2]{[TY]}.

\begin{prop}\label{prop:equality-only-alpha1-zero}
In the homogeneous cases, let $\alpha=(\alpha^{(1)},\alpha^{(2)})\in\Lambda_{G_1}\oplus\Lambda_{G_2}$. If
\[
\dim(\sO_{\Gr}^{0}\cap \sS_{\Gr}^{\alpha}) =\tfrac12d(-\alpha),
\]
then $\alpha^{(1)}=0$.
\end{prop}

\begin{proof}
By Corollary~\ref{lem:dimension-estimate-super}, the equality can hold only if $\alpha$ belongs to the negative geometric coroot monoid with respect to $B^-=B_1^+\times B_2^-$, and then $d(-\alpha)=\deg(-\alpha)$ by \eqref{eq:local-degree-normalization}. We denote by $\check\rho_{G_1}$, $\check\rho_{G_2}$, and $\check\rho_L$ the half-sums of the roots in $B_1^+$, $B_2^-$, and $B_2^-\cap L$, respectively. We should study when the dimension of the intersection $\sO_{\Gr}^{0}\cap\sS_{\Gr}^{\alpha}$ equals
\[
\langle\check\rho_{G_1},-\alpha^{(1)}\rangle +\langle\check\rho_{G_2},-\alpha^{(2)}\rangle =\tfrac12d(-\alpha).
\]

Via the second projection, the intersection ${\sO}_{\Gr}^{0}\cap \sS_{\Gr}^{\alpha}$ in $\Gr_{G_1\times G_2}$ is identified with the intersection
\begin{equation}\label{10.7}
\left( \sD^{-1}U_{G_1}^-(\bF)t^{\alpha^{(1)}}H_u(\bF)G_2(\bO)/G_2(\bO) \right) \cap \sS_{\Gr_{G_2}}^{\alpha^{(2)}}.
\end{equation}

To see when the equality holds, we decompose this intersection according to the $L$-coweight $\mu$. Namely, we consider
\begin{equation}
X_\mu:= \left( \sD^{-1}U_{G_1}^-(\bF)t^{\alpha^{(1)}}L(\bO)/L(\bO) \right) \cap \sS_{\Gr_L}^{\mu},
\end{equation}
where we regard $\alpha^{(1)}$ as a coweight of $L$ via $\iota$. Since the chosen $\sD$ belongs to the lower unipotent subgroup of $L$, $X_\mu$ is contained in $\sS_{\Gr_L}^{-,\alpha^{(1)}}\cap\sS_{\Gr_L}^{\mu}$. If $X_\mu\neq\varnothing$, then $\mu-\alpha^{(1)}$ is a sum of positive coroots for $B_L=L\cap B$, thus $\mu$ also belongs to the coroot lattice of $L$ (also of $G_2$, since $L$ is standard). These pieces $\pr_L^{-1}(X_\mu)\times_{\Gr_{G_2}}\sS_{\Gr_{G_2}}^{\alpha^{(2)}}$ form a locally finite stratification of \eqref{10.7}.

By Lemma \ref{lem:parabolic-product-decomposition} with $G'=G_2$ and $U_P^-=H_u$, we have
\begin{equation}
\dim\left( \pr_L^{-1}(X_\mu)\times_{\Gr_{G_2}}\sS_{\Gr_{G_2}}^{\alpha^{(2)}} \right) = \dim X_\mu + \dim\left( H_u(\bF)\cdot t^\mu G_2(\bO)/G_2(\bO) \cap \sS_{\Gr_{G_2}}^{\alpha^{(2)}} \right).
\end{equation}

Set-theoretically, $\pr_L^{-1}\bigl(\sS_{\Gr_L}^{\mu}\cap \sS_{\Gr_L}^{-,\alpha^{(1)}}\bigr) \times_{\Gr_{G_2}}\sS_{\Gr_{G_2}}^{\alpha^{(2)}}$ is contained in $\sS_{\Gr_{G_2}}^{-,\alpha^{(1)}}\cap \sS_{\Gr_{G_2}}^{\alpha^{(2)}}.$ Hence,
\begin{equation}
\begin{split}
\dim\left( H_u(\bF)\cdot t^\mu G_2(\bO)/G_2(\bO) \cap \sS_{\Gr_{G_2}}^{\alpha^{(2)}} \right)&\leq \langle \check\rho_{G_2},-\alpha^{(2)}+\alpha^{(1)}\rangle - \langle \check\rho_L,-\mu+\alpha^{(1)}\rangle\\
&= \langle \check\rho_{G_2},-\alpha^{(2)}\rangle +\langle\check\rho_L, \mu\rangle.
\end{split}
\end{equation}

Here, the last equality follows from the fact that $\check\rho_{G_2}-\check\rho_L$ vanishes on the coroot lattice of $L$.

On the other hand, note that $(G_1\times L)/\Delta(G_1)$ is a Borel-free spherical variety with spherical roots of type $T$, where $\Delta(G_1)$ is the image of $h\mapsto(h,\sD^{-1}\iota(h)\sD)$. Namely, the general linear and orthosymplectic cases are as in \cite[Remark 9.4.2]{[TY]}; in the $\on F(4)$ case, the flag variety of $G_1\times L$ is $(\BP^1)^3$ with the diagonal action of $\PGL_2$, and the stabilizer of three distinct points is trivial, while the stabilizer of two distinct points is a torus.

Note that $X_\mu$ is the central fiber of degree $(\alpha^{(1)},\mu)$ of the open SW Zastava space for this spherical variety. Since $G_1$ is reductive, the spherical variety has an invariant volume form, and the cocharacter $\rho$ is central in $G_1\times L$. Thus, according to the dimension estimate in \cite[Remark 5.4.4, Proposition 6.1.1, and Lemma 6.2.1]{[SW]}, one has
\[
\dim X_\mu \leq -\langle\check\rho_{G_1},\alpha^{(1)}\rangle -\langle\check\rho_L,\mu\rangle,
\]
and it is strict unless $\alpha^{(1)}=\mu=0$ (in which case $X_0=\{1\}$). Combining the last two inequalities, we see that every stratum of \eqref{10.7} (indexed by $\mu$) has dimension at most $\tfrac12d(-\alpha)$, and the inequality is strict unless $\alpha^{(1)}=\mu=0$. Taking the maximum over all the strata, we obtain the desired statement.
\end{proof}

\subsubsection{The ratio local system}\label{sec:ratio-local-system}

Let $\lambda$ be dominant with respect to $B^-$, and assume that the orbit $\sO_{\Gr}^{\lambda}$ is relevant. Then, $\cG^G$ has an $H(\bF)_c^\omega$-equivariant trivialization on $\sO_{\Gr}^{\lambda}$, and we denote it by $\tau_H^\lambda$. Furthermore, since $U$ is unipotent, the restriction of $\cG^G$ to any $U(\bF)_c^\omega$-orbit $\sS_{\Gr}^{\theta}$ has a $U(\bF)_c^\omega$-equivariant trivialization. Namely, for $i=1,2$, let $\tau_{i,\pm}^\nu$ be the unique (up to a scalar) $U_i^\pm(\bF)$-equivariant trivialization on $\sS_{G_i}^{\pm,\nu}$, as in \cite[\S 7.4.3]{[G7]}, and then the trivialization on $\sS_{\Gr}^{\theta}$ is $\tau_{1,-}^{\theta^{(1)}}\boxtimes\tau_{2,+}^{\theta^{(2)}}$. We denote by
\begin{equation}\label{eq:ratio-local-system}
\cL_{\lambda,\theta}
\end{equation}
the local system on $\sO_{\Gr}^{\lambda}\cap \sS_{\Gr}^{\theta}$ corresponding to the ratio of these two trivializations. We still denote by $\psi$ the corresponding Artin--Schreier local system (or the exponential $D$-module in the $D$-module setting).

\begin{lem}\label{lem:kummer-as-separation}
Let $Z$ be an irreducible component of $\sO_{\Gr}^{\lambda}\cap \sS_{\Gr}^{\theta}$. If $\cL_{\lambda,\theta}\otimes\psi$ is constant on a dense open subset of $Z$, then both $\cL_{\lambda,\theta}$ and $\psi$ are constant on that dense open subset.
\end{lem}

\begin{proof}
In characteristic $p>0$, the first factor is a Kummer local system, and the Kummer character has order prime to $p$, while the Artin--Schreier character has order dividing $p$. Thus, if their tensor product is trivial, then both of them are trivial. In the $D$-module setting, the Kummer D-module is regular singular, while the exponential D-module is irregular unless it is trivial, and hence both factors are trivial.
\end{proof}

\subsection{Perversity of the geometric-side factorization algebras}\label{sec:perversity} In this subsection, we prove that the factorization algebras in Definition \ref{def fact alg} are all perverse.

\subsubsection{The De Concini--Kac factorization algebra} We now prove the following perversity statement.

\begin{prop}\label{prop 7.0.1}\label{prop:Omega-geo-perverse}
For any quantum parameter $q$, the geometric factorization algebras $\Omega_{q}^{\DK}, \Omega_{q}^{\Lus}$ are perverse.
\end{prop}

\begin{proof}
We prove the claim for $\Omega_{q}^{\DK}$. The claim for $\Omega_{q}^{\Lus}$ follows by Verdier duality.

Recall that
\begin{equation}
F_{\Gr}^{\DK,\glob}: \SD_{\cG^G}^{H,\Hecke,\psi}(\cM)\longrightarrow \SD_{\cG^\Lambda}(\Conf)
\end{equation}
is given by
\begin{equation}
\cF\longmapsto \Pi^d v_{\Gr,*}\Bigl( p_{\Gr}^!(\cF)\overset{!}{\otimes} q_{\Gr}^!j_{*,\Gr}(\IC_{\Bun_B}) \Bigr)[\dim \Bun_G],
\end{equation}
where
\begin{equation}
p_{\Gr}:\overline{\cY}_{\Gr}\to \cM,\qquad q_{\Gr}:\overline{\cY}_{\Gr}\to \overline{\Bun}_B,\qquad v_{\Gr}:\overline{\cY}_{\Gr}\to \Conf
\end{equation}
are the natural maps. Thus, $\Omega_{q}^{\DK} = F_{\Gr}^{\DK,\glob}(\IC_{\Gr,\glob}^{0}).$

Since $\IC_{\Gr,\glob}^{0}$ is the clean extension from $\Bun_H^\omega\subset \cM$, by base change, the object $\Omega_{q}^{\DK}$ is obtained by applying $v_{\Gr,*}$ to the restriction of $\Pi^d \bigl(p_{\Gr}^!(\IC_{\Gr,\glob}^{0})\overset{!}{\otimes} q_{\Gr}^!j_{*,\Gr}(\IC_{\Bun_B})\bigr)[\dim \Bun_G]$ to the open locus $\Bun_H^\omega\times_{\Bun_G}'\Bun_B\subset \overline{\cY}_{\Gr}.$ Note that the restriction is a perverse local system given by the tensor product of $\psi$ and the ratio local system $\cL$, see \eqref{eq:ratio-local-system}. Note that the map $v_{\Gr}:\Bun_H^\omega\times_{\Bun_G}'\Bun_B\to \Conf$ is affine. Hence, $v_{\Gr,*}$ is right $t$-exact for the perverse $t$-structure, and we obtain $\Omega_{q}^{\DK}\in {}^p\SD^{\leq 0}(\Conf)$.

It remains to prove that $\Omega_{q}^{\DK}\in {}^p\SD^{\geq 0}(\Conf)$.

By the factorization property, it is enough to check the deepest diagonal strata (which are given by $\alpha\cdot C\subset \Conf^{\alpha}$). Thus, for every non-zero $\alpha=(\alpha^{(1)},\alpha^{(2)})\in \Lambda^{\on{neg}}$, and for every point $x\in C$, we need to show that the $!$-stalk $i_{\alpha\cdot x}^!(\Omega_{q}^{\DK})$ lands in degree no less than $1$.

By base change, this $!$-stalk is identified with the cohomology of the tensor product of the $\psi$-local system and the ratio local system on the central fiber
\begin{equation}\label{open int}
{\sO}_{\Gr}^{0}\cap {\sS}_{\Gr}^{\alpha}.
\end{equation}
The $!$-restriction of the tensor product above to this central fiber has $1$-dimensional $!$-stalks in degree $d(-\alpha)$, see Remark~\ref{rem:global-stratum-degree}.

By Corollary~\ref{lem:dimension-estimate-super}, every irreducible component of \eqref{open int} has dimension at most $\tfrac12d(-\alpha)$.

Therefore, only components of the critical dimension can contribute to perverse degree $0$. It suffices to prove that, on every such critical-dimensional component, the tensor product of the additive local system and $\cL$ is non-trivial on a dense smooth open subset (note that the Kummer local system $\cL$ cannot cancel with $\psi$).

By Corollary \ref{lem:dimension-estimate-super} and Proposition \ref{prop:equality-only-alpha1-zero}, the intersection \eqref{open int} can be of critical dimension only if $\alpha^{(1)}=0$, and $\alpha^{(2)}$ is an element in the negative coroot monoid for the $G_2$-factor.

Now, recall that $G=G_1\times G_2$, and $H_u$ is the unipotent radical of the parabolic subgroup $P^-\subset G_2$. Let $\mathcal O_{H_u}:=H_u(\bF)\cdot 1\subset \Gr_{G_2}$ be the unital $H_u(\bF)$-orbit, and consider the stratification indexed by $\mu$ in the proof of Proposition~\ref{prop:equality-only-alpha1-zero}. When $\alpha=(0,\alpha^{(2)})$, every piece with $\mu\neq0$ has dimension strictly smaller than the critical dimension, thus only the piece with $\mu=0$ can meet a critical-dimensional component densely, and this piece is
\begin{equation}\label{spec open int}
\mathcal O_{H_u}\cap \sS_{\Gr_{G_2}}^{\alpha^{(2)}}.
\end{equation}
First, we show that the character map $\psi$ is non-constant on every critical-dimensional irreducible component of $\mathcal O_{H_u}\cap\sS_{\Gr_{G_2}}^{\alpha^{(2)}}$.

Note that the intersection \eqref{spec open int} is contained in the usual semi-infinite intersection
\begin{equation}\label{usual sem int}
\sS_{\Gr_{G_2}}^{-,0}\cap \sS_{\Gr_{G_2}}^{\alpha^{(2)}}.
\end{equation}
Moreover, according to the dimension assumption, every critical-dimensional irreducible component of \eqref{spec open int} is dense in an irreducible component of \eqref{usual sem int}.

We denote by $\psi_{G_2}:\sS_{\Gr_{G_2}}^{-,0}\to\BG_a$ the usual Whittaker character map (see \cite[Section 7.1]{[FGV]}), and by $\psi_{\sD}$ the restriction of $\psi$ to \eqref{spec open int}. By the formulas in Section~\ref{trans D}, the coefficient of each negative simple geometric root coordinate outside $L$ in $\psi_{\sD}$ is non-zero. Since the torus $T_2$ acts on the root coordinates with different characters, the restriction of $\psi_{G_2}$ is a linear combination of the $T_2$-translates of $\psi_{\sD}$.

If $\psi_{\sD}$ is constant on a critical component, then so is $\psi_{G_2}$, since the connected group $B_L=L\cap B$ preserves each irreducible component of \eqref{spec open int} and $T_2\subset B_L$. Then, by density, $\psi_{G_2}$ is constant on the corresponding component of \eqref{usual sem int}, which contradicts Lemma~\ref{lem:Whittaker-nonconstant-on-MV-components} below.

Now, we choose a local coordinate $t$ and a trivialization of the $\omega$-twist. Let $V$ be a dense smooth open subset of a critical component, and choose $z=[u(t)]\in V$ with $u(t)\in H_u(\bF)$ and $\psi_{\sD}(z)\neq0$. Note that the loop rotation $a\mapsto[u(at)]$ preserves the component, and
\[
\psi_{\sD}([u(at)])=a^{-1}\psi_{\sD}(z),\qquad a\in\BG_m.
\]
Since this function has a simple pole at $a=0$, the additive local system is non-trivial on the preimage of $V$. Thus, by Lemma~\ref{lem:kummer-as-separation}, its tensor product with $\cL$ is also non-trivial on $V$.
\end{proof}

The following lemma is due to \cite[Proposition 7.1.7]{[FGV]} and \cite[Proposition 6.4]{[BFGM]}.

\begin{lem}\label{lem:Whittaker-nonconstant-on-MV-components}
Let $G'$ be a reductive group. For every non-zero coweight $\alpha$ in the negative coroot monoid and every irreducible component of $\sS_{\Gr_{G'}}^{-,0}\cap \sS_{\Gr_{G'}}^{\alpha},$ the Whittaker character map $\psi_{G'}:\sS_{\Gr_{G'}}^{-,0}\rightarrow \BG_a$ is non-constant on this component.
\end{lem}

\subsubsection{The small factorization algebra}

Recall the factorization algebra $\Omega_q^{\sm}$ of Definition \ref{def fact alg}. The following proposition is proved by the same method as in \cite[Proposition 18.7.4]{[GL]}.

\begin{prop}\label{prop 7.0.5}
Let $\cF\in \left( \SD_{\cG^G}^{H(\bF)_c^\omega,\psi}(\Gr_{G,c}^\omega) \right)^\heartsuit$, and let $\IC_{\bar{\sS}_{\Gr}^{\lambda}}$ denote the $!$-restriction to the fiber over $\lambda\cdot c$ of $j^{\sinf}_{\Gr,!*}(\omega_{\sS_{\Gr,\Conf_{\infty\cdot c}}})[d]$ in the construction over $\Conf_{\infty\cdot c}$ of Remark~\ref{rem 3.5.6}. Then:
\begin{enumerate}
\item The cohomology $H\bigl(\Gr_G^\omega,\ \cF\overset{!}{\otimes}\IC_{\bar{\sS}_{\Gr}^{\lambda}}\bigr)$ satisfies:
\begin{enumerate}
\item $H^i=0$ for $i<0$;
\item $H^0\simeq H^0\bigl(\sS_{\Gr}^{\lambda},\ \cF|_{\sS_{\Gr}^{\lambda}}[d(\lambda)]\bigr);$
\item

$H^1$ injects into $H^1\bigl(\sS_{\Gr}^{\lambda},\ \cF|_{\sS_{\Gr}^{\lambda}}[d(\lambda)]\bigr)$.
\end{enumerate}
\item The assertion $(a)$ also holds for
\begin{equation}
H\bigl(\Gr_G^\omega,\ \cF\overset{!}{\otimes}j^{\sinf}_{\Gr,\lambda,!}(\omega_{\sS_{\Gr}^{\lambda}})[d(\lambda)] \bigr),
\end{equation}
and its $H^0$ injects into $H^0\bigl(\sS_{\Gr}^{\lambda},\ \cF|_{\sS_{\Gr}^{\lambda}}[d(\lambda)]\bigr)$. Here, $j^{\sinf}_{\Gr,\lambda,!}$ denotes the $!$-extension from $\sS_{\Gr}^{\lambda}$.
\end{enumerate}
\end{prop}

\begin{proof}
We sketch again for completeness. The support of $\IC_{\bar{\sS}_{\Gr}^{\lambda}}$ is contained in
\begin{equation}
\overline{\sS}_{\Gr}^{\lambda} = \bigcup_{\lambda'\in \lambda+\Lambda_G^{\pos}} \sS_{\Gr}^{\lambda'}.
\end{equation}
Here, $\Lambda_G^{\pos}$ is the monoid generated by the positive coroots of the geometric group $G$ with respect to $B^-$.

Using the Cousin decomposition with respect to the stratification by semi-infinite orbits, it suffices to prove the following two vanishing statements:
\begin{equation}
H^i\bigl( \sS_{\Gr}^{\lambda'},\ \cF\overset{!}{\otimes}\IC_{\bar{\sS}_{\Gr}^{\lambda}}|_{\sS_{\Gr}^{\lambda'}} \bigr)=0 \qquad \text{for } i\leq 1
\end{equation}
for every $\lambda'\in (\lambda+\Lambda_G^{\pos})\setminus\{\lambda\}$, and
\begin{equation}
H^i\bigl( \sS_{\Gr}^{\lambda},\ \cF\overset{!}{\otimes}\IC_{\bar{\sS}_{\Gr}^{\lambda}}|_{\sS_{\Gr}^{\lambda}} \bigr)=0 \qquad \text{for } i<0.
\end{equation}

By \cite[Proposition 13.2.10]{[GL]}, one has
\begin{equation}
\IC_{\bar{\sS}_{\Gr}^{\lambda}}\big|_{\sS_{\Gr}^{\lambda'}} \simeq \omega_{\sS_{\Gr}^{\lambda'}}[d(\lambda')] \otimes K_{\lambda',c},
\end{equation}
where $K_{\lambda',c}$ lives in cohomological degrees $>1$ if $\lambda'\neq \lambda$. Hence, the claim reduces to showing that
\begin{equation}
H^i\bigl( \sS_{\Gr}^{\lambda'},\ \cF|_{\sS_{\Gr}^{\lambda'}} \bigr)=0 \qquad \text{for } i<d(\lambda')
\end{equation}
for all $\lambda'\in \lambda+\Lambda_G^{\pos}$.

Applying the Cousin decomposition with respect to the $H(\bF)_c^\omega$-orbits in $\Gr_{G,c}^\omega$, it is enough to show
\begin{equation}
H^i\bigl( \sO_{\Gr}^{\mu}\cap \sS_{\Gr}^{\lambda'},\ \cF|_{\sO_{\Gr}^{\mu}\cap \sS_{\Gr}^{\lambda'}} \bigr)=0 \qquad \text{for } i<d(\lambda').
\end{equation}

Note that $\cF$ is in the heart of the local Gaiotto category, the restriction of $\cF$ to the orbit $\sO_{\Gr}^{\mu}$ is lisse, and its $!$-stalks live in degrees no smaller than $d(\mu)$. Therefore, the restriction of $\cF$ to $\sO_{\Gr}^{\mu}\cap \sS_{\Gr}^{\lambda'}$ lives in perverse degree no less than $d(\mu) - \dim\bigl(\sO_{\Gr}^{\mu}\cap \sS_{\Gr}^{\lambda'}\bigr)$. Thus, the cohomologies live in the degrees no less than $d(\mu) - 2\dim\bigl(\sO_{\Gr}^{\mu}\cap \sS_{\Gr}^{\lambda'}\bigr)$.

By the dimension estimate for the corresponding SW Zastava space, see \cite[Corollary 8.3.1]{[TY]} in type A or Corollary \ref{lem:dimension-estimate-super}, one has
\begin{equation}
2\dim\bigl(\sO_{\Gr}^{\mu}\cap \sS_{\Gr}^{\lambda'}\bigr) \leq d(\mu) -d(\lambda')=\dim {\cY}_{\Gr,=\mu\cdot c}^{\lambda'}.
\end{equation}
Equivalently,
\begin{equation}
d(\mu) - 2\dim\bigl(\sO_{\Gr}^{\mu}\cap \sS_{\Gr}^{\lambda'}\bigr) \geq d(\lambda').
\end{equation}
This proves the required vanishing.

The same proof works for $(2)$. The only change is that one uses \cite[Theorem 1.5.5 and Proposition 3.3.8]{[G5]}, with the same argument as in \cite[proof of Proposition 13.2.10]{[GL]}, which gives
\begin{equation}
j^{\sinf}_{\Gr,\lambda,!}(\omega_{\sS_{\Gr}^{\lambda}})[d(\lambda)] \big|_{\sS_{\Gr}^{\lambda'}} \simeq \omega_{\sS_{\Gr}^{\lambda'}}[d(\lambda')] \otimes \widetilde K_{\lambda',c},
\end{equation}
where $\widetilde K_{\lambda',c}$ lives in cohomological degrees $>0$ for $\lambda'\neq \lambda$. Thus, by the Cousin exact sequence, we obtain $(a)$ and the asserted injection.
\end{proof}

As a corollary, we have the following.

\begin{prop}\label{prop 7.0.6}
For any quantum parameter $q$, the factorization algebra $\Omega_q^{\sm}$ is perverse.
\end{prop}

\begin{proof}
According to Corollary \ref{5.0.8}, one has $\BD(\Omega_q^{\sm})\simeq \Omega_{q^{-1}}^{\sm}.$ Thus, it is enough to show that $\Omega_q^{\sm}$ is concentrated in non-negative perverse degrees.

By factorization, it suffices to show that for any $\alpha\in \Lambda^{\on{neg}}\setminus\{0\}$, the $!$-stalk $i_{\alpha\cdot x}^!(\Omega_q^{\sm})$ is concentrated in degrees no less than $1$.

By Proposition \ref{prop 7.0.5}, this reduces to showing that
\begin{equation}
H^0\bigl( \sS_{\Gr}^{\alpha},\ \IC_{\Gr}^0|_{\sS_{\Gr}^{\alpha}}[d(\alpha)] \bigr)=0,
\end{equation}
which has already been proved in the proof of Proposition \ref{prop:Omega-geo-perverse}.
\end{proof}

\subsubsection{The $t$-exactness of the functors} Let
\[
j_{\cY}:\cY_{\Gr,\infty\cdot c}\hookrightarrow \overline{\cY}_{\Gr,\infty\cdot c}
\]
be the open embedding. We keep the notations $p_{\Gr},q_{\Gr},v_{\Gr}$ for the projections and their restrictions to $\cY_{\Gr,\infty\cdot c}$, and on a fixed component of $B$-degree, we denote $d_B=\dim\Bun_B$ and $d_G=\dim\Bun_G$. We define
\[
p_{\Gr}^\bullet\cF:=\Pi^d p_{\Gr}^!\cF[d_G-d_B].
\]

By base change and the projection formula, for a locally compact global Gaiotto object $\cF$, we have
\begin{equation}\label{eq:DK-open-kernel-comparison}
\Pi^d\left(p_{\Gr}^!\cF\overset!\otimes q_{\Gr}^!j_{*,\Gr}\IC_{\Bun_B}\right)[d_G] \xrightarrow{\sim}j_{\cY,*}p_{\Gr}^\bullet\cF.
\end{equation}
Here, we use $\IC_{\Bun_B}=\omega_{\Bun_B}[-d_B]$. Note that for sufficiently anti-dominant degree, $p_{\Gr}$ is smooth of relative dimension $d_B-d_G$, thus $p_{\Gr}^\bullet$ is $t$-exact by the standard factorization argument as in \cite[Proposition 9.1.1 and Remark 9.1.6]{[TY]}.

For any $\cF$ in the global Gaiotto category, we obtain the same formula with both $*$-extensions replaced by $!$-extensions or $!*$-extensions. In these cases, we first use the ULA property of the global semi-infinite sheaves over $\Bun_G$ for bounded defect and sufficiently anti-dominant degree, see \cite[Theorem 4.2.1]{[Cam]} and \cite[Lemma 8.2.1]{[Yan]}, and then apply the standard factorization argument as in \cite[Section 8.2]{[Yan]}.

Furthermore, we have:

\begin{thm}\label{thm:Grassmannian-functors-texact}
The restrictions of the functors
\begin{equation}
F_{\Gr}^{\Lus}, F_{\Gr}^{\DK}, F_{\Gr}^{\sm}: \SD_{\cG^G}^{H(\bF)_c^\omega,\psi,\locc}(\Gr^\omega_{G,c}) \longrightarrow \SD_{\cG^\Lambda}(\Conf_{\infty\cdot c})
\end{equation}
are $t$-exact. Their continuous extensions to $\SD_{\cG^G}^{H(\bF)_c^\omega,\psi,\ren}(\Gr^\omega_{G,c})$ are also $t$-exact, and we denote these extensions by the same notations.
\end{thm}

\begin{proof}
We first consider $F_{\Gr}^{\sm}$. Since $F_{\Gr}^{\sm}\circ \BD\simeq \BD\circ F_{\Gr}^{\sm}$, we only need to prove that $F_{\Gr}^{\sm}(\cF)$ lives in non-negative perverse degrees for $\cF$ in the heart. By factorization, this reduces to the assertion that the $!$-stalk of $F_{\Gr}^{\sm}(\cF)$ at every point $\lambda\cdot c\in \Conf_{\infty\cdot c}$ is concentrated in non-negative degrees, which follows from Proposition~\ref{prop 7.0.5}(1).

Then, by Verdier duality, it is sufficient to show the case for $F_{\Gr}^{\DK}$. By the dimension estimate in Corollary~\ref{lem:dimension-estimate-super}, the same cohomology calculation as in the proof of Proposition~\ref{prop 7.0.5} shows that $F_{\Gr}^{\DK}$ is left $t$-exact. To prove the right $t$-exactness, we use the global model. By Proposition~\ref{prop:local-global-functor-comparison} and \eqref{eq:DK-open-kernel-comparison}, we have
\[
F_{\Gr}^{\DK,\glob}(\cF) \simeq\bigl(v_{\Gr}|_{\cY_{\Gr,\infty\cdot c}}\bigr)_*p_{\Gr}^\bullet\cF.
\]
Since $p_{\Gr}^\bullet$ is $t$-exact and $v_{\Gr}|_{\cY_{\Gr,\infty\cdot c}}$ is affine, the right-hand side is right $t$-exact.

Finally, since the locally compact subcategory is preserved by the perverse truncations (as in Section~\ref{sec:global-duality}), the continuous extensions are also $t$-exact for the $t$-structures compatible with filtered colimits.
\end{proof}

\subsection{Comparison of the factorization algebras on the geometric and representation sides}\label{sec:factorization-comparison}

For completeness and the convenience of readers, we briefly sketch the proof in \cite[\S 6.2--6.6]{[G7]} and explain how the same argument applies in the Gaiotto setting. We temporarily distinguish the geometrically defined factorization algebras from the quantum supergroup ones by writing $\Omega_{q,\geo}^{?}, \Omega_{q,\quant}^{?}, ?\in\{\DK,\Lus,\sm\}.$ After the comparison theorem, we will omit the subscripts $\geo$ and $\quant$.

\begin{thm}\label{thm:geo-quant-comparison}
Assume the base field of schemes is of characteristic $0$, and $q$ satisfies condition $(*)$. Then, for the homogeneous data of Section~\ref{pol des}, there are canonical isomorphisms
\begin{equation}
\Omega_{q,\geo}^{\Lus}\simeq \Omega_{q,\quant}^{\Lus},\qquad \Omega_{q,\geo}^{\DK}\simeq \Omega_{q,\quant}^{\DK},\qquad \Omega_{q,\geo}^{\sm}\simeq \Omega_{q,\quant}^{\sm}.
\end{equation}
\end{thm}

\subsubsection{The Frobenius action on the geometric side}
\begin{prop}\label{prop:geometric-Frobenius-transfer}
Assume that $q$ is torsion valued on $Q_s$ and satisfies condition $(*)$. With the identification of Lemma \ref{lem:even-Frobenius-identification}, there is an action
\[
\Omega^{\#,\mathrm{cl}}\curvearrowright\Omega_{q,\geo}^{\Lus}
\]
for which the canonical map to $\Omega_{q,\geo}^{\sm}$ induces
\[
\operatorname{Bar}\bigl(\Omega^{\#,\mathrm{cl}}, \Omega_{q,\geo}^{\Lus}\bigr)\simeq\Omega_{q,\geo}^{\sm}.
\]
\end{prop}

\begin{proof}
As in \cite[Theorem 6.7.2 and Corollary 6.7.3]{[G7]}, we have the action and the bar isomorphism for the global semi-infinite sheaves, where the classical factorization algebra is $\Omega^{\#,\mathrm{cl}}$ by Lemma~\ref{lem:even-Frobenius-identification}. Taking the $!$-pullback along $q_{\Gr}$, the $!$-tensor product with $p_{\Gr}^!\IC_{\Gr,\glob}^{0}$, and then applying $\Pi^d v_{\Gr,*}[\dim\Bun_G]$, by Corollary~\ref{cor:omega-from-global-vacuum}, we obtain the required action and bar isomorphism, as in \cite[Section 6.7.4]{[G7]}.
\end{proof}

\subsubsection{The two Frobenius structures} Recall the open disjoint locus $\jmath:\Conf^\circ\hookrightarrow \Conf$, and we have canonical identifications $\jmath^!(\Omega_{q,\geo}^{?})\simeq \jmath^!(\Omega_{q,\quant}^{?}), ?\in\{\DK,\Lus,\sm\}.$ The problem is to extend these identifications across the diagonals. To repeat the proof in \cite{[G7]}, we essentially only need to add the following two facts.
\begin{enumerate}
\item First, the geometric objects $\Omega_{q,\geo}^{\DK}, \Omega_{q,\geo}^{\Lus}, \Omega_{q,\geo}^{\sm}$ are perverse by Propositions \ref{prop:Omega-geo-perverse} and \ref{prop 7.0.6}, and the factorization algebras arising from quantum supergroups are perverse by construction;
\item Second, the quantum supergroup side has the Frobenius structure, see Corollary \ref{cor:frob-action-factorization-super}:
\begin{equation}
\Omega^{\#,\mathrm{cl}}\curvearrowright \Omega_{q,\quant}^{\Lus},\qquad \Omega_{q,\quant}^{\Lus}\to \Omega_{q,\quant}^{\sm},
\end{equation}
and
\begin{equation}
\operatorname{Cone}\bigl((\Omega^{\#,\mathrm{cl}})_C\to(\Omega_{q,\quant}^{\Lus})_C\bigr) \otimes_{(\Omega^{\#,\mathrm{cl}})_C}\omega_C \xrightarrow{\sim} (\Omega_{q,\quant}^{\sm})_C.
\end{equation}

Also, we have the analogous geometric Frobenius statement for $\Omega_{q,\geo}^{\Lus}$,
\begin{equation}
\Omega^{\#,\mathrm{cl}}\curvearrowright \Omega_{q,\geo}^{\Lus},\qquad \Omega_{q,\geo}^{\Lus}\to \Omega_{q,\geo}^{\sm},
\end{equation}
and
\begin{equation}
\operatorname{Cone}\bigl((\Omega^{\#,\mathrm{cl}})_C\to(\Omega_{q,\geo}^{\Lus})_C\bigr) \otimes_{(\Omega^{\#,\mathrm{cl}})_C}\omega_C \xrightarrow{\sim} (\Omega_{q,\geo}^{\sm})_C.
\end{equation}
The geometric statements follow from Proposition \ref{prop:geometric-Frobenius-transfer} and the same description on the main diagonals as in \eqref{eq:cone-small-factorization-super}.
\end{enumerate}

\subsubsection{The generic comparison}
The following lemma is an analogue of the vanishing results in \cite[Proposition 9.3.1]{[TY]} (type $\GL(M|N)$), \cite[Proposition 3.9.1]{[FTY]} (type $\on{F}(4)$), and \cite[Sections 6.3--6.5]{[G]} (reductive case).

\begin{lem}\label{lem:orthosymplectic-vacuum-subtop}
For the orthosymplectic data of Section~\ref{pol des}, let $q$ be non-torsion valued on the even super roots. If $0\neq\alpha\in\Lambda^{\on{neg}}$ is not a negative simple super root, then
\begin{equation}\label{eq:orthosymplectic-vacuum-subtop}
H^1\left( \sO_{\Gr}^{0}\cap\sS_{\Gr}^{\alpha}, \IC_{\Gr}^{0}\overset{!}{\otimes} \omega_{\sS_{\Gr}^{\alpha}}[d(\alpha)] \right)=0.
\end{equation}
The assertion also holds after inverting $\hbar$ in the formal deformation.
\end{lem}

\begin{proof}
Let $D=-d(\alpha)=\deg(-\alpha)$. By Corollary~\ref{lem:dimension-estimate-super}, only the components of dimension $D/2$ or $(D-1)/2$ can contribute to $H^1$. For the components of dimension $(D-1)/2$, $H^1$ is the lowest possible cohomology group, and by the above dimension estimate, the argument with the divisors of determinant sections in \cite[Lemma 3.7.2]{[BFT1]} and \cite[Lemmas 9.1.8--9.1.9]{[TY]} shows that the ratio local system is non-trivial. Thus, the vanishing follows as in the proof of \cite[Proposition 9.3.1]{[TY]}.

For the components of dimension $D/2$, the strict bound in the proof of Proposition~\ref{prop:equality-only-alpha1-zero} reduces the calculation to $\alpha^{(1)}=\mu=0$, as in \cite[Remark 9.4.2]{[TY]}, thus only $H_u(\bF)\cdot1\cap\sS_{\Gr_{G_2}}^{\alpha^{(2)}}$ contributes. Then, by the same Fourier argument as in \cite[\S\S 6.4--6.5]{[G]}, applied to the residue map of $H_u(\bF)$, we only need to consider the case $\alpha^{(2)}=-n\beta$, where $\beta$ is a simple geometric coroot outside $L$. For $n>1$, the vanishing follows from the rank-one calculation in \cite[\S 6.5, Case (1)]{[G]}. 
\end{proof}

\begin{lem}\label{lem:generic-Gaiotto-comparison}
For the homogeneous data of Section~\ref{pol des}, let the parameter be non-torsion valued on the even super roots, or be the formal parameter $q((\hbar))$. Then, the three factorization algebras on the geometric side coincide and are the intermediate extension of their (common) restriction to the disjoint locus.

In characteristic zero, the identification on the disjoint locus canonically extends to
\begin{equation}\label{eq:generic-Gaiotto-comparison}
\Omega_{q((\hbar)),\geo}^{\DK}\simeq \Omega_{q((\hbar)),\quant}^{\DK},
\end{equation}
and similarly for a parameter which is non-torsion valued on the even super roots. Also, the three factorization algebras on the representation side coincide.
\end{lem}
\begin{proof}
The three factorization algebras on the geometric side coincide by the cleanness of the semi-infinite kernels, see \cite[Propositions 5.4.4 and 5.6.5]{[G7]}.

To show that the factorization algebra on the geometric side is the intermediate extension of its restriction to the disjoint locus, by perversity, factorization and Verdier duality, it suffices to prove
\[
H^1\left( \sO_{\Gr}^{0}\cap\sS_{\Gr}^{\alpha}, \IC_{\Gr}^{0}\overset{!}{\otimes} \omega_{\sS_{\Gr}^{\alpha}}[d(\alpha)] \right)=0
\]
for every non-zero $\alpha\in\Lambda^{\on{neg}}$ which is not a negative simple super root. This is \cite[Proposition 9.3.1]{[TY]} in type A, \cite[Propositions 3.9.1 and 3.9.5]{[FTY]} for $\on F(4)$, and Lemma~\ref{lem:orthosymplectic-vacuum-subtop} in the orthosymplectic case. This proves the assertion on the geometric side. 

In characteristic zero, the three factorization algebras on the representation side coincide by the PBW descriptions of the corresponding Hopf superalgebras, and they are the intermediate extension of the same local system on the disjoint locus.
\end{proof}

\subsubsection{The formal deformation} Following \cite[\S 6.5]{[G7]}, we consider the formal deformation $q[[\hbar]]$. Similar to \cite[(6.13)]{[G7]}, by Lemma~\ref{lem:generic-Gaiotto-comparison}, the open disjoint locus isomorphism and the isomorphism after inverting $\hbar$ give a canonical map
\begin{equation}\label{can fact dk}
\Omega_{q[[\hbar]],\geo}^{\DK}\longrightarrow \Omega_{q[[\hbar]],\quant}^{\DK}.
\end{equation}
Here, the target is the factorization algebra defined in \eqref{def h-deform}, corresponding to $U_{q[[\hbar]]}^{\DK}(N_s^+)$.

The following proposition can be proved with the same argument as in \cite[Proposition 6.5.5]{[G7]}.
\begin{prop}\label{prop:DK-small-formal-equivalence-Gaiotto}
If the characteristic of the base field is $0$ and the generic comparison \eqref{eq:generic-Gaiotto-comparison} holds, the following statements are equivalent:
\begin{enumerate}
\item $\Omega_{q[[\hbar]],\geo}^{\DK}\simeq \Omega_{q[[\hbar]],\quant}^{\DK};$
\item $\Omega_{q,\geo}^{\DK}\simeq \Omega_{q,\quant}^{\DK};$
\item $\Omega_{q,\geo}^{\sm}\simeq \Omega_{q,\quant}^{\sm}.$
\end{enumerate}
\end{prop}

\subsubsection{}
\begin{proof}[Proof of Theorem \ref{thm:geo-quant-comparison}]
We follow the proof in \cite[Section 6.6.5]{[G7]} and sketch.

Namely, for a parameter which is non-torsion valued on the even super roots, the assertion follows from Lemma~\ref{lem:generic-Gaiotto-comparison}.

For torsion valued $q$, the morphism \eqref{can fact dk} becomes an isomorphism after inverting $\hbar$ by Lemma~\ref{lem:generic-Gaiotto-comparison}, and the Frobenius statements on the geometric side and on the representation side are given by Proposition~\ref{prop:geometric-Frobenius-transfer} and Corollary~\ref{cor:frob-action-factorization-super}, respectively. Therefore, by the same induction on the configuration degree as in loc.\ cit., we obtain
\begin{equation}\label{sm iso}
\Omega_{q,\geo}^{\sm}\simeq\Omega_{q,\quant}^{\sm}.
\end{equation}
Then, the isomorphism $\Omega_{q,\geo}^{\DK}\simeq\Omega_{q,\quant}^{\DK}$ follows from Proposition~\ref{prop:DK-small-formal-equivalence-Gaiotto}, and the corresponding isomorphism for $\Omega_q^{\Lus}$ follows from the Verdier duality.
\end{proof}

\subsection{Comparison of characteristic zero and positive characteristic}\label{sec:positive-characteristic} Now, we compare the factorization algebras and the $!$-stalks of their factorization IC modules in characteristic zero and in sufficiently large positive characteristic. To define these objects over a localization of a number ring, we use the Kirillov model in \cite[Appendix A]{[GL]} to avoid the problem of lacking Whittaker/Gaiotto model.

\subsubsection{The Kirillov model}

Recall that the usual Whittaker model uses an exponential $D$-module in characteristic zero and an Artin--Schreier sheaf in positive characteristic, and the Kirillov model replaces this character sheaf by the averaging for $\BG_m\ltimes\BG_a$. Let $\cC$ be a sheaf category with an action of $\BG_m\ltimes\BG_a$. Following \cite[Section A.2]{[GL]}, let
\begin{equation}\label{eq:integral-Kirillov-model}
\begin{split}
\Kir_*(\cC)&:= \ker\left( \Av_*^{\BG_a}:\cC^{\BG_m} \longrightarrow\cC^{\BG_m\ltimes\BG_a} \right),\\
\Kir_!(\cC)&:= \ker\left( \Av_!^{\BG_a}:\cC^{\BG_m} \longrightarrow\cC^{\BG_m\ltimes\BG_a} \right).
\end{split}
\end{equation}
For the $!$-version, we work with constructible sheaves or holonomic $D$-modules. When $\cC=\SD(Y)$ (with the equivariant twisting understood), we denote them by $\Kir_*(Y)$ and $\Kir_!(Y)$, and denote the functors
\[
\begin{split}
\Pi_*(\cF)&= \operatorname{coFib}\left( \oblv\Av_*^{\BG_a}(\cF)\longrightarrow\cF \right),\\
\Pi_!(\cF)&= \operatorname{Fib}\left( \cF\longrightarrow\oblv\Av_!^{\BG_a}(\cF) \right).
\end{split}
\]

For the Gaiotto category with $\psi\neq0$, we first take the $\ker\psi$-invariants, and then consider the corresponding Kirillov model, where $\BG_m$ acts by the left multiplication through the integral cocharacter $\rho=(0,\rho^{(2)})$ in Section~\ref{super Borel}. In the $D$-module setting in characteristic zero and the $\ell$-adic setting in positive characteristic, it recovers the Gaiotto category by \cite[Proposition A.3.2]{[GL]}, and the same definition also works for Betti sheaves and for $\ell$-adic sheaves in characteristic zero.

The convolution with the $!$- and $*$-extensions of the constant sheaf along $\BG_m\hookrightarrow\BG_a$ gives rise to an equivalence $\Kir_*(Y)\simeq\Kir_!(Y)$, see \cite[Sections A.4.3--A.4.5]{[GL]}, and the Verdier duality exchanges the two categories on locally compact objects. Note that for equivariant maps, $f^!$ and $f_*$ preserve $\Kir_*$, while $f^*$ and $f_!$ preserve $\Kir_!$, as in \cite[Remark A.1.4]{[GL]}.

These averaging and convolution correspondences are defined over any integral base over which the group actions and the twisting are defined. Thus, they give an integral model for the constructions in the preceding sections, without choosing an additive character sheaf over that base.

\subsubsection{Spreading out}
Assume that $q$ has finite order $m$, and let us fix an $\ell$-adic coefficient field $\se$ which contains the values of $q$. We choose a number field $K$ which contains the $m$-th roots of unity, and we let
\[
S=\operatorname{Spec}\mathcal O_K[1/N_0], \qquad m\ell\mid N_0.
\]
Note that the generic point of $S$ has residue field $K$, and its closed points have finite residue fields of characteristic prime to $N_0$. In order to compare the constructions on the corresponding fibers, we should define the formal deformations of the factorization gerbes $\cG^G$ and $\cG^\Lambda$ over $S$. Here, the coefficient field $\se$ is fixed, and $\hbar$ is a formal parameter in the coefficients.

Consider the deformation $q_\hbar$ in \eqref{eq:formal-quadratic-bilinear-super}, i.e., $q_\hbar=q(1+\hbar)^{q_{\BZ}}$ in the $\ell$-adic setting. We need compatible deformations over $S$, whose restrictions to the geometric fibers are
\[
q_{\hbar,n}(\lambda)=q(\lambda)(1+\hbar)^{q_{\BZ}(\lambda)}\pmod{\hbar^{n+1}}, \qquad n\geq0.
\]

First, it is easier to see that the finite-order part $q$ is defined over $S$, namely, it is defined by the Kummer local systems associated with the covering $y\mapsto y^m$ of $\BG_{m,S}$. 

Let us consider the unipotent part over $S$. Below, we follow \cite[Sections 5.1--5.2 and the proof of Lemma 6.2.3]{[Mor]}. 

Note that we also need its monodromy to be compatible with the Galois action. For example, over a closed point with residue field $\BF_Q$, let us denote by $T$ the action of a generator of the tame inertia at $0$ on a local system on $\BG_m$, and by $F$ the action of the arithmetic Frobenius. Here, the tame inertia is described by the Kummer coverings $y\mapsto y^r$, where $r$ is prime to the residue characteristic. Since the arithmetic Frobenius sends $\zeta$ to $\zeta^Q$, we have
\[
FTF^{-1}=T^Q.
\]
Thus, if $T$ is unipotent, then its logarithm satisfies $F(\log T)F^{-1}=Q\log T$, and this is the reason why we use Tate twists in the coefficients.

More precisely, recall that
\[
\BZ_\ell(1)=\varprojlim_r\mu_{\ell^r}, \qquad \se(1)=\BZ_\ell(1)\otimes_{\BZ_\ell}\se,
\]
where the transition maps are given by $\zeta\mapsto\zeta^\ell$. We denote by $\se(j)=\se(1)^{\otimes j}$ the $j$-th Tate twist. Note that the arithmetic Frobenius at a closed point with residue field $\BF_Q$ acts on $\se(j)$ by $Q^j$. Then, we consider the following sheaves of coefficient algebras on $S$,
\[
\mathcal R_{S,n}=\bigoplus_{j=0}^{n}\se(j), \qquad \mathcal R_S=\varprojlim_n\mathcal R_{S,n},
\]
where the multiplication on $\mathcal R_{S,n}$ is given by the usual multiplication $\se(i)\otimes\se(j)\to\se(i+j)$ if $i+j\leq n$, and it is zero otherwise. Note that these are sheaves of $\ell$-adic coefficient algebras on $S$, rather than sheaves of regular functions. We denote by $\mathfrak m$ the kernel of the augmentation $\mathcal R_S\to\se$. At a geometric point $\bar s\to S$, a choice of a basis $u$ of $\se(1)_{\bar s}$ gives the identifications
\[
(\mathcal R_{S,n})_{\bar s}\simeq\se[u]/(u^{n+1}), \qquad \mathcal R_{\bar s}:=\varprojlim_n(\mathcal R_{S,n})_{\bar s}\simeq\se[[u]].
\]
If we set $\hbar=\exp(u)-1$, then $\mathcal R_{\bar s}$ is identified with $\se[[\hbar]]$, and $\mathfrak m_{\bar s}$ is identified with the ideal $(\hbar)$.  Note that at a closed point with residue field $\BF_Q$, the action $F(u)=Qu$ gives $F(1+\hbar)=(1+\hbar)^Q$, as required by the monodromy relation above.

On $\BG_{m,S}$, we take a rank-one local system with coefficients in $\mathcal R_{S,n}$ and tame monodromy $\exp(q_{\BZ}(\lambda)u)$ on each geometric fiber. This is compatible with the Galois action described above. Tensoring with the finite-order Kummer local system gives the parameter $q_{\hbar,n}(\lambda)$. They provide compatible deformations over $S$. The same construction applies to $b'_\hbar$.

In the following, when passing from the perverse normalization on the total space to that on a geometric fiber, we remove the shift $[\dim S]$ contributed by the base. By the \textit{generic perversity}, we mean the following well-known fact (see \cite[proof of Lemma 6.2.3]{[Mor]}): if the restriction of a bounded constructible complex on a scheme of finite type over $S$ to the geometric generic fiber is perverse, then, after shrinking $S$, its restrictions to all geometric fibers are perverse.

\subsubsection{} In order to spread out the De Concini--Kac factorization algebra on the representation side, we need the following lemma.
\begin{lem}\label{lem:formal-saturation-spreading}
Over a fixed finite union of components of the configuration space, after shrinking $S$, the intersection in \eqref{def h-deform} commutes with the restriction to geometric fibers.
\end{lem}
\begin{proof}
Let us denote by $\Omega_{\mathcal R}^{\circ}$ the deformation of $\Omega_q^\circ$ over $S$, and let
\[
A=\jmath_*\Omega_{\mathcal R}^{\circ},\qquad \mathcal{Q}=\operatorname{coker}(\jmath_!\Omega_{\mathcal R}^{\circ}\to A).
\]
We denote by $\mathcal{Q}[\mathfrak m^n]$ the largest subobject of $\mathcal{Q}$ killed by $\mathfrak m^n$. By the same argument as in \cite[Lemma A.2.8]{[G7]}, these subobjects stabilize, and we denote by $\mathcal{Q}_{\mathrm{tors}}$ their stable value. Note that, after inverting the formal parameter, the intermediate extension is the image of $\jmath_!\Omega_{\mathcal R}^{\circ}\to A$. So, the intersection defining $\Omega_{\mathcal R,\quant}^{\DK}$ in \eqref{def h-deform} is given by
\[
\Omega_{\mathcal R,\quant}^{\DK}=\ker(A\to \mathcal{Q}/\mathcal{Q}_{\mathrm{tors}}).
\]

Now, by the generic base change and the generic perversity, as in \cite[proof of Lemma 6.2.3]{[Mor]}, after shrinking $S$, the exact sequences defining $\mathcal{Q}$, $\mathcal{Q}_{\mathrm{tors}}$ and the above kernel remain exact on every geometric fiber, and $(\mathcal{Q}/\mathcal{Q}_{\mathrm{tors}})_{\bar s}$ remains flat. Thus, $(\mathcal{Q}_{\mathrm{tors}})_{\bar s}$ is exactly the $\mathfrak m_{\bar s}$-torsion subobject of $\mathcal{Q}_{\bar s}$, and hence, the above kernel gives the intersection on each geometric fiber, which proves the lemma.
\end{proof}

\subsubsection{Comparison of factorization algebras}
\begin{prop}\label{prop:spread-out-DK-comparison}
Assume that $q$ has finite order and satisfies condition $(*)$. For any $\alpha\in\Lambda^{\on{neg}}\setminus\{0\}$, there exists an integer $N(\alpha)$ such that, for every geometric point $\bar{s}$ of characteristic $p>N(\alpha)$, the comparison isomorphisms of Theorem \ref{thm:geo-quant-comparison} hold on
\[
\Conf_{\bar{s}}^{\geq\alpha} :=\bigsqcup_{\alpha\leq_{\super}\beta<_{\super}0} \Conf_{\bar{s}}^\beta.
\]
\end{prop}
\begin{proof}
We work over $\Conf_S^{\geq\alpha}$. After shrinking $S$, we may assume that the groups, affine Grassmannians, configuration spaces, gerbes and Zastava correspondences are all defined over $S$.

First, we consider the geometric side. Using the Kirillov model, the constructions in the preceding sections give $\Omega_{\mathcal R,\geo}^{\DK}$ over $S$. By Lemma~\ref{lem:proper-support-global-functors}, we may restrict to finite type spaces. We apply the generic base change with coefficients in $\se=\mathcal R_S/\mathfrak m$. Then, after shrinking $S$, the base change map for $\Omega_{\mathcal R,\geo}^{\DK}$ is an isomorphism modulo $\mathfrak m$, and hence, it is an isomorphism with formal coefficients.

Secondly, on the representation side, we define $\Omega_{\mathcal R,\quant}^{\DK}$ over $S$ by \eqref{def h-deform}, and by Lemma~\ref{lem:formal-saturation-spreading}, its geometric fibers are identified with the same construction over a field. We denote $A=\jmath_*\Omega_{\mathcal R}^{\circ}$, $B=\Omega_{\mathcal R,\quant}^{\DK}\subset A$. By Theorem~\ref{thm:geo-quant-comparison} and Proposition~\ref{prop:DK-small-formal-equivalence-Gaiotto}, the geometric generic fibers of $\Omega_{\mathcal R,\geo}^{\DK}$ and $B$ are identified. Furthermore, by the generic perversity modulo $\mathfrak m$, $\Omega_{\mathcal R,\geo}^{\DK}$ is perverse and flat after shrinking $S$.

Now, by adjunction, we have a map $\Omega_{\mathcal R,\geo}^{\DK}\to A$, and its image in $A/B$ has zero geometric generic fiber. The derived reduction of this image modulo $\mathfrak m$ is constructible and has zero geometric generic fiber. By Chevalley's theorem, the image of its support in $S$ is constructible, and hence consists of finitely many closed points. After removing these points, the reduction of this image is zero, and hence the image itself is zero, i.e., the map factors through
\[
c_{\mathcal R}:\Omega_{\mathcal R,\geo}^{\DK} \longrightarrow\Omega_{\mathcal R,\quant}^{\DK}.
\]

The cone of $c_{\mathcal R}$ has zero geometric generic fiber. By the same argument, it vanishes after shrinking $S$, and hence, $c_{\mathcal R}$ is an isomorphism.

Finally, by taking the geometric fibers and specializing at $\hbar=0$, we obtain $\Omega_{q,\geo}^{\DK}\simeq\Omega_{q,\quant}^{\DK}$. Then, by Verdier duality, we obtain the corresponding isomorphism for $\Omega_q^{\Lus}$, and by taking the Frobenius quotients, we obtain the isomorphism for $\Omega_q^{\sm}$. 
\end{proof}

\subsubsection{The $!$-stalks of factorization IC modules} Similarly, we can also compare the $!$-stalks of the factorization IC modules in characteristic zero and in sufficiently large positive characteristic.

\begin{prop}\label{prop:ICfact-costalk-dim-spread}
Assume that $q$ has finite order and satisfies condition $(*)$, and let $\theta\leq_{\super}\lambda$. Then, after removing finitely many primes, the dimensions $\dim H^k\!\left(i_{\theta}^!\IC_{\lambda}^{\Lus,\fact}\right)$ are independent of the characteristic of the base field. That is to say, there exists an integer $N(\lambda,\theta)$ such that, for every geometric point $\bar{s}$ of characteristic $p>N(\lambda,\theta)$,
\begin{equation}
\dim H^k\!\left(i_{\theta,\bar{s}}^!\IC_{\lambda,\bar{s}}^{\Lus,\fact}\right) =\dim H^k\!\left(i_{\theta,\eta}^!\IC_{\lambda,\eta}^{\Lus,\fact}\right),
\end{equation}
where $\eta$ denotes the geometric generic point of $S$.
\end{prop}
\begin{proof}
Since the assertion is \'{e}tale-local on the curve, we may assume that $C=\BA^1$ with the marked point $c$. Similarly to the proof of Proposition~\ref{prop:spread-out-DK-comparison}, we work over the component of $\Conf_{\leq\lambda\cdot c,S}$ which contains $\theta\cdot c$, together with the finitely many degrees needed for the factorization. The factorization algebra $\Omega_{q,S}^{\Lus}$ is defined over $S$ after a possible shrinking. Let
\[
j_{\lambda,S}:\Conf_{=\lambda\cdot c,S}\hookrightarrow \Conf_{\leq\lambda\cdot c,S},\qquad \IC_{\lambda,S}^{\Lus,\fact} =j_{\lambda,S,!*}\mathcal L_{\lambda,S}^{\fact},
\]
where $\mathcal L_{\lambda,S}^{\fact}$ is the perverse factorization object obtained from the constant even line along $\lambda\cdot c\simeq S$, as in \cite[Section 5.3.3 and Lemma 5.3.5]{[GL]}.

By the factorization property and the corresponding property of $\Omega_{q,S}^{\Lus}$, it is enough to consider the deepest stratum $\theta\cdot c\simeq S$. As in the proof of Lemma~\ref{lem:formal-saturation-spreading}, by the generic base change and the generic perversity applied to $j_{\lambda,S,!}\mathcal L_{\lambda,S}^{\fact}$, $j_{\lambda,S,*}\mathcal L_{\lambda,S}^{\fact}$ and the image of their canonical map, after shrinking $S$, the intermediate extension commutes with the restriction to geometric fibers, and so does its $!$-restriction along $i_{\theta,S}$. Since the latter is a bounded constructible complex on $S$, after shrinking $S$ once more, its cohomology sheaves are locally constant, and hence, their ranks agree with those at the geometric generic point.
\end{proof}

\begin{conj}
We expect that, if $q$ satisfies condition $(*)$, then the isomorphisms in Proposition~\ref{prop:DK-small-formal-equivalence-Gaiotto} still hold over any ground field.
\end{conj}

\section{The objects corresponding to irreducible representations}\label{sec:irreducibles}

In this section, we consider the functors $F_{\Gr}^{\Lus}$, $F_{\Gr}^{\sm}$, and $F_{\Fl}^{\Lus}$, and study the images of IC sheaves. The main result is Theorem~\ref{thm:image-of-ic}, which says that, if $\lambda$ is restricted dominant and relevant, then $F_{\Gr}^{\sm}$ and $F_{\Gr}^{\Lus}$ send $\IC_{\Gr}^{\lambda}$ to the corresponding irreducible factorization modules, up to the parity change $\Pi^{d(\lambda)}$. The proof is given in Sections~\ref{10.1}--\ref{sec:g1-direction}, and we explain the strategy at the end of Section~\ref{10.1}. Then, we study the metaplectic Hecke action and prove a geometric Steinberg isomorphism (Theorem~\ref{thm:metaplectic-hecke-action}). As an application, we discuss the spherical version and the graded small version of the Gaiotto conjecture.

\subsection{Images of IC sheaves}\label{10.1}

\subsubsection{Relevant orbits}

We consider the homogeneous data in Section~\ref{pol des}, and we use the notation of Section~\ref{sec:orbit-preliminaries}. Recall that the quadratic form $q$ in Section~\ref{3.2.5 sect} determines a symmetric bilinear form $b$, and here, we use the multiplicative convention. For convenience, we take the following notations.

\begin{defn}\label{def:mixed-pure-coroots}
Recall that we denote the simple coroots of $G_i$ by $\alpha_{i,j}$, for $i=1,2$. In the cotangent-type cases that we consider, i.e., $\GL(M|N)$, $\SOSp(2k|2l)$ and $\on F(4)$, there are two types of positive simple coroots of $G_1\times G_2$: we call $\alpha_{i,j}$ \emph{mixed} if it can be written as the sum of two positive simple odd isotropic super roots,
\[
\alpha_{i,j}=\alpha_{i,j}^-+\alpha_{i,j}^+,
\]
and we call it \emph{pure} if it is a positive simple even super root.
\end{defn}

\begin{prop}\label{prop:relevant-L-orbits}
In the three cases of Lemma~\ref{lem:L-orbit-indices}, the $(\cG^G,\psi)$-relevant $G_1(\bO)H_u(\bF)$-orbits are indexed by the subset of dominant coweights $\lambda$ of $G$ which further satisfy the condition: for any mixed simple coroot $\alpha_{i,j}$ of $G$, if
\[
\langle\lambda,\check\alpha_{i,j}\rangle=0,
\]
then
\[
b\bigl(\lambda,\alpha_{i,j}^{\pm}\bigr)=1.
\]
Here, $i=1,2$, and $\alpha_{i,j}=\alpha_{i,j}^-+\alpha_{i,j}^+$ for two odd simple super roots.
\end{prop}

\begin{proof}
Since the orbit is $\psi$-relevant, it remains to check that the twisting is trivial on the reductive quotient of the stabilizer $K_\lambda=H(\bF)\cap t^\lambda G(\bO)t^{-\lambda}$. Note that the projection of $K_\lambda$ to $G_1(\bF)$ has a pro-unipotent kernel. For the reductive quotient, in type A, we use \cite[proof of Proposition 3.3.4]{[TY]} and \cite[proof of Proposition 4.1.1]{[BFT1]}, see also \cite[Proposition 7]{Po}. For the orthosymplectic and $\on F(4)$ cases, see \cite[proof of Lemma 2.3.3]{[BFT]} and \cite[proof of Proposition 3.3.1(c)]{[FTY]}.

By these stabilizer calculations, we can reduce to the rank-one case for each mixed simple coroot $\alpha_{i,j}$ with $\langle\lambda,\check\alpha_{i,j}\rangle=0$. In this case, the twisting character on the corresponding one-dimensional torus is $b(\lambda,\alpha_{i,j}^{\pm})$ up to inversion, and
\[
b(\lambda,\alpha_{i,j}^-)b(\lambda,\alpha_{i,j}^+) =q(\alpha_{i,j})^{\langle\lambda,\check\alpha_{i,j}\rangle}=1.
\]
Thus, $b(\lambda,\alpha_{i,j}^-)=1$ if and only if $b(\lambda,\alpha_{i,j}^+)=1$.
\end{proof}

Assume that $\lambda$ is dominant with respect to $B^-$ and that its orbit in $\Gr_G^\omega$ is relevant, i.e.,
\[
\SD_{\cG^G}^{H^\omega(\bF),\psi} \bigl(\sO_{\Gr}^{\lambda}\bigr)\neq 0.
\]
We denote by $\IC_{\Gr,\glob}^{\lambda}$ the IC sheaf of the corresponding global stratum $\cM_{=\lambda\cdot c}$, and define $\IC_{\Gr}^{\lambda}$ to be its image under the local-global equivalence, with the normalization of Remark~\ref{rem:component-normalization}. Then, the restriction of $\IC_{\Gr}^{\lambda}$ to $\sO_{\Gr}^{\lambda}$ is the unique (since the stabilizer subgroup at $t^\lambda$ is connected) irreducible rank-one equivariant twisted local system with even coefficients, with the cohomological shift of Remark~\ref{rem:component-normalization}.

\subsubsection{Restricted coweights}

Recall that in \cite{[GL]}, a dominant coweight (on the geometric side) $\lambda\in\Lambda_G^+$ is called \emph{restricted} if
\[
0\leq \langle\lambda,\check\alpha_i\rangle <\operatorname{ord}\bigl(q(\alpha_i)\bigr)
\]
for every simple coroot $\alpha_i$ of $G$. On the dual side, these inequalities correspond to the usual restricted highest weights for the small quantum group.

In the super setting, the restricted condition is a bit different. Namely, for a mixed simple coroot, we also allow the upper bound. In this case, $b$ must satisfy an additional condition. For a simple coroot $\alpha_{i,j}$, let us denote
\[
q_{i,j}:=\operatorname{ord}\bigl(q(\alpha_{i,j})\bigr).
\]

\begin{defn}\label{def:restricted-coweight}
We call a dominant coweight $\lambda\in\Lambda_G^+$ \emph{restricted} if the following conditions hold.
\begin{enumerate}
\item If $\alpha_{i,j}$ is pure, then
\begin{equation}\label{eq:restricted-pure-inequality}
0\leq \langle\lambda,\check\alpha_{i,j}\rangle <q_{i,j}.
\end{equation}
\item If $\alpha_{i,j}$ is mixed, then
\begin{equation}\label{eq:restricted-mixed-inequality}
0\leq \langle\lambda,\check\alpha_{i,j}\rangle \leq q_{i,j},
\end{equation}
and if the equality holds in the upper bound, then
\begin{equation}\label{eq:restricted-wall-monodromy}
b\bigl(\lambda,\alpha_{i,j}^{\pm}\bigr)\neq 1.
\end{equation}
\end{enumerate}
\end{defn}

\subsubsection{The IC comparison and its reduction}

In this section, we will eventually prove the following theorem.

\begin{thm}\label{thm:image-of-ic}
Let $\lambda=(\lambda^{(1)},\lambda^{(2)})$ be restricted dominant and relevant. Then,
\[
F_{\Gr}^{\sm}\bigl(\IC_{\Gr}^{\lambda}\bigr) \simeq \Pi^{d(\lambda)}\IC_{\lambda}^{\sm,\fact}, \qquad F_{\Gr}^{\Lus}\bigl(\IC_{\Gr}^{\lambda}\bigr) \simeq \Pi^{d(\lambda)}\IC_{\lambda}^{\Lus,\fact}.
\]
Here, $\IC_{\lambda}^{\sm,\fact}$ and $\IC_{\lambda}^{\Lus,\fact}$ are the irreducible objects in the corresponding factorization module categories.
\end{thm}

For the calculations of local systems below, we work either in the $D$-module setting in characteristic zero, or in the $\ell$-adic setting in sufficiently large positive characteristic, as in Section~\ref{sec:positive-characteristic}.

First, we reduce the proof to a statement about the local systems $\cL_{\lambda,\theta}\otimes\psi$ of Section~\ref{sec:ratio-local-system} on the critical-dimensional irreducible components of $\sO_{\Gr}^{\lambda}\cap\sS_{\Gr}^{\theta}$.

Roughly speaking, by Proposition~\ref{prop 7.0.5} and Corollary~\ref{lem:dimension-estimate-super}, the $!$-stalks of the images of $\IC_{\Gr}^{\lambda}$ at $\theta\cdot c$ have $H^0\neq0$ only if $\cL_{\lambda,\theta}\otimes\psi$ is constant on a dense open subset of a critical-dimensional irreducible component of $\sO_{\Gr}^{\lambda}\cap\sS_{\Gr}^{\theta}$. The following theorem shows that this can happen only when $\theta=\lambda$.

\begin{thm}\label{thm:local-system-rigidity}
Let $\lambda=(\lambda^{(1)},\lambda^{(2)})$ be restricted dominant and relevant. Assume that
\[
\sO_{\Gr}^{(\lambda^{(1)},\lambda^{(2)})} \cap \sS_{\Gr}^{(\theta^{(1)},\theta^{(2)})}
\]
has a critical-dimensional irreducible component, such that the restriction of $\cL_{\lambda,\theta}\otimes\psi$ to a dense open subset of it is constant. Then,
\[
(\lambda^{(1)},\lambda^{(2)})=(\theta^{(1)},\theta^{(2)}).
\]
\end{thm}

\begin{proof}[Proof of Theorem~\ref{thm:image-of-ic}, assuming Theorem~\ref{thm:local-system-rigidity}]
Let $\dagger\in\{\sm,\Lus\}$. By Theorem~\ref{thm:Grassmannian-functors-texact}, $F_{\Gr}^{\dagger}(\IC_{\Gr}^{\lambda})$ is perverse. Thus, we should prove that, for $\theta\neq\lambda$, its $!$-stalk and $*$-stalk at $\theta\cdot c$ vanish in degree $0$.

First, we consider the $!$-stalk. Since $\IC_{\Gr}^{\lambda}$ is the $!*$-extension from $\sO_{\Gr}^{\lambda}$, by Proposition~\ref{prop 7.0.5} and Corollary~\ref{lem:dimension-estimate-super}, the orbits in the boundary of $\sO_{\Gr}^{\lambda}$ only contribute to the $!$-stalk in strictly positive degrees. Thus, a contribution in degree zero can only come from a critical-dimensional irreducible component of $\sO_{\Gr}^{\lambda}\cap\sS_{\Gr}^{\theta}$ on which $\cL_{\lambda,\theta}\otimes\psi$ is constant, and by Theorem~\ref{thm:local-system-rigidity}, it happens only if $\theta=\lambda$. Hence,
\begin{equation}\label{eq:restricted-IC-boundary-costalk}
H^0\bigl(i_\theta^!F_{\Gr}^{\dagger}(\IC_{\Gr}^{\lambda})\bigr)=0 \qquad(\theta\neq\lambda).
\end{equation}

Secondly, for the $*$-stalk, we use duality. Note that the small functor commutes with duality, and the dual of the Lusztig functor is $F_{\Gr}^{\DK}$. Since $\BD\IC_{\Gr}^{\lambda}$ is also the $!*$-extension, the same dimension estimate shows that the boundary orbits do not contribute in degree zero. Thus, by Theorem~\ref{thm:local-system-rigidity}, we obtain
\begin{equation}\label{eq:restricted-DK-boundary-costalk}
H^0\bigl(i_\theta^!F_{\Gr}^{\DK}(\BD\IC_{\Gr}^{\lambda})\bigr)=0 \qquad(\theta\neq\lambda).
\end{equation}
By taking the duals, we obtain the desired vanishing of the $*$-stalk. Note that, on the stratum $\Conf_{=\lambda\cdot c}$, the coefficient line is $\Pi^{d(\lambda)}\se$. Finally, by the standard factorization arguments and Remark~\ref{rem:target-parity}, we identify the resulting object with $\Pi^{d(\lambda)}\IC_{\lambda}^{\dagger,\fact}$.
\end{proof}

The proof of Theorem~\ref{thm:local-system-rigidity} is divided into two parts, corresponding to the $G_2$-part and the $G_1$-part of the coweights.

\begin{prop}\label{prop:part-q}
Under the assumptions of Theorem~\ref{thm:local-system-rigidity}, we have $\lambda^{(2)}=\theta^{(2)}$.
\end{prop}

\begin{prop}\label{prop:part-p}
Under the assumptions of Theorem~\ref{thm:local-system-rigidity}, we have $\lambda^{(1)}=\theta^{(1)}$.
\end{prop}

The proofs of Propositions~\ref{prop:part-q} and~\ref{prop:part-p} are the main technical parts of this section, and we first explain the idea. In the Whittaker case, one needs to consider the intersections of two semi-infinite orbits for opposite Borel subgroups, and Gaitsgory's argument in \cite[\S 7.4]{[G7]} reduces the problem to the rank-one case by using Kashiwara's crystal. In our case, $\sS_{\Gr}^{\theta}=\sS_{G_1}^{-,\theta^{(1)}}\times\sS_{G_2}^{+,\theta^{(2)}}$ is a product, but $\sO_{\Gr}^{\lambda}$ is not, since $G_1$ is embedded diagonally into $G_1\times G_2$. Thus, we project $\sO_{\Gr}^{\lambda}\cap\sS_{\Gr}^{\theta}$ to $\Gr_{G_1}$ and to $\Gr_{G_2}$, and estimate the dimensions of the strata of the images and of the fibers over them. It turns out that, for each projection, only one (!) stratum can have critical dimension, see Corollary~\ref{rem:unique-g2-critical-stratum} and Proposition~\ref{prop:g1-strata-dimension}. We will give these orbit and dimension calculations in Section~\ref{10.2}.

Then, in Section~\ref{sec:g2-direction}, we use these calculations to reduce Proposition~\ref{prop:part-q} to a statement about an intersection in $\Gr_{G_2}$. This intersection is contained in the intersection of two semi-infinite orbits of $G_2$ for opposite Borel subgroups. Furthermore, we estimate the codimension and reduce to the case where the similar argument for the Whittaker case applies. So, to prove Proposition~\ref{prop:part-q}, it is enough to consider the rank-one case $\lambda^{(2)}-\theta^{(2)}=q_{2,i}\alpha_{2,i}$. If the simple coroot is pure, this case is excluded by \cite[Theorem 9.4.9]{[GL]} and Lemma~\ref{lem:kummer-as-separation}; and if it is mixed, we will exclude this case by explicit matrix calculations, case by case. 

Finally, in Section~\ref{sec:g1-direction}, we prove Proposition~\ref{prop:part-p} by a similar general strategy (with some technical differences in details), with the roles of $G_1$ and $G_2$ interchanged. 

\ifdim\dimexpr\pagegoal-\pagetotal\relax<7\baselineskip\newpage\fi
\subsection{Common orbit and dimension calculations}\label{10.2}

As explained at the end of Section~\ref{10.1}, we study $\sO_{\Gr}^{\lambda}\cap\sS_{\Gr}^{\theta}$ by projecting it to $\Gr_{G_2}$ and to $\Gr_{G_1}$. For each projection, we stratify the image, and estimate the dimensions of the strata and of the fibers over them. Then, we determine which preimages of strata can have critical dimension. Here, we use the classification of orbits and the closure relations in Section~\ref{sec:orbit-preliminaries}.

\subsubsection{The projection to \texorpdfstring{$\Gr_{G_2}$}{Gr_G2} and its strata}\label{sec:g2-projection-strata}

Recall that $\sD$ is fixed in Section~\ref{trans D}. Consider the following equivariant product
\begin{equation}
\label{eq:Y-convolution}\fitmath{Y_{\theta^{(1)};\lambda^{(1)},\lambda^{(2)}} := \bigl(\sD^{-1}t^{\theta^{(1)}} U_1^-(\bF)G_1(\bO)\bigr) \times^{G_1(\bO)} \bigl(G_1(\bO)H_u(\bF)t^{-\lambda^{(1)}}\sD t^{\lambda^{(2)}} G_2(\bO)/G_2(\bO)\bigr).}
\end{equation}
Let
\begin{equation}\label{eq:g2-projection}
m_{\theta^{(1)};\lambda^{(1)},\lambda^{(2)}} \colon Y_{\theta^{(1)};\lambda^{(1)},\lambda^{(2)}}\longrightarrow \Gr_{G_2}
\end{equation}
be the map induced by the multiplication in $G_2(\bF)$. Then, we have
\begin{equation}\label{eq:intersection-as-g2-fiber-product}
\sO_{\Gr}^{(\lambda^{(1)},\lambda^{(2)})} \cap \sS_{\Gr}^{(\theta^{(1)},\theta^{(2)})} \simeq Y_{\theta^{(1)};\lambda^{(1)},\lambda^{(2)}} \times_{\Gr_{G_2}}\sS_{G_2}^{\theta^{(2)}}.
\end{equation}
Under this isomorphism, the projection $\sO_{\Gr}^{\lambda}\cap\sS_{\Gr}^{\theta}\to\sS_{G_2}^{\theta^{(2)}}$ is the base change of $m_{\theta^{(1)};\lambda^{(1)},\lambda^{(2)}}$.

The image of $Y_{\theta^{(1)};\lambda^{(1)},\lambda^{(2)}}$ under $m_{\theta^{(1)};\lambda^{(1)},\lambda^{(2)}}$ is
\[
\mathcal A_{\theta^{(1)};\lambda^{(1)},\lambda^{(2)}} := \sD^{-1}t^{\theta^{(1)}} U_1^-(\bF)G_1(\bO)H_u(\bF) t^{-\lambda^{(1)}}\sD t^{\lambda^{(2)}} G_2(\bO)/G_2(\bO).
\]
Clearly, $\mathcal A_{\theta^{(1)};\lambda^{(1)},\lambda^{(2)}}$ is stratified by the $\sD^{-1}t^{\theta^{(1)}}$-shifted $U_1^-(\bF)T_1(\bO)H_u(\bF)$-orbits
\begin{equation}\label{eq:g2-image-stratum}
\CS_{G_2}^{(\theta^{(1)};\lambda^{(1)\prime},\lambda^{(2)\prime})} := \sD^{-1}t^{\theta^{(1)}} U_1^-(\bF)T_1(\bO)H_u(\bF) t^{-\lambda^{(1)\prime}}\sD t^{\lambda^{(2)\prime}}G_2(\bO)/G_2(\bO).
\end{equation}

To classify these strata, we need the following lemma, which is an analogue of Lemma~\ref{lem:L-orbit-indices}, with $G_1(\bO)$ replaced by $U_1^-(\bF)T_1(\bO)$.

\begin{lem}\label{lem:semi-spherical-orbits}
The $U_1^-(\bF)T_1(\bO)H_u(\bF)$-orbits in $\Gr_{G_2}$ are indexed by pairs $\lambda=(\lambda^{(1)},\lambda^{(2)})$ satisfying the same conditions as in Lemma~\ref{lem:L-orbit-indices}, while we remove the requirement of the dominance condition for $\lambda^{(1)}$. A representative is given by $t^{-\lambda^{(1)}}\sD t^{\lambda^{(2)}}\in\Gr_{G_2}^\omega$. If the orbit indexed by $\lambda'$ is contained in the closure of the orbit indexed by $\lambda$, then $\lambda'\leq_{\super}\lambda$. We denote the index set by $\Lambda_G^L$.
\end{lem}

\begin{proof}
Let us compare the orbits of the following two groups,
\[
K=U_1^-(\bF)T_1(\bO)H_u(\bF),\qquad K_\eta=t^\eta G_1(\bO)t^{-\eta}H_u(\bF),
\]
where $\eta\in\Lambda_{G_1}$. We claim that, for fixed $x,y\in\Gr_{G_2}$ and sufficiently dominant $\eta$, we have
\[
y\in Kx\quad\Longleftrightarrow\quad y\in K_\eta x.
\]

Indeed, if $y=gx$ with $g\in K$, then $g\in K_\eta$ for sufficiently dominant $\eta$. Conversely, let $g\in K_\eta$ with $gx=y$, and write $g=ug_1$, where $u\in H_u(\bF)$ and $g_1\in t^\eta G_1(\bO)t^{-\eta}$. For fixed $x,y$ and sufficiently dominant $\eta$, we can write $g_1=u_-au_+$, where $u_-\in U_1^-(\bF)$, $a\in T_1(\bO)$, and $u_+\in t^\eta U_1^+(\bO)t^{-\eta}$. Note that $u_+$ is congruent to $1$ modulo a sufficiently high power of $t$, thus, it fixes $x$. Hence, $y=uu_-ax\in Kx$.

Now, up to the left translation by $t^{-\eta}$, the representatives and their uniqueness follow from Lemma~\ref{lem:L-orbit-indices}, where the indices become $(\lambda^{(1)}+\eta,\lambda^{(2)})$.

For the closure relation, we use the above comparison with $K_\eta$, and translate by $t^{-\eta}$ again. Then, by Lemma~\ref{lem:spherical-orbit-closure}, we have $(\lambda^{(1)\prime}+\eta,\lambda^{(2)\prime}) \leq_{\super}(\lambda^{(1)}+\eta,\lambda^{(2)})$, which implies $(\lambda^{(1)\prime},\lambda^{(2)\prime}) \leq_{\super}(\lambda^{(1)},\lambda^{(2)})$.
\end{proof}

\begin{lem}\label{lem:g2-strata-in-image}
A stratum $\CS_{G_2}^{(\theta^{(1)};\lambda^{(1)\prime},\lambda^{(2)\prime})}$ is contained in $\mathcal A_{\theta^{(1)};\lambda^{(1)},\lambda^{(2)}}$ (the image of \eqref{eq:g2-projection}) only if its index satisfies
\[
(\lambda^{(1)\prime},\lambda^{(2)\prime})\leq_{\super}(\lambda^{(1)},\lambda^{(2)}),
\]
and the difference $(\lambda^{(1)},\lambda^{(2)})-(\lambda^{(1)\prime},\lambda^{(2)\prime})$ is a non-negative integral linear combination of the positive odd simple super roots. Moreover, every such $\lambda^{(2)\prime}$ is dominant.
\end{lem}

\begin{proof}
Consider the convolution diagram
\begin{equation}\label{eq:dominant-translation-convolution}
t^\eta G_1(\bO)t^{-\eta}G_1(\bO) \overset{G_1(\bO)}{\times} G_1(\bO)H_u(\bF)t^{-\lambda^{(1)}}\sD t^{\lambda^{(2)}} G_2(\bO)/G_2(\bO)\longrightarrow\Gr_{G_2}.
\end{equation}
One can check that if we require $\lambda^{(1)},\lambda^{(2)},\eta$ to be dominant in $G$, then $t^\eta G_1(\bO)H_u(\bF)t^{-\eta-\lambda^{(1)}}\sD t^{\lambda^{(2)}}G_2(\bO)/G_2(\bO)$ is open dense in the image.

Now, assume that $\CS_{G_2}^{(\theta^{(1)};\lambda^{(1)\prime},\lambda^{(2)\prime})}$ is in $\mathcal A_{\theta^{(1)};\lambda^{(1)},\lambda^{(2)}}$. We translate by $t^{-\theta^{(1)}}\sD$. By Lemma~\ref{lem:semi-spherical-orbits}, there exists a very dominant $\eta$ such that $t^{-\lambda^{(1)\prime}}\sD t^{\lambda^{(2)\prime}}$ belongs to the image of \eqref{eq:dominant-translation-convolution}. After translating by $t^{-\eta}$, this implies $(\lambda^{(1)}+\eta,\lambda^{(2)})\geq_{\super}(\lambda^{(1)\prime}+\eta,\lambda^{(2)\prime})$, i.e., $(\lambda^{(1)},\lambda^{(2)})\geq_{\super}(\lambda^{(1)\prime},\lambda^{(2)\prime})$, see Lemma~\ref{lem:spherical-orbit-closure}.

Now, we claim that the difference $(\lambda^{(1)},\lambda^{(2)})-(\lambda^{(1)\prime},\lambda^{(2)\prime})$ is spanned by positive odd simple super roots.

Note that $\mathcal A_{\theta^{(1)};\lambda^{(1)},\lambda^{(2)}}$ is contained in
\[
\sD^{-1}G_1(\bF)\sD H_u(\bF)t^{\lambda^{(2)}}G_2(\bO)/G_2(\bO).
\]
Since $\CS_{G_2}^{(\theta^{(1)};\lambda^{(1)\prime},\lambda^{(2)\prime})}$ is contained in the above subspace, we have
\begin{itemize}
\item in the $\GL(M|N)$ case, $\lambda^{(2)\prime}_j=\lambda^{(2)}_j$ for $M+2\leq j\leq N$,
\item in the $\SOSp(2k|2l)$ case,
\begin{enumerate}
\item if $k\leq l$, i.e., $M$ is even, then $\lambda^{(2)\prime}_j=\lambda^{(2)}_j$ for $1\leq j\leq l-k$,
\item if $k>l$, i.e., $M$ is odd, then $\lambda^{(2)\prime}_j=\lambda^{(2)}_j$ for $1\leq j\leq k-l-1$,
\end{enumerate}
\item in the $\on F(4)$ case, $\lambda^{(2)\prime}_1+\lambda^{(2)\prime}_2=\lambda^{(2)}_1+\lambda^{(2)}_2$.
\end{itemize}
Thus, by the above equalities and the simple super roots in Section~\ref{super Borel}, the difference $(\lambda^{(1)},\lambda^{(2)})-(\lambda^{(1)\prime},\lambda^{(2)\prime})$ is a non-negative integral linear combination of positive odd simple super roots.

To prove that $\lambda^{(2)\prime}$ is dominant, we first notice that, by Lemma~\ref{lem:semi-spherical-orbits}, it satisfies the dominance inequalities for $L$. Thus, since $\lambda^{(2)}$ is dominant and the above equalities hold, we only need to show that
\begin{itemize}
\item in the $\GL(M|N)$ case, if $M+1<N$, then $\lambda^{(2)\prime}_{M+1}\leq\lambda^{(2)\prime}_{M+2}$,
\item in the $\SOSp(2k|2l)$ case,
\begin{enumerate}
\item if $k<l$, then $\lambda^{(2)\prime}_{l-k}\leq\lambda^{(2)\prime}_{l-k+1}$,
\item if $k>l+1$, then $\lambda^{(2)\prime}_{k-l-1}\leq\lambda^{(2)\prime}_{k-l}$,
\end{enumerate}
\item in the $\on F(4)$ case, $\lambda^{(2)\prime}_2-\lambda^{(2)\prime}_3\leq0$.
\end{itemize}
For $\on F(4)$, the odd simple super roots are $\alpha_1,\alpha_2,\alpha_3$, see Section~\ref{super Borel}. Subtracting $\alpha_1$ or $\alpha_2$ does not change $\lambda^{(2)}_2-\lambda^{(2)}_3$, while subtracting $\alpha_3$ decreases it by $1$, and thus, $\lambda^{(2)\prime}_2-\lambda^{(2)\prime}_3\leq\lambda^{(2)}_2-\lambda^{(2)}_3\leq0$. The first two cases are similar.
\end{proof}

\subsubsection{Critical strata for the projection to \texorpdfstring{$\Gr_{G_2}$}{Gr_G2}}

Now, we study the condition under which the intersections of the strata \eqref{eq:g2-image-stratum} with $\sS_{G_2}^{\theta^{(2)}}$ have preimages of critical dimension under \eqref{eq:g2-projection}. First, we estimate the dimensions of these intersections (Lemma~\ref{lem:g2-base-dimension} and Proposition~\ref{prop:g2-base-critical}), and then estimate the dimensions of the fibers \eqref{eq:g2-projection-fiber} (Proposition~\ref{prop:g2-fiber-dimension}). We will see in Corollary~\ref{rem:unique-g2-critical-stratum} that the critical dimension can only be attained for $(\lambda^{(1)\prime},\lambda^{(2)\prime})=(\theta^{(1)},\lambda^{(2)})$.

\begin{lem}\label{lem:g2-base-dimension}
For any $\lambda'=(\lambda^{(1)\prime},\lambda^{(2)\prime})\in\Lambda_G^L$ and $\theta=(\theta^{(1)},\theta^{(2)})\in\Lambda$, we have
\[
\begin{aligned}
&\dim\bigl( \CS_{G_2}^{(\theta^{(1)};\lambda^{(1)\prime},\lambda^{(2)\prime})} \cap\sS_{G_2}^{\theta^{(2)}} \bigr)\\
&\quad\leq \max_{\substack{m_{1,i},m_{2,j}\in\BZ_{\geq0}\\
\lambda' -\sum_i m_{1,i}\alpha_{1,i} -\sum_j m_{2,j}\alpha_{2,j} \geq_{\super}\theta}} \left(\sum_i m_{1,i}+\sum_j m_{2,j}\right).
\end{aligned}
\]
In particular,
\[
\dim\bigl( \CS_{G_2}^{(\theta^{(1)};\lambda^{(1)\prime},\lambda^{(2)\prime})} \cap\sS_{G_2}^{\theta^{(2)}} \bigr) \leq \frac12d(\lambda'-\theta).
\]
Furthermore, the dimension is critical only if
\[
(\lambda^{(1)\prime},\lambda^{(2)\prime})\geq_G(\theta^{(1)},\theta^{(2)}),
\]
i.e., the difference belongs to the positive coroot monoid of $G_1\times G_2$.
\end{lem}

\begin{proof}
Note that the stratum $\CS_{G_2}^{(\theta^{(1)};\lambda^{(1)\prime},\lambda^{(2)\prime})}$ only depends on $\theta^{(1)}-\lambda^{(1)\prime}$ and $\lambda^{(2)\prime}$. Thus, by translating $\lambda^{(1)\prime}$ and $\theta^{(1)}$ by the same coweight, we can assume that $\lambda^{(1)\prime}$ is dominant. Similarly, by translating $\lambda^{(2)\prime}$ and $\theta^{(2)}$ by the same sufficiently dominant coweight of $Z_L$, we can also assume that $\lambda^{(2)\prime}$ is dominant, and it changes neither the dimension of $\CS_{G_2}^{(\theta^{(1)};\lambda^{(1)\prime},\lambda^{(2)\prime})}\cap\sS_{G_2}^{\theta^{(2)}}$ nor the difference $\lambda'-\theta$. Now, the dimension of this intersection is bounded by the dimension of $\overline{\sO}_{\Gr}^{\lambda'}\cap\overline{\sS}_{\Gr}^{\theta}$, and if the former has critical dimension, then the latter is also critical. Hence, both the estimate and the condition for the equality follow from Corollary~\ref{lem:dimension-estimate-super}.
\end{proof}

\begin{cor}\label{cor:g2-boundary-codimension}
Assume that $\CS_{G_2}^{(\theta^{(1)};\lambda^{(1)\prime},\lambda^{(2)\prime})}$ is contained in $\mathcal A_{\theta^{(1)};\lambda^{(1)},\lambda^{(2)}}$, and that $\CS_{G_2}^{(\theta^{(1)};\lambda^{(1)\prime},\lambda^{(2)\prime})}\cap\sS_{G_2}^{\theta^{(2)}}$ has critical dimension $\frac12d(\lambda'-\theta)$. Then, the intersection of the boundary of $\CS_{G_2}^{(\theta^{(1)};\lambda^{(1)\prime},\lambda^{(2)\prime})}$ and the boundary of $\sS_{G_2}^{\theta^{(2)}}$ has dimension no more than
\[
\frac12d(\lambda'-\theta)-2.
\]
\end{cor}

\begin{proof}
By Lemma~\ref{lem:semi-spherical-orbits}, if a stratum is in the boundary of $\CS_{G_2}^{(\theta^{(1)};\lambda^{(1)\prime},\lambda^{(2)\prime})}$, then its index is strictly smaller than $\lambda'$ with respect to $\leq_{\super}$. Similarly, if a semi-infinite orbit is in the boundary of $\sS_{G_2}^{\theta^{(2)}}$, then its index is $\theta^{(2)}$ plus a non-zero sum of positive coroots. Thus, by Lemma~\ref{lem:g2-base-dimension}, the dimension bound for their intersection drops by at least $\frac12+1$. Since the original critical dimension is an integer, the intersection of the boundaries has dimension at most $\frac12d(\lambda'-\theta)-2$.
\end{proof}

\begin{prop}\label{prop:g2-base-critical}
Let $\lambda'=(\lambda^{(1)\prime},\lambda^{(2)\prime})\in\Lambda_G^L$. If $\CS_{G_2}^{(\theta^{(1)};\lambda^{(1)\prime},\lambda^{(2)\prime})}$ is contained in $\mathcal A_{\theta^{(1)};\lambda^{(1)},\lambda^{(2)}}$, then the intersection
\[
\CS_{G_2}^{(\theta^{(1)};\lambda^{(1)\prime},\lambda^{(2)\prime})} \cap\sS_{G_2}^{\theta^{(2)}}
\]
can have critical dimension $\frac12d(\lambda'-\theta)$ only if
\[
\lambda^{(1)\prime}=\theta^{(1)}, \qquad \lambda^{(2)\prime}-\theta^{(2)}\in\Lambda_{G_2}^{\on{pos}}.
\]
\end{prop}

\begin{proof}
By Lemma~\ref{lem:g2-base-dimension}, we have $\lambda'\geq_G\theta$.

Let $n=\frac12d(\lambda'-\theta)$. We consider the closure of $\CS_{G_2}^{(\theta^{(1)};\lambda^{(1)\prime},\lambda^{(2)\prime})}$ inside $\Gr_{G_2}$, and intersect it with $\overline{\sS}_{G_2}^{\theta^{(2)}}$. By the Mirkovi\'c--Vilonen hyperplane argument, the boundary of a semi-infinite orbit appears as a hyperplane in a projective embedding of the affine Grassmannian. Therefore, if an irreducible component of the intersection has critical dimension $n$, one can cut successively by boundary hyperplanes and obtain a sequence
\[
\sS_{G_2}^{\theta^{(2)}} \rightsquigarrow \sS_{G_2}^{\theta^{(2)}+\beta_1} \rightsquigarrow \cdots \rightsquigarrow \sS_{G_2}^{\theta^{(2)}+\beta_1+\cdots+\beta_n},
\]
where each $\beta_j$ is a positive simple coroot of $G_2$, such that after each cut, the corresponding intersection remains critical and its dimension drops by exactly one. This is the argument of \cite[\S 3]{[MV]} and is used in the same form in \cite[\S 1.5, Proposition 4.4.1]{[SW]}.

Indeed, let us consider the first cut. Every irreducible component of the hyperplane of the closure of this component has dimension at least $n-1$, and by Corollary~\ref{cor:g2-boundary-codimension}, it is not contained in the boundary of $\CS_{G_2}^{(\theta^{(1)};\lambda^{(1)\prime},\lambda^{(2)\prime})}$. Thus, a dense open subset of it is contained in $\CS_{G_2}^{(\theta^{(1)};\lambda^{(1)\prime},\lambda^{(2)\prime})}\cap\sS_{G_2}^{\theta^{(2)}+\beta}$ for some non-zero $\beta\in\Lambda_{G_2}^{\on{pos}}$, and by Lemma~\ref{lem:g2-base-dimension}, we have
\[
n-1\leq\tfrac12d\bigl(\lambda'-(\theta^{(1)},\theta^{(2)}+\beta)\bigr)=n-\tfrac12d(\beta).
\]
Since $d$ takes the value $2$ on every positive simple coroot, $\beta$ is a simple coroot, and the intersection $\CS_{G_2}^{(\theta^{(1)};\lambda^{(1)\prime},\lambda^{(2)\prime})}\cap\sS_{G_2}^{\theta^{(2)}+\beta}$ is again of critical dimension. Hence, we can repeat the argument.

Now, after the $n$-th cut, the intersection $\CS_{G_2}^{(\theta^{(1)};\lambda^{(1)\prime},\lambda^{(2)\prime})}\cap\sS_{G_2}^{\theta^{(2)}+\beta_1+\cdots+\beta_n}$ is non-empty and of critical dimension $0$. Thus, by Lemma~\ref{lem:g2-base-dimension} again, we have
\[
\lambda'\geq_G(\theta^{(1)},\theta^{(2)}+\beta_1+\cdots+\beta_n),
\]
and $\theta$ only differs from $\lambda'$ by exactly $n$-many simple coroots. It forces that $\lambda'=(\theta^{(1)},\theta^{(2)}+\beta_1+\cdots+\beta_n)$, since simple coroots are linearly independent, i.e., $\lambda^{(1)\prime}=\theta^{(1)}$ and $\lambda^{(2)\prime}-\theta^{(2)}=\beta_1+\cdots+\beta_n\in\Lambda_{G_2}^{\on{pos}}$.
\end{proof}

Next, we should bound the fibers. By the equivariant product \eqref{eq:Y-convolution}, the fiber of \eqref{eq:g2-projection} over the representative
\[
\sD^{-1}t^{\theta^{(1)}}t^{-\lambda^{(1)\prime}}\sD t^{\lambda^{(2)\prime}}G_2(\bO)/G_2(\bO)
\]
of the stratum \eqref{eq:g2-image-stratum} is identified with
\begin{equation}\label{eq:g2-projection-fiber}
\begin{split}
&\sD t^{\lambda^{(2)\prime}}G_2(\bO)t^{-\lambda^{(2)}}\sD^{-1}t^{\lambda^{(1)}} H_u(\bF)G_1(\bO)/G_1(\bO)\\
&\hspace{35mm}\cap U_1^-(\bF)t^{\lambda^{(1)\prime}}G_1(\bO)/G_1(\bO) \subset G_2(\bF)/G_1(\bO).
\end{split}
\end{equation}
Note that it also equals the fiber of $m_{\lambda^{(1)\prime};\lambda^{(1)},\lambda^{(2)}}$ over the point $t^{\lambda^{(2)\prime}}G_2(\bO)/G_2(\bO)$.

\begin{prop}\label{prop:g2-fiber-dimension}
For $(\lambda^{(1)},\lambda^{(2)})\in\Lambda_G^{L,+}$ and any $\lambda^{(1)\prime},\lambda^{(2)\prime},\theta^{(1)}$, the dimension of each fiber of the projection \eqref{eq:g2-projection} over a point of the stratum $\CS_{G_2}^{(\theta^{(1)};\lambda^{(1)\prime},\lambda^{(2)\prime})}$ is bounded by
\[
\frac12d\bigl((\lambda^{(1)}-\lambda^{(1)\prime},\lambda^{(2)}-\lambda^{(2)\prime})\bigr).
\]
Furthermore, the equality holds only if $\lambda^{(2)}=\lambda^{(2)\prime}$, and in this case, we have
\[
(\lambda^{(1)},\lambda^{(2)})\geq_G(\lambda^{(1)\prime},\lambda^{(2)\prime}).
\]
\end{prop}

\begin{proof}
Similar to the proof of Lemma~\ref{lem:g2-base-dimension}, by translating $\lambda^{(2)}$ and $\lambda^{(2)\prime}$ by the same coweight of $Z_L$, we can assume that $\lambda^{(2)}$ is dominant. Note that we have
\begin{equation}\label{eq:fiber-bound-by-full-intersection}
\sO_{\Gr}^{(\lambda^{(1)},\lambda^{(2)})} \cap\sS_{\Gr}^{(\lambda^{(1)\prime},\lambda^{(2)\prime})} \simeq Y_{\lambda^{(1)\prime};\lambda^{(1)},\lambda^{(2)}} \times_{\Gr_{G_2}}\sS_{G_2}^{\lambda^{(2)\prime}}.
\end{equation}
We denote $\lambda'=(\lambda^{(1)\prime},\lambda^{(2)\prime})$. According to Corollary~\ref{lem:dimension-estimate-super}, we have
\[
\dim\bigl(\overline{\sO}_{\Gr}^{\lambda}\cap\overline{\sS}_{\Gr}^{\lambda'}\bigr) \leq \frac12d(\lambda-\lambda').
\]
By \eqref{eq:fiber-bound-by-full-intersection}, the fiber \eqref{eq:g2-projection-fiber} is a subset of $\sO_{\Gr}^{\lambda}\cap\sS_{\Gr}^{\lambda'}$, and thus, we know that its dimension is no more than $\frac12d(\lambda-\lambda')$. Also, if the equality holds, then the same corollary implies that $\lambda\geq_G\lambda'$.

Then, we use the Mirkovi\'c--Vilonen hyperplane argument again to show that the dimension equals $\frac12d(\lambda-\lambda')$ only if $\lambda^{(2)}=\lambda^{(2)\prime}$. Assume that the equality holds, and denote $n=\frac12d(\lambda-\lambda')$. First, we note that the closure of an $n$-dimensional irreducible component of the fiber \eqref{eq:g2-projection-fiber} is actually contained in $\Gr_{G_1}\times\{t^{\lambda^{(2)\prime}}\}$. Thus, in the argument of Corollary~\ref{lem:dimension-estimate-super}, we have a sequence of $n$ cuts, which only occur in the $G_1$-direction, i.e., the $n$ changes of index are non-zero sums of positive simple coroots of $G_1$. In particular, the value of $d$ on their sum is at least $2n$. On the other hand, by Lemma~\ref{lem:orbit-S-nonempty-super-order}, their sum is at most $\lambda-\lambda'$ in the super root order. Since $d(\lambda-\lambda')=2n$, the difference $\lambda-\lambda'$ equals the sum of these changes, and hence, $\lambda^{(2)}=\lambda^{(2)\prime}$.
\end{proof}

\begin{cor}\label{rem:unique-g2-critical-stratum}
From the above analysis, for fixed $\lambda=(\lambda^{(1)},\lambda^{(2)})$ and $\theta=(\theta^{(1)},\theta^{(2)})$, the dimension of
\[
m_{\theta^{(1)};\lambda^{(1)},\lambda^{(2)}}^{-1} \bigl( \CS_{G_2}^{(\theta^{(1)};\lambda^{(1)\prime},\lambda^{(2)\prime})} \cap\sS_{G_2}^{\theta^{(2)}} \bigr)
\]
is no more than
\[
\frac12d(\lambda-\theta) =\frac12d(\lambda-\lambda')+\frac12d(\lambda'-\theta).
\]
The equality holds for at most one pair,
\[
(\lambda^{(1)\prime},\lambda^{(2)\prime})=(\theta^{(1)},\lambda^{(2)}).
\]
In this case, we also have
\[
(\lambda^{(1)},\lambda^{(2)})\geq_G(\theta^{(1)},\lambda^{(2)})\geq_G(\theta^{(1)},\theta^{(2)}).
\]
\end{cor}
\begin{proof}
The equality holds only if both the intersection $\CS_{G_2}^{(\theta^{(1)};\lambda^{(1)\prime},\lambda^{(2)\prime})}\cap\sS_{G_2}^{\theta^{(2)}}$ and the fiber \eqref{eq:g2-projection-fiber} have critical dimension, which implies that $\lambda^{(1)\prime}=\theta^{(1)}$ (see Proposition~\ref{prop:g2-base-critical}) and $\lambda^{(2)\prime}=\lambda^{(2)}$ (see Proposition~\ref{prop:g2-fiber-dimension}).
\end{proof}

\subsubsection{The projection to \texorpdfstring{$\Gr_{G_1}$}{Gr_G1} and its strata}\label{sec:g1-projection-strata}

Now, we consider the projection of $\sO_{\Gr}^{\lambda}\cap\sS_{\Gr}^{\theta}$ to $\Gr_{G_1}$, whose fibers are intersections in $\Gr_{G_2}$. Similar to \eqref{eq:Y-convolution}, we consider the following equivariant product
\begin{equation}
\label{eq:Z-convolution}\fitmath{Z_{\theta^{(2)};\lambda^{(2)},\lambda^{(1)}} := \bigl(\sD t^{\theta^{(2)}} U_2(\bF)G_2(\bO)\bigr) \times^{G_2(\bO)} \bigl(G_2(\bO)H_u(\bF)t^{-\lambda^{(2)}}\sD^{-1}t^{\lambda^{(1)}} G_1(\bO)/G_1(\bO)\bigr).}
\end{equation}
Let
\begin{equation}\label{eq:g1-projection}
m'_{\theta^{(2)};\lambda^{(2)},\lambda^{(1)}} \colon Z_{\theta^{(2)};\lambda^{(2)},\lambda^{(1)}} \longrightarrow G_2(\bF)/G_1(\bO)
\end{equation}
be the map induced by the multiplication in $G_2(\bF)$. Then, we have
\begin{equation}\label{eq:intersection-as-g1-fiber-product}
\sO_{\Gr}^{(\lambda^{(1)},\lambda^{(2)})} \cap\sS_{\Gr}^{(\theta^{(1)},\theta^{(2)})} \simeq Z_{\theta^{(2)};\lambda^{(2)},\lambda^{(1)}} \times_{G_2(\bF)/G_1(\bO)}\sS_{G_1}^{-,\theta^{(1)}}.
\end{equation}
Under the isomorphism \eqref{eq:intersection-as-g1-fiber-product}, the projection $\sO_{\Gr}^{\lambda}\cap\sS_{\Gr}^{\theta}\to\sS_{G_1}^{-,\theta^{(1)}}$ is the base change of $m'_{\theta^{(2)};\lambda^{(2)},\lambda^{(1)}}$. The image of \eqref{eq:g1-projection} is
\[
\mathcal B_{\theta^{(2)};\lambda^{(2)},\lambda^{(1)}} := \sD t^{\theta^{(2)}} U_2(\bF)G_2(\bO)H_u(\bF) t^{-\lambda^{(2)}}\sD^{-1}t^{\lambda^{(1)}} G_1(\bO)/G_1(\bO).
\]
We stratify $\mathcal B_{\theta^{(2)};\lambda^{(2)},\lambda^{(1)}}\cap\Gr_{G_1}$ by its intersections with the following subsets of $G_2(\bF)/G_1(\bO)$:
\begin{equation}\label{eq:g1-image-stratum}
\CS_{G_1}^{(\theta^{(2)};\lambda^{(2)\prime},\lambda^{(1)\prime})} := \sD t^{\theta^{(2)}} U_2(\bF)T_2(\bO)H_u(\bF) t^{-\lambda^{(2)\prime}}\sD^{-1}t^{\lambda^{(1)\prime}}G_1(\bO)/G_1(\bO).
\end{equation}

\begin{lem}\label{lem:g1-semi-spherical-indices}
\begin{enumerate}
\item The above strata \eqref{eq:g1-image-stratum} meeting $\Gr_{G_1}$ are indexed by a subset of the pairs $\lambda'=(\lambda^{(1)\prime},\lambda^{(2)\prime})\in\Lambda$ satisfying the same conditions as in Lemma~\ref{lem:L-orbit-indices}, while we remove the requirement of the dominance condition for $\lambda^{(2)\prime}$. A representative of $\CS_{G_1}^{(\theta^{(2)};\lambda^{(2)\prime},\lambda^{(1)\prime})}$ is given by
\[
\sD t^{\theta^{(2)}}t^{-\lambda^{(2)\prime}}\sD^{-1}t^{\lambda^{(1)\prime}} G_1(\bO)/G_1(\bO).
\]
\item If one of the strata defined in \eqref{eq:g1-image-stratum} intersects $\Gr_{G_1}$ and lies in the boundary of $\CS_{G_1}^{(\theta^{(2)};\lambda^{(2)\prime},\lambda^{(1)\prime})}$, then its index is strictly smaller than $\lambda'$ for $\leq_{\super}$.
\end{enumerate}
\end{lem}
\begin{proof}
{(1).} Let $U_L=U_2\cap L$. Up to the left translation by $t^{-\theta^{(2)}}\sD^{-1}$, we need consider the orbits in $L(\bF)$ for the left action of $U_L(\bF)T_2(\bO)$ and the right action of $G_1(\bO)$, which contain a point of the form
$x=t^{-\theta^{(2)}}\sD^{-1}a$, $a\in G_1(\bF)$.

We can apply the same argument as in the proof of Lemma~\ref{lem:semi-spherical-orbits} here, where we use the open orbit $(B_2\cap L)\sD^{-1}B_1^-$ in $L$, see Lemma~\ref{lem:trivial-H-B-intersection}. Namely, for fixed $x,y$ of this form and sufficiently dominant $\eta\in\Lambda_{G_2}$ (depending on $x,y$), we have
\[
y\in U_L(\bF)T_2(\bO)xG_1(\bO) \quad\Longleftrightarrow\quad y\in t^\eta L(\bO)t^{-\eta}xG_1(\bO).
\]
Then, we translate by $t^{-\eta}$ and take the inverses, and according to Lemma~\ref{lem:L-orbit-indices}, we obtain the representatives $t^{-\lambda^{(1)\prime}}\sD t^{\lambda^{(2)\prime}+\eta}$ and the uniqueness of their indices, with $\lambda^{(1)\prime}$ dominant. Finally, we take the inverses and translate back, which proves (1).

{(2).} Again, we apply the comparison with the left $t^\eta L(\bO)t^{-\eta}$-action in (1), and then translate by $t^{-\eta}$ and take the inverses. If $\CS_{G_1}^{(\theta^{(2)};\nu^{(2)},\nu^{(1)})}$ lies in the boundary of $\CS_{G_1}^{(\theta^{(2)};\lambda^{(2)\prime},\lambda^{(1)\prime})}$, then by Lemma~\ref{lem:spherical-orbit-closure}, we have $(\nu^{(1)},\nu^{(2)}+\eta) <_{\super}(\lambda^{(1)\prime},\lambda^{(2)\prime}+\eta)$, i.e., $(\nu^{(1)},\nu^{(2)}) <_{\super}(\lambda^{(1)\prime},\lambda^{(2)\prime})$.
\end{proof}

\subsubsection{Critical strata for the projection to \texorpdfstring{$\Gr_{G_1}$}{Gr_G1}}

We consider the intersections of the strata \eqref{eq:g1-image-stratum} with $\mathcal B_{\theta^{(2)};\lambda^{(2)},\lambda^{(1)}}\cap\Gr_{G_1}$. Similar to \eqref{eq:g2-projection-fiber}, the fiber of the map \eqref{eq:g1-projection} over the representative
\[
\sD t^{\theta^{(2)}} t^{-\lambda^{(2)\prime}}\sD^{-1}t^{\lambda^{(1)\prime}} G_1(\bO)/G_1(\bO)
\]
is
\begin{equation}\label{eq:g1-projection-fiber}
\begin{split}
&\sD^{-1}t^{\lambda^{(1)\prime}}G_1(\bO)H_u(\bF)t^{-\lambda^{(1)}} \sD t^{\lambda^{(2)}} G_2(\bO)/G_2(\bO)\\
&\hspace{35mm}\cap U_2(\bF)t^{\lambda^{(2)\prime}}G_2(\bO)/G_2(\bO).
\end{split}
\end{equation}

As in the case of projecting to $\Gr_{G_2}$, we have the following statement for the preimage of $\CS_{G_1}^{(\theta^{(2)};\lambda^{(2)\prime},\lambda^{(1)\prime})}\cap\sS_{G_1}^{-,\theta^{(1)}}$ under \eqref{eq:g1-projection}. We estimate its dimension by adding the dimension bounds for this intersection and for the fiber \eqref{eq:g1-projection-fiber}.

\begin{prop}\label{prop:g1-strata-dimension}
Let $\lambda'=(\lambda^{(1)\prime},\lambda^{(2)\prime})$, as in Lemma~\ref{lem:g1-semi-spherical-indices}. If $\CS_{G_1}^{(\theta^{(2)};\lambda^{(2)\prime},\lambda^{(1)\prime})}$ intersects with $\mathcal B_{\theta^{(2)};\lambda^{(2)},\lambda^{(1)}}\cap\Gr_{G_1}$, then we have
\[
\dim\bigl( \CS_{G_1}^{(\theta^{(2)};\lambda^{(2)\prime},\lambda^{(1)\prime})} \cap\sS_{G_1}^{-,\theta^{(1)}} \bigr) \leq \frac12d(\lambda'-\theta).
\]
The fiber \eqref{eq:g1-projection-fiber} has dimension at most
\[
\frac12d(\lambda-\lambda').
\]
By adding these two estimates, we obtain the following bound for the preimage:
\[
\dim (m'_{\theta^{(2)};\lambda^{(2)},\lambda^{(1)}})^{-1} \bigl( \CS_{G_1}^{(\theta^{(2)};\lambda^{(2)\prime},\lambda^{(1)\prime})} \cap\sS_{G_1}^{-,\theta^{(1)}} \bigr) \leq \frac12d(\lambda-\theta).
\]
Furthermore, the equality in the last bound holds only if
\[
(\lambda^{(1)\prime},\lambda^{(2)\prime})=(\lambda^{(1)},\theta^{(2)}).
\]
In this case, we have
\[
\lambda^{(2)}-\theta^{(2)}\in\Lambda_{G_2}^{\on{pos}}.
\]
\end{prop}

\begin{proof}
By Lemma~\ref{lem:g1-semi-spherical-indices}, $\lambda^{(1)\prime}$ is dominant. To prove the first inequality, we translate $\lambda^{(2)\prime}$ and $\theta^{(2)}$ by the same sufficiently dominant coweight $\eta\in X_*(T_2)$, so that $\lambda^{(2)\prime}+\eta$ is dominant. Since $t^\eta$ normalizes $U_2(\bF)T_2(\bO)H_u(\bF)$, it changes neither $\CS_{G_1}^{(\theta^{(2)};\lambda^{(2)\prime},\lambda^{(1)\prime})}\cap\sS_{G_1}^{-,\theta^{(1)}}$ nor $\frac12d(\lambda'-\theta)$. By the inclusion $T_2(\bO)\subset G_2(\bO)$, this intersection is contained in the image of
\[
\sO_{\Gr}^{\lambda'+(0,\eta)}\cap\sS_{\Gr}^{\theta+(0,\eta)}
\]
under the projection to $\Gr_{G_1}$, and thus, by Corollary~\ref{lem:dimension-estimate-super}, its dimension is bounded by $\frac12d(\lambda'-\theta)$.

Similarly, the fiber \eqref{eq:g1-projection-fiber} is contained in $\sO_{\Gr}^{\lambda}\cap\sS_{\Gr}^{\lambda'}$ (with the first coordinate $t^{\lambda^{(1)\prime}}$), thus, by the same corollary, its dimension is bounded by $\frac12d(\lambda-\lambda')$. Note that the sum of the two bounds is $\frac12d(\lambda-\theta)$.

Now, assume that the equality holds in the bound for the preimage under \eqref{eq:g1-projection}. Then, both $\CS_{G_1}^{(\theta^{(2)};\lambda^{(2)\prime},\lambda^{(1)\prime})}\cap\sS_{G_1}^{-,\theta^{(1)}}$ and the fiber \eqref{eq:g1-projection-fiber} have critical dimension. For the base intersection, we apply the same translation as above. By Lemma~\ref{lem:g1-semi-spherical-indices} and Corollary~\ref{lem:dimension-estimate-super}, with the same argument as in Corollary~\ref{cor:g2-boundary-codimension}, the intersection of the boundaries of $\CS_{G_1}^{(\theta^{(2)};\lambda^{(2)\prime},\lambda^{(1)\prime})}$ and $\sS_{G_1}^{-,\theta^{(1)}}$ has dimension at most $\frac12d(\lambda'-\theta)-2$. Then, the hyperplane arguments in Propositions~\ref{prop:g2-base-critical} and~\ref{prop:g2-fiber-dimension} apply, with $G_1$ and $G_2$ interchanged. The equality for the base implies that $\lambda^{(2)\prime}=\theta^{(2)}$, and the equality for the fiber implies that $\lambda^{(1)\prime}=\lambda^{(1)}$ and $\lambda^{(2)}-\lambda^{(2)\prime}\in\Lambda_{G_2}^{\on{pos}}$, which are the required conditions.
\end{proof}

Let us summarize the results of this subsection. By Corollary~\ref{rem:unique-g2-critical-stratum}, the preimage of $\CS_{G_2}^{(\theta^{(1)};\lambda^{(1)\prime},\lambda^{(2)\prime})}\cap\sS_{G_2}^{\theta^{(2)}}$ under \eqref{eq:g2-projection} has critical dimension only if $(\lambda^{(1)\prime},\lambda^{(2)\prime})=(\theta^{(1)},\lambda^{(2)})$, and by Proposition~\ref{prop:g1-strata-dimension}, the preimage of $\CS_{G_1}^{(\theta^{(2)};\lambda^{(2)\prime},\lambda^{(1)\prime})}\cap\sS_{G_1}^{-,\theta^{(1)}}$ under \eqref{eq:g1-projection} has critical dimension only if $(\lambda^{(1)\prime},\lambda^{(2)\prime})=(\lambda^{(1)},\theta^{(2)})$. In other words, for each of the two projections, a dense open subset of any critical-dimensional irreducible component of $\sO_{\Gr}^{\lambda}\cap\sS_{\Gr}^{\theta}$ maps into the single stratum given above. We will use both conclusions later, e.g., Lemma~\ref{lem:critical-G2-fiber}.

\subsection{The \texorpdfstring{$G_2$}{G2}-direction: proof that \texorpdfstring{$\lambda^{(2)}=\theta^{(2)}$}{lambda2=theta2}}\label{sec:g2-direction}

\subsubsection{The critical intersection and its boundary}

In this subsection, we prove Proposition~\ref{prop:part-q}. The idea is to reduce it to a statement about an intersection in $\Gr_{G_2}$, which is contained in the intersection of two opposite semi-infinite orbits of $G_2$. Namely, we first use the projection of $\sO_{\Gr}^{\lambda}\cap\sS_{\Gr}^{\theta}$ to $\Gr_{G_1}$ to find a critical-dimensional component of the fiber \eqref{eq:g1-projection-fiber} for $\lambda'=(\lambda^{(1)},\theta^{(2)})$. Then, we compare it with the following intersection in $\Gr_{G_2}$, to which we will apply the hyperplane argument.

We define
\[
K_{\lambda^{(1)}}^{(2)} := \sD^{-1}t^{\lambda^{(1)}} U_1^-(\bF)T_1(\bO)H_u(\bF)t^{-\lambda^{(1)}}\sD \cap G_2(\bF),
\]
and we consider the following intersection in $\Gr_{G_2}$,
\begin{equation}\label{eq:critical-G2-fiber}
\mathfrak F_2(\lambda^{(1)},\lambda^{(2)};\theta^{(2)}) := K_{\lambda^{(1)}}^{(2)}t^{\lambda^{(2)}} G_2(\bO)/G_2(\bO) \cap U_2(\bF)t^{\theta^{(2)}} G_2(\bO)/G_2(\bO).
\end{equation}
Note that $K_{\lambda^{(1)}}^{(2)}$ does not depend on $\lambda^{(1)}$, but we write $\lambda^{(1)}$ to emphasize that the equivariant trivialization below depends on $\lambda^{(1)}$.

By the explicit matrices in Section~\ref{trans D}, we have
\begin{equation}\label{eq:g2-orbit-in-semi-infinite}
\begin{aligned}
K_{\lambda^{(1)}}^{(2)}&\subset U_2^-(\bF)T_2(\bO),\\
\CT_{G_2}^{\lambda} &:= \CS_{G_2}^{(\lambda^{(1)};\lambda^{(1)},\lambda^{(2)})}\\
&=K_{\lambda^{(1)}}^{(2)}t^{\lambda^{(2)}}G_2(\bO)/G_2(\bO) \subset\sS_{G_2}^{-,\lambda^{(2)}}.
\end{aligned}
\end{equation}
The equivariant trivialization on this $K_{\lambda^{(1)}}^{(2)}$-orbit is induced from the original $H(\bF)$-orbit via the map
\begin{equation}\label{eq:g2-orbit-lift}
\begin{aligned}
r_2:\CT_{G_2}^{\lambda} &\longrightarrow\sO_{\Gr}^{\lambda},\\
x_2=\sD^{-1}g_1u\sD t^{\lambda^{(2)}}G_2(\bO) &\longmapsto(g_1t^{\lambda^{(1)}}G_1(\bO),x_2),
\end{aligned}
\end{equation}
where $g_1\in U_1^-(\bF)T_1(\bO)$ and $u\in H_u(\bF)$. Note that $r_2$ is well-defined. Indeed, since $\lambda^{(2)}$ is dominant, the $G_1$-component of the stabilizer of $t^{\lambda^{(2)}}$ in $K_{\lambda^{(1)}}^{(2)}$ lies in $U_1^-(\bO)T_1(\bO)$, which fixes $t^{\lambda^{(1)}}$.

We denote by $\tau_2^\lambda$ the induced $\sD^{-1}U_1^-(\bF)H_u(\bF)\sD$-equivariant trivialization. Namely, if we denote by $\pr_1$ the first projection, then it is given by
\begin{equation}\label{eq:g2-orbit-trivialization}
\tau_2^{\lambda} :=r_2^*\tau_H^\lambda\otimes \bigl((\pr_1\circ r_2)^*\tau_{1,-}^{\lambda^{(1)}}\bigr)^{-1}.
\end{equation}
On $\mathfrak F_2(\lambda^{(1)},\lambda^{(2)};\theta^{(2)})$, its ratio with $\tau_{2,+}^{\theta^{(2)}}$ is $r_2^*\cL_{\lambda,(\lambda^{(1)},\theta^{(2)})}$. With some abuse of notation, we use the same notation for this pullback, and denote the pullback $r_2^*\psi$ by $\psi$.

\begin{lem}\label{lem:critical-G2-fiber}
Under the assumptions of Theorem~\ref{thm:local-system-rigidity}, there is a critical-dimensional irreducible component of \eqref{eq:critical-G2-fiber}, such that the tensor product local system $\cL_{\lambda,(\lambda^{(1)},\theta^{(2)})}\otimes\psi$ is constant on a dense open subset of it. Here, the two trivializations defining this ratio local system $\cL_{\lambda,(\lambda^{(1)},\theta^{(2)})}$ are induced directly from the original $H(\bF)_c^\omega$-equivariant trivialization on $\sO_{\Gr}^{(\lambda^{(1)},\lambda^{(2)})}$ and the $U_2(\bF)$-equivariant trivialization on $\sS_{G_2}^{\theta^{(2)}}$.
\end{lem}

\begin{proof}
First, we find a critical-dimensional component of a fiber over $\Gr_{G_1}$. Let $C$ be a critical-dimensional component of $\sO_{\Gr}^{\lambda}\cap\sS_{\Gr}^{\theta}$ in Theorem~\ref{thm:local-system-rigidity}. By Proposition~\ref{prop:g1-strata-dimension}, a dense open subset of $C$ maps to $\CS_{G_1}^{(\theta^{(2)};\theta^{(2)},\lambda^{(1)})}\cap\sS_{G_1}^{-,\theta^{(1)}}$. Thus, the fiber of $C$ over a general point of its image has a component of dimension $\frac12d((0,\lambda^{(2)}-\theta^{(2)}))$, such that the restriction of the tensor product local system to a dense open subset of it is constant.

Secondly, we identify the induced local system on this component. Using the equivariant product \eqref{eq:Z-convolution}, we identify the component with a component of the fiber of $m'_{\theta^{(2)};\lambda^{(2)},\lambda^{(1)}}$ over $t^{\lambda^{(1)}}$. Since the first coordinate is fixed, the trivialization on $\sS_{G_1}^{-,\theta^{(1)}}$ only contributes a constant line, and thus, the induced local system is $\cL_{\lambda,(\lambda^{(1)},\theta^{(2)})}\otimes\psi$, up to a constant line.

Finally, we project to $\Gr_{G_2}$. Under the map $x_2\mapsto(t^{\lambda^{(1)}},x_2)$, the component above is a critical-dimensional component of $\sO_{\Gr}^{\lambda}\cap\sS_{\Gr}^{(\lambda^{(1)},\theta^{(2)})}$. According to Corollary~\ref{rem:unique-g2-critical-stratum}, this projection identifies a dense open subset of it with a dense open subset of a critical-dimensional component of \eqref{eq:critical-G2-fiber}. Since the two trivializations are the ones in the statement, the assertion follows.
\end{proof}

\begin{cor}\label{cor:critical-G2-fiber-boundary}
Let $C_2$ be a critical-dimensional irreducible component of \eqref{eq:critical-G2-fiber}, so that
\[
\dim C_2=\frac12d\bigl((0,\lambda^{(2)}-\theta^{(2)})\bigr).
\]
Let
\[
A=K_{\lambda^{(1)}}^{(2)}t^{\lambda^{(2)}}G_2(\bO)/G_2(\bO), \qquad B=U_2(\bF)t^{\theta^{(2)}}G_2(\bO)/G_2(\bO).
\]
Then, if we take all the closures in $\Gr_{G_2}$, we have
\[
\dim\bigl(\overline{C_2}\cap(\overline A\setminus A) \cap(\overline B\setminus B)\bigr) \leq\dim C_2-2.
\]
\end{cor}

\begin{proof}
This follows from Corollary~\ref{cor:g2-boundary-codimension}, with $\theta^{(1)}=\lambda^{(1)}$ and $\lambda'=\lambda$.
\end{proof}

From now on, let $C$ be the critical-dimensional component of \eqref{eq:critical-G2-fiber} given by Lemma~\ref{lem:critical-G2-fiber}. Note that it is dense in an irreducible component of $\sS_{G_2}^{-,\lambda^{(2)}}\cap\sS_{G_2}^{+,\theta^{(2)}}$ (their dimensions are the same).

\subsubsection{Reduction to rank one}

For $a=1,2$, let $\Lambda_{G_a}^{\sharp}:=\{\nu\in\Lambda_{G_a}\mid b(\nu,\Lambda_{G_a})=1\}$ be the metaplectic lattice of the pullback of $\cG^G$ to $\Gr_{G_a,\Ran}^{\omega}$.

In this subsection, we first show that $\lambda^{(2)}-\theta^{(2)}$ belongs to the metaplectic lattice of $G_2$, and then reduce the proof of Proposition~\ref{prop:part-q} to the rank-one case for a mixed simple coroot. The idea is as follows. If $\cL_{\lambda,(\lambda^{(1)},\theta^{(2)})}$ were the ratio of two semi-infinite trivializations, then Gaitsgory's argument in \cite[\S 7.4]{[G7]} would directly reduce the problem to rank one. In our case, one of the two trivializations is induced from the $H(\bF)$-orbit $\sO_{\Gr}^{\lambda}$, and thus, we compare both of them with the trivialization on $\sS_{G_2}^{-,\lambda^{(2)}}$, see \eqref{eq:g2-ratio-decomposition} below. The important observation is that the monodromy of $\cL_2^{(2)}$ along the boundary of $\sS_{G_2}^{+,\theta^{(2)}}$ cannot be cancelled by $\cL_1^{(2)}\otimes\psi$, see the proof of Lemma~\ref{lem:rank-one-fixed-orbit-comparison}. Hence, Gaitsgory's argument applies to $\cL_2^{(2)}$.

More precisely, we use the following three trivializations of $\cG^G$:
\begin{itemize}
\item $\tau_2^{\lambda}$ on $\CT_{G_2}^{\lambda}$, which is $\sD^{-1}U_1^-(\bF)H_u(\bF)\sD$-equivariant, see \eqref{eq:g2-orbit-trivialization};
\item $\tau_{2,-}^{\lambda^{(2)}}$ on $\sS_{G_2}^{-,\lambda^{(2)}}\supset\CT_{G_2}^{\lambda}$, which is $U_2^-(\bF)$-equivariant;
\item $\tau_{2,+}^{\theta^{(2)}}$ on $\sS_{G_2}^{+,\theta^{(2)}}$, which is $U_2(\bF)$-equivariant.
\end{itemize}
For two trivializations $\tau$ and $\tau'$ on the same space, we denote by $\tau\otimes\tau'^{-1}$ the local system corresponding to their ratio. In particular, $\cL_{\lambda,(\lambda^{(1)},\theta^{(2)})}=\tau_2^{\lambda}\otimes(\tau_{2,+}^{\theta^{(2)}})^{-1}$ on $\mathfrak F_2(\lambda^{(1)},\lambda^{(2)};\theta^{(2)})$. We define
\[
\cL_1^{(2)}:=\tau_2^{\lambda}\otimes\bigl(\tau_{2,-}^{\lambda^{(2)}}\bigr)^{-1}\ \text{ on }\CT_{G_2}^{\lambda},\qquad
\cL_2^{(2)}:=\tau_{2,+}^{\theta^{(2)}}\otimes\bigl(\tau_{2,-}^{\lambda^{(2)}}\bigr)^{-1}\ \text{ on }\sS_{G_2}^{-,\lambda^{(2)}}\cap\sS_{G_2}^{+,\theta^{(2)}}.
\]
Then, on $\mathfrak F_2(\lambda^{(1)},\lambda^{(2)};\theta^{(2)})$, we have
\begin{equation}\label{eq:g2-ratio-decomposition}
\cL_{\lambda,(\lambda^{(1)},\theta^{(2)})} \simeq\cL_1^{(2)}\otimes\bigl(\cL_2^{(2)}\bigr)^{-1}.
\end{equation}
Note that $\cL_1^{(2)}$ and $\psi$ are defined on the whole orbit $\CT_{G_2}^{\lambda}$, and they do not depend on $\theta^{(2)}$, while $\cL_2^{(2)}$ only depends on the two semi-infinite orbits of $G_2$.

Roughly speaking, by the following lemma, we can replace $\theta^{(2)}$ by $\theta^{(2)}+q_{2,i}\alpha_{2,i}$, such that there is still a critical-dimensional irreducible component on which the tensor product local system is constant, and by continuing this process, we reduce the problem to the rank-one case.

\begin{lem}\label{lem:rank-one-fixed-orbit-comparison}
Let $\theta'=(\lambda^{(1)},\theta^{(2)})$, and let $Z$ be a critical-dimensional irreducible component of the fiber of $\sO_{\Gr}^{\lambda}\cap\sS_{\Gr}^{\theta'}\to\Gr_{G_1}$ over $t^{\lambda^{(1)}}$, as in Lemma~\ref{lem:critical-G2-fiber}. Assume that $\cL_{\lambda,\theta'}\otimes\psi$ is constant on a dense open subset of $Z$. Let $\alpha_{2,i}$ be a simple geometric coroot, and let $m$ be the maximal number of successive applications of the crystal operator corresponding to $\alpha_{2,i}$ in \cite[\S\S 7.4.5--7.4.6]{[G7]} to the corresponding component of $\sS_{G_2}^{-,\lambda^{(2)}}\cap\sS_{G_2}^{+,\theta^{(2)}}$. Here, this operator increases the index of the positive semi-infinite orbit by $\alpha_{2,i}$. Then, $q(\alpha_{2,i})^m=1$. Moreover, if $q_i:=q_{2,i}=\ord(q(\alpha_{2,i}))<\infty$ and $m\geq q_i$, then the component corresponding to $q_i$ successive applications
\[
Z'\subset\sO_{\Gr}^{\lambda}\cap \sS_{\Gr}^{\theta'+q_i\alpha_{2,i}}
\]
has dimension $\dim Z-q_i$, and $\cL_{\lambda,\theta'+q_i\alpha_{2,i}}\otimes\psi$ is constant on a dense open subset of $Z'$.
\end{lem}
\begin{proof}
If $m=0$, there is nothing to prove. As in the proof of Lemma~\ref{lem:critical-G2-fiber}, a dense open subset of $Z$ is identified with a dense open subset of a critical-dimensional component $C'$ of $\mathfrak F_2(\lambda^{(1)},\lambda^{(2)};\theta^{(2)})$, which is dense in a component $Z_2$ of $\sS_{G_2}^{-,\lambda^{(2)}}\cap\sS_{G_2}^{+,\theta^{(2)}}$.

First, we describe the successive components by using the rank-one Levi subgroup of $G_2$ corresponding to $\alpha_{2,i}$. Let $Z_0$ be a smooth dense open subset of the component obtained after $m$ applications. By the product decomposition \cite[(7.5)]{[G7]} and the rank-one description in \cite[\S\S 7.3.2, 7.4.5]{[G7]}, we can choose coordinates $u_1,\ldots,u_m$ on $\BA^m$, such that dense open subsets of the successive components are identified with
\[
\{u_1=\cdots=u_r=0,\ u_{r+1}\neq0\}\times Z_0\quad(0\leq r<m),\qquad \{0\}\times Z_0\quad(r=m),
\]
inside the same product $\BA^m\times Z_0$ in $\sS_{G_2}^{-,\lambda^{(2)}}$. We call them the $r$-th strata, and we denote by $\cL_{2,r}^{(2)}$ the local system $\cL_2^{(2)}$ for $\theta^{(2)}+r\alpha_{2,i}$ on them, and by $\cL_0:=\cL_{2,m}^{(2)}$ the one on $Z_0$. Then, by \cite[Proposition 7.4.8 and \S 7.5.3]{[G7]}, we have
\[
\cL_{2,r}^{(2)}\simeq u_{r+1}^*\cK_r\otimes\pr_{Z_0}^*\cL_0, \qquad 0\leq r<m,
\]
where $\pr_{Z_0}:\BA^m\times Z_0\to Z_0$ is the projection, and $\cK_r$ is a Kummer local system on $\BG_m$ with monodromy $q(\alpha_{2,i})^{m-r}$ up to inversion. In particular, $\cK_r$ is trivial if and only if $q(\alpha_{2,i})^{m-r}=1$.

Now, by Corollary~\ref{cor:critical-G2-fiber-boundary}, the divisor $u_1=0$, which lies in the boundary of $\sS_{G_2}^{+,\theta^{(2)}}$, is not contained in the boundary of $\CT_{G_2}^{\lambda}$. Thus, it meets $\CT_{G_2}^{\lambda}$ in a dense open subset, where $\cL_1^{(2)}\otimes\psi$ is lisse. Since
\[
\cL_{\lambda,\theta'}\otimes\psi\simeq\bigl(\cL_1^{(2)}\otimes\psi\bigr)\otimes\bigl(\cL_{2,0}^{(2)}\bigr)^{-1}
\]
is constant on a dense open subset of $C'$, the monodromy of $\cL_{2,0}^{(2)}$ around $u_1=0$ cannot be cancelled by $\cL_1^{(2)}\otimes\psi$, and hence, it is trivial, i.e., $q(\alpha_{2,i})^m=1$.

Finally, assume that $m\geq q_i$. By the same boundary estimate at each step, the $q_i$-th stratum still meets $\CT_{G_2}^{\lambda}$ in a dense open subset, and as in the proof of Lemma~\ref{lem:critical-G2-fiber}, it gives the component $Z'$ in the statement, with the first coordinate $t^{\lambda^{(1)}}$ and dimension $\dim Z-q_i$. On the $0$-th stratum, since $q(\alpha_{2,i})^m=1$, the above isomorphism identifies $\cL_1^{(2)}\otimes\psi$ with $\pr_{Z_0}^*\cL_0$, up to a constant line. Both of them are lisse on the preimage of $\CT_{G_2}^{\lambda}$ in $\BA^m\times Z_0$, which is open, and hence smooth, so this identification extends there. Note that, for $0<r<q_i$, the Kummer factor is non-trivial, since $q(\alpha_{2,i})^{m-r}=q(\alpha_{2,i})^{-r}\neq1$, and thus, the tensor product local system is not constant on the intermediate strata. This is why we apply the crystal operator $q_i$ times. For $r=q_i$, the Kummer factor is trivial again, since $q(\alpha_{2,i})^{m-q_i}=1$ (if $q_i=m$, we have $\cL_{2,q_i}^{(2)}=\cL_0$ by definition). Therefore, by \eqref{eq:g2-ratio-decomposition} for $\theta^{(2)}+q_i\alpha_{2,i}$, the local system $\cL_{\lambda,\theta'+q_i\alpha_{2,i}}\otimes\psi\simeq(\cL_1^{(2)}\otimes\psi)\otimes(\cL_{2,q_i}^{(2)})^{-1}$ is constant on a dense open subset of $Z'$.
\end{proof}

By applying the above argument successively along the sequence in \cite[\S 7.4.6]{[G7]}, we also obtain the following corollary.

\begin{cor}\label{prop:g2-sharp-difference}
Under the assumptions of Theorem~\ref{thm:local-system-rigidity}, the local system $\cL_2^{(2)}$ is constant on a dense open subset of the component $C$ of \eqref{eq:critical-G2-fiber} constructed in Lemma~\ref{lem:critical-G2-fiber}. In particular,
\[
\lambda^{(2)}-\theta^{(2)}\in\Lambda_{G_2}^{\sharp}.
\]
\end{cor}

By Lemma~\ref{lem:rank-one-fixed-orbit-comparison}, Corollary~\ref{prop:g2-sharp-difference} and \cite[\S 7.4.6]{[G7]}, to prove Proposition~\ref{prop:part-q}, it is enough to consider the rank-one case
\begin{equation}\label{eq:g2-rank-one-difference}
\lambda^{(2)}-\theta^{(2)}=q_i\alpha_{2,i}, \qquad q_i=\ord(q(\alpha_{2,i}))<\infty.
\end{equation}
Here, both $\cL_{\lambda,(\lambda^{(1)},\theta^{(2)})}\otimes\psi$ and $\cL_2^{(2)}$ are constant on a dense open subset of a critical-dimensional component of \eqref{eq:critical-G2-fiber}. We should prove that this case cannot happen. If $\alpha_{2,i}$ is pure, then $\mathfrak F_2(\lambda^{(1)},\lambda^{(2)};\theta^{(2)})$ is just the usual intersection of opposite semi-infinite orbits, and it is known that the Whittaker local system on it is non-trivial, see \cite[Section 6]{[Lys]} or \cite[Theorem 9.4.9]{[GL]}, which contradicts Lemma~\ref{lem:kummer-as-separation}. Thus, it is enough to study the case when $\alpha_{2,i}$ is mixed, which we treat in the next subsection by explicit calculations.

\subsubsection{The mixed rank-one calculation}\label{10.3.3}

\begin{proof}[Completion of the proof of Proposition~\ref{prop:part-q}]
It remains to consider the mixed case in \eqref{eq:g2-rank-one-difference}. We denote $m=q_i$, and the intersection $\mathfrak F_2(\lambda^{(1)},\lambda^{(2)};\theta^{(2)})$ of \eqref{eq:critical-G2-fiber} can be rewritten as
\begin{equation}\label{eq:g2-mixed-intersection}
\sD^{-1}t^{\lambda^{(1)}}U_1^-(\bF)T_1(\bO)H_u(\bF) t^{-\lambda^{(1)}}\sD t^{\lambda^{(2)}}G_2(\bO)/G_2(\bO) \cap\sS_{G_2}^{\theta^{(2)}}.
\end{equation}
First, note that \eqref{eq:g2-mixed-intersection} is contained in $\sS_{G_2}^{-,\lambda^{(2)}}\cap\sS_{G_2}^{\theta^{(2)}} \simeq\BG_m\times\BA^{m-1}$, and our strategy is as follows. If the upper inequality in \eqref{eq:restricted-mixed-inequality} is strict, then we want to prove, case by case, that \eqref{eq:g2-mixed-intersection} is empty. If the equality holds in the upper bound, then we compute \eqref{eq:g2-mixed-intersection} explicitly, and we should prove that if the local system is constant, then $b(\lambda,\alpha_{2,i}^{\pm})=1$, which contradicts \eqref{eq:restricted-wall-monodromy}. Since the conjugation by $t^{\lambda^{(1)}}$ preserves $U_1^-(\bF)T_1(\bO)$, we omit it in the expressions below.

\case{The $\GL(M|N)$ case.}\par\nopagebreak[4]

For any point in the intersection \eqref{eq:g2-mixed-intersection}, we can find a representative of the form
\[
A=
\begin{pmatrix}
t^{\lambda^{(2)}_1}&&&&&&&\\
&\ddots&&&&&&\\
&&t^{\lambda^{(2)}_i}&&&&&\\
&&t^{\lambda^{(2)}_{i+1}-m}x&t^{\lambda^{(2)}_{i+1}}&&&&\\
&&&&\ddots&&&\\
&&&&&t^{\lambda^{(2)}_{M+1}}&&\\
&&&&&&\ddots&\\
&&&&&&&t^{\lambda^{(2)}_N}
\end{pmatrix},\qquad x\in(\bO/t^m\bO)^\times.
\]
If the point represented by $A$ lies in the intersection \eqref{eq:g2-mixed-intersection}, then there exist $u\in U_1^-(\bF)$, $a\in T_1(\bO)$ and $u'\in U_{M,N}^-(\bF)$ such that
\begin{equation}\label{eq:g2-gl-orbit-condition}
\sD^{-1}uau'\sD A\in t^{\lambda^{(2)}}\GL_N(\bO).
\end{equation}
Note that the left multiplication by $\sD$ only changes the $(M+1)$-st row, thus, for $1\leq j\leq M$, the $(j,i)$-entry of $uau'\sD A$ belongs to $t^{\lambda^{(2)}_j}\bO$. Then, by writing down the expression of $\sD^{-1}uau'\sD A$, we see that its $(M+1,i)$-entry belongs to
\[
t^{\lambda^{(2)}_i}\bO +t^{\lambda^{(2)}_{i+1}-m}\bO^\times +t^{\lambda^{(2)}_{i+1}}\bO+\cdots+t^{\lambda^{(2)}_M}\bO.
\]
If $0\leq\lambda^{(2)}_{i+1}-\lambda^{(2)}_i<m$, then the second term cannot be cancelled, and thus, this entry has valuation $\lambda^{(2)}_{i+1}-m<\lambda^{(2)}_{M+1}$, which contradicts \eqref{eq:g2-gl-orbit-condition}.

If $\lambda^{(2)}_{i+1}-\lambda^{(2)}_i=m$, then the intersection \eqref{eq:g2-mixed-intersection} is isomorphic to
\begin{equation}\label{eq:g2-gl-wall-cycle}
(\BG_m\setminus\{-1\})\times\BA^{m-1}.
\end{equation}
Indeed, for $i<M$, we consider the matrix
\[
a(s)=\begin{pmatrix}s&0\\1-s&1\end{pmatrix},\qquad s\in\bO^\times,
\]
in the positions $i,i+1$ of the first factor. Since its column sums are one, the conjugation by $\sD$ does not introduce any entry in the $(M+1)$-st row. For $i=M$, we take $a(s)$ to be the diagonal matrix with entry $s$ in position $M$ and $1$ elsewhere, and then the new $(M+1,M)$-entry is $1-s$. In both cases, we multiply $\sD^{-1}a(s)\sD t^{\lambda^{(2)}}$ on the right by the diagonal matrix with $s^{-1}$ in position $i$ and $1$ elsewhere, and we obtain $x=s^{-1}-1$, with $x(0)\neq0,-1$. The diagonal cocharacter is
\[
\nu=(e_i^{(1)},e_i^{(2)})=-\alpha_{2i-1}
\]
in \eqref{eq:typeA-duality-coordinates}, and thus, by pulling back $b_\lambda$ along $\nu$, the local monodromy around $x=-1$ is $b(\lambda,\nu)$ up to inversion. Since $\alpha_{2,i}=\alpha_{2i-1}+\alpha_{2i}$ and $b(\lambda,\alpha_{2,i})=1$, if the local system is constant, then we have $b(\lambda,\alpha_{2i-1})=b(\lambda,\alpha_{2i})=1$.

\case{The $\SOSp(2k|2l)$ case.}

\smallskip
\noindent\emph{(a) $M=2k$.} In this part, we write
\[
t^{\lambda^{(2)}}= \diag(t^{\lambda^{(2)}_1},\ldots,t^{\lambda^{(2)}_l},1, t^{-\lambda^{(2)}_l},\ldots,t^{-\lambda^{(2)}_1}),
\]
and use the analogous order in $G_1$. We denote $a_j=\lambda^{(2)}_{l-k+j}$, so that $a_1\leq\cdots\leq a_k\leq0$.

\smallskip
\noindent\emph{The short coroot case.} For any point in \eqref{eq:g2-mixed-intersection}, we can find a representative $A$, which is lower triangular and whose submatrix in positions $l-k+1,\ldots,l+1$ is
\[
A'=\begin{pmatrix} t^{a_1}&&&&&\\
&\ddots&&&&\\
&&t^{a_j}&&&\\
&&t^{a_{j+1}-m}x&t^{a_{j+1}}&&\\
&&&&\ddots&\\
&&&&&1
\end{pmatrix},\qquad x\in(\bO/t^m)^\times.
\]
Suppose that $\sD^{-1}uau'\sD A\in t^{\lambda^{(2)}}\SO(N)(\bO)$. As in type A, by writing down the expression of this matrix, we see that its entry in the middle row belongs to
\[
t^{a_j}\bO+t^{a_{j+1}-m}\bO^\times +t^{a_{j+1}}\bO+\cdots+t^{a_k}\bO.
\]
If $0\leq a_{j+1}-a_j<m$, then the second term cannot be cancelled, and this entry has valuation $a_{j+1}-m<0$, while it must belong to $\bO$, a contradiction. If $a_{j+1}-a_j=m$, then by the same matrix calculation as in type A, with the inverse transpose on the paired isotropic subspace, we obtain $x=s^{-1}-1$, and the intersection is isomorphic to $(\BG_m\setminus\{-1\})\times\BA^{m-1}$. Its diagonal cocharacter is
\[
\nu=(e_j^{(1)},e_{l-k+j}^{(2)}),
\]
which is the negative of one of the two odd summands of the mixed coroot. Thus, as in type A, if the local system is constant, then the characters corresponding to the two odd simple super roots are trivial.

\smallskip
\noindent\emph{The long coroot case.} Recall $a_k=\lambda^{(2)}_l\leq0$. For any point in \eqref{eq:g2-mixed-intersection}, we can find a representative $A$, which is lower triangular and whose submatrix in the middle three positions is
\[
A'=\begin{pmatrix} t^{a_k}&0&0\\
t^{-m}x&1&0\\
-\frac12t^{-a_k-2m}x^2&-t^{-a_k-m}x&t^{-a_k}
\end{pmatrix},\qquad x\in(\bO/t^m)^\times.
\]
In these positions, we have $\sD=u(1)$, where
\[
u(z)=\begin{pmatrix}1&0&0\\z&1&0\\-z^2/2&-z&1\end{pmatrix}.
\]
Since
\[
u(-1)\diag(c,1,c^{-1})u(1)u(z) =u\bigl(c^{-1}(1+z)-1\bigr)\diag(c,1,c^{-1}),
\]
the same submatrix of $\sD^{-1}uau'\sD A$ is
\[
\begin{pmatrix}c t^{a_k}&0&0\\d&1&0\\ *&*&c^{-1}t^{-a_k}\end{pmatrix},
\qquad d=t^{-m}x+t^{a_k}(1-c),\quad c\in\bO^\times.
\]
If $0\leq-a_k<m$, then $d$ has valuation $-m$, which contradicts $d\in\bO$. If $-a_k=m$, then the integrality condition is $c=1+x$ modulo $t^m$, and we obtain $(\BG_m\setminus\{-1\})\times\BA^{m-1}$. The diagonal cocharacter is
\[
\nu=(e_k^{(1)},e_l^{(2)}).
\]
Note that its negative is one of the odd super roots $(e_k^{(1)},-e_l^{(2)})$ and $(-e_k^{(1)},-e_l^{(2)})$, whose sum is $(0,-2e_l^{(2)})$. The local monodromy is $b(\lambda,\nu)$ up to inversion, thus, by \eqref{eq:restricted-wall-monodromy}, the local system cannot be constant.

\smallskip
\noindent\emph{(b) $M=2l+1$.} The preceding matrix calculation applies away from the end roots, using the sum of the two middle rows. We denote $a=\lambda^{(2)}_{k-1}$ and $b_0=\lambda^{(2)}_k$, so that $a\leq-|b_0|$.

For the coroot $-e_{k-1}^{(2)}+e_k^{(2)}$ and any point in \eqref{eq:g2-mixed-intersection}, we can find a representative $A$, which is lower triangular and whose submatrix in rows and columns $k-1,k,k+1,k+2$ is
\[
A'=\begin{pmatrix}
t^a&0&0&0\\
t^{b_0-m}x&t^{b_0}&0&0\\
0&0&t^{-b_0}&0\\
0&0&-t^{-a-m}x&t^{-a}
\end{pmatrix},\qquad x\in(\bO/t^m)^\times.
\]
If the point represented by $A$ lies in \eqref{eq:g2-mixed-intersection}, then there exists $g\in U_1^-(\bF)T_1(\bO)H_u(\bF)$ such that
\[
\sD^{-1}g\sD A\in t^{\lambda^{(2)}}\SO(2k)(\bO).
\]
By the explicit matrices in Sections~\ref{pol des} and \ref{trans D}, the same submatrix of $\sD^{-1}g\sD$ has the form
\begin{equation}\label{eq:SO-even-end-root-parameters}
g'=\begin{pmatrix}
s&0&0&0\\
sy_1&1&0&0\\
sy_2&0&1&0\\
-sy_1y_2&-y_2&-y_1&s^{-1}
\end{pmatrix},\qquad y_1+y_2=2(s^{-1}-1),\quad s\in\bO^\times.
\end{equation}
Thus,
\[
g'A'=\begin{pmatrix}
st^a&0&0&0\\
sy_1t^a+t^{b_0-m}x&t^{b_0}&0&0\\
sy_2t^a&0&t^{-b_0}&0\\
*&*&*&s^{-1}t^{-a}
\end{pmatrix}.
\]
By the $(2,1)$- and $(3,1)$-entries, we obtain
\[
sy_1+t^{b_0-a-m}x\in t^{b_0-a}\bO,\qquad sy_2\in t^{-a-b_0}\bO.
\]
Since $a+b_0\leq0$, by the second inclusion and $s(y_1+y_2)=2(1-s)$, we have $sy_1\in\bO$. If $0\leq b_0-a<m$, then $t^{b_0-a-m}x$ has negative valuation and cannot be cancelled, a contradiction.

If $b_0-a=m$ and $a+b_0<0$, then by the two inclusions, we have $x(0)=2(s(0)-1)\neq0,-2$. Conversely, we take $s=1+x/2$, $y_1=-x/s$, and $y_2=0$, and then $g'A'=\diag(st^a,t^{b_0},t^{-b_0},s^{-1}t^{-a})$. Extending $g'$ by the identity, we have $\sD g'\sD^{-1}\in U_1^-(\bO)T_1(\bO)$. Hence, the intersection \eqref{eq:g2-mixed-intersection} is $(\BG_m\setminus\{-2\})\times\BA^{m-1}$. Its diagonal cocharacter is
\[
\nu=(e_l^{(1)},e_{k-1}^{(2)}).
\]
Note that its negative is the common odd summand in
\[
\begin{split}
(0,-e_{k-1}^{(2)}+e_k^{(2)}) &=(-e_l^{(1)},-e_{k-1}^{(2)}) +(e_l^{(1)},e_k^{(2)}),\\
(0,-e_{k-1}^{(2)}-e_k^{(2)}) &=(-e_l^{(1)},-e_{k-1}^{(2)}) +(e_l^{(1)},-e_k^{(2)}).
\end{split}
\]
The local monodromy is $b(\lambda,\nu)$ up to inversion, thus, if the local system is constant, then $b(\lambda,\nu)=1$. Then, since the character of the even root is trivial, the character corresponding to the other odd simple super root is also trivial.

If $a+b_0=0$, then we take $y_1=0$ and $y_2=2(s^{-1}-1)$, with $s\in\BC^\times$. The matrix $g'$ stabilizes $t^{\lambda^{(2)}}G_2(\bO)$, and $\sD g'\sD^{-1}$ is integral lower triangular in $G_1$, thus, it fixes $t^{\lambda^{(1)}}$. Hence, by Proposition~\ref{prop:relevant-L-orbits}, the relevance condition implies that $b(\lambda,\nu)=1$, and we obtain the same conclusion. Finally, by interchanging the two middle rows and columns, we replace $b_0$ by $-b_0$, and the same argument works for the coroot $-e_{k-1}^{(2)}-e_k^{(2)}$, including the case $a-b_0=0$.

\case{The $\on{F}(4)$ case.}\par\nopagebreak[4]

On the $\on{SL}_2$-cover, the embedding of $G_1$ into $L$ is
\[
\iota(g)=\diag(g,g,Jg^{-t}J),\qquad J=\begin{pmatrix}0&1\\1&0\end{pmatrix}.
\]
Here, $H_u$ is the unipotent radical of the lower parabolic subgroup with block sizes $(2,2,2)$. We need to prove that
\begin{equation}\label{eq:g2-f4-mixed-intersection}
\sD^{-1}U_1^-(\bF)T_1(\bO)H_u(\bF)\sD t^{\lambda^{(2)}} \PSp(6)(\bO)/\PSp(6)(\bO) \cap\sS_{G_2}^{\theta^{(2)}} =\varnothing
\end{equation}
if the upper inequality in \eqref{eq:restricted-mixed-inequality} is strict, for either of the mixed coroots
\[
\alpha_{2,1}=-\varepsilon_1+\varepsilon_2, \qquad \alpha_{2,3}=-\varepsilon_3.
\]

\smallskip
\noindent\emph{The coroot $\alpha_{2,1}$.} In the case $\lambda^{(2)}-\theta^{(2)}=m\alpha_{2,1}$, for any point in \eqref{eq:g2-mixed-intersection}, we can find a representative $A$, which is lower triangular and whose upper-left $4\times4$ submatrix is
\[
A'=
\begin{pmatrix}
t^{\lambda^{(2)}_1}&&&\\
t^{\lambda^{(2)}_2-m}\bO^\times&t^{\lambda^{(2)}_2}&&\\
&&t^{\lambda^{(2)}_3}&\\
&&&t^{-\lambda^{(2)}_3}
\end{pmatrix}.
\]
If the point represented by $A$ lies in the intersection on the left-hand side of \eqref{eq:g2-f4-mixed-intersection}, then there exist
\[
u=\begin{pmatrix}1&0\\ f&1\end{pmatrix}\in U_1^-(\bF), \qquad a\in T_1(\bO),\qquad u'\in H_u(\bF),
\]
with
\[
\sD^{-1}uau'\sD A\in t^{\lambda^{(2)}}\PSp(6)(\bO).
\]
Since the conjugation by $\sD$ does not change the $(2,1)$-entry, we have $f\in t^{-\lambda^{(2)}_1-m+\lambda^{(2)}_2}\bO^\times$, and then, by writing down the expression of the above matrix, we see that its $(4,3)$-entry belongs to
\[
t^{\lambda^{(2)}_3}\bO +t^{-\lambda^{(2)}_1-m+\lambda^{(2)}_2+\lambda^{(2)}_3}\bO^\times +t^{-\lambda^{(2)}_3}\bO.
\]
Note that by the restricted dominant condition, we have $\lambda^{(2)}_2-\lambda^{(2)}_1<m$ in this case. Thus, the middle term cannot be cancelled, and its valuation is less than $-\lambda^{(2)}_3$, which contradicts the requirement that this entry belongs to $t^{-\lambda^{(2)}_3}\bO$.

If $\lambda^{(2)}_2-\lambda^{(2)}_1=m= \ord(\xi^8)$, recall that
\[
b(\lambda,\lambda') =\xi^{\,8\sum_{r=1}^3\lambda^{(2)}_r\lambda^{(2)\prime}_r-24\lambda^{(1)}\lambda^{(1)\prime}}.
\]
For $a^{-1}=\diag(c,c^{-1})$ on the $\on{SL}_2$ cover, the coordinate on the quotient torus in the projective groups is $s=c^2$. The corresponding cocharacter is
\[
\nu=\tfrac12(1;1,-1,1)=-\alpha_3.
\]
Thus, we have $b(\lambda,\nu)= \xi^{4\lambda^{(2)}_1-4\lambda^{(2)}_2+4\lambda^{(2)}_3-12\lambda^{(1)}}$, cf., \cite[Proposition 3.3.1(c)]{[FTY]} and the identification \eqref{eq:F4-duality-coordinates}. Then, we can write
\[
A=\sD^{-1}(u')^{-1}a^{-1}u^{-1}\sD t^{\lambda^{(2)}} g, \qquad g\in\PSp(6)(\bO),
\]
where, in fact, $g\in U_2^-(\bO)T_2(\bO)$, and the torus part is given by $\diag(c^{-1},c,c^{-1},c,c^{-1},c)$. The upper-left $4\times4$ submatrix of $\sD^{-1}(u')^{-1}a^{-1}u^{-1}\sD t^{\lambda^{(2)}}$ is given by
\[
\begin{pmatrix}
c t^{\lambda^{(2)}_1}&&&\\
-c^{-1}t^{\lambda^{(2)}_1}f&c^{-1}t^{\lambda^{(2)}_2}&&\\
*&*&c t^{\lambda^{(2)}_3}&\\
*&*&(c^{-1}-c^{-1}f-c)t^{\lambda^{(2)}_3}&c^{-1}t^{-\lambda^{(2)}_3}
\end{pmatrix}.
\]
If $\lambda^{(2)}_3<0$, then the coordinate $x=-c^{-2}f$ satisfies $x(0)=1-c(0)^{-2}$, and changing its sign, we obtain
\[
(\BG_m\setminus\{-1\})\times\BA^{m-1} \subset t^{\lambda^{(2)}_2-m}\bO^\times/t^{\lambda^{(2)}_2}\bO.
\]
The local monodromy is $b(\lambda,\nu)$ up to inversion, thus, if the local system is constant, then $b(\lambda,\alpha_{2,1}^{\pm})=1$. If $\lambda^{(2)}_3=0$, then the relevance condition for $\alpha_{2,3}$ (Proposition~\ref{prop:relevant-L-orbits}) implies that $b(\lambda,\alpha_{2,3}^{\pm})=1$, which implies that $b(\lambda,\alpha_{2,1}^{\pm})=1$.

\smallskip
\noindent\emph{The coroot $\alpha_{2,3}$.} In the case $\lambda^{(2)}-\theta^{(2)}=m\alpha_{2,3}$, for any point in \eqref{eq:g2-mixed-intersection}, we can find a representative $A$, which is lower triangular and whose upper-left $4\times4$ submatrix is
\[
A'=
\begin{pmatrix}
t^{\lambda^{(2)}_1}&&&\\
&t^{\lambda^{(2)}_2}&&\\
&&t^{\lambda^{(2)}_3}&\\
&&t^{-\lambda^{(2)}_3-m}\bO^\times&t^{-\lambda^{(2)}_3}
\end{pmatrix}.
\]
Now, by the $(2,1)$-entry, we have $f\in t^{-\lambda^{(2)}_1+\lambda^{(2)}_2}\bO$, and the $(4,3)$-entry belongs to
\[
t^{\lambda^{(2)}_3}\bO +t^{-\lambda^{(2)}_3-m}\bO^\times +t^{\lambda^{(2)}_3+\lambda^{(2)}_2-\lambda^{(2)}_1}\bO.
\]
Note that by the restricted dominant condition, we have $-2\lambda^{(2)}_3<m$ in this case. Thus, the second term cannot be cancelled, and this entry has valuation $-\lambda^{(2)}_3-m<-\lambda^{(2)}_3$, a contradiction. If $2\lambda^{(2)}_3+m=0$ and $m=\ord(\xi^4)$, then, by the same analysis with the same torus, the local monodromy is given by $b(\lambda,\nu)$ if $\lambda^{(2)}_1<\lambda^{(2)}_2$. If $\lambda^{(2)}_1=\lambda^{(2)}_2$, then the relevance condition for $\alpha_{2,1}$ (Proposition~\ref{prop:relevant-L-orbits}) implies that $b(\lambda,\alpha_{2,1}^{\pm})=1$. In either case, we obtain $b(\lambda,\alpha_{2,3}^{\pm})=1$, which contradicts the restricted assumption.

Thus, every non-zero difference in \eqref{eq:g2-rank-one-difference} is excluded, and hence, $\lambda^{(2)}=\theta^{(2)}$.
\end{proof}

\subsection{The \texorpdfstring{$G_1$}{G1}-direction: proof that \texorpdfstring{$\lambda^{(1)}=\theta^{(1)}$}{lambda1=theta1}}\label{sec:g1-direction}

In this subsection, we prove Proposition~\ref{prop:part-p}. By Proposition~\ref{prop:part-q}, we have $\lambda^{(2)}=\theta^{(2)}$. Now, we study the fiber of the projection \eqref{eq:g2-projection} over $t^{\lambda^{(2)}}$, which is an intersection in $\Gr_{G_1}$. The strategy is the same as in Section~\ref{sec:g2-direction}, with the roles of $G_1$ and $G_2$ interchanged. Namely, we first prove $\lambda^{(1)}-\theta^{(1)}\in\Lambda_{G_1}^{\sharp}$, and then reduce the proof of $\lambda^{(1)}=\theta^{(1)}$ to the rank-one case. However, the comparison with semi-infinite orbits of $G_1$ needs more work. First, we compare the critical-dimensional components of the fiber with those of the intersection of an orbit of a certain subgroup of $G_1(\bF)$ and a semi-infinite orbit, see Corollary~\ref{cor:g1-critical-component-comparison}. Secondly, in the orthosymplectic case, this subgroup is contained in $U_1(\bF)T_1(\bO)$ only after a conjugation,\footnote{we can also change the choice of $\sD$, but it will make the introduction of $\sD$ in Section \ref{trans D} more heavy} see Lemma~\ref{lem:orthogonal-auxiliary-Borel}.

\subsubsection{The critical fiber and its boundary}

Note that the fiber of the projection \eqref{eq:g2-projection} over the point $t^{\theta^{(2)}} G_2(\bO)/G_2(\bO)$ is given by
\begin{equation}
\label{eq:g1-fiber-intersection}\fitmath{\begin{aligned} &\sD t^{\theta^{(2)}} G_2(\bO)t^{-\lambda^{(2)}}\sD^{-1}t^{\lambda^{(1)}} H_u(\bF)G_1(\bO)/G_1(\bO) \cap U_1^-(\bF)t^{\theta^{(1)}} G_1(\bO)/G_1(\bO)\\
&\quad= \bigl( \sD t^{\theta^{(2)}} G_2(\bO)t^{-\lambda^{(2)}}\sD^{-1} \cap G_1(\bF) \bigr)t^{\lambda^{(1)}} G_1(\bO)/G_1(\bO) \cap U_1^-(\bF)t^{\theta^{(1)}} G_1(\bO)/G_1(\bO) \subset\Gr_{G_1}.
\end{aligned}}
\end{equation}

Under the assumptions of Theorem~\ref{thm:local-system-rigidity}, by Corollary~\ref{rem:unique-g2-critical-stratum} and the equivariant product \eqref{eq:Y-convolution}, the fiber \eqref{eq:g1-fiber-intersection} has a critical-dimensional irreducible component, such that the restriction of $\cL_{\lambda,\theta}\otimes\psi$ is constant on a dense open subset of it. 

\begin{cor}\label{cor:g1-fiber-dimension}
If $(\lambda^{(1)},\lambda^{(2)})\in\Lambda_G^{L,+}$, then
\[
\dim\eqref{eq:g1-fiber-intersection} \leq \frac12d(\lambda-\theta).
\]
\end{cor}

\begin{proof}
This follows from Proposition~\ref{prop:g2-fiber-dimension}, with $\lambda'=(\theta^{(1)},\theta^{(2)})$.
\end{proof}

\begin{lem}\label{lem:g1-orbit-closure-decomposition}
For $(\lambda^{(1)},\lambda^{(2)})\in\Lambda_G^{L,+}$, the closure of
\[
\sD t^{\theta^{(2)}} G_2(\bO)t^{-\lambda^{(2)}}\sD^{-1}t^{\lambda^{(1)}} H_u(\bF)G_1(\bO)/G_1(\bO)
\]
is contained in
\[
\bigcup_{ (\lambda^{(1)\prime},\lambda^{(2)\prime})\in\Lambda_G^{L,+}, \ (\lambda^{(1)\prime},\lambda^{(2)\prime})\leq_{\super}(\lambda^{(1)},\lambda^{(2)})} \sD t^{\theta^{(2)}} G_2(\bO)t^{-\lambda^{(2)\prime}}\sD^{-1}t^{\lambda^{(1)\prime}} H_u(\bF)G_1(\bO)/G_1(\bO).
\]
\end{lem}

\begin{proof}
This follows from Lemma~\ref{lem:L-orbit-indices} and Lemma~\ref{lem:spherical-orbit-closure}, after translating by $\sD t^{\theta^{(2)}}$.
\end{proof}

Here, we are only interested in the case $\lambda^{(2)}=\theta^{(2)}$, and then, by the preceding lemma, we obtain the following analogue of Corollary~\ref{cor:g2-boundary-codimension}.

\begin{cor}\label{cor:g1-boundary-codimension}
Assume that $\lambda^{(1)}\geq_{G_1}\theta^{(1)}$. Then, the intersection of the boundary of
\[
\sD t^{\lambda^{(2)}} G_2(\bO)t^{-\lambda^{(2)}}\sD^{-1}t^{\lambda^{(1)}} H_u(\bF)G_1(\bO)/G_1(\bO)
\]
in $G_2(\bF)/G_1(\bO)$ with the boundary of
\[
U_1^-(\bF)t^{\theta^{(1)}} G_1(\bO)/G_1(\bO)
\]
has dimension at most
\[
\frac12d\bigl((\lambda^{(1)}-\theta^{(1)},0)\bigr)-2.
\]
\end{cor}

\begin{proof}
The proof is the same as that of Corollary~\ref{cor:g2-boundary-codimension}, applied to the orbit decomposition in Lemma~\ref{lem:g1-orbit-closure-decomposition} and the dimension estimate in Corollary~\ref{cor:g1-fiber-dimension}.
\end{proof}

\subsubsection{Comparison with semi-infinite orbits}

\noindent\emph{Comparison of critical components.} In order to use the same reduction to the rank-one case as in the $G_2$-direction, we also consider the following intersection,
\begin{equation}
\label{eq:g1-limiting-intersection}\fitmath{\begin{aligned} &\sD t^{\theta^{(2)}} U_2(\bF)T_2(\bO)t^{-\lambda^{(2)}}\sD^{-1}t^{\lambda^{(1)}} H_u(\bF)G_1(\bO)/G_1(\bO) \cap U_1^-(\bF)t^{\theta^{(1)}} G_1(\bO)/G_1(\bO)\\
&\quad= \sD t^{\theta^{(2)}} U_L(\bF)T_2(\bO)t^{-\lambda^{(2)}}\sD^{-1}t^{\lambda^{(1)}} G_1(\bO)/G_1(\bO) \cap U_1^-(\bF)t^{\theta^{(1)}} G_1(\bO)/G_1(\bO)\\
&\quad= \bigl( \sD t^{\theta^{(2)}} U_L(\bF)T_2(\bO)t^{-\lambda^{(2)}}\sD^{-1} \cap G_1(\bF) \bigr)t^{\lambda^{(1)}} G_1(\bO)/G_1(\bO) \cap U_1^-(\bF)t^{\theta^{(1)}} G_1(\bO)/G_1(\bO).
\end{aligned}}
\end{equation}
 It is the analogue of \eqref{eq:critical-G2-fiber}, and we will compare its critical-dimensional components with those of \eqref{eq:g1-fiber-intersection}. Here, $U_L=U_2\cap L$, and since $T_2\subset L$, we take $T_2$ as the maximal torus of $L$.

\begin{cor}\label{cor:g1-limiting-dimension}
If $\lambda^{(1)}$ is dominant, then
\[
\dim\eqref{eq:g1-limiting-intersection} \leq \frac12d(\lambda-\theta).
\]
\end{cor}

\begin{proof}
The inequality follows from Proposition~\ref{prop:g1-strata-dimension} with $\lambda'=\lambda$.
\end{proof}

\begin{lem}\label{lem:g1-orbit-order}
For a given dominant $(\lambda^{(1)},\lambda^{(2)})\in\Lambda_G^{L,+}$, if
\[
\begin{split}
&G_2(\bO)t^{-\lambda^{(2)}}\sD^{-1}t^{\lambda^{(1)}} H_u(\bF) G_1(\bO)/G_1(\bO)\\
&\qquad\cap U_2(\bF)T_2(\bO)t^{-\lambda^{(2)\prime}}\sD^{-1}t^{\lambda^{(1)\prime}}H_u(\bF) G_1(\bO)/G_1(\bO)\neq\varnothing,
\end{split}
\]
then
\[
(\lambda^{(1)\prime},\lambda^{(2)\prime})\leq_{\super}(\lambda^{(1)},\lambda^{(2)}).
\]
\end{lem}

\begin{proof}
By \eqref{eq:g1-projection-fiber}, the intersection in the statement is non-empty if and only if the fiber of $\sO_{\Gr}^{\lambda}\cap\sS_{\Gr}^{\lambda'}\to\Gr_{G_1}$ over $t^{\lambda^{(1)\prime}}$ is non-empty, and then the assertion follows from Lemma~\ref{lem:orbit-S-nonempty-super-order}.
\end{proof}

\begin{lem}\label{lem:equal-index-intersection}
Let $\lambda\in\Lambda_G^{L,+}$, and let us take the representatives in Lemma~\ref{lem:L-orbit-indices}. Then, on the reduced spaces,
\[
\sO_{\Gr}^{\lambda}\cap\sS_{\Gr}^{\lambda}=\{t^\lambda\}.
\]
\end{lem}
\begin{proof}
First, we choose a local coordinate $t$ and a trivialization of the $\omega$-twist. By translating by a sufficiently dominant coweight of $Z_L$ in the second factor, we can assume that $\lambda$ is dominant for $B^-$, and then, by Corollary~\ref{lem:dimension-estimate-super}, we have $\dim(\sO_{\Gr}^{\lambda}\cap\sS_{\Gr}^{\lambda})\leq0$.

Note that the loop rotation $g(t)\mapsto g(a^{-1}t)$, $a\in\BG_m$, preserves both orbits, and it contracts $\sS_{\Gr}^{\lambda}$ to $t^\lambda$, since the root coordinates of $\sS_{\Gr}^{\lambda}$ have strictly positive weights. Thus, any other point of the intersection would have a one-dimensional orbit under the loop rotation, which contradicts the dimension bound. Hence, the reduced intersection is $\{t^\lambda\}$.
\end{proof}

\begin{lem}\label{lem:g1-open-stratum}
For $(\lambda^{(1)},\lambda^{(2)})\in\Lambda_G^{L,+}$, we have
\begin{equation}\label{eq:g1-open-stratum-intersection}
\begin{split}
&G_2(\bO)t^{-\lambda^{(2)}}\sD^{-1}t^{\lambda^{(1)}} H_u(\bF) G_1(\bO)/G_1(\bO)\\
&\qquad\cap U_2(\bF)T_2(\bO)t^{-\lambda^{(2)}}\sD^{-1}t^{\lambda^{(1)}} H_u(\bF) G_1(\bO)/G_1(\bO)\\
&\quad\simeq U_2(\bO)T_2(\bO)t^{-\lambda^{(2)}}\sD^{-1}t^{\lambda^{(1)}} H_u(\bF) G_1(\bO)/G_1(\bO).
\end{split}
\end{equation}
\end{lem}

\begin{proof}
It follows from the fact that, by Lemma~\ref{lem:equal-index-intersection}, the fiber \eqref{eq:g1-projection-fiber} for $\lambda'=\lambda$ with first coordinate $t^{\lambda^{(1)}}$ is just a point. We take the preimages of the fiber and of $t^{\lambda^{(2)}}$ under $G_2(\bF)\to\Gr_{G_2}$, multiply on the left by $t^{-\lambda^{(2)}}$, take inverses and intersect with $U_2(\bF)T_2(\bO)$, and the above fact implies that
\[
\begin{split}
&U_2(\bF)T_2(\bO)\cap G_2(\bO)t^{-\lambda^{(2)}}\sD^{-1}t^{\lambda^{(1)}}H_u(\bF) G_1(\bO)t^{-\lambda^{(1)}}\sD t^{\lambda^{(2)}}\\
&\hspace{35mm}=U_2(\bO)T_2(\bO).
\end{split}
\]
Then, multiplying on the right by $t^{-\lambda^{(2)}}\sD^{-1}t^{\lambda^{(1)}}H_u(\bF)G_1(\bO)$ and taking the quotient by $G_1(\bO)$, we obtain \eqref{eq:g1-open-stratum-intersection}.
\end{proof}

Combining Lemmas~\ref{lem:g1-orbit-order} and \ref{lem:g1-open-stratum}, we obtain the following corollary.

\begin{cor}\label{cor:g1-critical-component-comparison}
Assume that $\lambda^{(2)}=\theta^{(2)}$. Then, there is a bijection between the critical-dimensional irreducible components of \eqref{eq:g1-fiber-intersection} and the critical-dimensional irreducible components of
\begin{equation}
\label{eq:g1-integral-intersection}\fitmath{\begin{aligned} &\sD t^{\lambda^{(2)}} U_2(\bO)T_2(\bO)t^{-\lambda^{(2)}}\sD^{-1}t^{\lambda^{(1)}} H_u(\bF)G_1(\bO)/G_1(\bO) \cap U_1^-(\bF)t^{\theta^{(1)}} G_1(\bO)/G_1(\bO)\\
&\quad= \sD t^{\lambda^{(2)}} U_L(\bO)T_2(\bO)t^{-\lambda^{(2)}}\sD^{-1}t^{\lambda^{(1)}} G_1(\bO)/G_1(\bO) \cap U_1^-(\bF)t^{\theta^{(1)}} G_1(\bO)/G_1(\bO)\\
&\quad= \bigl( \sD t^{\lambda^{(2)}} U_L(\bO)T_2(\bO)t^{-\lambda^{(2)}}\sD^{-1} \cap G_1(\bF) \bigr)t^{\lambda^{(1)}} G_1(\bO)/G_1(\bO) \cap U_1^-(\bF)t^{\theta^{(1)}} G_1(\bO)/G_1(\bO).
\end{aligned}}
\end{equation}
Also, the set of critical-dimensional irreducible components of \eqref{eq:g1-integral-intersection} injects into the set of critical-dimensional irreducible components of \eqref{eq:g1-limiting-intersection}, and the equivariant trivializations and the corresponding ratio local systems match under this bijection and injection, on the common dense open subsets.
\end{cor}

\begin{proof}
We stratify the fiber \eqref{eq:g1-fiber-intersection} by its intersections with $\CS_{G_1}^{(\lambda^{(2)};\mu^{(2)},\mu^{(1)})}$, where $\mu=(\mu^{(1)},\mu^{(2)})$ is as in Lemma~\ref{lem:g1-semi-spherical-indices}. By Lemma~\ref{lem:g1-orbit-order}, we have $\mu\leq_{\super}\lambda$, and by Proposition~\ref{prop:g1-strata-dimension}, the strata with $\mu<_{\super}\lambda$ do not have critical dimension. By Lemma~\ref{lem:g1-open-stratum}, the intersection for $\mu=\lambda$ is \eqref{eq:g1-integral-intersection}, which proves the bijection. Since \eqref{eq:g1-integral-intersection} is locally closed and both dimensions are bounded by $\frac12d(\lambda-\theta)$, we obtain the injection by taking closures in \eqref{eq:g1-limiting-intersection}. The second claim is clear, since the trivializations are the same.
\end{proof}

\smallskip
\subsubsection{The subgroup $K_{\lambda^{(2)}}^{(1)}$}
Next, we compare \eqref{eq:g1-limiting-intersection} with $\sS_{G_1}^{+,\lambda^{(1)}}\cap\sS_{G_1}^{-,\theta^{(1)}}$. For this, we need the following group
\[
K_{\lambda^{(2)}}^{(1)} := \sD t^{\lambda^{(2)}} U_2(\bF)T_2(\bO)t^{-\lambda^{(2)}}\sD^{-1} \cap G_1(\bF),
\]
which plays a similar role to $K_{\lambda^{(1)}}^{(2)}$ in Section~\ref{sec:g2-direction}.

We expect an inclusion similar to \eqref{eq:g2-orbit-in-semi-infinite}. In type A and in the $\on F(4)$ case, we indeed have $K_{\lambda^{(2)}}^{(1)}\subset U_1(\bF)T_1(\bO)$, see \eqref{eq:g1-K-type-A} and \eqref{eq:g1-K-F4} below. However, in the $\SOSp$ case, this inclusion only holds after conjugation by a constant element $h\in B_1^-$, which is given in the following lemma.

\begin{lem}\label{lem:orthogonal-auxiliary-Borel}
For $M=2k$, we take
\[
h=\left(
\begin{array}{cccc|cccc}
i\sqrt{2}&&&&&&&\\
-i\sqrt{2}&i\sqrt{2}&&&&&&\\
&\ddots&\ddots&&&&&\\
&&-i\sqrt{2}&i\sqrt{2}&&&&\\ \hline
&&&&-\frac{i}{\sqrt{2}}&&&\\
&&&&-\frac{i}{\sqrt{2}}&-\frac{i}{\sqrt{2}}&&\\
&&&&\vdots&\vdots&\ddots&\\
&&&&-\frac{i}{\sqrt{2}}&-\frac{i}{\sqrt{2}}&\cdots&-\frac{i}{\sqrt{2}}
\end{array}\right),
\]
where the diagonal blocks have sizes $k,k$. For $M=2l+1$, we take
\[
h=\left(
\begin{array}{cccc|c|cccc}
\frac{i}{\sqrt{2}}&&&&&&&&\\
-\frac{i}{\sqrt{2}}&\frac{i}{\sqrt{2}}&&&&&&&\\
&\ddots&\ddots&&&&&&\\
&&-\frac{i}{\sqrt{2}}&\frac{i}{\sqrt{2}}&&&&&\\ \hline
&&&-1&1&&&&\\ \hline
&&&\frac{i}{\sqrt{2}}&-i\sqrt{2}&-i\sqrt{2}&&&\\
&&&\frac{i}{\sqrt{2}}&-i\sqrt{2}&-i\sqrt{2}&-i\sqrt{2}&&\\
&&&\vdots&\vdots&\vdots&\vdots&\ddots&\\
&&&\frac{i}{\sqrt{2}}&-i\sqrt{2}&-i\sqrt{2}&-i\sqrt{2}&\cdots&-i\sqrt{2}
\end{array}\right),
\]
where the diagonal blocks have sizes $l,1,l$. Then,
\[
h^{-1}K_{\lambda^{(2)}}^{(1)}h=\{g\in U_1(\bF)T_1(\bO)\mid gv=v\},
\]
where
\begin{equation}\label{eq:orthogonal-fixed-vector-normal-form}
v=\begin{cases}
(\underbrace{1,\ldots,1}_{k+1},\underbrace{0,\ldots,0}_{k-1})^t,&M=2k,\\
(\underbrace{1,\ldots,1}_{l+1},\underbrace{0,\ldots,0}_{l})^t,&M=2l+1.
\end{cases}
\end{equation}
Moreover, $ht^{\lambda^{(1)}}G_1(\bO)=t^{\lambda^{(1)}}G_1(\bO)$ for dominant $\lambda^{(1)}$, and $h$ preserves every $\sS_{G_1}^{-,\theta^{(1)}}$.
\end{lem}
\begin{proof}
It follows by a direct matrix calculation with $\sD$ in Section~\ref{trans D}.
\end{proof}

\subsubsection{Local systems and reduction to rank one}

In all three cases, we denote
\[
\begin{aligned}
\CT_{G_1}^{\lambda}&:=K_{\lambda^{(2)}}^{(1)}t^{\lambda^{(1)}}G_1(\bO)/G_1(\bO)\\
&=\CS_{G_1}^{(\lambda^{(2)};\lambda^{(2)},\lambda^{(1)})}\cap\Gr_{G_1}.
\end{aligned}
\]
Similar to \eqref{eq:g2-orbit-lift}, the equivariant trivialization on $\CT_{G_1}^{\lambda}$ is induced from the original $H(\bF)$-orbit via the map
\begin{equation}\label{eq:g1-orbit-lift}
\begin{aligned}
r_1:\CT_{G_1}^{\lambda}&\longrightarrow\sO_{\Gr}^{\lambda},\\
kt^{\lambda^{(1)}}G_1(\bO)&\longmapsto (kt^{\lambda^{(1)}}G_1(\bO),\sD^{-1}k\sD t^{\lambda^{(2)}}G_2(\bO)).
\end{aligned}
\end{equation}
Note that this map is well-defined. Indeed, if $k$ fixes $t^{\lambda^{(1)}}$, then by the dominance and the calculations of the subgroup, see \eqref{eq:g1-K-type-A}, Lemma~\ref{lem:orthogonal-auxiliary-Borel} and \eqref{eq:g1-K-F4}, we have $\sD^{-1}k\sD\in U_2(\bO)T_2(\bO)$, which fixes $t^{\lambda^{(2)}}$. Note also that the second coordinate lies in $\sS_{G_2}^{+,\lambda^{(2)}}$. If we denote by $\pr_2$ the second projection, then we define
\begin{equation}\label{eq:g1-orbit-trivialization}
\tau_1^{\lambda}:=r_1^*\tau_H^\lambda\otimes \bigl((\pr_2\circ r_1)^*\tau_{2,+}^{\lambda^{(2)}}\bigr)^{-1}.
\end{equation}

In the orthosymplectic case, we make the change of coordinates $x\mapsto h^{-1}x$, where $h$ is the matrix in Lemma~\ref{lem:orthogonal-auxiliary-Borel}. With some abuse of notation, we still denote by $K_{\lambda^{(2)}}^{(1)}$ and $\CT_{G_1}^{\lambda}$ the conjugate $h^{-1}K_{\lambda^{(2)}}^{(1)}h$ and the translate $h^{-1}\CT_{G_1}^{\lambda}$. We apply the same change of coordinates to \eqref{eq:g1-fiber-intersection}, \eqref{eq:g1-integral-intersection} and \eqref{eq:g1-limiting-intersection}, and keep the same notation for these intersections. Also, we pull back $r_1$, $\tau_1^\lambda$, $\tau_{1,-}^{\theta^{(1)}}$ and $\psi$ under $x\mapsto hx$. Thus, a matrix $k$ below represents $hkh^{-1}$ in the original coordinates. Note that, by Lemma~\ref{lem:orthogonal-auxiliary-Borel}, the point $t^{\lambda^{(1)}}G_1(\bO)$ and the negative semi-infinite orbits are preserved by $h$. Hence, in all three cases, we have
\[
\CT_{G_1}^{\lambda}\subset\sS_{G_1}^{+,\lambda^{(1)}}.
\]

Similar to \eqref{eq:critical-G2-fiber}, we denote the intersection \eqref{eq:g1-limiting-intersection} by
\[
\mathfrak F_1(\lambda^{(1)},\lambda^{(2)};\theta^{(1)}):=\CT_{G_1}^{\lambda}\cap\sS_{G_1}^{-,\theta^{(1)}}.
\]
It is contained in $\sS_{G_1}^{+,\lambda^{(1)}}\cap\sS_{G_1}^{-,\theta^{(1)}}$.

As in Section~\ref{sec:g2-direction}, we compare the trivializations with the ones on semi-infinite orbits. On $\CT_{G_1}^{\lambda}$, we denote
\[
\cL_1^{(1)}:=\tau_1^{\lambda}\otimes \bigl(\tau_{1,+}^{\lambda^{(1)}}|_{\CT_{G_1}^{\lambda}}\bigr)^{-1}.
\]
On $\sS_{G_1}^{+,\lambda^{(1)}}\cap\sS_{G_1}^{-,\theta^{(1)}}$, we denote
\[
\cL_2^{(1)}:=\tau_{1,-}^{\theta^{(1)}}\otimes \bigl(\tau_{1,+}^{\lambda^{(1)}}\bigr)^{-1},
\]
where both trivializations are restricted to this intersection. Since $\theta^{(2)}=\lambda^{(2)}$, on $\mathfrak F_1(\lambda^{(1)},\lambda^{(2)};\theta^{(1)})$, we have
\begin{equation}\label{eq:g1-ratio-decomposition}
r_1^*\cL_{\lambda,\theta}\simeq\cL_1^{(1)}\otimes\bigl(\cL_2^{(1)}\bigr)^{-1}.
\end{equation}
Note that, on \eqref{eq:g1-integral-intersection}, the second coordinate of $r_1$ is $t^{\lambda^{(2)}}$. Thus, as in Corollary~\ref{cor:g1-critical-component-comparison}, these local systems agree with the ones on the original fiber \eqref{eq:g1-fiber-intersection}, and we still denote the pullbacks under $r_1$ by $\cL_{\lambda,\theta}$ and $\psi$.

By Corollaries~\ref{cor:g1-boundary-codimension} and \ref{cor:g1-critical-component-comparison}, we can apply the argument of Lemma~\ref{lem:rank-one-fixed-orbit-comparison} to \eqref{eq:g1-ratio-decomposition}, with the two Borels interchanged. Then, as in Corollary~\ref{prop:g2-sharp-difference}, we obtain the following corollary.

\begin{cor}\label{prop:g1-sharp-difference}
Under the assumptions of Theorem~\ref{thm:local-system-rigidity}, we have
\[
\lambda^{(1)}-\theta^{(1)}\in\Lambda_{G_1}^{\sharp}.
\]
Furthermore, if the local system $\cL_{\lambda,\theta}\otimes\psi$ is constant on a dense open subset of a critical-dimensional component $Z$ of \eqref{eq:g1-fiber-intersection}, then the local system $\cL_2^{(1)}$ is constant on a dense open subset of the corresponding component of \eqref{eq:g1-limiting-intersection} in Corollary~\ref{cor:g1-critical-component-comparison}.
\end{cor}

\begin{lem}\label{lem:g1-hyperplane-block}
Let $Z$ be a critical-dimensional component of \eqref{eq:g1-limiting-intersection} corresponding to a component of the fiber \eqref{eq:g1-fiber-intersection}, as in Corollary~\ref{cor:g1-critical-component-comparison}. Suppose that $\cL_{\lambda,\theta}\otimes\psi$ and $\cL_2^{(1)}$ are constant on a dense open subset of $Z$. Let $q_i=\ord(q(\alpha_{1,i}))<\infty$. Here, the $i$-th operator is the one in Lemma~\ref{lem:rank-one-fixed-orbit-comparison}, with the two Borels interchanged, which increases the index of the negative semi-infinite orbit by $\alpha_{1,i}$. If we can apply it successively $q_i$ times, then the resulting component lies in $\sS_{G_1}^{-,\theta^{(1)}+q_i\alpha_{1,i}}$ and has dimension $\dim Z-q_i$, and both local systems remain constant on a dense open subset of this component.
\end{lem}
\begin{proof}
The proof is the same as that of Lemma~\ref{lem:rank-one-fixed-orbit-comparison}, with the two Borels interchanged. Here, by Corollaries~\ref{cor:g1-boundary-codimension} and \ref{cor:g1-critical-component-comparison}, we identify the resulting components with components which meet the fiber \eqref{eq:g1-fiber-intersection} in dense open subsets.
\end{proof}

Now, we reduce the proof of Proposition~\ref{prop:part-p} to the mixed rank-one case. Assume that $\lambda^{(1)}\neq\theta^{(1)}$. Then, by Corollary~\ref{prop:g1-sharp-difference} and Lemma~\ref{lem:g1-hyperplane-block}, we can repeat the argument as in Section~\ref{sec:g2-direction}, with the two Borels interchanged. Namely, we apply the operators successively, and in the last direction $\alpha_{1,i}$, we stop $\ord(q(\alpha_{1,i}))$ applications before reaching $t^{\lambda^{(1)}}$. Then, both local systems are still constant on a dense open subset of a component of the fiber \eqref{eq:g1-fiber-intersection}, and it is enough to consider the case
\begin{equation}\label{eq:g1-rank-one-difference}
\lambda^{(1)}-\theta^{(1)}=q_i\alpha_{1,i}, \qquad q_i=\ord(q(\alpha_{1,i}))<\infty.
\end{equation}
Here, $\theta$ denotes the index of the semi-infinite orbit after these cuts, and $\lambda$ is unchanged. If $\alpha_{1,i}$ is pure, then by \cite[Theorem 9.4.9]{[GL]}, we obtain a contradiction as in Section~\ref{sec:g2-direction}. Thus, it remains to consider the case that $\alpha_{1,i}$ is mixed. Note that for any simple coroot $\alpha_{1,i}$ of $G_1$, we always have that $\alpha_{1,i}$ is mixed, i.e., the pure case does not occur here.

\subsubsection{The mixed rank-one calculation}\label{10.4.5}

\begin{proof}[Completion of the proof of Proposition~\ref{prop:part-p}]
For a mixed coroot in \eqref{eq:g1-rank-one-difference}, we denote $m=q_i$, and consider
\begin{equation}\label{eq:g1-mixed-intersection}
\mathfrak F_1(\lambda^{(1)},\lambda^{(2)};\theta^{(1)})=\CT_{G_1}^{\lambda}\cap\sS_{G_1}^{-,\theta^{(1)}}.
\end{equation}
Note that \eqref{eq:g1-mixed-intersection} is contained in $\sS_{G_1}^{+,\lambda^{(1)}}\cap\sS_{G_1}^{-,\theta^{(1)}} \simeq\BG_m\times\BA^{m-1}$. As in Section~\ref{10.3.3}, if the upper inequality in \eqref{eq:restricted-mixed-inequality} is strict, then we only need to prove that \eqref{eq:g1-mixed-intersection} is empty. If equality holds in the upper bound, then we need to show that a constant local system would imply $b(\lambda,\alpha_{1,i}^{\pm})=1$, contrary to \eqref{eq:restricted-wall-monodromy}.

Note that, in the original coordinates, the families in the equality case have representatives $k$ satisfying $\sD^{-1}k\sD\in B_2(\bO)$, and hence they lie in the fiber \eqref{eq:g1-fiber-intersection}, since $t^{-\lambda^{(2)}}B_2(\bO)t^{\lambda^{(2)}}\subset G_2(\bO)$. Thus, their characters are those of the trivialization on $\sO_{\Gr}^{\lambda}$.

\case{The $\GL(M|N)$ case.}

In the upper-left $M\times M$ submatrix, we have
\begin{equation}\label{eq:g1-K-type-A}
K_{\lambda^{(2)}}^{(1)}= \{A\in B_M(\bF)\mid A_{jj}\in\bO^\times,\quad (1,\ldots,1)A=(1,\ldots,1)\}.
\end{equation}

For any point in \eqref{eq:g1-mixed-intersection}, we can find an upper triangular representative of the form
\[
A=\begin{pmatrix} t^{\lambda^{(1)}_1}&&&&&\\
&\ddots&&&&\\
&&t^{\lambda^{(1)}_i}&t^{\lambda^{(1)}_i-m}x&&\\
&&&t^{\lambda^{(1)}_{i+1}}&&\\
&&&&\ddots&\\
&&&&&t^{\lambda^{(1)}_M}
\end{pmatrix},\qquad x\in(\bO/t^m)^\times.
\]
Suppose that $kA\in t^{\lambda^{(1)}}\GL_M(\bO)$ for some $k\in K_{\lambda^{(2)}}^{(1)}$. By \eqref{eq:g1-K-type-A}, the column sums of $k$ are one, thus, the sum of the entries in the $(i+1)$-st column of $kA$ is
\[
t^{\lambda^{(1)}_{i+1}}+t^{\lambda^{(1)}_i-m}x.
\]
On the other hand, $kA$ is upper triangular, and its $(j,i+1)$-entry belongs to $t^{\lambda^{(1)}_j}\bO$. Since $\lambda^{(1)}$ is dominant, the displayed sum belongs to $t^{\lambda^{(1)}_{i+1}}\bO$. If $0\leq\lambda^{(1)}_i-\lambda^{(1)}_{i+1}<m$, then its valuation is $\lambda^{(1)}_i-m<\lambda^{(1)}_{i+1}$, a contradiction.

If $\lambda^{(1)}_i-\lambda^{(1)}_{i+1}=m$, then we take
\[
k(s)=\begin{pmatrix}1&1-s\\0&s\end{pmatrix}, \qquad s\in\bO^\times.
\]
We put this matrix in the positions $i,i+1$, and extend it by the identity. Multiplying $k(s)t^{\lambda^{(1)}}$ on the right by $\diag(1,s^{-1})$ in these positions, we obtain the upper-right entry $(s^{-1}-1)t^{\lambda^{(1)}_i-m}$. Hence, the intersection \eqref{eq:g1-mixed-intersection} is isomorphic to $(\BG_m\setminus\{-1\})\times\BA^{m-1}$, with coordinate $x=s^{-1}-1$. The diagonal cocharacter of these matrices is
\[
\nu=(e_{i+1}^{(1)},e_{i+1}^{(2)})=-\alpha_{2i+1}.
\]
The local monodromy around $x=-1$ is given by $b(\lambda,\nu)$, up to inversion. Since $\alpha_{1,i}=\alpha_{2i}+\alpha_{2i+1}$ and $b(\lambda,\alpha_{1,i})=1$, if the local system is constant, then we have $b(\lambda,\alpha_{2i})=b(\lambda,\alpha_{2i+1})=1$.

\case{The $\SOSp(2k|2l)$ case.}

By Lemma~\ref{lem:orthogonal-auxiliary-Borel}, $K_{\lambda^{(2)}}^{(1)}$ fixes the vector $v$ in \eqref{eq:orthogonal-fixed-vector-normal-form}. In the following, the representatives $A$ have the same entries as $t^{\lambda^{(1)}}$ outside the displayed submatrices. If a point represented by an upper triangular matrix $A$ lies in \eqref{eq:g1-mixed-intersection}, then there exists $g\in K_{\lambda^{(2)}}^{(1)}$ such that
\begin{equation}\label{eq:orthogonal-rank-one-orbit-condition}
gA\in t^{\lambda^{(1)}}\SO(M)(\bO).
\end{equation}

\smallskip
\noindent\emph{(a) $M=2k$.} We first consider the coroot $e_i^{(1)}-e_{i+1}^{(1)}$, with $1\leq i<k$. We denote $a=\lambda_i^{(1)}$ and $b=\lambda_{i+1}^{(1)}$. For any point in \eqref{eq:g1-mixed-intersection}, we can find an upper triangular representative $A$ whose submatrix in rows and columns $i,i+1,M-i,M+1-i$ is
\[
A'=\begin{pmatrix}
t^a&t^{a-m}x&0&0\\
0&t^b&0&0\\
0&0&t^{-b}&-t^{-b-m}x\\
0&0&0&t^{-a}
\end{pmatrix},\qquad x\in(\bO/t^m)^\times.
\]
Since $g$ fixes $v$ and preserves the bilinear form, by writing down the $(M+1-i)$-th column of $gA$, we have
\[
\sum_{r=k}^{M}(gA)_{r,M+1-i} =\sum_{r=k}^{M}A_{r,M+1-i} =t^{-a}-t^{-b-m}x.
\]
On the other hand, by \eqref{eq:orthogonal-rank-one-orbit-condition} and the dominance of $\lambda^{(1)}$, this expression belongs to $t^{-a}\bO$. If $0\leq a-b<m$, then its valuation is $-b-m<-a$, a contradiction.

If $a-b=m$, we first assume that $i<k-1$ or $a+b>0$. Then, by the dominance and \eqref{eq:orthogonal-rank-one-orbit-condition}, all the terms in the above sum except $(gA)_{M+1-i,M+1-i}$ belong to $t^{1-a}\bO$. Since the $(M+1-i,M+1-i)$-entry has valuation $-a$, we have $x(0)\neq1$. Conversely, we put the matrix
\[
k(s)=\begin{pmatrix}
s&1-s&0&0\\
0&1&0&0\\
0&0&1&1-s^{-1}\\
0&0&0&s^{-1}
\end{pmatrix},\qquad s\in\bO^\times,
\]
in the same four positions, and extend it by the identity. Note that this matrix fixes $v$. Multiplying $k(s)t^{\lambda^{(1)}}$ on the right by $\diag(s^{-1},1,1,s)$ in these positions, we obtain $A'$ with $x=1-s$. Hence, the intersection \eqref{eq:g1-mixed-intersection} is isomorphic to $(\BG_m\setminus\{1\})\times\BA^{m-1}$.

For the other end coroot $e_{k-1}^{(1)}+e_k^{(1)}$, we interchange the two middle rows and columns, and replace $b$ by $-b$. If $a+b=0$ and $a-b=m$, then we take the matrix obtained from $k(s)$ by this interchange, with $s\in\BC^\times$. Note that it fixes $t^{\lambda^{(1)}}$, and in the original coordinates, it stabilizes the pair of lattices, by the dominance argument before the $\GL(M|N)$ case above. Furthermore, its diagonal cocharacter is $(e_{k-1}^{(1)},e_l^{(2)})$, and by Proposition~\ref{prop:relevant-L-orbits}, the relevance condition implies that its character is trivial. The case $a-b=0$ and $a+b=m$ is similar.

\smallskip
\noindent\emph{(b) $M=2l+1$.} The preceding matrix calculation applies to $e_i^{(1)}-e_{i+1}^{(1)}$, for $1\leq i<l$. Here, the sum starts at the $(l+1)$-st row, and the $(l+1,M+1-i)$-entry of $gA$ is integral. Also, $a\geq m$ when $a-b=m$.

For the long coroot $2e_l^{(1)}$, we denote $a=\lambda_l^{(1)}\geq0$. For any point in \eqref{eq:g1-mixed-intersection}, we can find an upper triangular representative $A$ whose submatrix in rows and columns $l,l+1,l+2$ is
\[
A'=\begin{pmatrix}
t^a&t^{a-m}x&-\frac12t^{a-2m}x^2\\
0&1&-t^{-m}x\\
0&0&t^{-a}
\end{pmatrix},\qquad x\in(\bO/t^m)^\times.
\]
By writing down the $(l+2)$-th column of $gA$, we obtain
\[
\begin{aligned}
(gA)_{l+1,l+2}+(gA)_{l+2,l+2}
&=A_{l+1,l+2}+A_{l+2,l+2}\\
&=t^{-a}-t^{-m}x\in t^{-a}\bO.
\end{aligned}
\]
If $0\leq a<m$, then the second term cannot be cancelled, a contradiction. If $a=m$, then the $(l+1,l+2)$-entry of $gA$ is integral and its $(l+2,l+2)$-entry has valuation $-a$, thus, $x(0)\neq1$. Conversely, we put the matrix
\[
k(s)=\begin{pmatrix}
s&1-s&-\dfrac{(1-s)^2}{2s}\\
0&1&1-s^{-1}\\
0&0&s^{-1}
\end{pmatrix},\qquad s\in\bO^\times,
\]
in these three positions, and extend it by the identity. Note that it fixes $v$. Multiplying $k(s)t^{\lambda^{(1)}}$ on the right by $\diag(s^{-1},1,s)$, we obtain $A'$ with $x=1-s$. Hence, the intersection \eqref{eq:g1-mixed-intersection} is again isomorphic to $(\BG_m\setminus\{1\})\times\BA^{m-1}$.

Finally, we compute the monodromy. The matrices $k(s)$ above have the diagonal cocharacters
\begin{equation}\label{eq:orthogonal-G1-wall-cocharacter}
\nu_i=
\begin{cases}
(e_i^{(1)},e_{l-k+i+1}^{(2)}),&M=2k,\quad 1\leq i<k,\\
(e_i^{(1)},e_{k-l+i}^{(2)}),&M=2l+1,\quad 1\leq i\leq l.
\end{cases}
\end{equation}
Indeed, by the matrices in Lemma~\ref{lem:orthogonal-auxiliary-Borel}, the $(i+1)$-st diagonal entry in the middle $\SO(M+1)$-block of $\sD^{-1}hk(s)h^{-1}\sD$ is $s$. By \eqref{eq:osp-duality-coordinates}, $\nu_i$ is, up to sign, one of the two odd simple super roots defining the mixed coroot. Since $x=1-s$, the local monodromy at $x=1$ is given by $b(\lambda,\nu_i)$, up to inversion. Let us denote by $\alpha_{1,r}$ the coroot under consideration. Note that the requirement that $m=\ord(q(\alpha_{1,r}))$ implies that the ratio of the semi-infinite trivializations is trivial. Thus, we have $b(\lambda,\nu_i)=1$, either because the local system is constant, or by the relevance condition above. Together with the character of the even root, we obtain $b(\lambda,\alpha_{1,r}^{\pm})=1$.

\case{The $\on{F}(4)$ case.}

We have
\begin{equation}\label{eq:g1-K-F4}
K_{\lambda^{(2)}}^{(1)}= \left\{\begin{pmatrix}a&a^{-1}-a\\0&a^{-1}\end{pmatrix} :a\in\bO^\times\right\}/\{\pm1\}.
\end{equation}

Let $l=\lambda^{(1)}\geq0$. For any point in \eqref{eq:g1-mixed-intersection}, we can find a representative of the form
\[
A=\begin{pmatrix}t^l&t^{l-m}x\\0&t^{-l}\end{pmatrix},\qquad x\in(\bO/t^m)^\times.
\]
By \eqref{eq:g1-K-F4}, there exists $a\in\bO^\times$ such that
\[
\begin{pmatrix}a&a^{-1}-a\\0&a^{-1}\end{pmatrix}A\in t^lG_1(\bO).
\]
By writing down the upper-right entry, we obtain
\[
a t^{l-m}x+(a^{-1}-a)t^{-l}\in t^l\bO.
\]
If $2l<m$, then the first term has strictly smaller valuation than the other terms, a contradiction. If $2l=m$, then we obtain $x=1-a^{-2}$ modulo $t^m$. Conversely, if we multiply the inverse of the displayed matrix by $t^l$ and then by $\diag(a,a^{-1})$, we obtain $A$ with $x=1-a^{-2}$. Hence, the intersection \eqref{eq:g1-mixed-intersection} is isomorphic to $(\BG_m\setminus\{1\})\times\BA^{m-1}$. Note that the image torus in the projective groups has coordinate $s=a^2$, since $\diag(a,a^{-1})$ and its negative have the same image, and the cocharacter of this torus is $\frac12(1;1,-1,-1)=\alpha_2$. Thus, up to inversion, the monodromy of the local system around $x=1$ is given by
\[
b(\lambda,\alpha_2) =\xi^{4\lambda^{(2)}_1-4\lambda^{(2)}_2-4\lambda^{(2)}_3-12l}.
\]
Since $m=\ord(\xi^{12})$ and $2l=m$, we have $b(\lambda,\alpha_1+\alpha_2)=\xi^{-24l}=1$, and hence, if the local system is constant, then $b(\lambda,\alpha_1)=b(\lambda,\alpha_2)=1$.

Thus, the case \eqref{eq:g1-rank-one-difference} cannot happen, and hence, $\lambda^{(1)}=\theta^{(1)}$.
\end{proof}

\begin{rem}\label{rem:spherical-zastava}
For a mixed coroot, the calculations in the proofs of Propositions~\ref{prop:part-q} and~\ref{prop:part-p} show that \eqref{eq:g2-mixed-intersection} and \eqref{eq:g1-mixed-intersection} are empty if the upper inequalities in \eqref{eq:restricted-mixed-inequality} are strict. We expect that a similar result holds for a general Borel-free spherical homogeneous space $H\backslash G$, such that any spherical root is of type $T$. Namely, under the restricted condition with strict inequalities, we expect that the fiber of the corresponding rank-one Zastava space is empty.
\end{rem}

\subsection{Restricted coweights and the metaplectic Hecke action}

\subsubsection{The normalized Hecke action}

In this subsection, we study the metaplectic Hecke action on the Gaiotto category. Then, using the convolution with the Wakimoto objects $J_\gamma$, we construct objects whose images under $F_{\Fl}^{\Lus}$ are irreducible (Corollary~\ref{cor:iwahori-irreducible-image}). First, we recall the part of \cite[\S 9]{[GL]} which will be used below, where we replace the Whittaker category by the Gaiotto category.

Let $\Sph_{\cG^G}(G)$ be the $\cG^G$-twisted spherical Hecke category, and we denote by $G_{q}^\vee$ the metaplectic dual group. By the metaplectic geometric Satake equivalence, the heart of the spherical Hecke category is identified with the category of representations of $G_{q}^\vee$.

For finite-dimensional $V\in\Rep(G_{q}^\vee)^\heartsuit$, let $\Sat_{q}(V)$ be the corresponding spherical object with even coefficients, with the Satake tensor structure of \cite[Section 2.4]{[GL]}. If the weights of $V$ lie in one coset of the geometric coroot lattice, we define $\delta(V)=\delta(\eta)$ for any weight $\eta$ of $V$, which is well-defined since $\delta$ vanishes on geometric coroots. We consider the following Hecke action
\begin{equation}\label{eq:normalized-Gaiotto-Hecke-action}
\cF\star_{\rm H}\Sat_{q}(V) :=\bigl(\cF\star\Sat_{q}(V)\bigr)[-\delta(V)],
\end{equation}
where $\star$ is the geometric convolution. Note that, by the additivity of $\delta$, we obtain the associativity of this right Hecke action on $\SD_{\cG^G}^{H(\bF)_c^\omega,\psi}(\Gr_G^\omega)$.

Recall that $\Lambda^\sharp\subset\Lambda$ is the metaplectic lattice determined by the factorization gerbe. Note that $\Lambda^\sharp$ is also the weight lattice of $G_{q}^\vee$. We denote by
\[
(\Lambda^\sharp)^+ :=\Lambda^\sharp\cap\Lambda_G^+
\]
its dominant cone. For $\gamma\in(\Lambda^\sharp)^+$, we denote by $V^\gamma\in\Rep(G_{q}^\vee)$ the irreducible representation of highest weight $\gamma$.

For $\eta\in\Lambda^\sharp$, let $\operatorname{Tr}_\eta$ be the translation functor on factorization modules defined in \cite[\S\S 4.4, 5.2.1]{[GL]}. Namely, it is the pushforward under $D\mapsto D+\eta\cdot c$ on $\Conf_{\infty\cdot c}$, using the compatible identifications of $\cG^\Lambda$ under these translations, and these identifications use the same trivialization at $c$ as the metaplectic Satake equivalence, see \cite[\S 2.3.8]{[GL]}. In particular, it sends $\Delta_\lambda^{\sm,\fact}$ to $\Delta_{\lambda+\eta}^{\sm,\fact}$, and it is $t$-exact.

According to \cite[Theorem-Construction 16.2.4]{[GL]}, we have the following lemma.
\begin{lem}\label{lem:small-Gaiotto-Hecke-compatibility}
For $\cF\in\SD_{\cG^G}^{H(\bF)_c^\omega,\psi,\ren}(\Gr_G^\omega)$ and $V\in\Rep(G_{q}^\vee)^\heartsuit$, let $V=\bigoplus_{\eta\in\Lambda^\sharp}V_\eta$ be its weight decomposition. There is a canonical isomorphism of factorization modules
\begin{equation}\label{eq:small-Gaiotto-Hecke-compatibility}
F_{\Gr}^{\sm}\bigl(\cF\star_{\rm H}\Sat_{q}(V)\bigr) \simeq \bigoplus_{\eta\in\Lambda^\sharp} \Pi^{d(\eta)}V_\eta\otimes\operatorname{Tr}_\eta F_{\Gr}^{\sm}(\cF),
\end{equation}
compatible with tensor products in $V$. 
\end{lem}

\subsubsection{The Hecke action on IC sheaves}

The following statement is the Gaiotto analogue of \cite[Theorem 9.4.4]{[GL]}. The second claim should be understood as the geometric statement of the supergroup version of the Steinberg--Lusztig isomorphism, see Remark~\ref{rem:spherical-gaiotto-conjecture}.

\begin{thm}\label{thm:metaplectic-hecke-action}
\begin{enumerate}
\item \begin{samepage} The Hecke action of $\Rep(G_{q}^\vee)^\heartsuit$ on
\[
\SD_{\cG^G}^{H(\bF)_c^\omega,\psi,\ren}(\Gr_G^\omega)
\]
is $t$-exact.\par
\end{samepage}
\begin{samepage}
\noindent More precisely, for $V\in\Rep(G_{q}^\vee)^\heartsuit$, the functor
\[
\cF\longmapsto\cF\star_{\rm H}\Sat_{q}(V)
\]
is $t$-exact.\par
\end{samepage}

\item Let $\lambda$ be restricted dominant and relevant, and let $\gamma\in(\Lambda^\sharp)^+$. If $\lambda+\gamma$ is relevant, then
\[
\IC_{\Gr}^{\lambda} \star_{\rm H}\Sat_{q}(V^\gamma) \simeq \IC_{\Gr}^{\lambda+\gamma}.
\]
\end{enumerate}
\end{thm}

\begin{proof}
The proof of the $t$-exactness is the same as that of \cite[Proposition 9.3.2]{[GL]}. Namely, we compute the $!$-stalks of the convolution with $\Sat_{q}(V)$, stratify them by the $H(\bF)$-orbits $\sO_{\Gr}^{\nu}$ and the semi-infinite orbits $\sS_{\Gr}^{\eta}$, where $\nu\in\Lambda_G^{L,+}$ and $\eta\in\Lambda$, and then apply Corollary~\ref{lem:dimension-estimate-super} to estimate the degrees of these $!$-stalks. Here, we use the normalization \eqref{eq:normalized-Gaiotto-Hecke-action}. The other bound follows from the Verdier duality.

The second assertion is similar but more subtle. The same argument as in \cite[Theorem 9.4.4 and \S\S 9.4.6--9.4.8]{[GL]} with the dimension estimate reduces it to Theorem~\ref{thm:local-system-rigidity}. Namely, by Corollary~\ref{lem:dimension-estimate-super}, only the irreducible components of
\[
\sO_{\Gr}^\lambda\cap t^\theta\bigl(U(\bO)t^{-\gamma}G(\bO)/G(\bO)\bigr)
\]
of dimension $\tfrac12d(\lambda-\theta+\gamma)$ can contribute in degree zero. Note that these components are dense in critical-dimensional components of $\sO_{\Gr}^\lambda\cap\sS_{\Gr}^{\theta-\gamma}$, with the same local system up to a line. Thus, by Theorem~\ref{thm:local-system-rigidity}, we have $\theta=\lambda+\gamma$, and by Lemma~\ref{lem:equal-index-intersection}, we obtain the rank-one local system on $\sO_{\Gr}^{\lambda+\gamma}$. Finally, as in the proof of Theorem~\ref{thm:image-of-ic}, the duality excludes subobjects and quotients on the boundary, which proves (2).
\end{proof}

\subsubsection{The objects mapping to irreducibles}

Recall the projection $p\colon\Fl_G^\omega\to\Gr_G^\omega$ from \eqref{fl to gr}.

For the Wakimoto object $J_\gamma$ of \eqref{eq:twisted-bmw}, we consider the action defined by
\begin{equation}\label{eq:normalized-Wakimoto-action}
\cF\star_{\rm H}J_\gamma :=(\cF\star J_\gamma)[-\delta(\gamma)].
\end{equation}
Here, we use the normalization of Section~\ref{sec:gaiotto-standards} for $J_\gamma$. Then, for every Iwahori Gaiotto object $\cF$ and $\eta\in\Lambda^\sharp$, by the BMW relation \eqref{eq:bmw-convolution} and \eqref{5.1-colim-j!}, we have
\begin{equation}\label{eq:general-Wakimoto-compatibility}
F_{\Fl}^{\Lus}(\cF\star_{\rm H}J_\eta) \simeq\Pi^{d(\eta)}\operatorname{Tr}_\eta F_{\Fl}^{\Lus}(\cF).
\end{equation}

Now, suppose that we have chosen a decomposition
\[
\lambda =\lambda_{\res}+\gamma,
\]
where $\lambda_{\res}$ is restricted dominant and relevant, and $\gamma\in(\Lambda^\sharp)^+$ is dominant. We define an object of the Iwahori Gaiotto category by
\[
L_{\lambda}^{\Fl} :=p^!\bigl(\IC_{\Gr}^{\lambda_{\res}}\bigr) \star_{\rm H}J_\gamma,
\]
where $J_\gamma$ is the twisted Wakimoto sheaf attached to $\gamma$.

\begin{cor}\label{cor:iwahori-irreducible-image}
With the notation above,
\[
F_{\Fl}^{\Lus}\bigl(L_{\lambda}^{\Fl}\bigr) \simeq \Pi^{d(\lambda)}\IC_{\lambda}^{\Lus,\fact}.
\]
\end{cor}

\begin{proof}
Just note that, by \eqref{eq:general-Wakimoto-compatibility}, the convolution with $J_\gamma[-\delta(\gamma)]$ corresponds to $\Pi^{d(\gamma)}\operatorname{Tr}_\gamma$ on the configuration space. Thus, the assertion follows from \eqref{eq:flag-Grassmannian-functors}, Theorem~\ref{thm:image-of-ic}, and $d(\lambda)=d(\lambda_{\res})+d(\gamma)$.
\end{proof}

\subsubsection{Comparison with standard objects}

Combining Corollary~\ref{cor:iwahori-irreducible-image} with Proposition~\ref{5.0.1}, we can compare the $\RHom$ from standard objects to $L_\lambda^{\Fl}$ with the corresponding $\RHom$ in the category $\cO_{q}^{\mxd}$.

\begin{thm}\label{thm:iwahori-rhom-comparison}
We make the same assumptions as in Theorem~\ref{thm:geo-quant-comparison}. For any $\lambda$ and $L_\lambda^{\Fl}$ as in Corollary~\ref{cor:iwahori-irreducible-image} (here, $\lambda$ is regarded as a weight of the mixed quantum supergroup), and any weight $\nu$, the functor $F_{\Fl}^{\Lus}$ induces an isomorphism
\[
\RHom_{ \SD_{\cG^G}^{H(\bF)_c^\omega,\psi}(\Fl_G^\omega) } \bigl( \Delta_\nu,L_{\lambda}^{\Fl} \bigr) \overset{\sim}{\longrightarrow} \Pi^{d(\nu)+d(\lambda)}\RHom_{\cO_{q}^{\mxd}} \bigl( \Verma_\nu,\Irr_{\lambda} \bigr).
\]
\end{thm}

\begin{proof}
By Proposition~\ref{5.0.1} and Corollary~\ref{cor:iwahori-irreducible-image}, the left-hand side is $\Pi^{d(\nu)+d(\lambda)}i_\nu^!\IC_\lambda^{\Lus,\fact}$. Since the generator defining $\Delta_\nu^{\Lus,\fact}$ is even, it is isomorphic to
\[
\Pi^{d(\nu)+d(\lambda)}\RHom(\Delta_\nu^{\Lus,\fact},\IC_\lambda^{\Lus,\fact}).
\]
Then, the assertion follows, since the two objects correspond to $\Verma_\nu$ and $\Irr_\lambda$, respectively.
\end{proof}

\subsection{Behavior with respect to isogenies}
\label{sec:gaiotto-isogenies}

Note that the construction of $L_\lambda^{\Fl}$ uses a decomposition $\lambda=\lambda_{\res}+\gamma$, which may not exist for the original group $G$. In this subsection, following \cite[Sections 4.5 and 11.2--11.4]{[GL]}, we replace $G$ by a group $\widetilde G$ in the same isogeny class so that such decompositions exist, and compare the corresponding Hecke categories.

\subsubsection{The isogeny and the IC comparison}

Let
\begin{equation}\label{eq:gaiotto-isogeny}
1\longrightarrow G\longrightarrow\widetilde G \longrightarrow T_0\longrightarrow1
\end{equation}
be an exact sequence, where $T_0$ is a torus and $G$ is the group in one of the homogeneous data of Section~\ref{pol des}. We keep the subgroup $H\subset G\subset\widetilde G$ and its character $\psi$. Let $\widetilde T\subset\widetilde B$ be the maximal torus and the Borel subgroup containing $T\subset B$. We define the $\omega^\rho$-twist by the homomorphism $2\rho:\BG_m\to T\hookrightarrow\widetilde T$ of Section~\ref{sec:grading-cocharacter}, and we consider the homogeneous polarization $T_\psi^*(\widetilde G/H)$.

Since $\widetilde B\cap H=B\cap H=\{1\}$ and $\dim(\widetilde G/H)=\dim\widetilde B$, the homogeneous space $\widetilde G/H$ still has a free open $\widetilde B$-orbit. Let $\widetilde\Lambda=X_*(\widetilde T)$ and $\Lambda_0=X_*(T_0)$.

As in \cite[Section 11.4.4]{[GL]}, we choose a Zariski neighborhood $U\subset C$ of $c$, over which $\cG^G$ extends to a factorization gerbe $\cG^{\widetilde G}$ on $\Gr_{\widetilde G,\Ran_U}^{\omega}$, i.e., its pullback to $\Gr_{G,\Ran_U}^{\omega}$ is identified with the restriction of $\cG^G$. Moreover, its pullback to $H(\bF)_c^\omega$ has the multiplicative trivialization of Section~\ref{sec:factorization-gerbes}. As in Section~\ref{3.2.5 sect}, we denote by $\widetilde q$ and $\widetilde b$ the corresponding quadratic form and symmetric bilinear form.

Following \cite[Section 11.2.2]{[GL]}, we say that the isogeny \eqref{eq:gaiotto-isogeny} is strictly compatible with the geometric metaplectic data if
\begin{equation}\label{eq:gaiotto-isogeny-pushout}
\widetilde\Lambda =\Lambda\underset{\Lambda^\sharp}{\sqcup} \widetilde\Lambda^\sharp.
\end{equation}
Here, $\widetilde\Lambda^\sharp=\{\eta\in\widetilde\Lambda\mid\widetilde b(\eta,\widetilde\Lambda)=1\}$ is the metaplectic lattice. Note that the geometric coroots, the simple super roots and the orders $q_{a,i}=\ord(q(\alpha_{a,i}))$ are the same as for $G$ in Section~\ref{sec:irreducibles}, and the restricted coweights are defined as in Definition~\ref{def:restricted-coweight}, with $\lambda\in\widetilde\Lambda$ and $b$ replaced by $\widetilde b$.

The extension has natural quotients
\[
\widetilde G_1=\widetilde G/G_2,\qquad \widetilde G_2=\widetilde G/G_1, \qquad \widetilde G=\widetilde G_1\times_{T_0}\widetilde G_2.
\]
Let $\lambda=(\lambda^{(1)},\lambda^{(2)})$ be a coweight of $\widetilde G$, where $\lambda^{(a)}$ is its image in the coweight lattice of $\widetilde G_a$. The quotient maps $\widetilde G_a\to T_0$ send both components to the same $\xi\in\Lambda_0$, and we fix this $\xi$, which corresponds to $t^\xi\in\Gr_{T_0}$. Here, we consider the reduced ind-schemes, and then $(\Gr_{T_0})_{\mathrm{red}}=\Lambda_0$ by \cite[Remark~4.1.3]{[SW]}. We denote by $\widetilde{\bX}=\overline{H\backslash\widetilde G}^{\aff}$ the corresponding affine closure. For an object $(\cP_{\widetilde G},\sigma)$ of the global model in Section~\ref{sec:global-models}, the section is
\[
\sigma:C-\{c\}\longrightarrow \widetilde{\bX}_{\omega^\rho}\overset{\widetilde G}{\times} \cP_{\widetilde G}.
\]
Since $H\subset G$, the quotient map $H\backslash\widetilde G\to T_0$ extends to $\widetilde{\bX}\to T_0$. Since $2\rho$ factors through $G$, by composing $\sigma$ with this map, we obtain a section of $\cP_{T_0}:=\cP_{\widetilde G}\overset{\widetilde G}{\times}T_0$ over $C-\{c\}$, and hence a trivialization of $\cP_{T_0}$ there. We define the local functors for $\widetilde G$ by the correspondences and semi-infinite kernels of Section~\ref{sec:local-functors}, using $\cG^{\widetilde G}$ over $U$. For the construction over the configuration space of $C$, we use this datum near $c$ and the original datum for $G$ away from $c$, together with the factorization identifications of Section~\ref{3.3 fact prop}. Note that the factorization algebras are still $\Omega_q^{\Lus},\Omega_q^{\DK},\Omega_q^{\sm}$, since the configuration degrees away from $c$ belong to the same super root lattice $Q_s$ as in Section~\ref{super Borel}, and only the coweights at the marked point are enlarged.

In type A, we choose a section $T_0\to Z(\widetilde G)^0$ of $\widetilde G\to T_0$, which exists since $Z(G)$ is connected, see Lemma~\ref{lem:app-isogeny-components}, and we extend $d$ by zero on its coweight lattice. In the other cases, we let $d(\eta)=\langle2\check\rho_{\widetilde G},\eta\rangle$. Let
\begin{equation}\label{eq:app-isogeny-degree-correction}
\delta(\eta)=d(\eta) -\langle2\check\rho_{\widetilde G},\eta\rangle.
\end{equation}
Note that it agrees with the original correction on $\Lambda$, and it vanishes on geometric coroots. In type A, the translation by the lift of $t^\xi$ under the chosen section identifies $\Gr_G^\omega$ with the fiber of $\Gr_{\widetilde G}^\omega\to\Gr_{T_0}$ over $t^\xi$, and we use the corresponding local $t$-structure and the normalization of Remark~\ref{rem:component-normalization}.

The following two theorems are the analogues of Theorems~\ref{thm:image-of-ic} and \ref{thm:metaplectic-hecke-action}, and their proofs are given in Appendix~\ref{sec:gaiotto-isogeny-proofs}.

\begin{thmrestate}{thm:image-of-ic}\label{thm:isogeny-image-of-ic}
Suppose that the isogeny \eqref{eq:gaiotto-isogeny} is strictly compatible with the geometric metaplectic data. For every restricted dominant relevant coweight $\lambda\in\widetilde\Lambda$, we have
\[
F_{\Gr,\widetilde G}^{\sm}(\IC_{\Gr}^{\lambda}) \simeq\Pi^{d(\lambda)}\IC_{\lambda}^{\sm,\fact},\qquad F_{\Gr,\widetilde G}^{\Lus}(\IC_{\Gr}^{\lambda}) \simeq\Pi^{d(\lambda)}\IC_{\lambda}^{\Lus,\fact}.
\]
Moreover, if a component of $\sO_{\Gr}^{\lambda}\cap\sS_{\Gr}^{\theta}$ has dimension $\tfrac12d(\lambda-\theta)$ and $\cL_{\lambda,\theta}\otimes\psi$ is constant on a smooth dense open subset of that component, then $\theta=\lambda$.
\end{thmrestate}

\begin{thmrestate}{thm:metaplectic-hecke-action}\label{thm:isogeny-metaplectic-hecke-action}
Suppose that the isogeny \eqref{eq:gaiotto-isogeny} is strictly compatible with the geometric metaplectic data. The normalized Hecke action of $\Rep(\widetilde G_{\widetilde{q}}^\vee)^\heartsuit$ on
\[
\SD_{\cG^{\widetilde G}}^{H(\bF)_c^\omega,\psi,\ren} (\Gr_{\widetilde G}^\omega)
\]
is $t$-exact, with the correction $\delta$ of \eqref{eq:app-isogeny-degree-correction}. For $\lambda\in\widetilde\Lambda$ restricted dominant relevant, and $\gamma\in(\widetilde\Lambda^\sharp)^+$ such that $\lambda+\gamma$ is relevant, we have
\[
\IC_{\Gr}^{\lambda}\star_{\rm H} \Sat_{\widetilde{q}}(V^\gamma) \simeq\IC_{\Gr}^{\lambda+\gamma}.
\]
Here, $\Sat_{\widetilde{q}}$ is the Satake functor for $\widetilde G$, and $V^\gamma$ has highest weight $\gamma$.
\end{thmrestate}
\subsubsection{Construction of the isogeny}

First of all, we need the existence of a suitable isogeny.

\begin{prop}\label{prop:gaiotto-isogeny-existence}
For every homogeneous datum considered in this paper, there is an isogeny \eqref{eq:gaiotto-isogeny}, strictly compatible with the geometric metaplectic data, for which the metaplectic dual of $\widetilde G$ has simply connected derived group. Here, the extension of the geometric metaplectic datum is taken on an affine neighborhood of $c$, as in \cite[Section 11.4]{[GL]}.
\end{prop}
\begin{proof}
The existence of $\widetilde G$ and $\cG^{\widetilde G}$ follows from \cite[Proposition 11.3.6 and Section 11.4]{[GL]}, applied to the geometric metaplectic datum of $G$.

Moreover, to preserve Condition~\ref{cond:star}, we choose a section $s:\Lambda_0\to\widetilde\Lambda^\sharp$, which exists since $\Lambda_0$ is a free abelian group. Then, we have $\widetilde\Lambda=\Lambda\oplus s(\Lambda_0)$, and since $s(\Lambda_0)\subset\widetilde\Lambda^\sharp$, we can extend $b'$ to a bilinear form $\widetilde b'$ such that
\[
\widetilde b'(\nu,\nu)=\widetilde q(\nu),\qquad \widetilde b'(\lambda,s(\xi)) =\widetilde b'(s(\xi),\lambda)=1,
\]
for $\nu\in\widetilde\Lambda$, $\lambda\in\Lambda$ and $\xi\in\Lambda_0$. Since the super roots remain in $\Lambda$, Condition~\ref{cond:star} is preserved.
\end{proof}

\subsubsection{Restricted coweights}

By \cite[Lemma 9.4.2]{[GL]}, we obtain the following lemma. Here, for a mixed simple coroot $\alpha_{a,i}$ with $q_{a,i}<\infty$, if $\langle\lambda,\check\alpha_{a,i}\rangle$ is a positive multiple of $q_{a,i}$ and $\widetilde b(\lambda,\alpha_{a,i}^{\pm})\neq1$, then we choose the remainder to be $q_{a,i}$. In this way, we obtain the restricted and relevant conditions of Definition~\ref{def:restricted-coweight} and Proposition~\ref{prop:relevant-L-orbits}.

\begin{lem}\label{lem:restricted-relevant-decomposition}
Assume that the metaplectic dual of $\widetilde G$ has simply connected derived group. Every dominant relevant $\lambda\in\widetilde\Lambda$ can be written as
\[
\lambda=\lambda_{\res}+\gamma,
\]
where $\lambda_{\res}$ is restricted dominant and relevant and $\gamma\in(\widetilde\Lambda^\sharp)^+$.
\end{lem}

In particular, under the same assumptions, Corollary~\ref{cor:iwahori-irreducible-image} and Theorem~\ref{thm:iwahori-rhom-comparison} hold with $G,\Lambda$ replaced by $\widetilde G,\widetilde\Lambda$, for every dominant relevant $\lambda\in\widetilde\Lambda$.

\subsubsection{The Hecke categories}

The following proposition is from \cite[Proposition 11.2.5 and Corollary 11.2.6]{[GL]}.
\begin{prop}\label{prop:gaiotto-hecke-isogeny}
Assume that the isogeny \eqref{eq:gaiotto-isogeny} is strictly compatible with the geometric metaplectic data, and denote by $K$ and $\widetilde K$ the two metaplectic dual groups, with Cartan subgroups $T_K$ and $T_{\widetilde K}$, and by $\Sat_{\widetilde{q}}$ the Satake functor for $\widetilde G$. We denote
\[
\cC_G=\SD_{\cG^G}^{H^\omega(\bF),\psi}(\Gr_G^\omega), \qquad \cC_{\widetilde G} =\SD_{\cG^{\widetilde G}}^{H^\omega(\bF),\psi} (\Gr_{\widetilde G}^\omega).
\]
We consider the normalized right Satake actions, with the function $\delta$ fixed in \eqref{eq:app-isogeny-degree-correction}, which extends the correction in \eqref{eq:normalized-Gaiotto-Hecke-action} and vanishes on geometric coroots. Let
\[
\Hec(\cC_G) :=\cC_G\underset{\Rep(K)}\otimes \Vect, \qquad \DotHec(\cC_G) :=\cC_G\underset{\Rep(K)}\otimes \Rep(T_K),
\]
and similarly for $\widetilde G$. Let $i:\Gr_G^\omega\to\Gr_{\widetilde G}^\omega$ be the natural inclusion. The functor $\Rep(T_K)\to\Rep(T_{\widetilde K})$ is induced by $\Lambda^\sharp\hookrightarrow\widetilde\Lambda^\sharp$. Then, $i_*$ induces equivalences
\[
\Hec(\cC_G) \xrightarrow{\sim}\Hec(\cC_{\widetilde G})
\]
and
\[
\Rep(T_{\widetilde K}) \underset{\Rep(T_K)}\otimes \DotHec(\cC_G) \xrightarrow{\sim} \DotHec(\cC_{\widetilde G}).
\]
\end{prop}

\subsection{The spherical conjecture and the graded Hecke category}

\subsubsection{The spherical conjecture}

\begin{rem}\label{rem:spherical-gaiotto-conjecture}
There is an expected equivalence between the renormalized spherical Gaiotto category and the category of representations of the Lusztig quantum supergroup,
\begin{equation}\label{eq:spherical-gaiotto}
\SD_{\cG^G}^{H^\omega(\bF),\psi,\ren}(\Gr_G^\omega) \simeq\Rep_q^{\Lus}(G_s).
\end{equation}
Here, $\Rep_q^{\Lus}(G_s)=\Rep(U_q^{\Lus}(\fg_s))$ denotes $\Ind D^b(\Rep_q^{\Lus}(G_s)^{\mathrm{fd}})$, where $\mathrm{fd}$ means finite-dimensional representations with the prescribed $\Lambda$-weight decomposition.

However, the authors do not know the Lusztig quantum supergroup in the generality needed here. It is expected that there is a Lusztig quantum supergroup $U_{q}^{\Lus}(G_s)$, which is a Hopf superalgebra with the following properties. Let $\mathfrak g_{q}^\vee=\Lie(G_{q}^\vee)$, and let $T_{q}^\vee$ be its Cartan torus, whose character lattice is $\Lambda^\sharp$. We denote by $\Fr_{q,T}:T_s\to T_{q}^\vee$ the map dual to $\Lambda^\sharp\hookrightarrow\Lambda$, whose composite with $T_{q}^\vee\to T_s^\#$ is $\Fr_{q,T_s}$.

\begin{enumerate}
\item For any choice of super Borel, there is a Hopf sub-superalgebra $U_{q}^{\Lus}(B_s)$, which is the graded Hopf dual of $U_{q^{-1}}^{\DK}(B_s^-)$.

\item The Lusztig quantum supergroup should have a quantum Frobenius map
\[
\Fr_{q} \colon U_{q}^{\Lus}(\fg_s) \longrightarrow U(\mathfrak g_{q}^\vee),
\]
which is a super-cocentral map of Hopf superalgebras and fits into a short exact sequence of Hopf superalgebras
\begin{equation}\label{eq:whole-lusztig-frobenius-sequence}
u_{q}^{\sm}(G_s) \longrightarrow U_{q}^{\Lus}(G_s) \overset{\Fr_{q}}{\longrightarrow} \Fr_{q,T}^*U(\mathfrak g_{q}^\vee).
\end{equation}
Its restriction to either half coincides with \eqref{short ex seq Lus}.

\item The irreducible representations are parameterized by the same set as the simple objects in the $q$-twisted Gaiotto category, cf., \cite[Theorem 2]{Po} for the comparison with characteristic-$p$ representations in type A, when $M<N$ and $q^p=1$ with $p$ prime.

\item For $\lambda'$ relevant restricted dominant and $\lambda''\in(\Lambda^\sharp)^+$, the irreducible representation $\Irr_{\lambda'}^{\Lus}$ is also irreducible for the small quantum supergroup and the De Concini--Kac quantum supergroup, and there exists a Steinberg--Lusztig isomorphism
\begin{equation}\label{eq:super-steinberg-lusztig}
\Irr_{\lambda'+\lambda''}^{\Lus} \simeq \Irr_{\lambda'}^{\Lus} \otimes \Fr_{q}^*\bigl(V^{\lambda''}\bigr).
\end{equation}
\end{enumerate}

The conjectural equivalence \eqref{eq:spherical-gaiotto} should be $t$-exact, and should be compatible with the metaplectic Satake action in \eqref{eq:normalized-Gaiotto-Hecke-action}. At the heart level, the action of $\Rep(G_{q}^\vee)$ on $\Rep_q^{\Lus}(G_s)$ is given by tensoring with the pullback along $\Fr_{q}$, and the action of $\Sph_{\cG^G}(G)$ on $\SD_{\cG^G}^{H^\omega(\bF),\psi,\ren}(\Gr_G^\omega)$ is given by the normalized convolution \eqref{eq:normalized-Gaiotto-Hecke-action}. That explains why Theorem~\ref{thm:metaplectic-hecke-action} is the geometric version of the Steinberg--Lusztig isomorphism.

For $q$ generic, the equivalences \eqref{eq:spherical-gaiotto} in type A and in the $\on F(4)$ case have already been established in \cite{[G],[BFT1],[TY],[FTY]}, and for $\SOSp(2k|2l)$ with $k=l$ or $k=l+1$, the same proof as in \cite{[BFT1]} also works, see \cite[Section 1.5]{[BFT1]}. For the other orthosymplectic cases, the equivalence can be proved by the same method as in \cite{[TY]}. For the root of unity case, it is harder, and one can consider a weaker conjecture, i.e., an isomorphism of the corresponding Grothendieck groups.

\smallskip
\subsubsection{Ext multiplicities} Once there is a definition of the Lusztig quantum supergroup with the expected properties, the results of this section could be applied to compare the transition matrix coefficients of the canonical basis and the costandard basis on both sides, namely, to compare the Ext multiplicities in the Grothendieck groups:
\begin{equation}\label{eq:spherical-ext-comparison}
\begin{split}
&\bigl[ \RHom_{ \SD_{\cG^G}^{H^\omega(\bF),\psi}(\Gr_G^\omega) } (\Delta_{\nu,\Gr},\IC_{\Gr}^{\lambda}) \bigr]\\
&\qquad= \bigl[ \Pi^{d(\nu)+d(\lambda)}\RHom_{ \Rep_{q}^{\Lus}(G_s) } (\on{Kac}_\nu,\Irr_{\lambda}^{\Lus}) \bigr].
\end{split}
\end{equation}
Here, $\Delta_{\nu,\Fl}:=\Delta_\nu$ and $\Delta_{\nu,\Gr}:=p_*\Delta_{\nu,\Fl}$, with the induced normalization. For $B_s$, we take the same $\Ind D^b(\Rep^{\mathrm{fd}})$ convention. We assume that the restriction from $G_s$ to $B_s$ has a derived left adjoint $\Ind^L$ on the characters $\se_\nu$, with bounded finite-dimensional values. Then, with $\se_\nu$ as in Remark~\ref{rem:target-parity}, we let $\on{Kac}_\nu:=\Ind^L(\se_\nu)$, so that
\[
\RHom_{\Rep_{q}^{\Lus}(G_s)}(\on{Kac}_\nu,V) \simeq\RHom_{\Rep_{q}^{\Lus}(B_s)} (\se_\nu,\operatorname{Res}V).
\]
We call $\on{Kac}_\nu$ the Kac complex.\footnote{It need not belong to the heart for a non-standard super Borel, just as $\Delta_{\nu,\Gr}$ need not be perverse.}

By the adjunction for $p$, we have
\begin{equation}\label{eq:projection-formula-gr-fl}
\begin{split}
&\bigl[ \RHom_{ \SD_{\cG^G}^{H^\omega(\bF),\psi}(\Gr_G^\omega) } (\Delta_{\nu,\Gr},\IC_{\Gr}^{\lambda}) \bigr]\\
&\qquad= \bigl[ \RHom_{ \SD_{\cG^G}^{H^\omega(\bF),\psi}(\Fl_G^\omega) } (\Delta_{\nu,\Fl},p^!\IC_{\Gr}^{\lambda}) \bigr].
\end{split}
\end{equation}
Suppose that we have a decomposition
\[
\lambda =\lambda_{\res}+\gamma,
\]
where $\lambda_{\res}$ is relevant restricted dominant and $\gamma\in(\Lambda^\sharp)^+$. Then, according to Theorem~\ref{thm:metaplectic-hecke-action}, we have
\[
\IC_{\Gr}^{\lambda} \simeq \IC_{\Gr}^{\lambda_{\res}} \star_{\rm H}\Sat_{q}(V^{\gamma}).
\]
We denote by $Z_{q}:\Rep(G_{q}^\vee) \to\SD_{\cG^G}(I\backslash G(\bF)/I)$ the metaplectic version of Gaitsgory's central sheaf construction. Note that the pullback of the convolution with the Satake object is given by the convolution with the corresponding central sheaf, i.e.,
\[
p^!\bigl(\cF\star_{\rm H}\Sat_{q}(V)\bigr) \simeq p^!\cF\star_{\rm H}Z_{q}(V),
\]
and, for finite-dimensional $V$, $Z_{q}(V)$ has a finite filtration by the Wakimoto objects with associated graded $\bigoplus_\eta V_\eta\otimes J_\eta$, namely, the multiplicities are given by the dimensions of the corresponding weight spaces. For every weight $\gamma^{\prime}$ of $V^\gamma$, we have $\delta(\gamma^{\prime})=\delta(\gamma)$, since $\delta$ vanishes on geometric coroots, thus, the cohomological shifts in \eqref{eq:normalized-Gaiotto-Hecke-action} and \eqref{eq:normalized-Wakimoto-action} agree. By Theorem~\ref{thm:metaplectic-hecke-action}(2), the Wakimoto filtration, and the additivity of Grothendieck classes, we obtain
\begin{equation}\label{eq:central-sheaf-ext-chain}
\begin{split}
&\bigl[ \RHom_{\SD_{\cG^G}^{H^\omega(\bF),\psi}(\Fl_G^\omega)} (\Delta_{\nu,\Fl},p^!\IC_{\Gr}^{\lambda}) \bigr]\\
&\qquad= \sum_{\gamma'}\dim V^\gamma_{\gamma'} \bigl[ \RHom_{\SD_{\cG^G}^{H^\omega(\bF),\psi}(\Fl_G^\omega)} (\Delta_{\nu,\Fl},p^!\IC_{\Gr}^{\lambda_{\res}} \star_{\rm H}J_{\gamma'}) \bigr].
\end{split}
\end{equation}

For every weight $\gamma'$ of $V^\gamma$, we apply \eqref{eq:general-Wakimoto-compatibility}, \eqref{eq:flag-Grassmannian-functors} and Theorem~\ref{thm:image-of-ic}, and note that the translation by $\gamma'$ on the factorization side corresponds to tensoring by the line $\se_{\gamma'}$ on the representation side. Thus, by Proposition~\ref{5.0.1}, as in the proof of Theorem~\ref{thm:iwahori-rhom-comparison}, we have
\[
\begin{split}
&\RHom_{\SD_{\cG^G}^{H^\omega(\bF),\psi}(\Fl_G^\omega)} (\Delta_{\nu,\Fl},p^!\IC_{\Gr}^{\lambda_{\res}} \star_{\rm H}J_{\gamma'})\\
&\qquad\simeq \Pi^{d(\nu)+d(\lambda_{\res})+d(\gamma')} \RHom_{\cO_{q}^{\mxd}} (\Verma_\nu,\Irr_{\lambda_{\res}}\otimes\se_{\gamma'}).
\end{split}
\]
Also, we have $d(\gamma')\equiv d(\gamma)\pmod2$, since the difference of two weights of $V^\gamma$ lies in the geometric coroot lattice, and hence, the parity factor is $\Pi^{d(\nu)+d(\lambda)}$ in every summand.

Finally, by the adjunction for the Verma module, the Borel-stable filtration of $V^\gamma$, the adjunction for $\on{Kac}_\nu$, and the assumed Steinberg--Lusztig isomorphism \eqref{eq:super-steinberg-lusztig}, we obtain
\[
\begin{split}
&\sum_{\gamma'}\dim V^\gamma_{\gamma'} \bigl[\Pi^{d(\nu)+d(\lambda)} \RHom_{\cO_{q}^{\mxd}} (\Verma_\nu,\Irr_{\lambda_{\res}}\otimes\se_{\gamma'})\bigr]\\
&\quad= \bigl[\Pi^{d(\nu)+d(\lambda)} \RHom_{\Rep_{q}^{\Lus}(B_s)} (\se_\nu,\operatorname{Res}(\Irr_{\lambda_{\res}}^{\Lus}\otimes\Fr_{q}^*V^\gamma))\bigr]\\
&\quad= \bigl[\Pi^{d(\nu)+d(\lambda)} \RHom_{\Rep_{q}^{\Lus}(G_s)} (\on{Kac}_\nu,\Irr_{\lambda}^{\Lus})\bigr].
\end{split}
\]
\end{rem}

\subsubsection{The graded Hecke category}
\begin{rem}\label{rem:graded-small-quantum-group}
For completeness, we remark that the same proof (via the highest weight category argument) as in \cite[Section 19.3]{[GL]}, together with the results established in this paper, in particular Theorem~\ref{thm:metaplectic-hecke-action}, would imply that there is a $t$-exact equivalence of highest weight categories between
\[
\DotHec\bigl(\SD_{\cG^G}^{H^\omega(\bF),\psi,\ren}(\Gr_G^\omega)\bigr):=\SD_{\cG^G}^{H^\omega(\bF),\psi,\ren}(\Gr_G^\omega)\underset{\Rep(G_{q}^\vee)}{\otimes}\Rep(T_{q}^\vee)
\]
and the category $\Rep_q^\sm(\fg_s)$ of representations of the graded small quantum supergroup, defined in Section~\ref{sec:mixed-quantum-supergroups}. Here, the renormalization of both sides is given as follows: the geometric category is the ind-completion of the (non-complete) category of locally compact objects in $\SD_{\cG^G}^{H^\omega(\bF),\psi}(\Gr_G^\omega)$, and the representation side is the ind-completion of the (non-complete) category of finite-dimensional representations of the graded small quantum supergroup. The generic case follows as in Remark~\ref{rem:spherical-gaiotto-conjecture}, and we only need to consider the root of unity case. On factorization modules, the character of weight $\eta\in\Lambda^\sharp$ acts by $\Pi^{d(\eta)}\operatorname{Tr}_\eta$, as in Lemma~\ref{lem:small-Gaiotto-Hecke-compatibility}.

Note that, by Proposition~\ref{prop:gaiotto-isogeny-existence}, we may assume that the metaplectic dual has simply connected derived group. In this case, the IC comparison, the IC convolution formula, and the decomposition $\lambda=\lambda_{\res}+\gamma$ are given by Theorems~\ref{thm:isogeny-image-of-ic}, \ref{thm:isogeny-metaplectic-hecke-action} and Lemma~\ref{lem:restricted-relevant-decomposition}, respectively. Then, we recover the assertion for $G$ by Proposition~\ref{prop:gaiotto-hecke-isogeny} and \cite[Section 11.2]{[GL]}, by restricting to $\Lambda_0$-degree zero.

We follow the notation of \cite[\S\S 23--24 and 26.2--26.4]{[GL]}. For dominant $\lambda$, we define
\[
W_{\Fl}^{\lambda,!}:=\Delta_\lambda,\qquad W_{\Fl}^{\lambda,*}:=\nabla_\lambda^{\stan}[-2r],\qquad r=\dim(G/B),
\]
with the corresponding twisting and character, and let $W^{\lambda,?}=p_*W_{\Fl}^{\lambda,?}$. Here, the shifts come from the normalizations in Theorem~\ref{prop:global-local-strata} and \eqref{eq:bmw-normalized-duality}. We denote by $\overset{\bullet}{\cF}_{\infty/2}$ and $\overset{\bullet}{\IC}_{\infty/2}$ the graded Hecke liftings of the Iwahori and semi-infinite IC kernels in loc. cit. Note that, in this remark, the BMW kernels have the inverse gerbe and even coefficients. The right Hecke action on the kernels is given by
\[
K\longmapsto(K\star\Sat_{q^{-1}}(V))[\delta(V)],
\]
which is dual to \eqref{eq:normalized-Gaiotto-Hecke-action}. We consider the pairings of \cite[Sections 26.2--26.4]{[GL]}, which are given on locally compact objects by
\[
\Psi(\cF,K)=\Gamma(\Gr_G^\omega,\cF\overset!\otimes K),\qquad\Upsilon(\cF,\cF')=\RHom(\BD_{\Gr}\cF',\cF),
\]
and their graded Hecke versions and continuous extensions to the renormalized categories. Note that, in $\Psi$, the kernel has the inverse gerbe, and in $\Upsilon$, the Gaiotto categories have opposite gerbes and additive characters. Since $\delta$ vanishes on the geometric coroot lattice $Q_G$, we define
\[
(\mathcal T_\delta K)|_{(\Gr_G^\omega)^\gamma}:=K|_{(\Gr_G^\omega)^\gamma}[-2\delta(\gamma)],\qquad\gamma\in\Lambda/Q_G.
\]
Here, $(\Gr_G^\omega)^\gamma$ is the connected component indexed by $\gamma\in\Lambda/Q_G$. Note that $\mathcal T_\delta$ changes $[\delta(V)]$ in the Hecke action on the kernels into $[-\delta(V)]$.

The same strategy as in \cite{[GL]} works in the Gaiotto case, namely, we construct the costandard object
\[
\overset{\bullet}{M}_{\lambda}^*:=W_{\Fl}^{\lambda,*}\star_I\mathcal T_\delta(\overset{\bullet}{\cF}_{\infty/2}).
\]
For an arbitrary $\lambda$, we choose $\eta\in(Q_G\cap\Lambda^\sharp)^+$ such that $\lambda+\eta$ is dominant, and define $\overset{\bullet}{M}_{\lambda}^*$ by tensoring $\overset{\bullet}{M}_{\lambda+\eta}^*$ with the character of weight $-\eta$, as in \cite[Section 25.1.6]{[GL]}. Then, we define the standard object $\overset{\bullet}{M}_{\lambda}^!$ by applying the duality functor (change $q$ and $\psi$ to $q^{-1}$ and $-\psi$), and verify the orthogonality, the finiteness of the filtrations, and the fact that the functor preserves standards and costandards, as in \cite[Sections 25--26]{[GL]}. We only point out the following warnings.
\begin{enumerate}
\item For us, since we take $H$ transversal to $B$, the dominance means dominance with respect to $B^-$, as fixed in Section~\ref{trans D}, while in \cite{[GL]}, the $U$ used in the Whittaker condition is the same as the unipotent radical of the Borel used in the Iwahori subgroup.

\item The objects $W^{\lambda,!}$ and $W^{\lambda,*}$ used in the proof are obtained by pushforward along $p:\Fl_G^\omega\to\Gr_G^\omega$, and in the Whittaker case, they are just the $!$- and $*$-extensions of the rank-one equivariant local system on the corresponding orbit, which is used in the proof of orthogonality. In the Gaiotto case, it is no longer true in general that they are extensions of a rank-one local system. However, it does not harm the argument, namely, since one needs to consider limits and colimits, we can assume that $\lambda$ is dominant enough. For such $\lambda$, the same root-subgroup calculation as in the proof of Proposition~\ref{4.2.1} identifies the fiber of $\sO_{\Fl}^{\lambda}\to\sO_{\Gr}^{\lambda}$ with the unipotent radical $U^-$ of $B^-$, the restrictions of $W^{\lambda,!}$ and $W^{\lambda,*}$ to $\sO_{\Gr}^{\lambda}$ are of rank one, and also the stabilizer of $t^\lambda\in\Gr_G^\omega$ is pro-unipotent. Thus, the same orthogonality calculation and finite filtration argument in \cite[Section 25]{[GL]} apply, together with Theorem~\ref{thm:metaplectic-hecke-action}.

\item First, we take $\lambda$ dominant, and the general case follows from the character translation above. In the proof of the fact that $\overset{\bullet}{M}_{\lambda}^!$ corepresents $\Pi^{d(\lambda)}i_\lambda^!F_{\Gr}^{\sm}$, one needs to compare
\begin{equation}\label{eq:compare-two-psi}
\overset{\bullet}{\Psi}\bigl(-,t^\lambda\cdot\overset{\bullet}{\cF}_{\infty/2}\bigr)=\overset{\bullet}{\Psi}\bigl(-,t^\lambda\cdot\overset{\bullet}{\IC}_{\infty/2}\bigr).
\end{equation}
In \cite[Section 26.5]{[GL]}, this comparison is achieved by the fact that the Whittaker equivariance implies $U(\bO)$-equivariance. In the Gaiotto case, one does not have $U(\bO)$-equivariance anymore. Instead, for dominant $\lambda$ and $\lambda'$ with $\lambda-\lambda'\in Q_G\cap\Lambda^\sharp$, we have the convolution formula
\begin{equation}\label{eq:wakimoto-translation-property}
W_{\Fl}^{\lambda,*}=W_{\Fl}^{\lambda',*}\star_I J_{\lambda-\lambda'}^{\BD},
\end{equation}
which follows from \eqref{eq:normalized-standard-convolution} and \eqref{eq:bmw-normalized-duality}. Note that the functor $\mathcal T_\delta$ commutes with these BMW convolutions. As a result, by the averaging adjunction as in \cite[Lemma 26.4.3 and Corollary 26.4.4]{[GL]}, we have
\begin{equation}\label{eq:psi-translation-property}
\begin{split}
&\overset{\bullet}{\Psi}\bigl(-,t^\lambda\cF\bigr)[d(\lambda)]\\
&\quad=\overset{\bullet}{\Upsilon}\bigl(-,W_{\Fl}^{\lambda,*}\star_I\mathcal T_\delta(\cF)\bigr)\\
&\quad=\overset{\bullet}{\Upsilon}\bigl(-,W_{\Fl}^{\lambda',*}\star_I J_{\lambda-\lambda'}^{\BD}\star_I\mathcal T_\delta(\cF)\bigr)\\
&\quad=\overset{\bullet}{\Psi}\bigl(-,t^{\lambda'}\cdot(J_{\lambda-\lambda'}^{\BD}\star_I\cF)\bigr)[d(\lambda')].
\end{split}
\end{equation}
Now, we take the colimit over $\eta\in(Q_G\cap\Lambda^\sharp)^+$, ordered by dominance, and denote $U_\eta=t^\eta U(\bO)t^{-\eta}$. Note that these subgroups increase with $\eta$, and their union is $U(\bF)$. Also, on this lattice, $\delta=0$ and $d=\deg$ is even. Thus, we have
\begin{equation}\label{eq:averaging-colimit-comparison}
\begin{split}
&\overset{\bullet}{\Psi}\bigl(-,t^\lambda\cdot\Av_!^{U(\bF)}(\overset{\bullet}{\cF}_{\infty/2})\bigr)\\
&\quad=\colim_{\eta\in(Q_G\cap\Lambda^\sharp)^+}\overset{\bullet}{\Psi}\bigl(-,t^\lambda\cdot\Av_!^{U_\eta}(\overset{\bullet}{\cF}_{\infty/2})\bigr)\\
&\quad=\colim_{\eta\in(Q_G\cap\Lambda^\sharp)^+}\overset{\bullet}{\Psi}\bigl(-,t^{\lambda+\eta}\cdot(J_{-\eta}^{\BD}\star_I\overset{\bullet}{\cF}_{\infty/2})\bigr)[d(\eta)]\\
&\quad=\colim_{\eta\in(Q_G\cap\Lambda^\sharp)^+}\overset{\bullet}{\Psi}\bigl(-,t^\lambda\cdot\overset{\bullet}{\cF}_{\infty/2}\bigr)\\
&\quad=\overset{\bullet}{\Psi}\bigl(-,t^\lambda\cdot\overset{\bullet}{\cF}_{\infty/2}\bigr).
\end{split}
\end{equation}
Here, the second equality is the expression of the averaging by the Wakimoto convolution, and the third one follows from \eqref{eq:psi-translation-property}. Note that the transition maps are the identity under these identifications. Indeed, by adjunction, they are given by the maps between the two costandard objects in \eqref{eq:wakimoto-translation-property}, which are the identity on the defining local system. Furthermore, these maps commute with the right Hecke action.

By \cite[Corollary 24.3.4]{[GL]}, we have $\Av_!^{U(\bF)}\overset{\bullet}{\cF}_{\infty/2}\simeq\overset{\bullet}{\IC}_{\infty/2}$, which proves \eqref{eq:compare-two-psi}. Together with \eqref{eq:psi-translation-property}, it identifies $\RHom(\overset{\bullet}{M}_{\lambda}^!,-)$ with $\Pi^{d(\lambda)}i_\lambda^!F_{\Gr}^{\sm}$, thus, $\overset{\bullet}{M}_{\lambda}^!$ corresponds to $\Pi^{d(\lambda)}\Delta_\lambda^{\sm,\fact}$. By duality, the same argument applies to costandard objects. Note that the proof only relies on \eqref{eq:wakimoto-translation-property}, which is always true for spherical varieties, in place of the $U(\bO)$-equivariance, and it is completed as in \cite[Section 19.3]{[GL]}.
\end{enumerate}
\end{rem}

\section{Grothendieck groups and geometric \texorpdfstring{$!$}{!}-stalks}\label{sec:type-a-conventions} In the following sections, we are going to apply the results in the former part of this paper to study the category $\cO_{q}^{\mxd}$ of the mixed quantum supergroup in type A. The main goal is to give an explicit algorithm to calculate the coefficients of the Kazhdan--Lusztig polynomials, namely, the multiplicities of the derived Hom from the Verma module $\Verma_\mu$ to the irreducible module $\Irr_\lambda$. In this section, we first reduce this problem to the calculation of geometric $!$-stalks in the Grothendieck group, and then we fix the conventions on finite fields, twistings and Grothendieck groups. In Section~\ref{sec:finite-type-monodromic}, we study the finite type monodromic categories, and in Section~\ref{aff kl}, we use them to determine the canonical basis in the affine case.

In the finite and affine type A calculations, we denote by $D$ the sheaf category with coefficients in complexes of vector spaces. With super coefficients, the coefficient ring is extended by $\BZ[\pi]/(\pi^2-1)$, where $\pi$ is the class of the parity shift of the unit. Thus, we have $[K[1]]=-[K]$ and $[\Pi K]=\pi[K]$, and the action formulas extend by linearity. We take $K_0$ of bounded constructible Weil objects, with the convention on locally compact objects in Section~\ref{Gaiotto category}. For $\GL(M|N)$ with $M=N$, the possibly infinite standard expansions of these objects are interpreted in the completed module of Section~\ref{sec:balanced-affine-algorithm}.

\subsection{Reduction to geometric \texorpdfstring{$!$}{!}-stalks} In this subsection, we assume that $q$ is torsion valued on $Q_s$. Then, $\Lambda^\sharp$ contains a strictly dominant coweight. We fix a decomposition $\lambda=\lambda_{\res}+\gamma$ as in Corollary~\ref{cor:iwahori-irreducible-image}, where $\lambda_{\res}$ is restricted dominant and relevant, and $\gamma\in(\Lambda^\sharp)^+$ is dominant, and we let $\mu$ be an arbitrary weight.

\begin{rem}
In Sections \ref{sec:finite-type-monodromic} and \ref{aff kl}, we work with type A, where the metaplectic dual group has simply connected derived group. Thus, every dominant relevant $\lambda$ has a decomposition $\lambda=\lambda_{\res}+\gamma$ as above, and the results of Section~\ref{sec:irreducibles} suffice. In the other types, the metaplectic dual group need not have simply connected derived group, and we should first pass to the isogeny of Proposition~\ref{prop:gaiotto-isogeny-existence}, and then apply Lemma~\ref{lem:restricted-relevant-decomposition} and Corollary~\ref{cor:iwahori-irreducible-image} to $\widetilde G$, as in Section~\ref{sec:gaiotto-isogenies}.
\end{rem}

According to Theorem~\ref{thm:iwahori-rhom-comparison}, it is enough to calculate the multiplicity of the derived Hom from the standard object $\Delta_\mu$ with the index $\mu$ to the object $L_\lambda$. Moreover, for these indices, by Proposition~\ref{prop:ICfact-costalk-dim-spread}, the dimensions of its cohomology are identified with the corresponding $\Ext$ dimensions on the representation side, after reduction to sufficiently large positive characteristic.

Here, $L_\lambda:=L_\lambda^{\Fl}$ for the chosen decomposition, and for its translates by $\mu'\in\Lambda^\sharp$, we denote $L_{\lambda+\mu'}:=L_\lambda\star_{\rm H}J_{\mu'}$. Note that, if $\gamma+\mu'$ is dominant, then it is the same construction with $\gamma\in\Lambda^\sharp$ replaced by $\gamma+\mu'$.

For any $\mu'\in \Lambda^\sharp$, right convolving with $J_{\mu'}$, with the normalization above, gives an equivalence
\begin{equation}
\star_{\rm H} J_{\mu'}\colon \SD_{\kappa}^{H^\omega(\bF),\psi}(\Fl^\omega_G)\overset{\sim}{\rightarrow}\SD_{\kappa}^{H^\omega(\bF),\psi}(\Fl^\omega_G).
\end{equation}

In particular, we have
\begin{multline}\label{8.2}
\RHom_{\SD_{\kappa}^{H^\omega(\bF),\psi}(\Fl^\omega_G)}(\Delta_\mu, L_\lambda)\simeq \RHom_{\SD_{\kappa}^{H^\omega(\bF),\psi}(\Fl^\omega_G)}(\Delta_{\mu}\star_{\rm H} J_{\mu'}, L_{\lambda}\star_{\rm H} J_{\mu'})\\
\simeq \RHom_{\SD_{\kappa}^{H^\omega(\bF),\psi}(\Fl^\omega_G)}(\Delta_{\mu+\mu'}, L_{\lambda+\mu'}).
\end{multline}

If we choose $\mu'\in(\Lambda^\sharp)^+$ very dominant (such that $\mu+\mu'$ is strictly dominant), then the standard object $\Delta_{\mu+\mu'}$ is the $!$-extension of the rank-one equivariant perverse local system on the single orbit in $\Fl_G^\omega$ corresponding to $\mu+\mu'$, and \eqref{8.2} is the normalized $!$-stalk of the object $L_{\lambda+\mu'}$ at $t^{\mu+\mu'}$.

Thus, to understand the class of \eqref{8.2}, it is enough to determine the action of the affine Yokonuma algebra on the Grothendieck group of ${\SD_{\kappa}^{H^\omega(\bF),\psi}(\Fl^\omega_G)}$ and express the element $[L_{\lambda+\mu'}]$ in terms of the standard basis. Then, the Verdier duality matrix gives its expansion in the costandard basis, and the coefficient at $t^{\mu+\mu'}$ is the class of the normalized $!$-stalk.

In this way, we compute the class
\[
[\RHom(\Delta_\mu,L_\lambda)] =\sum_j(-1)^j[H^j\RHom(\Delta_\mu,L_\lambda)]
\]
in the Grothendieck group, and its weight specialization. For IC sheaves, the pointwise purity in Proposition~\ref{prop:pointwise-purity} determines the cohomological degrees, while for $L_\lambda$, which need not be pointwise pure, the calculation determines the corresponding elements in the Grothendieck group.

\subsection{Finite-field and twisting conventions} From now on, we assume that the supergroup is $\GL(M|N)$, and we choose the sheaf theory to be the Weil $\bar{\BQ}_\ell$-constructible sheaves on the spaces (such as $\Gr_G, \Fl_G, \Conf$, $G$, $G/B$ etc.) defined over the finite field $k=\BF_Q$. We denote by $\cF\langle n\rangle:=\cF[n](n/2)$ the normalized shift, where $(1/2)$ is the half Tate twist with Frobenius eigenvalue $Q^{-1/2}$. We take tensor products of mixed sheaf categories over $D_\mxd(\pt)$ (mixed sheaves on $\Spec\BF_Q$), and tensor products of their Grothendieck groups over $K_0^\mxd(\pt)$.

In order to agree with the conventions in \cite{[T1],[FGT]}, we first do some simplifications and set conventions.

\subsubsection{Matrix conventions} We choose the local coordinate $t$, and then the renormalization by $\omega$ does not change the category itself. Thus, for simplicity, we work with the non-renormalized category $D_\kappa^{H(\bF),\psi}(\Fl_G)$. From this section onward, we use \textbf{upper} triangular Borel subgroups. For $\GL_r$, we denote the upper triangular Borel subgroup by $B_r$, its unipotent radical by $U_r$, the diagonal torus by $T_r$, and the corresponding Iwahori subgroup by $I_r$, and we denote $\mathfrak i_r=\Lie(I_r)$. For $M<N$, we denote by $U_{M,N}$ the unipotent radical of the upper triangular parabolic subgroup of block sizes $(M+1,1,\ldots,1)$, which is the positive version of the unipotent subgroup used for the polarization in Section~\ref{pol des}.

Also, since we want to compare the Iwahori Gaiotto category with the finite type category, it is more convenient to work with the additive character which is given by the \textbf{constant} term, instead of taking the coefficient of the $t^{-1}$-term. It is conjugate to (the positive analogue of) the character in \eqref{chi} by $t^\zeta$, where $\zeta=(\underbrace{0,\ldots,0}_{M+1},1,\ldots,N-M-1)$. Namely,
\begin{equation}\label{eq:typeA-constant-character}
\psi(u)=\sum_{j=M+1}^{N-1}[t^0]u_{j,j+1}, \qquad u\in U_{M,N}(\bF).
\end{equation}
Here, $[t^0]$ means the constant coefficient, and an empty sum is zero. Note that the restriction of this character to the finite-dimensional group $U_{M,N}$ coincides with the additive character $\psi$ used for the finite type categories below. Since we will only work with this character later, it does no harm that we also denote the additive character for $U_{M,N}(\bF)$ by $\psi$.

\subsubsection{The determinant extension} In the affine setting, let $L_0=\bO^N$. We regard the punctured determinant line bundle on $\GL_N(\bF)$ as the Kac--Moody extension
\[
\widetilde \GL_N(\bF):= \left\{(g,\tau)\ \middle|\
\begin{gathered}
g\in \GL_N(\bF),\\
\tau\in \det(gL_0:L_0)^\times
\end{gathered}
\right\}
\]
of $\GL_N(\bF)$, \[1\rightarrow \BG_m\rightarrow \widetilde \GL_N(\bF)\rightarrow \GL_N(\bF)\rightarrow 1.\]

We define its doubled extension $\hat \GL_N(\bF)$, i.e., the Baer sum of $\widetilde \GL_N(\bF)$ with itself. Equivalently,
\[
\hat \GL_N(\bF) \simeq \widetilde \GL_N(\bF)\times^{\mathbb G_m}\mathbb G_m,
\]
where the map $\BG_m\rightarrow \BG_m$ is the square map, so that
\[
1\to \mathbb G_m \to \hat \GL_N(\bF) \to \GL_N(\bF)\to 1
\]
is the central extension underlying the metaplectic twisting used below. We denote by $\hat T_N(\bF)$ and $\hat I_N$ the preimages of $T_N(\bF)$ and $I_N$ in this extension.

For the diagonal torus $T_N\subset \GL_N$, with the cocharacter lattice identified naturally with $\BZ^N$, we use the notation of Section~\ref{deform algebra}. Then, the doubled extension above is encoded by the integral forms on the factor $\GL_N$,
\[
q_{\BZ}(\lambda)=\sum_{i=1}^N \lambda_i^2, \qquad b_{\BZ}(\lambda,\mu)=2\sum_{i=1}^N \lambda_i\mu_i,
\]
\[
b'_{\BZ}(\lambda,\mu)=\sum_{i=1}^N \lambda_i\mu_i, \qquad b_{\BZ}=b'_{\BZ}+(b'_{\BZ})^{\mathrm{op}}.
\]

Thus, the twisting $\kappa$ is identified with a multiplicative local system (i.e., a Kummer local system) $\cL_\kappa$ on $\BG_m$. On this factor, $q=\kappa^{q_{\BZ}}$, $b=\kappa^{b_{\BZ}}$ and $b'=\kappa^{b'_{\BZ}}$, in multiplicative notation. Using the canonical splitting over $I_N$, we denote by $D_\kappa(\Fl)$ the category of right $I_N$-equivariant sheaves on $\hat\GL_N(\bF)$ which are $(\BG_m,\cL_\kappa)$-equivariant for the central action. Under la correspondance faisceaux-fonctions, the Kummer local system corresponds to a multiplicative character \[\kappa\colon \BF_Q^\times\rightarrow \bar{\BQ}_\ell^\times\simeq \BC^\times,\] and by taking Frobenius traces, the objects of $D_\kappa(\Fl)$ with Weil structure correspond to $\bar{\BQ}_\ell$-valued functions on $(\hat\GL_N(\bF)/I_N)(\BF_Q)$ which are equivariant with respect to $\kappa$.

\subsubsection{The finite and affine categories} For $M<N$, we regard $\GL_M$ as a subgroup of $\GL_N$ by the naive embedding into the upper-left block, and denote by $B_M$ and $B_N$ the upper triangular Borel subgroups of $\GL_M$ and $\GL_N$, respectively. Given Kummer local systems $\cL_\chi$ on $T_M$ and $\cL_{\chi'}$ on $T_N$, we denote their coordinate characters by $\chi=(\chi_1,\ldots,\chi_M)$ and $\chi'=(\chi'_1,\ldots,\chi'_N)$. Here, $\cL_{\chi_i}$ is the pullback of $\cL_\chi$ to the $i$-th diagonal copy of $\BG_m$, and similarly for $\chi'_j$. We denote by
\[\cC^{\fin,\mxd}_{\chi,\chi'}(M|N):=D_{\mxd}^{(B_M,\chi)\ltimes (U_{M,N},\psi)}(\GL_N/(B_N,\chi')),\] \[\cC^{\aff,\mxd}_{\kappa,\chi,\chi'}(M|N):=D_{\mxd,\kappa}^{(I_M,\chi)\ltimes (U_{M,N}(\bF),\psi)}(\GL_N(\bF)/(I_N,\chi')).\] Similarly, in the case $M=N$, \[\cC^{\fin,\mxd}_{\chi,\chi'}(N|N):=D_{\mxd}^{(B_N,\chi)}(\GL_N/(B_N,\chi')\times \BA^N),\] \[\cC^{\aff,\mxd}_{\kappa,\chi,\chi'}(N|N):=D_{\mxd,\kappa}^{(I_N,\chi)}(\GL_N(\bF)/(I_N,\chi')\times\bF^N).\] The affine categories consist of locally compact objects as in Section~\ref{Gaiotto category}. For $M=N$, we use the mixed $D^{!*}$-version of \cite[\S 3.2]{[BFGT]} in the vector direction. For $M<N$, we sometimes include the additive character in the subscripts. Note that the Verdier duality inverts the twisting and the multiplicative characters, and replaces $\psi$ by $-\psi$.

We will use the same letter $\chi$ to denote the character and the corresponding character local system $\cL_\chi$.

\subsubsection{Specialization of Grothendieck groups} For any given $M,N$, we denote
\[
\bK_{\chi,\chi'}^{\fin,\mxd}(M|N) :=K_0\bigl(\mathcal C_{\chi,\chi'}^{\fin,\mxd}(M|N)\bigr),
\]
and
\[
\bK_{\kappa,\chi,\chi'}^{\aff,\mxd}(M|N) :=K_0\bigl(\mathcal C_{\kappa,\chi,\chi'}^{\aff,\mxd}(M|N)\bigr).
\]
We denote the graded Grothendieck group with a similar notation, \[\bK_{\chi,\chi'}^{\fin,\gr}(M|N):= \bK_{\chi,\chi'}^{\fin,\mxd}(M|N)\otimes_{K_0^\mxd(\pt)}\BZ[\bv,\bv^{-1}],\] \[\bK_{\kappa,\chi,\chi'}^{\aff,\gr}(M|N):= \bK_{\kappa,\chi,\chi'}^{\aff,\mxd}(M|N)\otimes_{K_0^\mxd(\pt)}\BZ[\bv,\bv^{-1}].\] Here, $K_0^\mxd(\pt)$ denotes the Grothendieck group of the mixed Weil sheaves at one point, which is the free abelian group generated by the classes $[\alpha]$ of the one-dimensional Weil sheaves with Frobenius eigenvalue $\alpha$, where $\alpha\in \bar{\BQ}_\ell^\times$ and its absolute value belongs to $\{Q^{\frac{m}{2}}|m\in \BZ\}$ under the identification $\bar{\BQ}_\ell\simeq \BC$ fixed below. The map $ K_0^\mxd(\pt)\rightarrow \BZ[\bv,\bv^{-1}]$ is the specialization map which just remembers the grading by weights. Namely,
\begin{equation}\label{9.67 specialization}
\begin{split}
\on{wt}\colon  K_0^\mxd(\pt)\rightarrow \BZ[\bv,\bv^{-1}],\\
\sum n_{\alpha}[\alpha]\mapsto \sum n_\alpha\bv^{2\log_{Q}(|\alpha|)}.
\end{split}
\end{equation}

\subsubsection{} Let $v=Q^{1/2}$ be the positive real square root. We fix an isomorphism $\bar{\BQ}_\ell\simeq\BC$, and we work with the complexified Grothendieck group, where we will put a subscript $\BC$, e.g.,
\[
\bK_{\chi,\chi'}^{\fin,\gr}(M|N)_\BC :=\bK_{\chi,\chi'}^{\fin,\gr}(M|N)\otimes_\BZ\BC,
\]
and similarly in the affine case. Note that we have another specialization map
\begin{equation}\label{9.67' specialization}
\begin{split}
\on{wt}'\colon  K_0^\mxd(\pt)\rightarrow \BC[\bv,\bv^{-1}],\\
\sum n_{\alpha}[\alpha]\mapsto \sum n_\alpha\cdot (\alpha/|\alpha|)\cdot\bv^{2\log_{Q}(|\alpha|)}.
\end{split}
\end{equation}

We denote the corresponding Grothendieck group by
\[
\bK_{\chi,\chi'}^{\fin,\gr,\mathrm{tr}}(M|N) :=\bK_{\chi,\chi'}^{\fin,\mxd}(M|N) \otimes_{K_0^\mxd(\pt),\,\on{wt}'}\BC[\bv,\bv^{-1}],
\]
and similarly in the affine case. Then, the evaluation at $\bv=v$ gives the Frobenius trace.

\subsubsection{} We denote by $\bK_{\chi,\chi'}^{\fin}(M|N)_\BC$ the space of $\BC$-valued left $((B_M,\chi)\ltimes (U_{M,N},\psi))$-equivariant and right $(B_N,\chi')$-equivariant functions on $\GL_N(\BF_Q)$ in the case $M<N$. Similarly for the mirabolic case and the affine case: $\bK_{\chi,\chi'}^{\fin}(N|N)_\BC$, $\bK_{\kappa,\chi,\chi'}^{\aff}(M|N)_\BC$. For $M<N$, the affine functions are supported on finitely many relevant double cosets. We have the specialization map
\begin{equation}\label{9.67'' specialization}
\begin{split}
\on{sp}\colon  K_0^\mxd(\pt)\rightarrow \BC,\\
\sum n_{\alpha}[\alpha]\mapsto \sum n_\alpha \alpha.
\end{split}
\end{equation}

Similarly, replacing mixed sheaves by Weil sheaves, we define the equivariant categories $\cC^{\fin,\Weil}$ and $\cC^{\aff,\Weil}$, with Grothendieck groups $\bK^{\fin,\Weil}$ and $\bK^{\aff,\Weil}$.

In summary, we have the following commutative diagram, where $\operatorname{ch}$ sends an object to its Grothendieck class.
\[
\begin{tikzcd}[column sep=2.2em,row sep=3.2em]
\cC^{\fin,\Weil}_{\chi,\chi'}(M|N) \arrow[r,"\operatorname{ch}"] & \bK^{\fin,\Weil}_{\chi,\chi'}(M|N) \arrow[r,hook] & \bK^{\aff,\Weil}_{\kappa,\chi,\chi'}(M|N) & \cC^{\aff,\Weil}_{\kappa,\chi,\chi'}(M|N) \arrow[l,"\operatorname{ch}"'] \\
\cC^{\fin,\mxd}_{\chi,\chi'}(M|N) \arrow[u,hook] \arrow[r,"\operatorname{ch}"] & \bK^{\fin,\mxd}_{\chi,\chi'}(M|N) \arrow[u,hook] \arrow[d,"\operatorname{wt}"'] \arrow[dd,bend right=74,"\operatorname{wt}'"'] \arrow[ddd,bend right=84,"\operatorname{sp}"'] \arrow[r,hook] & \bK^{\aff,\mxd}_{\kappa,\chi,\chi'}(M|N) \arrow[u,hook] \arrow[d,"\operatorname{wt}"] \arrow[dd,bend left=74,"\operatorname{wt}'"] \arrow[ddd,bend left=84,"\operatorname{sp}"'] & \cC^{\aff,\mxd}_{\kappa,\chi,\chi'}(M|N) \arrow[l,"\operatorname{ch}"'] \arrow[u,hook] \\
{} & \bK^{\fin,\gr}_{\chi,\chi'}(M|N) \arrow[r,hook] & \bK^{\aff,\gr}_{\kappa,\chi,\chi'}(M|N) & {} \\
{} & \bK^{\fin,\gr,\mathrm{tr}}_{\chi,\chi'}(M|N) \arrow[d,"{\bv=v}"'] \arrow[r,hook] & \bK^{\aff,\gr,\mathrm{tr}}_{\kappa,\chi,\chi'}(M|N) \arrow[d,"{\bv=v}"] & {} \\
{} & \bK^{\fin}_{\chi,\chi'}(M|N)_{\BC} \arrow[r,hook] & \bK^{\aff}_{\kappa,\chi,\chi'}(M|N)_{\BC} & {}
\end{tikzcd}
\]

\subsection{Grothendieck groups and bases} In the formulas below, we use even coefficient lines. Note that, for the Wakimoto action \eqref{eq:normalized-Wakimoto-action}, we have
\[
[\cF\star_{\rm H}J_\eta] =(-1)^{\delta(\eta)} [\cF\star J_\eta].
\]
The same formula holds for Satake objects supported on the connected component of $\Gr_G$ containing $t^\eta$.

We denote by
\[
\bH_{\kappa,\chi,\chi'}^{\aff,\mxd}(M),\qquad \bH_{\kappa,\chi,\chi'}^{\aff,\gr}(M),\qquad \bH_{\kappa,\chi,\chi'}^{\aff}(M)
\]
the mixed, graded and usual affine Yokonuma--Hecke blocks, respectively, see Section~\ref{aff kl}. Their direct sum over the character indices is an algebra, where a product is zero unless the right monodromy character of the first factor equals the left monodromy character of the second factor. We denote by $H_{w,\chi,\chi'}^{\mxd}$ and $\underline H_{w,\chi,\chi'}^{\mxd}$ the standard and canonical basis elements, and use the superscript $\gr$ or no superscript for their specializations. Similarly, for the module $\bK_{\kappa,\chi,\chi'}^{\aff,\mxd}(M|N)$, we use $H_{\tilde w,\chi,\chi'}^{\mxd}$ and $\underline H_{\tilde w,\chi,\chi'}^{\mxd}$, where $\tilde w=(w,\beta)$, with the same notation for specializations.

\section{Monodromic finite type category}\label{sec:finite-type-monodromic} In this section, we study the finite type monodromic categories $\cC^{\fin}_{\chi,\chi'}(M|N)$, which will be used in Section~\ref{aff kl} to calculate the canonical basis in the affine case, see Proposition~\ref{prop 9.4.17}. First, we show that these categories decompose into tensor products of categories with trivial characters (Theorem~\ref{Thm 9.2.5}). Then, we give the action formulas of the Yokonuma algebras, and determine the bar involution and the canonical basis in the non-monodromic case (Corollary~\ref{cor 9.3.12}). Finally, we make the block decomposition explicit on the level of Grothendieck groups, which gives an algorithm to calculate the canonical basis in the monodromic case. As in Section~\ref{sec:type-a-conventions}, we use upper triangular Borel subgroups, and write $W_N^{\fin}=S_N$. We denote by $\leq_{B_N\cap w^{-1}B_Nw}$ the partial order attached to $B_N\cap w^{-1}B_Nw$, i.e., $i\leq_{B_N\cap w^{-1}B_Nw} j$ if and only if $i\leq j$ and $w(i)\leq w(j)$.
\subsection{Block decomposition of finite type monodromic category} For the monodromic finite type mirabolic Hecke category $D(B_N,\chi\backslash \GL_N\times \BA^N/B_N,\chi')$, the following lemma is from \cite{[T1]}.
\begin{lem}\label{lem 9.2.1}
The $B_N$-orbits in $\GL_N/B_N\times \BA^N$ are indexed by $(w,\beta)$, where $w\in W_N^\fin$, and $\beta$ is a  subset of $\{1,2,\cdots, N\}$. We require that if $i\leq_{B_N\cap w^{-1}B_Nw} j$, then, $j\in \beta$ implies that $i\in \beta$. The representative is given by $(wB_N,v_{w(\beta)})$, where $v_{w(\beta)}=\sum_{i\in w(\beta)}e_i$. The non-zero entries of the permutation matrix $w$ are in positions $(w(j),j)$. The index set is denoted by $RB^\fin(N|N)$.
\end{lem}
\begin{rem}
By choosing the set of maximal elements $\sigma$ with respect to the partial order $\leq_{B_N\cap w^{-1}B_Nw}$, for each fixed $w$, there is a bijection between the subsets $\beta$ satisfying the condition above and the subsets $\sigma\subset[1,N]$ such that the map $\sigma\rightarrow \{1,2,\cdots,N\}$ given by $i\mapsto w(i)$ is strictly decreasing. Thus, the orbits are also indexed by $(w,\sigma)$, and the representative is taken as $(w, v_{w(\sigma)})\in \GL_N/B_N\times \BA^N$. The inverse map (sending $\sigma$ to $\beta$) is given by saturation by the partial order $\leq_{B_N\cap w^{-1}B_Nw}$.
\end{rem}

An immediate consequence of Lemma \ref{lem 9.2.1} is
\begin{cor}
A $B_N$-orbit in $\GL_N\times \BA^N/B_N$ indexed by $(w,\beta)$ is $(B_N,\chi)$, $(B_N,\chi')$-relevant if and only if for any $i\in \{1,2,\cdots, N\}\setminus w(\sigma)$, we have $\chi_{i}=\chi'_{w^{-1}(i)}$.
\end{cor}

\subsubsection{The case \texorpdfstring{$M<N$}{M<N}} A similar result also works for $M<N$. Namely, we have
\begin{lem}\label{lem 9.2.1+}
The $B_M\ltimes U_{M,N}$-orbits in $\GL_N/B_N$ are indexed by $(w,\beta)$, where $w\in W_N^\fin$, and $\beta$ is a subset of $w^{-1}([1,M])$. We require that if $i<_{B_N\cap w^{-1}B_Nw}j $, and $j\in \beta\sqcup \{w^{-1}(M+1)\}$, and $i\in w^{-1}([1,M])$, then $i\in \beta$.

Let $\sigma$ be the set of maximal elements of $\beta\cup\{w^{-1}(M+1)\}$ that belong to $w^{-1}([1,M])$. The orbit is $(B_M,\chi)\ltimes (U_{M,N},\psi)$, $(B_N,\chi')$-relevant if and only if the following conditions hold.
\begin{enumerate}
\item the sequence $\{M+1, M+2,\cdots, N\}\rightarrow \{1,2,\cdots, N\}$ given by $i\mapsto w^{-1}(i)$ is strictly decreasing,\\
\item for any $i\in \{1,2,\cdots, M\}\setminus w(\sigma)$, we have $\chi_{i}=\chi'_{w^{-1}(i)}$.
\end{enumerate}
The representative of the corresponding relevant orbit in $\GL_N$ is given by $g_{w(\beta)}w$, where $g_{w(\beta)}=\id+\sum_{i\in w(\beta)}e_i\otimes e_{M+1}^*$.
\end{lem}

We give an example in Figure~\ref{fig:finite-orbit-classification}. Note that, among rows $1,\ldots,M$, every dot above and to the left of a blue dot or the red dot must be blue.

\begin{figure}[htbp]
\centering
\[
M=4,\quad N=6,\quad (w(1),\ldots,w(6))=(4,5,2,1,6,3).
\]
\[
\beta=\{1,3,4\},\qquad \sigma=\{3,4\}.
\]
\resizebox{\linewidth}{!}{\begin{tikzpicture}[x=1cm,y=-1cm,font=\small]
\definecolor{orbitblue}{RGB}{32,94,176}
\definecolor{orbitred}{RGB}{184,40,47}

\fill[orbitblue!7] (0,0)--(3.76,0)--(3.76,0.94)--(2.82,0.94)--(2.82,1.88)--(0.94,1.88)--(0.94,3.76)--(0,3.76)--cycle; \draw[orbitblue!60,line width=.6pt] (3.76,0)--(3.76,0.94)--(2.82,0.94)--(2.82,1.88)--(0.94,1.88)--(0.94,3.76)--(0,3.76);
\begin{scope}[xshift=0cm]
\node[font=\normalsize] at (2.82,-1.08) {The pair $(w,\beta)$}; \draw[black!13,line width=.25pt] (0,0)--(0,5.64); \draw[black!13,line width=.25pt] (0,0)--(5.64,0); \draw[black!13,line width=.25pt] (0.94,0)--(0.94,5.64); \draw[black!13,line width=.25pt] (0,0.94)--(5.64,0.94); \draw[black!13,line width=.25pt] (1.88,0)--(1.88,5.64); \draw[black!13,line width=.25pt] (0,1.88)--(5.64,1.88); \draw[black!13,line width=.25pt] (2.82,0)--(2.82,5.64); \draw[black!13,line width=.25pt] (0,2.82)--(5.64,2.82); \draw[black!13,line width=.25pt] (3.76,0)--(3.76,5.64); \draw[black!13,line width=.25pt] (0,3.76)--(5.64,3.76); \draw[black!13,line width=.25pt] (4.7,0)--(4.7,5.64); \draw[black!13,line width=.25pt] (0,4.7)--(5.64,4.7); \draw[black!13,line width=.25pt] (5.64,0)--(5.64,5.64); \draw[black!13,line width=.25pt] (0,5.64)--(5.64,5.64); \node[font=\small,anchor=south] at (0.47,-0.05) {$1$}; \node[font=\small,anchor=east] at (-0.05,0.47) {$1$}; \node[font=\small,anchor=south] at (1.41,-0.05) {$2$}; \node[font=\small,anchor=east] at (-0.05,1.41) {$2$}; \node[font=\small,anchor=south] at (2.35,-0.05) {$3$}; \node[font=\small,anchor=east] at (-0.05,2.35) {$3$}; \node[font=\small,anchor=south] at (3.29,-0.05) {$4$}; \node[font=\small,anchor=east] at (-0.05,3.29) {$4$}; \node[font=\small,anchor=south] at (4.23,-0.05) {$5$}; \node[font=\small,anchor=east] at (-0.05,4.23) {$5$}; \node[font=\small,anchor=south] at (5.17,-0.05) {$6$}; \node[font=\small,anchor=east] at (-0.05,5.17) {$6$}; \draw[->,black!55] (0,-0.49)--(5.64,-0.49) node[right] {$i$}; \draw[->,black!55] (-0.58,0)--(-0.58,5.64) node[below] {$j$}; \draw[orbitred!65,densely dashed] (1.41,0)--(1.41,5.64); \draw[orbitred!65,densely dashed] (0,4.23)--(1.41,4.23); \fill[orbitblue] (0.47,3.29) circle (3.4pt); \fill[orbitred] (1.41,4.23) circle (3.4pt); \fill[orbitblue] (2.35,1.41) circle (3.4pt); \fill[orbitblue] (3.29,0.47) circle (3.4pt); \fill[black!85] (4.23,5.17) circle (3.4pt); \fill[black!85] (5.17,2.35) circle (3.4pt); \draw[orbitblue,line width=.7pt] (2.35,1.41) circle (5.6pt); \draw[orbitblue,line width=.7pt] (3.29,0.47) circle (5.6pt); \node[orbitred,font=\scriptsize,anchor=west,fill=white,inner sep=1pt] at (1.59,4.23) {$M+1=5$}; \node[orbitred,font=\scriptsize] at (1.41,6.11) {$w^{-1}(M+1)=2$};
\end{scope}
\begin{scope}[xshift=8.15cm]
\node[font=\normalsize] at (2.82,-1.08) {The representative $g_{w(\beta)}w$}; \draw[black!13,line width=.25pt] (0,0)--(0,5.64); \draw[black!13,line width=.25pt] (0,0)--(5.64,0); \draw[black!13,line width=.25pt] (0.94,0)--(0.94,5.64); \draw[black!13,line width=.25pt] (0,0.94)--(5.64,0.94); \draw[black!13,line width=.25pt] (1.88,0)--(1.88,5.64); \draw[black!13,line width=.25pt] (0,1.88)--(5.64,1.88); \draw[black!13,line width=.25pt] (2.82,0)--(2.82,5.64); \draw[black!13,line width=.25pt] (0,2.82)--(5.64,2.82); \draw[black!13,line width=.25pt] (3.76,0)--(3.76,5.64); \draw[black!13,line width=.25pt] (0,3.76)--(5.64,3.76); \draw[black!13,line width=.25pt] (4.7,0)--(4.7,5.64); \draw[black!13,line width=.25pt] (0,4.7)--(5.64,4.7); \draw[black!13,line width=.25pt] (5.64,0)--(5.64,5.64); \draw[black!13,line width=.25pt] (0,5.64)--(5.64,5.64); \node[font=\small,anchor=south] at (0.47,-0.05) {$1$}; \node[font=\small,anchor=east] at (-0.05,0.47) {$1$}; \node[font=\small,anchor=south] at (1.41,-0.05) {$2$}; \node[font=\small,anchor=east] at (-0.05,1.41) {$2$}; \node[font=\small,anchor=south] at (2.35,-0.05) {$3$}; \node[font=\small,anchor=east] at (-0.05,2.35) {$3$}; \node[font=\small,anchor=south] at (3.29,-0.05) {$4$}; \node[font=\small,anchor=east] at (-0.05,3.29) {$4$}; \node[font=\small,anchor=south] at (4.23,-0.05) {$5$}; \node[font=\small,anchor=east] at (-0.05,4.23) {$5$}; \node[font=\small,anchor=south] at (5.17,-0.05) {$6$}; \node[font=\small,anchor=east] at (-0.05,5.17) {$6$}; \draw[->,black!55] (0,-0.49)--(5.64,-0.49) node[right] {$i$}; \draw[->,black!55] (-0.58,0)--(-0.58,5.64) node[below] {$j$}; \fill[black!85] (0.47,3.29) circle (3.4pt); \fill[orbitred] (1.41,4.23) circle (3.4pt); \fill[black!85] (2.35,1.41) circle (3.4pt); \fill[black!85] (3.29,0.47) circle (3.4pt); \fill[black!85] (4.23,5.17) circle (3.4pt); \fill[black!85] (5.17,2.35) circle (3.4pt); \fill[orbitblue] (1.41,0.47) circle (3.1pt); \fill[orbitblue] (1.41,1.41) circle (3.1pt); \fill[orbitblue] (1.41,3.29) circle (3.1pt);
\end{scope}
\node[anchor=north west,align=left,text width=6.7cm,font=\footnotesize] at (0,6.79) { \textcolor{orbitblue}{Blue}: $i\in\beta$. Circled: $i\in\sigma=\{3,4\}$.\\[3pt]
\textcolor{orbitred}{Red}: row $5$, column $2$. }; \node[anchor=north west,align=left,text width=6.5cm,font=\footnotesize] at (8.15,6.79){\textcolor{orbitblue}{Blue}: the additional entries at\\[2pt]row $w(i)$, column $2$,\\[2pt]$i\in\beta$.};
\end{tikzpicture}
} \caption{An example of the orbit classification in Lemma~\ref{lem 9.2.1+}. Rows increase downwards and columns increase to the right. In the left diagram, all dots are entries of $w$. The blue dots correspond to $\beta$, and the circled dots correspond to $\sigma$. In the right diagram, every marked entry is $1$ and every other entry is $0$. In this example, the relevance conditions in Lemma~\ref{lem 9.2.1+} say $w^{-1}(5)>w^{-1}(6)$, $\chi_3=\chi'_6$, and $\chi_4=\chi'_1$. The first condition fails since $w^{-1}(5)=2<5=w^{-1}(6)$.}
\label{fig:finite-orbit-classification}
\end{figure}

We denote by $RB^\fin(M|N)$ the set of $\psi$-relevant indices, and by $RB^\fin_{\chi,\chi'}(M|N)$ its $(\chi,\chi')$-relevant subset. For $(w,\beta)\in RB^\fin_{\chi,\chi'}(M|N)$, we associate it with the subset $\sigma\subset \beta \cup\{w^{-1}(M+1)\}$ of the maximal elements with respect to the order $<_{B_N\cap w^{-1}B_Nw}$ that belong to $w^{-1}([1,M])$, as in Lemma~\ref{lem 9.2.1+}. Also, we have the reduced representative of the corresponding orbit $g_{w(\sigma)}w$, where $g_{w(\sigma)}=\id+\sum_{i\in w(\sigma)}e_i\otimes e_{M+1}^*$.

The following lemma is also easy to see.
\begin{lem}
The orbit indexed by $(w,\beta)$ is isomorphic to $\BG_a^{l(w)-M+|\beta\setminus \sigma|}\times \BG_m^{|\sigma|}$.
\end{lem}

\subsubsection{The monodromic category} We want to study the monodromic finite type categories $\cC^{\fin}_{\chi,\chi'}(M|N):=D(B_M\ltimes U_{M,N},\chi,\psi\backslash \GL_N/B_N,\chi')$ for $M<N$, and $\cC^{\fin}_{\chi,\chi'}(N|N):=D(B_N,\chi\backslash \GL_N\times\BA^N/B_N,\chi')$. Note that the equivariance character for $(b_L,b_R)\colon g\mapsto b_Lgb_R^{-1}$ is $\chi\boxtimes(\chi')^{-1}$. For $M>N$, we define
\[
\cC^{\fin}_{\chi,\chi'}(M|N;\psi) :=\cC^{\fin}_{(\chi')^{-1},\chi^{-1}}(N|M;-\psi)
\]
by inversion. For $M=N$, the involution $(g,v)\mapsto(g^{-1},g^{-1}v)$ induces a $t$-exact, weight-preserving equivalence
\begin{equation}\label{eq:balanced-inversion}
\cC^{\fin}_{\chi,\chi'}(N|N) \simeq\cC^{\fin}_{(\chi')^{-1},\chi^{-1}}(N|N).
\end{equation}
We choose the auxiliary Kummer local systems and their weight-zero Weil structures over a common finite extension, and we denote by $Q$ the cardinality of this finite field throughout the finite type and affine type algorithms below.

\subsubsection{The characters} Let us first simplify the multiplicative character $\chi$ (and $\chi'$).

\begin{lem}\label{lem:adjacent-character-swap}
If $\chi_{i}\neq \chi_{i+1}$, we have a $t$-exact equivalence
\begin{equation}\label{iso 9.3}
\cC^{\fin}_{\chi,\chi'}(M|N)\simeq \cC^{\fin}_{\tilde{\chi},\chi'}(M|N),
\end{equation}
which preserves weights. Here, $\tilde\chi$ is obtained from $\chi$ by interchanging $\chi_i$ and $\chi_{i+1}$. Similarly, if $\chi'_j\neq\chi'_{j+1}$, the right convolution induces
\begin{equation}\label{iso 9.3 right}
\cC^{\fin}_{\chi,\chi'}(M|N)\simeq\cC^{\fin}_{\chi,s_j\chi'}(M|N).
\end{equation}
Both equivalences are compatible with Verdier duality after inverting the multiplicative characters (and replacing $\psi$ by $-\psi$ when $M<N$).
\end{lem}
\begin{proof}
Let us only consider the case $M=N$; the general case is similar.

Consider the category $D(B_N,\tilde{\chi}\backslash \GL_N/B_N,\chi)$, and consider the $B_N$-orbit in $\GL_N/B_N$ corresponding to the simple reflection $s_i=(i,i+1)\in W_N^\fin$. According to the equivariance condition, it is actually a minimal relevant orbit. In particular, the $!$-extension of the rank-one perverse local system, pure of weight $0$, from this orbit is clean. Now, consider the convolution with this clean extension object; it is a $t$-exact equivalence. Indeed, taking convolution with this clean extension object is invertible with the inverse functor given by convolving with the $*$-extension (also clean) of the rank-one local system on the same orbit in $D(B_N,\chi\backslash \GL_N/B_N,\tilde{\chi})$, thus giving an equivalence. Furthermore, since the convolution diagram \[B_N s_i B_N\overset{B_N}{\times} \GL_N\times \BA^N/B_N\rightarrow \GL_N\times \BA^N/B_N\] is affine, and the convolution is given by $!$-pushforward, it is left $t$-exact. Let $P_i=B_N\cup B_Ns_iB_N$ be the corresponding minimal parabolic subgroup. We extend the convolution diagram by replacing $B_Ns_iB_N$ with $P_i$. The map
\[
P_i\overset{B_N}{\times}(\GL_N\times\BA^N)/B_N \longrightarrow(\GL_N\times\BA^N)/B_N
\]
is proper with fiber $P_i/B_N\simeq\BP^1$. Since the extension is clean, the convolution is also given by $*$-pushforward, and thus, it is also right $t$-exact. Similarly for the weight structure, Verdier duality and right convolution.
\end{proof}

As a corollary, we can arrange $\chi, \chi'$ such that
\begin{equation}\label{chi 1}
\begin{split}
\chi_{1}=\chi_{2}=\cdots=\chi_{k_1}=b_1,\\
\chi_{k_1+1}=\chi_{k_1+2}=\cdots=\chi_{k_1+k_2}=b_2,\\
\cdots\\
\chi_{(\sum_1^{l-1} k_i)+1}=\chi_{(\sum_1^{l-1} k_i)+2}=\cdots=\chi_{\sum_1^{l} k_i}=b_l,
\end{split}
\end{equation}
and
\begin{equation}\label{chi 2}
\begin{split}
\chi'_{1}=\chi'_{2}=\cdots=\chi'_{k'_1}=b_1,\\
\chi'_{k'_1+1}=\chi'_{k'_1+2}=\cdots=\chi'_{k'_1+k'_2}=b_2,\\
\cdots\\
\chi'_{(\sum_1^{l-1} k'_i)+1}=\chi'_{(\sum_1^{l-1} k'_i)+2}=\cdots=\chi'_{\sum_1^{l} k'_i}=b_l,
\end{split}
\end{equation}
where the $b_i$ are the pairwise different Kummer local systems occurring in $\chi$ and $\chi'$, the multiplicities $k_i,k_i'\geq0$ satisfy $\sum_{i=1}^l k_i=M$ and $\sum_{i=1}^l k_i'=N$, and a block with multiplicity zero is omitted.

\subsubsection{The block decomposition} We claim
\begin{thm}\label{Thm 9.2.5}
Let $\chi,\chi'$ be as in \eqref{chi 1} and \eqref{chi 2}. There is a $t$-exact equivalence
\begin{equation}\label{iso 9.6}
\cC^{\fin}_{\chi,\chi'}(M|N)\simeq \bigotimes_{i=1}^l \cC^{\fin}_{1,1}(k_i|k'_i).
\end{equation}
The equivalence is compatible with the Frobenius weight structure and with Verdier duality after inverting the multiplicative characters and replacing each additive character $\psi$ by $-\psi$.
\end{thm}
First, we note that

\begin{lem}\label{lem 9.2.6}
For any Kummer local system $\cL_b$, there is a canonical $t$-exact equivalence
\begin{equation}\label{eq 9.7}
\cC^{\fin}_{\chi,\chi'}(M|N)\simeq \cC^{\fin}_{b\chi,b\chi'}(M|N).
\end{equation}
\end{lem}
\begin{proof}
The equivalence is given by tensoring with the local system given by pulling back the Kummer local system $\cL_b$ along $\det\colon  \GL_N\rightarrow \BG_m$.
\end{proof}

\subsubsection{The Heisenberg group} Let us consider $B_M$, the upper triangular subgroup of $\GL_M$, and regard $B'_{M-1}$ as a subgroup of $B_M$ by embedding into the lower-right part. Note that we denote by $B_{M-1}$ the subgroup embedded into the upper-left part.

The key observation is the following.
\begin{lem}\label{9.2.7}
Assume $M\geq1$. Given a category $\cC$ with an action of $B_M\ltimes\BG_a^M$, we have a fully faithful embedding
\begin{equation}\label{9.8}
\cC^{B'_{M-1}\ltimes \BG_a^M,\chi,\psi'}\hookrightarrow \cC^{B_M,\chi},
\end{equation}
where the character $\psi'\colon \BG_a^M\rightarrow \BG_a$ is given by the projection to the first coordinate.

Furthermore, if the action of $B_M\ltimes \BG_a^M$ on the category comes from an action of $\GL_M\ltimes \BG_a^M$, we have
\begin{equation}\label{9.9}
\cC^{B_{M-1}\ltimes \BG_a^M,\chi',\psi}\hookrightarrow \cC^{B_M,\chi},
\end{equation}
where the character $\psi\colon \BG_a^M\rightarrow \BG_a$ is given by the projection to the $M$-th coordinate. Here, $\chi$ on $B'_{M-1}$ denotes the restriction to coordinates $2,\ldots,M$, and $\chi'$ on $B_{M-1}$ (not to be confused with the right character $\chi'$ elsewhere) has coordinates $\chi'_{j}=\chi_{j+1}$.

\end{lem}
\begin{proof}
Recall the proof of \cite[Lemma 5.3.1]{[TY1]}: for any category $\cD$ which admits an action of the Heisenberg group $\Hei=V\oplus V^*\oplus \BG_a$, there is an equivalence given by relative averaging
\begin{equation}
\cD^{V\oplus \BG_a,\psi'}\rightarrow\cD^{\BG_a,\psi'}\rightarrow \cD^{V^*\oplus \BG_a,\psi'}.
\end{equation}

Let $V$ be $(M-1)$-dimensional. We identify $V\oplus \BG_a$ as a subgroup $\begin{pmatrix}1&0&\BG_a\\
&\Id_{M-1}&V\\
&&1\end{pmatrix}$ of $B_M\ltimes \BG_a^M$. Similarly, we identify $V^*\oplus \BG_a$ as a subgroup of $B_M\ltimes \BG_a^M$ by $\begin{pmatrix}1&V^*&\BG_a\\
&\Id_{M-1}&0\\
&&1\end{pmatrix}$.

We have the following equivalence given by relative averaging
\begin{equation}
\cC^{(B'_{M-1},\chi)\ltimes (V\oplus \BG_a,\psi')}\simeq (\cC^{V\oplus \BG_a,\psi'})^{(B'_{M-1},\chi)}\simeq (\cC^{V^*\oplus \BG_a,\psi'})^{(B'_{M-1},\chi)}\simeq \cC^{(B'_{M-1},\chi)\ltimes (V^*\oplus \BG_a,\psi')}.
\end{equation}

The second step of the proof is the equivalence between the Kirillov model and the Whittaker model, see \cite[Appendix A]{[GL]}. Note that the category $\cC^{(B'_{M-1},\chi)\ltimes V^*}$ admits an action of $\BG_m\ltimes \BG_a= \begin{pmatrix} \BG_m&0&\BG_a\\
&1&0\\ &&1
\end{pmatrix}$, and we regard $\cC^{(B'_{M-1},\chi)\ltimes (V^*\oplus \BG_a,\psi')}$ as $(\cC^{(B'_{M-1},\chi)\ltimes V^*})^{\BG_a,\psi'}$.

Consider the full subcategory of $(\cC^{(B'_{M-1},\chi)\ltimes V^*})^{\BG_m,\chi_{1}}=\cC^{B_M,\chi}$ spanned by the kernel of the $*$-averaging functor
\begin{equation}
\Av_*^{\BG_a}\langle 1\rangle\colon (\cC^{(B'_{M-1},\chi)\ltimes V^*})^{\BG_m,\chi_{1}}\rightarrow (\cC^{(B'_{M-1},\chi)\ltimes V^*})^{(\BG_m,\chi_{1})\ltimes \BG_a}.
\end{equation}

By \cite[Proposition A.3.2]{[GL]}, with the normalization of Lemma~\ref{Whit and Kir}, this kernel is equivalent to $(\cC^{(B'_{M-1},\chi)\ltimes V^*})^{\BG_a,\psi'}$.

The second claim follows from the fact that left multiplying by the permutation element corresponding to $(1,2,3,\cdots,M)\in S_M$ induces an equivalence between the invariants with respect to $B_{M-1}$ and $B'_{M-1}$, i.e.,
\begin{equation}
\cC^{B_{M-1}\ltimes \BG_a^M,\chi',\psi}\simeq \cC^{B'_{M-1}\ltimes \BG_a^M,\chi,\psi'}.
\end{equation}

\end{proof}
\begin{rem}\label{rem 9.2.8}
According to the proof, the fully faithful functor \eqref{9.8} is an equivalence if and only if for any object $\cF\in \cC^{B_{M},\chi}$, we have $\Av_*^{\BG_a}(\cF)=0$. In particular, it is satisfied if $\cC^{(B_{M},\chi)\ltimes \BG_a}=0$.
\end{rem}

\begin{prop}\label{prop 9.2.9}
Assume $0<M\leq N$ and $\chi_1$ is different from $\chi'_j$ for any $1\leq j\leq N$. Let $\tau=(1,2,\ldots,M)$, and regard $\tau^{-1}\chi$ as the character $(\chi_2,\ldots,\chi_M)$ of $B_{M-1}$. Then, we have $\Av_*^{\BG_a}(\cF)=0$ for any object $\cF\in\cC^{\fin}_{\chi,\chi'}(M|N)$, where $\Av_*^{\BG_a}$ is the averaging functor of Remark~\ref{rem 9.2.8}. In particular, we have a $t$-exact equivalence
\begin{equation}\label{iso 9.14}
\cC^{\fin}_{\chi,\chi'}(M|N)\simeq \cC^{\fin}_{\tau^{-1}\chi,\chi'}(M-1|N),
\end{equation}
which is compatible with the Frobenius weight structure and with Verdier duality after inverting the multiplicative characters and replacing $\psi$ by $-\psi$.
\end{prop}
\begin{proof}
Assume that $M<N$. Then, the category $D(U_{M,N},\psi\backslash \GL_N/B_N,\chi')$ admits an action of $\GL_M\ltimes \BG_a^{M}$. Lemma~\ref{9.2.7} says that we have a fully faithful embedding
\begin{equation}
D(B_{M-1} \ltimes U_{M-1,N},\tau^{-1}\chi,\psi\backslash \GL_N/B_N,\chi') \hookrightarrow \cC^{\fin}_{\chi,\chi'}(M|N).
\end{equation}

We claim that it is an equivalence, and according to Remark \ref{rem 9.2.8}, we only need to show $D(((B_M,\chi)\ltimes\BG_a)\ltimes(U_{M,N},\psi)\backslash\GL_N/B_N,\chi')=0$. Here, $\BG_a$ acts by the matrices $\id+a e_1\otimes e_{M+1}^*$, which normalize $U_{M,N}$ and preserve $\psi$.

According to Lemmas~\ref{lem 9.2.1} and \ref{lem 9.2.1+} and the assumption, the index $(w,\beta)$ of any relevant orbit in $D(B_{M} \ltimes U_{M,N},\chi,\psi\backslash \GL_N/B_N,\chi')$ must satisfy the property that $1\in w(\sigma)$. On the other hand, if the orbit corresponding to $(w,\beta)$ belongs to a relevant orbit in $D(((B_M,\chi)\ltimes\BG_a)\ltimes(U_{M,N},\psi)\backslash\GL_N/B_N,\chi')$, then $(w,\beta)$ is relevant for $D(B_{M} \ltimes U_{M,N},\chi,\psi\backslash \GL_N/B_N,\chi')$, thus, $g_{w(\sigma)}=\id+e_{1}\otimes e_{M+1}^*+\cdots$. Furthermore, the left $B_{M} \ltimes U_{M,N}$, right $B_N$ orbit corresponding to $a\cdot (g_{w(\sigma)}w)$ is also relevant, where $a$ is an arbitrary element in $\BG_a$. However, it is impossible: taking $a=-1$, we remove the term $e_{1}\otimes e_{M+1}^*$ from $g_{w(\sigma)}$ and the corresponding orbit is not relevant, since $1$ no longer belongs to the corresponding new $w(\sigma)$.

For $M=N$, we apply Lemma~\ref{9.2.7} to $\cC=D((\GL_N\times\BA^N)/(B_N,\chi'))$, and the descent along the translation torsor identifies the source of \eqref{9.9} with $D(B_{N-1},\tau^{-1}\chi\backslash\GL_N/B_N,\chi')$. Then, the same argument applies. Namely, every relevant orbit has $1\in w(\sigma)$, while the translation by $-e_1$ removes this coordinate and gives an irrelevant orbit. Thus, the averaging along $\BG_a e_1$ vanishes.

Note that, for the descent along the translation torsor in the preceding paragraph, we use the normalized smooth pullback, and this descent and the equivalence in the proof of Lemma~\ref{9.2.7} preserve the perverse $t$-structure, weights and Verdier duality. For the Whittaker--Kirillov equivalence, we use Lemma~\ref{Whit and Kir}.

\end{proof}
\subsubsection{The Whittaker--Kirillov equivalence} Recall that the Whittaker--Kirillov equivalence of \cite[Proposition A.3.2]{[GL]} is given by $\Av_*^{\BG_m}\langle1\rangle$, with the inverse $\Av_!^{\BG_a,\psi}\langle-1\rangle$, and the same holds in the $\chi$-twisted case.
\begin{lem}\label{Whit and Kir}
\begin{enumerate}
\item Given a space $\cY$ over $\BF_Q$ which admits an action of $\BG_m\ltimes\BG_a$, the functor
\begin{equation}
D^{\BG_a,\psi}(\cY)\rightarrow D(\cY)\overset{\Av_*^{\BG_m,\chi}\langle 1\rangle}{\rightarrow} D^{\BG_m,\chi}(\cY)
\end{equation}
is $t$-exact.
\item The functor
\begin{equation}
D^{\BG_m,\chi}(\cY)\rightarrow D(\cY)\overset{\Av_!^{\BG_a,\psi}\langle -1\rangle}{\rightarrow}  D^{\BG_a,\psi}(\cY)
\end{equation}
is $t$-exact, preserves the weight structure, and is compatible with the duality.
\end{enumerate}
\end{lem}
\begin{proof}
For the first claim, we rewrite the functor above as the following composition $D^{\BG_a,\psi}(\cY)\simeq (D^{ \BG_a,\psi}(\cY\times \BG_m))^{\BG_m,\chi}\rightarrow (D^{\BG_a,\psi}(\cY\times \BG_a))^{\BG_m,\chi}\overset{FT}{\rightarrow} D^{\BG_m,\chi}(\cY\times \{\psi\})\simeq D^{\BG_m,\chi}(\cY)$.

Here, the action of $\BG_a$ on $\cY\times\BG_a$ is given by $(a,(y,a'))\mapsto(y+aa',a')$, and the second functor is the $*$-pushforward along the open embedding. The extension along this open embedding is clean, since $\BG_a$ acts trivially on $\cY\times\{0\}$ and $\psi$ is non-trivial. Hence, it is $t$-exact.

The second claim follows similarly; the only difference is that instead of using the pushforward along the open embedding, we use the fact that the pullback to an open subspace is $t$-exact and compatible with the duality.
\end{proof}

\begin{rem}
There is a more direct way to show Proposition~\ref{prop 9.2.9}. Namely, a fully faithful $t$-exact functor between bounded derived categories with finite-length hearts is an equivalence if it induces a bijection on simple objects.
\end{rem}
\subsubsection{Proof of the block decomposition}
\begin{proof}[Proof of Theorem \ref{Thm 9.2.5}]
By inversion, we can assume that $M\leq N$. Let $r_i=\max(k_i,k'_i)$ and $L=\sum_i r_i$. We choose additional characters, which are pairwise distinct and distinct from all $b_i$, and call them generic. First, according to Proposition~\ref{prop 9.2.9}, we extend the left character to rank $N$ by generic characters. Then, by Lemma~\ref{prop 9.2.11}, we add generic characters in pairs, since $\cC^{\fin,\mxd}_{\xi,\xi'}(1|1)\simeq D_\mxd(\pt)$ for $\xi\ne\xi'$. Thus, we obtain the case of rank $L|L$, where the extended characters are again denoted by $\chi,\chi'$, and the case $L=0$ is immediate.

Then, according to Lemma~\ref{lem:adjacent-character-swap}, we can switch adjacent rows and columns with generic values of characters, such that the characters $\chi$ and $\chi'$ satisfy
\begin{equation}
\begin{split}
\chi_{1}=\chi_{2}=\cdots=\chi_{k_1}=b_1,\\
\chi_{r_1+1}=\chi_{r_1+2}=\cdots=\chi_{r_1+k_2}=b_2,\\
\cdots\\
\chi_{(\sum_1^{l-1} r_i)+1}=\chi_{(\sum_1^{l-1} r_i)+2}=\cdots=\chi_{(\sum_1^{l-1} r_i)+k_l}=b_l,
\end{split}
\end{equation}
and
\begin{equation}
\begin{split}
\chi'_{1}=\chi'_{2}=\cdots=\chi'_{k'_1}=b_1,\\
\chi'_{r_1+1}=\chi'_{r_1+2}=\cdots=\chi'_{r_1+k'_2}=b_2,\\
\cdots\\
\chi'_{(\sum_1^{l-1} r_i)+1}=\chi'_{(\sum_1^{l-1} r_i)+2}=\cdots=\chi'_{(\sum_1^{l-1} r_i)+k'_l}=b_l,
\end{split}
\end{equation}
where the $b_i$ are pairwise different, $\sum_1^l k_i=M$, $\sum_1^l k'_i=N$, and all the other $\chi_{i}, \chi'_{j}$ are generic.

According to Lemma~\ref{prop 9.2.11}, this category is equivalent to the tensor product of the categories of ranks $r_i|r_i$. Since $r_i=\max(k_i,k'_i)$, the generic characters occur on at most one side of each factor. By Lemma~\ref{lem:adjacent-character-swap}, we move them to the first positions, and then remove them by Proposition~\ref{prop 9.2.9}, using inversion for the right characters. Then, we obtain $\cC^{\fin}_{b_i^{k_i},b_i^{k'_i}}(k_i|k'_i)$ in the $i$-th factor, where $b_i^k$ denotes the constant tuple of length $k$. Thus, Lemma~\ref{lem 9.2.6} with $b=b_i^{-1}$ implies the desired statement. The compatibility with weights and duality follows from that of these equivalences.
\end{proof}

\subsubsection{}
\begin{lem}\label{prop 9.2.11}
Assume $\sum_{a=1}^l r_a=L$, and write $I_a=[\sum_{b<a}r_b+1,\sum_{b\leq a}r_b]$, $\chi^{[a]}=\chi|_{I_a}$, and $(\chi')^{[a]}=\chi'|_{I_a}$. Suppose that for any $1\leq i,j\leq L$, if $\chi_{i}=\chi'_{j}$, then $i,j\in I_a$ for some $a$. Then, there is a $t$-exact equivalence
\begin{equation}\label{eq 9.20}
\cC^{\fin}_{\chi,\chi'}(L|L)\simeq \bigotimes_{a=1}^l\cC^{\fin}_{\chi^{[a]},(\chi')^{[a]}}(r_a|r_a).
\end{equation}
The equivalence preserves weights and is compatible with Verdier duality after inverting the characters. Empty blocks are omitted.
\end{lem}
\begin{proof}
By induction, it is enough to consider the case of two blocks, of sizes $r$ and $s$. First, using Lemma~\ref{lem:adjacent-character-swap}, we reverse the two blocks of $\chi'$. Let $\omega$ be the permutation which sends the first $s$ positions to the last $s$ positions and the last $r$ positions to the first $r$ positions, preserving the order within each block. For a relevant $(w,\sigma)$, if $w\notin(S_r\times S_s)\omega$, then there exist $a\leq s<b$ with $w(a)\leq r<w(b)$, and the assumption on the characters forces $a,b\in\sigma$, which contradicts the fact that $w|_\sigma$ is strictly decreasing.

Thus, all the relevant orbits lie in the affine open subset $P_{r,s}\omega B_L\times\BA^L$, where $P_{r,s}$ is the upper parabolic subgroup of $\GL_L$ corresponding to the partition $(r,s)$, and the restriction to this open subset is a clean, $t$-exact equivalence. After the right multiplication by $\omega^{-1}$, the left and right unipotent actions eliminate the off-diagonal blocks, and the remaining data are a matrix in $\GL_r\times\GL_s$ and a vector in $\BA^r\oplus\BA^s$. Then, the normalized smooth descent gives the two factors on the right-hand side of \eqref{eq 9.20}. These equivalences preserve weights and commute with Verdier duality after inverting the characters.
\end{proof}

\subsection{Action of Yokonuma algebra}
\subsubsection{Twisted Yokonuma algebra} Given two characters $\chi,\chi'$ of $T_N$, consider the space of $\BC$-valued functions on $\GL_N(\BF_Q)$ which are left $(B_N,\chi)$-invariant and right $(B_N,\chi')$-invariant, and denote it by $\bH_{\chi,\chi'}^\fin(N)_\BC$. Obviously, there is a convolution, normalized by $|B_N(\BF_Q)|^{-1}$,
\begin{equation}
\bH_{\chi,\chi'}^\fin(N)_\BC\otimes \bH_{\chi',\chi''}^\fin(N)_\BC\rightarrow \bH_{\chi,\chi''}^\fin(N)_\BC.
\end{equation}

Given $w\in W^\fin_N$, there is a non-zero $(\chi,\chi')$-invariant function on $B_NwB_N$ if and only if $\chi=w(\chi')$, cf., \cite[Observation 3.1]{[And]}. When $\chi=w(\chi')$, we denote by $T_{w,\chi,\chi'}$ the characteristic $(\chi,\chi')$-invariant function on the $B_N\times B_N$-orbit in $\GL_N$ corresponding to $w$, and normalized by the requirement that its value at $w\in \GL_N$ is $1$, and set $T_{w,\chi,\chi'}:=0$ otherwise.

The following properties appear in \cite[Proposition 3.3]{[And]}.
\begin{lem}
\begin{enumerate}
\item if $w_1,w_2\in W^{\fin}_N$ and $l(w_1)+l(w_2)=l(w_1w_2)$, then
\[T_{w_1, w_1w_2(\chi),w_2(\chi)}T_{w_2,w_2(\chi), \chi}=T_{w_1w_2, w_1w_2(\chi),\chi}.\]
\item if $s=s_i=(i,i+1)$, then we have
\[T_{s, \chi, s\chi}T_{s, s\chi,\chi}=Q\cdot T_{e,\chi,\chi}+(Q-1)\cdot\chi_{i+1}(-1)\cdot T_{s, \chi,\chi}.\]
\end{enumerate}
\end{lem}
\subsubsection{Convolution formulas} The convolution gives rise to a bimodule structure on the direct sum over the characters. Namely,
\begin{equation}
\bH_{\chi'',\chi}^\fin(M)_\BC\otimes \bK_{\chi,\chi'}^\fin(M|N)_\BC\otimes \bH_{\chi',\chi'''}^\fin(N)_\BC \rightarrow \bK_{\chi'',\chi'''}^\fin(M|N)_\BC.
\end{equation}
We need to understand the above bimodule structure. It is enough to determine the action of the basis elements $T_{s,\chi,\chi'}$ associated with simple reflections $s$.

Let $H=B_M\ltimes U_{M,N}$ for $M<N$, and $H=B_N$, $\psi=0$ for $M=N$. Following the convention in \cite{[And]}, we set \[\hat T_{\tilde{w},\chi,\chi'}:= \sum_{a\in H, b\in B_N} (\chi,\psi)(a)\cdot\chi'(b)\cdot\mathbb{1}_{a\dot{\tilde w}b}.\] Here, $\dot{\tilde w}$ is the chosen orbit representative, the sums are over finite-field points, and $(\chi,\psi)(a)$ and $\chi'(b)$ are the values of the character functions, respectively. We denote by $T_{\tilde w,\chi,\chi'}$ the equivariant function supported on this orbit with value $1$ at $\dot{\tilde w}$, and we set it to be $0$ if the orbit is irrelevant. Thus, we have \[\hat T_{\tilde{w},\chi,\chi'}=|\Stab_H(\dot{\tilde w}B_N/B_N)|\cdot T_{\tilde{w},\chi,\chi'}.\]

\subsubsection{} Let $\tilde{w}=(w,\beta)\in RB^\fin(M|N)$, and let $s=(i,i+1)\in S_N$, $i \in \{1,\cdots, N-1\}$. We denote $\tilde{w}s=(ws, s(\beta))$, and use $\Delta$ for the symmetric difference. In the left diagram of Figure~\ref{fig:finite-orbit-classification}, the operation $(w,\beta)\mapsto(ws,s\beta)$ exchanges the columns $i,i+1$ together with their blue dots, and the left multiplication by $s_i\in S_M$ exchanges the rows $i,i+1$. Replacing $\beta$ by $\beta\Delta\{i+1\}$ changes whether the dot in column $i+1$ is blue, and then we take $\sigma$ for the new pair. The following proposition on the right module structure of $\bK^\fin(N|N)_\BC$ over the finite Hecke algebra $\bH^{\fin}(N)_\BC$ was given in \cite[Proposition 2]{[T1]}, and its monodromic version was given in \cite[Theorem 3.1]{[And]}.
\begin{prop}
Fix $\tilde w = (w,\beta) \in RB^\fin(N|N)$ and $s = (i,i+1)$. Then, we have
\begin{equation}
\hat T_{w,\beta} T_s =
\begin{cases}
Q \hat T_{ws,s\beta}, & ws > w,\;  \iota\nsubseteq \sigma'\\

(Q-1)\hat T_{ws,s\beta} + \hat T_{ws,\beta\setminus\{i+1\}}, & ws > w,\; \iota\subseteq \sigma', \\

\hat T_{ws,s\beta} + (Q-1)\hat T_{w,\beta}, & ws < w,\; i \notin \beta, \\

\hat T_{ws,s\beta} + (Q-2)\hat T_{w,\beta} + \hat T_{w,\beta\setminus\{i+1\}}, & ws < w,\; i \in \beta,\; i+1 \in \sigma, \\

\hat T_{ws,(\beta\cup\{i+1\})} + (Q-1)\hat T_{w,\beta\cup\{i+1\}}, & ws < w,\; i \in \beta,\; i+1 \notin \sigma .
\end{cases}
\end{equation}
Here, $\iota:=\{i,i+1\}$, $\sigma=\sigma(w,\beta)$, and $\sigma'=\sigma(ws,s\beta)$.
\end{prop}

\begin{prop}
\begin{enumerate}
\item If $ws > w$ then
\[
\hat T_{\tilde w,\chi,\chi'} T_{s,\chi',s\chi'} =
\begin{cases}
Q\, \hat T_{\tilde{w}s,\chi,s\chi'}, &  \iota\nsubseteq \sigma', \\[6pt]
(Q-1)\hat T_{\tilde{w}s,\chi,s\chi'} + \hat T_{ws,\beta\setminus\{i+1\},\chi,s\chi'}, &  \iota\subseteq \sigma'.
\end{cases}
\]
\item If $ws < w$ and $i \notin \beta$ then
\[
\hat T_{\tilde w,\chi,\chi'} T_{s,\chi',s\chi'} = \hat T_{\tilde{w}s,\chi,s\chi'} + (Q-1)\chi'_i(-1)\hat T_{\tilde w,\chi,s\chi'}.
\]
\item If $ws < w$, $i \in \beta$ and $i+1 \in \sigma$ then
\[
\hat T_{\tilde w,\chi,\chi'} T_{s,\chi',s\chi'} = \hat T_{ws,s\beta,\chi,s\chi'} + \chi'_i(-1)\hat T_{w,\beta\setminus\{i+1\},\chi,s\chi'} + \chi'_i(-1) J\!\left( \frac{\chi'_{i+1}}{\chi'_{i}}, \frac{\chi'_i}{\chi_{w(i+1)}} \right) \hat T_{\tilde w,\chi,s\chi'}.
\]
\item If $ws < w$, $i \in \beta$, $i+1 \notin \sigma$ then
\[
\hat T_{\tilde w,\chi,\chi'} T_{s,\chi',s\chi'} = \hat T_{ws,s\beta\cup\{i\},\chi,s\chi'} + (Q-1)\chi'_{i+1}(-1) \hat T_{w,\beta\cup\{i+1\},\chi,s\chi'}.
\]
\end{enumerate}
Here, $J(\chi,\chi')=\sum_{a\in \BF_Q-\{0,1\}}\chi(a)\chi'(1-a)$ is the Jacobi sum.
\end{prop}
\begin{rem}
Note that our convention here differs from \cite{[And]} by changing the character $\chi$ to its inverse.
\end{rem}

\subsubsection{The case \texorpdfstring{$M<N$}{M<N}} One has the following generalization to the case $M<N$, and the proof is by straightforward calculation. Let $G(\psi,\chi)=\sum_{a\in \BF_Q^\times}\chi(a)\psi(a^{-1})$ be the Gauss sum, and denote this sum by $G(1,\chi)$ when the additive character is trivial.

\begin{prop}\label{9.3.4}
For $M<N$, we fix $\tilde w = (w,\beta) \in RB^\fin_{\chi,\chi'}(M|N)$ and $s = (i,i+1)$, and denote $\sigma^*=\sigma(w,\beta)\cup\{w^{-1}(M+1)\}$ and $\sigma^{',*}= \sigma(ws,s\beta)\cup\{(ws)^{-1}(M+1)\}$.

\begin{enumerate}
\item If $ws > w$ then
\[
\hat T_{\tilde w,\chi,\chi'} T_{s,\chi',s\chi'} =
\begin{cases}
Q\, \hat T_{\tilde{w}s,\chi,s\chi'}, & \iota \nsubseteq \sigma^{',*}, \\[6pt]
(Q-1)\hat T_{\tilde{w}s,\chi,s\chi'} + \hat T_{ws,s\beta\setminus\{i+1\},\chi,s\chi'}, & \iota\subseteq \sigma^{',*}.
\end{cases}
\]
\item If $ws < w$ and $i\in w^{-1}([1,M])\setminus \beta$ then
\[
\hat T_{\tilde w,\chi,\chi'} T_{s,\chi',s\chi'} = \hat T_{\tilde{w}s,\chi,s\chi'} + (Q-1)\chi'_i(-1)\hat T_{\tilde w,\chi,s\chi'}.\]
\item If $ws < w$, $i \in \beta\sqcup w^{-1}([M+1,N])$ and $i+1 \in \sigma\sqcup w^{-1}([M+1,N])$,
\begin{enumerate}
\item  if $M+1\leq w(i+1)=w(i)-1\leq N-1$, then
\[
\hat T_{\tilde w,\chi,\chi'} T_{s,\chi',s\chi'} = \chi'_{i+1}(-1)G(\psi,\chi'_i/\chi'_{i+1})\hat T_{\tilde w,\chi,s\chi'};\]
\item if $w(i)=M+1$, then \[
\begin{aligned}
\hat T_{\tilde w,\chi,\chi'} T_{s,\chi',s\chi'} &= \hat T_{ws,s\beta,\chi,s\chi'}\\
&{}+ \chi'_{i+1}(-1)\hat T_{w,\beta\setminus\{i+1\},\chi,s\chi'}\\
&{}+ \frac{\chi_{w(i+1)}\chi_i'}{\chi'_{i+1}}(-1) J\!\left(\frac{\chi_{w(i+1)}}{\chi_{i+1}'}, \frac{\chi_i'}{\chi_{w(i+1)}}\right) \hat T_{\tilde w,\chi,s\chi'}.
\end{aligned}
\]
\item otherwise,
\[
\begin{aligned}
\hat T_{\tilde w,\chi,\chi'} T_{s,\chi',s\chi'} &= \hat T_{ws,s\beta,\chi,s\chi'}\\
&{}+ \chi_{w(i+1)}(-1)\hat T_{w,\beta\setminus\{i+1\},\chi,s\chi'}\\
&{}+ \chi_{w(i+1)}(-1) J\!\left( \frac{\chi_{w(i+1)}}{\chi_{i+1}'}, \frac{\chi_i'}{\chi_{w(i+1)}}\right) \hat T_{\tilde w,\chi,s\chi'}.
\end{aligned}
\]
\end{enumerate}

\item If $ws < w$, $i \in \beta\sqcup w^{-1}([M+1,N])$, $i+1 \notin \sigma\sqcup w^{-1}([M+1,N])$,
\begin{enumerate}
\item if $w(i)\neq M+1$, or $w(i)=M+1$ and $i+1\in \beta\sqcup w^{-1}([M+1,N])$,
\[
\hat T_{\tilde w,\chi,\chi'} T_{s,\chi',s\chi'} = \hat T_{ws,s\beta\cup\{i\},\chi,s\chi'} + (Q-1)\chi'_{i+1}(-1) \hat T_{w,\beta\cup\{i+1\},\chi,s\chi'}.
\]
\item if $w(i)=M+1$ and $i+1\notin \beta\sqcup w^{-1}([M+1,N])$,
\[
\hat T_{\tilde w,\chi,\chi'} T_{s,\chi',s\chi'} = \hat T_{ws,s\beta\cup\{i\},\chi,s\chi'} + (Q-1)\chi'_{i}(-1) \hat T_{w,\beta\cup\{i+1\},\chi,s\chi'}.
\]
\end{enumerate}
\end{enumerate}
\end{prop}
\begin{proof}
We only sketch the calculation. Let $e_{\check\alpha_i}(a)=\id+a e_i\otimes e_{i+1}^*$, where $\check\alpha_i$ is the $i$-th simple root of $\GL_N$ for $B_N$. For a right $(B_N,\chi')$-equivariant function $f$, by the decomposition $B_NsB_N=\bigsqcup_{a\in\BF_Q}B_Ns e_{\check\alpha_i}(a)$, we have
\begin{equation}\label{eq:finite-rank-one-sum}
(fT_{s,\chi',s\chi'})(x) =\sum_{a\in\BF_Q}f\bigl(xe_{\check\alpha_i}(-a)s\bigr).
\end{equation}
\begin{itemize}
\item In the case (1), since $ws>w$ and $\tilde{w}$ is relevant, we have $w(i)<M+1$. The image of the convolution map contains only one orbit unless $\iota\subseteq \sigma^{',*}$. If $\iota\subseteq \sigma^{',*}$, there are two orbits $(ws, s\beta)$ and $(ws,s\beta\setminus\{i+1\})$. To determine the coefficients: if $w(i+1)\neq M+1$, then the coefficient of $\hat T_{ws,s\beta}$ is $\sum_{\BF_Q-\{1\}}\chi_{w(i)}((1-a)^{-1})(s\chi')_{i+1}(1-a)=Q-1$. Here, we use $i\notin \sigma$, which implies $\chi_{w(i)}=\chi'_i=(s\chi')_{i+1}$. If $w(i+1)= M+1$, then the coefficient of $\hat T_{ws,s\beta}$ is $\sum_{\BF_Q-\{-1\}}\chi_{w(i)}((1+a)^{-1})(s\chi')_{i+1}(1+a)=Q-1$. It is easy to see that the coefficient of $\hat T_{ws,s\beta\setminus\{i+1\},\chi,s\chi'}$ is $1$;
\item In the case (2), we have $w(i+1)<w(i)<M+1$, and the coefficient of $\hat T_{w,\beta}$ is $\sum_{\BF_Q-\{0\}}\chi_{w(i)}(a^{-1})(s\chi')_{i+1}(-a)=(Q-1) \chi'_i(-1)$. Here, we use $i\notin \sigma$, which implies $\chi_{w(i)}=\chi'_i=(s\chi')_{i+1}$;
\item In the case (3.a), the orbit $({ws,s\beta})$ is not relevant, and the coefficient of $\hat T_{w,\beta}$ is $\sum_{\BF_Q^\times}\psi(-a^{-1})(s\chi')_i(a^{-1})(s\chi')_{i+1}(-a)=\chi'_{i+1}(-1)G(\psi,\chi'_i/\chi'_{i+1})$;
\item In the case (3.b), if $w(i)=M+1$, then the coefficient of $\hat T_{ws,s\beta}$ is $1$, and the coefficient of $\hat T_{w,\beta\setminus\{i+1\}}$ is $\chi'_{i+1}(-1)$; the coefficient of the term $\hat T_{w,\beta}$ is $\sum_{\BF_Q-\{0,1\}}\chi_{i+1}'(-a^{-1})\chi_{w(i+1)}(a/(a-1))\chi_i'(a-1)=\chi_{w(i+1)}(-1)\chi_i'(-1)/\chi'_{i+1}(-1)\sum_{\BF_Q-\{0,1\}}(\chi_{w(i+1)}/\chi_{i+1}')(a) \cdot(\chi_i'/\chi_{w(i+1)})(1-a)=(\chi_{w(i+1)}\chi_i'/\chi'_{i+1})(-1) J(\chi_{w(i+1)}/\chi_{i+1}', \chi_i'/\chi_{w(i+1)})$;
\item In the case (3.c), according to the $\psi$-relevance condition, we have $w(i),w(i+1)\in [1,M]$. Let $A=\chi'_i$, $B=\chi'_{i+1}$, and $C=\chi_{w(i+1)}$. The coefficient of $\hat T_{w,\beta}$ is
\[
\fitmath{A(-1)\sum_{a\in\BF_Q\setminus\{0,1\}}(A/B)(a)(A/C)(1-a^{-1}) =A(-1)J(B/A,A/C) =C(-1)J(C/B,A/C).}
\]
The first equality follows from the substitution $a\mapsto a^{-1}$, and the second one from the standard identities for Jacobi sums.
\item In the case (4), we have $w(i+1)<M+1$;
\item In the case (4), if $w(i)\neq M+1$, and $i+1\notin \beta$, then the coefficient of $\hat T_{w,\beta\cup\{i+1\}}$ is $\sum_{\BF_Q-\{0\}}\chi_{w(i+1)}(-a)(s\chi')_{i}(a^{-1})=(Q-1) \chi'_{i+1}(-1)$. Here, we use $i+1\notin \sigma$, which implies $\chi_{w(i+1)}=\chi'_{i+1}$;
\item In the case (4), if $w(i)\neq M+1$, and $i+1\in \beta$, the coefficient of $\hat T_{w,\beta\cup\{i+1\}}=\hat T_{w,\beta}$ is also $\sum_{\BF_Q-\{0\}}\chi_{w(i+1)}(-a)(s\chi')_{i}(a^{-1})=(Q-1) \chi'_{i+1}(-1)$;
\item In the case (4), if $w(i)=M+1$, and $i+1\in \beta$, then the coefficient of $\hat T_{w,\beta\cup\{i+1\}}$ is $\sum_{\BF_Q-\{0\}}(s\chi')_i(a^{-1})\chi_{w(i+1)}(a)(s\chi')_{i+1}(-1)=(Q-1) \chi'_{i+1}(-1)$;
\item In the case (4), if $w(i)=M+1$, and $i+1\notin \beta$, then the coefficient of $\hat T_{w,\beta\cup\{i+1\}}$ is $\sum_{\BF_Q-\{0\}}(s\chi')_i(a^{-1})(s\chi')_{i+1}(-1)\chi_{w(i+1)}(a)=(Q-1) \chi'_{i}(-1)$.
\end{itemize}
\end{proof}

\subsubsection{} Let $S_{w,\beta}=|\Stab_H(g_{w(\beta)}wB_N)|$. By the orbit dimension formula,
\begin{equation}\label{eq:finite-stabilizer-normalization}
S_{w,\beta} =Q^{\binom N2-l(w)-|\beta|+|\sigma|}(Q-1)^{M-|\sigma|}, \qquad \hat T_{w,\beta}=S_{w,\beta}T_{w,\beta}.
\end{equation}
Thus, multiplying the coefficient of $\hat T_{\tilde v}$ in Proposition~\ref{9.3.4} by $S_{\tilde v}/S_{\tilde w}$, we obtain the following formulas for $T_{\tilde w}$.
\begin{cor}\label{9.3.4+}
Fix $\tilde w = (w,\beta) \in RB^\fin_{\chi,\chi'}(M|N)$ and $s = (i,i+1)$, and denote $\sigma = \sigma(\tilde w)$ and $\sigma' = \sigma(ws,s\beta)$. Let $\sigma^*=\sigma\cup\{w^{-1}(M+1)\}$ and $\sigma^{',*}=\sigma'\cup\{(ws)^{-1}(M+1)\}$.

\begin{enumerate}
\item If $ws > w$ then
\[
T_{\tilde w,\chi,\chi'} T_{s,\chi',s\chi'} =
\begin{cases}
T_{\tilde{w}s,\chi,s\chi'}, & \iota \nsubseteq  \sigma^{',*}, \\[6pt]
T_{\tilde{w}s,\chi,s\chi'} + T_{ws,s\beta\setminus\{i+1\},\chi,s\chi'}, & \iota\subseteq  \sigma^{',*}.
\end{cases}
\]
\item If $ws < w$ and $i \in w^{-1}([1,M])\setminus \beta$ then
\[
T_{\tilde w,\chi,\chi'} T_{s,\chi',s\chi'} = Q T_{\tilde{w}s,\chi,s\chi'} + (Q-1)\chi'_{i+1}(-1) T_{\tilde w,\chi,s\chi'}.
\]
\item If $ws < w$, $i \in \beta\sqcup w^{-1}([M+1,N])$ and $i+1 \in \sigma\sqcup w^{-1}([M+1,N])$,
\begin{enumerate}
\item  if $M+1\leq w(i+1)=w(i)-1\leq N-1$, then
\[
T_{\tilde w,\chi,\chi'} T_{s,\chi',s\chi'} = \chi'_{i+1}(-1)G(\psi,\chi'_i/\chi'_{i+1}) T_{\tilde w,\chi,s\chi'};\]
\item if $w(i)=M+1$, then \[
\begin{aligned}
T_{\tilde w,\chi,\chi'} T_{s,\chi',s\chi'} &= (Q-1) T_{ws,s\beta,\chi,s\chi'}\\
&{}+ (Q-1)\chi'_{i+1}(-1)T_{w,\beta\setminus\{i+1\},\chi,s\chi'}\\
&{}+ \frac{\chi_{w(i+1)}\chi_i'}{\chi'_{i+1}}(-1) J\!\left(\frac{\chi_{w(i+1)}}{\chi_{i+1}'}, \frac{\chi_i'}{\chi_{w(i+1)}}\right) T_{\tilde w,\chi,s\chi'}.
\end{aligned}
\]
\item otherwise,
\[
\begin{aligned}
T_{\tilde w,\chi,\chi'} T_{s,\chi',s\chi'} &= (Q-1) T_{ws,s\beta,\chi,s\chi'}\\
&{}+ (Q-1)\chi_{w(i+1)}(-1) T_{w,\beta\setminus\{i+1\},\chi,s\chi'}\\
&{}+ \chi_{w(i+1)}(-1) J\!\left( \frac{\chi_{w(i+1)}}{\chi_{i+1}'}, \frac{\chi_i'}{\chi_{w(i+1)}}\right) T_{\tilde w,\chi,s\chi'}.
\end{aligned}
\]
\end{enumerate}

\item If $ws < w$, $i \in \beta\sqcup w^{-1}([M+1,N])$ and $i+1 \notin \sigma\sqcup w^{-1}([M+1,N])$,
\begin{enumerate}
\item if $i+1\in \beta\sqcup w^{-1}([M+1,N])$,
\[
T_{\tilde w,\chi,\chi'} T_{s,\chi',s\chi'} = Q T_{ws,s\beta,\chi,s\chi'} + (Q-1)\chi'_{i+1}(-1) T_{w,\beta,\chi,s\chi'}.
\]
\item if $w(i)\neq M+1$ and $i+1\notin \beta\sqcup w^{-1}([M+1,N])$,
\[
T_{\tilde w,\chi,\chi'} T_{s,\chi',s\chi'} = T_{ws,s\beta\cup\{i\},\chi,s\chi'} + \chi'_{i+1}(-1) T_{w,\beta\cup\{i+1\},\chi,s\chi'}.
\]
\item if $w(i)=M+1$ and $i+1\notin \beta\sqcup w^{-1}([M+1,N])$,
\[
T_{\tilde w,\chi,\chi'} T_{s,\chi',s\chi'} = T_{ws,s\beta\cup\{i\},\chi,s\chi'} + \chi'_{i}(-1) T_{w,\beta\cup\{i+1\},\chi,s\chi'}.
\]
\end{enumerate}
\end{enumerate}
\end{cor}

\subsubsection{Generators in the non-monodromic case} Let us temporarily restrict our attention to the usual (non-monodromic) setting, and study $\bK^{\fin,\mxd}(M|N)$ and $\bK^{\fin,\gr}(M|N)$, in particular to determine the bar involution.

In the case of $M=N$, a collection of generators of the Grothendieck group was presented in \cite[Section 4.5]{[T1]}. Namely, consider the orbits corresponding to $(w,\beta)$, where $w=\id$, and $\beta=\{1,2,\cdots,i\}\subset \{1,2,\cdots,N\}$, where $i\in [0,N]$. The closure of the corresponding orbit in $\GL_N/B_N\times \BA^N$ is isomorphic to $\BA^i$, and the IC sheaves are constant. Thus, the corresponding canonical basis is easy to describe. By \cite[Corollary 2]{[T1]}, the classes of these IC sheaves form a set of generators of $\bK^{\fin,\mxd}(N|N)$, and thus, one can determine the involution of $\bK^{\fin,\mxd}(N|N)$.

In this section, we will study the generators of $\bK^{\fin,\mxd}(M|N)$ as a bimodule over $\bH^{\fin,\mxd}(M)\times \bH^{\fin,\mxd}(N)$, and the key observation is that it has a collection of clean generators.

\begin{prop}\label{9.3.5}
Given $M<N$, $\bK^\fin(M|N)_{\BC}$ is generated by $T_{\tilde w_k}$, $0\leq k\leq M$, as a bimodule over $\bH^\fin(M)_{\BC}\times\bH^\fin(N)_{\BC}$. The same holds for the graded and mixed versions, with generators $T^{\gr}_{\tilde w_k}$ and $T^{\mxd}_{\tilde w_k}$, respectively.
\begin{samepage}
Here, $\tilde{w}_k=(w,\beta)$, where $w$ is composed of three monotonic blocks:
\[
w = \big( \underbrace{1, 2, \dots, k}_{\text{Block 1}}, \quad \underbrace{N, N-1, \dots, M+1}_{\text{Block 2}}, \quad \underbrace{k+1, k+2, \dots, M}_{\text{Block 3}} \big),
\]
and $\beta=\{1, 2, \dots, k\}$.\par
\end{samepage}
\end{prop}
\begin{proof}
Let us do induction on the dimension of the corresponding orbit. Now, we assume that, for the orbits of dimension strictly smaller than $n$, we already know that $T_{\tilde{w}}$ belongs to the submodule generated by $\{T_{\tilde{w}_k}\}$.

Assume the dimension of the orbit corresponding to $(w,\beta)$ is $n$, and assume $l(ws_i)<l(w)$ for some $1\leq i\leq N-1$. If $\tilde{w}=(w,\beta)$ can be written as $\tilde{w}'s_i$, such that $\tilde{w}'=(ws_i,s_i\beta)$ is relevant, then one can generate $T_{\tilde{w}}$ by the formula $(1)$ in Corollary \ref{9.3.4+}. Indeed, the orbit corresponding to $(w,\beta\setminus \{i+1\})$ has dimension strictly smaller than the orbit corresponding to $\tilde{w}$. The bad cases are when $i+1\notin \beta$ and $i\in \beta\cup\{w^{-1}(M+1)\}$ (equivalently, $i\in \sigma\cup\{w^{-1}(M+1)\}$), or when both $w(i)$ and $w(i+1)$ belong to $\{M+1,\cdots,N\}$. In the first case, one cannot write $\tilde{w}$ as $\tilde{w}'s_i$, since $s_i\beta$ does not satisfy the saturation property; and in the second case, $\tilde{w}'$ does not satisfy the $\psi$-relevance condition, since the $\psi$-relevance condition requires  $w^{-1}(M+1) > w^{-1}(M+2) > \dots > w^{-1}(N)$. These are the only two cases in which we cannot further reduce $T_{\tilde{w}}$ to orbits of smaller dimension by formula $(1)$ via the right Hecke action in the non-monodromic case.

Similarly, for the left action, if we assume $1\leq i\leq M-1$ and $l(s_iw)<l(w)$, then the bad case is when $w^{-1}(i)\in \beta$ and $w^{-1}(i+1)\notin \beta$. Otherwise, $\tilde{w}=s_i\tilde{w}'$, and the corresponding left rank-one convolution calculation, with row indices in place of column indices, implies that we can generate $T_{\tilde{w}}$ by the functions corresponding to smaller orbits.

We need to classify all $(w,\beta)$ satisfying the following conditions.
\begin{enumerate}
\item For any $1 \leq i \leq M-1$, if $w^{-1}(i) > w^{-1}(i+1)$, then $w^{-1}(i) \in \beta$ and $w^{-1}(i+1) \notin \beta$. \item For any $1 \leq i \leq N-1$, if $w(i) > w(i+1)$ and $w(i+1) \in \{1, 2, \dots, M\}$, then $i \in \beta\cup\{w^{-1}(M+1)\}$ and $i+1 \notin \beta$. \end{enumerate}

The claim is that any $(w,\beta)$ of the above form is given by $\tilde{w}_k$, where $k\in \{0,1,\cdots,M\}$.

The condition $(1)$ and the saturation condition of $\beta$ imply that $w(\beta)=\{1,2,\cdots,k\}$, for some $k\in [0,M]$. The $\psi$-relevance condition and the condition (2) above imply that $w^{-1}(M+1), w^{-1}(M+2),\cdots, w^{-1}(N)$ is a decreasing sequence and each step decreases by $1$.

We claim that the set $\beta$ is an interval starting from $1$. If $j+1\in\beta$, but $j\notin \beta$, then the saturation condition $((i<j, w(i)<w(j), j\in\beta)\Rightarrow i\in \beta)$ implies that we have $w(j+1)<w(j)$. It contradicts the condition $(2)$, which requires that $j+1\notin \beta$. In particular, any element in $\{1,2,\cdots,N\}$ which is strictly bigger than $w^{-1}(M+1)$ does not belong to $\beta$.

Thus, $\beta=w(\beta)=[1,k]$, and $\{1,2,\cdots, N\}$ is divided into three segments $\{1,2,\cdots, w^{-1}(N)-1\}$, $\{w^{-1}(N), \cdots, w^{-1}(M+1)\}$, and $\{w^{-1}(M+1)+1,\ldots,N\}$, where the first one belongs to $\beta$, and the second and the third ones belong to the complement of $\beta$.

Finally, to avoid a contradiction with $(1)$ and $(2)$, we should have that $w$ is increasing on the segments $\{1,2,\cdots, w^{-1}(N)-1\}$ and $\{w^{-1}(M+1)+1, \cdots, N\}$. It implies that $k=w^{-1}(N)-1$ and our $(w,\beta)=\tilde{w}_k$.
\end{proof}

\begin{prop}\label{prop 9.3.6}
Every $\tilde{w}_{k}$ is minimal relevant. In other words, there is no other relevant orbit in the closure of the orbit corresponding to $\tilde{w}_{k}$.
\end{prop}
Before proving Proposition \ref{prop 9.3.6}, let us study the closure relation of the $B_M\ltimes U_{M,N}$-orbits in $\GL_N/B_N$.

\subsubsection{The closure relation} The closure relation of $B_N$-orbits in $\GL_N/B_N\times\BA^N$ is studied in \cite{[Mag]}, which we briefly recall here.

\begin{enumerate}
\item  Let $w^{[i,j]}$ be the $i \times j$ submatrix of the permutation matrix $w$ consisting of the bottom $i$ rows and the first $j$ columns, and denote its rank by
\[ r_{i,j}(w) = \text{rank} \left( w_{p,s} \right)_{N-i+1 \le p \le N, 1 \le s \le j} \]

\item  Let $v_{w(\beta)}^{[i]}$ be the sub-vector consisting of the last $i$ components of $v_{w(\beta)}$. We define
\[ \tilde{r}_{i,j}(w, \beta) = \text{rank} \left[ w^{[i,j]} \mid v_{w(\beta)}^{[i]} \right] \] where the matrix is augmented by the partial vector $v_{w(\beta)}^{[i]}$ as an additional column.
\end{enumerate}
\begin{thm}
The orbit corresponding to $(w,\beta)$ lands in the closure of the orbit corresponding to $(w',\beta')$ if and only if
\begin{enumerate}
\item \[r_{i,j}(w)\leq r_{i,j}(w'), \forall 1\leq i,j\leq N;\]
\item \[ \tilde{r}_{i,j}(w, \beta) \le \tilde{r}_{i,j}(w', \beta') \quad \forall \ 0 \le i \le N, \ 0 \le j \le N. \]
\end{enumerate}
\end{thm}

The closure relation above gives rise to the closure relation for the general $M,N$ case.

Given $M<N$, consider the subspace
\begin{equation}
\GL_N/B_N\times (\BG_m\times \BA^{M})\subset \GL_N/B_N\times \BA^N.
\end{equation}
Here, we regard $\BG_m\times \BA^M$ as a subspace of $\BA^N$ by embedding into the first $M+1$ coordinates, and $\BG_m=\BA^1\setminus \{0\}$ is embedded into the $(M+1)$-st coordinate.

From the definition, the above embedding is an embedding of $B_N$-varieties. Furthermore, since $B_N$ acts transitively on $\BG_m\times \BA^{M}$, we have the following lemma.
\begin{lem}
The embedding
\begin{equation}
\GL_N/B_N\simeq \GL_N/B_N\times e_{M+1}\hookrightarrow    \GL_N/B_N\times (\BG_m\times \BA^{M})
\end{equation}
induces a bijection between the $B_N$-orbits in $\GL_N/B_N\times (\BG_m\times \BA^{M})$ and the $\Stab_{B_N}(e_{M+1})\simeq B_M\ltimes B_{M,N}$-orbits in $\GL_N/B_N$, and preserves the order relation.

Here, $B_{M,N}=T_{N-M-1}\ltimes U_{M,N}$. In particular, the torus of $B_M\ltimes B_{M,N}$ is $T_M\times T_{N-M-1}\subset T_N$, which consists of all the diagonal coordinates except the $(M+1)$-st one.
\end{lem}
Furthermore, according to the classification of $B_M\ltimes U_{M,N}$-orbits in $\GL_N/B_N$ in Lemma~\ref{lem 9.2.1+}, we have the following lemma.
\begin{lem}
Any $B_M\ltimes U_{M,N}$-orbit in $\GL_N/B_N$ is also a $B_M\ltimes B_{M,N}$-orbit.
\end{lem}

Thus, we have the following description of the closure relation of $B_M\ltimes U_{M,N}$-orbits in $\GL_N/B_N$.
\begin{prop}
In the case of $M<N$, the $B_M\ltimes U_{M,N}$-orbit corresponding to $(w,\beta)$ lands in the closure of the orbit corresponding to $(w',\beta')$ if and only if
\begin{enumerate}
\item \[r_{i,j}(w)\leq r_{i,j}(w'), \forall 1\leq i,j\leq N;\]
\item \[ \tilde{r}_{i,j}(w, \beta) \le \tilde{r}_{i,j}(w', \beta') \quad \forall \ 0 \le i \le N, \ 0 \le j \le N. \]
\end{enumerate}
Here, we regard $(w,\beta)\in RB^\fin(M|N)$ as an element in $RB^\fin(N|N)$ by adding $w^{-1}(M+1)$ to $\beta$.
\end{prop}
\subsubsection{The minimal relevant orbits} Now, we can prove Proposition~\ref{prop 9.3.6}, using the augmented matrix of $\tilde w_k=(w_k,\beta_k)$. Let $\beta_k^*=\beta_k\cup\{w_k^{-1}(M+1)\}$, and then the augmented matrix is
\[
(w_k \mid v_{w_k(\beta_k^*)}) = \left(
\begin{array}{ccc|c}
\begin{matrix}
1 & & \\
& \ddots & \\
& & 1
\end{matrix} & \mathbf{0} & \mathbf{0} &
\begin{matrix}
1 \\
\vdots \\
1
\end{matrix} \\
\hline \mathbf{0} & \mathbf{0} &
\begin{matrix}
1 & & \\
& \ddots & \\
& & 1
\end{matrix} &
\begin{matrix}
0 \\
\vdots \\
0
\end{matrix} \\
\hline \mathbf{0} &
\begin{matrix}
& & 1 \\
& \reflectbox{$\ddots$} & \\
1 & &
\end{matrix} & \mathbf{0} &
\begin{matrix}
1 \\
0 \\
\vdots \\
0
\end{matrix}
\end{array}
\right)
\begin{array}{l}
\text{Rows } 1 \dots k \\ \\
\text{Rows } k+1 \dots M \\ \\
\text{Rows } M+1 \dots N.
\end{array}
\]
\begin{proof}
Suppose that a relevant orbit $(w,\beta)$ is contained in the closure of the orbit corresponding to $\tilde w_k$. Let $d=N-M$ and $c=k+d$. For $j\leq k$, by the rank inequality $r_{N-j,j}(w)\leq r_{N-j,j}(w_k)=0$, we have $w(j)=j$. Then, we have $r_{N-k,c}(w)=d$, while $\widetilde r_{N-k,c}(w_k,\beta_k^*)=d$. Hence, the vector in rows $k+1,\ldots,N$ lies in the span of the first $c$ columns, and since its $(M+1)$-st coordinate is non-zero, we have $w^{-1}(M+1)\leq c$. By the relevance condition, we have
\[
k<w^{-1}(N)<\cdots<w^{-1}(M+1)\leq c.
\]
These $d$ positions fill $[k+1,c]$. Then, the remaining rank inequalities force the remaining values $k+1,\ldots,M$ to occur increasingly, thus, $w=w_k$. Finally, the saturation condition with $w^{-1}(M+1)=c$ forces $[1,k]\subset\beta$, and the same augmented rank inequality excludes any additional vector coordinate in rows $k+1,\ldots,M$. Thus, $\beta=[1,k]$.
\end{proof}
\subsubsection{Involution on the finite non-monodromic Grothendieck group} Similar to $\bK^{\fin,\mxd}(M|N)$, $\bK^{\fin,\gr}(M|N)$ and $\bK^{\fin,\gr}(M|N)_{\BC}$, one can also define the mixed and graded versions of the Grothendieck group of the finite Hecke category, $\bH^{\fin,\mxd}(N)$, $\bH^{\fin,\gr}(N)$ and $\bH^{\fin,\gr}(N)_{\BC}$. They are modules over $\bR^{\mathrm{mxd}}=K^\mxd_0(\pt)$, $\BZ[\bv,\bv^{-1}]$, and $\BC[\bv,\bv^{-1}]$, respectively. The latter two are (different kinds of) specializations of $\bH^{\fin,\mxd}(N)$, which is the Grothendieck group of $\cH^{\fin,\mxd}(N):= D^{B_N}_{\mxd}(\GL_N/B_N)$.

Let $\bar{\BQ}_\ell(1)$ denote the Tate twist over $\Spec(\BF_Q)$, and let \[ [v]=[\sqrt{Q}]:=\big[\bar{\BQ}_\ell\!\left(-\tfrac12\right)\big]\in \bR^{\mathrm{mxd}}.
\]

For a $B$-orbit corresponding to ${w}$, we denote by $\Delta_{{w}}$ the $!$-extension of the constant sheaf shifted by $\langle l({w})\rangle$ (such that the local system is pure of weight $0$ and perverse), and we define $\IC_{{w}}$ and $\nabla_{{w}}$ similarly.

Under the characteristic map, (the specialization of) the $!$-extension of the constant sheaf corresponds to $T_{w}$, and we denote the corresponding elements in the other Grothendieck groups by $T^\mxd_{w}$ and $T^\gr_{w}$.

Denote by $H^\mxd_w$ and $\underline{H}^\mxd_w$ the corresponding elements of $\Delta_{{w}}$ and $\IC_{{w}}$ in the mixed Grothendieck group. Similarly, their specializations are denoted by $H^\gr_w$, $\underline{H}^\gr_w$, $H_w$, $\underline{H}_w$.

In $\bR^{\mathrm{mxd}}$, a shift by $[1]$ induces a multiplication by $-1$ due to the alternating sum of cohomological degrees, and a twist by $(1/2)$ corresponds to multiplying the Frobenius eigenvalues by $Q^{-1/2}$, which maps to $v^{-1}$ under the specialization map to $\BC$. Thus, we have \[H^\mxd_w=(-[v])^{-l(w)}T^\mxd_w,\] and similar relations for the other versions. The finite mixed Hecke algebra is generated by $H^\mxd_s$ for simple reflections $s\in W^\fin$, and satisfies the relations $H^\mxd_w H^\mxd_{w'}=H^\mxd_{ww'}$, if $l(w)+l(w')=l(ww')$, and $H^\mxd_s H^\mxd_s= H^\mxd_e+([v^{-1}]-[v])H^\mxd_s$.
\subsubsection{The Hecke bar involution}\label{10.2.7} The Verdier duality on $D^{B_N}_{\mxd}(\GL_N/B_N)$ induces a bar involution
\begin{equation}
\bar{\cdot}\colon \bH^{\fin,\mxd}(N)\rightarrow \bH^{\fin,\mxd}(N).
\end{equation}
The bar involution changes $[\alpha]\in \bR^\mxd$ to $[\alpha^{-1}]\in \bR^\mxd$.

Let $\bR^{\mxd}_{<0}\subset \bR^{\mxd}$ denote the additive subgroup spanned by classes of strictly negative weight, i.e., those $[\alpha]$, where $\alpha\in \bar{\BQ}_\ell^\times$ and its absolute value belongs to $\{Q^{\frac{m}{2}}|m\in \BZ_{<0}\}$. Then, the canonical basis $\underline{H}_w^\mxd$ is characterized by the bar-invariance and the triangularity with respect to the standard basis ${H}_w^\mxd$, i.e., $\underline{H}^\mxd_{w}$ is the unique element satisfying
\begin{equation}\label{bar h mxd}
\overline{\underline{H}^\mxd_{w}} = \underline{H}^\mxd_{w}, \qquad \underline{H}^\mxd_{w} \equiv H^\mxd_{w} \pmod{\sum_{y<w} \bR^{\mathrm{mxd}}_{<0}\, H^\mxd_{y}}.
\end{equation}

In particular, for a simple reflection $s=(i,i+1)\in W_N^{\fin}$, the canonical basis element in the mixed Grothendieck group is given by
\[
\underline{H}^\mxd_{s} = H^\mxd_{s}-[v^{-1}]\cdot H^\mxd_{e}.
\]

By taking specialization, we obtain the bar involution
\begin{equation}
\bar{\cdot}\colon \bH^{\fin,\gr}(N)\rightarrow \bH^{\fin,\gr}(N),
\end{equation}
and the canonical basis $\underline{H}_w^\gr$ for the graded Grothendieck group, satisfying
\begin{equation}\label{bar h gr}
\overline{\underline{H}^\gr_{w}} = \underline{H}^\gr_{w}, \qquad \underline{H}^\gr_{w} \equiv H^\gr_{w} \pmod{\sum_{y<w} \bv^{-1}\BZ[\bv^{-1}]\, H^\gr_{y}}.
\end{equation}

\subsubsection{The mixed action}\label{sec 10.2.9} The convolution gives rise to a functor
\[\cH^{\fin,\mxd}(M)\otimes\cC^{\fin,\mxd}(M|N) \otimes\cH^{\fin,\mxd}(N)\rightarrow \cC^{\fin,\mxd}(M|N),\] and passing to Grothendieck groups, the above functor gives rise to a bimodule structure
\begin{equation}\label{10.26}
\bH^{\fin,\mxd}(M)\otimes\bK^{\fin,\mxd}(M|N) \otimes\bH^{\fin,\mxd}(N)\rightarrow \bK^{\fin,\mxd}(M|N).
\end{equation}

Furthermore, if one takes the specialization $\mathrm{sp}$, it gives rise to the (non-monodromic case of the) formulas in Corollary~\ref{9.3.4+},
\begin{equation}
\bH^{\fin}(M)_\BC\otimes\bK^{\fin}(M|N)_\BC \otimes\bH^{\fin}(N)_\BC\rightarrow \bK^{\fin}(M|N)_\BC.
\end{equation}
Note that the coefficients in the above formula are given by the sum of the eigenvalues of the Frobenius that appear in the coefficients of the formula for the mixed Grothendieck group.

By considering the above formulas for $\BF_{Q^n}$ for all $n\geq 1$, we recover the formula at the mixed Grothendieck group level, i.e., \eqref{10.26}. Indeed, if we have a $k$-dimensional space with eigenvalues $\alpha_1,\cdots, \alpha_k$, then the trace of the $n$-th power of Frobenius will be the sum of the $n$-th powers of $\alpha_1,\cdots, \alpha_k$. From this information for all $n$, one can extract all the eigenvalues, with their multiplicities.

\subsubsection{The canonical basis} Similar to the Hecke algebra, we also take the (renormalized) standard bases $H_{\tilde{w}}^\mxd$, $H_{\tilde{w}}^\gr$, and $H_{\tilde{w}}$. Under \funshf, they are given by (the specialization of) the $!$-extension of the rank-one equivariant perverse local system, pure of weight $0$. In particular, we have $H_{\tilde{w}}^\mxd=(-[v])^{-l(\tilde{w})}T_{\tilde{w}}^\mxd$, and similarly for $H_{\tilde{w}}^\gr$ and $H_{\tilde{w}}$. Here, $l(\tilde w)$ denotes the dimension of the corresponding finite type orbit.

Consider the Verdier duality on $\cC^{\fin,\mxd}(M|N)$. In this non-monodromic case, we identify the Grothendieck groups for $\psi$ and $-\psi$ by their normalized standard bases. Then, the Verdier duality induces a bar involution
\begin{equation}
\bar{}\colon \bK^{\fin,\mxd}(M|N)\rightarrow \bK^{\fin,\mxd}(M|N).
\end{equation}

Also, we have the canonical bases $\underline{H}_{\tilde{w}}^\mxd$ and $\underline{H}_{\tilde{w}}^\gr$, which are characterized by the requirements similar to \eqref{bar h mxd} and \eqref{bar h gr}. We should give a combinatorial way to determine the bar involution and the canonical basis in terms of the standard basis.

Such a description is actually a direct corollary of the generating property, Proposition~\ref{9.3.5}, and the minimal relevance condition of the orbits corresponding to $\tilde{w}_k$, Proposition~\ref{prop 9.3.6}.

\begin{cor}\label{cor 9.3.12}
\begin{enumerate}
\item There is a unique involution on $\bK^{\fin,\mxd}(M|N)$, semilinear for the involution $[\alpha]\mapsto[\alpha^{-1}]$, such that, for any $0\leq k\leq M$, we have
\begin{equation}
\begin{split}
\overline{\underline{H}}^\mxd_{(w_k, \beta_k)}={\underline{H}}^\mxd_{(w_k, \beta_k)},
\end{split}
\end{equation}
and \begin{equation}
\begin{split}
\overline{h_M\cdot r}=\overline{h_M}\cdot \overline{r},\\
\overline{r\cdot h_N}=\overline{r}\cdot \overline{h_N},
\end{split}
\end{equation}
for any $r\in \bK^{\fin,\mxd}(M|N)$, $h_M\in \bH^{\fin,\mxd}(M)$, and $h_N\in \bH^{\fin,\mxd}(N)$.
\item The involution in (1) is induced by the Grothendieck--Verdier duality on $\GL_N/B_N$.
\end{enumerate}
\end{cor}
\begin{proof}
It follows immediately from Propositions~\ref{9.3.5} and \ref{prop 9.3.6}.
\end{proof}

\subsection{Determining the isomorphism in Theorem \ref{Thm 9.2.5}} On the level of Grothendieck groups, Theorem~\ref{Thm 9.2.5} gives the following corollary.
\begin{cor}\label{cor 10.2.14}
There is an isomorphism
\begin{equation}
\bK^{\fin,\mxd}_{\chi,\chi'}(M|N)\simeq \bigotimes_{i=1}^l\bK^{\fin,\mxd}_{1,1}(k_i|k_i').
\end{equation}
Furthermore, the Verdier duality functor induces a bar involution
\begin{equation}
\bar{\cdot}\colon \bK^{\fin,\mxd}_{\chi,\chi'}(M|N)\rightarrow \bK^{\fin,\mxd}_{(\chi)^{-1},(\chi')^{-1}}(M|N),
\end{equation}
and it is compatible with the bar involution of the non-monodromic finite type Grothendieck groups, up to changing the characters to their inverses,
\[\xymatrix{ \bK^{\fin,\mxd}_{\chi,\chi'}(M|N)\ar[r]\ar[d]&\bK^{\fin,\mxd}_{(\chi)^{-1},(\chi')^{-1}}(M|N)\ar[d]\\
\bigotimes_{i=1}^l\bK^{\fin,\mxd}_{1,1}(k_i|k_i')\ar[r]& \bigotimes_{i=1}^l\bK^{\fin,\mxd}_{1,1}(k_i|k_i'). }\] Here, for the inverse characters, we take the dual kernels and auxiliary Weil structures, and the additive character $-\psi$. The isomorphism sends each canonical basis element to a canonical basis element tensored with a one-dimensional Weil representation pure of weight $0$.
\end{cor}
\begin{proof}
This follows from Theorem~\ref{Thm 9.2.5}, since two Weil structures on a geometric simple sheaf differ by a one-dimensional Weil representation.
\end{proof}

In this subsection, we want to give a combinatorial description of this isomorphism of Grothendieck groups, which gives the bar involution and the canonical basis of $\bK^{\fin,\mxd}_{\chi,\chi'}(M|N)$ in terms of the standard basis. According to the proof of Theorem~\ref{Thm 9.2.5}, we have to present the following isomorphism explicitly
\begin{equation}
\bK^{\fin,\mxd}_{\chi,\chi'}(M|N) \simeq \bK^{\fin,\mxd}_{\tau^{-1}\chi,\chi'}(M-1|N)
\end{equation}
given by Proposition~\ref{prop 9.2.9}, for $0<M\leq N$ and $\chi_1\neq\chi'_j$ for every $1\leq j\leq N$. For $M<N$ ($M=N$ is similar), this is the composition
\begin{equation}\label{9.28}
\begin{split}
K_0(D^{B_M,\chi\ltimes U_{M,N},\psi\times (B_N,\chi')}_\mxd(\GL_N))\\
\overset{\Av_!^{\BG_a,\psi'}\langle-1\rangle}{\simeq} K_0(D^{(B'_{M-1},\chi)\ltimes (V^*\oplus \BG_a,\psi')\ltimes U_{M,N},\psi\times (B_N,\chi')}_\mxd(\GL_N))\\
\overset{\FT}{\simeq}K_0(D^{(B'_{M-1},\chi)\ltimes (V\oplus \BG_a,\psi')\ltimes U_{M,N},\psi\times (B_N,\chi')}_\mxd(\GL_N))\\
\overset{}{\simeq}K_0(D^{(B_{M-1},\tau^{-1}\chi)\ltimes U_{M-1,N},\psi\times (B_N,\chi')}_\mxd(\GL_N)).
\end{split}
\end{equation}

Denote by $H^{(1)}$, $H^{(2)}$, $H^{(3)}$, and $H^{(4)}$ the standard bases of the above four Grothendieck groups, corresponding to the $!$-extensions of the equivariant pure of weight $0$ perverse local systems. For the same reason as in Section~\ref{sec 10.2.9}, it is enough to consider the isomorphism after taking the specialization $\mathrm{sp}$, namely, at the $\BC$-valued function level, over all finite extensions.

Before performing the isomorphisms, let us index the bases of these four Grothendieck groups. The relevant orbits in the first and the fourth categories have already appeared in Lemmas~\ref{lem 9.2.1} and \ref{lem 9.2.1+}. Similarly, we have the following lemma.
\begin{lem}
\begin{enumerate}
\item The relevant $(B'_{M-1},\chi)\ltimes (V^*\oplus \BG_a,\psi')\ltimes (U_{M,N},\psi)\times (B_N,\chi')$-orbits in $\GL_N$ are indexed by $(w,\beta)$, where $w\in W^\fin_N$, and $\beta$ is a subset of $w^{-1}([2,M])$. We require that if $i<_{B_N\cap w^{-1}B_Nw}j$, and $j\in \beta \cup\{w^{-1}(M+1)\}$, and $i\in w^{-1}([2,M])$, then $i\in \beta$. Also, if $i\in \{2,3,\cdots, M\}\setminus w(\sigma)$ ($\sigma$ is the set of elements in $\beta$ which are maximal with respect to the partial order $<_{B_N\cap w^{-1}B_Nw}$ in $\beta \cup\{w^{-1}(M+1)\}$), then $\chi_i=\chi'_{w^{-1}(i)}$. If $1\leq_{B_N\cap wB_Nw^{-1}}i$, then $i\notin w(\beta)$. Also, $w^{-1}(1)>w^{-1}(M+1)>w^{-1}(M+2)>\cdots>w^{-1}(N)$. The representative is given by $g_{w(\beta)}w=(\id+\sum_{a\in w(\beta)}e_a\otimes e_{M+1}^*)w$.
\item The relevant $(B'_{M-1},\chi)\ltimes (V\oplus \BG_a,\psi')\ltimes (U_{M,N},\psi)\times (B_N,\chi')$-orbits in $\GL_N$ are indexed by $(w,\beta)$, where $w\in W^\fin_N$, and $\beta$ is a subset of $[2,M]$.  We require that if $j\in \beta\cup \{1\}$, $i>_{B_N\cap wB_Nw^{-1}}j$, and $i\in [2,M]$, then $i\in \beta$. Also, if $i\in \{2,3,\cdots, M\}\setminus \sigma$ ($\sigma$ is the set of elements in $\beta$ which are minimal with respect to the partial order $<_{B_N\cap wB_Nw^{-1}}$ in $\beta\sqcup \{1\}$), then $\chi_i=\chi'_{w^{-1}(i)}$. If $i\leq_{B_N\cap wB_Nw^{-1}}M+1$, then $i\notin \beta$. Also, $w^{-1}(1)>w^{-1}(M+1)>w^{-1}(M+2)>\cdots>w^{-1}(N)$. The representative is given by $(\id+\sum_{a\in \beta}e_{1}\otimes e_{a}^*)w$.
\end{enumerate}
\end{lem}
\begin{proof}
The additive characters require $w^{-1}(1)>w^{-1}(M+1)>\cdots>w^{-1}(N)$. The remaining conditions follow from row operations and the torus stabilizer, as in Lemma~\ref{lem 9.2.1+}.
\end{proof}
\subsubsection{The first isomorphism in \eqref{9.28}}\label{sec 10.3.1} In the formulas below, $\psi$ is the fixed non-trivial additive character used in the central character $\psi'$. Recall that the relevant $(B_M,\chi)\ltimes (U_{M,N},\psi)\times (B_N,\chi')$-orbits in $\GL_N$ are classified by Lemmas~\ref{lem 9.2.1} and \ref{lem 9.2.1+}. We take the representative $(\id+\sum_{a\in w(\beta)}e_a\otimes e_{M+1}^*)w\in \GL_N$. Removing $\{1\}$ from $w(\beta)$ induces a bijection between the corresponding relevant orbits. Indeed, the assumption on the characters forces $1\in w(\sigma)$, and the inverse map adds $1$. By a direct calculation, we have the following lemma.
\begin{lem} $\Av_!^{\BG_a,\psi'}\langle-1\rangle(H^{(1)}_{(w,\beta)})=-1/\sqrt{Q}\cdot G(-\psi, \chi'_{w^{-1}(1)}/\chi_1) H^{(2)}_{(w,\beta\setminus\{w^{-1}(1)\})}$.
\end{lem}
\begin{proof}
Let $r=\chi_1/\chi'_{w^{-1}(1)}$. By the equivariance, the value of the characteristic function at the coefficient $a\neq0$ of $e_1\otimes e_{M+1}^*$ is $r(a)$. Thus, by averaging, we obtain
\[
-Q^{-1/2}\sum_{a\in\BF_Q^\times}r(a)\psi(-a) =-Q^{-1/2}G(-\psi,r^{-1}).
\]
Since the two orbits have the same dimension, the same coefficient occurs for the normalized standard objects.
\end{proof}

\subsubsection{The second isomorphism in \eqref{9.28}}\label{sec 10.3.2} It is a bit cumbersome to do the Fourier transform. To carry out the calculation, we use the Heisenberg model. Recall that we denote $\Hei=V\oplus V^*\oplus\BG_a$ in the proof of Lemma~\ref{9.2.7}, and identify it with an upper-left matrix subgroup $\begin{pmatrix}1&V^*&\BG_a\\
&\Id_{M-1}&V\\
&&1\end{pmatrix}$ of $\GL_N$.

The restriction along \[e\rightarrow \Hei\] induces an isomorphism
\begin{equation}
\begin{split}
\BC^{(B'_{M-1},\chi)\ltimes (V^*\oplus \BG_a,\psi')\ltimes (U_{M,N},\psi)\times (B_N,\chi')}(\GL_N)\\ \simeq\\
\BC^{(B'_{M-1},\chi)\ltimes \Hei\ltimes (U_{M,N},\psi)\times (V^*\oplus \BG_a,-\psi')}(\GL_N/(B_N,\chi')\times \Hei)\\ \simeq \\  \BC^{(B'_{M-1},\chi)\ltimes \Hei\ltimes (U_{M,N},\psi)}(\GL_N/(B_N,\chi')\times \Hei/(V^*\oplus \BG_a,-\psi')).
\end{split}
\end{equation}
Here, in the second function space, $V^*\oplus \BG_a$ acts on $\Hei$ from the right.

The standard object corresponding to $(\id+\sum_{a\in w(\beta)}e_a\otimes e_{M+1}^*)w$ goes to the standard object corresponding to the orbit passing through $(w, \id-\sum_{a\in w(\beta)}e_a\otimes e_{M+1}^*)\in \GL_N\times \Hei$.

Also,
\begin{equation}
\begin{split}
\BC^{(B'_{M-1},\chi)\ltimes (V\oplus \BG_a,\psi')\ltimes U_{M,N},\psi\times (B_N,\chi')}(\GL_N)\\ \simeq\\
\BC^{(B'_{M-1},\chi)\ltimes \Hei\ltimes (U_{M,N},\psi)}(\GL_N/(B_N,\chi')\times \Hei/(V\oplus \BG_a,-\psi')).
\end{split}
\end{equation}

The standard object corresponding to the orbit passing through $(\id+\sum_{a\in \beta}e_{1}\otimes e_{a}^*)w$ goes to the standard object corresponding to the orbit passing through $(w, \id-\sum_{a\in \beta}e_{1}\otimes e_{a}^*)\in \GL_N\times \Hei$.

With the Heisenberg models, the Fourier equivalence in \eqref{9.28} is the usual Fourier equivalence over $\GL_N/B_N$. Thus, it is enough to calculate the image of the standard object under the isomorphism induced by the Fourier equivalence \[D(\Hei/V^*\oplus \BG_a,-\psi')\simeq D(\Hei/V\oplus \BG_a,-\psi'),\] i.e., \[D(V)\simeq D(V^*).\]

Consider the fiber over $w\in \GL_N/B_N$. We should study
\begin{equation}
\begin{split}
\BC^{(T'_{M-1},\chi/w\chi')\ltimes (U_{M+1}\cap wU_Nw^{-1})}(V)\\
\overset{\FT}{\simeq} \\
\BC^{(T'_{M-1},\chi/w\chi')\ltimes (U_{M+1}\cap wU_Nw^{-1})}(V^*).
\end{split}
\end{equation}

We normalize the one-dimensional Fourier transform so that $\mathbb 1_{\{0\}}$ maps to $-Q^{-1/2}\mathbb 1_{\BA^1}$ and $\mathbb 1_{\BA^1}$ maps to $-Q^{1/2}\mathbb 1_{\{0\}}$. The restriction of the characteristic function corresponding to $(w,\beta)$ to $V$ (resp. $V^*$) is clear. Namely, under the natural identification \[V\simeq \BG_a\times \BG_a\times\cdots\times \BG_a,\] the characteristic function is the product of functions $T_i$ on each factor. To be more precise,
\begin{enumerate}
\item if $i\in w(\sigma)$, then $T_i$ is the $(\BG_m,\chi_i/\chi'_{w^{-1}(i)})$-invariant characteristic function (renormalized to have value $1$ at $-1\in\BG_a$) on $\BG_m$;
\item if $i\in w(\beta\setminus \sigma)$, then $T_i$ is the constant function on $\BG_a$;
\item if $i\notin w(\beta)$, then $T_i$ is the characteristic function on $\{0\}\subset \BG_a$.
\end{enumerate}
Note that when $i\notin w(\sigma)$, we have that $\chi_i/\chi'_{w^{-1}(i)}$ is trivial.

For $V^*\simeq\Hei/(V\oplus\BG_a)$, the torus acts through the inverse characters, and we replace $w(\sigma)$ and $w(\beta\setminus\sigma)$ by $\sigma$ and $\beta\setminus\sigma$ in the above description.

For each $\BG_a$, we denote by $T_b$ the equivariant characteristic function on (the big $\BG_m$-orbit) $\BA^1\setminus \{0\}$, renormalized to have value $1$ at $-1$, and by $T_{\{0\}}$ the equivariant characteristic function at (the small $\BG_m$-orbit) $\{0\}$ if it is relevant, and otherwise, it is set to be $0$.
\begin{lem}
\begin{enumerate}
\item if $i\in w(\sigma)$, then the Fourier transform of $T_i$ is $-{Q}^{-1/2}\cdot (G(-\psi, \chi'_{w^{-1}(i)}/\chi_i)T_b+ G(1, \chi'_{w^{-1}(i)}/\chi_i)T_{\{0\}})$;
\item if $i\in w(\beta\setminus \sigma)$, then the Fourier transform of $T_i$ is $-{Q}^{1/2} T_{\{0\}}$;
\item if $i\notin w(\beta)$, then the Fourier transform of $T_i$ is $-{Q}^{-1/2}\cdot(T_{\{0\}}+T_b)$.
\end{enumerate}
\end{lem}
\begin{proof}
Let $h(r,v,z)=\begin{pmatrix}1&r&z\\&\Id_{M-1}&v\\&&1\end{pmatrix}$. Then, the function on $V$ is $F_w(v)=f(h(0,-v,0)w)$, and the identity
\[
h(0,-v,0)h(-y,0,0)=h(-y,0,-yv)h(0,-v,0)
\]
gives the Fourier kernel $\psi(-yv)$. Let $r=\chi_i/\chi'_{w^{-1}(i)}$. For $x,y\neq0$, we have $T_i(x)=r(-x)$ and $T_b(y)=r^{-1}(-y)$, thus,
\[
\sum_{x\in\BF_Q^\times}r(-x)\psi(-xy) =G(-\psi,r^{-1})r^{-1}(-y).
\]
At $y=0$, the sum is $G(1,r^{-1})$. The other cases follow from the Fourier transforms of $\mathbb 1_{\BA^1}$ and $\mathbb 1_{\{0\}}$.
\end{proof}

To summarize, we have the following description of the Fourier isomorphism
\begin{prop}
$\FT(H^{(2)}_{w,\beta})=\sum_{\beta'}a_{\beta'} H^{(3)}_{w,\beta'}$, where $\beta'\subset[2,M]$ indexes the relevant orbits in the third category, and $a_{\beta'}=\prod_{i=2}^M b_i$, with
\[
b_i=\begin{cases} -Q^{-1/2}G(-\psi,\chi'_{w^{-1}(i)}/\chi_i), & i\in w(\sigma)\cap\beta',\\
Q^{-1}G(1,\chi'_{w^{-1}(i)}/\chi_i), & i\in w(\sigma),\ i\notin\beta',\\
1, & i\in w(\beta\setminus\sigma),\ i\notin\beta',\\
0, & i\in w(\beta\setminus\sigma)\cap\beta',\\
1, & i\notin w(\beta),\ i\in\beta',\\
-Q^{-1/2}, & i\notin w(\beta),\ i\notin\beta'.
\end{cases}
\]
\end{prop}

\subsubsection{The third isomorphism in \eqref{9.28}}\label{sec 10.3.3}

Let us denote $\tau=(1,2,\cdots,M)\in S_N$. As we mentioned, the left translation by $\tau^{-1}$ induces an isomorphism
\begin{equation}\label{10.42}
\BC^{(B'_{M-1},\chi)\ltimes (V\oplus \BG_a,\psi')\ltimes U_{M,N},\psi\times (B_N,\chi')}(\GL_N) \overset{}{\simeq}\BC^{(B_{M-1},\tau^{-1}\chi)\ltimes U_{M-1,N},\psi\times (B_N,\chi')}(\GL_N).
\end{equation}

The orbit corresponding to $(\id+\sum_{a\in \beta}e_{1}\otimes e_{a}^*)w$ goes to the orbit corresponding to $\tau^{-1}\cdot(\id+\sum_{a\in \beta}e_{1}\otimes e_{a}^*)w= \tau^{-1}\cdot(\id+\sum_{a\in \beta}e_{1}\otimes e_{a}^*)\tau \tau^{-1}w= (\id+\sum_{a\in \beta}e_{M}\otimes e_{\tau^{-1}(a)}^*)\tau^{-1}w$.

However, it is different from the representative that we took in Lemma \ref{lem 9.2.1+}, which is of the form $(\id+\sum_{a\in w(\beta)}e_{a}\otimes e_{M}^*)w$.

Consider the $(B'_{M-1},\chi)\ltimes (V\oplus \BG_a,\psi')\ltimes (U_{M,N},\psi)\times (B_N,\chi')$-orbit in $\GL_N$ corresponding to $(w,\beta)$. Denote by $\sigma\subset \{2,3,\cdots,M\}$ the set of minimal elements. In particular, we have that the sequence $w^{-1}\colon \sigma\rightarrow\{1,\ldots,N\}$ is decreasing, and $w^{-1}(\sigma)<w^{-1}(1)$. For convenience, we take the reduced representative $(\id+\sum_{a\in \sigma}e_{1}\otimes e_{a}^*)w$. The invariance condition implies that any equivariant function in $\BC^{(B'_{M-1},\chi)\ltimes (V\oplus \BG_a,\psi')\ltimes U_{M,N},\psi\times (B_N,\chi')}(\GL_N)$ has the same value at the representative $(\id+\sum_{a\in\beta}e_1\otimes e_a^*)w$ and the reduced representative here. Also, the reduced representative goes to the reduced representative $(\id+\sum_{a\in \sigma}e_{M}\otimes e_{\tau^{-1}(a)}^*)\tau^{-1}w$ under \eqref{10.42}.

\begin{lem}
Let $\sigma=\{\sigma_1<\sigma_2<\cdots<\sigma_k\}$, where, if $k=0$, we interpret $w_\sigma$ as the identity and the product of character values as $1$.

Under the third isomorphism, $H^{(3)}_{(w,\beta)}$ goes to \[\prod_{1\leq i\leq k}(w_\sigma\tau^{-1}w\chi')_{\tau^{-1}(\sigma_i)}(-1) H^{(4)}_{(w_\sigma)\tau^{-1}w, ((w_\sigma)\tau^{-1}w)^{-1}(\gamma)}.\] Here, $w_\sigma=(\tau^{-1}(\sigma_1),\tau^{-1}(\sigma_2),\cdots, \tau^{-1}(\sigma_k),M)\in S_M$, and, if we denote $q=w_\sigma\tau^{-1}w$, then we let
\[
\gamma=\{i\in[1,M-1]:\ i\leq j,\ q^{-1}(i)\leq q^{-1}(j) \text{ for some }j\in\tau^{-1}(\sigma)\cup\{M\}\}.
\]
Note that $q^{-1}(\gamma)$ is the index specified by Lemma~\ref{lem 9.2.1+} for the final representative.
\end{lem}
\begin{proof}
Let $p_0=M$, $p_a=\tau^{-1}(\sigma_a)$ and $w_0=\tau^{-1}w$. Then, $p_1<\cdots<p_k<M$ and $w_0^{-1}(p_0)>w_0^{-1}(p_1)>\cdots>w_0^{-1}(p_k)$. Let $E_{ab}=e_a\otimes e_b^*$, and let $D$ be the diagonal matrix with entries $-1$ at $p_0,\ldots,p_{k-1}$ and $1$ elsewhere. Then, by direct row and column operations, we have
\begin{equation}\label{eq:cyclic-orbit-normalization}
\begin{split}
&\left(\id+\sum_{1\leq a<b\leq k}E_{p_ap_b}\right)
\left(\id+\sum_{a=1}^k E_{p_0p_a}\right)
\left(\id-\sum_{a=0}^{k-1}E_{p_{a+1}p_a}\right)\\
&\qquad=\left(\id+\sum_{a=1}^k E_{p_aM}\right)w_\sigma D.
\end{split}
\end{equation}
Note that the left multiplier lies in $U_{M-1}$, and the right multiplier lies in $w_0U_Nw_0^{-1}$, thus, both of them have character value $1$. Moving $D$ to the right of $w_0$, we obtain
\[
\prod_{a=0}^{k-1}\chi'_{w_0^{-1}(p_a)}(-1) =\prod_{a=1}^k(w_\sigma w_0\chi')_{p_a}(-1),
\]
as required. Finally, applying Lemma~\ref{lem 9.2.1+} to $(\id+\sum_{a=1}^k E_{p_aM})w_\sigma w_0$, we obtain $\gamma$ as defined above, and the translation preserves the orbit dimensions and their normalizations.
\end{proof}
\begin{lem}\label{lem:balanced-Fourier-coefficients}
Assume that $M=N$ and $\chi_{1}\neq\chi'_{j}$ for every $j$. Then, the equivalence \eqref{iso 9.14}, which is given by the normalized Fourier transform with the Artin--Schreier kernel associated with $\psi(-\langle v,\xi\rangle)$, the normalized restriction to $\xi=e_1^*$ and the conjugation by $\tau^{-1}$, has the trace formula (we denote by $\mathcal D_N$ the induced map on functions)
\begin{equation}\label{eq:balanced-Fourier-sum}
(\mathcal D_Nf)(g)=\sum_{v\in\BF_Q^N}\psi(-v_1)f(\tau g,v).
\end{equation}
Moreover, the coefficient of $H_z$ in $\mathcal D_N(H_u)$ is
\begin{equation}\label{eq:balanced-Fourier-standard-coefficient}
(-\sqrt Q)^{d_z-d_u} \sum_{v\in\BF_Q^N}\psi(-v_1)T_u(\tau g_z,v),
\end{equation}
where $g_z$ is the chosen representative, $d_u,d_z$ are the orbit dimensions, and $T_u$ has value $1$ at its representative.
\end{lem}
\begin{proof}
By Proposition~\ref{prop 9.2.9}, the equivalence is given by the Fourier transform and the descent from the open orbit of $e_1^*$. The factors $(-1)^NQ^{-N/2}$ from the Fourier transform and $(-1)^NQ^{N/2}$ from the restriction with $\langle-N\rangle$ cancel, which gives \eqref{eq:balanced-Fourier-sum}. The second formula follows from $H_u=(-\sqrt Q)^{-d_u}T_u$.
\end{proof}

\subsubsection{The bar involution and canonical basis in the monodromic case} The formulas in Sections~\ref{sec 10.3.1}, \ref{sec 10.3.2}, and \ref{sec 10.3.3} provide an explicit way to calculate the images of $H_{(w,\beta)}$ and $\underline{H}_{(w,\beta)}$ in $\bK^\fin_{\chi,\chi'}(M|N)$ (resp. $\bK^{\fin,\mxd}_{\chi,\chi'}(M|N)$) under the isomorphism of Grothendieck groups induced by \eqref{iso 9.6}.

Namely, first, we extend the characters as in the proof of Theorem~\ref{Thm 9.2.5} (we can extend the base finite field $\BF_Q$ to $\BF_{Q^n}$ to allow additional characters), using \eqref{9.28} for unequal ranks, Lemma~\ref{lem:balanced-Fourier-coefficients} for equal ranks, and Lemma~\ref{prop 9.2.11} for the added rank $1|1$ factors. Then, we obtain
\begin{equation}
\bK^{\fin,\mxd}_{\chi,\chi'}(M|N) \overset{}{\simeq}K_0(D^{(B_{L},\chi)\times (B_L,\chi')}_\mxd(\GL_L\times \BA^L)).
\end{equation}

Secondly, according to \eqref{iso 9.3}, we permute the characters $\chi$ and $\chi'$ such that they satisfy the conditions in Lemma~\ref{prop 9.2.11}. Now, according to \eqref{eq 9.20}, we have
\begin{equation}
K_0(D^{(B_{L},\chi)\times (B_L,\chi')}_\mxd(\GL_L\times \BA^L)) \overset{}{\simeq}\bigotimes_{i=1}^l K_0(D^{(B_{r_i},\chi^{[i]})\times (B_{r_i},(\chi')^{[i]})}_\mxd(\GL_{r_i}\times \BA^{r_i})).
\end{equation}
Note that the permutation of the right characters in the proof of Lemma~\ref{prop 9.2.11} is given by the right Hecke action, and the equivalences in that proof identify the standard, costandard and IC objects with the corresponding products, with their induced Weil structures.

Thirdly, using Lemma~\ref{lem:balanced-Fourier-coefficients}, \eqref{eq:balanced-inversion}, and then \eqref{9.28} for unequal ranks, we have an explicit formula
\begin{equation}
K_0(D^{(B_{r_i},\chi^{[i]})\times (B_{r_i},(\chi')^{[i]})}_\mxd(\GL_{r_i}\times \BA^{r_i}))\simeq  \bK^{\fin,\mxd}_{b_i^{k_i},b_i^{k_i'}}(k_i|k_i'),
\end{equation}
for any $i$. Also, each coordinate character in the restrictions of $\chi$ and $\chi'$ to $T_{k_i}$ and $T_{k'_i}$ equals $b_i$, as in \eqref{chi 1} and \eqref{chi 2}.

Finally, with \eqref{eq 9.7}, we have
\begin{equation}
\bK^{\fin,\mxd}_{b_i^{k_i},b_i^{k_i'}}(k_i|k_i')\simeq \bK^{\fin,\mxd}_{1,1}(k_i|k_i').
\end{equation}
This step does not involve any combinatorial difficulty. Geometrically, it is given by tensoring with a local system defined on the whole space.

To summarize, we have an explicit formula
\begin{equation}
\bK^{\fin,\mxd}_{\chi,\chi'}(M|N) \overset{}{\simeq}\bigotimes_{i=1}^l \bK^{\fin,\mxd}_{1,1}(k_i|k_i').
\end{equation}
Let $A$ be the matrix of this isomorphism with respect to the standard bases. For a canonical basis element on the right-hand side, let $b$ be the column vector of its coefficients in the standard basis, and let $x=A^{-1}b$. By Corollary~\ref{cor 10.2.14}, the coefficient $a$ of $x$ on the unique maximal orbit $u$ is invertible, and we have
\begin{equation}\label{eq:finite-IC-normalization}
\underline H_u^\mxd=a^{-1}x.
\end{equation}
Let $S$ be the resulting matrix of the canonical basis, and let $S^D$ be the corresponding matrix for the inverse multiplicative characters and $-\psi$, with the same ordering of orbits. Then, the Verdier duality gives
\begin{equation}\label{eq:finite-bar-matrix}
B=S^D(\overline S)^{-1},\qquad c\longmapsto B\overline c,
\end{equation}
where $\overline{[\alpha]}=[\alpha^{-1}]$. Taking specializations, we obtain the graded and graded trace formulas, where, on the latter coefficient ring, we have $\overline{z\bv^n}=\overline z\,\bv^{-n}$, and $\overline z$ is the complex conjugation.

\section{Kazhdan--Lusztig polynomial for the twisted Iwahori Gaiotto category}\label{aff kl}

In this section, we want to study the affine Kazhdan--Lusztig polynomials of the (twisted) Iwahori Gaiotto category in the $\GL(M|N)$ case. Note that Chen and Yi study affine Kazhdan--Lusztig--Vogan polynomials for Iwahori subgroup orbits in loop spaces of symmetric varieties, see \cite[Sections 8--9]{[CY]}.

Roughly speaking, our strategy is as follows. First, we classify the relevant orbits, and identify the category of sheaves on the preimage of the unit $H(\bF)$-orbit in $\Gr_G$ with the finite type category of Section~\ref{sec:finite-type-monodromic}, and show that there is no relevant orbit in the boundary of this preimage (Lemma~\ref{13.2.1} and Proposition~\ref{prop:finite-subcategory-clean}). Then, we compute the action of the affine Yokonuma algebras on the Grothendieck group, and show that every affine standard object is obtained from a standard object of the finite type category by left and right Hecke convolutions whose convolution diagrams are isomorphisms (see Proposition~\ref{prop 9.4.17} and the remark following it). Thus, the bar involution, and hence the canonical basis, is determined by the finite type calculation in Section~\ref{sec:finite-type-monodromic} and the action formulas. The case $M=N$ is treated in Section~\ref{sec:balanced-affine-algorithm}.

\subsection{Classification of relevant orbits} Let us denote by $\Fl_M$ (resp. $\Fl_N$) the affine flags associated with $\GL_M$ (resp. $\GL_N$). We denote by $W_N^{\ext}=S_N\ltimes\BZ^N$ the extended affine Weyl group, and by $W_N^{\aff}=S_N\ltimes\{\lambda\in\BZ^N:\sum_i\lambda_i=0\}$ its affine Weyl subgroup. Similar to \cite{[FGT]}, we regard $W^\ext_N$ as the set of permutations of $\BZ$ such that $w(i + N ) = w(i) + N$, i.e., $N$-periodic permutations. Also, we use the following notation in loc. cit.: choose a basis $\{e_1,\cdots,e_N\}$ of an $N$-dimensional vector space and set $e_{Nk+j}:=t^{-k}e_j$. For a subset $A\subset\BZ$ bounded above, we let $e_A:=\sum_{i\in A}e_i\in\bF^N$. We denote $[a,b]_N=N\BZ+\{a,a+1,\ldots,b\}$ and $[a]_N=a+N\BZ$. For the completeness of the introduction of the twisted Iwahori Gaiotto category, we also include the classification of irreducible objects in the case of $\GL(M|N)$. The orthosymplectic and $\on{F}(4)$ cases are considered in Appendix~\ref{sec:orthogonal-exceptional-Iwahori}.

For $M=N$, we consider the diagonal $I_N$-orbits in $\Fl_N\times\bF^N$ with non-zero vector component. They are indexed by the pairs $(w,\beta)$ of \cite[\S\S 2.1--2.3]{[FGT]}. Here, $w\in W_N^{\ext}$, and $\beta$ is bounded above, contains every sufficiently negative integer, and satisfies
\[
j\in\beta,\quad i\leq j,\quad w(i)\leq w(j) \quad\Longrightarrow\quad i\in\beta.
\]
The representative is $(w,e_{w(\beta)})$. We treat the case $M=N$ in Section~\ref{sec:balanced-affine-algorithm}, see also \cite[\S\S 2.5, 2.9]{[FGT]}, and in the calculations below, we assume that $M<N$.

For $M<N$, it is enough to index the $I_M\ltimes U_{M,N}(\bF)$-orbits in $\Fl_N=\GL_N(\bF)/I_N$, with the conventions of Section~\ref{sec:type-a-conventions}.

When $M=0$, the auxiliary group denoted below by $I_M\ltimes U_{M-1,N}(\bF)$ is understood to be $U_N(\bF)$, and all left Hecke actions are omitted.

In the Gaiotto case, it is known that the stabilizer subgroup of any orbit is connected; thus, classifying all irreducible objects is equivalent to classifying all relevant orbits.

The following lemma is straightforward.
\begin{lem}\label{lem 3.2.2}
If $1\leq M<N$, then the set of $I_M\ltimes U_{M-1,N}(\bF)$-orbits in $\Fl_N$ is indexed by the extended Weyl group $W^\ext_N$ of $\GL_N$.
\end{lem}

\subsubsection{Orbit indices} For a periodic matrix group $K$, we write $i\leq_K j$ if $e_i\otimes e_j^*\in\Lie(K)$, and $i>_K j$ if $j\leq_K i$ and $i\neq j$. For $K=I_N$, this is the usual order on integers.

To classify the $I_M\ltimes U_{M,N}(\bF)$-orbits in $\Fl_N$, let us introduce the following partial order on $[1,M]_N:=N\BZ+\{1,2,\cdots, M\}$, the set of integers which belong to $\{1,2,\cdots, M\}$ modulo $N$.
\begin{defn}
Regard $I_M$ as a subgroup of $I_N$ by the upper-left embedding, and fix an element $w\in W^\ext_{N}$. For $i,j\in [1,M]_N$, we say $i\leq_{I_M\cap wI_Nw^{-1}}j$, if:
\begin{enumerate}
\item $i\leq j$, and
\item $w^{-1}(i)\leq w^{-1}(j)$.
\end{enumerate}
\end{defn}

\begin{rem}
As the notation suggests, $i\leq_{I_M\cap wI_Nw^{-1}} j$ if and only if $e_i\otimes e_j^*\in \Lie(I_M)\cap w\Lie(I_N)w^{-1}$. Here, we regard elements in $\GL_N(\bF)$ as periodic $\BZ\times \BZ$-matrices.
\end{rem}

\begin{prop}\label{prop 3.2.2}
The $I_M\ltimes U_{M,N}(\bF)$-orbits in $\Fl_N$ are indexed by pairs $(w,\beta)$, where $w\in W^{\ext}_N$ and $\beta\subset [1,M]_N$, such that
\begin{enumerate}
\item $\beta$ is bounded above;
\item if $i\leq_{I_M\cap wI_Nw^{-1}} j\in [1,M]_N$, and $j\in \beta$, then $i\in \beta$;
\item if $i\in [1,M]_N$, and $w^{-1}(i)<w^{-1}(M+1)$, then $i$ belongs to $\beta$,
\end{enumerate}
and the representative is given by $g_{\beta}w$, where \[g_{\beta}:=\Id+\sum_{i\in\beta}e_{i}\otimes e^*_{M+1}\in \GL_N(\bF).\]
\end{prop}

\begin{proof}
For $M=0$, this is the usual semi-infinite Iwahori decomposition with $\beta=\varnothing$. In the following, we assume that $M\geq1$.

According to Lemma \ref{lem 3.2.2}, any $I_M\ltimes U_{M,N}(\bF)$-orbit in $\Fl_N$ admits a representative $g_{\beta}w$, and two such representatives $g_{\beta}w$ and $g_{\beta'}w'$ belong to the same orbit only if $w=w'$. One needs to figure out the condition on $\beta, \beta'$ such that the orbit corresponding to $g_{\beta}w$ is the same as the orbit corresponding to $g_{\beta'}w$.

Assume that two such representatives belong to the same orbit. That is to say, there exists $h\in I_M\ltimes U_{M,N}(\bF)$, such that
\begin{equation}
(g_{\beta'})^{-1}h g_\beta\in wI_Nw^{-1}.
\end{equation}

Note that $g_{\beta}=\Id+\sum_{i\in\beta}e_{i}\otimes e^*_{M+1}$. The above is just equivalent to the existence of $h'\in I_M\cap wI_N w^{-1}$, such that
\begin{equation}
e_\beta\otimes e_{M+1}^*-h'\cdot e_{\beta'}\otimes e_{M+1}^*\in \on{Lie}(w I_N w^{-1}).
\end{equation}

It is equivalent to the existence of $h''\in I_M\cap wI_N w^{-1}$, such that
\begin{equation}
e_{\tilde{\beta}}\otimes e_{M+1}^*=h''\cdot e_{\tilde{\beta}'}\otimes e_{M+1}^*,
\end{equation}
or, just, $ e_{\tilde{\beta}}=h''\cdot e_{\tilde{\beta}'}$, where $\tilde\beta=\beta\cup\{i\in[1,M]_N:w^{-1}(i)<w^{-1}(M+1)\}$, and similarly for $\tilde\beta'$. From this description, the classification follows immediately.
\end{proof}

The orbit corresponding to $(w,\beta)$ here corresponds to the orbit indexed by $(w,w^{-1}(\beta))$ in \cite{[FGT]}. We adopt the convention in loc. cit., which is described as follows.

\begin{proprestate}{prop 3.2.2}\label{prop 3.2.7}
The $I_M\ltimes U_{M,N}(\bF)$-orbits in $\Fl_N$ are indexed by $(w,\beta)$, such that $w\in W^\ext_N$, $\beta\subset w^{-1}([1,M]_N)$, and
\begin{enumerate}
\item $\beta$ is bounded above;
\item if $i,j\in w^{-1}([1,M]_N)$, $i\leq j$, $w(i)\leq w(j)$ (we denote the order by $i\leq_{w^{-1}I_M w\cap I_N} j$), and $j\in \beta$, then $i\in \beta$;
\item if $i\in w^{-1}([1,M]_N)$, and $i<w^{-1}(M+1)$, then $i\in \beta$.
\end{enumerate}
The representative is given by $g_{w(\beta)}w$.
\end{proprestate}

For an index $(w,\beta)$, we denote by $\sigma$ the set of maximal elements of $\beta$ with respect to the order $\leq_{w^{-1}I_Mw\cap I_N}$ which are bigger than $w^{-1}(M+1)$, and by $T_{w,\sigma}\subset T_M$ the coordinate subtorus on the residues which are not represented in $w(\sigma)$. For $w=t^\lambda w^{\fin}$, we denote $n_i=\lambda_i$, so that $w^{-1}(e_i)=t^{-n_i}e_{(w^{\fin})^{-1}(i)}$.

Figures~\ref{fig:affine-orbit-classification} and \ref{fig:affine-orbit-classification-dispersed} below give the periodic version of Figure~\ref{fig:finite-orbit-classification}. Note that, among rows $[1,M]_N$, every dot above and to the left of a blue dot must be blue, and every dot to the left of the column $w^{-1}(M+1)$ must be blue.

\begin{figure}[htbp]
\centering
\[
M=4,\qquad N=6,\qquad w(i+6)=w(i)+6.
\]
\resizebox{\linewidth}{!}{\begin{tikzpicture}[x=1cm,y=-1cm,font=\small]
\definecolor{orbitblue}{RGB}{32,94,176}
\definecolor{orbitred}{RGB}{184,40,47}

\begin{scope}[xshift=0cm]
\node[font=\normalsize] at (2.880,-2.200) {The pair $(w,\beta)$}; \node at (2.880,-1.750) {$(w(1),\ldots,w(6))=(4,5,2,1,6,3)$}; \node at (2.880,-1.350) {$\beta=\{i\in w^{-1}([1,4]_6):i<2\}\cup\{3,4\}$}; \node at (2.880,-0.950) {$\sigma=\{3,4\}$}; \draw[black!13,line width=.25pt] (0,0)--(0,5.76); \draw[black!13,line width=.25pt] (0,0)--(5.76,0); \draw[black!13,line width=.25pt] (0.32,0)--(0.32,5.76); \draw[black!13,line width=.25pt] (0,0.32)--(5.76,0.32); \draw[black!13,line width=.25pt] (0.64,0)--(0.64,5.76); \draw[black!13,line width=.25pt] (0,0.64)--(5.76,0.64); \draw[black!13,line width=.25pt] (0.96,0)--(0.96,5.76); \draw[black!13,line width=.25pt] (0,0.96)--(5.76,0.96); \draw[black!13,line width=.25pt] (1.28,0)--(1.28,5.76); \draw[black!13,line width=.25pt] (0,1.28)--(5.76,1.28); \draw[black!13,line width=.25pt] (1.6,0)--(1.6,5.76); \draw[black!13,line width=.25pt] (0,1.6)--(5.76,1.6); \draw[black!13,line width=.25pt] (1.92,0)--(1.92,5.76); \draw[black!13,line width=.25pt] (0,1.92)--(5.76,1.92); \draw[black!13,line width=.25pt] (2.24,0)--(2.24,5.76); \draw[black!13,line width=.25pt] (0,2.24)--(5.76,2.24); \draw[black!13,line width=.25pt] (2.56,0)--(2.56,5.76); \draw[black!13,line width=.25pt] (0,2.56)--(5.76,2.56); \draw[black!13,line width=.25pt] (2.88,0)--(2.88,5.76); \draw[black!13,line width=.25pt] (0,2.88)--(5.76,2.88); \draw[black!13,line width=.25pt] (3.2,0)--(3.2,5.76); \draw[black!13,line width=.25pt] (0,3.2)--(5.76,3.2); \draw[black!13,line width=.25pt] (3.52,0)--(3.52,5.76); \draw[black!13,line width=.25pt] (0,3.52)--(5.76,3.52); \draw[black!13,line width=.25pt] (3.84,0)--(3.84,5.76); \draw[black!13,line width=.25pt] (0,3.84)--(5.76,3.84); \draw[black!13,line width=.25pt] (4.16,0)--(4.16,5.76); \draw[black!13,line width=.25pt] (0,4.16)--(5.76,4.16); \draw[black!13,line width=.25pt] (4.48,0)--(4.48,5.76); \draw[black!13,line width=.25pt] (0,4.48)--(5.76,4.48); \draw[black!13,line width=.25pt] (4.8,0)--(4.8,5.76); \draw[black!13,line width=.25pt] (0,4.8)--(5.76,4.8); \draw[black!13,line width=.25pt] (5.12,0)--(5.12,5.76); \draw[black!13,line width=.25pt] (0,5.12)--(5.76,5.12); \draw[black!13,line width=.25pt] (5.44,0)--(5.44,5.76); \draw[black!13,line width=.25pt] (0,5.44)--(5.76,5.44); \draw[black!13,line width=.25pt] (5.76,0)--(5.76,5.76); \draw[black!13,line width=.25pt] (0,5.76)--(5.76,5.76); \node[font=\scriptsize,anchor=south] at (0.16,-0.05) {$-5$}; \node[font=\scriptsize,anchor=east] at (-0.05,0.16) {$-5$}; \node[font=\scriptsize,anchor=south] at (0.48,-0.05) {$-4$}; \node[font=\scriptsize,anchor=east] at (-0.05,0.48) {$-4$}; \node[font=\scriptsize,anchor=south] at (0.8,-0.05) {$-3$}; \node[font=\scriptsize,anchor=east] at (-0.05,0.8) {$-3$}; \node[font=\scriptsize,anchor=south] at (1.12,-0.05) {$-2$}; \node[font=\scriptsize,anchor=east] at (-0.05,1.12) {$-2$}; \node[font=\scriptsize,anchor=south] at (1.44,-0.05) {$-1$}; \node[font=\scriptsize,anchor=east] at (-0.05,1.44) {$-1$}; \node[font=\scriptsize,anchor=south] at (1.76,-0.05) {$0$}; \node[font=\scriptsize,anchor=east] at (-0.05,1.76) {$0$}; \node[font=\scriptsize,anchor=south] at (2.08,-0.05) {$1$}; \node[font=\scriptsize,anchor=east] at (-0.05,2.08) {$1$}; \node[font=\scriptsize,anchor=south] at (2.4,-0.05) {$2$}; \node[font=\scriptsize,anchor=east] at (-0.05,2.4) {$2$}; \node[font=\scriptsize,anchor=south] at (2.72,-0.05) {$3$}; \node[font=\scriptsize,anchor=east] at (-0.05,2.72) {$3$}; \node[font=\scriptsize,anchor=south] at (3.04,-0.05) {$4$}; \node[font=\scriptsize,anchor=east] at (-0.05,3.04) {$4$}; \node[font=\scriptsize,anchor=south] at (3.36,-0.05) {$5$}; \node[font=\scriptsize,anchor=east] at (-0.05,3.36) {$5$}; \node[font=\scriptsize,anchor=south] at (3.68,-0.05) {$6$}; \node[font=\scriptsize,anchor=east] at (-0.05,3.68) {$6$}; \node[font=\scriptsize,anchor=south] at (4,-0.05) {$7$}; \node[font=\scriptsize,anchor=east] at (-0.05,4) {$7$}; \node[font=\scriptsize,anchor=south] at (4.32,-0.05) {$8$}; \node[font=\scriptsize,anchor=east] at (-0.05,4.32) {$8$}; \node[font=\scriptsize,anchor=south] at (4.64,-0.05) {$9$}; \node[font=\scriptsize,anchor=east] at (-0.05,4.64) {$9$}; \node[font=\scriptsize,anchor=south] at (4.96,-0.05) {$10$}; \node[font=\scriptsize,anchor=east] at (-0.05,4.96) {$10$}; \node[font=\scriptsize,anchor=south] at (5.28,-0.05) {$11$}; \node[font=\scriptsize,anchor=east] at (-0.05,5.28) {$11$}; \node[font=\scriptsize,anchor=south] at (5.6,-0.05) {$12$}; \node[font=\scriptsize,anchor=east] at (-0.05,5.6) {$12$}; \draw[black!40,dashed,line width=.6pt] (0,0) rectangle (1.92,1.92); \draw[black!70,line width=.8pt] (1.92,1.92) rectangle (3.84,3.84); \draw[black!40,dashed,line width=.6pt] (3.84,3.84) rectangle (5.76,5.76); \node[black!45] at (-0.36,-0.36) {$\ddots$}; \node[black!45] at (5.98,5.96) {$\ddots$}; \draw[->,black!55] (0,-0.490)--(5.760,-0.490) node[right] {$i$}; \draw[->,black!55] (-0.580,0)--(-0.580,5.760) node[below] {$j$}; \draw[orbitred!65,densely dashed] (2.400,0)--(2.400,5.760); \draw[orbitred!65,densely dashed] (0,3.360)--(2.400,3.360); \fill[orbitblue] (0.16,1.12) circle (2.8pt); \fill[black!85] (0.48,1.44) circle (2.8pt); \fill[orbitblue] (0.8,0.48) circle (2.8pt); \fill[orbitblue] (1.12,0.16) circle (2.8pt); \fill[black!85] (1.44,1.76) circle (2.8pt); \fill[orbitblue] (1.76,0.8) circle (2.8pt); \fill[orbitblue] (2.08,3.04) circle (2.8pt); \fill[orbitred] (2.4,3.36) circle (2.8pt); \fill[orbitblue] (2.72,2.4) circle (2.8pt); \fill[orbitblue] (3.04,2.08) circle (2.8pt); \fill[black!85] (3.36,3.68) circle (2.8pt); \fill[black!85] (3.68,2.72) circle (2.8pt); \fill[black!85] (4,4.96) circle (2.8pt); \fill[black!85] (4.32,5.28) circle (2.8pt); \fill[black!85] (4.64,4.32) circle (2.8pt); \fill[black!85] (4.96,4) circle (2.8pt); \fill[black!85] (5.28,5.6) circle (2.8pt); \fill[black!85] (5.6,4.64) circle (2.8pt); \draw[orbitblue,line width=.7pt] (2.72,2.4) circle (4.5pt); \draw[orbitblue,line width=.7pt] (3.04,2.08) circle (4.5pt); \node[orbitred,font=\scriptsize,anchor=west,fill=white,inner sep=1pt] at (2.58,3.36) {$M+1=5$}; \node[orbitred,font=\scriptsize] at (2.4,6.23) {$w^{-1}(M+1)=2$};
\end{scope}
\begin{scope}[xshift=8.15cm]
\node[font=\normalsize] at (2.880,-2.200) {The representative $g_{w(\beta)}w$}; \node at (2.880,-1.750) {$(w(1),\ldots,w(6))=(4,5,2,1,6,3)$}; \node at (2.880,-1.350) {$\beta=\{i\in w^{-1}([1,4]_6):i<2\}\cup\{3,4\}$}; \node at (2.880,-0.950) {$\sigma=\{3,4\}$}; \draw[black!13,line width=.25pt] (0,0)--(0,5.76); \draw[black!13,line width=.25pt] (0,0)--(5.76,0); \draw[black!13,line width=.25pt] (0.32,0)--(0.32,5.76); \draw[black!13,line width=.25pt] (0,0.32)--(5.76,0.32); \draw[black!13,line width=.25pt] (0.64,0)--(0.64,5.76); \draw[black!13,line width=.25pt] (0,0.64)--(5.76,0.64); \draw[black!13,line width=.25pt] (0.96,0)--(0.96,5.76); \draw[black!13,line width=.25pt] (0,0.96)--(5.76,0.96); \draw[black!13,line width=.25pt] (1.28,0)--(1.28,5.76); \draw[black!13,line width=.25pt] (0,1.28)--(5.76,1.28); \draw[black!13,line width=.25pt] (1.6,0)--(1.6,5.76); \draw[black!13,line width=.25pt] (0,1.6)--(5.76,1.6); \draw[black!13,line width=.25pt] (1.92,0)--(1.92,5.76); \draw[black!13,line width=.25pt] (0,1.92)--(5.76,1.92); \draw[black!13,line width=.25pt] (2.24,0)--(2.24,5.76); \draw[black!13,line width=.25pt] (0,2.24)--(5.76,2.24); \draw[black!13,line width=.25pt] (2.56,0)--(2.56,5.76); \draw[black!13,line width=.25pt] (0,2.56)--(5.76,2.56); \draw[black!13,line width=.25pt] (2.88,0)--(2.88,5.76); \draw[black!13,line width=.25pt] (0,2.88)--(5.76,2.88); \draw[black!13,line width=.25pt] (3.2,0)--(3.2,5.76); \draw[black!13,line width=.25pt] (0,3.2)--(5.76,3.2); \draw[black!13,line width=.25pt] (3.52,0)--(3.52,5.76); \draw[black!13,line width=.25pt] (0,3.52)--(5.76,3.52); \draw[black!13,line width=.25pt] (3.84,0)--(3.84,5.76); \draw[black!13,line width=.25pt] (0,3.84)--(5.76,3.84); \draw[black!13,line width=.25pt] (4.16,0)--(4.16,5.76); \draw[black!13,line width=.25pt] (0,4.16)--(5.76,4.16); \draw[black!13,line width=.25pt] (4.48,0)--(4.48,5.76); \draw[black!13,line width=.25pt] (0,4.48)--(5.76,4.48); \draw[black!13,line width=.25pt] (4.8,0)--(4.8,5.76); \draw[black!13,line width=.25pt] (0,4.8)--(5.76,4.8); \draw[black!13,line width=.25pt] (5.12,0)--(5.12,5.76); \draw[black!13,line width=.25pt] (0,5.12)--(5.76,5.12); \draw[black!13,line width=.25pt] (5.44,0)--(5.44,5.76); \draw[black!13,line width=.25pt] (0,5.44)--(5.76,5.44); \draw[black!13,line width=.25pt] (5.76,0)--(5.76,5.76); \draw[black!13,line width=.25pt] (0,5.76)--(5.76,5.76); \node[font=\scriptsize,anchor=south] at (0.16,-0.05) {$-5$}; \node[font=\scriptsize,anchor=east] at (-0.05,0.16) {$-5$}; \node[font=\scriptsize,anchor=south] at (0.48,-0.05) {$-4$}; \node[font=\scriptsize,anchor=east] at (-0.05,0.48) {$-4$}; \node[font=\scriptsize,anchor=south] at (0.8,-0.05) {$-3$}; \node[font=\scriptsize,anchor=east] at (-0.05,0.8) {$-3$}; \node[font=\scriptsize,anchor=south] at (1.12,-0.05) {$-2$}; \node[font=\scriptsize,anchor=east] at (-0.05,1.12) {$-2$}; \node[font=\scriptsize,anchor=south] at (1.44,-0.05) {$-1$}; \node[font=\scriptsize,anchor=east] at (-0.05,1.44) {$-1$}; \node[font=\scriptsize,anchor=south] at (1.76,-0.05) {$0$}; \node[font=\scriptsize,anchor=east] at (-0.05,1.76) {$0$}; \node[font=\scriptsize,anchor=south] at (2.08,-0.05) {$1$}; \node[font=\scriptsize,anchor=east] at (-0.05,2.08) {$1$}; \node[font=\scriptsize,anchor=south] at (2.4,-0.05) {$2$}; \node[font=\scriptsize,anchor=east] at (-0.05,2.4) {$2$}; \node[font=\scriptsize,anchor=south] at (2.72,-0.05) {$3$}; \node[font=\scriptsize,anchor=east] at (-0.05,2.72) {$3$}; \node[font=\scriptsize,anchor=south] at (3.04,-0.05) {$4$}; \node[font=\scriptsize,anchor=east] at (-0.05,3.04) {$4$}; \node[font=\scriptsize,anchor=south] at (3.36,-0.05) {$5$}; \node[font=\scriptsize,anchor=east] at (-0.05,3.36) {$5$}; \node[font=\scriptsize,anchor=south] at (3.68,-0.05) {$6$}; \node[font=\scriptsize,anchor=east] at (-0.05,3.68) {$6$}; \node[font=\scriptsize,anchor=south] at (4,-0.05) {$7$}; \node[font=\scriptsize,anchor=east] at (-0.05,4) {$7$}; \node[font=\scriptsize,anchor=south] at (4.32,-0.05) {$8$}; \node[font=\scriptsize,anchor=east] at (-0.05,4.32) {$8$}; \node[font=\scriptsize,anchor=south] at (4.64,-0.05) {$9$}; \node[font=\scriptsize,anchor=east] at (-0.05,4.64) {$9$}; \node[font=\scriptsize,anchor=south] at (4.96,-0.05) {$10$}; \node[font=\scriptsize,anchor=east] at (-0.05,4.96) {$10$}; \node[font=\scriptsize,anchor=south] at (5.28,-0.05) {$11$}; \node[font=\scriptsize,anchor=east] at (-0.05,5.28) {$11$}; \node[font=\scriptsize,anchor=south] at (5.6,-0.05) {$12$}; \node[font=\scriptsize,anchor=east] at (-0.05,5.6) {$12$}; \draw[black!40,dashed,line width=.6pt] (0,0) rectangle (1.92,1.92); \draw[black!70,line width=.8pt] (1.92,1.92) rectangle (3.84,3.84); \draw[black!40,dashed,line width=.6pt] (3.84,3.84) rectangle (5.76,5.76); \node[black!45] at (-0.36,-0.36) {$\ddots$}; \node[black!45] at (5.98,5.96) {$\ddots$}; \draw[->,black!55] (0,-0.490)--(5.760,-0.490) node[right] {$i$}; \draw[->,black!55] (-0.580,0)--(-0.580,5.760) node[below] {$j$}; \draw[orbitred!65,densely dashed] (2.400,0)--(2.400,5.760); \draw[orbitred!65,densely dashed] (0,3.360)--(2.400,3.360); \fill[black!85] (0.16,1.12) circle (2.8pt); \fill[black!85] (0.48,1.44) circle (2.8pt); \fill[black!85] (0.8,0.48) circle (2.8pt); \fill[black!85] (1.12,0.16) circle (2.8pt); \fill[black!85] (1.44,1.76) circle (2.8pt); \fill[black!85] (1.76,0.8) circle (2.8pt); \fill[black!85] (2.08,3.04) circle (2.8pt); \fill[orbitred] (2.4,3.36) circle (2.8pt); \fill[black!85] (2.72,2.4) circle (2.8pt); \fill[black!85] (3.04,2.08) circle (2.8pt); \fill[black!85] (3.36,3.68) circle (2.8pt); \fill[black!85] (3.68,2.72) circle (2.8pt); \fill[black!85] (4,4.96) circle (2.8pt); \fill[black!85] (4.32,5.28) circle (2.8pt); \fill[black!85] (4.64,4.32) circle (2.8pt); \fill[black!85] (4.96,4) circle (2.8pt); \fill[black!85] (5.28,5.6) circle (2.8pt); \fill[black!85] (5.6,4.64) circle (2.8pt); \fill[orbitblue] (0.48,0.16) circle (2.8pt); \fill[orbitblue] (2.4,0.16) circle (2.8pt); \fill[orbitblue] (4.32,0.16) circle (2.8pt); \fill[orbitblue] (0.48,0.48) circle (2.8pt); \fill[orbitblue] (2.4,0.48) circle (2.8pt); \fill[orbitblue] (4.32,0.48) circle (2.8pt); \fill[orbitblue] (2.4,0.8) circle (2.8pt); \fill[orbitblue] (4.32,0.8) circle (2.8pt); \fill[orbitblue] (0.48,1.12) circle (2.8pt); \fill[orbitblue] (2.4,1.12) circle (2.8pt); \fill[orbitblue] (4.32,1.12) circle (2.8pt); \fill[orbitblue] (2.4,2.08) circle (2.8pt); \fill[orbitblue] (4.32,2.08) circle (2.8pt); \fill[orbitblue] (2.4,2.4) circle (2.8pt); \fill[orbitblue] (4.32,2.4) circle (2.8pt); \fill[orbitblue] (4.32,2.72) circle (2.8pt); \fill[orbitblue] (2.4,3.04) circle (2.8pt); \fill[orbitblue] (4.32,3.04) circle (2.8pt); \fill[orbitblue] (4.32,4) circle (2.8pt); \fill[orbitblue] (4.32,4.32) circle (2.8pt); \fill[orbitblue] (4.32,4.96) circle (2.8pt);
\end{scope}
\node[anchor=north west,align=left,text width=6.7cm,font=\footnotesize] at (0,6.91) { \textcolor{orbitblue}{Blue}: $i\in\beta$. Circled: $i\in\sigma=\{3,4\}$.\\[3pt]
\textcolor{orbitred}{Red}: row $5$, column $2$. }; \node[anchor=north west,align=left,text width=6.5cm,font=\footnotesize] at (8.15,6.91){\textcolor{orbitblue}{Blue}: the additional entries at\\[2pt]row $w(i)+6k$, column $2+6k$,\\[2pt]$i\in\beta$, $k\in\mathbb Z$.};
\end{tikzpicture}
} \caption{An example of the orbit classification in Proposition~\ref{prop 3.2.7}. This is a very special case: $w$ is obtained by repeating the finite permutation matrix in Figure~\ref{fig:finite-orbit-classification} along the diagonal. In the left diagram, all dots are entries of $w$. The blue dots correspond to $\beta$, and the circled dots correspond to $\sigma$. The right diagram is the corresponding representative as a periodic matrix, where every marked entry is $1$ and every other entry is $0$. The matrices $w$ and $g_{w(\beta)}w$ are $6$-periodic; $\beta$ is bounded above. The red dot has row $M+1$ and column $w^{-1}(M+1)$. The solid square marks rows and columns $1,\ldots,N$.}
\label{fig:affine-orbit-classification}
\end{figure}

\begin{figure}[htbp]
\centering
\[
M=4,\qquad N=6,\qquad w(i+6)=w(i)+6.
\]
\resizebox{\linewidth}{!}{\begin{tikzpicture}[x=1cm,y=-1cm,font=\small]
\definecolor{orbitblue}{RGB}{32,94,176}
\definecolor{orbitred}{RGB}{184,40,47}

\begin{scope}[xshift=0cm]
\node[font=\normalsize] at (2.880,-2.200) {The pair $(w,\beta)$}; \node at (2.880,-1.750) {$(w(1),\ldots,w(6))=(4,-1,2,1,12,3)$}; \node at (2.880,-1.350) {$\beta=\{i\in w^{-1}([1,4]_6):i<8\}\cup\{9,10\}$}; \node at (2.880,-0.950) {$\sigma=\{9,10\}$}; \draw[black!13,line width=.25pt] (0,0)--(0,5.76); \draw[black!13,line width=.25pt] (0,0)--(5.76,0); \draw[black!13,line width=.25pt] (0.32,0)--(0.32,5.76); \draw[black!13,line width=.25pt] (0,0.32)--(5.76,0.32); \draw[black!13,line width=.25pt] (0.64,0)--(0.64,5.76); \draw[black!13,line width=.25pt] (0,0.64)--(5.76,0.64); \draw[black!13,line width=.25pt] (0.96,0)--(0.96,5.76); \draw[black!13,line width=.25pt] (0,0.96)--(5.76,0.96); \draw[black!13,line width=.25pt] (1.28,0)--(1.28,5.76); \draw[black!13,line width=.25pt] (0,1.28)--(5.76,1.28); \draw[black!13,line width=.25pt] (1.6,0)--(1.6,5.76); \draw[black!13,line width=.25pt] (0,1.6)--(5.76,1.6); \draw[black!13,line width=.25pt] (1.92,0)--(1.92,5.76); \draw[black!13,line width=.25pt] (0,1.92)--(5.76,1.92); \draw[black!13,line width=.25pt] (2.24,0)--(2.24,5.76); \draw[black!13,line width=.25pt] (0,2.24)--(5.76,2.24); \draw[black!13,line width=.25pt] (2.56,0)--(2.56,5.76); \draw[black!13,line width=.25pt] (0,2.56)--(5.76,2.56); \draw[black!13,line width=.25pt] (2.88,0)--(2.88,5.76); \draw[black!13,line width=.25pt] (0,2.88)--(5.76,2.88); \draw[black!13,line width=.25pt] (3.2,0)--(3.2,5.76); \draw[black!13,line width=.25pt] (0,3.2)--(5.76,3.2); \draw[black!13,line width=.25pt] (3.52,0)--(3.52,5.76); \draw[black!13,line width=.25pt] (0,3.52)--(5.76,3.52); \draw[black!13,line width=.25pt] (3.84,0)--(3.84,5.76); \draw[black!13,line width=.25pt] (0,3.84)--(5.76,3.84); \draw[black!13,line width=.25pt] (4.16,0)--(4.16,5.76); \draw[black!13,line width=.25pt] (0,4.16)--(5.76,4.16); \draw[black!13,line width=.25pt] (4.48,0)--(4.48,5.76); \draw[black!13,line width=.25pt] (0,4.48)--(5.76,4.48); \draw[black!13,line width=.25pt] (4.8,0)--(4.8,5.76); \draw[black!13,line width=.25pt] (0,4.8)--(5.76,4.8); \draw[black!13,line width=.25pt] (5.12,0)--(5.12,5.76); \draw[black!13,line width=.25pt] (0,5.12)--(5.76,5.12); \draw[black!13,line width=.25pt] (5.44,0)--(5.44,5.76); \draw[black!13,line width=.25pt] (0,5.44)--(5.76,5.44); \draw[black!13,line width=.25pt] (5.76,0)--(5.76,5.76); \draw[black!13,line width=.25pt] (0,5.76)--(5.76,5.76); \node[font=\scriptsize,anchor=south] at (0.16,-0.05) {$-5$}; \node[font=\scriptsize,anchor=east] at (-0.05,0.16) {$-5$}; \node[font=\scriptsize,anchor=south] at (0.48,-0.05) {$-4$}; \node[font=\scriptsize,anchor=east] at (-0.05,0.48) {$-4$}; \node[font=\scriptsize,anchor=south] at (0.8,-0.05) {$-3$}; \node[font=\scriptsize,anchor=east] at (-0.05,0.8) {$-3$}; \node[font=\scriptsize,anchor=south] at (1.12,-0.05) {$-2$}; \node[font=\scriptsize,anchor=east] at (-0.05,1.12) {$-2$}; \node[font=\scriptsize,anchor=south] at (1.44,-0.05) {$-1$}; \node[font=\scriptsize,anchor=east] at (-0.05,1.44) {$-1$}; \node[font=\scriptsize,anchor=south] at (1.76,-0.05) {$0$}; \node[font=\scriptsize,anchor=east] at (-0.05,1.76) {$0$}; \node[font=\scriptsize,anchor=south] at (2.08,-0.05) {$1$}; \node[font=\scriptsize,anchor=east] at (-0.05,2.08) {$1$}; \node[font=\scriptsize,anchor=south] at (2.4,-0.05) {$2$}; \node[font=\scriptsize,anchor=east] at (-0.05,2.4) {$2$}; \node[font=\scriptsize,anchor=south] at (2.72,-0.05) {$3$}; \node[font=\scriptsize,anchor=east] at (-0.05,2.72) {$3$}; \node[font=\scriptsize,anchor=south] at (3.04,-0.05) {$4$}; \node[font=\scriptsize,anchor=east] at (-0.05,3.04) {$4$}; \node[font=\scriptsize,anchor=south] at (3.36,-0.05) {$5$}; \node[font=\scriptsize,anchor=east] at (-0.05,3.36) {$5$}; \node[font=\scriptsize,anchor=south] at (3.68,-0.05) {$6$}; \node[font=\scriptsize,anchor=east] at (-0.05,3.68) {$6$}; \node[font=\scriptsize,anchor=south] at (4,-0.05) {$7$}; \node[font=\scriptsize,anchor=east] at (-0.05,4) {$7$}; \node[font=\scriptsize,anchor=south] at (4.32,-0.05) {$8$}; \node[font=\scriptsize,anchor=east] at (-0.05,4.32) {$8$}; \node[font=\scriptsize,anchor=south] at (4.64,-0.05) {$9$}; \node[font=\scriptsize,anchor=east] at (-0.05,4.64) {$9$}; \node[font=\scriptsize,anchor=south] at (4.96,-0.05) {$10$}; \node[font=\scriptsize,anchor=east] at (-0.05,4.96) {$10$}; \node[font=\scriptsize,anchor=south] at (5.28,-0.05) {$11$}; \node[font=\scriptsize,anchor=east] at (-0.05,5.28) {$11$}; \node[font=\scriptsize,anchor=south] at (5.6,-0.05) {$12$}; \node[font=\scriptsize,anchor=east] at (-0.05,5.6) {$12$}; \draw[black!40,dashed,line width=.6pt] (0,0) rectangle (1.92,1.92); \draw[black!70,line width=.8pt] (1.92,1.92) rectangle (3.84,3.84); \draw[black!40,dashed,line width=.6pt] (3.84,3.84) rectangle (5.76,5.76); \node[black!45] at (-0.36,-0.36) {$\ddots$}; \node[black!45] at (5.98,5.96) {$\ddots$}; \draw[->,black!55] (0,-0.490)--(5.760,-0.490) node[right] {$i$}; \draw[->,black!55] (-0.580,0)--(-0.580,5.760) node[below] {$j$}; \draw[orbitred!65,densely dashed] (4.320,0)--(4.320,5.760); \draw[orbitred!65,densely dashed] (0,3.360)--(4.320,3.360); \fill[orbitblue] (0.16,1.12) circle (2.8pt); \fill[orbitblue] (0.8,0.48) circle (2.8pt); \fill[orbitblue] (1.12,0.16) circle (2.8pt); \fill[black!85] (1.44,3.68) circle (2.8pt); \fill[orbitblue] (1.76,0.8) circle (2.8pt); \fill[orbitblue] (2.08,3.04) circle (2.8pt); \fill[black!85] (2.4,1.44) circle (2.8pt); \fill[orbitblue] (2.72,2.4) circle (2.8pt); \fill[orbitblue] (3.04,2.08) circle (2.8pt); \fill[black!85] (3.36,5.6) circle (2.8pt); \fill[orbitblue] (3.68,2.72) circle (2.8pt); \fill[orbitblue] (4,4.96) circle (2.8pt); \fill[orbitred] (4.32,3.36) circle (2.8pt); \fill[orbitblue] (4.64,4.32) circle (2.8pt); \fill[orbitblue] (4.96,4) circle (2.8pt); \fill[black!85] (5.6,4.64) circle (2.8pt); \draw[orbitblue,line width=.7pt] (4.64,4.32) circle (4.5pt); \draw[orbitblue,line width=.7pt] (4.96,4) circle (4.5pt); \node[orbitred,font=\scriptsize,anchor=west,fill=white,inner sep=1pt] at (4.5,3.36) {$M+1=5$}; \node[orbitred,font=\scriptsize] at (4.32,6.23) {$w^{-1}(M+1)=8$};
\end{scope}
\begin{scope}[xshift=8.15cm]
\node[font=\normalsize] at (2.880,-2.200) {The representative $g_{w(\beta)}w$}; \node at (2.880,-1.750) {$(w(1),\ldots,w(6))=(4,-1,2,1,12,3)$}; \node at (2.880,-1.350) {$\beta=\{i\in w^{-1}([1,4]_6):i<8\}\cup\{9,10\}$}; \node at (2.880,-0.950) {$\sigma=\{9,10\}$}; \draw[black!13,line width=.25pt] (0,0)--(0,5.76); \draw[black!13,line width=.25pt] (0,0)--(5.76,0); \draw[black!13,line width=.25pt] (0.32,0)--(0.32,5.76); \draw[black!13,line width=.25pt] (0,0.32)--(5.76,0.32); \draw[black!13,line width=.25pt] (0.64,0)--(0.64,5.76); \draw[black!13,line width=.25pt] (0,0.64)--(5.76,0.64); \draw[black!13,line width=.25pt] (0.96,0)--(0.96,5.76); \draw[black!13,line width=.25pt] (0,0.96)--(5.76,0.96); \draw[black!13,line width=.25pt] (1.28,0)--(1.28,5.76); \draw[black!13,line width=.25pt] (0,1.28)--(5.76,1.28); \draw[black!13,line width=.25pt] (1.6,0)--(1.6,5.76); \draw[black!13,line width=.25pt] (0,1.6)--(5.76,1.6); \draw[black!13,line width=.25pt] (1.92,0)--(1.92,5.76); \draw[black!13,line width=.25pt] (0,1.92)--(5.76,1.92); \draw[black!13,line width=.25pt] (2.24,0)--(2.24,5.76); \draw[black!13,line width=.25pt] (0,2.24)--(5.76,2.24); \draw[black!13,line width=.25pt] (2.56,0)--(2.56,5.76); \draw[black!13,line width=.25pt] (0,2.56)--(5.76,2.56); \draw[black!13,line width=.25pt] (2.88,0)--(2.88,5.76); \draw[black!13,line width=.25pt] (0,2.88)--(5.76,2.88); \draw[black!13,line width=.25pt] (3.2,0)--(3.2,5.76); \draw[black!13,line width=.25pt] (0,3.2)--(5.76,3.2); \draw[black!13,line width=.25pt] (3.52,0)--(3.52,5.76); \draw[black!13,line width=.25pt] (0,3.52)--(5.76,3.52); \draw[black!13,line width=.25pt] (3.84,0)--(3.84,5.76); \draw[black!13,line width=.25pt] (0,3.84)--(5.76,3.84); \draw[black!13,line width=.25pt] (4.16,0)--(4.16,5.76); \draw[black!13,line width=.25pt] (0,4.16)--(5.76,4.16); \draw[black!13,line width=.25pt] (4.48,0)--(4.48,5.76); \draw[black!13,line width=.25pt] (0,4.48)--(5.76,4.48); \draw[black!13,line width=.25pt] (4.8,0)--(4.8,5.76); \draw[black!13,line width=.25pt] (0,4.8)--(5.76,4.8); \draw[black!13,line width=.25pt] (5.12,0)--(5.12,5.76); \draw[black!13,line width=.25pt] (0,5.12)--(5.76,5.12); \draw[black!13,line width=.25pt] (5.44,0)--(5.44,5.76); \draw[black!13,line width=.25pt] (0,5.44)--(5.76,5.44); \draw[black!13,line width=.25pt] (5.76,0)--(5.76,5.76); \draw[black!13,line width=.25pt] (0,5.76)--(5.76,5.76); \node[font=\scriptsize,anchor=south] at (0.16,-0.05) {$-5$}; \node[font=\scriptsize,anchor=east] at (-0.05,0.16) {$-5$}; \node[font=\scriptsize,anchor=south] at (0.48,-0.05) {$-4$}; \node[font=\scriptsize,anchor=east] at (-0.05,0.48) {$-4$}; \node[font=\scriptsize,anchor=south] at (0.8,-0.05) {$-3$}; \node[font=\scriptsize,anchor=east] at (-0.05,0.8) {$-3$}; \node[font=\scriptsize,anchor=south] at (1.12,-0.05) {$-2$}; \node[font=\scriptsize,anchor=east] at (-0.05,1.12) {$-2$}; \node[font=\scriptsize,anchor=south] at (1.44,-0.05) {$-1$}; \node[font=\scriptsize,anchor=east] at (-0.05,1.44) {$-1$}; \node[font=\scriptsize,anchor=south] at (1.76,-0.05) {$0$}; \node[font=\scriptsize,anchor=east] at (-0.05,1.76) {$0$}; \node[font=\scriptsize,anchor=south] at (2.08,-0.05) {$1$}; \node[font=\scriptsize,anchor=east] at (-0.05,2.08) {$1$}; \node[font=\scriptsize,anchor=south] at (2.4,-0.05) {$2$}; \node[font=\scriptsize,anchor=east] at (-0.05,2.4) {$2$}; \node[font=\scriptsize,anchor=south] at (2.72,-0.05) {$3$}; \node[font=\scriptsize,anchor=east] at (-0.05,2.72) {$3$}; \node[font=\scriptsize,anchor=south] at (3.04,-0.05) {$4$}; \node[font=\scriptsize,anchor=east] at (-0.05,3.04) {$4$}; \node[font=\scriptsize,anchor=south] at (3.36,-0.05) {$5$}; \node[font=\scriptsize,anchor=east] at (-0.05,3.36) {$5$}; \node[font=\scriptsize,anchor=south] at (3.68,-0.05) {$6$}; \node[font=\scriptsize,anchor=east] at (-0.05,3.68) {$6$}; \node[font=\scriptsize,anchor=south] at (4,-0.05) {$7$}; \node[font=\scriptsize,anchor=east] at (-0.05,4) {$7$}; \node[font=\scriptsize,anchor=south] at (4.32,-0.05) {$8$}; \node[font=\scriptsize,anchor=east] at (-0.05,4.32) {$8$}; \node[font=\scriptsize,anchor=south] at (4.64,-0.05) {$9$}; \node[font=\scriptsize,anchor=east] at (-0.05,4.64) {$9$}; \node[font=\scriptsize,anchor=south] at (4.96,-0.05) {$10$}; \node[font=\scriptsize,anchor=east] at (-0.05,4.96) {$10$}; \node[font=\scriptsize,anchor=south] at (5.28,-0.05) {$11$}; \node[font=\scriptsize,anchor=east] at (-0.05,5.28) {$11$}; \node[font=\scriptsize,anchor=south] at (5.6,-0.05) {$12$}; \node[font=\scriptsize,anchor=east] at (-0.05,5.6) {$12$}; \draw[black!40,dashed,line width=.6pt] (0,0) rectangle (1.92,1.92); \draw[black!70,line width=.8pt] (1.92,1.92) rectangle (3.84,3.84); \draw[black!40,dashed,line width=.6pt] (3.84,3.84) rectangle (5.76,5.76); \node[black!45] at (-0.36,-0.36) {$\ddots$}; \node[black!45] at (5.98,5.96) {$\ddots$}; \draw[->,black!55] (0,-0.490)--(5.760,-0.490) node[right] {$i$}; \draw[->,black!55] (-0.580,0)--(-0.580,5.760) node[below] {$j$}; \draw[orbitred!65,densely dashed] (4.320,0)--(4.320,5.760); \draw[orbitred!65,densely dashed] (0,3.360)--(4.320,3.360); \fill[black!85] (0.16,1.12) circle (2.8pt); \fill[black!85] (0.8,0.48) circle (2.8pt); \fill[black!85] (1.12,0.16) circle (2.8pt); \fill[black!85] (1.44,3.68) circle (2.8pt); \fill[black!85] (1.76,0.8) circle (2.8pt); \fill[black!85] (2.08,3.04) circle (2.8pt); \fill[black!85] (2.4,1.44) circle (2.8pt); \fill[black!85] (2.72,2.4) circle (2.8pt); \fill[black!85] (3.04,2.08) circle (2.8pt); \fill[black!85] (3.36,5.6) circle (2.8pt); \fill[black!85] (3.68,2.72) circle (2.8pt); \fill[black!85] (4,4.96) circle (2.8pt); \fill[orbitred] (4.32,3.36) circle (2.8pt); \fill[black!85] (4.64,4.32) circle (2.8pt); \fill[black!85] (4.96,4) circle (2.8pt); \fill[black!85] (5.6,4.64) circle (2.8pt); \fill[orbitblue] (0.48,0.16) circle (2.8pt); \fill[orbitblue] (2.4,0.16) circle (2.8pt); \fill[orbitblue] (4.32,0.16) circle (2.8pt); \fill[orbitblue] (0.48,0.48) circle (2.8pt); \fill[orbitblue] (2.4,0.48) circle (2.8pt); \fill[orbitblue] (4.32,0.48) circle (2.8pt); \fill[orbitblue] (2.4,0.8) circle (2.8pt); \fill[orbitblue] (4.32,0.8) circle (2.8pt); \fill[orbitblue] (0.48,1.12) circle (2.8pt); \fill[orbitblue] (2.4,1.12) circle (2.8pt); \fill[orbitblue] (4.32,1.12) circle (2.8pt); \fill[orbitblue] (2.4,2.08) circle (2.8pt); \fill[orbitblue] (4.32,2.08) circle (2.8pt); \fill[orbitblue] (2.4,2.4) circle (2.8pt); \fill[orbitblue] (4.32,2.4) circle (2.8pt); \fill[orbitblue] (4.32,2.72) circle (2.8pt); \fill[orbitblue] (2.4,3.04) circle (2.8pt); \fill[orbitblue] (4.32,3.04) circle (2.8pt); \fill[orbitblue] (4.32,4) circle (2.8pt); \fill[orbitblue] (4.32,4.32) circle (2.8pt); \fill[orbitblue] (4.32,4.96) circle (2.8pt);
\end{scope}
\node[anchor=north west,align=left,text width=6.7cm,font=\footnotesize] at (0,6.91) { \textcolor{orbitblue}{Blue}: $i\in\beta$. Circled: $i\in\sigma=\{9,10\}$.\\[3pt]
\textcolor{orbitred}{Red}: row $5$, column $8$. }; \node[anchor=north west,align=left,text width=6.5cm,font=\footnotesize] at (8.15,6.91){\textcolor{orbitblue}{Blue}: the additional entries at\\[2pt]row $w(i)+6k$, column $8+6k$,\\[2pt]$i\in\beta$, $k\in\mathbb Z$.};
\end{tikzpicture}
} \caption{Another example of Proposition~\ref{prop 3.2.7}. The red dot has row $M+1$ and column $w^{-1}(M+1)$. The solid square marks rows and columns $1,\ldots,N$.}
\label{fig:affine-orbit-classification-dispersed}
\end{figure}

\FloatBarrier
\subsubsection{Relevance} Similar to \cite{Po}, we have the following proposition.
\begin{prop}
An orbit corresponding to $(w, \beta)$ is relevant if and only if,
\begin{enumerate}
\item for any $k\in\{M+2,M+3,\cdots,N\}$, we have $w^{-1}(k-1)>w^{-1}(k)$;
\item for any $i\in \{1,2,\cdots, M\}$, if $N\BZ+i$ does not contain any element of $w(\sigma)$, then $b(\lambda,-)$ restricts trivially to the $i$-th coordinate torus. If $b(\lambda,\mu)=\eta^{\sum_{j=1}^N\lambda_j\mu_j}$, where $\eta$ has order $c$, this condition is $\eta^{n_i}=1$, or equivalently $c\mid n_i$.
\end{enumerate}
\end{prop}
\begin{proof}
We need to prove that the conditions in the statement are equivalent to the following two conditions on the stabilizer subgroup of $I_M\ltimes U_{M,N}(\bF)$ at $g_{w(\beta)}\cdot w\in\Fl_N$:
\begin{enumerate}
\item it belongs to the kernel of $\psi$;
\item the character sheaf induced on it by the twisting gerbe is trivial.
\end{enumerate}
The first condition is equivalent to the condition (1) in the statement.

For the second condition, let us consider the reductive quotient of the stabilizer subgroup of $I_M\ltimes U_{M,N}(\bF)$, which is $T_{w,\sigma}\simeq\BG_m^{M-|\sigma|}$. Thus, the second condition is the triviality of $b(\lambda,-)$ on its coordinate tori, which is (2).
\end{proof}
For example, the orbit in Figure~\ref{fig:affine-orbit-classification} is not relevant, since $w^{-1}(5)=2<5=w^{-1}(6)$, while the orbit in Figure~\ref{fig:affine-orbit-classification-dispersed} is relevant. Indeed, we have $w^{-1}(5)=8>-1=w^{-1}(6)$, $\sigma=\{9,10\}$, and $\lambda=(0,0,0,0,1,-1)$. Thus, $b(\lambda,-)$ is trivial on $T_{w,\sigma}$, which is the product of the third and fourth coordinate tori.

\begin{rem}\label{rem 3.2.9}
Let $\cB_\lambda$ be the character sheaf on $T_N$ determined by $b(\lambda,-)$, see \cite[\S 9.5.1]{[GL]}. The twisting condition is $\cB_\lambda|_{T_{w,\sigma}}\simeq\mathbb 1$.
\end{rem}
For future use, we also need a monodromic version of the above proposition, and the proof is the same.

\begin{defn}
For $w=t^\lambda w^{\fin}\in W_N^{\ext}$ and a multiplicative rank-one character sheaf $\chi$ on $T_N$, we define
\[
w(\chi):=\cB_\lambda\otimes w^{\fin}(\chi).
\]
\end{defn}

\begin{prop}\label{prop 3.2.10}
An orbit corresponding to $(w, \beta)$ is relevant in $D_{\kappa}^{(I_M,\chi)\ltimes (U_{M,N}(\bF),\psi)}(\GL_N(\bF)/I_N,\chi')$ if and only if,
\begin{enumerate}
\item for any $k\in\{M+2,M+3,\cdots,N\}$, we have $w^{-1}(k-1)>w^{-1}(k)$;
\item $w(\chi')$ and $\chi$ have the same restriction to $T_{w,\sigma}$.
\end{enumerate}
\end{prop}
\subsubsection{Periodic characters}\label{sec:periodic-characters} We can regard any element in $\GL_N(\bF)$ as an $N$-periodic $\BZ\times\BZ$ matrix. In particular, $w$ is the permutation matrix with non-zero entries at $(w(i),i)$, $i\in\BZ$, and $g_\beta$ is regarded as $\Id+\sum_{i\in\beta,k\in\BZ}e_{i+kN}\otimes e^*_{M+1+kN}$. Also, we can regard the characters $\chi'$ and $\chi$ as sequences $\widetilde\chi'$ on $\BZ$ and $\widetilde\chi$ on $[1,M]_N$, such that $\widetilde\chi'(i+N)=\eta\widetilde\chi'(i)$ and $\widetilde\chi(i+N)=\eta\widetilde\chi(i)$, where $(\cB_\lambda)_i=\eta^{\lambda_i}$. Note that, for our doubled extension, we have $\eta=\kappa^2$ by Lemma~\ref{lem:affine-torus-lift-convention}. Then, the affine Weyl action is given by
\[
(w\widetilde\chi')_i=\widetilde\chi'(w^{-1}(i)) =\eta^{\lambda_i}\chi'_{(w^{\fin})^{-1}(i)}, \qquad w=t^\lambda w^{\fin}.
\]
In the action formulas, we omit the tildes. In this notation, condition (2) of Proposition~\ref{prop 3.2.10} becomes $\chi_{w(i)}=\chi'_i$ whenever $i\in w^{-1}([1,M]_N)$ and $(i+N\BZ)\cap\sigma=\varnothing$.
\subsubsection{The semi-infinite version} We include this subsection for completeness, and it is not used in the rest of this paper.

Assume $M<N$. Let $G=\GL_M\times\GL_N$ and let $H\subset G$ be the constant Gaiotto subgroup. Let $\Fl_G^{\seminf}$ be the semi-infinite flag variety associated with $G$, i.e., $\Fl_G^{\seminf}=G(\bF)/U(\bF)T(\bO)$. We choose representatives $z_r$ of the finitely many $H$-orbits on $G/B$, and let
\[
S_r=z_r^{-1}Hz_r\cap B,\qquad \Gamma_r=\operatorname{im}(S_r\longrightarrow T),\qquad \Lambda_r=\operatorname{im}\bigl(X_*(\Gamma_r)\longrightarrow\Lambda\bigr).
\]
\begin{prop}\label{prop:semi-infinite-orbit-classes}
The $H(\bF)$-orbits on $\Fl_G^{\seminf}$ have representatives $z_rt^\lambda U(\bF)T(\bO)$, indexed by $r$ and $\bar\lambda\in\Lambda/\Lambda_r$, where $\lambda$ is any lift.
\end{prop}
\begin{proof}
Each $S_r$ is an extension of a split torus by a split unipotent group, thus, $H^1(\bF,S_r)=1$. Hence, the elements $z_r$ represent all the double cosets in $H(\bF)\backslash G(\bF)/B(\bF)$. Since $B(\bF)/U(\bF)T(\bO)\simeq\Lambda$ and $S_r(\bF)\to\Gamma_r(\bF)$ is surjective, by taking valuations, we obtain
\[
S_r(\bF)\backslash B(\bF)/U(\bF)T(\bO)\simeq\Lambda/\Lambda_r.
\]
\end{proof}
The relevance of these orbits is determined by restricting the additive character and the twisting character to the stabilizer of $z_rt^\lambda U(\bF)T(\bO)$.

\subsection{From finite type category to Iwahori Gaiotto category} In this subsection, we want to study the Grothendieck group of the Iwahori Gaiotto category with the help of the results about the finite type Gaiotto category.

The starting point is the following general lemma. In our case, $G=\GL_M\times \GL_N$, $H=\GL_M\ltimes U_{M,N}$. Recall that we make the convention that the character $\psi$ which appears in the calculation of the canonical basis is given by taking the coefficient of the constant term, rather than the coefficient of $t^{-1}$, see Section~\ref{sec:type-a-conventions}.
\begin{lem}\label{13.2.1}
For any $\kappa$ and $\chi$, there is an equivalence
\begin{equation}
D(H,\psi\backslash G/B,\chi)\simeq D_{\kappa}(H(\bF),\psi\backslash H(\bF)G(\bO)/I,\chi).
\end{equation}
\end{lem}
\begin{proof}

The lemma follows from the following equivalences
\begin{multline}
D_{\kappa}(H(\bF),\psi\backslash H(\bF)G(\bO)/I,\chi)\simeq D_{\kappa}(H(\bO),\psi\backslash G(\bO)/I,\chi)\\ \simeq D(H,\psi\backslash (H_1\backslash G(\bO)/I),\chi) \simeq D(H,\psi\backslash G/B,\chi).
\end{multline}

Here, the second equivalence follows from the fact that the trivializations on $H(\bO)$ and $I$ are the restrictions of the canonical trivialization over $G(\bO)$, and the third equivalence follows from the fact that the group $H_1:=\ker(H(\bO)\to H)$ is pro-unipotent, acts trivially on $G(\bO)/I\simeq G/B$, and lies in $\ker\psi$.
\end{proof}

We denote by $\sO_{0}:= H(\bF)G(\bO)/I\subset \Fl_G$ the preimage of the unit $H(\bF)$-orbit in $\Gr_G$ under the natural projection $\Fl_G\rightarrow \Gr_G$.

\begin{prop}\label{prop:finite-subcategory-clean}
There is no relevant orbit in the boundary of $\sO_{0}$ for any $\chi, \chi'$, and $\kappa$.
\end{prop}
\begin{proof}
We consider the case when $\kappa=1$ and $\chi,\chi'$ are trivial, and the general case is similar, since we only use the $\psi$-relevance condition below.

Note that, in the upper triangular convention, we use decreasing orbit parameters $\xi^{(1)},\xi^{(2)}$. Namely, the $\GL_M(\bO)\ltimes U_{M,N}(\bF)$-orbits in $\Gr_{N}$ are indexed by two sequences $\xi^{(1)}=(\xi^{(1)}_1\geq\xi^{(1)}_2\geq\cdots\geq\xi^{(1)}_M)$ and $\xi^{(2)}=(\xi^{(2)}_1\geq\xi^{(2)}_2\geq\cdots\geq \xi^{(2)}_{M+1}, \xi^{(2)}_{M+2},\cdots, \xi^{(2)}_N)$, see \cite[Lemma 3.3.3]{[TY]}, after inverse transpose and changing the signs of the orbit parameters. If the preimage of such an orbit in $\Fl_N$ contains an $(I_M\ltimes U_{M,N}(\bF),\psi)$-relevant orbit, then
\begin{equation}\label{rel gi}
\xi^{(2)}_{M+1}\geq \xi^{(2)}_{M+2}\geq \xi^{(2)}_{M+3}\geq\cdots\geq \xi^{(2)}_N.
\end{equation}
(We remark here that if one uses the additive character in loc. cit., then the above condition will become $\xi^{(2)}_{M+1}\geq \xi^{(2)}_{M+2}-1\geq \xi^{(2)}_{M+3}-2\geq\cdots\geq \xi^{(2)}_N-(N-M-1)$.)

However, we claim that the boundary of the unit $\GL_M(\bO)\ltimes U_{M,N}(\bF)$-orbit does not contain any orbit corresponding to $(\xi^{(1)},\xi^{(2)})$ satisfying the above condition.

Indeed, note that the above condition implies that the $\GL_M(\bO)\ltimes U_{M,N}(\bF)$-orbit corresponding to $(\xi^{(1)},\xi^{(2)})$ in $\Gr_N$ is relevant with respect to the additive character in loc. cit., which is given by taking the coefficient of $t^{-1}$, while the unit $\GL_M(\bO)\ltimes U_{M,N}(\bF)$-orbit is minimal relevant.

For completeness, we include here the proof that the unit orbit is minimal relevant.

We combine $\xi^{(1)}$ and $\xi^{(2)}$ into the sequence
\[
\xi^{(2)}_1,\xi^{(1)}_1,\ldots, \xi^{(2)}_M,\xi^{(1)}_M, \xi^{(2)}_{M+1},\ldots,\xi^{(2)}_N.
\]
The closure relation of $\GL_M(\bO)\ltimes U_{M,N}(\bF)$-orbits is studied in \cite[Proposition 8.2.1]{[TY]}: if the orbit corresponding to $(\xi^{(1)},\xi^{(2)})$ lands in the closure of the unit orbit, then the partial sum $S_s$ of the first $s$ terms is no bigger than $0$ for any $1\leq s\leq M+N$, and is $0$ if $s=M+N$. Furthermore, the condition \eqref{rel gi} says $\xi^{(2)}_{M+1}\geq \xi^{(2)}_{M+2}\geq \xi^{(2)}_{M+3}\geq\cdots\geq \xi^{(2)}_N$.

The above conditions imply that $(\xi^{(1)}, \xi^{(2)})=(0,0)$. Indeed, since $S_1=\xi^{(2)}_1\leq0$ and $\xi^{(2)}$ is non-increasing, all its entries are non-positive. On the other hand, $\xi^{(2)}_N=S_{M+N}-S_{M+N-1}\geq0$. Thus, $\xi^{(2)}=0$. For $M>0$, the remaining partial sums are $S_{2k}=\sum_{j=1}^k\xi^{(1)}_j$. In particular, $\xi^{(1)}_1=S_2\leq0$, so all entries of $\xi^{(1)}$ are non-positive by monotonicity. Their sum is $S_{M+N}=0$, and hence $\xi^{(1)}=0$ as well.

Thus, there is no orbit in the boundary of the unit $\GL_M(\bO)\ltimes U_{M,N}(\bF)$-orbit in $\Gr_N$ whose preimage in $\Fl_N$ contains a relevant $(I_M\ltimes U_{M,N}(\bF),\psi)$-orbit.
\end{proof}

\subsection{Action of affine Yokonuma algebra} Recall that in the $\kappa$-twisted setting, our sheaves and functions are over the (doubled) Kac--Moody extension $\hat \GL_N(\bF)$. Before we state the action formula in the affine case, we need to introduce some sections/liftings of certain elements in the loop group to $\hat \GL_N(\bF)$.

We use the following notation, cf. \cite{[GGK]}. For $1\leq j\leq N$ and $a\in\BG_m$, we denote by $\pi_j(a)$ the canonical lift of the diagonal matrix with $a$ in position $j$ and $1$ elsewhere, using the splitting over $T_N(\bO)$, and write $\pi_j=\pi_j(-1)$.

For $v\in\bF^M$, let $g(v):=\Id+v e_{M+1}^*$, and denote by $\underline g(v)$ its canonical lift. Note that $g_A=g(e_A)$ for $A\subset[1,M]_N$ bounded above.

For any real root $\check\alpha$, consider the corresponding $1$-dimensional subgroup in $\GL_N(\bF)$, which has a canonical (and unique multiplicative) section in $\hat \GL_N(\bF)$. We denote by $e_{\check\alpha}(a)$ the root subgroup element, and by $\underline e_{\check\alpha}(a)$ its canonical lift. Here, for an affine real root, the parameter $a$ belongs to the ground field, and the power of $t$ is included in the root.

We define
\begin{equation}
w_{\check\alpha}(a) := \underline e_{\check\alpha}(a)\cdot \underline e_{-\check\alpha}(-a^{-1})\cdot \underline e_{\check\alpha}(a)
\end{equation}
We denote $w_{\check\alpha}:=w_{\check\alpha}(1)$, and define
\begin{equation}\label{h sect}
\underline{h}_{\check\alpha}(a) := w_{\check\alpha}(a)\cdot w_{\check\alpha}(-1)\in \hat T_N(\bF).
\end{equation}

For $N\geq2$ and $0\leq i\leq N-1$, we take $\check\alpha_i$ to be a simple (possibly affine) root and set
\[
\underline{s_i}:= {\pi}_{i+1} w_{\check\alpha_i}(1).
\]

We extend the liftings $\underline{s_i}$ to the whole $W_N^\ext$, just as points, i.e., we do not assume that they preserve the group structure. Also, we take the lifting of $g_{w(\beta)}\in U_N(\bF)$ to $\hat\GL_N(\bF)$ for any $\beta$ using the canonical splitting over $U_N(\bF)$. We denote by $\underline{\tilde{w}}=\underline{g_{w(\beta)}}\underline{w}$ the product of the corresponding sections $\underline{w}$ and $\underline{g_{w(\beta)}}$.

For $w_1,w_2\in W_N^{\ext}$, we define $z(w_1,w_2)$ and $c(w_1,w_2)$ by
\[
\underline{w_1}\,\underline{w_2} =z(w_1,w_2)\,\underline{w_1w_2},\qquad c(w_1,w_2)=\kappa(z(w_1,w_2))^{-1}.
\]
Here, we have $f(zg)=\kappa(z)f(g)$ for our trace functions.

The torus element $h_{\check\alpha_i}(a)$, $a\in\BG_m$, has two lifts, namely, $\underline h_i^{\mathrm{Iw}}(a)$ given by the splitting over the Iwahori subgroup, and $\underline h_i(a):=\underline h_{\check\alpha_i}(a)$ defined in \eqref{h sect}. Their ratio is computed by the tame symbol as follows.
\begin{lem}\label{lem:affine-torus-lift-convention}
We use the tame-symbol convention
\[
(f,g)_{\mathrm{tame}} =(-1)^{\operatorname{val}(f)\operatorname{val}(g)} \left.\frac{f^{\operatorname{val}(g)}} {g^{\operatorname{val}(f)}}\right|_{t=0}, \qquad (t,a)_{\mathrm{tame}}=a^{-1}.
\]
Then, for the doubled extension with integral bilinear form $b_{\BZ}(\lambda,\mu)=2\sum_i\lambda_i\mu_i$, we have
\[
[\underline{t^\lambda},\pi_i(a)]=a^{-2\lambda_i}, \qquad [x,y]:=xyx^{-1}y^{-1}.
\]
Consequently, the coordinate monodromy is $\eta=\kappa^2$. Moreover, we have
\[
\underline h_0(a)=a^{-2}\underline h_0^{\mathrm{Iw}}(a), \qquad \underline h_i(a)=\underline h_i^{\mathrm{Iw}}(a) \quad(1\leq i\leq N-1).
\]
\end{lem}
\begin{proof}
Note that the torus commutator is given by the tame symbol and $b_{\BZ}$. Thus, for $w=t^\lambda w^{\fin}$, we have
\[
\underline w^{-1}\pi_i(a)\underline w =a^{2\lambda_i}\pi_{(w^{\fin})^{-1}(i)}(a),
\]
and hence the left character is $\kappa^{2\lambda_i}\chi'_{(w^{\fin})^{-1}(i)}$. For the affine root, let $d=\diag(1,\ldots,1,t)$. Then, we have
\[
\underline h_0(a) =\underline d\,\underline h_{-\check\theta}(a)\,\underline d^{-1} =a^{-2}\underline h_0^{\mathrm{Iw}}(a),
\]
where $\check\theta$ is the highest root and $h_{-\check\theta}(a)=\diag(a^{-1},1,\ldots,1,a)$. For finite simple roots, both lifts agree. The same computation applies to the embedded $\GL_M$.
\end{proof}

We extend the notation for torus lifts to periodic indices by
\begin{equation}\label{eq:periodic-torus-lifts}
\pi_{j+N}(a)=a^2\pi_j(a).
\end{equation}
By Lemma~\ref{lem:affine-torus-lift-convention}, we have
\[
\underline w\pi_j(a)\underline w^{-1}=\pi_{w(j)}(a), \qquad \underline h_i(a)=\pi_i(a)\pi_{i+1}(a^{-1}),
\]
also for $i=0$. For the left affine simple root, we take the successor in $[1,M]_N$ in the second identity.

\subsubsection{Twisted affine Yokonuma algebra} Similarly to the finite case, given two characters $\chi,\chi'$ of $T_N$, we consider the space of $\BC$-valued functions on $\hat\GL_N(\bF)$ which are left $(\BG_m,\kappa)\times(I_N,\chi)$-invariant and right $(I_N,\chi')$-invariant, with compact support modulo the two Iwahori subgroups and the center, and denote it by $\bH_{\kappa,\chi,\chi'}^{\aff}(N)_\BC$. We have the convolution
\begin{equation}\label{10.54}
\bH_{\kappa,\chi,\chi'}^{\aff}(N)_\BC\otimes \bH_{\kappa,\chi',\chi''}^\aff(N)_\BC\rightarrow \bH_{\kappa,\chi,\chi''}^\aff(N)_\BC.
\end{equation}

Given $w\in W^\ext_N$, there is a non-zero $(\kappa,\chi,\chi')$-invariant function on $\hat I_N\underline{w}I_N$ if and only if $\chi=w(\chi')$. When $\chi=w(\chi')$, we denote by $T_{w,\kappa,\chi,\chi'}$ the characteristic $(\kappa,\chi,\chi')$-invariant function on the $\hat I_N\times I_N$-orbit in $\hat\GL_N(\bF)$ corresponding to $\underline{w}$, and normalized by the requirement that its value at $\underline{w}\in \hat\GL_N(\bF)$ is $1$, and set $T_{w,\kappa,\chi,\chi'}:=0$ otherwise. We omit $\kappa$ in the notation as the twisting $\kappa$ is clear. We also denote $T_{w,\chi}:=T_{w,w(\chi),\chi}$, and use the same convention for $H_{w,\chi}^{\mxd}$ and $H_{w,\chi}^{\gr}$.

\begin{lem}
\begin{enumerate}
\item if $w_1,w_2\in W^{\ext}_N$ and $l(w_1)+l(w_2)=l(w_1w_2)$, then we have
\[T_{w_1, w_1(w_2(\chi)),w_2(\chi)}T_{w_2,w_2(\chi), \chi}=c(w_1,w_2)\cdot T_{w_1w_2, w_1(w_2(\chi)),\chi}.\]
\item if $s=s_i=(i,i+1)$, then we have
\[T_{s, \chi, s\chi}T_{s, s\chi,\chi}=Q\cdot c(s,s)\cdot T_{e,\chi,\chi}+(Q-1)\cdot c(s,s)\cdot\chi_{i}(-1)\cdot T_{s, \chi,\chi}.\]
\end{enumerate}
\end{lem}
\begin{proof}
We only show the second one, and calculate the coefficient of $T_{s,\chi,\chi}$. Note that, if $a\neq 0$, we have \[ \underline{s_i}\underline e_{\check\alpha_i}(a)\underline{s_i}=z(s_i,s_i)\cdot \underline e_{\check\alpha_i} (a^{-1})\cdot \pi_{i+1}\cdot \underline{s_i}\cdot \underline h_i(-a)\cdot \underline e_{\check\alpha_i}(a^{-1}),
\]
and the coefficient of $T_{s,\chi,\chi}$ is given by $\sum_{a\in \BF_Q\setminus\{0\}} c(s_i,s_i)\cdot \chi_{i}(a^{-1})\chi_{i+1}(a)\chi_{i}(-1)=(Q-1)\cdot c(s_i,s_i)\cdot \chi_{i}(-1)$. Here, we use the fact that if $T_{s,\chi,\chi}$ is relevant, then $\chi_i=\chi_{i+1}$.
\end{proof}

\subsubsection{The left and right actions} The convolution gives rise to a bimodule structure on the direct sum over the characters. Namely,
\begin{equation}
\bH_{\kappa,\chi'',\chi}^\aff(M)_\BC\otimes \bK_{\kappa,\chi,\chi'}^\aff(M|N)_\BC\otimes \bH_{\kappa,\chi',\chi'''}^\aff(N)_\BC \rightarrow \bK_{\kappa,\chi'',\chi'''}^\aff(M|N)_\BC.
\end{equation}
The action of length-zero elements is given by left and right translation, with the scalars from the chosen lifts. It remains to determine the action of the basis elements $T_{s,\chi,\chi'}$ associated with simple reflections $s$.

\subsubsection{Standard basis} Assume $M<N$. Let $H=I_M\ltimes U_{M,N}(\bF)$, and let $\widehat H$ be its preimage in the central extension. Here, all groups are taken on $\BF_Q$-points. We normalize the Haar measure by giving $I_M\ltimes U_{M,N}(\bO)$ volume $1$, and the central kernel $\BF_Q^\times$ total mass $1$. Let $\theta$ be the character determined by $\kappa,\chi,\psi$, and for $g\in\hat\GL_N(\bF)$, let $\delta_{g,\chi'}$ be the function supported on $gI_N$, with value $\chi'(b)$ at $gb$. We define
\[
\hat T_{g,\chi,\chi'}(x) :=\int_{\widehat H}\theta(a)\delta_{g,\chi'}(a^{-1}x)\,da.
\]
Note that the stabilizer of $gI_N$ in $\widehat H$ is compact open, and $\hat T_{g,\chi,\chi'}$ is non-zero if and only if the orbit corresponding to $g$ is relevant.

If the orbit corresponding to $g\in \hat \GL_N(\bF)$ is relevant, we let $T_{g,\chi,\chi'}$ be the left $(\BG_m,\kappa)\times (I_M,\chi)\ltimes (U_{M,N}(\bF),\psi)$-invariant and right $(I_N,\chi')$-invariant function on the $\BG_m\times I_M\ltimes U_{M,N}(\bF)\times I_N$-orbit in $\hat\GL_N(\bF)$, which is normalized such that its value at $g$ is $1$. If the corresponding orbit is irrelevant, we let $T_{g,\chi,\chi'}=0$. Thus, \[\hat T_{g,\chi,\chi'}=\on{vol}(\Stab_{\widehat H}(gI_N))\cdot T_{g,\chi,\chi'}.\]

Obviously, we can regard $\hat T_{g,\chi,\chi'}$ and $T_{g,\chi,\chi'}$ as right $(T_N,\chi')$-equivariant functions on the punctured square of the determinant line over $\GL_N(\bF)/I_N^0$, with central character $\kappa$, where $I_N^0=\ker(I_N\to T_N)$.

\subsubsection{Relative dimensions} Given an orbit $\tilde{w}=(w,\beta)$, we denote by $l(\tilde{w})$ the relative dimension of the $I_M\ltimes U_{M,N}(\bF)$-orbit in $\Fl_N$ corresponding to $\tilde{w}$ with respect to the unit orbit, i.e., the orbit corresponding to $(\id, [1,M]_N\cap \{i|i\leq M+1\})$. By a direct calculation, we have the following lemma.
\begin{lem}\label{lem:affine-orbit-relative-dimension}
We have $l(\tilde{w})=l^{I_M\ltimes U_{M,N}(\bF)}(w)+\#\{i\in \beta|i> w^{-1}(M+1)\}$. Here,  $l^{I_M\ltimes U_{M,N}(\bF)}(w)$ is the relative dimension of the $I_M\ltimes U_{M,N}(\bF)$-orbit through $wI_N\in\Fl_N$.     Explicitly, $l^{I_M\ltimes U_{M,N}(\bF)}(w)$ equals
\begin{multline}
\#\{(i,j)\in \BZ^2|1\leq j\leq N, i\leq_{I_M\ltimes U_{M,N}(\bF)} j, w^{-1}(i)>_{I_N}w^{-1}(j), i\leq_{I_N}j\}\\
-\#\{(i,j)\in \BZ^2|1\leq j\leq N, i\leq_{I_M\ltimes U_{M,N}(\bF)} j, w^{-1}(i)\leq_{I_N}w^{-1}(j), i>_{I_N}j\}.
\end{multline}

Furthermore, the corresponding orbit is isomorphic to the product of $\BG_m^{|\sigma|}$ and an (infinite-dimensional, if $M<N-1$) affine space.
\end{lem}

\subsubsection{Convolution formulas} We denote by $\hat g_{\tilde v}$ the chosen representative of the orbit indexed by $\tilde v$. Then, the right and left convolution coefficients are given respectively by
\begin{equation}\label{eq:affine-exact-convolution-sum}
\begin{split}
(T_{\tilde u}T_{s_i})(\hat g_{\tilde v}) &=\sum_{a\in\BF_Q} T_{\tilde u}\bigl(\hat g_{\tilde v}\underline e_{\check\alpha_i}(-a)\underline s_i^{-1}\bigr),\\
(T_{s_i}T_{\tilde u})(\hat g_{\tilde v}) &=\sum_{a\in\BF_Q} T_{\tilde u}\bigl(\underline s_i^{-1}\underline e_{\check\alpha_i}(-a)\hat g_{\tilde v}\bigr).
\end{split}
\end{equation}
These formulas follow from the Iwahori coset decompositions of $I_N\underline s_iI_N$ and $I_M\underline s_iI_M$, where each Iwahori subgroup has volume $1$. Here, the character indices are as in \eqref{10.54}, and the torus factors are given by \eqref{eq:periodic-torus-lifts}.

In the comparisons below, $ws>_{I_M\ltimes U_{M-1,N}(\bF)}w$ means that the relative dimension of the $I_M\ltimes U_{M-1,N}(\bF)$-orbit through $wsI_N\in\Fl_N$ is one greater than that of the orbit through $wI_N$, and similarly on the left.

Now, we give the action formula for the affine Yokonuma algebras on $\bK_{\kappa,\chi,\chi'}^{\aff}(M|N)_{\BC}$. In the case when $M=N$ and $\kappa,\chi,\chi'$ are trivial, it is given in \cite[Proposition 2]{[FGT]}. For the right action, let $s=s_i$ and $\iota=\{i,i+1\}$ in periodic notation.

In the left diagram of Figure~\ref{fig:affine-orbit-classification}, the operation $(w,\beta)\mapsto(ws_i,s_i\beta)$ exchanges the columns $i+kN$ and $i+1+kN$ for every $k\in\BZ$, together with their blue dots, and then we take $\sigma$ for the new pair. We denote by $\tilde w s$ the index of the new orbit, which satisfies Proposition~\ref{prop 3.2.7}, and put $\sigma'=\sigma(\tilde w s)$. For $i\in w^{-1}([1,M+1]_N)$, let us denote by $\sigma_i$ the maximal integer in $[i]_N\cap(\beta\cup\{w^{-1}(M+1)\})$ whenever this set is non-empty.

\begin{prop}
Let $\tilde w$ be relevant. Then, we have the following right action formulas.

\begin{enumerate}
\item If $ws >_{I_{M}\ltimes U_{M-1,N}(\bF)} w$,
\begin{enumerate}
\item if there is no $k$ such that $\iota+kN \subseteq \sigma'\cup\{(ws)^{-1}(M+1)\}$, then
\[\hat T_{\tilde w,\chi,\chi'} T_{s,\chi',s\chi'} =Q\cdot c(w,s)\cdot \hat T_{\tilde{w}s,\chi,s\chi'},\]
\item if there is $k$ such that $\iota+kN\subseteq \sigma'\cup\{(ws)^{-1}(M+1)\}$,
\[\hat T_{\tilde w,\chi,\chi'} T_{s,\chi',s\chi'}=(Q-1)\cdot c(w,s)\cdot\hat T_{\tilde{w}s,\chi,s\chi'} + c(w,s)\cdot\hat T_{ws,s\beta\setminus\{i+1+kN\},\chi,s\chi'},\]
\end{enumerate}
\item If $ws<_{I_{M}\ltimes U_{M-1,N}(\bF)} w$, and $w(i),w(i+1)\in[1,M]_N$,
\begin{enumerate}
\item If $\sigma_i+1\in \beta\setminus\sigma$ then
\[
\hat T_{\tilde w,\chi,\chi'} T_{s,\chi',s\chi'} = c(w,s)\cdot\hat T_{\tilde{w}s,\chi,s\chi'} + \chi'_i(-1)G(\mathbb{1},\chi'_{i}/\chi'_{i+1})\cdot\hat T_{\tilde w,\chi,s\chi'}.
\]
\item If $\sigma_i+1\notin \beta$,
\[
\hat T_{\tilde w,\chi,\chi'} T_{s,\chi',s\chi'} = c(w,s)\cdot\hat T_{ws,s\beta\cup\{\sigma_i\}, \chi,s\chi'} + (Q-1)\cdot\chi_{{w}(i+1)}(-1)\cdot\hat T_{w,\beta\cup \{\sigma_i+1\},\chi,s\chi'}.
\]
\item If $\sigma_i+1\in \sigma$,
\begin{equation*}
\begin{aligned}
\hat T_{\tilde w,\chi,\chi'} T_{s,\chi',s\chi'} &= c(w,s)\cdot\hat T_{ws,s\beta,\chi,s\chi'}\\
&{}+ \chi_{w(i+1)}(-1)\cdot\hat T_{w,\beta\setminus\{\sigma_i+1\},\chi,s\chi'}\\
&{}+ \chi_{w(i+1)}(-1) J( {\chi_{w(i+1)}}/{\chi_{i+1}'}, {\chi_i'}/{\chi_{w(i+1)}}) \cdot\hat T_{\tilde w,\chi,s\chi'}.
\end{aligned}
\end{equation*}
\end{enumerate}
\item If $ws<_{I_{M}\ltimes U_{M-1,N}(\bF)} w$, and $w(i)\in[M+2,N]_N$,
\begin{enumerate}
\item if $w(i+1)=w(i)-1$, then
\[
\hat T_{\tilde w,\chi,\chi'} T_{s,\chi',s\chi'} = \chi'_{i+1}(-1)G(\psi,\chi'_i/\chi'_{i+1})\hat T_{\tilde w,\chi,s\chi'};\]
\item otherwise \[ \hat T_{\tilde w,\chi,\chi'} T_{s,\chi',s\chi'} = c(w,s)\cdot\hat T_{ws,s\beta,\chi,s\chi'} + \chi'_{i+1}(-1)G(\mathbb{1},\chi'_i/\chi'_{i+1}) \hat T_{\tilde{w},\chi,s\chi'}.
\]
\end{enumerate}

\item If $ws<_{I_{M}\ltimes U_{M-1,N}(\bF)} w$, and $w(i)\in [M+1]_N$,
\begin{enumerate}
\item if $\sigma_i+1=w^{-1}(M+1)+1\in \sigma$,
\begin{equation*}
\begin{aligned}
\hat T_{\tilde w,\chi,\chi'} T_{s,\chi',s\chi'} &= c(w,s)\cdot\hat T_{ws,s\beta,\chi,s\chi'}\\
&{}+ \chi'_{i+1}(-1)\hat T_{w,\beta\setminus\{\sigma_i+1\},\chi,s\chi'}\\
&{}+ \frac{\chi_{w(i+1)}\chi_i'}{\chi'_{i+1}}(-1) J( {\chi_{w(i+1)}}/{\chi_{i+1}'}, {\chi_i'}/{\chi_{w(i+1)}}) \cdot\hat T_{\tilde w,\chi,s\chi'}
\end{aligned}
\end{equation*}
\item if $w^{-1}(M+1)+1\notin \beta$
\[
\hat T_{\tilde w,\chi,\chi'} T_{s,\chi',s\chi'} = c(w,s)\cdot \hat T_{ws,s\beta\cup\{\sigma_i\},\chi,s\chi'} + (Q-1)\cdot\chi'_{i}(-1) \hat T_{w,\beta\cup\{\sigma_i+1\},\chi,s\chi'}.
\]
\item if $w^{-1}(M+1)+1\in \beta\setminus \sigma$,
\[
\hat T_{\tilde w,\chi,\chi'} T_{s,\chi',s\chi'} = c(w,s)\cdot \hat T_{ws,s\beta,\chi,s\chi'} + \chi'_{i}(-1)G(\mathbb{1},\chi'_{i}/\chi'_{i+1}) \hat T_{w,\beta,\chi,s\chi'}.
\]
\end{enumerate}
\end{enumerate}
\end{prop}

\begin{proof}
On the two periodic rows and columns, the calculation is the same as in the proof of Proposition~\ref{9.3.4}, using \eqref{eq:periodic-torus-lifts}. Note that multiplying the chosen Weyl lifts contributes $c(w,s)$ exactly to the terms with Weyl element $ws$. Thus, we obtain the formulas in the case (1), where, in (1.b), the removed index $i+1+kN$ is maximal in $s\beta$.

Now, we consider the case when $ws <_{I_{M}\ltimes U_{M-1,N}(\bF)} w$, and we should calculate the product $\underline{g}_{w(\beta)}\underline{w} \underline{s_{i}} \underline{e}_{-\check\alpha_i}(a)$, for $a\in \BG_a$.

If $a=0$, then the product $\underline{g}_{w(\beta)}\underline{w}\underline{s_i}$ contributes to the term $\hat T_{\tilde{w}s,\chi,s\chi'}$, provided that the corresponding orbit is relevant. Note that $(ws,s\beta)$ is not necessarily saturated, and we should modify it by saturation. Namely, if $(ws,s\beta)$ does not satisfy the conditions in Proposition~\ref{prop 3.2.7}, then we let $q=ws$ and replace $s\beta$ by
\[
\begin{split}
&\{i\in q^{-1}([1,M]_N):i\leq j,\ q(i)\leq q(j) \text{ for some }j\in s\beta\}\\
&\qquad\cup\{i\in q^{-1}([1,M]_N):i<q^{-1}(M+1)\}.
\end{split}
\]
This is the index of the same orbit in that proposition.

Assume $a\neq 0$. Using $\underline{h}_{\check\alpha_i}(a)\underline{h}_{-{\check\alpha_i}}(a)=1$ and the rank-one matrix identity, we have
\begin{equation}\label{10.57}
\begin{split}
\underline{g}_{w(\beta)}\underline{w} \underline{s_{i}} \underline{e}_{-\check\alpha_i}(a)\\
\simeq \underline{g}_{w(\beta)} \underline{e}_{w(-\check\alpha_i)}(a^{-1})\underline{w} \pi_{i+1} \underline{h}_{\check\alpha_i}(a)\underline{e}_{\check\alpha_i}(a^{-1})\\
\simeq \underline{e}_{w(-\check\alpha_i)}(a^{-1})\underline g(e_{w(\beta)}-a^{-1}e_{w(\beta_{i+1})})\cdot \underline{w} \pi_{i+1}\cdot  \underline{h}_{\check\alpha_i}(a)\underline{e}_{\check\alpha_i}(a^{-1}).
\end{split}
\end{equation}
Here, the second $\simeq$ allows left multiplication by an element of $U_{M,N}(\bF)\cap\ker(\psi)$, with its canonical lift. We let $\beta_{i+1}:=(\beta\cap [i]_N)+1$, which equals $\sigma_i+1-N\BN$ if $w(i)\in[1,M]_N$, and is empty otherwise.

A.) If $w(i)\notin [M+1]_N$, we have ${e}_{w(-{\check\alpha_i})}(a^{-1})\in I_M\ltimes U_{M,N}(\bF)$.

I.) Assume $w(i)\in[1,M]_N$ and $\sigma_i+1\in \beta\setminus \sigma$, or $w(i)\in [M+2,N]_N$. Then, the value of a left $(\BG_m,\kappa)\times (I_M,\chi)\ltimes (U_{M,N}(\bF),\psi)$-invariant and right $(I_N,s\chi')$-invariant function at the point $\underline{e}_{w(-{\check\alpha_i})}(a^{-1})\underline g(e_{w(\beta)}-a^{-1}e_{w(\beta_{i+1})})\cdot \underline{w} \pi_{i+1}\cdot  \underline{h}_{\check\alpha_i}(a)\underline{e}_{{\check\alpha_i}}(a^{-1})$ is the same as the value at the point $\underline{e}_{w(-{\check\alpha_i})}(a^{-1})\underline{g}_{w(\beta)}\cdot \underline{w} \pi_{i+1}\cdot  \underline{h}_{\check\alpha_i}(a)\underline{e}_{{\check\alpha_i}}(a^{-1})$.  Thus, if ${e}_{w(-{\check\alpha_i})}(a^{-1})$ is in the kernel of $\psi$, the coefficient of $\hat T_{\tilde{w},\chi,s\chi'}$ equals $\sum_{a\in \BF_Q\setminus\{0\}} (s\chi')_{i+1}(-1)\cdot (s\chi')_{i}(a^{-1})\cdot (s\chi')_{i+1}(a)$ if $i\neq 0$, and $\sum_{a\in \BF_Q\setminus\{0\}} (s\chi')_{i+1}(-1)\cdot \kappa^2(a)\cdot (s\chi')_{i+N}(a^{-1})\cdot (s\chi')_{i+1}(a)$ in the case $i=0$,  which also equals $\sum_{a\in \BF_Q\setminus\{0\}} (s\chi')_{i+1}(-1)\cdot (s\chi')_{i}(a^{-1})\cdot (s\chi')_{i+1}(a)$. This gives the formulas for (2.a) and (3.b). If ${e}_{w(-{\check\alpha_i})}(a^{-1})$ is not in the kernel of $\psi$, i.e., the case of (3.a), the coefficient of the term $\hat T_{\tilde{w},\chi,s\chi'}$ is $\sum_{a\in \BF_Q\setminus\{0\}}\psi(-a^{-1})(s\chi')_{i+1}(-1)(s\chi')_i(a^{-1})(s\chi')_{i+1}(a)=\chi'_{i+1}(-1)G(\psi,\chi'_i/\chi'_{i+1})$.

II.) If $\sigma_i+1\notin \beta$ and $w(i)\in [1,M]_N$, then every function $f$ with the equivariance conditions in I.) satisfies
\begin{multline}
f\bigl(\underline{e}_{w(-{\check\alpha_i})}(a^{-1})\underline g(e_{w(\beta)}-a^{-1}e_{w(\beta_{i+1})})\cdot \underline{w} \pi_{i+1}\cdot  \underline{h}_{\check\alpha_i}(a)\underline{e}_{{\check\alpha_i}}(a^{-1})\bigr)\\= f\bigl(\underline{e}_{w(-{\check\alpha_i})}(a^{-1})\pi_{w(i+1)}(-a^{-1})\underline{g}_{w(\beta\cup \{\sigma_i+1\})}\cdot \pi_{w(i+1)}(-a)\underline{w} \pi_{i+1}\underline{h}_{{\check\alpha_i}}(a)\underline{e}_{{\check\alpha_i}}(a^{-1})\bigr).
\end{multline}
Here, the matrix identity uses $t-a^{-1}=-a^{-1}(1-at)$ in place of $-a^{-1}$, and the extra factors $(1-at)^{\pm1}$ belong to the Iwahori congruence subgroup, on which the characters and the tame symbols are trivial.

From the above, we obtain that the coefficient of $\hat T_{w,\beta\cup \{\sigma_i+1\},\chi,s\chi'}$ is given by
\[
\begin{aligned}
&\sum_{a\in \BF_Q-\{0\}} \chi_{w(i+1)}(-a)(s\chi')_{i+1}(a^{-1}) (s\chi')_i(a^{-1})(s\chi')_{i+1}(a)\\
&\qquad=\chi_{{w}(i+1)}(-1) \sum_{a\in \BF_Q-\{0\}}(\chi_{w(i+1)}/\chi'_{i+1})(a).
\end{aligned}
\]

III.) Assume $\sigma_i+1\in\sigma$ and $w(i)\in[1,M]_N$. If $a=1$, it contributes to a term $\chi'_i(-1)\hat T_{w,\beta\setminus\{\sigma_i+1\},\chi,s\chi'}$. If this orbit is relevant, then $\chi_{w(i+1)}=\chi'_i$, and hence the coefficient is $\chi_{w(i+1)}(-1)$. If $a\neq0,1$, we have
\begin{multline}
\underline{e}_{w(-{\check\alpha_i})}(a^{-1})\underline g(e_{w(\beta)}-a^{-1}e_{w(\beta_{i+1})})\cdot \underline{w} \pi_{i+1}\cdot  \underline{h}_{\check\alpha_i}(a)\underline{e}_{{\check\alpha_i}}(a^{-1})\\= \underline{e}_{w(-{\check\alpha_i})}(a^{-1})\pi_{w(i+1)}(1-a^{-1})\underline{g}_{w(\beta\cup \{\sigma_i+1\})}\cdot \pi_{w(i+1)}((1-a^{-1})^{-1})\underline{w} \pi_{i+1}\underline{h}_{{\check\alpha_i}}(a)\underline{e}_{{\check\alpha_i}}(a^{-1}).
\end{multline}

From the above, we obtain that the coefficient of the term $\hat T_{\tilde w,\chi,s\chi'}$ is given by $\sum_{a\in \BF_Q\setminus\{0,1\}}\chi_{w(i+1)}((1-a^{-1})^{-1})(s\chi')_{i+1}(1-a^{-1})(s\chi')_{i+1}(-1)(s\chi')_{i}(a^{-1})(s\chi')_{i+1}(a)=\chi_{w(i+1)}(-1) J\!\left( \frac{\chi_{w(i+1)}}{\chi_{i+1}'}, \frac{\chi_i'}{\chi_{w(i+1)}}\right)$.

B.) Assume $w(i)\in [M+1]_N$.

I.) In the case of $\sigma_i+1\in \sigma$, if $a=-1$, then $\underline{g}_{w(\beta)} \underline{e}_{w(-\check\alpha_i)}(a^{-1})\underline{w} \pi_{i+1} \underline{h}_{\check\alpha_i}(a)\underline{e}_{{\check\alpha_i}}(a^{-1})$ gives $\chi'_{i+1}(-1)\hat T_{w,\beta\setminus\{\sigma_i+1\},\chi,s\chi'}$. If $a\neq0,-1$, the invariance conditions imply that the value at $\underline{g}_{w(\beta)} \underline{e}_{w(-\check\alpha_i)}(a^{-1})\underline{w} \pi_{i+1} \underline{h}_{\check\alpha_i}(a)\underline{e}_{{\check\alpha_i}}(a^{-1})$ is the same as the value at $ \pi_{w(i+1)}(1+a^{-1})\underline{g}_{w(\beta)}\cdot \pi_{w(i+1)}((1+a^{-1})^{-1})\underline{w} \pi_{i+1}\underline{h}_{{\check\alpha_i}}(a)\underline{e}_{{\check\alpha_i}}(a^{-1})$. Namely,
\begin{multline}
\underline{g}_{w(\beta)} \underline{e}_{w(-\check\alpha_i)}(a^{-1})\underline{w} \pi_{i+1} \underline{h}_{\check\alpha_i}(a)\underline{e}_{{\check\alpha_i}}(a^{-1})\\
=\pi_{w(i+1)}(1+a^{-1}+b)\underline{g}_{w(\beta)}\cdot \pi_{w(i+1)}((1+a^{-1}+b')^{-1})\underline{w} \pi_{i+1}\underline{h}_{{\check\alpha_i}}(a)\underline{e}_{{\check\alpha_i}}(a^{-1}),
\end{multline}
where $b,b'\in t\bO$.

Based on this fact, we obtain that the coefficient of $\hat T_{\tilde w,\chi,s\chi'}$ is given by $\sum_{a\in \BF_Q-\{0,-1\}}\chi_{w(i+1)}((1+a^{-1})^{-1})(s\chi')_{i+1}(-1-a^{-1})(s\chi')_i(a^{-1})(s\chi')_{i+1}(a)=\frac{\chi_{w(i+1)}\chi_i'}{\chi'_{i+1}}(-1) J\!\left( \frac{\chi_{w(i+1)}}{\chi_{i+1}'}, \frac{\chi_i'}{\chi_{w(i+1)}}\right)$.

II.) In the case $\sigma_i+1\notin \beta$, the invariance conditions imply that the value at $\underline{g}_{w(\beta)} \underline{e}_{w(-\check\alpha_i)}(a^{-1})\underline{w} \pi_{i+1} \underline{h}_{\check\alpha_i}(a)\underline{e}_{{\check\alpha_i}}(a^{-1})$ is the same as the value at $ \pi_{w(i+1)}(a^{-1})\underline{g}_{w(\beta\cup\{\sigma_i+1\})}\cdot \pi_{w(i+1)}(a)\underline{w} \pi_{i+1}\underline{h}_{{\check\alpha_i}}(a)\underline{e}_{{\check\alpha_i}}(a^{-1})$.

Thus, the coefficient of $\hat T_{w,\beta\cup\{\sigma_i+1\},\chi,s\chi'}$ is given by
\begin{multline*}
\sum_{a\in \BF_Q\setminus\{0\}}\chi_{w(i+1)}(a)(s\chi')_{i+1}(-a^{-1}) (s\chi')_i(a^{-1})(s\chi')_{i+1}(a)\\
=\chi'_{i}(-1)\sum_{a\in \BF_Q\setminus\{0\}}(\chi_{w(i+1)}/\chi'_{i+1})(a).
\end{multline*}

III.) In the case $\sigma_i+1\in \beta\setminus \sigma$, if $a\neq 0$, the invariance conditions imply that the value at $\underline{g}_{w(\beta)} \underline{e}_{w(-\check\alpha_i)}(a^{-1})\underline{w} \pi_{i+1} \underline{h}_{\check\alpha_i}(a)\underline{e}_{{\check\alpha_i}}(a^{-1})$ is the same as the value at $\underline{g}_{w(\beta)}\cdot \underline{w} \pi_{i+1}\underline{h}_{{\check\alpha_i}}(a)\underline{e}_{{\check\alpha_i}}(a^{-1})$. Thus, the coefficient of $\hat T_{w,\beta,\chi,s\chi'}$ is given by $\sum_{a\in \BF_Q\setminus\{0\}}(s\chi')_{i+1}(-1)(s\chi')_{i}(a^{-1})(s\chi')_{i+1}(a)=\chi'_{i}(-1)\sum_{a\in \BF_Q\setminus\{0\}}(\chi'_{i}/\chi'_{i+1})(a)=\chi'_{i}(-1)G(\mathbb{1},\chi'_{i}/\chi'_{i+1})$.
\end{proof}

\subsubsection{}
Since $\hat T_{g,\chi,\chi'}=\on{vol}(\Stab_{\widehat H}(gI_N))\cdot T_{g,\chi,\chi'}$, as in Corollary~\ref{9.3.4+}, we obtain the following formulas for the functions $T_{\tilde w,\chi,\chi'}$ from the above proposition.
\begin{prop}\label{right yoko}
We have the following module structure:

\begin{enumerate}
\item If $ws >_{I_{M}\ltimes U_{M-1,N}(\bF)} w$,
\begin{enumerate}
\item if there is no $k$ such that $\iota+kN \subseteq \sigma'\cup\{(ws)^{-1}(M+1)\}$, then
\[ T_{\tilde w,\chi,\chi'} T_{s,\chi',s\chi'} = c(w,s)\cdot T_{\tilde{w}s,\chi,s\chi'},\]
\item if there is $k$ such that $\iota+kN\subseteq \sigma'\cup\{(ws)^{-1}(M+1)\}$,
\[ T_{\tilde w,\chi,\chi'} T_{s,\chi',s\chi'}= c(w,s)\cdot T_{\tilde{w}s,\chi,s\chi'} + c(w,s)\cdot T_{ws,s\beta\setminus\{i+1+kN\},\chi,s\chi'},\]
\end{enumerate}
\item If $ws<_{I_{M}\ltimes U_{M-1,N}(\bF)} w$, and $w(i),w(i+1)\in[1,M]_N$,
\begin{enumerate}
\item If $\sigma_i+1\in \beta\setminus \sigma$ then
\[ T_{\tilde w,\chi,\chi'} T_{s,\chi',s\chi'} = Q\cdot c(w,s)\cdot T_{\tilde{w}s,\chi,s\chi'} + \chi'_i(-1)G(\mathbb{1},\chi'_{i}/\chi'_{i+1})\cdot T_{\tilde w,\chi,s\chi'}.
\]
\item If $\sigma_i+1\notin \beta$,
\[
T_{\tilde w,\chi,\chi'} T_{s,\chi',s\chi'} = c(w,s)\cdot T_{ws,s\beta\cup\{\sigma_i\},\chi,s\chi'} + \chi_{w(i+1)}(-1)\cdot T_{w,\beta\cup \{\sigma_i+1\},\chi,s\chi'}.
\]
\item If $\sigma_i+1\in \sigma$,
\begin{multline*}
T_{\tilde w,\chi,\chi'} T_{s,\chi',s\chi'} = (Q-1)\cdot c(w,s)\cdot T_{ws,s\beta,\chi,s\chi'} + (Q-1)\cdot\chi_{w(i+1)}(-1)\cdot T_{w,\beta\setminus\{\sigma_i+1\},\chi,s\chi'} \\ + \chi_{w(i+1)}(-1) J( {\chi_{w(i+1)}}/{\chi_{i+1}'}, {\chi_i'}/{\chi_{w(i+1)}}) \cdot T_{\tilde w,\chi,s\chi'}.
\end{multline*}
\end{enumerate}
\item If $ws<_{I_{M}\ltimes U_{M-1,N}(\bF)} w$, and $w(i)\in[M+2,N]_N$,
\begin{enumerate}
\item if $w(i+1)=w(i)-1$, then
\[
T_{\tilde w,\chi,\chi'} T_{s,\chi',s\chi'} = \chi'_{i+1}(-1)G(\psi,\chi'_i/\chi'_{i+1}) T_{\tilde w,\chi,s\chi'};\]
\item otherwise
\begin{multline*}
T_{\tilde w,\chi,\chi'} T_{s,\chi',s\chi'} = Q\cdot c(w,s)\cdot T_{ws,s\beta,\chi,s\chi'} + \chi'_{i+1}(-1)G(\mathbb{1},\chi'_i/\chi'_{i+1}) T_{\tilde{w},\chi,s\chi'}
\end{multline*}
\end{enumerate}

\item If $ws<_{I_{M}\ltimes U_{M-1,N}(\bF)} w$, and $w(i)\in [M+1]_N$,
\begin{enumerate}
\item if $\sigma_i+1=w^{-1}(M+1)+1\in \sigma$,
\begin{multline*}
T_{\tilde w,\chi,\chi'} T_{s,\chi',s\chi'} = (Q-1)\cdot c(w,s)\cdot T_{ws,s\beta,\chi,s\chi'} + (Q-1)\cdot \chi'_{i+1}(-1) T_{w,\beta\setminus\{\sigma_i+1\},\chi,s\chi'} \\+ \frac{\chi_{w(i+1)}\chi_i'}{\chi'_{i+1}}(-1) J( {\chi_{w(i+1)}}/{\chi_{i+1}'}, {\chi_i'}/{\chi_{w(i+1)}}) \cdot T_{\tilde w,\chi,s\chi'}
\end{multline*}
\item if $w^{-1}(M+1)+1\notin \beta$
\begin{multline*}
T_{\tilde w,\chi,\chi'} T_{s,\chi',s\chi'} = c(w,s)\cdot  T_{ws,s\beta\cup\{\sigma_i\},\chi,s\chi'} + \chi'_{i}(-1) T_{w,\beta\cup\{\sigma_i+1\},\chi,s\chi'}.
\end{multline*}
\item if $w^{-1}(M+1)+1\in \beta\setminus \sigma$,
\[
T_{\tilde w,\chi,\chi'} T_{s,\chi',s\chi'} = Q\cdot c(w,s)\cdot T_{ws,s\beta,\chi,s\chi'} + \chi'_{i}(-1)G(\mathbb{1},\chi'_{i}/\chi'_{i+1}) T_{w,\beta,\chi,s\chi'}.
\]
\end{enumerate}
\end{enumerate}
\end{prop}

\subsubsection{The left action} Similarly, consider the left action of the twisted affine Yokonuma algebra. Below, let us denote by $j^+$ the successor of $j$ in $[1,M]_N$, and let $\iota=\{i,i^+\}$, with $i\in\{1,\ldots,M-1\}$ for the finite simple reflections and $i=M-N$, $i^+=1$ for the affine simple reflection. We denote by $\check\alpha$ the corresponding geometric affine root, and set $\pi_i'=w_{\check\alpha}^{-1}\pi_{i^+}w_{\check\alpha}$. Then, in the calculation of the left action, \eqref{10.57} is replaced by
\[
\underline e_{-\check\alpha}(a)\pi_{i^+}w_{\check\alpha} \underline g_{w(\beta)}\underline w =\underline e_{\check\alpha}(a^{-1}) \underline h_{\check\alpha}(a^{-1})\pi_i' \underline e_{-\check\alpha}(a^{-1}) \underline g_{w(\beta)}\underline w,
\]
for $a\neq0$, where $\pi_i'=\pi_i$ by \eqref{eq:periodic-torus-lifts}.

For $\tilde w=(w,\beta)$, we denote $s\tilde w:=(sw,\beta)$.

\begin{prop}\label{left yoko}
Let $\sigma'':=\sigma(s\tilde{w})$, and let $w(\sigma)_i$ be the maximal element in $[i]_N\cap w(\beta)$. Then, we have the following left action formulas.

\begin{enumerate}
\item If $sw >_{I_{M}\ltimes U_{M-1,N}(\bF)} w$,
\begin{enumerate}
\item if there is no $k$ such that $\iota+kN \subseteq sw\sigma''$, then
\[ T_{s,\chi} T_{\tilde w,\chi,\chi'} = c(s,w)\cdot T_{s\tilde w,s\chi,\chi'},\]
\item if there is $k$ such that $\iota+kN\subseteq sw\sigma''$,
\[  T_{s,\chi} T_{\tilde w,\chi,\chi'}=  c(s,w)\cdot T_{s\tilde w,s\chi,\chi'} + c(s,w)\cdot T_{sw,\beta\setminus\{w^{-1}(i+kN)\},s\chi,\chi'},\]
\end{enumerate}
\item If $sw<_{I_{M}\ltimes U_{M-1,N}(\bF)} w$,
\begin{enumerate}
\item If $(w(\sigma)_i)^+\in w(\beta\setminus\sigma)$ then
\[T_{s,\chi} T_{\tilde w,\chi,\chi'} = Q\cdot c(s,w)\cdot T_{s\tilde{w},s\chi,\chi'} + \chi_{i^+}(-1)G(\mathbb{1},\chi_{i^+}/\chi_{i})\cdot T_{\tilde w,s\chi,\chi'}.
\]
\item If $(w(\sigma)_i)^+\notin w(\beta)$,
\[
T_{s,\chi} T_{\tilde w,\chi,\chi'} = c(s,w)\cdot T_{{sw},\beta\cup \{w^{-1}((w(\sigma)_i)^+)\}, s\chi,\chi'} + \chi'_{w^{-1}(i^+)}(-1)\cdot T_{w,\beta\cup \{w^{-1}((w(\sigma)_i)^+)\},s\chi,\chi'}.
\]
\item If $(w(\sigma)_i)^+\in w(\sigma)$,
\begin{multline*}
T_{s,\chi} T_{\tilde w,\chi,\chi'} = (Q-1)\cdot c(s,w)\cdot T_{sw,\beta,s\chi,\chi'} + (Q-1)\cdot\chi'_{w^{-1}(i^+)}(-1)\cdot T_{w,\beta\setminus\{w^{-1}((w(\sigma)_i)^+)\},s\chi,\chi'} \\ + \chi'_{w^{-1}(i^+)}(-1) J( {\chi'_{w^{-1}(i^+)}}/{\chi_{i}}, {\chi_{i^+}}/{\chi'_{w^{-1}(i^+)}}) \cdot T_{\tilde w,s\chi,\chi'}.
\end{multline*}
\end{enumerate}
\end{enumerate}
\end{prop}
\begin{proof}
The calculation is the same as in the proof of Proposition~\ref{9.3.4}, with rows and columns interchanged, using the preceding identity and \eqref{eq:periodic-torus-lifts}. Note that multiplying the chosen lifts contributes $c(s,w)$ to the terms with Weyl element $sw$.
\end{proof}
\subsubsection{The mixed affine Yokonuma algebra} We denote by $\cH_{\kappa,\chi,\chi'}^{\aff,\mxd}(N)$ the category of $\kappa$-twisted left $(I_N,\chi)$-equivariant and right $(I_N,\chi')$-equivariant mixed sheaves on $\GL_N(\bF)$, and by $\bH_{\kappa,\chi,\chi'}^{\aff,\mxd}(N)$ its Grothendieck group. The convolution induces the multiplicative structures
\begin{equation}\label{yoko M}
\bH^{\aff,\mxd}_{\kappa,\chi',\chi}(M)\otimes \bH^{\aff,\mxd}_{\kappa,\chi,\chi''}(M)\rightarrow \bH^{\aff,\mxd}_{\kappa,\chi',\chi''}(M),
\end{equation}
and
\begin{equation}\label{yoko N}
\bH^{\aff,\mxd}_{\kappa,\chi',\chi}(N)\otimes \bH^{\aff,\mxd}_{\kappa,\chi,\chi''}(N)\rightarrow \bH^{\aff,\mxd}_{\kappa,\chi',\chi''}(N).
\end{equation}

We denote by $H^\mxd_{w,\chi,\chi'}$ the standard basis element in the mixed affine Yokonuma algebra, corresponding to the $!$-extension of the $\kappa$-twisted $(I_N,\chi)$-left and $(I_N,\chi')$-right equivariant perverse pure of weight $0$ local system supported on the orbit indexed by $w$, i.e., the usual $!$-extension of the local system shifted by $\langle l(w)\rangle$. Namely, we have\begin{equation} H^\mxd_{w,\chi,\chi'} = (-[v])^{-l(w)} \cdot T^\mxd_{w,\chi,\chi'}.
\end{equation}
Its specialization to $\bH_{\kappa,\chi,\chi'}^{\aff}(N)_\BC$ under $\on{sp}$ corresponds to $H_{w,\chi,\chi'}=(-v)^{-l(w)}T_{w,\chi,\chi'}$.

The Verdier duality functor induces a bar involution
\begin{equation}
\bar{\cdot}\colon \bH^{\aff,\mxd}_{\kappa,\chi,\chi'}(N)\rightarrow \bH^{\aff,\mxd}_{(\kappa)^{-1},\chi^{-1},(\chi')^{-1}}(N),
\end{equation}
which is uniquely characterized by the requirement that $\bar{[a]}=[a^{-1}]\in \bR^\mxd$, and

\begin{equation}
\bar{H}^{\mxd}_{w,\chi,\chi'} H^\mxd_{w^{-1},(\chi')^{-1},\chi^{-1}} =[c(w,w^{-1})^{-1}]\,H^\mxd_{e,\chi^{-1},\chi^{-1}}.
\end{equation}
Here, $H^\mxd_{w^{-1},(\chi')^{-1},\chi^{-1}}\in \bH^{\aff,\mxd}_{(\kappa)^{-1},(\chi')^{-1},\chi^{-1}}(N)$, i.e., the twisting is $\kappa^{-1}$, and thus, the multiplication formula uses $c(w,w^{-1})^{-1}$.

As in the finite type case, the Grothendieck groups for inverse twisting and characters are identified by their normalized standard bases. Consider the pure of weight $0$ IC sheaf on the orbit corresponding to $w$. The canonical basis $\underline{H}^{\mxd}_{w,\chi,\chi'}$ and its graded version after taking the specialization $\on{wt}$, i.e., the corresponding elements of the pure of weight $0$ IC sheaf in the Grothendieck group, are characterized by requirements similar to those in \eqref{bar h mxd} and \eqref{bar h gr}.

For example, for a simple reflection $s=(i,i+1)\in W_N^{\mathrm{ext}}$ with $\chi=s\chi'$, the canonical basis element in the mixed Grothendieck group is given by
\[
\underline{H}^\mxd_{s,\chi,\chi'} =
\begin{cases}
H^\mxd_{s,\chi,\chi'}-[v^{-1}][\chi_i(-1)]\cdot H^\mxd_{e,\chi,\chi'}, & \text{if } \chi_i=\chi_{i+1},\\[4pt]
H^\mxd_{s,\chi,\chi'}, & \text{if } \chi_i\neq \chi_{i+1}.
\end{cases}
\]

For the affine simple reflection, we use the convention of Section~\ref{sec:periodic-characters}.

By taking the specialization with respect to $\on{wt}$, we obtain
\begin{equation}
\bar{\cdot}\colon \bH^{\aff,\gr}_{\kappa,\chi,\chi'}(N)\rightarrow \bH^{\aff,\gr}_{(\kappa)^{-1},\chi^{-1},(\chi')^{-1}}(N),
\end{equation}
and the canonical basis $\underline{H}^\gr_{w,\chi,\chi'}$.
\begin{rem}
It does not matter whether we put the superscripts ``aff'' or ``fin'' for the basis here, since the basis that is defined in the finite setting matches the affine one by the canonical embedding of the finite Hecke/Yokonuma algebra into the affine Hecke/Yokonuma algebra.
\end{rem}

\subsubsection{The mixed action}

For the same reason as in Section~\ref{sec 10.2.9}, by considering the action formulas over all $\BF_{Q^n}$, one recovers the action formula
\begin{equation}\label{10.61}
\bH_{\kappa,\chi,\chi'}^{\aff,\mxd}(M)\otimes \bK_{\kappa,\chi',\chi''}^{\aff,\mxd}(M|N) \longrightarrow\bK_{\kappa,\chi,\chi''}^{\aff,\mxd}(M|N),
\end{equation}
from Propositions~\ref{right yoko} and \ref{left yoko}. Here, the formula is given by the convolution functor of the category of mixed sheaves
\begin{equation}
\cH_{\kappa,\chi,\chi'}^{\aff,\mxd}(M)\otimes \cC_{\kappa,\chi',\chi''}^{\aff,\mxd}(M|N) \longrightarrow\cC_{\kappa,\chi,\chi''}^{\aff,\mxd}(M|N),
\end{equation}
under the characteristic map. The right action is recovered in the same way.

To be more precise, to recover the formula \eqref{10.61}, we just need to replace the Jacobi sums and Gauss sums in the original formulas (i.e., on the level of $\BC$-valued functions) by the classes of the corresponding compactly supported cohomology complexes in $\bR^\mxd$.

For characters $\chi, \chi'$, we replace $J(\chi,\chi')$ by

\begin{equation}\label{9.53}
J(\chi,\chi')\mapsto  \begin{cases}        [Q]-2\cdot[1] \text{, if $\chi=\chi'=\id$;}\\
-[1] \text{, if exactly one of $\chi$ and $\chi'$ is $\id$;}\\
-[\chi(-1)]\text{, if $\chi\chi'=\id$ and $\chi,\chi'\neq \id$;}\\
-[-J(\chi,\chi')]\text{, if none of $\chi,\chi', \chi\chi'$ is $ \id$.}
\end{cases}
\end{equation}

For the Gauss sum, we replace $G(\psi, \chi)$ by
\begin{equation}\label{9.54}
G(\psi, \chi)\mapsto
\begin{cases}
-[1]\text{, if $\psi\neq \id$, $\chi=\id$;}\\
-[-G(\psi,\chi)]\text{, if $\psi\neq \id$, $\chi\neq\id$.}\\
[Q]-[1]\text{, if $\psi=\chi=\id$;}\\
0\text{, if $\chi\neq\id,\psi=\id$;}\\
\end{cases}
\end{equation}

\subsubsection{Weight specialization} We apply the specialization $\on{wt}$, and just remember the grading by weights. We denote by $\bJ, \bG\in \BZ[\bv,\bv^{-1}]$ the images of the Jacobi sum and Gauss sum in \eqref{9.53} and \eqref{9.54} under the map \eqref{9.67 specialization}.

Applying the above map to the formulas \eqref{yoko M}, \eqref{yoko N}, and \eqref{10.61}, we rewrite the above propositions using the graded standard basis $H=(-\bv)^{-l(\tilde{w})}T$, and only remember the coefficients in $\BZ[\bv,\bv^{-1}]$, which contain the weight information.
\begin{cor}
\begin{enumerate}
\item If $l(w_1)+l(w_2)=l(w_1w_2)$, then
\begin{equation}
H^\gr_{w_1,w_1(w_2(\chi)),w_2(\chi)}H^\gr_{w_2, w_2(\chi),\chi}=H^\gr_{w_1w_2, w_1(w_2(\chi)), \chi};
\end{equation}
\item
\begin{equation}
H^\gr_{s, \chi,s\chi}H^\gr_{s,s\chi,\chi}= H^\gr_{e,\chi,\chi}+(\bv^{-1}-\bv)\cdot H^\gr_{s,\chi,\chi}.
\end{equation}
Note that, if $\chi_i\neq\chi_{i+1}$, then we have $T^\gr_{s,\chi,\chi}=H^\gr_{s,\chi,\chi}=0$.
\end{enumerate}
\end{cor}

More explicitly, if $T_{\tilde w}T_s=\sum_{\tilde v}a_{\tilde v}T_{\tilde v}$, then the coefficient of $H^{\gr}_{\tilde v}$ in $H^{\gr}_{\tilde w}H^{\gr}_s$ is
\[
\operatorname{wt}(a_{\tilde v}) (-\bv)^{l(\tilde v)-l(\tilde w)-1}.
\]
The same rule applies to the left action. Here, the character values and $c(w_1,w_2)$ have absolute value $1$, thus, their classes are sent to $1$ under $\on{wt}$.
\begin{cor}\label{right yoko wt}
We have the following module structure:
\begin{enumerate}
\item If $ws >_{I_{M}\ltimes U_{M-1,N}(\bF)} w$,
\begin{enumerate}
\item if there is no $k$ such that $\iota+kN \subseteq \sigma'\cup\{(ws)^{-1}(M+1)\}$, then
\[ H^\gr_{\tilde w,\chi,\chi'} H^\gr_{s,\chi',s\chi'} = H^\gr_{\tilde{w}s,\chi,s\chi'},\]
\item if there is $k$ such that $\iota+kN\subseteq\sigma'\cup\{(ws)^{-1}(M+1)\}$,
\[ H^\gr_{\tilde w,\chi,\chi'} H^\gr_{s,\chi',s\chi'}= H^\gr_{\tilde{w}s,\chi,s\chi'} -\bv^{-1}\cdot H^\gr_{ws,s\beta\setminus\{i+1+kN\},\chi,s\chi'},\]
\end{enumerate}
\item If $ws<_{I_{M}\ltimes U_{M-1,N}(\bF)} w$, and $w(i),w(i+1)\in[1,M]_N$,
\begin{enumerate}
\item If $\sigma_i+1\in \beta\setminus\sigma$ then
\[ H^\gr_{\tilde w,\chi,\chi'} H^\gr_{s,\chi',s\chi'} = H^\gr_{\tilde{w}s,\chi,s\chi'} -\bv^{-1}\cdot \bG(\mathbb{1},\chi'_{i}/\chi'_{i+1})\cdot H^\gr_{\tilde w,\chi,s\chi'}.
\]
\item If $\sigma_i+1\notin \beta$,
\[
H^\gr_{\tilde w,\chi,\chi'} H^\gr_{s,\chi',s\chi'} = -\bv^{-1}\cdot H^\gr_{ws,s\beta\cup\{\sigma_i\},\chi,s\chi'} + H^\gr_{w,\beta\cup \{\sigma_i+1\},\chi,s\chi'}.
\]
\item If $\sigma_i+1\in \sigma$,
\begin{multline*}
H^\gr_{\tilde w,\chi,\chi'} H^\gr_{s,\chi',s\chi'} = (1-\bv^{-2})\cdot H^\gr_{ws,s\beta,\chi,s\chi'} + (1-\bv^{-2})\cdot H^\gr_{w,\beta\setminus\{\sigma_i+1\},\chi,s\chi'}\\
-\bv^{-1} \bJ( {\chi_{w(i+1)}}/{\chi_{i+1}'}, {\chi_i'}/{\chi_{w(i+1)}}) \cdot H^\gr_{\tilde w,\chi,s\chi'}.
\end{multline*}
\end{enumerate}
\item If $ws<_{I_{M}\ltimes U_{M-1,N}(\bF)} w$, and $w(i)\in[M+2,N]_N$,
\begin{enumerate}
\item if $w(i+1)=w(i)-1$, then
\[
H^\gr_{\tilde w,\chi,\chi'} H^\gr_{s,\chi',s\chi'} = -\bv^{-1}\bG(\psi,\chi'_i/\chi'_{i+1}) H^\gr_{\tilde w,\chi,s\chi'};\]
\item otherwise,
\begin{equation*}
H^\gr_{\tilde w,\chi,\chi'} H^\gr_{s,\chi',s\chi'}= H^\gr_{ws,s\beta,\chi,s\chi'} -\bv^{-1}\cdot \bG(\mathbb{1},\chi'_i/\chi'_{i+1})\cdot H^\gr_{\tilde{w},\chi,s\chi'}.
\end{equation*}
\end{enumerate}

\item If $ws<_{I_{M}\ltimes U_{M-1,N}(\bF)} w$, and $w(i)\in [M+1]_N$,
\begin{enumerate}
\item if $\sigma_i+1=w^{-1}(M+1)+1\in \sigma$,
\begin{multline*}
H^\gr_{\tilde w,\chi,\chi'} H^\gr_{s,\chi',s\chi'} = (1-\bv^{-2})\cdot  H^\gr_{ws,s\beta,\chi,s\chi'} + (1-\bv^{-2})\cdot  H^\gr_{w,\beta\setminus\{\sigma_i+1\},\chi,s\chi'}\\-\bv^{-1} \bJ({\chi_{w(i+1)}}/{\chi_{i+1}'}, {\chi_i'}/{\chi_{w(i+1)}}) \cdot H^\gr_{\tilde w,\chi,s\chi'}
\end{multline*}
\item if $w^{-1}(M+1)+1\notin \beta$
\[
H^\gr_{\tilde w,\chi,\chi'} H^\gr_{s,\chi',s\chi'} = -\bv^{-1} H^\gr_{ws,s\beta\cup\{\sigma_i\},\chi,s\chi'} +H^\gr_{w,\beta\cup\{\sigma_i+1\},\chi,s\chi'}.\]
\item if $w^{-1}(M+1)+1\in \beta\setminus \sigma$,
\[
H^\gr_{\tilde w,\chi,\chi'} H^\gr_{s,\chi',s\chi'} = H^\gr_{ws,s\beta,\chi,s\chi'} -\bv^{-1}\bG(\mathbb{1},\chi'_{i}/\chi'_{i+1}) H^\gr_{w,\beta,\chi,s\chi'}.
\]
\end{enumerate}
\end{enumerate}
\end{cor}

\ifdim\dimexpr\pagegoal-\pagetotal\relax<7\baselineskip\newpage\fi
\begin{cor}\label{yoko left wt}
We have the following module structure:

\begin{enumerate}
\item If $sw >_{I_{M}\ltimes U_{M-1,N}(\bF)} w$,
\begin{enumerate}
\item if there is no $k$ such that $\iota+kN \subseteq sw\sigma''$, then
\[H^\gr_{s,s\chi,\chi} H^\gr_{\tilde w,\chi,\chi'} = H^\gr_{s\tilde w,s\chi,\chi'},\]
\item if there is $k$ such that $\iota+kN\subseteq sw\sigma''$,
\[  H^\gr_{s,s\chi,\chi} H^\gr_{\tilde w,\chi,\chi'}=  H^\gr_{s\tilde w,s\chi,\chi'} -\bv^{-1} H^\gr_{sw,\beta\setminus\{w^{-1}(i+kN)\},s\chi,\chi'},\]
\end{enumerate}
\item If $sw<_{I_{M}\ltimes U_{M-1,N}(\bF)} w$,
\begin{enumerate}
\item If $(w(\sigma)_i)^+\in w(\beta\setminus\sigma)$ then
\[H^\gr_{s,s\chi,\chi} H^\gr_{\tilde w,\chi,\chi'} = H^\gr_{s\tilde w,s\chi,\chi'} -\bv^{-1}\bG(\mathbb{1},\chi_{i^+}/\chi_{i})\cdot H^\gr_{\tilde w,s\chi,\chi'}.
\]
\item If $(w(\sigma)_i)^+\notin w(\beta)$,
\[
H^\gr_{s,s\chi,\chi} H^\gr_{\tilde w,\chi,\chi'} = -\bv^{-1} H^\gr_{sw,\beta\cup\{w^{-1}((w(\sigma)_i)^+)\},s\chi,\chi'} + H^\gr_{w,\beta\cup \{w^{-1}((w(\sigma)_i)^+)\},s\chi,\chi'}.
\]
\item If $(w(\sigma)_i)^+\in w(\sigma)$,
\begin{multline*}
H^\gr_{s,s\chi,\chi} H^\gr_{\tilde w,\chi,\chi'} = (1-\bv^{-2})\cdot H^\gr_{sw,\beta,s\chi,\chi'} +(1-\bv^{-2})\cdot H^\gr_{w,\beta\setminus\{w^{-1}((w(\sigma)_i)^+)\},s\chi,\chi'} \\-\bv^{-1} \bJ( {\chi'_{w^{-1}(i^+)}}/{\chi_{i}}, {\chi_{i^+}}/{\chi'_{w^{-1}(i^+)}}) \cdot H^\gr_{\tilde w,s\chi,\chi'}.
\end{multline*}
\end{enumerate}
\end{enumerate}
\end{cor}

\ifdim\dimexpr\pagegoal-\pagetotal\relax<7\baselineskip\newpage\fi In particular,\par\nopagebreak[4]
\begin{cor}
\begin{enumerate}
\item  If $ws >_{I_{M}\ltimes U_{M-1,N}(\bF)} w$, and there is no $k$ such that $\iota+kN \subseteq \sigma'\cup\{(ws)^{-1}(M+1)\}$, then
\[H^\gr_{\tilde w,\chi,\chi'} \underline{H}^\gr_{s,\chi',s\chi'} = \begin{cases} H^\gr_{\tilde{w}s,\chi,s\chi'}, & \text{if } \chi'_i\neq\chi'_{i+1},\\
H^\gr_{\tilde{w}s,\chi,s\chi'}-\bv^{-1}H^\gr_{\tilde{w},\chi,s\chi'}, & \text{if } \chi'_i=\chi'_{i+1}.
\end{cases}\]
\item If $sw >_{I_{M}\ltimes U_{M-1,N}(\bF)} w$, and there is no $k$ such that $\iota+kN \subseteq sw\sigma''$, then
\[\underline{H}^\gr_{s,s\chi,\chi} H^\gr_{\tilde w,\chi,\chi'} = \begin{cases} H^\gr_{s\tilde w,s\chi,\chi'}, & \text{if } \chi_i\neq\chi_{i^+},\\
H^\gr_{s\tilde w,s\chi,\chi'}-\bv^{-1} H^\gr_{\tilde w,s\chi,\chi'}, & \text{if } \chi_i=\chi_{i^+}.
\end{cases}\]
\end{enumerate}
\end{cor}

\subsection{Generating property} For $M<N$, we denote by $RB^\aff(M|N)$ the set of indices $(w,\beta)$ in Proposition~\ref{prop 3.2.7} satisfying $w^{-1}(M+1)>\cdots>w^{-1}(N)$, i.e., the $\psi$-relevant indices. For fixed twisting and multiplicative characters, by Proposition~\ref{prop 3.2.10}(2), the relevance also requires that $\chi$ and $w(\chi')$ have the same restriction to $T_{w,\sigma}$. The idea of this subsection is to reduce every relevant orbit, by convolutions whose convolution diagrams are isomorphisms, to a reduced one, which comes from the finite type category up to translations, see Proposition~\ref{prop 9.4.17}. Now, let us consider a partial order on $RB^\aff(M|N)$.
\begin{defn}
We denote by $\preccurlyeq$ the partial order generated by:
\begin{enumerate}
\item $\tilde{w}\preccurlyeq s\tilde{w}$, if the corresponding convolution diagram
\[I_M s I_M\overset{I_M}{\times}(I_M\ltimes U_{M,N}(\bF))\cdot g_{\tilde{w}}I_N\rightarrow \GL_N(\bF)\] is an isomorphism to the orbit corresponding to $s\tilde{w}$. It is the case (1.a) in Corollary~\ref{yoko left wt}, i.e., we have $sw>_{I_M\ltimes U_{M-1,N}(\bF)}w$, and there is no $k$ such that $\iota+kN\subseteq sw\sigma''$;
\item $\tilde{w}\preccurlyeq \tilde{w}s$, if the corresponding convolution diagram
\[(I_M\ltimes U_{M,N}(\bF))\cdot g_{\tilde{w}}I_N\overset{I_N}{\times}I_NsI_N\rightarrow \GL_N(\bF)\] is an isomorphism to the orbit corresponding to $\tilde{w}s$. It is the case (1.a) in Corollary~\ref{right yoko wt}, i.e., we have $ws>_{I_M\ltimes U_{M-1,N}(\bF)}w$, and there is no $k$ such that $\iota+kN\subseteq\sigma'\cup\{(ws)^{-1}(M+1)\}$.
\end{enumerate}
\end{defn}

We have $\tilde{w}\preccurlyeq \tilde{w}'$ if and only if there exist $w_M\in W^\aff_M$ and $w_N\in W^\aff_N$, such that the convolution diagram corresponding to $T_{w_M} T_{\tilde{w}} T_{w_N}$ is an isomorphism to the support of $\tilde{w}'$. In particular, we have $|\sigma(\tilde w)|=|\sigma(\tilde w')|$, an isomorphism \[\BO_{\tilde{w}'}\simeq \BA^{l(w_M)}\times \BO_{\tilde{w}}\times \BA^{l(w_N)},\] and \[H^\gr_{w_M,w_M(\chi),\chi} H^\gr_{\tilde w,\chi,\chi'} H^\gr_{w_N,\chi',w_N^{-1}(\chi')}=H^\gr_{\tilde w',w_M(\chi),w_N^{-1}(\chi')},\] for any $\chi,\chi'$.

\begin{lem}\label{lem:affine-reduction-termination}
For $\tilde w=(w,\beta)\in RB^\aff(M|N)$, let $d=N-M$. For $1\leq r\leq M$, let $r+Nm_r$ be the largest element of $w(\beta)$ congruent to $r$ modulo $N$, and let $S(\tilde w)=\sum_{r=1}^M m_r$. Then, $S$ is constant along the relations defining $\preccurlyeq$, and
\[
l(\tilde w)+(d-1)S(\tilde w)\geq 0.
\]
In particular, every strictly descending chain with respect to $\preccurlyeq$ is finite.
\end{lem}
\begin{proof}
For the second relation defining $\preccurlyeq$, the indices change from $(w,\beta)$ to $(ws,s\beta)$, and thus, $(ws)(s\beta)=w(\beta)$. For the first relation, the left finite simple reflections permute the $m_r$, and the left affine simple reflection sends $(m_1,m_M)$ to $(m_M+1,m_1-1)$ when $M\geq2$. Thus, $S$ has the same value before and after each of these operations. Let $p_a=w^{-1}(a)$ and
\[
C=\#\{j\in\beta:j>p_{M+1}\} =\sum_{r=1}^M\left(m_r+\left\lceil\frac{p_r-p_{M+1}}N\right\rceil\right).
\]
By the formula for the relative dimension in Lemma~\ref{lem:affine-orbit-relative-dimension}, we have
\begin{equation}\label{eq:affine-reduction-length-bound}
\begin{split}
l(\tilde w)={}&L_M(w)+C +\sum_{r=1}^M\sum_{a=M+2}^N \left\lceil\frac{p_r-p_a}{N}\right\rceil\\
&+\sum_{M+1\leq a<b\leq N} \left\lceil\frac{p_a-p_b}{N}\right\rceil .
\end{split}
\end{equation}
Here
\[
L_M(w):=\dim\bigl(I_M/(I_M\cap wI_Nw^{-1})\bigr)
\]
is the dimension of the $I_M$-orbit through $wI_N\in\Fl_N$, and hence it is non-negative. Since $p_{M+1}>\cdots>p_N$ by the $\psi$-relevance, we obtain
\[
l(\tilde w)\geq C+(d-1)(C-S) =dC-(d-1)S\geq-(d-1)S.
\]
Finally, along a strictly descending chain, $S$ is constant and $l$ strictly decreases. Since $l\geq -(d-1)S$, the chain is finite.
\end{proof}

\subsubsection{Reduced subset of $RB^\aff(M|N)$}\label{sec:reduced-affine-indices} We consider the following subset of $RB^{\aff}(M|N)$:
\begin{enumerate}
\item if $i,i+1\in w^{-1}([1,M]_N)$, $w(i)>w(i+1)$, then there exists $k\in \BZ$ such that $i+kN\in \sigma$, and $i+1+kN\in \sigma$ or $\notin \beta$;
\item if $w(i)=M+1$, $w(i+1)\in [1,M]_N$, then $i+1\in \sigma$ or $\notin \beta$;
\item if $i\in [1,M-1]_N$, and $w^{-1}(i)>w^{-1}(i+1)$, then there exists $k\in \BZ$ such that $w^{-1}(i)+kN\in \sigma$, and $w^{-1}(i+1)+kN\in \sigma$ or $\notin \beta$;
\item if $w^{-1}(M)>w^{-1}(N+1)$, then there exists $k\in \BZ$ such that $w^{-1}(M)+kN\in \sigma$ and $w^{-1}(N+1)+kN\in \sigma$ or $\notin \beta$;
\item $w^{-1}(N)=w^{-1}(N-1)-1=w^{-1}(N-2)-2=\cdots =w^{-1}(M+1)-N+M+1$.
\end{enumerate}
We call it the reduced subset of $RB^\aff(M|N)$, and denote it by $RB^\aff_{\on{red}}(M|N)$. Note that the convolution isomorphisms transport the twisted equivariant local systems. The following proposition says that every element can be reduced to this subset.
\begin{prop}\label{prop:affine-reduction-to-reduced}
Any minimal element in $RB^\aff(M|N)$, with respect to $\preccurlyeq$, belongs to $RB^\aff_{\on{red}}(M|N)$, and for every $\tilde w\in RB^\aff(M|N)$, there exists $\tilde v\in RB^\aff_{\on{red}}(M|N)$ such that $\tilde v\preccurlyeq\tilde w$.

In particular, the bar involution on the direct sum of $\bK^{\aff,\mxd}_{\kappa,\chi,\chi'}(M|N)$ over the twistings, the characters and the additive characters $\psi,-\psi$ is determined by its values on the standard classes indexed by $RB^\aff_{\on{red}}(M|N)$ and its compatibility with the Hecke actions. The canonical basis is then determined by the triangularity condition in \eqref{bar h mxd}, the normalization of the leading coefficient in \eqref{eq:finite-IC-normalization}, and compatibility with Verdier duality. Here, duality sends each basis element to the corresponding element for the inverse twisting and characters.
\end{prop}
\begin{proof}
By Lemma~\ref{lem:affine-reduction-termination}, the successive reductions reach a minimal element in finitely many steps. Thus, it is enough to verify the five conditions for an element $\tilde w\in RB^\aff(M|N)$ which is minimal with respect to $\preccurlyeq$.

First, we prove (5). Condition (5) is automatic when $N=M+1$. In proving it below, we assume $N\geq M+2$. If $i\in w^{-1}([M+2,N]_N)$ and $i+1\in w^{-1}([1,M]_N)$, then we switch the $i$-th and $(i+1)$-th columns. Since $\tilde{w}s=(ws, s\beta)$ still satisfies Proposition~\ref{prop 3.2.7} and the $\psi$-relevance condition, and the convolution diagram associated with $H_{\tilde{w}s,\chi,s\chi'}H_{s, s\chi',\chi'}$ is an isomorphism to the support of $H_{\tilde w,\chi,\chi'}$, we obtain $\tilde{w}s\preccurlyeq \tilde{w}$. Thus, for a minimal element $\tilde{w}\in RB^\aff(M|N)$, the congruence classes $w^{-1}([M+1,N]_N)$ are adjacent to each other. Then, the situation is similar to the Whittaker case, and we claim that any such minimal element satisfies the condition (5).

Indeed, for the $N$-periodic interval $w^{-1}([M+1,N]_N)$, we note that if $i\in w^{-1}([M+2,N]_N)$, then $i+1\notin w^{-1}([1,M]_N)$ according to the minimality assumption, and thus, $w^{-1}([M+1]_N)$ is the rightmost congruence class of $w^{-1}([M+1,N]_N)$. Therefore, $w^{-1}(M+1)-1\in w^{-1}([M+2,N]_N)$. According to the $\psi$-relevance condition, we have \[w(w^{-1}(M+1)-1)>M+1.\]

Furthermore, according to the minimality assumption, \[w(w^{-1}(M+1)-1)<M+1+N.\] Otherwise, we switch the columns $w^{-1}(M+1)$ and $w^{-1}(M+1)-1$, and get an orbit which is relevant and strictly smaller than the original orbit.

Thus, $w(w^{-1}(M+1)-1)\in [M+2,N]$.

Note that, according to the $\psi$-relevance condition, $w^{-1}(M+2)<w^{-1}(M+1)$, and $w^{-1}(M+2)\geq w^{-1}(w(w^{-1}(M+1)-1))=w^{-1}(M+1)-1$. It implies that $w^{-1}(M+2)=w^{-1}(M+1)-1$. Repeating the argument, one obtains the condition (5).

Then, we prove (1) as follows. If $(w,\beta)$ appears as $(w's,s\beta')$ for some $(w',\beta')\in RB^\aff(M|N)$, where $w'=ws$ and $\beta'=s\beta$, then the bad case is when there exists $k$ such that $\iota+kN\subseteq \sigma$, as in the case (1.b) of Proposition~\ref{right yoko}. If $(w',\beta')$ is not in $RB^\aff(M|N)$, it means there exists $k\in \BZ$ such that $i+1+kN\in \beta'$ but $i+kN\notin \beta'$, which breaks the saturation requirement in condition (2) of Proposition~\ref{prop 3.2.7}. In this case, we have $i+kN\in \sigma=\sigma(w's,s\beta')=\sigma(\tilde{w})$ and $i+1+kN\notin \beta$. The conditions (2)--(4) are obtained in the same way.
\end{proof}

\begin{prop}\label{prop 9.4.16}
For any reduced element $(w,\beta)$,
\begin{enumerate}
\item  If $i+1, i\in w^{-1}([1,M]_N)$, and $i+1\in\beta$, then $i\in \beta$;
\item If $i+1\in w^{-1}([M+1]_N)+1$, and $i+1\in \beta$, then $i-N+M\in \beta$;
\item If $i+1\in w(\beta)$, then $i\in w(\beta)$ if $i+1\in[2,M]_N$, and $i-N+M\in w(\beta)$ if $i+1\in[1]_N$.
\end{enumerate}
\end{prop}
\begin{proof}
We prove (1) and (2), and the claim (3) follows in the same way by considering the left action.

(1). If $i\notin \beta$, it implies that $w(i)>w(i+1)$, since otherwise, by condition (2) of Proposition~\ref{prop 3.2.7}, $i+1\in\beta$ would imply $i\in\beta$.

According to the requirement (1) in the definition of the reduced element, $i+1\in \beta$ implies that $i\in \beta$. Otherwise, $i+kN\in\sigma$ implies $k< 0$. If $i+1+kN\in \sigma$ or $\notin \beta$, it implies that $i+1+(k+1)N\notin \beta$, so $i+1\notin \beta$, a contradiction.

(2). If $i+1\leq w^{-1}(M+1)+1$, then the claim just follows from the requirement on $\beta$, see Proposition~\ref{prop 3.2.7}(3). Also, by the requirement (2) in the definition of the reduced element, $w^{-1}(M+1)+1+N\notin \beta$.
\end{proof}

We denote
\[
\Theta_{\tilde w} =\beta\sqcup\bigcup_{a=M+1}^N \bigl(w^{-1}(a)-N\BZ_{\geq0}\bigr).
\]

In the left matrix of Figure~\ref{fig:affine-orbit-classification}, $\Theta_{\tilde w}$ consists of the column indices of the blue dots and of the dots in rows $a-kN$, where $M+1\leq a\leq N$ and $k\geq0$.

The above proposition and conditions (2,5) in the definition of the reduced subset $RB^\aff_{\on{red}}(M|N)$ say that $\Theta_{\tilde w}$ and $w(\beta)$ do not have holes as subsets of $\BZ$ and $[1,M]_N$, respectively. It implies that we can translate any reduced element to an element coming from the finite type Gaiotto category, which we have already studied in the former section.

Figures~\ref{fig:reduced-affine-indices} and~\ref{fig:nonreduced-dispersed-indices} compare a reduced index with indices which are not reduced.

\begin{figure}[htbp]
\centering
\[
M=1,\qquad N=3,\qquad w(i+3)=w(i)+3.
\]
\resizebox{\linewidth}{!}{\begin{tikzpicture}[x=1cm,y=-1cm,font=\small]
\definecolor{orbitblue}{RGB}{32,94,176}
\definecolor{orbitred}{RGB}{184,40,47}

\begin{scope}[xshift=0cm]
\node[font=\normalsize] at (2.205,-2.200) {Reduced}; \node at (2.205,-1.750) {$(w(1),\ldots,w(3))=(3,2,1)$}; \node at (2.205,-1.350) {$\beta=3\BZ_{\leq0}$}; \node at (2.205,-0.950) {$\sigma=\varnothing$}; \draw[black!13,line width=.25pt] (0.000,0.000)--(0.000,4.410); \draw[black!13,line width=.25pt] (0.000,0.000)--(4.410,0.000); \draw[black!13,line width=.25pt] (0.490,0.000)--(0.490,4.410); \draw[black!13,line width=.25pt] (0.000,0.490)--(4.410,0.490); \draw[black!13,line width=.25pt] (0.980,0.000)--(0.980,4.410); \draw[black!13,line width=.25pt] (0.000,0.980)--(4.410,0.980); \draw[black!13,line width=.25pt] (1.470,0.000)--(1.470,4.410); \draw[black!13,line width=.25pt] (0.000,1.470)--(4.410,1.470); \draw[black!13,line width=.25pt] (1.960,0.000)--(1.960,4.410); \draw[black!13,line width=.25pt] (0.000,1.960)--(4.410,1.960); \draw[black!13,line width=.25pt] (2.450,0.000)--(2.450,4.410); \draw[black!13,line width=.25pt] (0.000,2.450)--(4.410,2.450); \draw[black!13,line width=.25pt] (2.940,0.000)--(2.940,4.410); \draw[black!13,line width=.25pt] (0.000,2.940)--(4.410,2.940); \draw[black!13,line width=.25pt] (3.430,0.000)--(3.430,4.410); \draw[black!13,line width=.25pt] (0.000,3.430)--(4.410,3.430); \draw[black!13,line width=.25pt] (3.920,0.000)--(3.920,4.410); \draw[black!13,line width=.25pt] (0.000,3.920)--(4.410,3.920); \draw[black!13,line width=.25pt] (4.410,0.000)--(4.410,4.410); \draw[black!13,line width=.25pt] (0.000,4.410)--(4.410,4.410); \node[anchor=south,font=\scriptsize] at (0.245,-0.030) {$-2$}; \node[anchor=east,font=\scriptsize] at (-0.030,0.245) {$-2$}; \node[anchor=south,font=\scriptsize] at (0.735,-0.030) {$-1$}; \node[anchor=east,font=\scriptsize] at (-0.030,0.735) {$-1$}; \node[anchor=south,font=\scriptsize] at (1.225,-0.030) {$0$}; \node[anchor=east,font=\scriptsize] at (-0.030,1.225) {$0$}; \node[anchor=south,font=\scriptsize] at (1.715,-0.030) {$1$}; \node[anchor=east,font=\scriptsize] at (-0.030,1.715) {$1$}; \node[anchor=south,font=\scriptsize] at (2.205,-0.030) {$2$}; \node[anchor=east,font=\scriptsize] at (-0.030,2.205) {$2$}; \node[anchor=south,font=\scriptsize] at (2.695,-0.030) {$3$}; \node[anchor=east,font=\scriptsize] at (-0.030,2.695) {$3$}; \node[anchor=south,font=\scriptsize] at (3.185,-0.030) {$4$}; \node[anchor=east,font=\scriptsize] at (-0.030,3.185) {$4$}; \node[anchor=south,font=\scriptsize] at (3.675,-0.030) {$5$}; \node[anchor=east,font=\scriptsize] at (-0.030,3.675) {$5$}; \node[anchor=south,font=\scriptsize] at (4.165,-0.030) {$6$}; \node[anchor=east,font=\scriptsize] at (-0.030,4.165) {$6$}; \draw[black!40,dashed,line width=.6pt] (0.000,0.000) rectangle (1.470,1.470); \draw[black!40,dashed,line width=.6pt] (2.940,2.940) rectangle (4.410,4.410); \draw[black!70,line width=.8pt] (1.470,1.470) rectangle (2.940,2.940); \draw[->,black!55] (0,-0.490)--(4.410,-0.490) node[right] {$i$}; \draw[->,black!55] (-0.580,0)--(-0.580,4.410) node[below] {$j$}; \draw[orbitred!65,densely dashed] (2.205,0)--(2.205,4.410); \draw[orbitred!65,densely dashed] (0,2.205)--(2.205,2.205); \fill[black!85] (0.245,1.225) circle (2.8pt); \fill[black!85] (0.735,0.735) circle (2.8pt); \fill[orbitblue] (1.225,0.245) circle (2.8pt); \fill[black!85] (1.715,2.695) circle (2.8pt); \fill[orbitred] (2.205,2.205) circle (2.8pt); \fill[black!85] (2.695,1.715) circle (2.8pt); \fill[black!85] (3.185,4.165) circle (2.8pt); \fill[black!85] (3.675,3.675) circle (2.8pt); \fill[black!85] (4.165,3.185) circle (2.8pt); \node[black!45] at (-0.230,-0.240) {$\ddots$}; \node[black!45] at (4.590,4.590) {$\ddots$}; \draw[black!35] (0.000,5.060)--(4.410,5.060); \node[anchor=east] at (-0.160,5.060) {$\Theta_{\tilde w}$}; \fill[black!80] (0.245,5.060) circle (2.2pt); \fill[black!80] (0.735,5.060) circle (2.2pt); \fill[black!80] (1.225,5.060) circle (2.2pt); \fill[black!80] (1.715,5.060) circle (2.2pt); \fill[black!80] (2.205,5.060) circle (2.2pt); \draw[black!55,fill=white,line width=.6pt] (2.695,5.060) circle (2.2pt); \draw[black!55,fill=white,line width=.6pt] (3.185,5.060) circle (2.2pt); \draw[black!55,fill=white,line width=.6pt] (3.675,5.060) circle (2.2pt); \draw[black!55,fill=white,line width=.6pt] (4.165,5.060) circle (2.2pt); \node at (2.205,5.610) {$\Theta_{\tilde w}=(-\infty,2]$};
\end{scope}
\begin{scope}[xshift=7.1cm]
\node[font=\normalsize] at (2.205,-2.200) {Not reduced}; \node at (2.205,-1.750) {$(w(1),\ldots,w(3))=(3,2,1)$}; \node at (2.205,-1.350) {$\beta=3\BZ_{\leq2}$}; \node at (2.205,-0.950) {$\sigma=\{6\}$}; \draw[black!13,line width=.25pt] (0.000,0.000)--(0.000,4.410); \draw[black!13,line width=.25pt] (0.000,0.000)--(4.410,0.000); \draw[black!13,line width=.25pt] (0.490,0.000)--(0.490,4.410); \draw[black!13,line width=.25pt] (0.000,0.490)--(4.410,0.490); \draw[black!13,line width=.25pt] (0.980,0.000)--(0.980,4.410); \draw[black!13,line width=.25pt] (0.000,0.980)--(4.410,0.980); \draw[black!13,line width=.25pt] (1.470,0.000)--(1.470,4.410); \draw[black!13,line width=.25pt] (0.000,1.470)--(4.410,1.470); \draw[black!13,line width=.25pt] (1.960,0.000)--(1.960,4.410); \draw[black!13,line width=.25pt] (0.000,1.960)--(4.410,1.960); \draw[black!13,line width=.25pt] (2.450,0.000)--(2.450,4.410); \draw[black!13,line width=.25pt] (0.000,2.450)--(4.410,2.450); \draw[black!13,line width=.25pt] (2.940,0.000)--(2.940,4.410); \draw[black!13,line width=.25pt] (0.000,2.940)--(4.410,2.940); \draw[black!13,line width=.25pt] (3.430,0.000)--(3.430,4.410); \draw[black!13,line width=.25pt] (0.000,3.430)--(4.410,3.430); \draw[black!13,line width=.25pt] (3.920,0.000)--(3.920,4.410); \draw[black!13,line width=.25pt] (0.000,3.920)--(4.410,3.920); \draw[black!13,line width=.25pt] (4.410,0.000)--(4.410,4.410); \draw[black!13,line width=.25pt] (0.000,4.410)--(4.410,4.410); \node[anchor=south,font=\scriptsize] at (0.245,-0.030) {$-2$}; \node[anchor=east,font=\scriptsize] at (-0.030,0.245) {$-2$}; \node[anchor=south,font=\scriptsize] at (0.735,-0.030) {$-1$}; \node[anchor=east,font=\scriptsize] at (-0.030,0.735) {$-1$}; \node[anchor=south,font=\scriptsize] at (1.225,-0.030) {$0$}; \node[anchor=east,font=\scriptsize] at (-0.030,1.225) {$0$}; \node[anchor=south,font=\scriptsize] at (1.715,-0.030) {$1$}; \node[anchor=east,font=\scriptsize] at (-0.030,1.715) {$1$}; \node[anchor=south,font=\scriptsize] at (2.205,-0.030) {$2$}; \node[anchor=east,font=\scriptsize] at (-0.030,2.205) {$2$}; \node[anchor=south,font=\scriptsize] at (2.695,-0.030) {$3$}; \node[anchor=east,font=\scriptsize] at (-0.030,2.695) {$3$}; \node[anchor=south,font=\scriptsize] at (3.185,-0.030) {$4$}; \node[anchor=east,font=\scriptsize] at (-0.030,3.185) {$4$}; \node[anchor=south,font=\scriptsize] at (3.675,-0.030) {$5$}; \node[anchor=east,font=\scriptsize] at (-0.030,3.675) {$5$}; \node[anchor=south,font=\scriptsize] at (4.165,-0.030) {$6$}; \node[anchor=east,font=\scriptsize] at (-0.030,4.165) {$6$}; \draw[black!40,dashed,line width=.6pt] (0.000,0.000) rectangle (1.470,1.470); \draw[black!40,dashed,line width=.6pt] (2.940,2.940) rectangle (4.410,4.410); \draw[black!70,line width=.8pt] (1.470,1.470) rectangle (2.940,2.940); \draw[->,black!55] (0,-0.490)--(4.410,-0.490) node[right] {$i$}; \draw[->,black!55] (-0.580,0)--(-0.580,4.410) node[below] {$j$}; \draw[orbitred!65,densely dashed] (2.205,0)--(2.205,4.410); \draw[orbitred!65,densely dashed] (0,2.205)--(2.205,2.205); \fill[black!85] (0.245,1.225) circle (2.8pt); \fill[black!85] (0.735,0.735) circle (2.8pt); \fill[orbitblue] (1.225,0.245) circle (2.8pt); \fill[black!85] (1.715,2.695) circle (2.8pt); \fill[orbitred] (2.205,2.205) circle (2.8pt); \fill[orbitblue] (2.695,1.715) circle (2.8pt); \fill[black!85] (3.185,4.165) circle (2.8pt); \fill[black!85] (3.675,3.675) circle (2.8pt); \fill[orbitblue] (4.165,3.185) circle (2.8pt); \draw[orbitblue,line width=.7pt] (4.165,3.185) circle (4.5pt); \node[black!45] at (-0.230,-0.240) {$\ddots$}; \node[black!45] at (4.590,4.590) {$\ddots$}; \draw[black!35] (0.000,5.060)--(4.410,5.060); \node[anchor=east] at (-0.160,5.060) {$\Theta_{\tilde w}$}; \fill[black!80] (0.245,5.060) circle (2.2pt); \fill[black!80] (0.735,5.060) circle (2.2pt); \fill[black!80] (1.225,5.060) circle (2.2pt); \fill[black!80] (1.715,5.060) circle (2.2pt); \fill[black!80] (2.205,5.060) circle (2.2pt); \fill[black!80] (2.695,5.060) circle (2.2pt); \draw[black!55,fill=white,line width=.6pt] (3.185,5.060) circle (2.2pt); \draw[black!55,fill=white,line width=.6pt] (3.675,5.060) circle (2.2pt); \fill[black!80] (4.165,5.060) circle (2.2pt); \node at (2.205,5.610) {$\Theta_{\tilde w}=(-\infty,3]\cup\{6\}$};
\end{scope}
\end{tikzpicture}
} \caption{A reduced index and an index which is not reduced. Here, blue dots represent $\beta$ and circled dots represent $\sigma$, as in Figure~\ref{fig:affine-orbit-classification}. Both indices are $\psi$-relevant. The red dot has row $M+1$ and column $w^{-1}(M+1)$. The solid square marks rows and columns $1,\ldots,N$. The right one is not reduced, since $3\in\beta$ but $3\notin\sigma$, and hence condition~(2) in Section~\ref{sec:reduced-affine-indices} fails.}
\label{fig:reduced-affine-indices}
\end{figure}

\begin{figure}[htbp]
\centering
\[
M=2,\qquad N=4,\qquad w(i+4)=w(i)+4.
\]
\resizebox{\linewidth}{!}{\begin{tikzpicture}[x=1cm,y=-1cm,font=\small]
\definecolor{orbitblue}{RGB}{32,94,176}
\definecolor{orbitred}{RGB}{184,40,47}

\begin{scope}[xshift=0cm]
\node[font=\normalsize] at (2.925,-2.200) {Not reduced}; \node at (2.925,-1.750) {$(w(1),\ldots,w(4))=(8,3,1,2)$}; \node at (2.925,-1.350) {$\beta=\{4k-1,4k:k\leq0\}$}; \node at (2.925,-0.950) {$\sigma=\varnothing$}; \draw[black!13,line width=.25pt] (0.000,0.000)--(0.000,9.450); \draw[black!13,line width=.25pt] (0.450,0.000)--(0.450,9.450); \draw[black!13,line width=.25pt] (0.900,0.000)--(0.900,9.450); \draw[black!13,line width=.25pt] (1.350,0.000)--(1.350,9.450); \draw[black!13,line width=.25pt] (1.800,0.000)--(1.800,9.450); \draw[black!13,line width=.25pt] (2.250,0.000)--(2.250,9.450); \draw[black!13,line width=.25pt] (2.700,0.000)--(2.700,9.450); \draw[black!13,line width=.25pt] (3.150,0.000)--(3.150,9.450); \draw[black!13,line width=.25pt] (3.600,0.000)--(3.600,9.450); \draw[black!13,line width=.25pt] (4.050,0.000)--(4.050,9.450); \draw[black!13,line width=.25pt] (4.500,0.000)--(4.500,9.450); \draw[black!13,line width=.25pt] (4.950,0.000)--(4.950,9.450); \draw[black!13,line width=.25pt] (5.400,0.000)--(5.400,9.450); \draw[black!13,line width=.25pt] (5.850,0.000)--(5.850,9.450); \draw[black!13,line width=.25pt] (0.000,0.000)--(5.850,0.000); \draw[black!13,line width=.25pt] (0.000,0.450)--(5.850,0.450); \draw[black!13,line width=.25pt] (0.000,0.900)--(5.850,0.900); \draw[black!13,line width=.25pt] (0.000,1.350)--(5.850,1.350); \draw[black!13,line width=.25pt] (0.000,1.800)--(5.850,1.800); \draw[black!13,line width=.25pt] (0.000,2.250)--(5.850,2.250); \draw[black!13,line width=.25pt] (0.000,2.700)--(5.850,2.700); \draw[black!13,line width=.25pt] (0.000,3.150)--(5.850,3.150); \draw[black!13,line width=.25pt] (0.000,3.600)--(5.850,3.600); \draw[black!13,line width=.25pt] (0.000,4.050)--(5.850,4.050); \draw[black!13,line width=.25pt] (0.000,4.500)--(5.850,4.500); \draw[black!13,line width=.25pt] (0.000,4.950)--(5.850,4.950); \draw[black!13,line width=.25pt] (0.000,5.400)--(5.850,5.400); \draw[black!13,line width=.25pt] (0.000,5.850)--(5.850,5.850); \draw[black!13,line width=.25pt] (0.000,6.300)--(5.850,6.300); \draw[black!13,line width=.25pt] (0.000,6.750)--(5.850,6.750); \draw[black!13,line width=.25pt] (0.000,7.200)--(5.850,7.200); \draw[black!13,line width=.25pt] (0.000,7.650)--(5.850,7.650); \draw[black!13,line width=.25pt] (0.000,8.100)--(5.850,8.100); \draw[black!13,line width=.25pt] (0.000,8.550)--(5.850,8.550); \draw[black!13,line width=.25pt] (0.000,9.000)--(5.850,9.000); \draw[black!13,line width=.25pt] (0.000,9.450)--(5.850,9.450); \node[anchor=south,font=\scriptsize] at (0.225,-0.040) {$-4$}; \node[anchor=south,font=\scriptsize] at (0.675,-0.040) {$-3$}; \node[anchor=south,font=\scriptsize] at (1.125,-0.040) {$-2$}; \node[anchor=south,font=\scriptsize] at (1.575,-0.040) {$-1$}; \node[anchor=south,font=\scriptsize] at (2.025,-0.040) {$0$}; \node[anchor=south,font=\scriptsize] at (2.475,-0.040) {$1$}; \node[anchor=south,font=\scriptsize] at (2.925,-0.040) {$2$}; \node[anchor=south,font=\scriptsize] at (3.375,-0.040) {$3$}; \node[anchor=south,font=\scriptsize] at (3.825,-0.040) {$4$}; \node[anchor=south,font=\scriptsize] at (4.275,-0.040) {$5$}; \node[anchor=south,font=\scriptsize] at (4.725,-0.040) {$6$}; \node[anchor=south,font=\scriptsize] at (5.175,-0.040) {$7$}; \node[anchor=south,font=\scriptsize] at (5.625,-0.040) {$8$}; \node[anchor=east,font=\scriptsize] at (-0.040,0.225) {$-8$}; \node[anchor=east,font=\scriptsize] at (-0.040,0.675) {$-7$}; \node[anchor=east,font=\scriptsize] at (-0.040,1.125) {$-6$}; \node[anchor=east,font=\scriptsize] at (-0.040,1.575) {$-5$}; \node[anchor=east,font=\scriptsize] at (-0.040,2.025) {$-4$}; \node[anchor=east,font=\scriptsize] at (-0.040,2.475) {$-3$}; \node[anchor=east,font=\scriptsize] at (-0.040,2.925) {$-2$}; \node[anchor=east,font=\scriptsize] at (-0.040,3.375) {$-1$}; \node[anchor=east,font=\scriptsize] at (-0.040,3.825) {$0$}; \node[anchor=east,font=\scriptsize] at (-0.040,4.275) {$1$}; \node[anchor=east,font=\scriptsize] at (-0.040,4.725) {$2$}; \node[anchor=east,font=\scriptsize] at (-0.040,5.175) {$3$}; \node[anchor=east,font=\scriptsize] at (-0.040,5.625) {$4$}; \node[anchor=east,font=\scriptsize] at (-0.040,6.075) {$5$}; \node[anchor=east,font=\scriptsize] at (-0.040,6.525) {$6$}; \node[anchor=east,font=\scriptsize] at (-0.040,6.975) {$7$}; \node[anchor=east,font=\scriptsize] at (-0.040,7.425) {$8$}; \node[anchor=east,font=\scriptsize] at (-0.040,7.875) {$9$}; \node[anchor=east,font=\scriptsize] at (-0.040,8.325) {$10$}; \node[anchor=east,font=\scriptsize] at (-0.040,8.775) {$11$}; \node[anchor=east,font=\scriptsize] at (-0.040,9.225) {$12$}; \draw[black!40,dashed,line width=.6pt] (0.450,2.250) rectangle (2.250,4.050); \draw[black!40,dashed,line width=.6pt] (4.050,5.850) rectangle (5.850,7.650); \draw[black!70,line width=.8pt] (2.250,4.050) rectangle (4.050,5.850); \draw[->,black!55] (0,-0.490)--(5.850,-0.490) node[right] {$i$}; \draw[->,black!55] (-0.580,0)--(-0.580,9.450) node[below] {$j$}; \draw[orbitred!65,densely dashed] (2.925,0)--(2.925,9.450); \draw[orbitred!65,densely dashed] (0,5.175)--(2.925,5.175); \fill[orbitblue] (0.225,1.125) circle (2.8pt); \fill[black!85] (0.675,5.625) circle (2.8pt); \fill[black!85] (1.125,3.375) circle (2.8pt); \fill[orbitblue] (1.575,2.475) circle (2.8pt); \fill[orbitblue] (2.025,2.925) circle (2.8pt); \fill[black!85] (2.475,7.425) circle (2.8pt); \fill[orbitred] (2.925,5.175) circle (2.8pt); \fill[black!85] (3.375,4.275) circle (2.8pt); \fill[black!85] (3.825,4.725) circle (2.8pt); \fill[black!85] (4.275,9.225) circle (2.8pt); \fill[black!85] (4.725,6.975) circle (2.8pt); \fill[black!85] (5.175,6.075) circle (2.8pt); \fill[black!85] (5.625,6.525) circle (2.8pt); \node[black!45] at (-0.230,-0.240) {$\ddots$}; \node[black!45] at (6.030,9.630) {$\ddots$}; \draw[black!35] (0.000,10.100)--(5.850,10.100); \node[anchor=east] at (-0.160,10.100) {$\Theta_{\tilde w}$}; \fill[black!80] (0.225,10.100) circle (2.2pt); \fill[black!80] (0.675,10.100) circle (2.2pt); \fill[black!80] (1.125,10.100) circle (2.2pt); \fill[black!80] (1.575,10.100) circle (2.2pt); \fill[black!80] (2.025,10.100) circle (2.2pt); \draw[black!55,fill=white,line width=.6pt] (2.475,10.100) circle (2.2pt); \fill[black!80] (2.925,10.100) circle (2.2pt); \draw[black!55,fill=white,line width=.6pt] (3.375,10.100) circle (2.2pt); \draw[black!55,fill=white,line width=.6pt] (3.825,10.100) circle (2.2pt); \draw[black!55,fill=white,line width=.6pt] (4.275,10.100) circle (2.2pt); \draw[black!55,fill=white,line width=.6pt] (4.725,10.100) circle (2.2pt); \draw[black!55,fill=white,line width=.6pt] (5.175,10.100) circle (2.2pt); \draw[black!55,fill=white,line width=.6pt] (5.625,10.100) circle (2.2pt); \node at (2.925,10.650) {$\Theta_{\tilde w}=(-\infty,0]\cup\{2\}$};
\end{scope}
\begin{scope}[xshift=8cm]
\node[font=\normalsize] at (2.925,-2.200) {Not reduced}; \node at (2.925,-1.750) {$(w(1),\ldots,w(4))=(4,3,-3,6)$}; \node at (2.925,-1.350) {$\beta=\{4k-1,4k:k\leq0\}$}; \node at (2.925,-0.950) {$\sigma=\varnothing$}; \draw[black!13,line width=.25pt] (0.000,0.000)--(0.000,9.450); \draw[black!13,line width=.25pt] (0.450,0.000)--(0.450,9.450); \draw[black!13,line width=.25pt] (0.900,0.000)--(0.900,9.450); \draw[black!13,line width=.25pt] (1.350,0.000)--(1.350,9.450); \draw[black!13,line width=.25pt] (1.800,0.000)--(1.800,9.450); \draw[black!13,line width=.25pt] (2.250,0.000)--(2.250,9.450); \draw[black!13,line width=.25pt] (2.700,0.000)--(2.700,9.450); \draw[black!13,line width=.25pt] (3.150,0.000)--(3.150,9.450); \draw[black!13,line width=.25pt] (3.600,0.000)--(3.600,9.450); \draw[black!13,line width=.25pt] (4.050,0.000)--(4.050,9.450); \draw[black!13,line width=.25pt] (4.500,0.000)--(4.500,9.450); \draw[black!13,line width=.25pt] (4.950,0.000)--(4.950,9.450); \draw[black!13,line width=.25pt] (5.400,0.000)--(5.400,9.450); \draw[black!13,line width=.25pt] (5.850,0.000)--(5.850,9.450); \draw[black!13,line width=.25pt] (0.000,0.000)--(5.850,0.000); \draw[black!13,line width=.25pt] (0.000,0.450)--(5.850,0.450); \draw[black!13,line width=.25pt] (0.000,0.900)--(5.850,0.900); \draw[black!13,line width=.25pt] (0.000,1.350)--(5.850,1.350); \draw[black!13,line width=.25pt] (0.000,1.800)--(5.850,1.800); \draw[black!13,line width=.25pt] (0.000,2.250)--(5.850,2.250); \draw[black!13,line width=.25pt] (0.000,2.700)--(5.850,2.700); \draw[black!13,line width=.25pt] (0.000,3.150)--(5.850,3.150); \draw[black!13,line width=.25pt] (0.000,3.600)--(5.850,3.600); \draw[black!13,line width=.25pt] (0.000,4.050)--(5.850,4.050); \draw[black!13,line width=.25pt] (0.000,4.500)--(5.850,4.500); \draw[black!13,line width=.25pt] (0.000,4.950)--(5.850,4.950); \draw[black!13,line width=.25pt] (0.000,5.400)--(5.850,5.400); \draw[black!13,line width=.25pt] (0.000,5.850)--(5.850,5.850); \draw[black!13,line width=.25pt] (0.000,6.300)--(5.850,6.300); \draw[black!13,line width=.25pt] (0.000,6.750)--(5.850,6.750); \draw[black!13,line width=.25pt] (0.000,7.200)--(5.850,7.200); \draw[black!13,line width=.25pt] (0.000,7.650)--(5.850,7.650); \draw[black!13,line width=.25pt] (0.000,8.100)--(5.850,8.100); \draw[black!13,line width=.25pt] (0.000,8.550)--(5.850,8.550); \draw[black!13,line width=.25pt] (0.000,9.000)--(5.850,9.000); \draw[black!13,line width=.25pt] (0.000,9.450)--(5.850,9.450); \node[anchor=south,font=\scriptsize] at (0.225,-0.040) {$-4$}; \node[anchor=south,font=\scriptsize] at (0.675,-0.040) {$-3$}; \node[anchor=south,font=\scriptsize] at (1.125,-0.040) {$-2$}; \node[anchor=south,font=\scriptsize] at (1.575,-0.040) {$-1$}; \node[anchor=south,font=\scriptsize] at (2.025,-0.040) {$0$}; \node[anchor=south,font=\scriptsize] at (2.475,-0.040) {$1$}; \node[anchor=south,font=\scriptsize] at (2.925,-0.040) {$2$}; \node[anchor=south,font=\scriptsize] at (3.375,-0.040) {$3$}; \node[anchor=south,font=\scriptsize] at (3.825,-0.040) {$4$}; \node[anchor=south,font=\scriptsize] at (4.275,-0.040) {$5$}; \node[anchor=south,font=\scriptsize] at (4.725,-0.040) {$6$}; \node[anchor=south,font=\scriptsize] at (5.175,-0.040) {$7$}; \node[anchor=south,font=\scriptsize] at (5.625,-0.040) {$8$}; \node[anchor=east,font=\scriptsize] at (-0.040,0.225) {$-8$}; \node[anchor=east,font=\scriptsize] at (-0.040,0.675) {$-7$}; \node[anchor=east,font=\scriptsize] at (-0.040,1.125) {$-6$}; \node[anchor=east,font=\scriptsize] at (-0.040,1.575) {$-5$}; \node[anchor=east,font=\scriptsize] at (-0.040,2.025) {$-4$}; \node[anchor=east,font=\scriptsize] at (-0.040,2.475) {$-3$}; \node[anchor=east,font=\scriptsize] at (-0.040,2.925) {$-2$}; \node[anchor=east,font=\scriptsize] at (-0.040,3.375) {$-1$}; \node[anchor=east,font=\scriptsize] at (-0.040,3.825) {$0$}; \node[anchor=east,font=\scriptsize] at (-0.040,4.275) {$1$}; \node[anchor=east,font=\scriptsize] at (-0.040,4.725) {$2$}; \node[anchor=east,font=\scriptsize] at (-0.040,5.175) {$3$}; \node[anchor=east,font=\scriptsize] at (-0.040,5.625) {$4$}; \node[anchor=east,font=\scriptsize] at (-0.040,6.075) {$5$}; \node[anchor=east,font=\scriptsize] at (-0.040,6.525) {$6$}; \node[anchor=east,font=\scriptsize] at (-0.040,6.975) {$7$}; \node[anchor=east,font=\scriptsize] at (-0.040,7.425) {$8$}; \node[anchor=east,font=\scriptsize] at (-0.040,7.875) {$9$}; \node[anchor=east,font=\scriptsize] at (-0.040,8.325) {$10$}; \node[anchor=east,font=\scriptsize] at (-0.040,8.775) {$11$}; \node[anchor=east,font=\scriptsize] at (-0.040,9.225) {$12$}; \draw[black!40,dashed,line width=.6pt] (0.450,2.250) rectangle (2.250,4.050); \draw[black!40,dashed,line width=.6pt] (4.050,5.850) rectangle (5.850,7.650); \draw[black!70,line width=.8pt] (2.250,4.050) rectangle (4.050,5.850); \draw[->,black!55] (0,-0.490)--(5.850,-0.490) node[right] {$i$}; \draw[->,black!55] (-0.580,0)--(-0.580,9.450) node[below] {$j$}; \draw[orbitred!65,densely dashed] (2.925,0)--(2.925,9.450); \draw[orbitred!65,densely dashed] (0,5.175)--(2.925,5.175); \fill[orbitblue] (0.225,2.925) circle (2.8pt); \fill[black!85] (0.675,3.825) circle (2.8pt); \fill[black!85] (1.125,3.375) circle (2.8pt); \fill[orbitblue] (1.575,0.675) circle (2.8pt); \fill[orbitblue] (2.025,4.725) circle (2.8pt); \fill[black!85] (2.475,5.625) circle (2.8pt); \fill[orbitred] (2.925,5.175) circle (2.8pt); \fill[black!85] (3.375,2.475) circle (2.8pt); \fill[black!85] (3.825,6.525) circle (2.8pt); \fill[black!85] (4.275,7.425) circle (2.8pt); \fill[black!85] (4.725,6.975) circle (2.8pt); \fill[black!85] (5.175,4.275) circle (2.8pt); \fill[black!85] (5.625,8.325) circle (2.8pt); \node[black!45] at (-0.230,-0.240) {$\ddots$}; \node[black!45] at (6.030,9.630) {$\ddots$}; \draw[black!35] (0.000,10.100)--(5.850,10.100); \node[anchor=east] at (-0.160,10.100) {$\Theta_{\tilde w}$}; \fill[black!80] (0.225,10.100) circle (2.2pt); \fill[black!80] (0.675,10.100) circle (2.2pt); \fill[black!80] (1.125,10.100) circle (2.2pt); \fill[black!80] (1.575,10.100) circle (2.2pt); \fill[black!80] (2.025,10.100) circle (2.2pt); \fill[black!80] (2.475,10.100) circle (2.2pt); \fill[black!80] (2.925,10.100) circle (2.2pt); \draw[black!55,fill=white,line width=.6pt] (3.375,10.100) circle (2.2pt); \draw[black!55,fill=white,line width=.6pt] (3.825,10.100) circle (2.2pt); \draw[black!55,fill=white,line width=.6pt] (4.275,10.100) circle (2.2pt); \draw[black!55,fill=white,line width=.6pt] (4.725,10.100) circle (2.2pt); \draw[black!55,fill=white,line width=.6pt] (5.175,10.100) circle (2.2pt); \draw[black!55,fill=white,line width=.6pt] (5.625,10.100) circle (2.2pt); \node at (2.925,10.650) {$\Theta_{\tilde w}=(-\infty,2]$};
\end{scope}
\end{tikzpicture}
} \caption{Two more indices which are not reduced. Both indices are still $\psi$-relevant. For the left one, $w^{-1}(3)=2$ and $w^{-1}(4)=-3$, thus, condition~(5) in Section~\ref{sec:reduced-affine-indices} fails, and for the right one, $w^{-1}(1)=7>w^{-1}(2)=0$ and $\sigma=\varnothing$, thus, condition~(3) fails.}
\label{fig:nonreduced-dispersed-indices}
\end{figure}
\FloatBarrier

\subsubsection{Translation to the finite flag variety}
\begin{prop}\label{prop 9.4.17}
Up to left and right (length-zero) translations, every orbit indexed by $RB^\aff_{\on{red}}(M|N)$ meets the finite flag variety $\GL_N/B_N\subset\Fl_N$. Its standard object is the $!$-extension of the image of the corresponding finite type standard object under the equivalence in Lemma \ref{13.2.1}.
\end{prop}
\begin{proof}
In terms of the figures above, a right translation moves every dot and $\beta$ by the same number of columns, and a left translation moves the dots in rows $[1,M]_N$ by the same number of steps in this set (i.e., $1\mapsto 2\mapsto\cdots\mapsto M\mapsto N+1\mapsto\cdots$), keeping the other rows and $\beta$ fixed. By a right translation, we obtain $\tilde w_1=(w_1,\beta_1)$ such that $\Theta_{\tilde w_1}=(-\infty,N]$. Furthermore, by a left translation, we obtain $\tilde w_2=(w_2,\beta_2)$ such that
\[
\beta_2=\beta_1,\qquad w_2(\beta_2)=(-\infty,N]\cap[1,M]_N.
\]
Thus,
\[
w_2\bigl((-\infty,N]\bigr) =w_2(\beta_2)\sqcup \bigcup_{a=M+1}^N(a-N\BZ_{\geq0}) =(-\infty,N],
\]
which implies $w_2\in S_N$. Note that replacing each entry of $g_{w_2(\beta_2)}$ by its constant term does not change the corresponding point of $\Fl_N$, and hence this point lies in $\GL_N/B_N$. Finally, the assertion about standard objects follows from Lemma~\ref{13.2.1}.
\end{proof}
\begin{rem}
Combined with the reduction to $RB^\aff_{\on{red}}(M|N)$ above, this shows that every affine standard object is obtained from a standard object of the finite type category by left and right Hecke convolutions, and the convolution diagrams on these orbits are isomorphisms.
\end{rem}

\begin{lem}\label{lem:affine-finite-IC-expansion}
Assume $M<N$. Each affine IC class has a finite expansion in the standard classes. Each affine standard class also has a finite expansion in the IC classes.
\end{lem}
\begin{proof}
This follows from the finite number of relevant orbits in each orbit closure and the triangularity of the transition matrix between the standard and IC classes.
\end{proof}

\subsubsection{The closure relation for $I_M\ltimes U_{M,N}(\bF)$-orbits in $\Fl_N$} In this subsection, we study the closure relation of the Gaiotto orbits in the affine flags.

We denote by $\leq_{I,(N,N)}$ the closure relation of the $I_N$-orbits in $\Fl_N\times \bF^N$ given in \cite[Proposition 5]{[FGT]}.

Let
\[
\rho^+:=(\underbrace{0,\dots,0}_{M},1,2,\dots,N-M),
\]
and
\[
J=(I_M\times T_{N-M}(\bO))\ltimes U_{M-1,N}(\bF), \qquad v_0=e_{M+1},\qquad K=\Stab_J(v_0),
\]
and then
\[
K=(I_M\ltimes U_{M,N}(\bF))\rtimes T_{\{M+2,\ldots,N\}}(\bO).
\]
\begin{lem}\label{lem:fixed-vector-closure-comparison}
For $x,y\in\Fl_N$, we have $y\in\overline{K\cdot x}$ if and only if, for all sufficiently large integers $n$,
\begin{equation}\label{nn cloure}
(t^{-n\rho^+}y,t^{-n\rho^+}v_0) \leq_{I,(N,N)} (t^{-n\rho^+}x,t^{-n\rho^+}v_0).
\end{equation}
\end{lem}
\begin{proof}
Let $J_n=\Ad_{t^{n\rho^+}}I_N$ and $K_n=\Stab_{J_n}(v_0)$. Note that the closure of the $K$-orbit passing through $x$ is the union of the closures of the $(K\cap J_n)$-orbits. When $n$ is big enough, we can write each element of $K_n$ as a product of an element of $K\cap J_n$ and an element of the unipotent subgroup opposite to $U_{M-1,N}$, and the latter fixes $x$, thus, $(K\cap J_n)x=K_nx$. Hence, $y$ lies in the closure of the $K$-orbit passing through $x$ if and only if $y\in\overline{K_nx}$ for all sufficiently large $n$.

Since $K_n=\Stab_{J_n}(v_0)$, we have
\[
y\in\overline{K_nx} \quad\Longleftrightarrow\quad (y,v_0)\in\overline{J_n\cdot(x,v_0)}.
\]
Indeed, for $v\in J_nv_0$, the matrix $\Id+(v-v_0)e_{M+1}^*$ belongs to $J_n$, depends regularly on $v$, and sends $v_0$ to $v$. After translation by $t^{-n\rho^+}$, we obtain \eqref{nn cloure}.
\end{proof}

As a corollary, we have
\begin{cor}
For two $I_M\ltimes U_{M,N}(\bF)$-orbits in $\Fl_N$ indexed by $\tilde{w}$ and $\tilde{w}'$, we have
\begin{equation}
\BO_{\tilde{w}}\subset \overline{\BO}_{\tilde{w}'}
\end{equation}
if and only if, for all sufficiently large $n$,
\[
t^{-n\rho^+}(w,e_{w(B)}) \leq_{I,(N,N)} t^{-n\rho^+}(w',e_{w'(B')}),
\]
where
\[
B=\beta\cup w^{-1}(M+1-N\BZ_{\geq0}),\qquad B'=\beta'\cup (w')^{-1}(M+1-N\BZ_{\geq0}).
\]
\end{cor}
\begin{proof}
Note that there is an order-preserving bijection between $I_M\ltimes U_{M,N}(\bF)$-orbits in $\Fl_N$ and $I_M\ltimes U_{M-1,N}(\bF)$-orbits in $\Fl_N\times\bF^M$, where $\bF^M$ is embedded in $\bF^N$ via $(f_1,\cdots,f_M)\mapsto(f_1,\cdots,f_M,1,0,\cdots,0)$. After translation by $t^{-n\rho^+}$, the corresponding $I_N$-orbit is represented by $t^{-n\rho^+}(w,e_{w(B)})$ for sufficiently large $n$. The assertion follows from Lemma~\ref{lem:fixed-vector-closure-comparison}.
\end{proof}

By Corollaries~\ref{cor 9.3.12} and \ref{cor 10.2.14}, with the normalization \eqref{eq:finite-IC-normalization}, and Propositions~\ref{prop:affine-reduction-to-reduced}, \ref{prop 9.4.16}, and \ref{prop 9.4.17}, we have

\begin{thm}
Assume $M<N$.
\begin{enumerate}
\item There exists a unique involution
\[
\bar{\cdot}\colon r \longmapsto \overline{r}
\]
on the direct sum over $\kappa,\chi,\chi'$ and the additive characters $\psi,-\psi$, whose component is
\[
\bar{\cdot}\colon \bK^{\aff,\mxd}_{\kappa,\chi,\chi',\psi}(M|N) \longrightarrow \bK^{\aff,\mxd}_{\kappa^{-1},\chi^{-1},(\chi')^{-1},-\psi}(M|N),
\]
such that its restriction to the finite type subcategory is the bar involution of Corollary~\ref{cor 10.2.14}, described in \eqref{eq:finite-bar-matrix}, and
\[
\overline{h_M r}=\overline{h_M}\,\overline{r}, \qquad \overline{rh_N}=\overline{r}\,\overline{h_N}
\]
for any $h_M\in \bH^{\aff,\mxd}_{\kappa,\chi'',\chi}(M)$, $h_N\in \bH^{\aff,\mxd}_{\kappa,\chi',\chi'''}(N)$, and $r\in \bK^{\aff,\mxd}_{\kappa,\chi,\chi',\psi}(M|N)$.

\item The involution in \textup{(1)} is induced by the Grothendieck--Verdier duality, with additive character $-\psi$ on the target.
\end{enumerate}
\end{thm}

\subsection{The case $M=N$}\label{sec:balanced-affine-algorithm} In this subsection, we consider the category $\cC^{\aff,\mxd}_{\kappa,\chi,\chi'}(N|N)$ of Section~\ref{sec:type-a-conventions}. For an index $(w,\beta)$, we denote by $\sigma$ the set of maximal elements of $\beta$ for the order given by $i\leq j$ and $w(i)\leq w(j)$. The left and right convolution formulas are obtained in the same way as in Section~\ref{sec:finite-type-monodromic}, using the periodic characters of Section~\ref{sec:periodic-characters} and the torus relations \eqref{eq:periodic-torus-lifts}. For trivial $\kappa,\chi,\chi'$, these are the formulas of \cite[Proposition 2]{[FGT]}.

\subsubsection{Reduction to the finite type category} Let $\Fl_N^{\fin}=\GL_N(\bO)/I_N\simeq\GL_N/B_N$.
\begin{lem}\label{lem:balanced-affine-reduction}
Every relevant orbit is obtained from an orbit meeting $\Fl_N^{\fin}\times k^N$ by left and right convolutions, and the corresponding convolution diagram is an isomorphism.
\end{lem}
\begin{proof}
The reduction in Proposition~\ref{prop:affine-reduction-to-reduced} and the argument in Proposition~\ref{prop 9.4.16}(1,3) apply here. Namely, by convolutions decreasing $l(w)$, we obtain a pair for which each of $\beta$ and $w(\beta)$ contains every integer smaller than one of its elements. As in the proof of Proposition~\ref{prop 9.4.17}, by right and left length-zero translations, we can make both of them equal to $(-\infty,N]$. Then, we have $w\in S_N$, and
\[
e_{w(\beta)}=\sum_{j=1}^N\sum_{m\geq0}t^m e_j =(1-t)^{-1}(e_1+\cdots+e_N).
\]
Note that the scalar matrix $(1-t)\Id\in I_N$ fixes $wI_N$ and sends $e_{w(\beta)}$ to $e_1+\cdots+e_N$. Thus, the pair belongs to the same orbit as $(w,e_1+\cdots+e_N)$, see also \cite[\S 2.5]{[FGT]}.
\end{proof}

For the factor $\bO^N$, we use the mixed $D^{!*}$-category of \cite[\S 3.2]{[BFGT]}, as in Section~\ref{sec:type-a-conventions}. Let $\mathrm{red}:\Fl_N^{\fin}\times\bO^N\to\Fl_N^{\fin}\times k^N$ be the reduction modulo $t$, and we denote by $\mathrm{red}^{\bullet}$ its normalized pullback. For $r\geq1$, we denote by $\mathrm{red}_r$ the reduction modulo $t$ from $\Fl_N^{\fin}\times(\bO/t^r\bO)^N$ to $\Fl_N^{\fin}\times k^N$. On this quotient, for an object $\mathcal A$ of the finite type category, the normalized pullback $\mathrm{red}^{\bullet}\mathcal A$ is represented by
\[
\mathrm{red}_r^*\mathcal A\langle N(r-1)\rangle.
\]
Here, the twisting and the equivariance are as in Section~\ref{sec:type-a-conventions}.

Let $j:\Fl_N^{\fin}\times\bO^N\hookrightarrow\Fl_N\times\bF^N$ be the closed embedding, and let
\begin{equation}\label{eq:balanced-affine-seed-functor}
\mathcal E(\mathcal A)=j_*\mathrm{red}^{\bullet}\mathcal A.
\end{equation}
Note that the first congruence subgroup fixes $\Fl_N^{\fin}$, and acts transitively on the lifts of each non-zero vector in $k^N$. Thus, by the smooth pullback and the closed pushforward, we obtain the following lemma, see also \cite[\S 2.6]{[FGT]}.
\begin{lem}\label{lem:balanced-affine-seeds}
For every relevant finite type orbit $u=(w,\beta_{\fin})$ with $\beta_{\fin}\ne\varnothing$, its preimage under $\mathrm{red}$ is one affine orbit. Its standard and IC objects are obtained by applying $\mathcal E$ to the corresponding finite type objects. Moreover, $\mathcal E$ commutes with Verdier duality.
\end{lem}

\subsubsection{Standard expansions and the bar involution} For $\tilde w=(w,\beta)$, let
\[
l(\tilde w)=l(w)+\#(\beta\setminus\BZ_{\leq0}) -\#(\BZ_{\leq0}\setminus\beta),
\]
as in \cite[\S 2.4]{[FGT]}. This is the relative dimension of the orbit with respect to the orbit indexed by $(e,\BZ_{\leq0})$. We denote by $RB_c$ the set of relevant indices for $(\kappa,\chi,\chi')$ in the connected component $c=N^{-1}\sum_{i=1}^N(w(i)-i)$. As in \cite[\S 2.7]{[FGT]}, we allow formal sums whose lengths are bounded above on each connected component, i.e., we consider
\[
\widehat{\bK}^{\aff,\mxd}_{\kappa,\chi,\chi'}(N|N) =\bigoplus_{c\in\BZ}\varinjlim_{L\in\BZ} \prod_{\substack{\tilde w\in RB_c\\l(\tilde w)\leq L}} \bR^{\mxd}H^{\aff,\mxd}_{\tilde w}.
\]
Here, each product has the product topology, where $\bR^{\mxd}$ is discrete, the transition maps extend the coefficients by zero, and we use the direct-limit topology and the direct-sum topology.

In order to see that the Verdier duality and the convolutions are continuous for this topology, note that, for a fixed index $\tilde u$ and an integer $L$, there are only finitely many $\tilde w$ in its component with $\tilde u\leq_{I,(N,N)}\tilde w$ and $l(\tilde w)\leq L$. Indeed, the closed conditions that the vector belongs to a given lattice in either flag bound $\max\beta$ and $\max w(\beta)$ from below. We write $w(i)=\pi(i)+N\nu_i$, with $\pi\in S_N$ and $\sum_i\nu_i=c$, and let $i+Nk_i$ be the largest element of $\beta$ in the residue class $i$. Let $x_i=k_i+1$ and $y_i=\nu_i+x_i$, and then we have $l(\tilde w)=l(w)+\sum_i x_i$.

We reorder the indices so that $\nu_1\geq\cdots\geq\nu_N$. The defining condition on $\beta$ forbids $x_i\geq x_j+2$ and $y_i\geq y_j+2$ simultaneously. It follows that $x_i\leq x_j+1$ and $y_i\geq y_j-1$ for $i<j$. Thus, $x_N,y_1$ are bounded below, and the adjacent differences $x_{i+1}-x_i$ and $y_i-y_{i+1}$ are at least $-1$.

By the affine length formula, we have $l(w)\geq\sum_{i<j}(\nu_i-\nu_j)-N(N-1)/2$, and hence
\[
\begin{split}
2l(\tilde w)\geq{}& N(x_N+y_1)-c-N(N-1)\\
&+\sum_{i=1}^{N-1}i(2N-2i-1)(x_{i+1}-x_i)\\
&+\sum_{i=1}^{N-1}(N-i)(2i-1)(y_i-y_{i+1}).
\end{split}
\]
Since all coefficients are positive, the upper bound on $l(\tilde w)$ bounds $x_N,y_1$ and these differences from above, and hence there are only finitely many possibilities for $(w,\beta)$. Note that the standard expansion of the Verdier dual of a standard class is supported on its orbit closure, thus, the duality preserves the upper length bound, and by the finiteness above, the coefficient of each standard class only depends on finitely many coefficients of the original sum. The convolution formulas have the same property, and they change the length by a bounded amount. Thus, these maps are continuous for the topology above. Note also that the finite sums of standard classes are dense.

For a relevant finite type orbit $(w,\beta_{\fin})$, we have
\[
\mathcal E(H^{\fin,\mxd}_{w,\beta_{\fin}})=
\begin{cases}
H^{\aff,\mxd}_{w,\,\beta_{\fin}\cup\BZ_{\leq0}}, &\beta_{\fin}\ne\varnothing,\\[2pt]
\displaystyle\sum_{\beta\subset\BZ_{\leq0}} (-[v])^{-\#(\BZ_{\leq0}\setminus\beta)}H^{\aff,\mxd}_{w,\beta}, &\beta_{\fin}=\varnothing.
\end{cases}
\]
Here, the sum runs over the relevant orbit indices $(w,\beta)$, and we omit the characters from the notation. The second formula is the standard expansion of the normalized constant sheaf on $t\bO^N$ in the vector factor.

By Lemmas~\ref{lem:balanced-affine-reduction} and \ref{lem:balanced-affine-seeds}, and \eqref{eq:finite-bar-matrix}, with the same argument as in the preceding theorem, we obtain the following proposition, see also \cite[Proposition 6]{[FGT]}.
\begin{prop}\label{prop:balanced-affine-bar}
There exists a unique continuous semilinear involution on the direct sum of the completed modules, whose component is
\[
\widehat{\bK}^{\aff,\mxd}_{\kappa,\chi,\chi'}(N|N) \longrightarrow \widehat{\bK}^{\aff,\mxd}_{\kappa^{-1},\chi^{-1},(\chi')^{-1}}(N|N), \qquad r\longmapsto\overline r,
\]
such that
\[
\overline{\mathcal E(\mathcal A)}=\mathcal E(\overline{\mathcal A}), \qquad \overline{h_1r}=\overline{h_1}\,\overline r, \qquad \overline{rh_2}=\overline r\,\overline{h_2},
\]
for finite type classes $\mathcal A$ and
\[
h_1\in\bH^{\aff,\mxd}_{\kappa,\chi'',\chi}(N), \qquad h_2\in\bH^{\aff,\mxd}_{\kappa,\chi',\chi'''}(N).
\]
On standard classes it is induced by Grothendieck--Verdier duality.
\end{prop}

\section{Pointwise purity of IC sheaves}\label{sec:purity}
\subsection{Contracting actions} In this section, we are going to show that the IC sheaves arising from the type A twisted Iwahori Gaiotto category are pointwise pure. Roughly speaking, for each relevant orbit, we construct a transversal slice of finite type with a $\BG_m$-action contracting it to a point of the orbit, and then the purity of the $*$-stalks and $!$-stalks at this point follows from the contraction argument, see Corollary~\ref{cor:contracting-slice-purity}. Here, we follow the conventions of Section~\ref{aff kl}, i.e., upper triangular Iwahori subgroups and unipotent groups, and the additive character on the constant term.

Note that there are several different equivalent models that we can work with.
\begin{enumerate}
\item ($\Fl_N$-model) we consider left $(I_M\ltimes U_{M,N}(\bF),\psi)$-equivariant sheaves on $\Fl_N$ in the case $M<N$, and we consider left $I_N$-equivariant sheaves on $\Fl_N\times \bF^N$ in the case $M=N$;
\item ($\Fl_G$-model) we consider left $(\GL_M(\bF)\ltimes U_{M,N}(\bF),\psi)$-equivariant sheaves on $\Fl_M\times \Fl_N$ in the case $M<N$, and we consider left $\GL_N(\bF)$-equivariant sheaves on $\Fl_N\times \Fl_N\times \bF^N$ in the case $M=N$;
\item ($\bX(\bF)$-model) in the case $M<N$, we consider right $(I_M\times I_N)$-equivariant sheaves on $(\GL_M(\bF)\ltimes U_{M,N}(\bF),\psi)\backslash (\GL_M(\bF)\times \GL_N(\bF))$, and in the case $M=N$, we consider right $(I_N\times I_N)$-equivariant sheaves on the mirabolic space $\GL_N(\bF)\backslash\bigl((\GL_N(\bF)\times \GL_N(\bF))\times \bF^N\bigr)$.
\end{enumerate}

In order to simplify the notation used in the main calculation, and due to the fact that the normal space of the $I_M\ltimes U_{M,N}(\bF)$-orbit in $\Fl_N$ is identified with the normal space of the $\GL_M(\bF)\ltimes U_{M,N}(\bF)$-orbit in $\Fl_M\times \Fl_N$, we use the $\Fl_N$-model in this subsection to construct the contracting cocharacters and calculate the normal weights. Then, in Section~\ref{14.2}, we construct finite type transversal slices from the global model for $M<N$ (Proposition~\ref{prop:global-model}) and from the local model for $M=N$ (Proposition~\ref{prop:balanced-finite-transversal}). Here, we denote by $\BG_m^{\mathrm{rot}}$ the group acting by loop rotation.

\begin{prop}\label{prop 9.5.1}

In the case $M<N$, we consider a relevant $(I_M\ltimes U_{M,N}(\bF),\psi)$-orbit in $\Fl_N$ corresponding to $\tilde{w}=(w,\beta)$, where $w$ is an extended Weyl group element of $\GL_N$, viewed as an $N$-periodic permutation matrix.

We denote by $E_{a,b}$ the elementary $N$-periodic matrix whose non-zero entries are $1$ in positions $(a+jN,b+jN)$, $j\in\BZ$. Put
\[
u_\sigma:=1+\sum_{a\in w(\sigma)}E_{a,M+1}, \qquad x_{\tilde w}:=u_\sigma wI_N\in \Fl_N.
\]

In the case $M=N$ and a non-zero mirabolic vector, put $u_\sigma:=\sum_{a\in w(\sigma)}e_a$, $x_{\tilde{w}}:= (w, u_\sigma)\in \Fl_N\times \bF^N$. We denote the orbit of $x_{\tilde w}$ in either case by $\BO_{\tilde w}$.

Then, there exists
\[
\lambda:\mathbb G_m\to T_N\times \mathbb G_m^{\mathrm{rot}}
\]
such that:
\begin{enumerate}
\item $\lambda$ fixes $x_{\tilde w}$;
\item $\lambda$ preserves any orbit of the acting group in the indicated model and the Whittaker character $\psi$;
\item if $M<N$, all $\lambda$-weights on
\[
T_{x_{\tilde w}}\Fl_N\big/T_{x_{\tilde w}}\BO_{\tilde w}
\]
\begin{samepage}
(resp. if $M=N$, all $\lambda$-weights on
\[
T_{x_{\tilde w}}(\Fl_N\times \bF^N)\big/T_{x_{\tilde w}}\BO_{\tilde w} \]) are strictly positive.\par
\end{samepage}
\end{enumerate}

Furthermore, let $X=\Fl_N$ for $M<N$ and $X=\Fl_N\times\bF^N$ for $M=N$, and let
\[
\pi:\overset{\circ}{\cP}_{\det}\to X
\]
be the punctured determinant line bundle (pulled back from $\Fl_N$ in the mirabolic case). Then, $\lambda$ admits a lift
\[
\widetilde\lambda:\mathbb G_m\to \operatorname{Aut}(\overset{\circ}{\cP}_{\det})
\]
whose induced action on the fiber $\pi^{-1}(x_{\tilde w})$ is trivial and which preserves the preimages of the orbits. Let $\widetilde{\BO}_{\tilde w}=\pi^{-1}(\BO_{\tilde w})$. For $\tilde x\in\pi^{-1}(x_{\tilde w})$, the natural map identifies
\[
T_{\tilde x}\overset{\circ}{\cP}_{\det}/ T_{\tilde x}\widetilde{\BO}_{\tilde w} \simeq T_{x_{\tilde w}}X/T_{x_{\tilde w}}\BO_{\tilde w}
\]
equivariantly, and thus, its weights are strictly positive.
\end{prop}

\begin{proof}
We work with the $N$-periodic basis $\{e_i\}_{i\in\mathbb Z}$ satisfying
\[
e_{i+N}=t^{-1}e_i.
\]

We want a one-parameter subgroup of $T_N\times \BG_m^{\mathrm{rot}}$ encoded by a function \[ \gamma:\mathbb Z\to \mathbb Z,\qquad \gamma_{i+N}=\gamma_i+(M+1), \] via \[ \lambda(z)\cdot e_i=z^{\gamma_i}e_i. \]

Here, we define a rational weight function
\[
\gamma':\mathbb Z\to \frac{1}{M+1}\BZ, \qquad \gamma'_{i+N}=\gamma'_i+1,
\]
which will determine a rational grading on the $N$-periodic basis $\{e_i\}_{i\in\mathbb Z}$, and hence on the elementary matrices by
\[
\mathrm{wt}_{\gamma'}(E_{p,q})=\gamma'_p-\gamma'_q.
\]

Then, we will take $\gamma=(M+1)\gamma'$. Note that the sign of every weight, as well as the vanishing of every weight, is unchanged by passing from $\gamma$ to $\gamma'$, and thus, in the calculation below, we denote $\mathrm{wt}=\mathrm{wt}_{\gamma'}$.

We define $\gamma'$ as follows.

\smallskip

\noindent \emph{(a) On the classes $M+1,\dots,N$.} Set
\[
\gamma'_k:=0\qquad (k=M+1,\dots,N).
\]
Then, for every $k=M+2,\dots,N$,
\[
\mathrm{wt}(E_{k-1,k}) = \gamma'_{k-1}-\gamma'_k = 0.
\]
Hence, $\lambda$ preserves the character $\psi$.

\smallskip
\noindent \emph{(b) On the classes $\{1,\dots,M\}$.} We choose the value of the weight similarly to \cite{[T1]}. Namely, let $B_0=\beta\setminus(\beta-N)$ (note that $\beta-N\subset\beta$). We claim that it has exactly one element in each of the residue classes $w^{-1}(1),\ldots,w^{-1}(M)$ modulo $N$. Indeed, for $M<N$, Proposition~\ref{prop 3.2.7} implies that $\beta$ contains all sufficiently negative indices in these classes. In the mirabolic case, we choose $j\in\beta$, which exists since the vector is non-zero. Then, for every sufficiently negative $i$, we have $i<j$ and $w(i)<w(j)$, and thus, the saturation condition on $\beta$ implies the same assertion. Since $\beta$ is bounded above, it has a unique maximal element in each class. Now, we consider the partial order on $B_0$, where $i\prec j$ means that $i<j$ and $w(i)<w(j)$. For $M<N$, we adjoin a maximal element $*$, with $i\prec *$ whenever $i<w^{-1}(M+1)$, in order to include the third condition of Proposition~\ref{prop 3.2.7}, and in the mirabolic case, we add no element. We denote $h(i)=\max\{n\mid i=i_0\prec i_1\prec\cdots\prec i_n\}$, where the $i_j$ are in the above partially ordered set, and define
\[
\gamma'_{w(i)}=-\frac{h(i)}{M+1}\quad(i\in B_0)
\]
and extend it by $\gamma'_{j+N}=\gamma'_j+1$. Note that the elements of $B_0$ which are maximal in this order are exactly the elements of $\sigma$. Hence,
\[
\gamma'_{w(i)}=0\ (i\in\sigma),\qquad -1<\gamma'_{w(i)}<0\ (i\in B_0\setminus\sigma).
\]
Here, the strict lower bound holds because of the denominator $M+1$, also for chains ending at $*$. Furthermore, for $i\prec j$ in $B_0$, we have $\gamma'_{w(i)}<\gamma'_{w(j)}$.

In particular, if $M<N$, then for any $a\in w(\sigma)$, we have $\gamma'_a=0=\gamma'_{M+1}$. Thus,
\[
\mathrm{wt}(E_{a,M+1})=0 \qquad (\forall a\in w(\sigma)),
\]
and $u_\sigma$ is fixed by $\lambda$, and so $x_{\tilde w}=u_\sigma wI_N$ is fixed. Similarly for the mirabolic case.

Next, we show that $\lambda$ preserves every orbit. Note that the action of $\lambda(z)$ is the composition of a constant diagonal element of $T_M$ and $r_z:t\mapsto z^{-(M+1)}t$, and both factors normalize the acting group. Thus, it is enough to show that one can left and right multiply the representative of an orbit by diagonal elements, such that the product equals the representative after the loop rotation. Namely, for any orbit index $\tilde w'=(w',\beta')$, the elements of $w'(\sigma')$ lie in distinct residue classes modulo $N$, and if $M<N$, then we choose $a(z)\in T_M$ such that $r_z(u_{\sigma'})=a(z)u_{\sigma'}a(z)^{-1}$. Since $w'$ is monomial, $b(z):=(w')^{-1}a(z)^{-1}r_z(w')$ belongs to $T_N$, and
\[
r_z(u_{\sigma'}w') =a(z)u_{\sigma'}w'b(z).
\]
In the mirabolic case, we choose $a(z)\in T_N$ with $r_z(u_{\sigma'})=a(z)u_{\sigma'}$, and with the same definition of $b(z)$, we have $r_z(w')=a(z)w'b(z)$. Since $T_M$ lies in the acting group and $T_N\subset I_N$, the action preserves every orbit.

Now, let us first prove the positivity of the weights in the case $M<N$. We denote $\mathfrak u_{M,N}=\Lie(U_{M,N})$, $\mathfrak i_M=\Lie(I_M)$ and $\mathfrak i_N=\Lie(I_N)$, and set
\[
\mathfrak h:=\Lie(I_M\ltimes U_{M,N}(\bF)),\qquad L:=\mathfrak h+\Ad_{u_\sigma}(w\mathfrak i_Nw^{-1}).
\]
We will show that every non-positive weight vector in $\mathfrak{gl}_N(\bF)$ lies in $L$, which implies the positivity, since the normal space $T_{x_{\tilde w}}\Fl_N/T_{x_{\tilde w}}\BO_{\tilde w}$ is identified with $\mathfrak{gl}_N(\bF)/L$, see the end of this case. Namely, we consider the root vectors in the following five groups.

\smallskip

\noindent \emph{(1) Root directions involving $w(\sigma)$ and $M+1$.} For $a,a'\in w(\sigma)$ and $r\le 0$, we have
\[
\mathrm{wt}(E_{a+rN,M+1})=r,\qquad \mathrm{wt}(E_{a+rN,a'})=r,\qquad \mathrm{wt}(E_{M+1+rN,a})=r.
\]
Now, we have
\[
\Ad_{u_\sigma}(E_{a+rN,a})=E_{a+rN,a}-E_{a+rN,M+1},
\]
and $E_{a+rN,a}\in \mathfrak i_M\subset \mathfrak h$, so
\[
E_{a+rN,M+1}\in L.
\]
Likewise,
\[
\Ad_{u_\sigma}(E_{a+rN,a'})=E_{a+rN,a'}-E_{a+rN,M+1}.
\]
If $a+rN\le a'$, then $E_{a+rN,a'}\in \mathfrak i_M$. If $a+rN>a'$, then $a>a'$, so
\[
w^{-1}(a+rN)=w^{-1}(a)+rN<w^{-1}(a'),
\]
by the requirement of $w(\sigma)$. Hence,
\[
E_{a+rN,a'}\in w\mathfrak i_Nw^{-1}.
\]
Thus,
\[
E_{a+rN,a'}\in L.
\]

Finally, for $r\le 0$ and $a\in w(\sigma)$, we have
\[
w^{-1}(M+1+rN)=w^{-1}(M+1)+rN<w^{-1}(a),
\]
because $w^{-1}(a)\in \sigma$ and every element of $\sigma$ is $>w^{-1}(M+1)$. Hence,
\[
E_{M+1+rN,a}\in w\mathfrak i_Nw^{-1}.
\]
Furthermore, we have
\[
\Ad_{u_\sigma}(E_{M+1+rN,a}) = E_{M+1+rN,a} +\sum_{a'\in w(\sigma)}E_{a'+rN,a} - E_{M+1+rN,M+1}-\sum_{a'\in w(\sigma)}E_{a'+rN,M+1},
\]
and the two sums are already in $L$, while
\[
E_{M+1+rN,M+1}\in w\mathfrak i_Nw^{-1} \qquad (r\le 0).
\]
Note that its conjugate by $u_\sigma$ is
\[
E_{M+1+rN,M+1}+\sum_{a'\in w(\sigma)}E_{a'+rN,M+1}.
\]
Since the sum already lies in $L$, we have $E_{M+1+rN,M+1}\in L$, and hence,
\[
E_{M+1+rN,a}\in L.
\]

\smallskip

\noindent \emph{(2) Root directions from the tail $\{M+2,\dots,N\}$.} For $M+1\le c\le b\le N$, the $\psi$-relevance condition implies
\[
w^{-1}(b)\le w^{-1}(c).
\]
Now, let $b,b'\in \{M+2,\dots,N\}$, $a'\in w(\beta\setminus((\beta-N)\cup \sigma))$ and $a\in w(\sigma)$. Then, we have
\[
\mathrm{wt}(E_{b+rN,M+1})=r, \qquad \mathrm{wt}(E_{b+rN,a})=r, \qquad \mathrm{wt}(E_{b+rN,b'})=r,\]
\[
r<\mathrm{wt}(E_{b+rN,a'})<r+1.
\]

If $\mathrm{wt}(E_{b+rN,M+1})\le 0$, then $r\le 0$, so
\[
w^{-1}(b+rN)=w^{-1}(b)+rN\le w^{-1}(M+1),
\]
and thus,
\[
E_{b+rN,M+1}\in w\mathfrak i_Nw^{-1}.
\]
It implies
\[
E_{b+rN,M+1}\in \Ad_{u_\sigma}(w\mathfrak i_Nw^{-1}),
\]
since $\Ad_{u_\sigma}(E_{b+rN,M+1})=E_{b+rN,M+1}$.

If $\mathrm{wt}(E_{b+rN,a})\le 0$, then again $r\le 0$, hence
\[
w^{-1}(b+rN)\le w^{-1}(M+1)<w^{-1}(a),
\]
and therefore
\[
E_{b+rN,a}\in w\mathfrak i_Nw^{-1}.
\]
Also, we have
\[
\Ad_{u_\sigma}(E_{b+rN,a})=E_{b+rN,a}-E_{b+rN,M+1},
\]
so $E_{b+rN,a}\in L$.

If $\mathrm{wt}(E_{b+rN,b'})\leq0$, then $r\leq0$. Note that, for $b<b'$, the vector lies in $\mathfrak u_{M,N}(\bF)\subset\mathfrak h$ for every $r$, and for $b\geq b'$, the relevance condition implies
\[
w^{-1}(b+rN)=w^{-1}(b)+rN\le w^{-1}(b'),
\]
hence
\[
E_{b+rN,b'}\in w\mathfrak i_Nw^{-1}.
\]
This vector is fixed by the conjugation by $u_\sigma$, and thus, in either case, we have
\[
E_{b+rN,b'}\in L.
\]

If $\mathrm{wt}(E_{b+rN,a'})\le 0$, then $r< 0$, so
\[
w^{-1}(b+rN)=w^{-1}(b)+rN\le w^{-1}(M+1)-N,
\]
while \[w^{-1}(a')+N>w^{-1}(M+1).\] Thus,
\[
E_{b+rN,a'}\in w\mathfrak i_Nw^{-1}.
\]
It implies
\[
E_{b+rN,a'}\in \Ad_{u_\sigma}(w\mathfrak i_Nw^{-1}).
\]
\smallskip

\noindent \emph{(3) Root directions from the first block to the tail.} Every root vector of type
\[
E_{i+rN,b}\qquad (1\le i\le M+1,\ b\in\{M+2,\dots,N\})
\]
belongs to
\[
\mathfrak u_{M,N}(\bF)\subset \mathfrak h.
\]

\smallskip

\noindent \emph{(4) Root directions from $w(\beta\setminus((\beta-N)\cup \sigma))$ to $[1,M+1]_N$.} For $a',a''\in w(\beta\setminus ((\beta-N)\cup \sigma) )$ and $a\in w(\sigma)$, we have
\[
\begin{gathered}
r-1<\mathrm{wt}(E_{a'+rN,M+1})<r,\\
r-1<\mathrm{wt}(E_{a'+rN,a})<r,\\
r-1<\mathrm{wt}(E_{a'+rN,a''})<r+1.
\end{gathered}
\]

For a non-positive weight, we necessarily have $r\leq 0$. Since $a'\in w(\beta)$, either (i) there exists $\tilde{a}\in w(\sigma)$ such that $a'\leq \tilde{a}$ and $w^{-1}(a')\leq w^{-1}(\tilde{a})$, or (ii) $w^{-1}(a')< w^{-1}(M+1)$. Either of these two conditions implies that $E_{a'+rN, M+1}\in L$ and $E_{a'+rN, a}\in L$.

Namely, in the case (i),
\[
\Ad_{u_\sigma}(E_{a'+rN,\tilde{a}})=E_{a'+rN,\tilde{a}}-E_{a'+rN,M+1},
\]
and $E_{a'+rN,\tilde{a}}\in \mathfrak i_M \cap w\mathfrak i_Nw^{-1}$, so
\[
E_{a'+rN,M+1}\in L.
\]
Note that
\[
\Ad_{u_\sigma}(E_{a'+rN,a})=E_{a'+rN,a}-E_{a'+rN,M+1},
\]
and we already have $E_{a'+rN,M+1}\in L$. Either $E_{a'+rN,a}\in w\mathfrak i_Nw^{-1}$ or $E_{a'+rN,a}\in \mathfrak i_M$ implies that $E_{a'+rN,a}\in L$. If not, we have $w^{-1}(a'+rN)>w^{-1}(a)$ and $a'+rN>a$, which contradicts the maximality of $a\in w(\sigma)$.

Likewise, in the case (ii),
\[
\Ad_{u_\sigma}(E_{a'+rN,M+1})=E_{a'+rN,M+1},
\]
and \[w^{-1}(a'+rN)=w^{-1}(a')+rN< w^{-1}(M+1),\] so $E_{a'+rN,M+1}\in L$. Similarly to the argument in the case (i), by the maximality of $a\in w(\sigma)$, we have $E_{a'+rN,a}\in L$.

Finally, we have
\[
\Ad_{u_\sigma}(E_{a'+rN,a''})=E_{a'+rN,a''}.
\]
If $E_{a'+rN,a''}\notin L$, then $E_{a'+rN,a''}\notin \mathfrak i_M$ and $E_{a'+rN,a''}\notin w\mathfrak i_Nw^{-1}$, which implies that $a'+rN>a''$ and $w^{-1}(a'+rN)>w^{-1}(a'')$. By the requirement $a''\in w(\beta\setminus (\beta-N))$, we have $r=0$, since otherwise, $r<0$ implies that $w^{-1}(a'')+N$ also belongs to $\beta$. Now, $a'>a''$ and $w^{-1}(a')>w^{-1}(a'')$ imply that $\gamma'_{a'}>\gamma'_{a''}$, and $\mathrm{wt}(E_{a'+rN,a''})>0$.

\smallskip\noindent \emph{(5) Root directions to $w(\beta\setminus((\beta-N)\cup \sigma))$ from $[1,M+1]_N$.} Note that the root directions from $w(\beta\setminus((\beta-N)\cup \sigma))$ have already been considered in (4). For $a\in w(\sigma)$ and $a'\in w(\beta\setminus((\beta-N)\cup \sigma))$, we have
\[
r<\mathrm{wt}(E_{a+rN,a'})<r+1,\qquad r<\mathrm{wt}(E_{M+1+rN,a'})<r+1.
\]

If $\mathrm{wt}(E_{a+rN,a'})\leq 0$, then $r<0$. Note that
\[
\Ad_{u_\sigma}(E_{a+rN,a'})=E_{a+rN,a'}.
\]
Either $E_{a+rN,a'}\in \mathfrak i_M$ or $E_{a+rN,a'}\in w\mathfrak i_Nw^{-1}$ implies that $E_{a+rN,a'}\in L$. Otherwise, $a+rN>a'$ and $w^{-1}(a+rN)>w^{-1}(a')$, so $w^{-1}(a')+N\in \beta$, which contradicts the requirement of $a'$.

If $\mathrm{wt}(E_{M+1+rN,a'})\leq 0$, then $r<0$. Note that
\[
\Ad_{u_\sigma}(E_{M+1+rN,a'})=\sum_{a\in w(\sigma)}E_{a+rN,a'}+E_{M+1+rN,a'},
\]
where $E_{a+rN,a'}$ is already in $L$, and thus, we need to show that $E_{M+1+rN,a'}\in w\mathfrak i_Nw^{-1}$. Otherwise, $w^{-1}(M+1)-N\geq w^{-1}(M+1+rN)>w^{-1}(a')$, so $w^{-1}(a')+N\in \beta$, which contradicts the requirement of $a'$.

Note that, if $\sigma$ is empty, then the conjugation formula for $E_{M+1+rN,M+1}$ has an empty sum, and the same calculation applies for $r\leq0$. The above five cases cover all the cases. Thus, every non-positive weight vector of $\mathfrak{gl}_N(\bF)$ lies in
\[
L=\mathfrak h+\Ad_{u_\sigma}(w\mathfrak i_Nw^{-1}).
\]
Therefore, all weights on
\[
T_{x_{\tilde w}}\Fl_N\big/T_{x_{\tilde w}}\BO_{\tilde w} \simeq \mathfrak{gl}_N(\bF)\big/L
\]
are strictly positive.

\smallskip

Now, we consider the case $M=N$ and a non-zero vector. In this case, the acting group is $I_N$, there is no Whittaker character, and we use the rational weight function $\gamma'$ constructed above.

To keep the notation similar, we identify $\GL_N(\bF)\ltimes \bF^N$ with the mirabolic loop subgroup of $\GL_{N+1}(\bF)$, and thus,
\[
\Fl_N\times \bF^N \simeq (\GL_N(\bF)\ltimes \bF^N)/(I_N\ltimes 0),
\]
and we let $\widetilde u_\sigma:=1+\sum_{a\in w(\sigma)}E_{a,N+1}$ be the corresponding unipotent element in $\GL_{N+1}(\bF)$. Here, the last column refers to the auxiliary basis vector fixed by $\lambda$, not to the periodic vector $e_{N+1}=t^{-1}e_1$ in the original $N$-dimensional representation.

Set
\[
\mathfrak h:=\mathfrak i_N\ltimes 0, \qquad L:=\mathfrak h+\Ad_{\widetilde u_\sigma}\bigl(w\mathfrak i_Nw^{-1}\ltimes 0\bigr) \subset \mathfrak{gl}_N(\bF)\ltimes \bF^N.
\]
Then,
\[
T_{x_{\tilde w}}(\Fl_N\times \bF^N)\big/T_{x_{\tilde w}}\BO_{\tilde w} \simeq (\mathfrak{gl}_N(\bF)\ltimes \bF^N)\big/L.
\]
We will show that every non-positive weight vector of $\mathfrak{gl}_N(\bF)\ltimes \bF^N$ lies in $L$.

\smallskip

\noindent \emph{(1) Flag directions.} The flag directions are of the following four types:
\[
E_{a+rN,a'}, E_{b+rN,a}, E_{a+rN,b}, E_{b+rN,b'}.
\]
Here, $a,a'\in w(\sigma)$, and $b,b'\in w(\beta\setminus((\beta-N)\cup \sigma))$.

These are exactly the flag directions considered above in the cases (1), (4) and (5), with the tail removed, but their conjugation formulas also involve the vector directions $E_{p,N+1}$. Thus, we first treat these vector directions in the next two cases, and then, the same matrix calculation as above applies to the flag directions.

\smallskip

\noindent \emph{(2) Vector directions indexed by $w(\sigma)$.} We identify the vector directions with the matrix units $E_{p,N+1}$. For $a\in w(\sigma)$ and $r\in \BZ$, we have
\[
\mathrm{wt}(E_{a+rN,N+1})=r.
\]
Assume $r\le 0$. Then,
\[
\Ad_{\widetilde u_\sigma}(E_{a+rN,a}) = E_{a+rN,a}-E_{a+rN,N+1}.
\]
Note that we have
\[
E_{a+rN,a}\in \mathfrak i_N\cap w\mathfrak i_Nw^{-1}.
\]
Therefore,
\[
(E_{a+rN,a},0)\in \mathfrak h, \qquad \Ad_{\widetilde u_\sigma}(E_{a+rN,a},0)\in L,
\]
and hence
\[
E_{a+rN,N+1}\in L.
\]

\smallskip

\noindent \emph{(3) Vector directions indexed by $w(\beta\setminus((\beta-N)\cup \sigma))$.} Let $b\in w(\beta\setminus((\beta-N)\cup \sigma))$. Then, we have $-1<\gamma'_b<0$, and
\[
r-1<\mathrm{wt}(E_{b+rN,N+1})<r.
\]
If $\mathrm{wt}(E_{b+rN,N+1})\le 0$, then necessarily $r\le 0$.

Since $\sigma$ is the set of $\prec$-maximal elements of $\beta\setminus(\beta-N)$, there exists $a\in w(\sigma)$ such that
\[
b<a, \qquad w^{-1}(b)<w^{-1}(a).
\]
Hence,
\[
b+rN<a, \qquad w^{-1}(b+rN)=w^{-1}(b)+rN<w^{-1}(a),
\]
so
\[
E_{b+rN,a}\in \mathfrak i_N\cap w\mathfrak i_Nw^{-1}.
\]
Furthermore,
\[
\Ad_{\widetilde u_\sigma}(E_{b+rN,a}) = E_{b+rN,a}-E_{b+rN,N+1}.
\]
Therefore,
\[
(E_{b+rN,a},0)\in \mathfrak h, \qquad \Ad_{\widetilde u_\sigma}(E_{b+rN,a},0)\in L,
\]
and hence,
\[
E_{b+rN,N+1}\in L.
\]

Combining the flag directions and the vector directions, we conclude that every non-positive weight vector of $\mathfrak{gl}_N(\bF)\ltimes \bF^N$ lies in $L$. Therefore, all $\lambda$-weights on the normal space are strictly positive.

Finally, let
\[
\widetilde{\GL}_N(\bF)\rightarrow \GL_N(\bF)
\]
be the Kac--Moody central extension corresponding to the determinant line bundle, and
\[
\pi:\overset{\circ}{\cP}_{\det}\to \Fl_N
\]
the associated punctured line bundle. Recall that, as in Section~\ref{sec:type-a-conventions}, the doubled extension $\hat\GL_N(\bF)$ corresponds to $\cP_{\det}^{\otimes2}$. Since $\lambda$ lies in
\[
T_N\times \mathbb G_m^{\mathrm{rot}},
\]
it lifts to an action on $\overset{\circ}{\cP}_{\det}$. If the induced action on the fiber $\pi^{-1}(x_{\tilde w})$ is given by the character $z\mapsto z^c$, then composing the lift with the inverse central cocharacter $z\mapsto z^{-c}$ produces a new lift, still denoted $\widetilde\lambda$, whose action on $\pi^{-1}(x_{\tilde w})$ is trivial. This central modification does not change the induced action on $\Fl_N$, hence does not change any of the above positivity properties. Furthermore, the new lift also fixes the fiber of $\cP_{\det}^{\otimes2}$ at $x_{\tilde w}$, and hence, we also obtain a lift for the doubled extension and its central character $\kappa$.
\end{proof}

\subsection{Transversal slices}\label{14.2}
\label{sec:global-transversal-slices} For a relevant index $\tilde w'$, we denote by $\IC_{\tilde w'}$ the normalized IC sheaf of the rank-one twisted equivariant local system on $\BO_{\tilde w'}$ in the $\Fl_N$-model, with the conventions of Section~\ref{sec:type-a-conventions}. Let $x=x_{\tilde w}$ be the representative of a relevant orbit $\BO_{\tilde w}$ in the support of $\IC_{\tilde w'}$, and let $\lambda$ be the cocharacter of Proposition~\ref{prop 9.5.1}. In order to prove the pointwise purity, we construct a transversal slice of finite type to the orbit of $x$, which is contracted to a point over $x$, and then we compare the $*$-stalk and $!$-stalk of $\IC_{\tilde w'}$ at $x$ with those of IC sheaves on this slice at the contraction point, up to normalized shifts. Namely, for $M<N$, we work on the global model and use a correspondence of generic Hecke modifications, and for $M=N$, we work on finite type truncations of the local model and use the action of a finite type quotient of $I_N$. In both cases, the slice is provided by the general Lemma~\ref{lem:finite-type-transversal} below. We refer readers to \cite[\S6.3]{[CY]} for contracting transversal slices to Iwahori orbits in loop spaces of symmetric varieties.

In Sections~\ref{14.2} and~\ref{sec:pointwise-purity-sub}, we denote by $k$ an algebraic closure of the finite field of Section~\ref{sec:type-a-conventions}, and we assume that the twisting and monodromy have finite order. Also, we will need the normalized shifts $\langle n\rangle=[n](n/2)$ of that section, and recall that, for a smooth map $f$ of relative dimension $r$, the normalized pullback $f^*\langle r\rangle$ preserves perverse sheaves and purity of weight zero. Here, by an IC complex, we mean an IC sheaf up to a shift $\langle n\rangle$, and note that these shifts preserve purity of weight zero and pointwise purity.

\subsubsection{A lemma about equivariant transversal slices}\label{sec:transversal-correspondences}
In this subsection, we prove a general lemma about contracting transversal slices on algebraic stacks, which will be applied to the global model for $M<N$ and to the local model for $M=N$ below. We say that a $\BG_m$-action on a scheme $P$ contracts $P$ to a point $p$ if the action map extends to a morphism $\BA^1\times P\to P$, whose restriction to $\{0\}\times P$ is constant with value $p$. Note that, in this case, every invariant open subset of $P$ containing $p$ is equal to $P$. Given a point $a$ of an algebraic stack $\mathcal B$ with a $\BG_m$-action, an equivariant structure on $a$ consists of isomorphisms $s\cdot a\simeq a$ compatible with the multiplication in $\BG_m$. We denote by $T_a\mathcal B$ the vector space of first-order deformations of $a$ up to isomorphism, and note that an equivariant structure on $a$ induces an action of $\BG_m$ on it.

\begin{lem}\label{lem:finite-type-transversal}
Let $\BG_m$ act on algebraic stacks $\mathcal B$ and $\mathcal R$ locally of finite type, where $\mathcal B$ has affine diagonal, and let
\[
\mathcal B\xleftarrow{s_0}\mathcal R\xrightarrow{t_0}\mathcal B
\]
be equivariant, representable, smooth maps. Let $a\in\mathcal B(k)$ and $y_0\in\mathcal R(k)$ be $\BG_m$-invariant points with equivariant structures, and with equivariant identifications $s_0(y_0)\simeq a\simeq t_0(y_0)$. Assume that $\operatorname{Aut}_{\mathcal B}(a)$ is smooth. Let
\[
\mathcal R_a=\mathcal R\times_{t_0,\mathcal B}\{a\},\qquad D=\operatorname{im}\bigl(T_{y_0}\mathcal R_a\xrightarrow{ds_0}T_a\mathcal B\bigr),
\]
and assume that all weights of $\BG_m$ on $T_a\mathcal B/D$ are positive. Then, there are an affine scheme $P$ of finite type, contracted to a point $p$ by $\BG_m$, and an equivariant map $f\colon(P,p)\to(\mathcal B,a)$ such that the composite
\[
T_pP\xrightarrow{df_p}T_a\mathcal B\longrightarrow T_a\mathcal B/D
\]
is an isomorphism. Furthermore, there is an invariant open substack $\mathcal R_P^\circ$ of $P\times_{f,\mathcal B,s_0}\mathcal R$, containing the point over $(p,y_0)$, such that the projection $\pi$ and the map $\tau$ induced by $t_0$,
\[
P\xleftarrow{\pi}\mathcal R_P^\circ\xrightarrow{\tau}\mathcal B,
\]
are representable and smooth, and $\pi$ is surjective. If $\mathcal B$ is a scheme and $a$ has an invariant affine open neighborhood, we may take $f$ to be a closed immersion into this neighborhood.
\end{lem}
Roughly speaking, $P$ is a transversal slice to the directions $D$ at $a$, and it is contracted to the point $p$ over $a$. Note that $f$ need not be an immersion, and the smooth correspondence $P\xleftarrow{\pi}\mathcal R_P^\circ\xrightarrow{\tau}\mathcal B$ will be used to compare sheaves on $\mathcal B$ and on $P$, see Lemma~\ref{lem:transversal-IC-restriction}.
\begin{proof}
The proof consists of two steps. First, we construct a transversal slice to $D$ in an equivariant smooth chart of $\mathcal B$. Note that the equivariant structure at $a$ splits the extension of $\BG_m$ by $\operatorname{Aut}_{\mathcal B}(a)$, which is defined by the stabilizer of $a$ in $[\mathcal B/\BG_m]$. Since $\operatorname{Aut}_{\mathcal B}(a)$ is smooth, according to \cite[Theorem 4.1, p.~692]{[AHR]}, there is an equivariant smooth map $q\colon(U,u)\to(\mathcal B,a)$, with $U$ affine of finite type, which induces an isomorphism of the stabilizers of $u$ and $a$ in $\BG_m$, and in particular, $u$ is fixed. We replace $U$ by the connected component of $u$. In the scheme case of the last assertion, we instead take $U$ to be the connected component of $a$ in the given affine open neighborhood, and $q$ to be the inclusion. Let $D_U=(dq_u)^{-1}(D)$ and $d=\dim D_U$. Since $q$ is smooth, the differential $dq_u$ is surjective, and we have
\[
T_uU/D_U\simeq T_a\mathcal B/D.
\]
We choose homogeneous functions $f_1,\ldots,f_d$ vanishing at $u$, whose differentials restrict to a basis of $D_U^*$, and set $P_0=V(f_1,\ldots,f_d)$. Then, $T_uP_0$ is a complement to $D_U$, and it has positive weights, since $T_uU/D_U\simeq T_a\mathcal B/D$.

Note that $P_0$ need not be contracted to $u$, and we need to cut it further as follows. Namely, let $A=k[P_0]=\bigoplus_nA_n$, with the convention $g(s\cdot y)=s^ng(y)$ for $g\in A_n$ and $s\in\BG_m$, and let $\mathfrak m$ be the ideal of $u$. We claim that the ideal
\[
J=(A_n\ (n<0),\ A_0\cap\mathfrak m)
\]
maps to zero in every $A/\mathfrak m^j$. Indeed, with this convention, the differential at $u$ of an element of $A_n$ is non-zero only on tangent vectors of weight $n$, and thus, $\mathfrak m/\mathfrak m^2$ is concentrated in positive degrees. Hence, each $\mathfrak m^i/\mathfrak m^{i+1}$, which is a quotient of $\operatorname{Sym}^i(\mathfrak m/\mathfrak m^2)$, is concentrated in degrees at least $i$, and $A/\mathfrak m^j$ is concentrated in non-negative degrees, with degree-zero part $k$. Now, by Krull's intersection theorem, we have $J_{\mathfrak m}=0$, see \cite[Tag 00IP]{[Stacks]}, and thus, $P=\operatorname{Spec}(A/J)$ agrees with $P_0$ near $u$. Furthermore, the grading on its coordinate ring is non-negative, with degree-zero part $k$, which defines the required contraction. Then, we take $p=u$ and $f=q|_P$.

Secondly, we construct the two smooth maps. Consider $\mathcal Z=U\times_{q,\mathcal B,s_0}\mathcal R$ with the projections
\[
\xymatrix{
& \mathcal Z \ar[dl]_{\sigma} \ar[dr]^{\operatorname{pr}_{\mathcal R}} & \\
U & & \mathcal R.
}
\]
We denote $v=t_0\circ\operatorname{pr}_{\mathcal R}\colon \mathcal Z\to\mathcal B$. Then, at the invariant point $z_0$ over $(u,y_0)$, we have
\begin{equation}\label{eq:slice-relative-tangent}
d\sigma\bigl(T_{z_0}(\mathcal Z\times_{v,\mathcal B}\{a\})\bigr)=D_U.
\end{equation}
Indeed, we have $\mathcal Z\times_{v,\mathcal B}\{a\}\simeq U\times_{q,\mathcal B,s_0}\mathcal R_a$, and a tangent vector $\xi\in T_uU$ lifts to this fiber product if and only if $dq_u(\xi)$ lies in the image $D$ of $T_{y_0}\mathcal R_a$.

Let $U^\circ=U\setminus\operatorname{Supp}(J)$, where the support of $J$ on $P_0$ is viewed as a closed subset of $U$. Note that, over $U^\circ$, we have $P=P_0$, and hence, the substack $P\times_U \mathcal Z=P\times_{f,\mathcal B,s_0}\mathcal R$ of $\mathcal Z$ is defined there by $f_1\circ\sigma=\cdots=f_d\circ\sigma=0$. We define $\mathcal R_P^\circ$ to be the open locus in $P\times_U\mathcal Z$, where the differentials of these functions are linearly independent on the tangent spaces of the fibers of $v$. This locus is invariant and contains $z_0$ by \eqref{eq:slice-relative-tangent}. Here, we denote by $\pi$ and $\tau$ the restrictions of $\sigma$ and $v$ to this open locus, respectively.

Finally, the map $\pi$ is smooth, since it is the restriction of the base change of $s_0$ to an open substack, and its image is an invariant open subset of $P$ containing $p$, which equals $P$ by contraction. On the other hand, the map $v$ is smooth, since $\operatorname{pr}_{\mathcal R}$ is a base change of $q$ and $t_0$ is smooth. Hence, by the relative Jacobian criterion, $\tau$ is smooth.
\end{proof}

\subsubsection{Restriction of IC sheaves}
Now, let us use the smooth correspondence of Lemma~\ref{lem:finite-type-transversal} to compare the IC complexes on $\mathcal B$ and on $P$.

\begin{lem}\label{lem:transversal-IC-restriction}
In the situation of Lemma~\ref{lem:finite-type-transversal}, let $K$ be an IC sheaf on $\mathcal B$ such that
\[
s_0^*K\simeq t_0^*K\otimes\mathcal L
\]
for a rank-one lisse sheaf $\mathcal L$ on $\mathcal R$. Then, $f^*K$ and $f^!K$ are IC complexes on $P$, and the $*$-stalk of the first and the $!$-stalk of the second at $p$ are $a^*K$ and $a^!K$, respectively. Here, $a$ also denotes the map $\Spec k\to\mathcal B$. If the data are defined over a finite field and $K$ and $\mathcal L$ are pure of weight zero, then so are $f^*K$ and $f^!K$.
\end{lem}
\begin{proof}
Since $f\circ\pi\simeq s_0\circ\operatorname{pr}_{\mathcal R}$ on $\mathcal R_P^\circ$, the given isomorphism implies
\begin{equation}\label{eq:transversal-IC-pullback}
\pi^*f^*K\simeq\tau^*K\otimes\operatorname{pr}_{\mathcal R}^*\mathcal L.
\end{equation}
Note that the right-hand side is an IC complex, since $\tau$ is smooth. Since $\pi$ is smooth and surjective, $f^*K$ is an IC complex by smooth descent, see \cite[Section 4.2]{[BBD]}. Applying the same argument to $\BD K$ and then Verdier duality gives the assertion for $f^!K$. Here, we use the isomorphism $s_0^*\BD K\simeq t_0^*\BD K\otimes\mathcal L^{-1}$, which holds since $s_0$ and $t_0$ have the same relative dimension.

Since purity can be checked after a smooth surjective pullback, the assertion about weights follows from \eqref{eq:transversal-IC-pullback} and the same argument for $\BD K$.
\end{proof}

\subsubsection{Generic Hecke modifications}\label{sec:generic-Hecke-modifications}
Now, we construct a correspondence on the global models of Section~\ref{sec:global-models}, to which Lemma~\ref{lem:finite-type-transversal} applies. Roughly speaking, the correspondence modifies the objects of the global model at a varying divisor away from $c$ by rational $H$-valued functions, and its tangent image at the point $a$ fixed below is the space of locally trivial deformations, see Proposition~\ref{prop:general-H-Hecke-correspondence}. We will work in the following general setting.

Let $H\subset G$ be a smooth affine subgroup. We fix an affine $G$-variety of finite type containing $H\backslash G$ as an open orbit, and consider its global model as in Section~\ref{sec:global-models}, with fixed twisting data. Let $C=\BP^1$ and $c=0$, and let $z=t^{-1}$ be the coordinate on $C-c$. We denote by $\infty$ the point $z=0$. Furthermore, we fix an $H^\omega$-bundle $\cP_H^\omega$ with a trivialization on $C-c$, where, in our application, we have $\cP_H^\omega=\omega^\rho\overset{Z_L}{\times}H$ as in Section~\ref{sec:global-models}.

Let $\cB\subset\cM'_{\infty\cdot c}$ be the substack, locally of finite type, defined by a fixed bound for the poles at $c$. Note that it has affine diagonal, since the isomorphisms of the underlying $G$-bundles form an affine scheme, and preserving the section and the flag at $c$ are closed conditions. Also, it is preserved by the generic Hecke modifications away from $c$. We fix a point $a\in\cB(k)$ with no defect on $C-c$, and then, its section defines an $H^\omega$-reduction of the underlying $G$-bundle on $C-c$. Here, we also fix an identification of this reduction with $\cP_H^\omega|_{C-c}$.

We assume that $\BG_m$ acts on $\cB$ through an action on $C$ fixing $c$ and $\infty$ and a compatible action on the twisted bundles, and that $a$ is $\BG_m$-fixed, compatibly with the identification above. In terms of the trivialization $\cP_H^\omega|_{C-c}\simeq(C-c)\times H$, we also assume that the action is given by
\[
\lambda(s)(z,h)=(s^mz,\alpha_s(h)),\qquad (z,h)\in(C-c)\times H,\quad s\in\BG_m,
\]
for some $m\in\BZ$ and an algebraic action $\alpha\colon\BG_m\to\operatorname{Aut}(H)$ by group automorphisms, and we denote $\alpha_s=\alpha(s)$.

We denote by $\fh_a$ the sheaf of infinitesimal automorphisms of $a$. Namely, for an open subset $V\subset C$, the space $\fh_a(V)$ consists of the sections over $V$ of the adjoint bundle of the underlying $G$-bundle which preserve the section over $V-c$ and, if $c\in V$, the flag at $c$. Note that, on $C-c$, it is the adjoint bundle of the $H^\omega$-reduction, and $\fh_a$ is a vector bundle, since it is a coherent subsheaf of a vector bundle on a smooth curve.

We denote by $D_a\subset T_a\cB$ the subspace of first-order deformations $a_\epsilon$ of $a$ over $k[\epsilon]/(\epsilon^2)$ which are locally trivial, i.e., isomorphic to $a\times\Spec k[\epsilon]/(\epsilon^2)$ on the members of an open covering of $C$, up to isomorphisms which reduce to the identity at $\epsilon=0$. By comparing the local trivializations on the overlaps, we obtain a \v{C}ech $1$-cocycle with values in $\fh_a$, and this construction identifies
\[
D_a\simeq H^1(C,\fh_a).
\]
Note that, if a locally closed substack $\cS\subset\cB$ classifies the families which are locally isomorphic to $a$, then $T_a\cS=D_a$, and in particular, this applies to the global stratum in Proposition~\ref{prop:global-model}.

In order to construct the modifications, we need the following coordinates on $H$. The following lemma is a consequence of \cite[Theorems 13.33 and 13.41, Corollaries 17.39 and 16.56]{[Milne]}.

\begin{lem}\label{lem:equivariant-group-coordinates}
Let $H^\circ$ be the identity component of a smooth affine group $H$. For a $\BG_m$-action on $H$ by automorphisms, there are non-negative integers $u,l$ and an equivariant open immersion
\[
\vartheta:\BA^u\times(\BG_m)^l\longrightarrow H^\circ,\qquad \vartheta(0,1)=1,
\]
where the action on $\BA^u$ is diagonal and that on $(\BG_m)^l$ is trivial.
\end{lem}

\begin{proof}
We denote $\widetilde H=H^\circ\rtimes\BG_m$, and choose a maximal torus $T_{\widetilde H}$ of $\widetilde H$ containing the subgroup $\{1\}\rtimes\BG_m$. We denote by $C_{\widetilde H}(T_{\widetilde H})$ the centralizer of $T_{\widetilde H}$ in $\widetilde H$, and choose a cocharacter $\mu:\BG_m\to T_{\widetilde H}$ such that $C_{\widetilde H}(\mu(\BG_m))=C_{\widetilde H}(T_{\widetilde H})$. We denote
\[
U_\pm=\{g\in\widetilde H\mid\lim_{s\to0}\mu(s)^{\pm1}g\mu(s)^{\mp1}=1\}.
\]
Note that these are smooth connected unipotent subgroups, and by \cite[Theorem 13.33]{[Milne]}, the multiplication map induces an open immersion
\[
U_-\times C_{\widetilde H}(T_{\widetilde H})\times U_+\longrightarrow \widetilde H.
\]
According to \cite[Theorem 13.41]{[Milne]}, we can take $T_{\widetilde H}$-equivariant affine coordinates on $U_\pm$. Also, as a variety, $C_{\widetilde H}(T_{\widetilde H})$ is a product of $T_{\widetilde H}$ and an affine space, and $\{1\}\rtimes\BG_m$ acts trivially on it, see \cite[Corollaries 17.39 and 16.56]{[Milne]}. Finally, by taking the fiber over $1$ of the projection $\widetilde H\to\BG_m$, we obtain the required coordinates.
\end{proof}

\begin{prop}\label{prop:general-H-Hecke-correspondence}
We can choose a connected algebraic stack $\cR$ of finite type
with a $\BG_m$-action and morphisms
\[
\cB\xleftarrow{s_0}\cR\xrightarrow{t_0}\cB.
\]
Both $s_0$ and $t_0$ are $\BG_m$-equivariant, representable, and smooth.

The identity modification of $a$ gives a $\BG_m$-fixed point
$y_0\in\cR(k)$. We have equivariant isomorphisms
\[
s_0(y_0)\simeq a\simeq t_0(y_0).
\]
On the $t_0$-fiber over $a$, the image of the differential of $s_0$
at $y_0$ is exactly $D_a$. In other words,
\[
\operatorname{im}\Bigl(
T_{y_0}\bigl(\cR\underset{t_0,\cB}{\times}\{a\}\bigr)
\xrightarrow{ds_0}T_a\cB
\Bigr)=D_a.
\]
\end{prop}

\begin{proof}
First, we construct rational $H$-valued functions whose first-order variations produce all the deformations in $D_a$. Namely, we choose $n\geq1$ such that $H^1(C,\fh_a(n\infty))=0$, and then, the boundary map
\begin{equation}\label{eq:transversal-principal-parts}
\partial\colon H^0\bigl(C,\fh_a(n\infty)/\fh_a\bigr)\twoheadrightarrow H^1(C,\fh_a)
\end{equation}
is surjective. Since $a$ has no defect near $\infty$, the trivialization identifies $\fh_a$ with $\Lie(H)\otimes\cO_C$ there. Hence, we have
\[
H^0\bigl(C,\fh_a(n\infty)/\fh_a\bigr)\simeq\bigoplus_{j=1}^n\Lie(H)\,z^{-j}.
\]

Let $\vartheta$ be as in Lemma~\ref{lem:equivariant-group-coordinates}, applied to the action $\alpha$, and let $r_1,\ldots,r_u$ be the weights of $\BG_m$ on the coordinates $x_1,\ldots,x_u$ of $\BA^u$. We denote by $y_1,\ldots,y_l$ the coordinates on $(\BG_m)^l$, and let $\BA$ be the affine space with coordinates $a_{ij}$ and $b_{qj}$, for $1\leq i\leq u$, $1\leq q\leq l$ and $1\leq j\leq n$. For $\xi\in\BA$, we define
\[
x_i(z)=\sum_{j=1}^na_{ij}z^{-j},\qquad y_q(z)=1+\sum_{j=1}^nb_{qj}z^{-j},\qquad h_\xi(z)=\vartheta\bigl(x(z);y(z)\bigr).
\]
Here, $h_\xi$ is a rational map from $C-c$ to $H^\circ$, and since $d\vartheta_{(0,1)}$ is an isomorphism onto $\Lie(H)$, the map
\[
T_0\BA\longrightarrow H^0\bigl(C,\fh_a(n\infty)/\fh_a\bigr),\qquad \xi\longmapsto\xi(z):=d\vartheta_{(0,1)}\Bigl(\sum_{j=1}^n(a_{\bullet j},b_{\bullet j})\,z^{-j}\Bigr),
\]
is an isomorphism. Note that $h_0=1$ and $h_{\epsilon\xi}(z)=1+\epsilon\xi(z)$ for $\epsilon^2=0$. Now, we let $\BG_m$ act on $\BA$ with weights $r_i+jm$ on $a_{ij}$ and $jm$ on $b_{qj}$. Then, $x_i$ and $y_q$ transform as the coordinates of $\BA^u\times(\BG_m)^l$, and by the equivariance of $\vartheta$, we have
\[
h_{\lambda(s)\xi}(s^mz)=\alpha_s\bigl(h_\xi(z)\bigr).
\]

Secondly, we choose a divisor $E$ outside which $h_\xi$ and $h_\xi^{-1}$ are regular. Namely, we set $p_q(z):=z^ny_q(z)=z^n+\sum_{j=1}^nb_{qj}z^{n-j}$, and fix a faithful representation $H\hookrightarrow\GL(V)$. Note that the matrix entries of $\vartheta(x,y)$ and of $\vartheta(x,y)^{-1}$ are polynomials in $x_i$, $y_q$ and $y_q^{-1}=z^n/p_q(z)$, and hence, the matrix entries of $h_\xi$ and $h_\xi^{-1}$ are rational functions whose denominators are products of powers of $z$ and of the $p_q(z)$. Then, we choose an integer $N_E\geq1$ such that the polynomial
\[
e(z):=\Bigl(z\prod_{q=1}^lp_q(z)\Bigr)^{N_E}\in k[\BA][z]
\]
clears these denominators, i.e., $e(z)h_\xi(z)$ and $e(z)h_\xi(z)^{-1}$ have matrix entries in $k[\BA][z]$. We denote by $E\subset(C-c)\times\BA$ the divisor $e(z)=0$, and by $E_\xi$ its fiber over $\xi\in\BA$. Since $e(z)$ is monic in $z$ of degree $N_E(1+nl)$, the divisor $E$ is finite and flat over $\BA$, and we have $E_0=N_E(1+nl)\cdot\infty$. By construction, the map $h_\xi$ is a regular map to $H^\circ$ outside $E$, and since $e_{\lambda(s)\xi}(s^mz)=s^{mN_E(1+nl)}e_\xi(z)$, the divisor $E$ is $\BG_m$-invariant. Also, we denote by $3E$ the subscheme defined by $e(z)^3=0$.

In order to keep the parameters of finite type, we consider trivializations along the finite subscheme $3E$ instead of formal trivializations. We denote by $\operatorname{Res}_{3E/\BA}H$ the group scheme over $\BA$ of $H$-valued maps on $3E$. Namely, for a scheme $\BA'$ over $\BA$, its $\BA'$-points are the morphisms $3E\times_\BA\BA'\to H$, with pointwise multiplication. Note that it is the Weil restriction of $H\times3E$ along $3E\to\BA$, see \cite[Tag 05Y8]{[Stacks]}, and it is smooth over $\BA$, since $H$ is smooth and $3E\to\BA$ is finite flat.

Thirdly, we construct the correspondence. Let $(\BA\times\cB)^\circ\subset\BA\times\cB$ be the open substack of pairs $(\xi,a')$ such that $a'$ has no defect along $E_\xi$. We denote by $\widetilde\cR$ the stack classifying $\xi\in\BA$, an object $a'\in\cB$ with $(\xi,a')\in(\BA\times\cB)^\circ$, and an isomorphism $\varphi$ of the $H^\omega$-reduction of $a'$ on $3E_\xi$ with $\cP_H^\omega|_{3E_\xi}$, and by $\operatorname{pr}_{\BA}$ and $s_0$ its projections to $\BA$ and $\cB$, respectively. Note that, by using the trivialization of $\cP_H^\omega$ on $C-c$, the choices of $\varphi$ form a torsor under $\operatorname{Res}_{3E/\BA}H$, and hence, $(\operatorname{pr}_{\BA},s_0)\colon\widetilde\cR\to(\BA\times\cB)^\circ$ is representable, smooth and surjective.

Now, let us define the map $t_0$. Let $(\xi,a',\varphi)$ be a point of $\widetilde\cR$ over a test scheme. Since $H$ is smooth, locally on the test scheme, $\varphi$ extends to an isomorphism $\widehat\varphi$ over the formal neighborhood of $E_\xi$. Then, we define $t_0(\xi,a',\varphi)$ by gluing the restriction of $a'$ to the complement of $E_\xi$ with the trivial $H^\omega$-bundle on the formal neighborhood of $E_\xi$, via $h_\xi\circ\widehat\varphi$ on the punctured formal neighborhood, see \cite{[BL]} and \cite[Section 1.4]{[Zhu]}. Since the resulting object has the same poles at $c$, it lies in $\cB$.

Note that the target does not depend on the choice of $\widehat\varphi$. Indeed, two extensions agree on $3E_\xi$, and thus, they differ by an element $g$ with $g=1+e(z)^3X$ in the faithful representation, where $X$ is a matrix which is regular on the formal neighborhood of $E_\xi$. By the choice of $e(z)$,
\[
h_\xi gh_\xi^{-1}-1=e(z)\bigl(e(z)h_\xi\bigr)X\bigl(e(z)h_\xi^{-1}\bigr)
\]
is regular and divisible by $e(z)$, and the same holds with $g$ replaced by $g^{-1}$. Since $H$ is closed in $\GL(V)$ and $e(z)$ is a non-zero-divisor, $h_\xi gh_\xi^{-1}$ is an automorphism of the trivial $H^\omega$-bundle on the formal neighborhood, which identifies the two glued objects canonically and compatibly with compositions. Thus, $t_0\colon\widetilde\cR\to\cB$ is well-defined, and it is $\BG_m$-equivariant by the identity for $h_\xi$ above. Furthermore, both $s_0$ and $t_0$ are representable. Namely, for $s_0$, it follows from the torsor description, and for $t_0$, it follows from the fact that an automorphism of $a'$ which induces the identity on the glued object is the identity away from $E_\xi$, hence everywhere.

Next, we check the smoothness near the identity modification. We denote by $y_0\in\widetilde\cR(k)$ the point consisting of $\xi=0$, the object $a$, and the restriction to $3E_0$ of the fixed identification of the reduction of $a$, where this restriction is possible since $E_0$ is supported at $\infty\in C-c$. Since $h_0=1$, we have canonical equivariant isomorphisms $s_0(y_0)\simeq a\simeq t_0(y_0)$, and $y_0$ is $\BG_m$-fixed.

Note that, since $(\operatorname{pr}_{\BA},s_0)$ is smooth and $(\BA\times\cB)^\circ$ is flat over $\BA$, the stack $\widetilde\cR$ is flat over $\BA$. On the other hand, since $h_0=1$, the maps $(\operatorname{pr}_{\BA},s_0)$ and $(\operatorname{pr}_{\BA},t_0)$ agree over $0\in\BA$, and in particular, the restriction of $(\operatorname{pr}_{\BA},t_0)$ to the fiber $\widetilde\cR_0$ is smooth over $\{0\}\times\cB$. Thus, by the fiberwise criteria of flatness and of smoothness, $(\operatorname{pr}_{\BA},t_0)$ is smooth at every point of $\widetilde\cR_0$, see \cite[Tags 039C and 01V4]{[Stacks]}. Since the projection $\BA\times\cB\to\cB$ is smooth, it follows that $s_0$ and $t_0$ are smooth on an open neighborhood of $\widetilde\cR_0$, and at the points of $\widetilde\cR_0$, they have the same relative dimension. Now, we choose a connected quasi-compact open neighborhood of $y_0$ in this locus, on which the relative dimensions are constant, and let $\cR$ be the union of its $\BG_m$-translates. Then, $\cR$ is a connected, $\BG_m$-invariant open substack of finite type, on which $s_0$ and $t_0$ are smooth of the same constant relative dimension.

Finally, let us compute the tangent image. Given $\xi\in T_0\BA$, by modifying the trivial deformation of $a$ by $h_{\epsilon\xi}^{-1}=1-\epsilon\xi(z)$ along $E_{\epsilon\xi}$, we obtain a first-order deformation $a'_\epsilon$ of $a$, and $(\epsilon\xi,a'_\epsilon)$, with the induced trivialization on $3E_{\epsilon\xi}$, is a tangent vector to $\cR\times_{t_0,\cB}\{a\}$ at $y_0$. By the definition of the boundary map, the class of $a'_\epsilon$ is $-\partial(\xi(z))$, and hence, by \eqref{eq:transversal-principal-parts}, the tangent image contains $D_a$. Conversely, a first-order modification whose target is the trivial deformation of $a$ changes only the gluing near $E_{\epsilon\xi}\subset C-c$, where $a$ has no defect, and hence, its source deformation is locally trivial, i.e., it belongs to $D_a$. Thus, the tangent image is $D_a$.
\end{proof}

Now, we show that a generic Hecke equivariant object $K$ on $\cB$ satisfies
\[
s_0^*K\simeq t_0^*K\otimes\mathcal L,
\]
where $\mathcal L$ is the prescribed rank-one character sheaf. Let $r_E=\deg(E/\BA)=N_E(1+nl)$, and let $\rho\colon\cR'\to\cR$ be the pullback of the finite flat cover of $\BA$ which orders the $r_E$ zeros of $e(z)$, counted with multiplicity. Over $\cR'$, the modification at the ordered zeros defines a map to the universal Hecke stack \eqref{univ heck}, and hence, \eqref{univ heck iso} provides the comparison after the $!$-pullback along $\rho$. Note that this comparison is defined over the Ran space, and thus, it is invariant under the symmetric group $S_{r_E}$ permuting the zeros, including the locus where zeros coincide. By Verdier duality, it turns into a comparison of the $*$-pullbacks along $\rho$, which descends to $\cR$ by the projection formula, since $(\rho_*\bar\BQ_\ell)^{S_{r_E}}=\bar\BQ_\ell$. Then, by applying Verdier duality again, we obtain the $!$-pullback comparison on $\cR$, which is the displayed isomorphism without a shift, since $s_0$ and $t_0$ are smooth of the same relative dimension. Finally, we include the twisting by using the corresponding torsors, as in Remark~\ref{rem 9.5.3}.

\subsubsection{The global model}
Now, we go back to the case $M<N$, and apply the above constructions to the global stratum corresponding to the orbit of $x$. Here, we transport $\pi_\Fl$ and the action $\lambda$ by the change of Borel subgroup and the conjugation in Section~\ref{sec:type-a-conventions}. Let
\[
x_G=(t^0,x)\in\Fl_M\times\Fl_N,\qquad a=\pi_\Fl(x_G)\in\cM'_{\infty\cdot c}(k).
\]
Here, $t^0\in\Fl_M$ denotes the base point, and we denote by $\mathcal S=\overset{\circ}{\cM}_{\Fl,\tilde w}$ the global stratum containing $a$, and by $\sO_G$ the $H(\bF)$-orbit of $x_G$, so that $\pi_\Fl^{-1}(\mathcal S)=\sO_G$. Note that, since $\mathcal S$ is a good-elsewhere stratum, $a$ has no defect on $C-c$, as required in Section~\ref{sec:generic-Hecke-modifications}. Furthermore, the action $\lambda$ lifts to $\Fl_M\times\Fl_N$ by the same loop rotation and the upper-left $M\times M$ part of its constant torus factor, and the lift fixes $x_G$ and extends to the global model, where, on the trivialized $H^\omega$-bundle over $C-c$, it has the form used in Section~\ref{sec:generic-Hecke-modifications}. Finally, we denote by $K$ the global IC sheaf corresponding to $\IC_{\tilde w'}$ under Theorem~\ref{prop:global-local-strata}, and by $\mathcal B\subset\cM'_{\infty\cdot c}$ the substack of Section~\ref{sec:generic-Hecke-modifications}, for a pole bound at $c$ such that $\mathcal B$ contains the support of $K$.

\begin{prop}\label{prop 9.5.2}\label{prop:global-model}
There are an equivariant smooth map $q\colon(U,u)\to(\mathcal B,a)$, with $U$ affine of finite type (so that $q(U)$ is of finite type), and an invariant closed subscheme $i\colon P\hookrightarrow U$ contracted to $u$, such that
\[
T_uU=T_uP\oplus T_u\mathcal S_U,\qquad \mathcal S_U:=U\times_{\cM'_{\infty\cdot c}}\mathcal S.
\]
Let $K_U$ be the normalized smooth pullback of $K$ to $U$. Then, $i^*K_U$ and $i^!K_U$ are IC complexes on $P$. The $*$-stalk of the first and the $!$-stalk of the second at $u$ are the $*$-stalk and $!$-stalk of $\IC_{\tilde w'}$ at $x$, respectively, up to shifts $\langle n\rangle$.
\end{prop}
\begin{proof}
Note that, for $K$, the hypothesis of Lemma~\ref{lem:transversal-IC-restriction} is the generic Hecke equivariance shown after Proposition~\ref{prop:general-H-Hecke-correspondence}. Thus, in order to apply Lemma~\ref{lem:finite-type-transversal} to the correspondence of Proposition~\ref{prop:general-H-Hecke-correspondence}, we only need to check that $\operatorname{Aut}(a)$ is smooth and that the weights on $T_a\mathcal B/D_a$ are positive.

First of all, in type A, the group $H$ is the group of invertible elements of the affine subspace $1+\Lie(H)$, which is closed under multiplication, see \eqref{emd gl}, and the flag conditions at $c$ are also linear in the matrix entries. Thus, $\operatorname{Aut}(a)$ is an open subscheme of $1+H^0(C,\mathfrak h_a)$, and hence, it is smooth. Furthermore, the global stratum $\mathcal S$ classifies the families which are locally isomorphic to $a$, and it is smooth at $a$, since its gluing obstructions lie in $H^2(C,\mathfrak h_a)=0$, and we have $T_a\mathcal S=D_a$, see \cite[Lecture 2]{[Hei]}.

Secondly, for the positivity, recall that $\pi_\Fl$ is formally smooth, and since $\pi_\Fl^{-1}(\mathcal S)=\sO_G$, we have an equivariant isomorphism
\[
T_{x_G}(\Fl_M\times\Fl_N)/T_{x_G}\sO_G\simeq T_a\cM'_{\infty\cdot c}/T_a\mathcal S.
\]
Since the subgroup $\GL_M(\bF)\subset H(\bF)$ acts transitively on $\Fl_M$, and the fiber of $\sO_G$ over $t^0$ is $\BO_{\tilde w}$, the left-hand side is identified with $T_x\Fl_N/T_x\BO_{\tilde w}$, whose weights are positive by Proposition~\ref{prop 9.5.1}, and thus, the weights on $T_a\mathcal B/D_a$ are also positive.

Therefore, we can apply the construction in the proof of Lemma~\ref{lem:finite-type-transversal}, and we obtain $q$, $P$, $\pi$ and $\tau$. Since $\mathcal S$ is smooth at $a$ and $\mathcal S_U\to\mathcal S$ is a base change of $q$, the stack $\mathcal S_U$ is smooth at $u$, and we have $T_u\mathcal S_U=(dq_u)^{-1}(T_a\mathcal S)=(dq_u)^{-1}(D_a)$, which implies $T_uU=T_uP\oplus T_u\mathcal S_U$. Note that the image $q(U)$ is an invariant quasi-compact open substack of $\mathcal B$, and we restrict $\mathcal R_P^\circ$ to $\tau^{-1}(q(U))$, on which the map $\pi$ remains surjective by contraction. Since $f=q\circ i$, the assertion about IC complexes follows from Lemma~\ref{lem:transversal-IC-restriction}. Finally, according to Theorem~\ref{prop:global-local-strata}, the $*$-stalk and $!$-stalk of $K$ at $a$ are identified with the normalized $*$-stalk and $!$-stalk of $\IC_{\tilde w'}$ at $x$, respectively, and we include the twisting as in Remark~\ref{rem 9.5.3}.
\end{proof}

\subsubsection{The case \texorpdfstring{$M=N$}{M=N}}\label{sec:slices-balanced}
In this case, we will apply Lemmas~\ref{lem:finite-type-transversal} and~\ref{lem:transversal-IC-restriction} to finite type truncations of the local model, with the action of a finite type quotient of $I_N$.

We go back to the point $x=x_{\tilde w}$ and the cocharacter $\lambda$ of Proposition~\ref{prop 9.5.1}, where $x=(F_\bullet,v)\in\Fl_N\times\bF^N$, with $v\neq0$. We choose an $I_N$-stable Schubert variety $Z\subset\Fl_N$ and an integer $A\geq0$ such that $x$ and the support of $\IC_{\tilde w'}$ are contained in $Z\times t^{-A}\bO^N$. For $n\geq1$, we denote
\[
Y=Z\times t^{-A}\bO^N,\qquad X_n=Z\times\bigl(t^{-A}\bO^N/t^n\bO^N\bigr),
\]
and let $\pi_n\colon Y\to X_n$ be the reduction modulo $t^n\bO^N$. Also, we set $x_n=\pi_n(x)$ and $S_n=I_N\cdot x_n$. For $n'\geq n$, let $p_{n',n}\colon X_{n'}\to X_n$ be the reduction modulo $t^n\bO^N$, which is smooth of relative dimension $N(n'-n)$. Note that the group $I_N$ and the action $\lambda$ preserve $Z$, $t^{-A}\bO^N$ and $t^n\bO^N$, and thus, they act on $X_n$, and $\lambda$ fixes $x_n$.

\begin{prop}\label{prop:balanced-finite-transversal}
For all sufficiently large $n$, the following hold.
\begin{enumerate}
\item The scheme $X_n$ carries a normalized IC sheaf $\IC_n$ whose normalized smooth pullbacks give the local IC object $\IC_{\tilde w'}$ near $x$. For $n'\geq n$,
\begin{equation}\label{eq:balanced-IC-transition}
\IC_{n'}\simeq p_{n',n}^*\IC_n\langle N(n'-n)\rangle.
\end{equation}
\item There are a $\lambda$-invariant affine open neighborhood $U$ of $x_n$ in $X_n$ and a $\lambda$-invariant closed subscheme $i\colon P\hookrightarrow U$ contracted to $x_n$, such that $T_{x_n}U=T_{x_n}P\oplus T_{x_n}S_n$.
\item The complexes $i^*\IC_n$ and $i^!\IC_n$ are IC complexes on $P$.
\end{enumerate}
\end{prop}

\begin{proof}
Note that the transition functors in the vector direction are normalized pullbacks, and hence, Assertion (1), including \eqref{eq:balanced-IC-transition}, follows from \cite[Section 3.2]{[BFGT]}.

For (2) and (3), we first replace $I_N$ by a finite type quotient. Namely, for $e\geq1$, let $I_N^e=\ker(\GL_N(\bO)\to\GL_N(\bO/t^e\bO))$, and we choose $e_0\geq1$ such that $I_N^{e_0}$ acts trivially on $Z$. Then, for $e\geq\max(e_0,n+A)$, the group $I_N^e$ acts trivially on $X_n$, and thus, $I_N$ acts on $X_n$ through $\Gamma_e=I_N/I_N^e$, which is the group of invertible elements of the finite-dimensional algebra $\mathfrak i_N/t^e\Mat_N(\bO)$, where $\mathfrak i_N=\Lie(I_N)$, and hence, it is smooth and connected. We also take $e$ large enough such that the character sheaves factor through $\Gamma_e$. Note that the action $\lambda$ normalizes $I_N$ and $I_N^e$, and thus, it acts on $\Gamma_e$ by conjugation.

Next, we compare the orbit of $x$ with its truncations. We choose a row vector $\ell\in\bF^{1\times N}$ with $\ell v=1$, and let $n$ be large enough such that $1+\Delta v\,\ell\in I_N^{e_0}$ for every $\Delta v\in t^n\bO^N$. Then, $1+\Delta v\,\ell$ fixes $F_\bullet$ and sends $v$ to $v+\Delta v$, and hence, we have $x+(0,\Delta v)\in I_N\cdot x$. By applying this to $h^{-1}\Delta v$, for $h\in I_N$, we also obtain $h\cdot x+(0,\Delta v)\in I_N\cdot x$, which implies that the orbit of $x$ in $Y$ is $\pi_n^{-1}(S_n)$, and $p_{n',n}^{-1}(S_n)=S_{n'}$ for $n'\geq n$. Note that the stabilizer of $x_n$ in $\Gamma_e$ is defined by the linear conditions that $g$ preserves $F_\bullet$ and $(g-1)v\in t^n\bO^N$, together with the invertibility, and hence, it is smooth and connected, and the orbit map $\Gamma_e\to S_n$ is smooth. Since $\pi_n$ is the inverse limit of the smooth projections $p_{n',n}$, these equalities imply
\[
T_{x_n}X_n/T_{x_n}S_n\simeq T_xY/T_x(I_N\cdot x)\subset T_x(\Fl_N\times\bF^N)/T_x\BO_{\tilde w}.
\]
Here, the weights on the right-hand side are positive by Proposition~\ref{prop 9.5.1}.

Finally, we denote by $W$ the intersection of $Z$ with a $\lambda$-invariant affine Pl\"ucker chart around the $T_N$-fixed flag $F_\bullet$, so that $W\times(t^{-A}\bO^N/t^n\bO^N)$ is a $\lambda$-invariant affine open neighborhood of $x_n$, and we apply Lemma~\ref{lem:finite-type-transversal}, in the scheme case, to
\[
X_n\xleftarrow{s_0}\Gamma_e\times X_n\xrightarrow{t_0}X_n,\qquad s_0(g,y)=y,\quad t_0(g,y)=g\cdot y,
\]
with $a=x_n$ and $y_0=(1,x_n)$. Note that both maps are smooth of the same relative dimension, and the map $g\mapsto(g,g^{-1}x_n)$ identifies $\Gamma_e$ with the fiber of $t_0$ over $x_n$, and hence, we have $D=T_{x_n}S_n$ by the smoothness of the orbit map. Thus, we obtain (2), with $U=W\times(t^{-A}\bO^N/t^n\bO^N)$. Furthermore, by the twisted equivariance of $\IC_n$, we have $s_0^*\IC_n\simeq t_0^*\IC_n\otimes\mathcal L$, with $\mathcal L$ pulled back from $\Gamma_e$, and hence, Lemma~\ref{lem:transversal-IC-restriction} implies (3). Also, we include the torsors as in Remark~\ref{rem 9.5.3}.
\end{proof}

\subsubsection{Twisting and monodromy}\label{sec:slices-twisting}
Given an action $a_\lambda:\BG_m\times P\to P$, we say that a complex $K$ is weakly $\lambda$-equivariant if $a_\lambda^*K\simeq\operatorname{pr}_P^*K$, where $\operatorname{pr}_P:\BG_m\times P\to P$ is the projection.

\begin{rem}\label{rem 9.5.3}
We include the twisting and monodromy by using the punctured line bundle $\cP_{\det}^{\otimes2}$ and the Cartan torsors $\GL_r(\bF)/I_r^0\to\Fl_r$, where $I_r^0=\ker(I_r\to T_r)$, with the conventions of Section~\ref{sec:type-a-conventions}. On the global model, we also need the $T$-torsor of trivializations of the $T$-bundle induced by the $B$-reduction at $c$, and note that the proof of Theorem~\ref{prop:global-local-strata} applies with these torsors, since the generic Hecke modifications occur away from $c$. On $X_n$, we pull back the torsors from $Z$.

According to Proposition~\ref{prop 9.5.1}, $\lambda$ lifts to the punctured determinant line bundle, and hence to its square. Also, $\lambda$ lifts naturally to the Cartan torsors and to the torsor of trivializations of the $T$-bundle at $c$. For each torsor, we compose the lift with the inverse of the cocharacter of its structure torus by which it acts on the fiber over the fixed point of $P$. Then, the new lift acts trivially on this fiber and induces the same action on $P$, and we claim that the torsors become equivariantly trivial over $P$. Indeed, since $P$ is affine, bases of the fibers of the associated line bundles at the fixed point lift to invariant sections, which are nowhere zero by contraction. Thus, with the weight-zero normalization of the character sheaves, the IC comparisons in Propositions~\ref{prop:global-model} and~\ref{prop:balanced-finite-transversal} provide the IC complexes $K_P=i^*K_U$ and $K_P^!=i^!K_U$ on $P$, which are pure of weight zero, and in the case $M=N$, they are $i^*\IC_n$ and $i^!\IC_n$, respectively.

By the calculation in Proposition~\ref{prop 9.5.1} and intermediate extension, we have
\[
a_\lambda^*K_P\simeq\mathcal E\boxtimes K_P,
\]
where $a_\lambda:\BG_m\times P\to P$ is the action map and $\mathcal E$ is the rank-one local system on $\BG_m$ induced by the prescribed multiplicative characters. Note that, for finite-order twisting and monodromy, $\mathcal E$ has finite order, and hence, if we replace $\lambda(s)$ and its lift by $\lambda(s^\nu)$ and the corresponding lift, for a positive multiple $\nu$ of this order, then $\mathcal E$ becomes constant, while the positivity of the weights is preserved. Therefore, $K_P$ is weakly equivariant for the contracting action, and by Verdier duality, so is $K_P^!$.
\end{rem}

\subsection{Pointwise purity}\label{sec:pointwise-purity-sub}
In this subsection, we prove the pointwise purity of the IC sheaves. With the preparation above, it suffices to repeat the argument in \cite[Proposition 7]{[T1]} on the slices constructed in Section~\ref{14.2}. Recall that a Weil complex $K$ which is pure of weight zero is pointwise pure if, for every point $y$ defined over a finite field, $H^j(i_y^*K)$ and $H^j(i_y^!K)$ are pure of weight $j$ for all $j$, see \cite[Appendix A.8]{[CY]}. Here, for the local IC object $\IC_{\tilde w'}$, we consider its $*$-stalks and $!$-stalks on the global model if $M<N$, and on the finite type quotients $X_n$ if $M=N$, which are well-defined up to normalized shifts by Theorem~\ref{prop:global-local-strata} and \eqref{eq:balanced-IC-transition}, respectively. Note that these shifts preserve pointwise purity.

\begin{cor}\label{cor:contracting-slice-purity}
Let $P$ be a scheme of finite type with a $\BG_m$-action contracting $P$ to a point $p$, and let $K_P$ be a complex on $P$ which is pure of weight zero and weakly equivariant for this action, all defined over a finite field. Then, the $*$-stalk and the $!$-stalk of $K_P$ at $p$ are pure of weight zero.
\end{cor}
\begin{proof}
Since $K_P$ is weakly equivariant for a contracting action, the restriction to $p$ induces a Frobenius-compatible isomorphism $H(P,K_P)\simeq(K_P)_p$, see \cite[proof of Proposition 7]{[T1]} and \cite[Appendix A.8, proof of Proposition 113]{[CY]}. On the other hand, for each $j$, since $K_P$ is pure of weight zero, the weights of $H^j(P,K_P)$ are at least $j$, and the weights of $\mathcal H^j(K_P)_p$ are at most $j$, and hence, both are pure of weight $j$. The assertion for the $!$-stalk follows by applying the same argument to the Verdier dual.
\end{proof}

\begin{prop}\label{prop:pointwise-purity}
For a relevant orbit $\BO_{\tilde w'}$, with finite-order twisting and monodromy data defined over a finite field, the IC sheaf $\IC_{\tilde w'}$ is pointwise pure.
\end{prop}
\begin{proof}
First of all, note that the $*$- and $!$-restrictions of $\IC_{\tilde w'}$ to an irrelevant orbit vanish, since such an orbit carries no non-zero twisted equivariant local system. Thus, by the twisted equivariance, it is enough to consider the representatives $x=x_{\tilde w}$ of relevant orbits, since the rank-one character sheaves have weight zero, and hence, they do not change the weights of the $*$-stalks and $!$-stalks along an orbit. Also, by replacing the base field by a finite extension, we may assume that all the finite type data below are defined over it.

In the case $M<N$, we consider the slice of Proposition~\ref{prop:global-model}. According to Remark~\ref{rem 9.5.3}, the IC complexes $K_P$ and $K_P^!$ on $P$ are pure of weight zero and weakly equivariant for the contracting action, and hence, by Corollary~\ref{cor:contracting-slice-purity}, the $*$-stalk of $K_P$ and the $!$-stalk of $K_P^!$ at the invariant point are pure of weight zero. On the other hand, by Proposition~\ref{prop:global-model}, they are the $*$-stalk and $!$-stalk of $\IC_{\tilde w'}$ at $x$, up to shifts $\langle n\rangle$, which proves the assertion in this case.

In the case $M=N$, for a point $x=(F_\bullet,v)$ with $v\ne0$, we take $n$ sufficiently large and apply the same argument by using Proposition~\ref{prop:balanced-finite-transversal}. Finally, let us consider the points of $X_n$ with zero vector coordinate. Such a point is of the form $\pi_n(y)$, where $y=(F'_\bullet,v')\in Y$ with $v'\in t^n\bO^N$ and $v'\ne0$, defined over a finite extension, and we take $n'\geq n$ large enough such that the previous argument applies to $\pi_{n'}(y)$. Then, by \eqref{eq:balanced-IC-transition}, the $*$-stalk and $!$-stalk of $\IC_n$ at $\pi_n(y)$ are identified with those of $\IC_{n'}$ at $\pi_{n'}(y)$, up to normalized shifts, and hence, the purity at $\pi_{n'}(y)$ implies the purity at $\pi_n(y)$. This completes the proof.
\end{proof}

\appendix
\section{The IC comparison after a change of $G$ within the same isogeny class}
\label{sec:gaiotto-isogeny-proofs}

In this appendix, we prove Theorems~\ref{thm:isogeny-image-of-ic} and \ref{thm:isogeny-metaplectic-hecke-action}. The proofs are analogous to those of Theorems~\ref{thm:image-of-ic} and \ref{thm:metaplectic-hecke-action}, and we only explain the necessary modifications. We use the notation of Section~\ref{sec:gaiotto-isogenies}, and in particular, we have
\[
\widetilde G_1=\widetilde G/G_2,\qquad \widetilde G_2=\widetilde G/G_1,\qquad \widetilde G=\widetilde G_1\times_{T_0}\widetilde G_2.
\]
Here, the first and second components refer to the projections to $\widetilde G_1$ and $\widetilde G_2$, respectively. For a coweight $\lambda\in\widetilde\Lambda$, we denote its components by $\lambda=(\lambda^{(1)},\lambda^{(2)})$, where $\lambda^{(a)}$ denotes its image in the coweight lattice of $\widetilde G_a$. Note that both components have the same image in $\Lambda_0$, and two coweights with the same image in $\Lambda_0$ differ by an element of $\Lambda$.

Roughly speaking, the difference is as follows. For $\widetilde\eta\in\widetilde\Lambda$ over $\xi\in\Lambda_0$, the right multiplication by $t^{\widetilde\eta}$ identifies $G(\bF)/K$, where $K=G(\bF)\cap t^{\widetilde\eta}\widetilde G(\bO)t^{-\widetilde\eta}$, with the fiber of $\Gr_{\widetilde G}\to\Gr_{T_0}$ over $t^\xi$, and hence we have to check the arguments of Sections~\ref{sec:orbit-preliminaries} and~\ref{sec:irreducibles} with $G(\bO)$ replaced by $K$. If $\widetilde\eta$ can be chosen central, then $K=G(\bO)$, and this covers all the cases except one orthogonal case, in which $K$ involves the group $K_\eta$ defined below, whose root subgroups are described in \eqref{eq:enlarged-K-root-filtration}.

\subsection{The components and the orbit representatives}

\begin{lem}\label{lem:app-isogeny-components}
If $Z(G)$ is connected, then for every $\xi\in\Lambda_0$, there exists $\eta_0\in X_*(Z(\widetilde G)^0)$ mapping to $\xi$ under the map induced by $Z(\widetilde G)^0\hookrightarrow\widetilde G\to T_0$ in \eqref{eq:gaiotto-isogeny}. Multiplication by $t^{\eta_0}$ gives an isomorphism
\[
\begin{gathered}
\Gr_G\xrightarrow{\sim} \Gr_{\widetilde G}\times_{\Gr_{T_0}}\{t^\xi\},\\
gG(\bO)/G(\bO)\longmapsto t^{\eta_0}g\widetilde G(\bO)/\widetilde G(\bO).
\end{gathered}
\]
\end{lem}
\begin{proof}
We denote $Z=Z(\widetilde G)^0$. Since $\widetilde G/G$ is a torus, we have $\widetilde G=G\cdot Z=(G\times Z)/(G\cap Z)$ and $G\cap Z\subset Z(G)$. If $Z(G)$ is connected, then $Z(G)\subset Z$, and hence $1\to Z(G)\to Z\to T_0\to1$ is an exact sequence of tori, and $X_*(Z)\to\Lambda_0$ is surjective. Note that the fiber over $t^0$ is $G(\bF)\widetilde G(\bO)/\widetilde G(\bO)\simeq\Gr_G$, and the multiplication by the central element $t^{\eta_0}$ maps it isomorphically onto the fiber over $t^\xi$.
\end{proof}

In particular, if $Z(G)$ is connected, e.g., in type A and in the $\on F(4)$ case, the multiplication by $t^{\eta_0}$ in Lemma~\ref{lem:app-isogeny-components} identifies $\Gr_G$ with the fiber over $t^\xi$, together with the $H(\bF)$-orbits, their stabilizers and their intersections with semi-infinite orbits, and the pairings with geometric roots are the same. Under this map, a coweight $\mu\in\Lambda$ becomes $\mu+\eta_0$, so the twisting is computed with $\widetilde b(\mu+\eta_0,\nu)$, for $\nu\in\widetilde\Lambda$. Thus, the proofs of Theorems~\ref{thm:image-of-ic} and \ref{thm:metaplectic-hecke-action} apply with these coweights and with $\widetilde b$ in place of $b$, without any other change. Note that the same argument applies when $\widetilde G\simeq G\times T_0$, with $\eta_0$ taken in the factor $T_0$.

\begin{rem}\label{rem:app-isogeny-orthogonal}
Suppose that $\{G_1,G_2\}=\{\SO(2k),\SO(2l+1)\}$ with $k\geq2$ (for $k=1$, the center of $G$ is connected). We denote $Z=Z(\widetilde G)^0$. As above, $\widetilde G=(G\times Z)/(G\cap Z)$, where $G\cap Z$ is either $1$ or $\mu_2=Z(G)$. In the first case, $\widetilde G\simeq G\times T_0$. In the second case, we choose coordinates $Z\simeq\BG_m\times T'$ with $\mu_2\subset\BG_m\times\{1\}$. Then $(\SO(2k)\times\BG_m)/\mu_2\simeq\operatorname{GSO}(2k)$ via $(g,c)\mapsto cg$, and
\[
\widetilde G\simeq\operatorname{GSO}(2k)\times\SO(2l+1)\times T',
\]
with the first two factors interchanged if $G_1=\SO(2l+1)$. Here, $\operatorname{GSO}(2k)$ consists of the matrices $g$ satisfying $g^tJ_{2k}g=s(g)J_{2k}$ and $\det g=s(g)^k$, with $J_{2k}$ as in Section~\ref{pol des} and $s(g)\in\BG_m$. Note that this defines a homomorphism $s\colon\operatorname{GSO}(2k)\to\BG_m$, and $T_0\simeq\BG_m\times T'$, with the first factor given by $s$. For $g\in\operatorname{GSO}(2k)(\bF)$, let $\epsilon=\ord_t s(g)$. Since the multiplication by the central elements $t^mI_{2k}$, $m\in\BZ$, and $t^{\eta'}$, $\eta'\in X_*(T')$, changes $\xi$ by $(2m,0)$ and $(0,\eta')$, respectively, it suffices to consider $\xi=(\epsilon,0)$ with $\epsilon\in\{0,1\}$. For $\epsilon=0$, the fiber over $t^0$ of the map $\Gr_{\operatorname{GSO}(2k)}\to\Gr_{\BG_m}$ induced by $s$ is $\Gr_{\SO(2k)}$, and hence the calculations are the same as in Sections~\ref{sec:orbit-preliminaries} and~\ref{sec:irreducibles}.
\end{rem}

\subsubsection{The orthogonal case with \texorpdfstring{$\epsilon=1$}{epsilon=1}}
We assume that the characteristic is not two. We denote by $a_0\in\{1,2\}$ the index with $G_{a_0}=\SO(2k)$, and by $a_1$ the other index, so that $a_0=1$ if and only if $k\leq l$. In the second case of Remark~\ref{rem:app-isogeny-orthogonal}, we have $\widetilde G_{a_0}\simeq\operatorname{GSO}(2k)\times T'$ and $\widetilde G_{a_1}\simeq G_{a_1}\times T_0$. We omit the $T'$-coordinate from the matrices below. Note that the diagonal torus of $\operatorname{GSO}(2k)$ consists of the matrices $\diag(z_1,\ldots,z_k,cz_k^{-1},\ldots,cz_1^{-1})$. For $\lambda\in\widetilde\Lambda$, let $g_\lambda\in\operatorname{GSO}(2k)(\bF)$ be the image of $t^{\lambda^{(a_0)}}$. Let $t^{r_1},\ldots,t^{r_k}$ be its first $k$ diagonal entries, with $r_i\in\BZ$, and let $\epsilon=\ord_ts(g_\lambda)$. We denote
\begin{equation}\label{eq:app-isogeny-orthogonal-coordinates}
\bar\lambda^{(a_0)}=(r_1-\epsilon/2,\ldots,r_k-\epsilon/2)\in\tfrac12\BZ^k.
\end{equation}

Assume $\epsilon=1$, and let
\[
\eta=\bigl(\tfrac12,\ldots,\tfrac12\bigr),\qquad g_\eta=\begin{pmatrix}tI_k&0\\0&I_k\end{pmatrix},\qquad K_\eta=g_\eta\SO(2k)(\bO)g_\eta^{-1}.
\]
Then, we have $g_\eta\in\operatorname{GSO}(2k)(\bF)$ and $s(g_\eta)=t$. We denote $\mu_\lambda^{(a_0)}=(r_1-1,\ldots,r_k-1)\in X_*(T_{a_0})$, so that $g_\lambda=t^{\mu_\lambda^{(a_0)}}g_\eta$ and $\bar\lambda^{(a_0)}=\mu_\lambda^{(a_0)}+\eta$, and we denote by $\mu_\lambda^{(a_1)}\in X_*(T_{a_1})$ the $G_{a_1}$-component of $\lambda^{(a_1)}$. Note that the map $hK_\eta\mapsto hg_\eta\operatorname{GSO}(2k)(\bO)$ identifies $\SO(2k)(\bF)/K_\eta$ with the fiber of $\Gr_{\operatorname{GSO}(2k)}\to\Gr_{\BG_m}$ over $t$, and $g_\lambda\operatorname{GSO}(2k)(\bO)$ corresponds to $t^{\mu_\lambda^{(a_0)}}K_\eta$.

Let $\check\varepsilon_1,\ldots,\check\varepsilon_k$ be the basis of $X^*(T_{a_0})$ dual to $e_1^{(a_0)},\ldots,e_k^{(a_0)}$. By \eqref{eq:osp-duality-coordinates}, the super weight $\varepsilon_j$ corresponds to $-e_j^{(1)}$ if $a_0=1$, and to $e_j^{(2)}$ if $a_0=2$. For a geometric root $\check\beta$ of $\SO(2k)$, we have
\begin{equation}\label{eq:enlarged-K-root-filtration}
K_\eta\cap U_{\check\beta}(\bF) =U_{\check\beta} \bigl(t^{\langle\eta,\check\beta\rangle}\bO\bigr), \qquad K_\eta\cap T_{a_0}(\bF)=T_{a_0}(\bO).
\end{equation}
Note that the pairings $\langle\eta,\check\beta\rangle$ are integers. Namely, they are zero for $\check\varepsilon_i-\check\varepsilon_j$ and $\pm1$ for $\pm(\check\varepsilon_i+\check\varepsilon_j)$. Also, in $k\times k$ blocks, the elements of $K_\eta$ have the form
\begin{equation}\label{eq:enlarged-K-block-matrix}
\begin{pmatrix}A&tB\\t^{-1}C&D\end{pmatrix},
\qquad
\begin{pmatrix}A&B\\C&D\end{pmatrix}\in\SO(J_{2k},\bO).
\end{equation}
The stabilizer of $t^{\mu_\lambda^{(a_0)}}K_\eta$ in $\SO(2k)(\bF)$ is $t^{\mu_\lambda^{(a_0)}}K_\eta t^{-\mu_\lambda^{(a_0)}}$, and its intersection with $U_{\check\beta}(\bF)$ is $U_{\check\beta}(t^{\langle\bar\lambda^{(a_0)},\check\beta\rangle}\bO)$. Thus, in all the calculations in root subgroups, $\lambda^{(a_0)}$ is replaced by $\bar\lambda^{(a_0)}$.

\begin{lem}\label{lem:enlarged-orthogonal-normalforms}
Assume $\epsilon=1$. If $a_0=1$, the $K_\eta H_u(\bF)$-orbits in $\Gr_{G_2}$ have representatives
\begin{equation}\label{eq:enlarged-orthogonal-normalform}
t^{-\mu_\lambda^{(1)}}\sD t^{\mu_\lambda^{(2)}}G_2(\bO)/G_2(\bO),
\end{equation}
where $\lambda\in\widetilde\Lambda$ runs over the coweights over $\xi$ such that $(\bar\lambda^{(1)},\mu_\lambda^{(2)})$ satisfies the inequalities of Lemma~\ref{lem:L-orbit-indices}. In particular, the first inequalities are
\[
\bar\lambda^{(1)}_1\geq\cdots\geq \bar\lambda^{(1)}_{k-1}\geq|\bar\lambda^{(1)}_k|.
\]
If $a_0=2$, the $G_1(\bO)H_u(\bF)$-orbits in $G_2(\bF)/K_\eta$ have representatives $t^{-\mu_\lambda^{(1)}}\sD t^{\mu_\lambda^{(2)}}K_\eta/K_\eta$, where the inequalities of Lemma~\ref{lem:L-orbit-indices} are imposed on $(\mu_\lambda^{(1)},\bar\lambda^{(2)})$. In both cases, the stabilizers are connected, and if the orbit indexed by $\lambda'$ is contained in the closure of the orbit indexed by $\lambda$, then $\lambda'\leq_{\super}\lambda$.
\end{lem}

\begin{proof}
Note that, for $a_0=2$, the element $g_\eta$ normalizes $P^-$, and hence $G_2(\bF)=P^-(\bF)K_\eta$. In both cases, the parabolic projection of Lemma~\ref{lem:L-orbit-indices} reduces the assertions to the case $\SO(M)\subset\SO(M+1)$, with $M$ as in Section~\ref{pol des}, and then they follow by the same proofs as in \cite[Lemmas 2.3.2--2.3.3 and the ``only if'' part of Theorem 3.3.5(a,b)]{[BFT]}, using \eqref{eq:enlarged-K-root-filtration}. Note that, in the even orthogonal factor, the lattices $\mathcal L=hg_\eta\bO^{2k}$, $h\in\SO(2k)(\bF)$, satisfy $\mathcal L^\perp=t^{-1}\mathcal L$ for $J_{2k}$, instead of $\mathcal L^\perp=\mathcal L$. For instance, for $h=t^{\mu_\lambda^{(a_0)}}$, the exponents on the paired standard basis vectors $e_i,e_{2k+1-i}$ are $r_i,1-r_i$, and by subtracting $\tfrac12$, we obtain $\pm\bar\lambda^{(a_0)}_i$.
\end{proof}

We denote by $\sO_{\Gr}^{\lambda}$ the corresponding $H(\bF)$-orbit in the fiber of $\Gr_{\widetilde G}$ over $t^\xi$.

\subsubsection{Relevance}

The relevance condition is the same as in Proposition~\ref{prop:relevant-L-orbits}, with $b$ replaced by $\widetilde b$. Indeed, by \eqref{eq:enlarged-K-root-filtration}, we have the same residue conditions and the same equalities among the coordinates of $\bar\lambda^{(a_0)}$ and $\mu_\lambda^{(a_1)}$, and hence the stabilizer has the same reductive factors. Also, on its one-dimensional tori, the twisting is $\widetilde b(\lambda,\alpha_{a,i}^{\pm})$, where $\alpha_{a,i}^{\pm}$ are the simple super roots of Definition~\ref{def:mixed-pure-coroots}. Since the stabilizers are connected, each relevant orbit carries a unique irreducible equivariant local system, which defines $\IC_{\Gr}^{\lambda}$.

\subsection{Intersection dimensions and the two projections}
\label{sec:app-isogeny-projections}

By the above discussion, it remains to consider the second case of Remark~\ref{rem:app-isogeny-orthogonal} with $\xi=(1,0)$. In this subsection, we explain how the dimension estimates of Section~\ref{sec:orbit-preliminaries} and the two projections of Section~\ref{sec:irreducibles} are modified in this case. We choose $\widetilde\eta\in\widetilde\Lambda$ over $\xi$ such that the $\operatorname{GSO}(2k)$-component of $t^{\widetilde\eta}$ is $g_\eta$ and its $G_{a_1}$- and $T^{\prime}$-components in the product of Remark~\ref{rem:app-isogeny-orthogonal} are trivial. For $a=1,2$ and $\nu\in\widetilde\Lambda$ over $\xi$, let
\[
K_a=G_a(\bF)\cap t^{\widetilde\eta^{(a)}}\widetilde G_a(\bO)t^{-\widetilde\eta^{(a)}},\qquad \mu_\nu^{(a)}=\nu^{(a)}-\widetilde\eta^{(a)}\in\Lambda_{G_a}.
\]
Then, we have $K_{a_0}=K_\eta$, $K_{a_1}=G_{a_1}(\bO)$, and $\mu_\lambda^{(a)}$ is as above. Note that the multiplication on the right by $t^{\widetilde\eta^{(a)}}$ identifies $G_a(\bF)/K_a$ with the fiber of $\Gr_{\widetilde G_a}$ over $t^\xi$, and $t^{\mu_\nu^{(a)}}K_a$ corresponds to $t^{\nu^{(a)}}$.

We use the statements and proofs of Sections~\ref{sec:orbit-preliminaries} and~\ref{sec:irreducibles} with the following replacements:
\begin{itemize}
\item $G_a(\bO)$ is replaced by $K_a$;
\item in the formulas for points and subsets, $\lambda^{(a)},\theta^{(a)},\lambda^{(a)\prime}$ are replaced by $\mu_\lambda^{(a)},\mu_\theta^{(a)},\mu_{\lambda'}^{(a)}$;
\item the indices are $\lambda,\theta,\lambda'\in\widetilde\Lambda$ over $\xi$, the orbit indices satisfy the inequalities of Lemma~\ref{lem:enlarged-orthogonal-normalforms}, and the super order, $d$ and $\widetilde b$ are evaluated on these coweights.
\end{itemize}
Note that, since $\lambda-\theta\in\Lambda$, the super order and $d(\lambda-\theta)$ are the same as for $G$. The proofs remain valid for the following reasons. First, $K_{a_0}=g_\eta\SO(2k)(\bO)g_\eta^{-1}$ with $g_\eta$ in the maximal torus of $\operatorname{GSO}(2k)$, so the Iwasawa and Gauss decompositions hold for $K_a$. Secondly, we have $K_a\supset\ker(G_a(\bO)\to G_a(\bO/t^2\bO))$. Finally, the fiber of $\Gr_{\widetilde G}$ over $t^\xi$ is a union of connected components of $\Gr_{\widetilde G}$ preserved by loop rotation, so the hyperplane argument of \cite[Section 3]{[MV]} and loop rotation can be used in $\Gr_{\widetilde G}$.

\begin{lem}\label{lem:app-isogeny-full-dimension}
Let $\lambda,\theta\in\widetilde\Lambda$ have the same image in $\Lambda_0$, with $\lambda$ dominant. If $\overline\sO_{\Gr}^{\lambda}\cap \overline\sS_{\Gr}^{\theta}\neq\varnothing$, then $\lambda-\theta\in \Lambda^{\mathrm{pos}}$. Moreover,
\begin{equation}\label{eq:app-isogeny-full-dimension}
\begin{split}
&\dim\bigl(\overline\sO_{\Gr}^{\lambda} \cap\overline\sS_{\Gr}^{\theta}\bigr)\\
&\quad\leq \max_{\substack{m_{a,i}\in\BZ_{\geq0}\\
\lambda-\theta-\sum_{a,i}m_{a,i}\alpha_{a,i} \in \Lambda^{\mathrm{pos}}}} \sum_{a,i}m_{a,i} \leq\tfrac12d(\lambda-\theta).
\end{split}
\end{equation}
Equality in the last bound is possible only if $\lambda-\theta$ is a sum of positive geometric simple coroots. For $\mu\in\widetilde\Lambda$ over $\xi$, the difference between the degree of the global stratum $\cM_{=\mu\cdot c}$ and that of the semi-infinite kernel of index $\theta$ is $d(\mu-\theta)$, with the normalization of Remark~\ref{rem:global-stratum-degree}.
\end{lem}
\begin{proof}
With the representatives and closure relations of Lemma~\ref{lem:enlarged-orthogonal-normalforms}, the proofs of Lemma~\ref{lem:orbit-S-nonempty-super-order} and Corollary~\ref{lem:dimension-estimate-super} apply, since the geometric coroots and the values of $d$ on them are the same as for $G$. Also, the relative shift follows from the same normalization as in Remark~\ref{rem:global-stratum-degree}.
\end{proof}

\subsubsection{The two projections}
\label{sec:app-isogeny-two-projections}

Let $\lambda,\theta\in\widetilde\Lambda$ over $\xi$. Note that, for the embedding before the conjugation by $\sD$, the orbit representative is $(t^{\mu_\lambda^{(1)}}K_1/K_1,\sD t^{\mu_\lambda^{(2)}}K_2/K_2)$. As in \eqref{eq:Y-convolution} and \eqref{eq:Z-convolution}, we consider
\begin{equation}\label{eq:app-isogeny-based-products}
\begin{split}
Y_{\widetilde\eta}={}&(\sD^{-1}t^{\mu_\theta^{(1)}}U_1^-(\bF)K_1) \mathbin{\times^{K_1}}\\
&\quad(K_1H_u(\bF)t^{-\mu_\lambda^{(1)}} \sD t^{\mu_\lambda^{(2)}}K_2/K_2),\\
Z_{\widetilde\eta}={}&(\sD t^{\mu_\theta^{(2)}}U_2(\bF)K_2) \mathbin{\times^{K_2}}\\
&\quad(K_2H_u(\bF)t^{-\mu_\lambda^{(2)}} \sD^{-1}t^{\mu_\lambda^{(1)}}K_1/K_1).
\end{split}
\end{equation}
Then, their multiplication maps take values in $G_2(\bF)/K_2$ and $G_2(\bF)/K_1$, respectively, and their fiber products with $U_2(\bF)t^{\mu_\theta^{(2)}}K_2/K_2$ over $G_2(\bF)/K_2$ and with $U_1^-(\bF)t^{\mu_\theta^{(1)}}K_1/K_1$ over $G_2(\bF)/K_1$, respectively, are both identified with $\sO_{\Gr}^{\lambda}\cap\sS_{\Gr}^{\theta}$.

Let $\lambda'\in\widetilde\Lambda$ over $\xi$. As in \eqref{eq:g2-projection-fiber} and \eqref{eq:g1-projection-fiber}, the fibers of the projections of $\sO_{\Gr}^{\lambda}\cap\sS_{\Gr}^{\lambda'}$ to $\Gr_{\widetilde G_2}$ and $\Gr_{\widetilde G_1}$ over $t^{\lambda^{(2)\prime}}$ and $t^{\lambda^{(1)\prime}}$, respectively, are identified with
\begin{equation}\label{eq:app-isogeny-canonical-fibers}
\begin{split}
&\sD t^{\mu_{\lambda'}^{(2)}}K_2t^{-\mu_\lambda^{(2)}}\sD^{-1} t^{\mu_\lambda^{(1)}}H_u(\bF)K_1/K_1 \ \cap U_1^-(\bF)t^{\mu_{\lambda'}^{(1)}}K_1/K_1,\\
&\sD^{-1}t^{\mu_{\lambda'}^{(1)}}K_1H_u(\bF)t^{-\mu_\lambda^{(1)}} \sD t^{\mu_\lambda^{(2)}}K_2/K_2 \ \cap U_2(\bF)t^{\mu_{\lambda'}^{(2)}}K_2/K_2.
\end{split}
\end{equation}

\subsubsection{The strata of the two projections}

The analogues of Lemmas~\ref{lem:semi-spherical-orbits}, \ref{lem:g2-strata-in-image} and \ref{lem:g1-semi-spherical-indices} hold with the same proofs, in which the groups $t^\nu G_1(\bO)t^{-\nu}$ and $t^\nu L(\bO)t^{-\nu}$ are replaced by $t^\nu K_1t^{-\nu}$ and $t^\nu(L(\bF)\cap K_2)t^{-\nu}$. Here, $\nu$ denotes the sufficiently dominant coweight which is denoted by $\eta$ in those proofs, since $\eta$ already has another meaning in this appendix.

\begin{lem}\label{lem:app-isogeny-critical-strata}
Let $\lambda,\theta\in\widetilde\Lambda$ over $\xi$, with $\lambda$ dominant.
\begin{enumerate}
\item A stratum of $\sO_{\Gr}^{\lambda}\cap\sS_{\Gr}^{\theta}$ given by either multiplication map in \eqref{eq:app-isogeny-based-products} can have the critical dimension $\tfrac12d(\lambda-\theta)$ only at the index $(\theta^{(1)},\lambda^{(2)})$ for $Y_{\widetilde\eta}$, resp. $(\lambda^{(1)},\theta^{(2)})$ for $Z_{\widetilde\eta}$.
\item Let $C$ be a critical-dimensional irreducible component of the fiber over $t^{\lambda^{(1)}}$ of the projection $\sO_{\Gr}^{\lambda}\cap \sS_{\Gr}^{(\lambda^{(1)},\theta^{(2)})} \to\Gr_{\widetilde G_1}$. Then, the projection to $\Gr_{\widetilde G_2}$ identifies a dense open subset of $C$ with a dense open subset of a critical-dimensional irreducible component of the analogue of \eqref{eq:critical-G2-fiber}, where the critical dimension is $\tfrac12d\bigl((0,\lambda^{(2)}-\theta^{(2)})\bigr)$. Moreover, the analogue of Lemma~\ref{lem:critical-G2-fiber} holds.
\item If $\lambda^{(2)}=\theta^{(2)}$, the critical-dimensional irreducible components of the fiber of $\sO_{\Gr}^{\lambda}\cap\sS_{\Gr}^{\theta}\to\Gr_{\widetilde G_2}$ over $t^{\lambda^{(2)}}$ correspond to those of
\begin{equation}\label{eq:app-isogeny-first-integral}
\begin{split}
&\sD t^{\mu_\lambda^{(2)}}\bigl(K_2\cap U_2(\bF)T_2(\bO)\bigr) t^{-\mu_\lambda^{(2)}}\sD^{-1}t^{\mu_\lambda^{(1)}}H_u(\bF)K_1/K_1\\
&\hspace{25mm}\cap U_1^-(\bF)t^{\mu_\theta^{(1)}}K_1/K_1,
\end{split}
\end{equation}
and the inclusion of \eqref{eq:app-isogeny-first-integral} into the analogue of \eqref{eq:g1-limiting-intersection} induces an injection on irreducible components of critical dimension. These identifications preserve the restrictions of the two trivializations defining $\cL_{\lambda,\theta}$.
\item For a critical-dimensional irreducible component $Z$ of the analogue of \eqref{eq:critical-G2-fiber} as in (2), the intersection of $\overline Z$ with the two boundaries specified in Corollary~\ref{cor:critical-G2-fiber-boundary} has dimension at most $\dim Z-2$. The same assertion holds for a critical-dimensional irreducible component $Z$ as in (3), with the boundaries specified in Corollary~\ref{cor:g1-boundary-codimension}.
\end{enumerate}
\end{lem}
\begin{proof}
By Lemma~\ref{lem:app-isogeny-full-dimension}, the proofs of Propositions~\ref{prop:g2-base-critical}, \ref{prop:g2-fiber-dimension} and \ref{prop:g1-strata-dimension} and of Corollary~\ref{rem:unique-g2-critical-stratum} apply, and we obtain (1). Then, (2) follows as in Lemma~\ref{lem:critical-G2-fiber}, (3) as in Corollary~\ref{cor:g1-critical-component-comparison}, and (4) as in Corollaries~\ref{cor:critical-G2-fiber-boundary} and \ref{cor:g1-boundary-codimension}. Note that, in the proof of Lemma~\ref{lem:g1-open-stratum}, the group $U_2(\bO)T_2(\bO)=G_2(\bO)\cap U_2(\bF)T_2(\bO)$ is replaced by $K_2\cap U_2(\bF)T_2(\bO)$.
\end{proof}

\subsection{The local system in the rank-one case}
\label{sec:app-isogeny-cuts}

Let $Z$ be an irreducible component of $\sO_{\Gr}^{\lambda}\cap\sS_{\Gr}^{\theta}$ as in the last assertion of Theorem~\ref{thm:isogeny-image-of-ic}, i.e., $\dim Z=\tfrac12d(\lambda-\theta)$, and $\cL_{\lambda,\theta}\otimes\psi$ is constant on a smooth dense open subset of $Z$. As in Sections~\ref{sec:g2-direction} and \ref{sec:g1-direction}, we show that $\theta=\lambda$. Note that the unipotent radicals and the rank-one products are the same as there, and their intersections with $K_a$ are given by \eqref{eq:enlarged-K-root-filtration}. By Lemma~\ref{lem:app-isogeny-critical-strata}, the cuts have critical dimension and meet the original Gaiotto orbit in dense open subsets, and hence the proofs of Lemmas~\ref{lem:rank-one-fixed-orbit-comparison} and \ref{lem:g1-hyperplane-block} apply without change, and give the analogues of Corollaries~\ref{prop:g2-sharp-difference} and \ref{prop:g1-sharp-difference}. Namely, first for $a=2$, and then, by using $\lambda^{(2)}=\theta^{(2)}$, for $a=1$, a non-zero difference is reduced by $q_{a,i}$ successive cuts at a time to
\begin{equation}\label{eq:app-isogeny-rank-one-difference}
\lambda^{(a)}-\theta^{(a)}=q_{a,i}\alpha_{a,i},
\end{equation}
with $q_{a,i}$ as in Section~\ref{sec:gaiotto-isogenies}, and $\cL_{\lambda,\theta}\otimes\psi$ and the ratio local system of the two semi-infinite trivializations are constant on a dense open subset of the resulting component. Here, $\theta$ denotes the index of the semi-infinite orbit after the cuts. Finally, it remains to exclude \eqref{eq:app-isogeny-rank-one-difference}. For a pure coroot in the sense of Definition~\ref{def:mixed-pure-coroots}, this follows from \cite[Theorem 9.4.9 and Remark 9.5.4]{[GL]} and Lemma~\ref{lem:kummer-as-separation}. For a mixed coroot, it follows from the rank-one calculations in the proofs of Propositions~\ref{prop:part-q} and \ref{prop:part-p}, with $b$ replaced by $\widetilde b$ and with the valuations in the $\SO(2k)$-factor given by \eqref{eq:enlarged-K-root-filtration}. Indeed, for the roots $\pm(\check\varepsilon_i-\check\varepsilon_j)$, the differences $\mu_{\lambda,i}^{(a_0)}-\mu_{\lambda,j}^{(a_0)}=\bar\lambda_i^{(a_0)}-\bar\lambda_j^{(a_0)}$ are unchanged, and for $\pm(\check\varepsilon_i+\check\varepsilon_j)$, the factors $t^{\pm1}$ in the off-diagonal blocks of \eqref{eq:enlarged-K-block-matrix} replace $\mu_{\lambda,i}^{(a_0)}+\mu_{\lambda,j}^{(a_0)}$ by $\mu_{\lambda,i}^{(a_0)}+\mu_{\lambda,j}^{(a_0)}+1=\bar\lambda_i^{(a_0)}+\bar\lambda_j^{(a_0)}$ in the corresponding valuation bounds. Thus, the end-root calculations give the inequalities of Section~\ref{sec:irreducibles}, written in terms of $\bar\lambda^{(a_0)}$. For instance, if $a_0=2$, the inclusion $t^{-\lambda^{(2)}}B_2(\bO)t^{\lambda^{(2)}}\subset G_2(\bO)$ used in Section~\ref{10.4.5} becomes $t^{-\mu_\lambda^{(2)}}B_2(\bO)t^{\mu_\lambda^{(2)}}\subset K_\eta$, which holds since $\langle\bar\lambda^{(2)},\check\beta\rangle\leq0$ for every positive root $\check\beta$.

\par \ifdim\dimexpr\pagegoal-\pagetotal\relax<7\baselineskip\newpage\fi
\subsection{Proofs of the statements}

\begin{proof}[Proof of Theorem~\ref{thm:isogeny-image-of-ic}]
By Lemma~\ref{lem:app-isogeny-components} and Remark~\ref{rem:app-isogeny-orthogonal}, we only need to consider the orthogonal case with $\epsilon=1$. In this case, the assertion $\theta=\lambda$ follows from Section~\ref{sec:app-isogeny-cuts}, and it replaces Theorem~\ref{thm:local-system-rigidity} in the proof of Theorem~\ref{thm:image-of-ic}, which then applies with the dimension estimates and relative shifts of Lemma~\ref{lem:app-isogeny-full-dimension}.
\end{proof}

\begin{proof}[Proof of Theorem~\ref{thm:isogeny-metaplectic-hecke-action}]
It follows from the same proof as that of Theorem~\ref{thm:metaplectic-hecke-action}, using Lemma~\ref{lem:app-isogeny-full-dimension} and the last assertion of Theorem~\ref{thm:isogeny-image-of-ic}.
\end{proof}

\section{Relevant Iwahori orbits in the orthosymplectic and exceptional cases}
\label{sec:orthogonal-exceptional-Iwahori}

In this appendix, we classify the relevant orbits, and hence the irreducible objects, in the Iwahori Gaiotto category in the orthosymplectic and $\on{F}(4)$ cases. Note that this classification is not needed in the proofs of the main theorems, and we include it for completeness. Here, we work over $\BC$, with $\bF=\BC((t))$ and $\bO=\BC[[t]]$.

\subsection{Classification of orbits}
\label{sec:Iwahori-Levi-reduction}

In this subsection, we first reduce the classification of orbits and stabilizers to the Levi subgroup. Let $G_1,G_2,L$ and $\iota\colon G_1\hookrightarrow L$ be as in Section~\ref{pol des}. We use upper triangular conventions analogous to those of Section~\ref{sec:type-a-conventions}, and take $P=L\ltimes H_u$ to be the upper triangular parabolic subgroup of $G_2$, and $H=G_1\ltimes H_u$, with the diagonal embedding of $G_1$. For $a=1,2$, let $I_a\subset G_a(\bO)$ be the upper triangular Iwahori subgroup and $\Fl_{G_a}=G_a(\bF)/I_a$. Let $B_L=L\cap B_{G_2}$, and note that $\iota^{-1}(B_L)=B_{G_1}$. We denote $I_L=L(\bF)\cap I_2$ and $I_H=I_1\ltimes H_u(\bF)$. As in Section~\ref{aff kl}, it is enough to classify the $I_H$-orbits in $\Fl_{G_2}$. We denote by $\Phi_2^+,\Phi_L^+$ the positive geometric roots for $B_{G_2},B_L$, by $W_2,W_L$ the finite Weyl groups, and let
\[
{}^LW_2=\{w_2\in W_2:w_2^{-1}\Phi_L^+\subset\Phi_2^+\}.
\]
We choose representatives of $W_2$ in $N_{G_2}(T_2)(\BC)$. By the Iwasawa decomposition $G_2(\bF)=P(\bF)G_2(\bO)$ and the Bruhat decomposition of $G_2(\BC)$ with respect to $P(\BC)$ and $B_{G_2}(\BC)$, we have
\begin{equation}\label{eq:Iwahori-parabolic-cells}
G_2(\bF)=\coprod_{w_2\in{}^LW_2}P(\bF)w_2I_2, \qquad L(\bF)\cap w_2I_2w_2^{-1}=I_L.
\end{equation}
The second equality follows from $w_2^{-1}\Phi_L^+\subset\Phi_2^+$, since, in both groups, the positive and negative root subgroups of $L$ have parameters in $\bO$ and $t\bO$, respectively. We denote $U_{w_2}=H_u(\bF)\cap w_2I_2w_2^{-1}$. Then, we have $P(\bF)\cap w_2I_2w_2^{-1}=I_L\ltimes U_{w_2}$, so the $H_u(\bF)$-orbits in $P(\bF)w_2I_2/I_2$ are identified with $L(\bF)/I_L$. Thus, the $I_H$-orbits in $\Fl_{G_2}$ are indexed by
\begin{equation}\label{eq:Iwahori-Levi-indices}
w_2\in{}^LW_2,\qquad I_1\text{-orbits in }L(\bF)/I_L.
\end{equation}
Here, the representative is $lw_2I_2/I_2$, with $l\in L(\bF)$, and the corresponding $H(\bF)$-orbit in $\Fl_{G_1}\times\Fl_{G_2}$ has the representative $(t^0,lw_2I_2/I_2)$, where $t^0$ is the point represented by the identity matrix. In the following subsections, we make this classification explicit.

\subsubsection{Stabilizers} We denote by $R_l=I_1\cap lI_Ll^{-1}$ the stabilizer of $lI_L$ in $I_1$. Then, we have
\begin{equation}\label{eq:Iwahori-stabilizer-splitting}
\operatorname{Stab}_{I_H}(lw_2I_2/I_2) =R_l\ltimes\operatorname{Ad}_l(U_{w_2}).
\end{equation}
For a geometric root $\check\alpha$ of $H_u$, let
\begin{equation}\label{eq:Iwahori-root-level}
\epsilon_{w_2}(\check\alpha)=
\begin{cases}
0,&w_2^{-1}\check\alpha\in\Phi_2^+,\\
1,&w_2^{-1}\check\alpha\in-\Phi_2^+.
\end{cases}
\end{equation}
Note that the corresponding root subgroup of $U_{w_2}$ has parameter in $t^{\epsilon_{w_2}(\check\alpha)}\bO$, and $U_{w_2}$ is the product of these root subgroups. Let $\psi^0\colon H\to\BG_a$ be the additive character given in \eqref{eq:orthogonal-Iwahori-character} and Section~\ref{sec:F4-Iwahori-relevance}. As in Section~\ref{aff kl}, we use
\[
\psi(u)=[t^0]\psi^0(u),
\]
where $[t^0]$ denotes the constant coefficient. We choose a non-zero root vector $X_{\check\alpha}\in\operatorname{Lie}(H_u)_{\check\alpha}$ over $\BC$. Then, the stabilizer \eqref{eq:Iwahori-stabilizer-splitting} lies in the kernel of $\psi$ if and only if, for every geometric root $\check\alpha$ of $H_u$,
\begin{equation}\label{eq:Iwahori-general-additive-test}
\operatorname{ord}_t \bigl(d\psi^0(\operatorname{Ad}_lX_{\check\alpha})\bigr) +\epsilon_{w_2}(\check\alpha)\geq1.
\end{equation}
Here, $\operatorname{ord}_t(0)=+\infty$. Indeed, $\psi^0$ vanishes on $R_l$ and factors through the abelianization of $H_u$, and $[t^0](ax)=0$ for every $x\in t^\epsilon\bO$ if and only if $\operatorname{ord}_ta+\epsilon\geq1$.

\subsection{The orthosymplectic orbits}

In this subsection, we classify the $I_1$-orbits in $L(\bF)/I_L$ in the orthosymplectic case.

\subsubsection{The orthogonal affine flag variety} We first describe the orthogonal affine flag variety inside the general linear one. Let $M,N,m$ be as in Section~\ref{pol des}, so that $N=M+1+2m$, and let $n=M+1$. We use the matrices $J_M,J_n$ of Section~\ref{pol des}, and put $i^*=M+1-i$ for $1\leq i\leq M$. We consider $V=\bF^M\oplus\bF e_n$, where $e_n$ is the last standard basis vector. The symmetric bilinear form is given by $J=\diag(J_M,1)$. We denote by $\sC\in\GL_n(\BC)$ the matrix identifying the ordered basis of the orthogonal complement of the anisotropic vector in Section~\ref{pol des} with the standard basis of $\bF^M$, and the anisotropic vector with $e_n$. Then $\sC^tJ\sC=J_n$. Let $I_M,I_n$ be the upper triangular Iwahori subgroups of $\GL_M(\bF),\GL_n(\bF)$, and let
\[
I_M^{\rm orth}=\SO(J_M,\bF)\cap I_M,\qquad I_n^{\rm orth}=\SO(J,\bF)\cap \sC I_n\sC^{-1}.
\]
The map $gI_n^{\rm orth}\mapsto g\sC I_n$ embeds $\SO(J,\bF)/I_n^{\rm orth}$ into $\GL_n(\bF)/I_n$, and its image consists of the cosets $xI_n$ with a representative satisfying $x^tJx=J_n$ and $\det x=\det\sC$.

We denote by $\mathcal Q_n$ the set of symmetric $n\times n$ matrices $Q$ over $\bO$ satisfying
\begin{equation}\label{eq:orthogonal-compatible-Gram}
Q_{ij}\in t\bO\quad(i+j\leq n),\qquad Q_{i,n+1-i}\in\bO^\times.
\end{equation}
Equivalently, $Q=b^tJ_nb$ for some $b\in I_n$. Indeed, modulo $t$, the conditions \eqref{eq:orthogonal-compatible-Gram} imply $\bar Q=\bar b^{\,t}J_n\bar b$ for an invertible upper triangular matrix $\bar b$, and then the coefficients of $t,t^2,\ldots$ in $b$ can be chosen successively. We define $\mathcal Q_M$ similarly.

\subsubsection{Classification of orbits} In this subsubsection, we determine which $I_M$-orbits in $\GL_n(\bF)/I_n$ meet the orthogonal affine flag variety. By Proposition~\ref{prop 3.2.7}, with $N$ replaced by $n$, these orbits have representatives $r=g_{w(\beta)}w$. We will show in the proof of Proposition~\ref{prop:adjacent-orthogonal-Iwahori-indices} below that the orbit of $r$ meets the orthogonal affine flag variety if and only if there exists a matrix $Q$ satisfying
\begin{equation}\label{eq:orthogonal-Gram-filter}
Q\in\mathcal Q_M,\qquad r^t\operatorname{diag}(Q,1)r\in\mathcal Q_n.
\end{equation}
Note that, if $r\in t^{-A}\operatorname{Mat}_n(\bO)$ with $A\geq0$, then \eqref{eq:orthogonal-Gram-filter} depends only on $Q$ modulo $t^{2A+1}$. Thus, it suffices to take $Q=\sum_{a=0}^{2A}Q_at^a$ with $Q_a^t=Q_a$, in which case \eqref{eq:orthogonal-Gram-filter} consists of finitely many linear equations and non-vanishing conditions on the coefficients of $Q$.

\begin{prop}\label{prop:adjacent-orthogonal-Iwahori-indices}
The $I_M^{\mathrm{orth}}$-orbits on $\SO(J,\bF)/I_n^{\mathrm{orth}}$ are indexed by the type A pairs $(w,\beta)$ for which there exists $Q$ satisfying \eqref{eq:orthogonal-Gram-filter}. Every stabilizer is connected.
\end{prop}

\begin{proof}
The proof consists of three steps. First, we show that the general linear orbit indexed by $(w,\beta)$ meets the orthogonal affine flag variety if and only if there is $Q$ satisfying \eqref{eq:orthogonal-Gram-filter}. If this orbit contains an orthogonal flag, the flag has a representative $x=\diag(h,1)rb^{-1}$ with $h\in I_M$, $b\in I_n$, and $x^tJx=J_n$. This condition is equivalent to
\[
r^t\diag(h^tJ_Mh,1)r=b^tJ_nb,
\]
so $Q=h^tJ_Mh$ satisfies \eqref{eq:orthogonal-Gram-filter}. Conversely, if $Q$ satisfies \eqref{eq:orthogonal-Gram-filter}, then there exist $h\in I_M$ and $b\in I_n$ with
\[
h^tJ_Mh=Q,\qquad b^tJ_nb=r^t\diag(Q,1)r.
\]
Then, $x=\diag(h,1)rb^{-1}$ satisfies $x^tJx=J_n$. After multiplying $x$ on the right by $-\Id_n$ if $n$ is odd, or on the left by $\diag(-\Id_M,1)$ if $M$ is odd, we may assume that $\det x=\det\sC$. The corresponding point of $\SO(J,\bF)/I_n^{\rm orth}$ is represented by $g=x\sC^{-1}$.

Secondly, we show that the intersection is a single $I_M^{\rm orth}$-orbit. We consider the involutions
\[
\tau(a)=J_M^{-1}a^{-t}J_M\quad(a\in\GL_M(\bF)),\qquad xI_n\longmapsto J^{-1}x^{-t}J_nI_n.
\]
The first one preserves $I_M$, the second one is well-defined on $\GL_n(\bF)/I_n$, and they are compatible with the action of $\GL_M(\bF)$ through $a\mapsto\diag(a,1)$ and fix the orthogonal flags. We denote by $S\subset I_M$ the stabilizer of an orthogonal flag $xI_n$. Then, $S$ is $\tau$-stable. By the stabilizer calculation following Proposition~\ref{prop 3.2.7}, and since conjugation by $h\in I_M$ does not change the image in the diagonal quotient, $S$ is an extension of $T_{w,\sigma}$ by a pro-unipotent group $S_u$, where $S\to T_{w,\sigma}$ is the restriction of the projection of $I_M$ to the diagonal torus $T_M(\BC)$. Note that this projection is $\tau$-equivariant, with $\tau(\diag(z_i))=\diag(z_{i^*}^{-1})$, so $T_{w,\sigma}$ and $S_u$ are $\tau$-stable.

We claim that $H^1(\langle\tau\rangle,S)=1$, where $\langle\tau\rangle=\{1,\tau\}$, i.e., that every $s\in S$ with $s\tau(s)=1$ can be written as $s=a^{-1}\tau(a)$ with $a\in S$. We first solve this equation in $T_{w,\sigma}$. Namely, on a pair of coordinates $(z_i,z_{i^*})$, $i\neq i^*$, we have $\tau(u,v)=(v^{-1},u^{-1})$, so $s\tau(s)=1$ gives $s=(z,z)$, and $a=(1,z^{-1})$ is a solution. On a middle coordinate, $\tau(z)=z^{-1}$, and the equation is $a^{-2}=s$, which has a solution over $\BC$. Let $a_0\in S$ be a lift of the solution in $T_{w,\sigma}$. Then, $s_0=a_0s\tau(a_0)^{-1}$ belongs to $S_u$ and satisfies $s_0\tau(s_0)=1$. For $S_u$, we use $\tau$-stable central filtrations of its unipotent quotients. Note that, on each additive subquotient, if $v+\tau(v)=0$, then $v=\tau(u)-u$ with $u=-v/2$. Thus, by correcting successively in the kernels of the quotients, we obtain $c\in S_u$ with $s_0=c^{-1}\tau(c)$, and then $s=(ca_0)^{-1}\tau(ca_0)$.

Now, let $xI_n,yI_n$ be two orthogonal flags in the same general linear orbit, and let $k\in I_M$ with $yI_n=\diag(k,1)xI_n$. By applying the involution, we obtain $s=k^{-1}\tau(k)\in S$ with $s\tau(s)=1$, and we choose $a\in S$ with $s=a^{-1}\tau(a)$. Then, $ka^{-1}$ is fixed by $\tau$ and still sends $xI_n$ to $yI_n$, and hence the two flags are in the same $I_M\cap\mathrm O(J_M,\bF)$-orbit. Note that, for even $M$, this group has determinant one, and for odd $M$, the group $I_n\cap\mathrm O(J_n,\bF)$ has determinant one, since $n$ is even, so the condition $\det x=\det y=\det\sC$ forces the element of $I_M\cap\mathrm O(J_M,\bF)$ to have determinant one as well.

Finally, we show that the stabilizer is connected. By the calculation on $S_u$ above, $S^\tau\to T_{w,\sigma}^\tau$ is surjective, with connected pro-unipotent kernel $S_u^\tau$. Note that the coordinates of $T_{w,\sigma}^\tau$ satisfy $z_{i^*}=z_i^{-1}$, and a middle coordinate is $\pm1$, which is forced to be $1$ by the determinant-one condition. Thus, the stabilizer in $I_M^{\rm orth}$ is an extension of a torus by a connected pro-unipotent group, and hence it is connected.
\end{proof}

\subsubsection{The orbit representatives} Note that the $\SO(n)$ factor of $L$ in Section~\ref{pol des} is defined by $J_n$, whereas the matrix $g$ above preserves $J=\diag(J_M,1)$. Therefore, in \eqref{eq:orthogonal-Iwahori-Levi-matrix} below, we use
\[
g_{\mathrm{adj}}=\sC^{-1}g\sC\in\SO(J_n,\bF).
\]
Since $L=\SO(n)\times\BG_m^m$, by \eqref{eq:Iwahori-Levi-indices} and Proposition~\ref{prop:adjacent-orthogonal-Iwahori-indices}, the $I_H$-orbits in $\Fl_{G_2}$ are indexed by
\begin{equation}\label{eq:orthogonal-Iwahori-all-indices}
w_2\in{}^LW_2,\qquad (c_1,\ldots,c_m)\in\BZ^m,\qquad (w,\beta)\text{ as in Proposition~\ref{prop:adjacent-orthogonal-Iwahori-indices}},
\end{equation}
with no dominance condition on the $c_j$, and the representative is $lw_2I_2/I_2$, where
\begin{equation}\label{eq:orthogonal-Iwahori-Levi-matrix}
l=\operatorname{diag} (t^{c_1},\ldots,t^{c_m},g_{\mathrm{adj}},t^{-c_m},\ldots,t^{-c_1}).
\end{equation}
Here, $w_2$ is a signed permutation satisfying $w_2^{-1}\Phi_L^+\subset\Phi_2^+$, with an even number of sign changes in the even orthogonal case.

\subsection{The orthosymplectic relevance conditions}
\label{sec:orthogonal-Iwahori-relevance}

In this subsection, we determine which orthosymplectic orbits are relevant.

\subsubsection{The additive character} Let
\[
e_1,\ldots,e_m, f_1,\ldots,f_n,e_m^*,\ldots,e_1^*
\]
be an ordered basis of the defining representation $\bF^N$ of $G_2=\SO(N)$, with $(e_i,e_j^*)=\delta_{ij}$ and $(f_s,f_{s'})=\delta_{s,n+1-s'}$. The $\SO(n)$ factor of $L$ acts on the span of $f_1,\ldots,f_n$; these vectors are denoted by $v_1,\ldots,v_n$ in Section~\ref{pol des}. As in that section, we denote by $v$ the anisotropic vector fixed by $G_1=\SO(M)$, i.e.,
\[
v=f_{(n+1)/2}\quad(n\text{ odd}),\qquad v=(f_{n/2}+f_{n/2+1})/\sqrt2\quad(n\text{ even}).
\]
Note that the matrix $\sC$ above satisfies $\sC v=e_n$. For $m>0$ and $u\in H_u$, let $X(u)=(u_{m,m+1},\ldots,u_{m,m+n})$. Then
\begin{equation}\label{eq:orthogonal-Iwahori-character}
\psi^0(u)=\sum_{j=1}^{m-1}u_{j,j+1}+X(u)v, \qquad \psi(u)=[t^0]\psi^0(u).
\end{equation}
For $m=0$, the character is zero. Since $\SO(M)$ fixes $v$, $\psi^0$ extends to an additive character of $H=\SO(M)\ltimes H_u$ by $\psi^0(g,u)=\psi^0(u)$, which is the positive version of the character in \cite[Section 3.2]{[BFT]}.

Let $\check\varepsilon^{(2)}_j$, $1\leq j\leq\lfloor N/2\rfloor$, be the basis of $X^*(T_2)$ dual to the coordinate coweights $e_j^{(2)}$ in \eqref{eq:osp-duality-coordinates}, so that they are dual to $\delta_j$ if $k\leq l$, and to $\varepsilon_j$ if $k>l$. Note that the weight of $e_j$ is $\check\varepsilon^{(2)}_j$ for $1\leq j\leq m$. We denote by $\check\eta_s$ the weight of $f_s$. Then, we have
\[
\check\eta_s=\check\varepsilon^{(2)}_{m+s},\qquad \check\eta_{n+1-s}=-\check\varepsilon^{(2)}_{m+s} \quad(1\leq s\leq\lfloor n/2\rfloor),
\]
and $\check\eta_{(n+1)/2}=0$ for odd $n$. For $m>0$, let
\[
\check\alpha_j=\check\varepsilon^{(2)}_j-\check\varepsilon^{(2)}_{j+1}, \qquad \check\beta_s=\check\varepsilon^{(2)}_m-\check\eta_s,
\]
the weights of the entries $u_{j,j+1}$ and $u_{m,m+s}$. Note that, for odd $n$, since the middle weight is zero, $\check\beta_{(n+1)/2}=\check\varepsilon^{(2)}_m$ is a short geometric root.

\subsubsection{The twisting} Let $\sigma$ and $T_{w,\sigma}$ be as defined after Proposition~\ref{prop 3.2.7}, for $\GL_M\subset\GL_n$, and let
\[
\mathcal J=\{i\in\{1,\ldots,M\}:(i+n\BZ)\cap w(\sigma)=\varnothing\}.
\]
Then
\[
T_{w,\sigma}=\left\{\diag(z_1,\ldots,z_M)\in T_M: z_i=1\text{ for }i\notin\mathcal J\right\}.
\]
By the proof of Proposition~\ref{prop:adjacent-orthogonal-Iwahori-indices}, $\mathcal J$ is stable under $i\mapsto i^*$, and the reductive quotient of the orthogonal stabilizer is the subtorus of $T_{w,\sigma}$ given by
\[
z_{i^*}=z_i^{-1}\qquad(i\in\mathcal J),
\]
with middle diagonal entry $1$ if $M$ is odd. Let
\[
w=\operatorname{diag}(t^{n_1},\ldots,t^{n_n})w^{\rm fin}, \qquad w^{\rm fin}\in S_n.
\]
Note that the torus $T_{w,\sigma}$ acts on the fiber of the relative determinant line at the orbit representative by the character $\prod_{i\in\mathcal J}z_i^{n_i}$, up to inversion. At an orthogonal flag, the involution $\tau$ of the proof of Proposition~\ref{prop:adjacent-orthogonal-Iwahori-indices} fixes the flag and preserves this line, and hence the character is invariant under $\tau$, which sends $z_i$ to $z_{i^*}^{-1}$. By comparing the exponents, we obtain
\begin{equation}\label{eq:orthogonal-free-pair-exponents}
n_{i^*}=-n_i\qquad(i\in\mathcal J).
\end{equation}
We denote $\zeta=b\bigl((0,e_1^{(2)}),(0,e_1^{(2)})\bigr)$, with $b$ as in Section~\ref{3.2.5 sect}. Note that, by the Weyl invariance of $b$, the same value is obtained for every coordinate cocharacter $e_j^{(2)}$ of $T_2$.

\begin{prop}\label{prop:orthogonal-Iwahori-scalar-relevance}
For the twisting given by the determinant line bundle in Remark~\ref{rem:determinant-twisting}, an orbit indexed by \eqref{eq:orthogonal-Iwahori-all-indices} is relevant if and only if the following hold.
\begin{enumerate}
\item If $m>0$, then
\begin{align}
c_j-c_{j+1}+\epsilon_{w_2}(\check\alpha_j)&\geq1 &&(1\leq j<m),\label{eq:orthogonal-Iwahori-tail-test}\\
c_m+\operatorname{ord}_t\bigl((g_{\mathrm{adj}}^{-1}v)_s\bigr)+ \epsilon_{w_2}(\check\beta_s)&\geq1 &&(1\leq s\leq n),
\label{eq:orthogonal-Iwahori-boundary-test}
\end{align}
where $(g_{\mathrm{adj}}^{-1}v)_s$ is the $s$-th coordinate in the basis $f_1,\ldots,f_n$. There is no condition if $m=0$.
\item For $i\in\mathcal J$ with $i<i^*$, we have
\begin{equation}\label{eq:orthogonal-Iwahori-scalar-relevance}
\zeta^{n_i}=1.
\end{equation}
\end{enumerate}
\end{prop}
\begin{proof}
As in Section~\ref{aff kl}, the orbit is relevant if and only if the stabilizer \eqref{eq:Iwahori-stabilizer-splitting} lies in the kernel of $\psi$, and the character sheaf induced by the twisting on it is trivial. Note that the conjugation by $l$ multiplies the entry $u_{j,j+1}$ by $t^{c_j-c_{j+1}}$ and sends the row $X(u)$ to $t^{c_m}X(u)g_{\mathrm{adj}}^{-1}$, and hence, by \eqref{eq:Iwahori-general-additive-test}, the first condition is equivalent to (1).

For the second condition, we may pass to the reductive quotient of the stabilizer, which, by Proposition~\ref{prop:adjacent-orthogonal-Iwahori-indices}, is the product of the tori $(z_i,z_{i^*})=(s,s^{-1})$, for $i\in\mathcal J$ with $i<i^*$. Since the bilinear form corresponding to $\cP_{\det}$ is half the trace (Remark~\ref{rem:determinant-twisting}), such a torus acts on the fiber of $\cP_{\det}$ at the orbit representative by $s\mapsto s^{n_i}$, by \eqref{eq:orthogonal-free-pair-exponents}, and hence the condition is $\zeta^{n_i}=1$.
\end{proof}

In other words, condition (2) says that $\ord(\zeta)$ divides $n_i$ if $\zeta$ has finite order, and $n_i=0$ if $\zeta$ has infinite order.

\begin{rem}\label{rem:orthogonal-Iwahori-general-twisting}
For a general factorization gerbe, condition (2) of Proposition~\ref{prop:orthogonal-Iwahori-scalar-relevance} is replaced by the condition that the character sheaf induced by the gerbe on the stabilizer \eqref{eq:Iwahori-stabilizer-splitting} is trivial. The additive condition (1) is the same.
\end{rem}

\subsection{The exceptional orbits}
\label{sec:F4-Iwahori-orbits}

In this subsection, we classify the orbits in the $\on{F}(4)$ case and give their representatives.

\subsubsection{The Levi subgroup} In the $\on{F}(4)$ case, we have $G_1=\PGL_2$ and $G_2=\PSp(6)$, and the Levi subgroup is
\begin{equation}\label{eq:F4-Iwahori-Levi}
L=(\GL_2\times\GL_2)/\BG_m, \qquad z\longmapsto(z\Id_2,z\Id_2).
\end{equation}
Note that this is the Levi subgroup $(\GL_2\times\on{SL}_2)/\mu_2$ of \cite[Section 3.1]{[FTY]}. For $A,B\in\GL_2(\bF)$, we denote by $l(A,B)\in L(\bF)$ the element represented by $(A,B)$. Under the embedding $L\subset G_2$, a matrix representative of $l(A,B)$ is
\begin{equation}\label{eq:F4-Iwahori-Levi-matrix}
g(A,B)=\operatorname{diag} \bigl(A,B,\det(B)JA^{-t}J\bigr), \qquad J=\begin{pmatrix}0&1\\1&0\end{pmatrix},
\end{equation}
where the symplectic form on $\bF^6$ is given by the matrix
\[
\Omega=
\begin{pmatrix}
0&0&0&0&0&1\\
0&0&0&0&1&0\\
0&0&0&1&0&0\\
0&0&-1&0&0&0\\
0&-1&0&0&0&0\\
-1&0&0&0&0&0
\end{pmatrix}.
\]
Note that $g(A,B)^t\Omega g(A,B)=\det(B)\Omega$, and multiplying $A,B$ by the same scalar does not change the resulting element of $L$. The embedding of $G_1$ is given by $g\mapsto l(g,g)$.

We denote by $\check\varepsilon^{(2)}_i$, $1\leq i\leq3$, the coordinate characters of the diagonal torus of $\Sp(6)$, so that, under \eqref{eq:F4-duality-coordinates}, we have $\langle\varepsilon_j,\check\varepsilon^{(2)}_i\rangle=\delta_{ij}$. Since the positive roots of $L$ are $\check\varepsilon^{(2)}_1-\check\varepsilon^{(2)}_2$ and $2\check\varepsilon^{(2)}_3$, and $I_L$ is upper triangular in both $\GL_2$ factors, the set ${}^LW_2$ of Section~\ref{sec:Iwahori-Levi-reduction} is
\begin{equation}\label{eq:F4-Iwahori-finite-Weyl}
{}^LW_2=\{w_2:\ w_2^{-1}(\check\varepsilon^{(2)}_1-\check\varepsilon^{(2)}_2)>0, \quad w_2^{-1}(2\check\varepsilon^{(2)}_3)>0\}.
\end{equation}
Note that it has $|W_2|/|W_L|=48/4=12$ elements, which are listed in Section~\ref{sec:F4-Iwahori-relevance}.

\subsubsection{The orbit indices and representatives} In this subsubsection, we describe the $I_1$-orbits in $L(\bF)/I_L$ by pairs of oriented edges in the Bruhat--Tits tree of $\PGL_2$. Recall that its vertices are the $\bO$-lattices in $\bF^2$ up to multiplication by $\bF^\times$, and two vertices are adjacent if they have representatives with one contained in the other with one-dimensional quotient; see \cite[Chapter II]{[Serre]}. By an oriented edge, we mean an ordered pair $(x,y)$ of adjacent vertices, and we do not require that the lattice representing $x$ be contained in the one representing $y$. The group $\PGL_2(\bF)$ acts on the tree by $g[\mathcal L]=[g\mathcal L]$, where $[\mathcal L]$ denotes the vertex represented by $\mathcal L$.

Let $e_1,e_2$ be the standard basis of $\bF^2$. For $a\in\BZ$, let $v_a$ be the vertex represented by $t^a\bO e_1+\bO e_2$. The vertices $v_0$ and $v_{-1}$ are adjacent, and we put $e_0=(v_0,v_{-1})$. Note that the upper triangular Iwahori subgroup $I_1$ is the stabilizer of $e_0$, so $\PGL_2(\bF)/I_1$ is identified with the set of oriented edges. For an oriented edge $e=(x,y)$, let $\partial_0(e)=x$. For a vertex $v=[g\bO^2]$, with $g\in\GL_2(\bF)$, let
\[
\operatorname{par}(v)=\operatorname{ord}_t\det g\pmod2.
\]
This is well-defined, since changing the representative changes $\operatorname{ord}_t\det g$ by an even integer. For example, $\operatorname{par}(v_a)\equiv a\pmod2$.

For a point $l(A,B)I_L/I_L\in L(\bF)/I_L$, let
\[
e_A=Ae_0,\qquad e_B=Be_0,\qquad c=\operatorname{ord}_t\det A-\operatorname{ord}_t\det B.
\]
Since $\partial_0(e_A)=[A\bO^2]$ and $\partial_0(e_B)=[B\bO^2]$, we have
\begin{equation}\label{eq:F4-Iwahori-central-parity}
c\equiv\operatorname{par}(\partial_0(e_A))-\operatorname{par}(\partial_0(e_B))\pmod2.
\end{equation}
Conversely, we fix two oriented edges $e_A,e_B$, and choose $A,B$ with $Ae_0=e_A$ and $Be_0=e_B$. Then, the points of $L(\bF)/I_L$ giving these two edges are represented by $(t^rA,B)$, $r\in\BZ$. Indeed, multiplying both matrices by the same scalar gives the same element of $L$, and multiplying either matrix by a scalar in $\bO^\times$ does not change the coset modulo $I_L$. Note that the corresponding integers are $\operatorname{ord}_t\det A-\operatorname{ord}_t\det B+2r$. Thus, the points of $L(\bF)/I_L$ are indexed by the triples $(e_A,e_B,c)$ satisfying \eqref{eq:F4-Iwahori-central-parity}, and $g\in I_1$ acts by
\[
g\cdot(e_A,e_B,c)=(ge_A,ge_B,c),
\]
since multiplying $A$ and $B$ on the left by the same matrix does not change $c$. Therefore, it remains to classify the pairs $(e_A,e_B)$ under $I_1$. There are two cases, according to the smallest subtree containing $e_0,e_A,e_B$, as illustrated in the following picture.

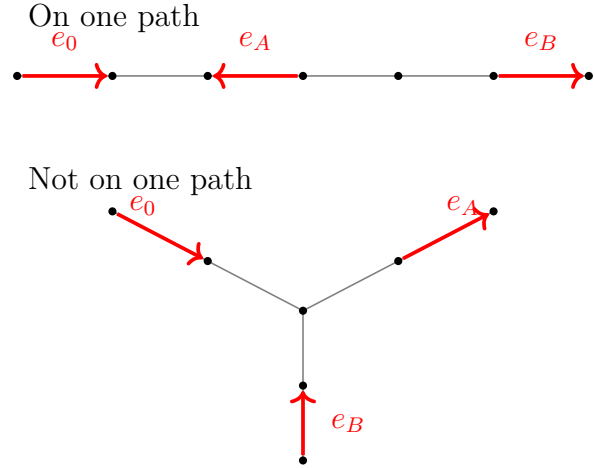
\begin{center}
\begin{tikzpicture}[
scale=0.9, selected/.style={ red, line width=1.4pt, ->, shorten >=2pt, shorten <=2pt }, ordinary/.style={gray, line width=0.6pt} ] \node[anchor=west] at (0,0.85) {On one path}; \draw[ordinary] (0,0)--(8.4,0); \draw[selected] (0,0)--(1.4,0) node[midway,above=5pt] {$e_0$}; \draw[selected] (4.2,0)--(2.8,0) node[midway,above=5pt] {$e_A$}; \draw[selected] (7,0)--(8.4,0) node[midway,above=5pt] {$e_B$}; \foreach \x in {0,1.4,2.8,4.2,5.6,7,8.4} \fill (\x,0) circle (1.7pt);

\node[anchor=west] at (0,-1.55) {Not on one path}; \coordinate (o)  at (4.2,-3.45); \coordinate (a1) at (2.8,-2.72); \coordinate (a2) at (1.4,-1.99); \coordinate (b1) at (5.6,-2.72); \coordinate (b2) at (7,-1.99); \coordinate (c1) at (4.2,-4.55); \coordinate (c2) at (4.2,-5.65);

\draw[ordinary] (a2)--(o)--(b2); \draw[ordinary] (o)--(c2); \draw[selected] (a2)--(a1); \draw[selected] (b1)--(b2); \draw[selected] (c2)--(c1) node[midway,right=6pt] {$e_B$}; \node[red] at (1.85,-1.9) {$e_0$}; \node[red] at (6.55,-1.9) {$e_A$};

\foreach \p in {o,a1,a2,b1,b2,c1,c2} \fill (\p) circle (1.7pt);
\end{tikzpicture}
\end{center}

\subsubsection{} I. First, suppose that $e_0,e_A,e_B$ lie on a simple path, i.e., a path with no repeated vertices, where we ignore the orientations. After applying an element of $I_1$, we may assume that they lie on the path with vertices $v_a$, $a\in\BZ$; see the proof of Proposition~\ref{prop:F4-Iwahori-orbit-indices}. We write
\[
e_A=(v_{r_A},v_{r_A+\sigma_A}),\qquad e_B=(v_{r_B},v_{r_B+\sigma_B}).
\]
These positions and orientations relative to $e_0$ determine the $I_1$-orbit of $(e_A,e_B)$, and hence the $I_1$-orbits in this case are indexed by
\begin{equation}\label{eq:F4-Iwahori-path-indices}
r_A,r_B,c\in\BZ,\quad \sigma_A,\sigma_B\in\{+1,-1\}, \qquad c-r_A+r_B\in2\BZ.
\end{equation}
Let
\[
s=\begin{pmatrix}0&1\\1&0\end{pmatrix},\qquad n_{+1}=s,\quad n_{-1}=\Id_2,\quad h=(c-r_A+r_B)/2.
\]
Then, the representative is given by
\begin{equation}\label{eq:F4-Iwahori-path-matrices}
A=t^h\begin{pmatrix}t^{r_A}&0\\0&1\end{pmatrix}n_{\sigma_A}, \qquad B=\begin{pmatrix}t^{r_B}&0\\0&1\end{pmatrix}n_{\sigma_B}.
\end{equation}

\subsubsection{} II. Now, suppose that the three edges do not lie on one path. Then, the smallest subtree containing them consists of three paths meeting at one vertex, with $e_0,e_A,e_B$ as their respective last edges. For $i=0,A,B$, we denote by $\ell_i$ the length of the path from the common vertex to the far endpoint of $e_i$, and we set $\sigma_i=1$ if $e_i$ is oriented away from the common vertex, and $\sigma_i=-1$ otherwise. Then, $a_i=\ell_i-(\sigma_i+1)/2$ is the distance from the common vertex to $\partial_0(e_i)$, and the indices satisfy
\begin{equation}\label{eq:F4-Iwahori-tripod-indices}
\ell_0,\ell_A,\ell_B\geq1,\qquad \sigma_0,\sigma_A,\sigma_B\in\{+1,-1\},\qquad c\in\BZ,\quad c-a_A+a_B\in2\BZ.
\end{equation}
Let $h=(c-a_A+a_B)/2$, let $P_0=\Id_2$, $P_A=s$, $P_B=\begin{psmallmatrix}1&1\\0&1\end{psmallmatrix}$, and let $M_i=P_i\diag(t^{a_i},1)n_{\sigma_i}$. Then, the representative is given by
\begin{equation}\label{eq:F4-Iwahori-tripod-matrices}
A=t^h M_0^{-1}M_A,\qquad B=M_0^{-1}M_B.
\end{equation}
Indeed, the three paths $P_iv_0,P_iv_1,P_iv_2,\ldots$ meet only at $v_0$, the edge $M_ie_0$ lies on the $i$-th path with the indices $(\ell_i,\sigma_i)$, and $M_0^{-1}$ sends $M_0e_0$ to $e_0$.

\subsubsection{}
\begin{prop}\label{prop:F4-Iwahori-orbit-indices}
The $H(\bF)$-orbits in $\Fl_{\PGL_2}\times\Fl_{\PSp(6)}$ are indexed by the following two sets of parameters:
\begin{enumerate}
\item[I.] $(w_2,r_A,r_B,\sigma_A,\sigma_B,c)$, where $w_2\in{}^LW_2$ and $(r_A,r_B,\sigma_A,\sigma_B,c)$ satisfies \eqref{eq:F4-Iwahori-path-indices}; these correspond to the case where $e_0,e_A,e_B$ lie on one simple path;
\item[II.] $(w_2,\ell_0,\ell_A,\ell_B,\sigma_0,\sigma_A,\sigma_B,c)$, where $w_2\in{}^LW_2$ and $(\ell_0,\ell_A,\ell_B,\sigma_0,\sigma_A,\sigma_B,c)$ satisfies \eqref{eq:F4-Iwahori-tripod-indices}; these correspond to the case where the three edges do not lie on one simple path.
\end{enumerate}
The representative is $(t^0,l(A,B)w_2I_2/I_2)$, where $A,B$ are given by \eqref{eq:F4-Iwahori-path-matrices} in case~I and by \eqref{eq:F4-Iwahori-tripod-matrices} in case~II. Every stabilizer is connected.
\end{prop}

\begin{proof}
By \eqref{eq:Iwahori-Levi-indices} and the discussion above, it suffices to classify the $I_1$-orbits of pairs $(e_A,e_B)$. In case~I, we choose a path containing $e_0,e_A,e_B$ and infinite in both directions, with its two ends ordered so that $e_0$ points towards the first one. In case~II, we choose three infinite paths starting at the common vertex and containing $e_0,e_A,e_B$, respectively. Recall that the infinite simple paths starting at a fixed vertex correspond to the points of $\BP^1(\bF)$, see \cite[Chapter II, \S 1.1, p.~72]{[Serre]}. Thus, the two ends of the path in case~I give two distinct points of $\BP^1(\bF)$, and the three paths in case~II give three distinct points. Now, let $(e_A,e_B)$ and $(e_A',e_B')$ have the same indices. Since $\PGL_2(\bF)$ acts transitively on ordered pairs and triples of distinct points of $\BP^1(\bF)$, there exists $g\in\PGL_2(\bF)$ mapping the chosen ends for the first configuration to the corresponding ends for the second. In case~II, $g$ then maps the common vertex and the three paths to the corresponding ones, so it fixes $e_0$, i.e., $g\in I_1$, and it maps $(e_A,e_B)$ to $(e_A',e_B')$. In case~I, the same holds after composing $g$ with an element of the stabilizer of the two ends, which acts on the path by translations. Conversely, since $I_1$ fixes $e_0$ and preserves the two components of the tree with the interior of $e_0$ removed, it preserves the indices. The representatives are given by \eqref{eq:F4-Iwahori-path-matrices} and \eqref{eq:F4-Iwahori-tripod-matrices}.

Next, we compute the stabilizers. The stabilizer of $l(A,B)I_L/I_L$ in $I_1$ is
\begin{equation}\label{eq:F4-Iwahori-Levi-stabilizer}
R_l=I_1\cap A I_1A^{-1}\cap B I_1B^{-1},
\end{equation}
since the action of $I_1$ preserves $c$. In case~I, we denote by $v_{r_-},v_{r_+}$ the endpoints of the smallest path containing the three edges, with $r_-\leq-1<0\leq r_+$. Then, we have
\[
R_l\simeq \left\{
\begin{pmatrix}u&x\\y&v\end{pmatrix}:
u,v\in\bO^\times,\quad x\in t^{r_+}\bO,\quad y\in t^{-r_-}\bO \right\}/\bO^\times,
\]
where $\bO^\times$ acts by scalar multiplication. Note that the homomorphism
\[
R_l\longrightarrow\BG_m,\qquad \left[\begin{pmatrix}u&x\\y&v\end{pmatrix}\right] \longmapsto \frac{u(0)}{v(0)},
\]
is surjective, and its kernel is pro-unipotent. In case~II, we denote by $\mathcal L\subset\bF^2$ a lattice representing the common vertex. Then, the reduction modulo $t$ gives a homomorphism
\[
\Stab_{\PGL_2(\bF)}([\mathcal L]) \longrightarrow \PGL(\mathcal L/t\mathcal L)\simeq\PGL_2(\BC)
\]
with pro-unipotent kernel. Since the first edges of the three paths starting at the common vertex correspond to three distinct points of $\BP(\mathcal L/t\mathcal L)$, see \cite[Chapter II, \S 1.1, p.~72]{[Serre]}, and every element of $R_l$ fixes these three points, its image is the identity, and hence $R_l$ is pro-unipotent. Thus, in both cases, $R_l$ is connected, and hence so is the stabilizer \eqref{eq:Iwahori-stabilizer-splitting}.
\end{proof}

\subsection{The exceptional relevance conditions}
\label{sec:F4-Iwahori-relevance}

In this subsection, we determine which $\on{F}(4)$ orbits are relevant.

\subsubsection{The additive character} For $u\in H_u$, let $X=(u_{i,j+2})_{1\leq i,j\leq2}$. Then $\psi^0(u)=u_{1,3}+u_{2,4}=\operatorname{tr}(X)$ and $\psi(u)=[t^0](u_{1,3}+u_{2,4})$. Note that the conjugation by $l(A,B)$ sends $X$ to $AXB^{-1}$, and hence, for $C=B^{-1}A$, we have
\begin{equation}\label{eq:F4-Iwahori-character-coefficients}
\psi(\operatorname{Ad}_{l(A,B)}u)=[t^0]\operatorname{tr}(CX).
\end{equation}
Also, the geometric roots corresponding to the entries $X_{ij}$ are
\[
\begin{aligned}
\check\beta_{11}&=\check\varepsilon^{(2)}_1-\check\varepsilon^{(2)}_3, &\check\beta_{12}&=\check\varepsilon^{(2)}_1+\check\varepsilon^{(2)}_3,\\
\check\beta_{21}&=\check\varepsilon^{(2)}_2-\check\varepsilon^{(2)}_3, &\check\beta_{22}&=\check\varepsilon^{(2)}_2+\check\varepsilon^{(2)}_3.
\end{aligned}
\]
Then, for the representatives \eqref{eq:F4-Iwahori-path-matrices} and \eqref{eq:F4-Iwahori-tripod-matrices}, we have
\begin{align}
C&=t^h n_{\sigma_B}^{-1}
\begin{pmatrix}t^{r_A-r_B}&0\\0&1\end{pmatrix}n_{\sigma_A}
&&\text{in case I},
\label{eq:F4-Iwahori-path-coefficients}\\
C&=t^h n_{\sigma_B}^{-1}
\begin{pmatrix}-t^{a_A-a_B}&t^{-a_B}\\t^{a_A}&0\end{pmatrix}n_{\sigma_A}
&&\text{in case II}.
\label{eq:F4-Iwahori-tripod-coefficients}
\end{align}

We denote by $(a_1,a_2,a_3)$ the element $w_2$ satisfying $w_2(\check\varepsilon^{(2)}_i)=\operatorname{sgn}(a_i) \check\varepsilon^{(2)}_{|a_i|}$. The twelve elements of ${}^LW_2$ and the values $\epsilon_{w_2}(\check\beta_{ij})$ of \eqref{eq:Iwahori-root-level} are as follows:
\[
\bigl(\epsilon_{w_2}(\check\beta_{11}),\epsilon_{w_2}(\check\beta_{12}),\epsilon_{w_2}(\check\beta_{21}),\epsilon_{w_2}(\check\beta_{22})\bigr)=
\begin{cases}
(0,0,0,0),&w_2=(1,2,3),\\
(1,1,1,1),&w_2=(-2,-1,3),\\
(0,0,1,1),&w_2=(1,-2,3),\ (-2,1,3),\\
(0,0,1,0),&w_2=(1,3,\pm2),\\
(1,0,1,1),&w_2=(-2,3,\pm1),\\
(1,0,1,0),&w_2=(3,1,\pm2),\ (3,-2,\pm1).
\end{cases}
\]

\subsubsection{The twisting} Note that, in case~II, the stabilizer is pro-unipotent, so there is no twisting condition. In case~I, in the coordinates of \eqref{eq:F4-duality-coordinates}, the cocharacter of $l(A,B)$ is
\begin{equation}\label{eq:F4-Iwahori-translation}
\mu^{(2)}= \left(h+r_A-\frac{r_B}{2},\ h-\frac{r_B}{2},\ \frac{r_B}{2}\right) \in X_*(T_{\PSp(6)}),
\end{equation}
and the reductive quotient $\BG_m$ of $R_l$ in the proof of Proposition~\ref{prop:F4-Iwahori-orbit-indices} is the image of the cocharacter $z\mapsto\diag(z,1)$ of $G_1$. Under the diagonal embedding of $G_1$, this cocharacter corresponds to
\begin{equation}\label{eq:F4-Iwahori-primitive-torus}
\nu=(\nu^{(1)},\nu^{(2)}) =\left(\frac12,\left(\frac12,-\frac12,\frac12\right)\right).
\end{equation}

\begin{prop}\label{prop:F4-Iwahori-relevance}
An orbit of Proposition~\ref{prop:F4-Iwahori-orbit-indices} is relevant if and only if
\begin{enumerate}
\item for every $1\leq i,j\leq2$, we have
\begin{equation}\label{eq:F4-Iwahori-four-tests}
\operatorname{ord}_t C_{ji} +\epsilon_{w_2}(\check\beta_{ij})\geq1;
\end{equation}
\item in case~I, we have
\begin{equation}\label{eq:F4-Iwahori-twisting-test}
b\bigl((0,\mu^{(2)}),\nu\bigr)=1.
\end{equation}
\end{enumerate}
\end{prop}
\begin{proof}
As in the proof of Proposition~\ref{prop:orthogonal-Iwahori-scalar-relevance}, we check the additive character and the twisting on the stabilizer \eqref{eq:Iwahori-stabilizer-splitting}. By \eqref{eq:F4-Iwahori-character-coefficients}, $d\psi^0(\operatorname{Ad}_{l(A,B)}X_{\check\beta_{ij}})$ is a non-zero constant multiple of $C_{ji}$, so condition \eqref{eq:Iwahori-general-additive-test} is equivalent to (1). In case~II, $R_l$ is pro-unipotent, so there is no further condition. In case~I, the reductive quotient of the stabilizer is the one-dimensional torus given by $\nu$, and the point $(t^0,l(A,B)w_2I_2/I_2)$ has translation part $(0,\mu^{(2)})$, and hence, as in Proposition~\ref{prop:relevant-L-orbits}, the twisting on this torus is given by $b((0,\mu^{(2)}),\nu)$, which gives (2).
\end{proof}

\begin{cor}\label{cor:orthogonal-exceptional-Iwahori-IC}
In the orthosymplectic case, the irreducible objects in the perverse heart of the Iwahori Gaiotto category, up to parity shift, are indexed by the relevant orbits of Proposition~\ref{prop:orthogonal-Iwahori-scalar-relevance}. In the $\on{F}(4)$ case, they are indexed by the relevant orbits of Proposition~\ref{prop:F4-Iwahori-relevance}. For each such index, the object is the IC extension of the irreducible equivariant twisted local system with even coefficients on its orbit, with the normalizations used in Section~\ref{aff kl}.
\end{cor}

\begin{rem}\label{rem:Iwahori-isogeny-components}
For the isogeny \eqref{eq:gaiotto-isogeny}, let $\widetilde I$ be the Iwahori subgroup of $\widetilde G(\bO)$ with $\widetilde I\cap G(\bF)=I$. We choose $\widetilde\eta\in\widetilde\Lambda$ over $\xi\in\Lambda_0$ as in Appendix~\ref{sec:gaiotto-isogeny-proofs}, namely, central if possible, and as in Section~\ref{sec:app-isogeny-projections} otherwise. We denote $I_{\widetilde\eta}=t^{\widetilde\eta}It^{-\widetilde\eta}$, so that the right multiplication by $t^{\widetilde\eta}$ gives an $H(\bF)$-equivariant isomorphism
\[
G(\bF)/I_{\widetilde\eta} \xrightarrow{\sim} \Fl_{\widetilde G}\times_{\Gr_{T_0}}\{t^\xi\}, \qquad gI_{\widetilde\eta}\longmapsto gt^{\widetilde\eta}\widetilde I.
\]
If $\widetilde\eta$ is central, then $I_{\widetilde\eta}=I$, and the representatives above are multiplied on the right by $t^{\widetilde\eta}$. In the other orthogonal cases, we use the same calculations with $I$ replaced by $I_{\widetilde\eta}$. For a geometric root $\check\beta$, we have
\[
I_{\widetilde\eta}\cap U_{\check\beta}(\bF)=
\begin{cases}
U_{\check\beta}(t^{\langle\widetilde\eta,\check\beta\rangle}\bO),&\check\beta>0,\\
U_{\check\beta}(t^{1+\langle\widetilde\eta,\check\beta\rangle}\bO),&\check\beta<0,
\end{cases}
\]
where positivity is taken with respect to the Borel subgroup defining $I$. Then, we check the restriction of $\psi$ to each stabilizer, and calculate the twisting at $gt^{\widetilde\eta}$, with $\widetilde b$ in place of $b$.
\end{rem}

\section*{Acknowledgments} We would like to express our sincere gratitude to Michael Finkelberg and Alexander Braverman for their contributions to this paper, their helpful discussions, and their generous assistance. We are especially grateful to Michael Finkelberg for his constant support throughout the project and for his comments, which helped improve the manuscript.

We also thank Alexander Braverman for bringing D. Gaiotto's conjecture to our attention.

We worked on this paper during 2022--2026. During this time, we got many ideas and learned a lot related to this project from discussions and email correspondence with Lin Chen, Michael Finkelberg, Alexander Braverman, Gurbir Dhillon, Ivan Losev, Sam Raskin, Dennis Gaitsgory, Lingfei Yi, Tsao-Hsien Chen, Weiqiang Wang, Hiroyuki Yamane, Jinkui Wan, Jialiang Zou, Zhiwei Yun, Daxin Xu, Zicheng Qian, Lei Zhang, Jize Yu, Toan Pham, Quan Situ, and Valentin Buciumas.

This paper is part of a long series of papers with Alexander Braverman, Victor Ginzburg and Michael Finkelberg on mirabolic Hecke algebras and, later, their categorification. We are grateful to them for sharing many ideas. Finally, we thank Dennis Gaitsgory for teaching us many techniques in geometric Langlands theory.

\section*{Disclosure of AI use}
We used ChatGPT for the following purposes:
\begin{itemize}
    \item polishing the exposition and making the notation consistent. Due to the fact that different parts of the work were finished at different times, we used several different systems of notations and conventions in the first draft written by the authors and we unified the notation with AI's help;
    \item proofreading. ChatGPT identified many minor mistakes and the following three errors in the first version written by the authors:
    \begin{itemize}
        \item In the first version of Section~\ref{sec:g1-direction}, we claimed that $K_{\lambda^{(2)}}^{(1)}\subset U_1(\bF)T_1(\bO)$ and calculated only the $\GL(M|N)$ and $\on F(4)$ cases. ChatGPT pointed out that this inclusion is false in the orthosymplectic case for the choice of $\sD$ in Section~\ref{trans D}. The argument is corrected directly by conjugating $K_{\lambda^{(2)}}^{(1)}$ by the element $h$ in Lemma~\ref{lem:orthogonal-auxiliary-Borel}.
        \item In the original proof of the block decomposition in Lemma~\ref{prop 9.2.11}, we made a basic mistake by claiming that all relevant orbits lie in $P\times\BA^L$, where $P\subset\GL_L$ is the corresponding parabolic subgroup. In the revised proof, we first permute the right monodromy characters and then restrict to an appropriate open subset.
        \item In the original proof of pointwise purity, our argument for finite typeness in Section \ref{14.2} was incorrect. In the revised version, we used ChatGPT to draft a new Section~\ref{14.2} by suggesting the strategy that there should be a general lemma about the existence of transversal slices for stacks, as stated in Lemma \ref{lem:finite-type-transversal}, and that, for any global model, there should be a Hecke correspondence diagram satisfying the conditions of that lemma, i.e., Proposition \ref{prop:general-H-Hecke-correspondence}. We checked and revised Section \ref{14.2}.
    \end{itemize}
    \item completing the rank-one calculations for $\SOSp(2k|2l)$ in Sections \ref{10.3.3} and \ref{10.4.5}: in the first version, we calculated the cases of $\GL(M|N),\on{F}(4)$ and $\SOSp(2k|2l)$ when $k\leq l$ in Sections \ref{10.3.3}; and we calculated the cases of $\GL(M|N),\on{F}(4)$ in Section \ref{10.4.5}. We stated that the other case was similar.
    \item drafting the appendices. For Appendix~\ref{sec:gaiotto-isogeny-proofs}, we asked AI to prove analogues of Theorems~\ref{thm:image-of-ic} and \ref{thm:metaplectic-hecke-action} for $\widetilde G$ by the same method, which we believed to work but did not write down in the first version. For Appendix~\ref{sec:orthogonal-exceptional-Iwahori}, we asked AI to classify the relevant orbits in the orthosymplectic and $\on{F}(4)$ cases for completeness (we do not need the classification in the proof of the main theorems). We checked and revised the Appendices.
    \item reorganizing the structure of Section \ref{sec:irreducibles}. In the first version, the exposition was not linear, and we asked AI to reorganize Sections \ref{10.1}--\ref{sec:g1-direction}.
    \item generating some of the figures in this paper.
\end{itemize}
We also used Claude to polish the exposition, to proofread the manuscript, and to reorganize Sections~\ref{sec:factorization-identification} and \ref{sec:irreducibles}. 

The authors are responsible for the content of this paper.

\end{document}